\chardef\bslash=`\\

\newcommand{\phibar}{\overline{\phi}}
\def\Rone{R_1}
\def\Stwo{S_2}
\def\Ob{\mathrm{Ob}}

\def\Cous{\mathrm{Cous}}
\def\Kos{\mathrm{Kos}}
 \def\wm{\widetilde{\mathfrak{m}}}
 \def\RRR{\widetilde{R}}

\def\met{m\text{-}et}
\def\pad{p\text{-}ad}
\def\Parmet{\Par\text{-}m\text{-}et}

\documentclass{amsart}
 
\usepackage[bookmarksopen,bookmarksdepth=2]{hyperref}
\usepackage{etex} 
\usepackage{phoenician}
\usepackage{linearb}
\usepackage{linearA}
\usepackage{scalerel}
\usepackage{graphicx}

\usepackage{diagrams}
\newarrow{Equals}{=}{=}{=}{=}{=}

\usepackage[utf8]{inputenc}
\usepackage[T1]{fontenc}
\usepackage{fixltx2e}
\usepackage{mathrsfs}

\newcommand{\bY}{\mathbf{Y}}
\newcommand{\Kot}{\operatorname{Kott}}
\def\Km{\mathrm{Km}}
\def\NS{\mathrm{NS}}
\def\AFF{AF\kern-0.12em{F}}

\usepackage[OT2,T1]{fontenc}
\DeclareSymbolFont{cyrletters}{OT2}{wncyr}{m}{n}
\DeclareMathSymbol{\Sha}{\mathalpha}{cyrletters}{"58}

\newcommand{\whalingship}{\C}
\def\wtheta{\widetilde{\theta}}
\def\rrr{\widetilde{r}}

\def\HH{H}
\def\comp{\mathrm{comp}}
\def\fflat{\mathrm{flat}}
\renewcommand{\star}{*}
\def\RGamma{\mathrm{R}\Gamma}
\def\varrhobar{\overline{\varrho}}
\def\Asai{\mathrm{As}}
\def\PSp{\mathrm{PSp}}
\def\PGSp{\mathrm{PGSp}}
\def\wA{\widetilde{A}}
\newcommand{\oscr}{\mathcal{O}}
\newcommand{\Ha}{\mathrm{Ha}}
\newcommand{\an}{\textrm{an}}
\newcommand{\tm}{\widetilde{\m}}
\newcommand{\tT}{\widetilde{T}}

\newcommand{\ext}{\operatorname{ext}}
\def\SheafHom{\mathop{\mathcal{H}\!{\it om}}\nolimits}

\def\XXX{\mathcal{X}}
\def\doubleslash{/\kern-0.2em{/}}

\def\normU{U}
\def\Unaive{U^{\mathrm{naive}}}

\def\tK{\widetilde{K}}
\def\fX{\mathfrak{X}}
\def\fY{\mathfrak{Y}}
\def\fC{\mathfrak{C}}
\def\sgn{\mathrm{sgn}}
\def\MM{\widetilde{M}}

\def\gk{\mathfrak{k}}

\def\ppp{\mathfrak{p}}

\def\af{a \kern-0.05em{f}}
\def\PPP{\wp}
\def\HomFR{\Hom(F,\R)}
\newcommand{\tLambda}{\widetilde{\Lambda}}
\newcommand{\alphat}{\tilde{\alpha}}
\newcommand{\betat}{\tilde{\beta}}
\newcommand{\etat}{\tilde{\eta}}

\def\can{\mathrm{can}}

\def\cusp{\mathrm{cusp}}

\def\OL{\mathcal{O}}

\def\bbb{\mathfrak{b}}

\newcommand{\Jac}{\operatorname{Jac}}
\newcommand{\Favoid}{F^{(\mathrm{avoid})}} 
\newcommand{\Iw}{\mathrm{Iw}}
\newcommand{\Sph}{\mathrm{Sph}}
\newcommand{\Kli}{\mathrm{Kli}}
\newcommand{\Si}{\mathrm{Si}}
\newcommand{\Par}{\mathrm{Par}}

\newcommand{\Spa}{\mathrm{Spa}}

\newcommand{\fee}[1]{
\begin{tikzpicture}[#1]
\draw (0,1ex) -- (0.9ex,1ex);
\draw (0.9ex,1ex) -- (0.9ex,0);
\draw (0.9ex,0) -- (0,1ex);
\end{tikzpicture}}
\def\tri{\fee{scale=1.1}}

\renewcommand{\mathbb}{\mathbf}

\newcommand{\alphabeta}{\scaleobj{0.7}{\textlinb{\LinearAC}}}
\newcommand{\alphabar}{\overline{\alpha}}
\newcommand{\betabar}{\overline{\beta}}
\newcommand{\alphabetabar}{\overline{\alphabeta}}
\newcommand{\recGT}{\rec_{\operatorname{GT}}}
\newcommand{\recGTp}{\rec_{\operatorname{GT},p}}
\newcommand{\Tor}{\mathrm{Tor}}

\newcommand{\psibar}{{\overline{\psi}}}

\newcommand{\chibar}{\overline{\chi}}
\newcommand{\rbar}{\overline{r}}

\newcommand{\sbar}{\overline{s}}

\newcommand{\zbar}{\overline{z}}
\newcommand{\gr}{\operatorname{gr}}

\usepackage{tikz}
\usepackage{tikz-cd}
\usetikzlibrary{arrows, matrix}
\usetikzlibrary{arrows.meta,
                positioning}

\usepackage[shortlabels]{enumitem}\usetikzlibrary{decorations.pathmorphing}

\newcommand{\Lie}{{\operatorname{Lie}\,}}
\newcommand{\semis}{{\operatorname{ss}}}
\newcommand{\Fsemis}{{\operatorname{F-ss}}}

\newcommand{\fg}{{\mathfrak{g}}}
\newcommand{\fb}{{\mathfrak{b}}}
\newcommand{\mf}{\mathfrak}
\newcommand{\BC}{{\operatorname{BC}}}
\newcommand{\nInd}{{\operatorname{n-Ind}}}

\newcommand{\ndiv}{{\mbox{$\not| $}}}
\newcommand{\Br}{{\operatorname{Br}}}
\newcommand{\rec}{{\operatorname{rec}}}

\newcommand{\lra}{\longrightarrow}

\newcommand{\CA}{{\mathcal{A}}}
\newcommand{\CO}{{\mathcal{O}}}
\newcommand{\CR}{{\mathcal{R}}}
\newcommand{\CS}{{\mathcal{S}}}
\newcommand{\CX}{{\mathcal{X}}}

\newcommand{\barE}{\overline{{E}}}
\newcommand{\barK}{\overline{{K}}}
\newcommand{\barM}{\overline{{M}}}
\newcommand{\Ebar}{\overline{{E}}}
\newcommand{\barr}{\overline{{r}}}

\def\PSL{\mathrm{PSL}}
\def\PGL{\mathrm{PGL}}

\usepackage{amssymb,amsmath,amsfonts,amsthm,epsfig,amscd,stmaryrd}
\usepackage{stmaryrd}
\renewcommand{\llbracket}{[[}
\renewcommand{\rrbracket}{]]}

\usepackage[all,cmtip,poly]{xy}
\usepackage{color}

\newcommand{\loc}{\operatorname{loc}}
\newcommand{\ad}{\operatorname{ad}}
\newcommand{\diag}{\operatorname{diag}}
\newcommand{\tr}{\operatorname{tr}}
\newcommand{\Tr}{\operatorname{Tr}}

\newcommand{\CNL}{\operatorname{CNL}}

\newcommand{\red}{\operatorname{red}}
\newcommand{\colim}{\varinjlim}

\newcommand{\univ}{{\operatorname{univ}}}

\newcommand{\To}{\longrightarrow}
\newcommand{\isoto}{\stackrel{\sim}{\To}}
\newcommand{\wotimes}{\widehat{\otimes}}

\newcommand{\codim}{\operatorname{codim}}

\newcommand{\gP}{\mathfrak{P}}

\def\iso{\buildrel \sim \over \longrightarrow}

\newcommand{\id}{\operatorname{id}}

\newcommand{\RHom}{{\mathrm{RHom}}}

\newtheorem{theorem}[subsubsection]{Theorem}
\newtheorem{thm}[subsubsection]{Theorem}
\newtheorem{lemma}[subsubsection]{Lemma}
\newtheorem{lem}[subsubsection]{Lemma}
\newtheorem{xlem}[subsubsection]{Expected Lemma}

\newtheorem{cor}[subsubsection]{Corollary}
\newtheorem{conj}[subsubsection]{Conjecture}

\newtheorem{prop}[subsubsection]{Proposition}

\theoremstyle{definition}
\newtheorem{df}[subsubsection]{Definition}
\newtheorem{defn}[subsubsection]{Definition}

\theoremstyle{remark}
\newtheorem{remark}[subsubsection]{Remark}
\newtheorem{rem}[subsubsection]{Remark}

\newtheorem{example}[subsubsection]{Example}

\newtheorem{hypothesis}[subsubsection]{Hypothesis}

\def\numequation{\addtocounter{subsubsection}{1}\begin{equation}}
\def\nummultline{\addtocounter{subsubsection}{1}\begin{multline}}
\def\anumequation{\addtocounter{subsection}{1}\begin{equation}}
\def\anummultline{\addtocounter{subsection}{1}\begin{multline}}

\newif\iffinalrun
\iffinalrun
\else
 \fi

\iffinalrun
  \newcommand{\need}[1]{}
  \newcommand{\mar}[1]{}
\else
  \newcommand{\need}[1]{{\tiny *** #1}}
\newcommand{\mar}[1]{\marginpar{\raggedright\tiny FIXME: #1 }}\fi

\newcommand{\A}{\mathbf{A}}
\newcommand{\C}{\CC}
\newcommand{\F}{\FF}

\newcommand{\Q}{\QQ}
\newcommand{\R}{\RR}
\newcommand{\Z}{\ZZ}

\newcommand{\G}{\cG}
\renewcommand{\L}{\cL}

\renewcommand{\O}{\cO}

\newcommand{\m}{\frakm}

\newcommand{\p}{\frakp}
\newcommand{\q}{\frakq}

\newcommand{\CC}{{\mathbb C}}
\newcommand{\FF}{{\mathbb F}}
\newcommand{\GG}{{\mathbb G}}
\newcommand{\RR}{{\mathbb R}}
\newcommand{\QQ}{{\mathbb Q}}
\newcommand{\TT}{{\mathbb T}}
\newcommand{\TTt}{\widetilde{\TT}}

\newcommand{\TTT}{\widetilde{\TT}}
\newcommand{\ZZ}{{\mathbb Z}}

\newcommand{\ba}{\ensuremath{\mathbf{a}}}

\newcommand{\cC}{{\mathcal C}}
\newcommand{\cD}{{\mathcal D}}

\newcommand{\cF}{{\mathcal F}}
\newcommand{\cG}{{\mathcal G}}
\newcommand{\cH}{{\mathcal H}}
\newcommand{\cI}{{\mathcal I}}

\newcommand{\cL}{{\mathcal L}}
\newcommand{\CL}{{\mathcal{L}}}
\newcommand{\cM}{{\mathcal M}}
\newcommand{\cN}{{\mathcal N}}
\newcommand{\cO}{{\mathcal O}}
\newcommand{\ocal}{{\mathcal O}}
\newcommand{\cP}{{\mathcal P}}

\newcommand{\cR}{{\mathcal R}}
\newcommand{\cS}{{\mathcal S}}
\newcommand{\cT}{{\mathcal T}}
\newcommand{\cU}{{\mathcal U}}

\newcommand{\cX}{{\mathcal X}}
\newcommand{\cY}{{\mathcal Y}}
\newcommand{\cZ}{{\mathcal Z}}
\newcommand{\frakm}{\mathfrak{m}}

\newcommand{\frakp}{\mathfrak{p}}
\newcommand{\frakq}{\mathfrak{q}}

\newcommand{\fra}{{\mathfrak a}}

\newcommand{\Fbar}{\overline{\F}}
\newcommand{\Qbar}{\overline{\Q}}

\newcommand{\Fp}{\F_p}
\newcommand{\Fq}{\F_q}
\newcommand{\Fl}{\F_{l}}

\newcommand{\Fpbar}{\Fbar_p}
\newcommand{\Fqbar}{\Fbar_q}

\newcommand{\Zp}{\Z_p}

\newcommand{\Ql}{\Q_{l}}
\newcommand{\Qp}{\Q_p}

\newcommand{\Qlbar}{\Qbar_{l}}
\newcommand{\Qlbartimes}{\Qlbar^{\times}}
\newcommand{\Qpbar}{\Qbar_p}

\newcommand{\Qpbartimes}{\Qpbar^{\times}}

\DeclareMathOperator{\Aut}{Aut}
\DeclareMathOperator{\coker}{coker}

\DeclareMathOperator{\End}{End}
\DeclareMathOperator{\Ext}{Ext}
\DeclareMathOperator{\Fil}{Fil}
\DeclareMathOperator{\Gal}{Gal}
\DeclareMathOperator{\GL}{GL}
\DeclareMathOperator{\GO}{GO}
\DeclareMathOperator{\GSp}{GSp}
\DeclareMathOperator{\GSpin}{GSpin}

\DeclareMathOperator{\Hom}{Hom}
\DeclareMathOperator{\im}{im}
\DeclareMathOperator{\Ind}{Ind}

\DeclareMathOperator{\Max}{Max}

\DeclareMathOperator{\ord}{ord}
\DeclareMathOperator{\Pic}{Pic}
\DeclareMathOperator{\Proj}{Proj}

\DeclareMathOperator{\SL}{SL}
\DeclareMathOperator{\SO}{SO}
\DeclareMathOperator{\Sp}{Sp}
\DeclareMathOperator{\Spec}{Spec}

\DeclareMathOperator{\Spf}{Spf}

\DeclareMathOperator{\Sym}{Sym}

\DeclareMathOperator{\WD}{WD}
\DeclareMathOperator{\len}{len}

\newcommand{\mult}{\mathrm{mult}}
\newcommand{\ab}{\mathrm{ab}}

\newcommand{\dR}{\mathrm{dR}}

\newcommand{\Frob}{\mathrm{Frob}}

\newcommand{\HT}{\mathrm{HT}}
\newcommand{\Id}{\mathrm{Id}}
\newcommand{\nr}{\mathrm{nr}}

\newcommand{\rhobar}{\overline{\rho}}

                               \newcommand{\into}{\hookrightarrow}
\newcommand{\onto}{\twoheadrightarrow}

\newcommand{\Gm}{\GG_m}

\newcommand{\Art}{{\operatorname{Art}}}

\newcommand{\varepsilonbar}{\overline{\varepsilon}}
\def\vareps{\varepsilon}
\def\varepsbar{\varepsilonbar}
\newcommand{\Res}{\operatorname{Res}}

\title[Potential modularity of abelian surfaces]{Abelian Surfaces over totally real fields are Potentially Modular}

\author[G. Boxer]{George Boxer}  \email{gboxer@math.uchicago.edu} \address{The University of Chicago,
5734 S University Ave,
Chicago, IL 60637, USA}

\author[F. Calegari]{Frank Calegari}  \email{fcale@math.uchicago.edu} \address{The University of Chicago,
5734 S University Ave,
Chicago, IL 60637, USA}

\author[T. Gee]{Toby Gee} \email{toby.gee@imperial.ac.uk} \address{Department of
  Mathematics, Imperial College London,
  London SW7 2AZ,~UK}
  
\author[V. Pilloni]{Vincent Pilloni}\email{vincent.pilloni@ens-lyon.fr}\address{Unit\'e de Math\'ematiques pures et appliqu\'ees,
Ecole normale sup\'erieure de Lyon,
46 all\'ee d'Italie,
69 364 Lyon Cedex 07,
France}
  
\thanks{G.B. was supported in
part by NSF postdoctoral fellowship DMS-1503047.
F.C. was supported in part by NSF Grants DMS-1404620,
  DMS-1701703, and DMS-2001097. T.G.\ was
  supported in part by a Leverhulme Prize, EPSRC grant EP/L025485/1,
  ERC Starting Grant 306326, and a Royal Society Wolfson Research
  Merit Award. V.P was supported in part by the ANR-14-CE25-0002-01 Percolator, and the ERC-2018-COG-818856-HiCoShiVa. }
  
\begin{document}

\begin{abstract}
  We show that abelian surfaces (and consequently curves of genus~2)
  over totally real fields are potentially modular. As a consequence,
  we obtain the expected meromorphic continuation and functional
  equations of their Hasse--Weil zeta functions. We furthermore show
  the modularity of infinitely many abelian surfaces~$A$ over~$\Q$
  with~$\End_{\C}A=\Z$. We also deduce modularity and potential
  modularity results for genus one curves over (not necessarily CM)
  quadratic extensions of totally real fields.
\end{abstract}

\maketitle

\setcounter{tocdepth}{2}
{\footnotesize
\tableofcontents
}

\section{Introduction}
\subsection{Our main theorems}

Let $X$ be a smooth, projective variety of dimension~$m$ over a number field~$F$  with good reduction outside a finite set of primes~$S$.
Associated to~$X$, one may write down a global Hasse--Weil zeta function:
$$\zeta_X(s) = \prod \frac{1}{1 - N(x)^{-s}},$$
where the product runs over all the closed points~$x$ of some (any) smooth proper integral model~$\XXX/\OL_F[1/S]$ for~$X$.
(We suppress~$S$ from the notation --- different choices of~$S$ only change~$\zeta_X(s)$ by a finite number of Euler factors.)
The function~$\zeta_X(s)$ is absolutely convergent for~$\mathrm{Re}(s)
> 1 + m$. We have the following:
\begin{conj}[Hasse--Weil Conjecture, cf.~\cite{Serre1969-1970}, in particular Conj.\ C9] \label{conj:serre} The function~$\zeta_X(s)$ extends to a meromorphic function of~$\C$. There exists a positive real number~$A \in \R^{>0}$,  non-zero rational functions~$P_v(T)$ for~$v|S$, and infinite Gamma factors~$\Gamma_v(s)$ for~$v| \infty$ such that:
$$\xi(s) = \zeta_X(s) \cdot  A^{s/2}  \cdot \prod_{v | \infty} \Gamma_v(s)  \cdot \prod_{v|S} P_v(N(v)^{-s})$$
satisfies the functional equation~$\xi(s) = w \cdot \xi(m+1-s)$ with~$w = \pm 1$.
\end{conj}
(In Serre's formulation of the conjecture, the Gamma factors are also given explicitly in terms of the Archimedean Hodge structures of $X$.)
This conjecture appears to be first formulated in print (albeit in a less precise form and only
for curves) on the final page of~\cite{MR0045416}.
If~$F = \Q$ and~$X$ is a point, then~$\zeta_X(s)$ is  the Riemann zeta function, 
and Conjecture~\ref{conj:serre} follows from Riemann's functional equation~\cite{Riemann}. 
If~$F$ is a general number field but~$X$ is still a point,
then~$\zeta_X(s)$ is the Dedekind zeta function~$\zeta_F(s)$, and Conjecture~\ref{conj:serre} is a theorem of Hecke~\cite{Hecke}.
If~$X$ is a curve of genus zero, then (up to bad Euler factors) $\zeta_X(s) = \zeta_F(s) \zeta_F(s-1)$, and Conjecture~\ref{conj:serre} follows immediately. 
More generally, if~$X$ is any smooth projective variety  whose cohomology is generated by algebraic cycles over~$\overline{F}$, then~$\zeta_X(s)$
is a finite product of Artin~$L$-functions (up to translation), and  Conjecture~\ref{conj:serre} in this case is a consequence of Brauer's theorem~\cite{Brauer}.
In the case when the Galois representations associated to the~$l$-adic
cohomology of~$X$ are \emph{potentially} abelian (e.g.\ an abelian
variety with CM), Conjecture~\ref{conj:serre}
is also a consequence of the results of Hecke and Brauer. 

The fundamental work of Wiles~\cite{MR1333035,MR1333036} and the subsequent
work of Breuil, Conrad, Diamond, and Taylor~\cite{MR1639612,MR1839918} proved Conjecture~\ref{conj:serre} for curves~$X/\Q$ of genus one,
since (again up to a finite number of Euler factors) $\zeta_X(s) = \zeta_{\Q}(s) \zeta_{\Q}(s-1)/L(E,s)$ (where~$E = \Jac(X)$), and the modularity of~$E$ implies the holomorphy and functional
equation for~$L(E,s)$. 
More generally, the potential modularity results of~\cite{MR1954941} imply Conjecture~\ref{conj:serre} for curves~$X/F$ of genus one over any totally real field.
The methods used in these papers have been vastly generalized over the past~25 years
due to the enormous efforts of many people. On the other hand, these
methods have until recently been extremely reliant on the assumption
that the Hodge numbers~$h^{p,q} = \dim H^{p,q}_{\dR}(X) = \dim
H^p(X,\Omega^q)$ of~$X$ are at most~$1$ for all~$p$ and~$q$, or at least that such an inequality holds (suitably interpreted) for the irreducible motives occurring in the cohomology of~$X$. While
many such motives exist inside the cohomology of Shimura varieties, there is a paucity of natural geometric examples satisfying this condition. For example, if~$X$ is a curve of
genus~$g$, then~$h^{1,0}=h^{0,1} = g$, and so the original Taylor--Wiles method only
applies when~$g = 0$ or~$1$. 
For genus two curves, we prove the following theorem.

\begin{theorem}\label{theorem:puppy} Let~$X$ be either a genus two curve or an abelian surface over a totally real field~$F$. Then Conjecture~\ref{conj:serre}
holds for~$X$.
\end{theorem}

We prove Theorem~\ref{theorem:puppy} as a  corollary of the following
theorem. 

\begin{theorem} \label{theorem:main} Let~$X$ be either a genus two curve or an abelian surface over a totally real field~$F$. Then~$X$
is potentially automorphic.
\end{theorem}Here  by \emph{potentially automorphic} we mean that there exists a finite Galois
extension~$L/F$ such that
the compatible system of Galois representations~$\CR$ attached to~$H^1(X_{\Qbar},\Q_p)$
(as~$p$ varies) over~$L$
 is automorphic in a precisely circumscribed sense which we make
 explicit in Definition~\ref{defn: automorphic compatible system general number
    field}. (See also Remark~\ref{rem: relationship between modular
    and automorphic} for a discussion of how we distinguish between
  \emph{automorphic} and~\emph{modular} in this paper; this
  distinction is made purely for technical convenience, and can safely
  be ignored while reading this introduction.) In particular, an immediate consequence is that the~$L$-function
    of~$H^1(X_{\Qbar},\Q_p)$ \emph{as a~$G_{L}$-representation} extends to a holomorphic function on
    all of~$\C$.
 Theorem~\ref{theorem:puppy} follows from
Theorem~\ref{theorem:main} via a standard argument with Brauer's
theorem and base change, together (in the case of abelian surfaces)
with known functorialities in small rank.
(Some care must be taken in this deduction if
the $p$-adic Galois representations associated to~$X$ become reducible
after restriction to~$L$; 
this issue does not arise in the most
interesting cases of Theorem~\ref{theorem:main}, in particular the
case of an abelian surface~$X$ with~$\End_{\C}(X)=\Z$.)

Theorem~\ref{theorem:main} (and thus also Theorem~\ref{theorem:puppy})
is a consequence of Theorem~\ref{thm: potential modularity implies meromorphic continuation} and Corollary~\ref{cor: genus 2 curves},
which in turn are deduced from our main modularity lifting theorem, Theorem~\ref{thm: final R equals T}.
As a consequence of Theorem~\ref{theorem:main}, we also deduce the following potential modularity result for genus one 
curves (see Theorem~\ref{thm: potential modularity of elliptic curves}):

\begin{theorem} \label{thm:genusone} Let~$X$ be a genus one curve over a quadratic extension~$K/F$ of a totally real field~$F$.
Then~$X$ is potentially modular.
\end{theorem}

When~$K/F$ is totally real, this result has been known for some time (\cite{MR1954941}). When~$K/F$ is totally imaginary,
however, the result was only recently proved in~\cite{10author}. For all other quadratic extensions (such as~$F = \Q(\sqrt{2})$ and~$K = \Q(\sqrt[4]{2})$), the result is new.  (See the remarks in~\S\ref{section:comparison} for a comparison between
the methods of this paper with those of~\cite{10author}.)

Just as elliptic
curves over~$\Q$ can be associated (via the modularity theorem) to  modular forms of weight~$2$,
the Langlands program predicts that abelian surfaces over~$\Q$ should
be modular in the sense that they correspond to
certain weight 2 Siegel modular forms. 
This is because (due to the existence of polarizations)
the Galois representations associated to
the  $p$-adic Tate modules of abelian surfaces
are naturally valued in~$\GSp_4(\Qp)$, and~$\GSp_4$ is its own
Langlands dual group.
A consideration of the Hodge--Tate weights 
then suggests
that the corresponding automorphic forms on~$\GSp_4$
should be of weight~$2$ (see \S\ref{subsec: paramodular} for a
more detailed discussion of this).

Our methods also have implications for the \emph{modularity} (as opposed
to potential modularity) of
abelian surfaces over totally real fields.  Here is an example of what
can be proven by our methods.

\begin{thm} \label{thm:bcgpthmintro}  There exist infinitely many modular
  abelian surfaces~$A/\Q$ up to twist with~$\End_{\C}A=\Z$.\end{thm}

As a consequence, one deduces that the~$L$-function
associated to ~$A$ in Theorem~\ref{thm:bcgpthmintro} (that is, the $L$-function associated to the Galois
representation~$H^1(A_{\Qbar},\Q_p)$ for any prime~$p$)
  has a holomorphic continuation to the entire complex plane. Note that Theorem~\ref{theorem:main} only implies
  that this $L$-function has a  meromorphic continuation, with no control over any
  possible poles. (This is for essentially the same reason that
  Brauer's theorem proves the meromorphic continuation of Artin
  $L$-functions, but not the holomorphic continuation.) In fact, we can also prove an analogous theorem for any totally real field~$F$ in which~$3$ splits completely; see
 Theorem~\ref{thm:bcgpthm}. 

To put Theorem~\ref{thm:bcgpthmintro} into
 context, note firstly that if~$\End_{\C}(A) \ne \Z$, then the Galois representations associated to~$A$ become reducible over some finite extension, and hence
one may use (or prove) special cases of functoriality to reduce the problem to the modularity of representations of dimensions~$2$ or~$1$.
Results of this kind appear in the
papers~\cite{YoshidaOriginal,Yoshida,MR2318628,MR2887605,MR3466871,BergerBianchi}. (Several
of these arguments could now be redone more systematically in light of
the monumental work of Arthur~\cite{MR2058604,MR3135650}.)  

In the ``typical'' case that~$\End_{\C}(A) = \Z$, Brumer and
Kramer~\cite{MR3165645} formulated the \emph{paramodular conjecture},
which gives a precise prescription for the ``optimal'' level structure
for an automorphic form corresponding to a given abelian surface; in
particular, this in principle reduces the conjecture for a given~$A$
to an explicit computation of a (finite-dimensional) space of Siegel
modular forms. They furthermore showed that the smallest prime
conductor of an abelian surface  is~$277$;
in combination with the computations of~ \cite{MR3315514}, this
demonstrates that the conjecture is true in prime conductor less
than~$277$ (because there are neither any abelian surfaces nor suitable Siegel
modular forms).

These considerations are taken further in the recent
papers~\cite{BPVY,Berger}. In particular, these papers succeed in
establishing for the first time the modularity of (finitely many, up
to twist) abelian surfaces~$A$ with~$\End_{\C}(A) = \Z$.  (The
explicit examples in~\cite{BPVY} are conductors~$277$, $353$,
and~$587$, and the example in~\cite{Berger} is of conductor~$731$. It
should be noted that the abelian surfaces considered in
Theorem~\ref{thm:bcgpthmintro} do not include any of these examples;
as explained below, Theorem~\ref{thm:bcgpthmintro} is proved by
proving the existence of infinitely many abelian surfaces to which our
modularity lifting theorems apply, rather than by starting with
explicit examples of small conductor.)  These papers ultimately rely
on elaborate explicit computations of low weight Siegel modular forms,
developed in part by Poor and
Yuen~\cite{MR3315514,MR3713095,MR3498287}.

\subsubsection{Our modularity lifting theorem}
We now state our main modularity lifting theorem as it applies to abelian
surfaces.  The following theorem is proved
in~\S\ref{sec:modularityapplications}, see
Proposition~\ref{prop:infinitetwo}. (It
is possible to slightly weaken the hypothesis at~$v|p$ to deal with
certain abelian surfaces which have semistable reduction
at~$v|p$.)

\begin{theorem} \label{theorem:infinite} Let~$F$ be a totally real
  field in which~$p>2$ splits completely. Let~$A/F$ be an abelian
  surface with good ordinary reduction at all places~$v|p$, 
and suppose that, at each~$v|p$, the unit root crystalline eigenvalues are distinct modulo~$p$. 
Assume that~$A$ admits a polarization of degree prime to~$p$.  Let
$$\rhobar_{A,p}: G_F \rightarrow \GSp_4(\F_p)$$
denote the dual of the  mod-$p$ Galois representation associated
to~$A[p]$, and assume that~$\rhobar_{A,p}$ is vast and tidy in the
sense of Definitions~\ref{defn:vast} and~\ref{defn:tidy}.
Assume that~$\rhobar_{A,p}$ is  ordinarily modular, in the sense that
there exists an automorphic representation $\pi$ of $\GSp_4/F$ of parallel weight~$2$ and central
character~$|\cdot|^2$ which is ordinary at all~$v|p$, such
    that~$\rhobar_{\pi,p}\cong\rhobar_{A,p}$,
    and~$\rho_{\pi,p}|_{G_{F_v}}$ is pure for all finite
    places~$v$ of~$F$.
 Then~$A$ is modular, corresponding to a Hilbert--Siegel eigenform of
 parallel weight two.
\end{theorem}

Moreover, Proposition~\ref{prop:applications} shows that the modularity hypotheses on~$\rhobar_{A,p}$ can be omitted 
 in the following situations:
 \vspace{-1.5mm}
\begin{enumerate}
\item $p = 3$, and~$\rhobar_{A,3}$ is induced from a~$2$-dimensional
  representation with inverse cyclotomic determinant  defined over a totally
  real quadratic extension~$E/F$ in which~$3$ is unramified.
\item $p = 5$, and~$\rhobar_{A,5}$ is induced from a~$2$-dimensional
  representation  valued in~$\GL_2(\F_5)$  with inverse cyclotomic character defined over a totally real quadratic
extension~$E/F$ in which ~$5$ is unramified.
\item $\rhobar_{A,p}$ is induced from a character of a quartic CM
  field~$H/F$ in which~$p$ splits completely. 
\end{enumerate}

Theorem~\ref{theorem:infinite}  may be viewed as the genus two  analogue of~\cite[Thm.\
0.2]{MR1333035}, which is the main modularity lifting result proved in
that paper.  Proposition~\ref{prop:applications} is then
 the analogue  of~\cite[Thm.\ 0.6]{MR1333035}, which is a modularity  result for residually projectively dihedral representations.
The reason one cannot prove an analogue of ~\cite[Thm.\ 0.3]{MR1333035}  (which proves that all ordinary semistable elliptic curves over~$\Q$ with~$\rhobar_{E,3}$ absolutely irreducible are modular)
is that there is no argument to reduce the residual modularity of a surjective mod-$3$ representation~$\rhobar_3:G_F \rightarrow \GSp_4(\F_3)$
 (as in~\S5 of \emph{ibid}) to special cases of the Artin Conjecture (proved by Langlands--Tunnell).
Note that the difficulty is not simply that~$\GSp_4(\F_3)$ is not solvable (some of
the indicated representations above for~$p = 3$ and~$5$ are
non-solvable), but also that Artin representations do not contribute to the coherent cohomology of Shimura varieties in any  setting 
other than holomorphic (Hilbert) modular forms of weight one.

For~$E/F$ a totally real quadratic extension, the inductions of (modular) representations~$\varrhobar: G_{E}
\rightarrow \GL_2(\F_3)$ with
determinant~$\varepsbar^{-1}$ to~$G_{F}$ provide a large source of residually modular~$\rhobar$. We then
show that any such~$\rhobar: G_F \rightarrow \GSp_4(\F_3)$
with suitable determinant and local conditions at places~$v|3$ is equal to~$\rhobar_{A,3}$ for infinitely many  abelian surfaces~$A/F$ with~$\End_{\C}(A) = \Z$ and with good ordinary reduction at~$v|3$ (see Theorem~\ref{boxershorts}).  Theorem~\ref{theorem:infinite}
then implies that all such~$A$ are modular, and hence implies Theorem~\ref{thm:bcgpthmintro}.

\subsection{An overview of our argument}

Let $A$ be an abelian surface over a totally real field $F$. We may
assume that $\mathrm{End}_{F}(A) = \Z$ as otherwise, $A$ is of
$\mathrm{GL}_2$-type, in which case it is known that $A$ is
potentially modular. If $\mathrm{End}_{F}(A) = \Z$, a generalization
of the paramodular conjecture  predicts the existence of a holomorphic weight~ $2$ Hilbert--Siegel modular cuspidal eigenform $f$ (for the group $\mathrm{GSp}_4/F$) associated to $A$ in the sense that we have an equality of $L$-functions $L(f,s) = L (H^1(A),s)$.    
If such an equality holds,  we say that $A$ is modular. 

In this paper, we establish that  (under some mild further restrictions on~$A$),
 after possibly replacing the field $F$ by a finite totally real extension $F'$, the conjecture is true.

\begin{rem} There are situations where we don't prove
(even potentially) the paramodular conjecture for $A$. This is due to the presence of non-trivial endomorphisms of~$A$ over $\Qbar$. Nevertheless, we always express the $L$-function  of $A$ using automorphic forms on  groups $\mathrm{GL}_i/K$ for $i \in \{1,2,4\}$ and $K$ a number field, and thus establish Conjecture~\ref{conj:serre}. 
\end{rem}

On the surface, the modularity conjecture for abelian surfaces appears to be a generalization of the modularity conjecture for  elliptic curves. However, this analogy is somewhat
misleading.   Elliptic curves are regular motives with weights
$(0,1)$, whereas abelian surfaces are  irregular motives with weights
$(0,0,1,1)$. On the automorphic side, weight $2$ Hilbert modular
cuspforms occur in a single degree of the Betti and coherent
cohomology of the Hilbert modular varieties. Under mild assumptions,
there is an  elliptic curve associated to any Hilbert modular cuspidal
eigenform with rational Hecke field.

In contrast, weight $2$ Hilbert--Siegel modular cuspforms only occur
in the coherent cohomology of the Hilbert--Siegel modular
variety. More precisely, a holomorphic weight $2$ Hilbert--Siegel
modular cuspidal eigenform  can be viewed as a section of a line
bundle $\omega^2$ over the Hilbert--Siegel modular variety $X$;
here~$X$ is a smooth algebraic variety defined over~$\Q$ of
dimension~$3[F:\Q]$ which parametrizes abelian schemes of
dimension~$2[F:\Q]$ equipped with an action of~$\OL_F$, a level
structure, and a polarization.  Moreover, in the ``generic case'',
such an eigenform contributes to cohomology in degrees~$0$ to~$[F:\Q]$. Since the Hecke eigenvalues associated to such modular forms are not realized in the \'etale cohomology of a Shimura variety,  we don't know how to associate a  ``motive'' to a weight $2$ Hilbert--Siegel modular cuspidal eigenform, but only a compatible system of Galois representations which should correspond to the system of $\ell$-adic realizations of this motive. These Galois representations are constructed by using congruences. 

From a technical point of view, it turns out that the modularity
conjecture for abelian surfaces over a totally real field~$F$ is closely
related
to the~$2$-dimensional
odd Artin conjecture for~$F$  (now a theorem), which is the existence of a bijection preserving $L$-functions between the following  objects:

\begin{itemize}  
\item Irreducible, totally odd, two dimensional complex representations
  of the absolute Galois group of $F$, and

\item Hilbert modular cuspidal eigenforms (newforms) of weight one. 
\end{itemize}

$2$-dimensional odd Artin representations have irregular Hodge--Tate weights~$(0,0)$, and Hilbert modular forms of weight one only occur in the  coherent cohomology of the Hilbert modular variety, where
they contribute  in degrees~$0$ to~$[F:\Q]$.

We now review some of the strategies employed in the proof of Artin's
conjecture, as they have served as an inspiration for our current
work. As with almost all  modularity theorems, one  proceeds by combining
a modularity lifting theorem with residual modularity (that is, the modularity of the mod~$p$ representation).
In the case of Artin's conjecture, residual modularity ultimately (if quite indirectly) comes
from the Langlands--Tunnell theorem, 
 whereas in our
setting, the
residual potential modularity comes from a straightforward application
of Taylor's method~\cite{MR1954941} using a theorem of
Moret-Bailly. Accordingly, we ignore the question of residual
modularity for the rest of this introduction, and concentrate on
explaining the modularity lifting theorems.

The first modularity (lifting) theorems which applied to  two dimensional odd Artin  representations~$\rho$  over $\Q$ were obtained by 
 Buzzard--Taylor and Buzzard~\cite{MR1709306,MR1937198}. There is  an
 obstruction to generalizing the Taylor--Wiles method (which was
 originally applied in the regular case of Hodge--Tate
 weights~$(0,1)$ and weight two modular
 forms~\cite{MR1333035,MR1333036}) to the irregular case of
 weights~$(0,0)$ and weight one modular forms. This obstruction lies
 in the fact that weight one forms occur in degrees~$0$ and~$1$ of the
 coherent cohomology and  that there exist non-liftable mod~$p$
 weight one eigenforms. (There is also a reflection of this
 obstruction on the Galois theoretic side --- the corresponding local
 deformation ring at~$p$ has dimension one less in the irregular
 weight case.) Instead, Buzzard and Taylor proceed quite differently.
 
 Choose a prime~$p$ and view~$\rho$ as a~$p$-adic
 representation with finite image. We also assume that~$\rho$ is
 unramified at~$p$ and let~$\alpha$, $\beta$ denote the Frobenius
 eigenvalues. For simplicity, we also assume that $\alphabar \neq
 \betabar$ (where the bar denotes reduction modulo~$p$). We have that  $$\rho\vert_{G_{\Q_p}}  \simeq \begin{pmatrix}
    \lambda_\alpha &0\\0&\lambda_\beta
  \end{pmatrix}$$
  for the unramified characters $\lambda_\alpha$ and $\lambda_\beta$
  taking a Frobenius element to~$\alpha, \beta$ respectively.

 The strategy of Buzzard and Taylor is to first replace the space of
 classical weight one modular forms by a bigger space of ordinary
 $p$-adic modular forms of weight one. On the Galois side, classical weight one eigenforms (of level prime to $p$) have associated Galois  representations which are  unramified at $p$, while an  ordinary $p$-adic modular form $f$ of weight one has an associated Galois representation  which may be ramified at $p$ of the form:
 
 $$\rho_f\vert_{I_{\Q_p}}  \simeq \begin{pmatrix}
 1 &\star\\0& 1
  \end{pmatrix}$$
 
 Moreover, $f$ should be  classical if and only if $\star = 0$. 
 A key advantage of working with ordinary~$p$-adic modular forms
 is that they are defined as sections of a line bundle over the ordinary locus, which is affine,
 and thus only occur in cohomological degree~$0$. It follows that
ordinary $p$-adic modular forms of weight one are unobstructed for
congruences 
and one can (assuming residual modularity) apply the Taylor--Wiles method 
in this setting
to  deduce the existence of two $p$-adic ordinary weight one modular forms $f_\alpha$ and $f_\beta$ such that $\rho_{f_\alpha} = \rho_{f_\beta} = \rho$ and $U_p f_\alpha = \alpha f_\alpha$, $U_p f_\beta = \beta f_\beta$.

We observe that the existence of both $f_\alpha$ and $f_\beta$ witnesses the fact that $\rho$ is unramified at $p$. In order to show that $f_\alpha$ and $f_\beta$ are classical forms of weight one, one forms the linear combinations $h = (\alpha f_\alpha - \beta f_\beta)/(\alpha - \beta)$ and $g =  (f_\alpha -  f_\beta)/(\alpha - \beta)$. The property that  $\rho_{f_\alpha} = \rho_{f_\beta} = \rho$ and the explicit relation between $q$-expansions and Hecke eigenvalues translates into the geometric property that $\Frob (h) = g$.  Using rigid analytic techniques, one can show that this property implies that $f_\alpha, f_\beta$ are classical forms of weight one. 
This strategy has been successfully generalized to any totally real
field~\cite{MR3104551,MR3294624,MR3581174,MR3581178,MR3639600}.

From a different direction, the paper~\cite{CG} introduced an alternate method for proving modularity lifting results in 
weight one, by modifying the method of Taylor--Wiles and exploiting the Galois representations associated to coherent cohomology
classes in  all  degrees.  This method eliminates the delicate
classicality theorem in weight one because one only works with classical
(but possibly higher degree) cohomology. 
This method allows in principle to deal with any obstructed situation,
but requires some non-trivial input. For  $2$-dimensional odd Artin
representations over a totally real fields, one needs to prove that
(after suitable localization at a maximal ideal of the Hecke algebra)
the cohomology in weight one is supported in degrees $0$ to $[F: \Q]$
(this  is actually automatic here for cohomological dimension
reasons), and that the Galois representations in all cohomological
degrees satisfy a form of local-global compatibility (at places above
$p$). This last property has been proved when $F= \Q$ where one can
reduce to studying degree $0$ torsion cohomology classes and use the
``doubling method'' described below, but has not yet been proved for all primes~$p$ over a general totally real field
(though see~\cite{MRX} for some partial results).

 After this discussion of Artin's conjecture, we return to the paramodular conjecture.  We first assume that $F = \Q$ and fix a prime $p$. We assume that $A$ has ordinary good reduction at $p$ so that 
 
$$\rho_{A,p} \vert_{G_{\Q_p}} \simeq \begin{pmatrix}
 \lambda_\alpha &0 & \star & \star \\0& \lambda_\beta & \star & \star \\
 0 & 0 & \lambda_{\beta}^{-1} \varepsilon^{-1} & 0 \\
 0 &0 & 0 & \lambda_\alpha^{-1} \varepsilon^{-1}
  \end{pmatrix},$$
  where, additionally, we assume that $\alphabar \ne \betabar$.
  (The Weil bounds together with the Cebotarev density theorem guarantee
 an ample source of such primes~$p$.) Tilouine and his
  collaborators~\cite{Til1,Til2,Til3,Til4,MR2234862, MR2264659, Til5} developed
  modularity lifting results for~$\GSp_4/\Q$ in regular weight.  In
  the case of Hodge--Tate weights $(0,0,1,1)$,  the paper
  ~\cite{MR2920881}  
  applied these techniques to ordinary $p$-adic
  modular  forms of weight $2$ to produce (under technical
  assumptions) two $p$-adic eigenforms $f_\alpha$ and $f_\beta$ associated to $A$ (see also \cite{MR2264659, MR3025747}, where the case of certain $\mathrm{GSp}_4$-type abelian varieties is treated).  
 
 Similarly to  the case of $\mathrm{GL}_2/\Q$, an ordinary $p$-adic modular form of weight $2$ has a Galois representation whose restriction to inertia at $p$ has the shape:

 $$ \begin{pmatrix}
1 & \star_1 & \star & \star \\0& 1 & \star & \star \\
 0 & 0 &  \varepsilon^{-1} &  \star_2 \\
 0 &0 & 0 &  \varepsilon^{-1}
  \end{pmatrix}.$$
Such a form should be classical if and only if  its Galois representation is de Rham --- equivalently: $\star_1= \star_2 = 0$ (because of the symplectic structure, the vanishing of~$\star_1$
is equivalent to the vanishing of~$\star_2$).
 
As before, the existence of both $f_\alpha$ and $f_\beta$ witnesses the property that $A$ is de Rham at $p$. 
  One difficulty, however, is that the Fourier expansions of Siegel modular forms
are \emph{not} explicitly determined by the Hecke eigenvalues (although we often have an abstract multiplicity one theorem). In particular,
 one doesn't know how to deduce geometrically from~$\rho_{f_\alpha} = \rho_{f_\beta}  = \rho_{A,p}$ 
  that there exist suitable linear combinations of $f_\alpha$ and $f_\beta$ giving rise to the desired
  form~$f$ by   mimicking the Buzzard--Taylor argument.

   In another direction, in~\cite{CGGSp4} the modified Taylor--Wiles method  was applied to low weight Siegel modular forms over~$\Q$.
There were a number of serious difficulties which prevented the
authors  from deducing any unconditional modularity lifting for abelian surfaces.
The idea of the method is to consider (a suitable localization of) the full cohomology complex $\mathrm{R}\Gamma(X, \omega^2)$ where $X$ is an integral model over $\Z_p$ of the Siegel threefold. The required inputs are: 

\begin{enumerate}
\item to prove that the cohomology is only supported in degrees~$0$
  and $1$, and
\item to prove local-global compatibility for the cohomology classes. 
\end{enumerate} 

The first point  is subtle in the weight
 of interest, because the cohomology groups will not generally vanish
 before localization at some non-Eisenstein maximal ideal~$\m$ (and
 indeed this point was not established in weight~$2$ in~\cite{CGGSp4}).
The paper \cite{CGGSp4} proved the second point for torsion degree $0$
cohomology classes, using a ``doubling'' argument that we will return
to below.

One crucial new
ingredient which allows us to proceed in the symplectic case and deal with $(1)$  is the higher Hida theory developed
for~$\GSp_4$ over~$\Q$ in~\cite{pilloniHidacomplexes}.  The idea of~\cite{pilloniHidacomplexes} is (loosely speaking) to work over the larger space
which is the complement of the  supersingular locus (the rank~$\ge  1$ strata),
 which is now no longer affine.
(Since we are working in mixed characteristic, one should imagine this taking place in the category of formal schemes, as in classical Hida theory.)  Since the cohomological dimension of these spaces is 
one (more precisely, the image of these spaces in the minimal compactification has
cohomological dimension one, which is sufficient for our purposes), there should exist complexes of 
amplitude~$[0,1]$
computing the coherent cohomology of all the relevant vector bundles. The main result of~\cite{pilloniHidacomplexes} is that suitably constructed Hida idempotents cut down such a complex to  a perfect complex, and moreover that the cohomology of this perfect complex is computed in characteristic zero by the space of weight $2$ automorphic forms of interest. 
 A crucial ingredient  in order to study the coherent cohomology is therefore the introduction of  Hecke operators at~$p$ and their associated projectors. 

 A version over~$\Q$ of our modularity lifting theorem could be proved
 by applying the patching method of~\cite{CG} to the higher Hida
 complexes of~\cite{pilloniHidacomplexes}.  It should nevertheless be
 noted that, even if we were only interested in theorems over~$\Q$,
 we are forced to prove a modularity lifting theorem
 for any totally real field $F$ (and prime $p$ which splits completely
 in it). This is  because we need to employ Taylor's Ihara avoidance
 technique~\cite{tay} to deal with issues of level raising and
 lowering at places away from~$p$, and this step crucially relies on
 using solvable base change.  We can then combine this
 modularity lifting result with base change techniques and the Moret-Bailly argument to
 achieve residual potential modularity, in order to prove our main potential
 modularity
 theorem. 

In the light of the above discussion, in order to prove a modularity
lifting theorem for Hilbert--Siegel modular forms it is natural to
consider  (a suitable localization of) either the cohomology complex
$\mathrm{R}\Gamma(X, \omega^2)$ where $X$ is an integral model over
$\Z_p$ of the Hilbert--Siegel space, or of the ordinary part of the
cohomology complex for a subspace of~$X$ obtained from the $p$-rank stratification. The required inputs for the modified Taylor--Wiles method are now:

\begin{enumerate}
\item to prove that the cohomology is only supported in degrees~$0$ to
  $[F:\Q]$, and
\item to prove  local-global compatibility for the cohomology classes. 
\end{enumerate} 

It is to some extent possible to solve $(1)$ using higher Hida theory
(although there are some issues), but $(2)$ seems to be a more serious
problem because we only know how to prove that the Galois representations associated to torsion classes in~$H^i$  satisfy
the right local-global compatibility condition at~$v|p$
if~$i=0$. Accordingly, we are unable to argue directly with such complexes.

Let the number of non-zero degrees of cohomology of the spaces we are
considering be~$l_0+1$; we refer to~$l_0$ as the \emph{defect}. (The
original Taylor--Wiles method only applies if~$l_0=0$, while
if~$l_0>0$ we use the method of~\cite{CG}. As mentioned above, $l_0$
also has a Galois-theoretic interpretation: the sum of the dimensions
of the local deformation rings is~$l_0$ less than the corresponding
dimension in the defect~$0$ case.) One key trick we employ in this
paper is to reduce to situations where we only have to consider
cohomology in at most two degrees (so the defect is at most one), i.e.\ it suffices to work with complexes consisting of at most two
terms.  This is where we take advantage of the product situation at
$p$ (because $p$ splits in the totally real field). (Implicitly, what
happens in this case is that any cohomology occurring in~$H^1$ can
also be seen via the Bockstein homomorphism as coming from~$H^0$,
\emph{provided} that the characteristic zero classes in~$H^1$ are also
seen by the characteristic zero classes in~$H^0$, and this can be
established by automorphic considerations; so we only have to prove
local-global compatibility for~$H^0$.) We now explain how we do this
in slightly more detail.

We assume that  $A$ has ordinary good reduction at all places   $v|p$, so that 
 
$$\rho_{A,p} \vert_{G_{F_v}} \simeq \begin{pmatrix}
 \lambda_{\alpha_v} &0 & \star & \star \\0& \lambda_{\beta_v} & \star & \star \\
 0 & 0 & \lambda_{\beta_v}^{-1} \varepsilon^{-1} & 0 \\
 0 &0 & 0 & \lambda_{\alpha_v}^{-1} \varepsilon^{-1}
  \end{pmatrix},$$ where we furthermore assume that~$\alphabar_v\ne\betabar_v$.

Although we expect that there should be a weight $2$ eigenform
associated to $A$ of spherical level at $p$ (because $A$ has good
reduction at~$p$), it turns out that because~$A$ is ordinary at~$p$,
it is more natural to look for an eigenform $f$ associated to $A$ of Klingen level at $p$. The Klingen level structure is given by choosing a subgroup of order $p$ inside $A[v]$ for all $v|p$. At Klingen level at~$v$,
there is a Hecke operator~$U_{\Kli(v),1}$ whose eigenvalue on $f$ should be~$\alpha_v + \beta_v$,
and a second Hecke operator~$U_{\Kli(v),2}$ whose eigenvalue should be~$\alpha_v \beta_v$.
We observe that the second operator has an invertible eigenvalue (we say that $f$ is Klingen ordinary) and this corresponds to the fact  that the Galois representation $\rho_{A,p}|_{G_{F_v}}$  is ordinary.

There is another level structure that plays a role: the  Iwahori level structure given by choosing a complete self dual flag of subgroups inside  $A[v]$. 
For each $v|p$, there are two degeneracy maps  from Iwahori level to Klingen level, and there  are  Hecke
 operators~$U_{\Iw(v),1}, U_{\Iw(v),2} = U_{\Kli(v),2}$  at Iwahori level.  Pulling back the expected form $f$ by the degeneracy maps should yield eigenforms at Iwahori level which have 
  eigenvalues~$\alpha_v$
 and~$\beta_v$ for $U_{\Iw(v),1}$ (we call them Iwahori ordinary).

We now return to the question of using modularity lifting theorems to
find~$f$. First of all, modularity lifting theorems with $p$-adic
ordinary modular forms (i.e.\ with $l_0=0$) allow us to construct $2^{[F:\Q]}$ Iwahori ordinary $p$-adic modular forms 
whose eigenvalue for $U_{\Iw(v),1}$ is~$\alpha_v $ or $\beta_v$,
and  whose eigenvalue for $U_{\Kli(v),2}$ is~$\alpha_v \beta_v$. We suspect that these forms are classical, but as explained before, we don't know how to establish any geometric relation between them. 

As a second step we apply a modularity lifting theorem in the case
that the defect $l_0$ equals one. Let us isolate a place
$v|p$. Using higher Hida theory, we construct a perfect complex
of amplitude $[0,1]$ which is obtained by taking the ordinary (more precisely Iwahori ordinary  at $w \neq v$, Klingen ordinary  at $v$)
cohomology of the open subspace of the Hilbert--Siegel Shimura variety which is ordinary and carries an Iwahori level structure at all places $w \neq v$, and has $p$-rank at least one at $v$ and carries a Klingen level structure. 

We manage to prove that this cohomology carries a Galois
representation which  has the following type of local-global
compatibility property: \begin{enumerate}

\item For all places $w|p$, $w \neq v$: 
$$ \rho_{A,p} \vert_{I_{F_w}} \simeq \begin{pmatrix}
1 & \star & \star & \star \\0& 1 & \star & \star \\
 0 & 0 & \varepsilon^{-1} & \star \\
 0 &0 & 0 & \varepsilon^{-1}
  \end{pmatrix}.$$
  \item For $v$: $$ \rho_{A,p} \vert_{I_{F_v}} \simeq\begin{pmatrix}
1 & 0 & \star & \star \\0& 1 & \star & \star \\
 0 & 0 & \varepsilon^{-1} & 0\\
 0 &0 & 0 & \varepsilon^{-1}
  \end{pmatrix}.$$
\end{enumerate}

Using the methods of~\cite{CG}, we can prove a modularity lifting theorem, and produce 
$2^{[F:\Q]-1}$  $p$-adic modular forms (which converge a lot more in the $v$ direction)
 whose eigenvalue for $U_{\Iw(w),1}$ is~$\alpha_w $ or $\beta_w$ if $w
 \neq v$,  and whose eigenvalue for $U_{\Kli(v),1}$ is~$\alpha_v  +
 \beta_v$, and  whose eigenvalue for $U_{\Kli(w),2}=U_{\Iw(w),2}$ is~$\alpha_w \beta_w$ for all $w|p$. 
 
 Our last step is to prove lots of linear relations between all these forms  we have constructed. This step ultimately relies upon an abstract multiplicity one result
 which we prove using the Taylor--Wiles method. Exploiting these linear relations  and using \'etale descent  techniques, we first manage to construct a Klingen ordinary weight $2$ modular form defined  on the open subspace of the Hilbert--Siegel Shimura variety which  has $p$-rank at least one at all $v|p$ and carries a Klingen level structure. We then manage, using analytic continuation techniques, to prove that this form extends to the full Shimura variety and is therefore classical.

\subsection{An outline of the paper}We briefly explain the outline of
the paper; we refer the reader to the introductions to the individual
sections for a further explanation of their contents, and for some
elaborations on the overview of our arguments above.

In \S\ref{sec:background} we recall some more or less standard
background material on Galois representations, the local Langlands
correspondence, local representation theory, and related topics. \S\ref{sec: integral
  Shimura varieties} discusses the Shimura varieties which we use, and
some properties of their integral models and compactifications, and
recalls the approach to the normalization of Hecke operators on coherent cohomology via
cohomological correspondences which was introduced
in~\cite{pilloniHidacomplexes}. 

In \S\ref{sec: Hida complexes} we construct the Hida complexes
that we work with, and prove some of their basic properties (in particular, we
prove that they are perfect complexes). In \S\ref{sec:doubling}
we establish the ``doubling'' results that we will later use to prove
local--global compatibility for Hilbert--Siegel modular forms over torsion rings.
The basic strategy (employed
in a number of other places, see~\cite{MR1074305,edix:weights,MR3247800,CG,CGGSp4})
is to show that we can embed (via degeneracy maps)
two copies of our space of  ordinary modular forms at Klingen level into a space
of ordinary modular forms of Iwahori level. This allows us to show that the corresponding
Galois representations are ordinary (in the Iwahori sense) in two different ways, namely,
with~$\alpha_w$ and~$\beta_w$ as unramified subspaces. Then the genericity assumption~$\alphabar_w \ne \betabar_w$ forces there to be a~$2$-dimensional unramified summand of our representation.
The key technical difficulty is proving that the direct sum of the degeneracy maps does
indeed give an embedding. All previous incarnations of the doubling phenomenon ultimately
relied on the~$q$-expansion principle, but our argument is more
geometric, and ultimately rests on analyzing the effect of the Hecke operator~$Z_w=U_{\Kli(w),1} - U_{\Iw(w),1}$ along the~$w$-non-ordinary locus. 

In  \S\ref{sec: higher Coleman
  theory} we prove that a characteristic zero classicality result for
the~$H^0$ of our Hida complexes, using Coleman theory. We also show that the complexes we
consider are balanced, in the sense that they have Euler
characteristic zero, 
using a somewhat intricate interplay between three objects --- the complex of classical
forms, the complex of overconvergent forms, and our
complex of (Klingen) ordinary forms.

In~\S\ref{sec:CG} we carry out our main Taylor--Wiles patching arguments
in the cases that~$l_0=0$ and~$l_0=1$. We then prove our main
modularity lifting theorem in \S\ref{sec:gluing}, using analytic
continuation, \'etale descent, and linear algebra arguments based on
the doubling results of~\S\ref{sec:doubling} to reduce to the
classicality results of~\S\ref{sec: higher Coleman
  theory}.

In~\S\ref{section: potential modularity of abelian surfaces}, we apply our main automorphy lifting theorem
to prove the potential automorphy of abelian surfaces. The basic idea is to
use a version of the~$p$-$q$ trick (first employed by Wiles as the
$3$-$5$ trick), together with an application of a theorem
of Moret-Bailly, to connect general abelian surfaces via a chain of congruences
to the restriction of scalars of an elliptic curve over a totally real quadratic extension of~$F$,
which we know already by~\cite{MR1954941} to be potentially modular.
We are also left to deal directly with some cases
 of abelian surfaces with small Mumford--Tate groups,  which can mostly
 be done immediately with an appeal to the theory
 of Grossencharacters. We also include a number of applications as mentioned in the introduction,
 including elliptic curves over quadratic extensions of~$F$.
 
 In~\S\ref{sec:modularityapplications}, we give applications to the automorphy of abelian
 surfaces. We show that, given any mod~$3$ representation~$\rhobar: G_{\Q}
 \rightarrow \GSp_4(\F_3)$ with (inverse) cyclotomic similitude character, 
  it can be realized (in infinitely many ways) as the ~$3$-torsion of an abelian surface over~$\Q$.
 Here we exploit some classical geometry related to the Burkhardt quartic, which is isomorphic
 to a compactification of~$\mathcal{A}_2(3)$. The key point is to show that the variety
 given by the twist of~$\mathcal{A}_2(3)$ by~$\rhobar$ has sufficiently many rational points.
 We do this by proving it is unirational over~$\Q$ via a map of degree at most~$6$.
 The argument is similar to that of~\cite{MR1415322}, except that it is applied not to the
 twist of~$\mathcal{A}_2(3)$ itself but to a twist of a degree~$6$ rational cover, which has
 the pleasing property (unlike the Burkhardt quartic itself) that the birational map 
 to~$\mathbf{P}^3$ over~$\Q$ can be made equivariant with respect to the action of the
 automorphism group~$\PSp_4(\F_3)$. Finally, we conclude with a discussion of
 the paramodular conjecture and its relationship to the standard conjectures,
 and explain why the original formulation of this conjecture requires a minor modification.

\subsection{Some further remarks}

For length reasons, we did not try to optimize all of our theorems --- for example, our arguments would surely extend
to prove the potential automorphy of some~$\GSp_4$-type abelian varieties, but sticking with abelian surfaces makes the Moret-Bailly arguments
somewhat simpler, and (by using a trick) we manage to avoid any character
building whatsoever. However, we have gone to some lengths to treat
the case~$p=3$, and to use a weaker notion of $p$-distinguishedness
than in~\cite{CGGSp4}; while this is not necessary for our
applications to potential modularity, it significantly increases the
applicability of our theorems to actual modularity problems.

\subsubsection{The work of Arthur} \label{arthur}
It should be noted that we use Arthur's multiplicity formula for the
discrete spectrum of~$\GSp_4$, as announced in~\cite{MR2058604}. A
proof of this (relying on Arthur's work for symplectic and orthogonal
groups in~\cite{MR3135650}) was given in~\cite{GeeTaibi}, but this
proof is only as unconditional as the results of~\cite{MR3135650}
and~\cite{SFTT1,SFTT2}. In particular, it depends on cases of the
twisted weighted fundamental lemma that were announced
in~\cite{MR2735371}, but whose proofs have not yet appeared, as well as on the references [A24], [A25], [A26] and~[A27]
in~\cite{MR3135650}, which at the time of writing have not appeared publicly.  

\subsubsection{Curves of higher genus}

One may well ask whether the methods of this paper could be used to prove (potential) modularity of curves of genus~$g \ge 3$  whose Jacobians have  trivial endomorphism rings. 
At the moment, this seems exceedingly unlikely without some substantial new idea. All generalizations of the Taylor--Wiles method to this point require that the automorphic representations
in question  are associated to the Betti cohomology groups of locally
symmetric spaces, or the coherent cohomology groups of Shimura varieties, which have integral structures and hence allow one to talk about
\emph{congruences} between automorphic forms. Symplectic motives of rank~$2g$ over~$\Q$ are conjecturally associated to automorphic
representations for the (split) orthogonal group~$\SO_{2g+1}$ (when~$g=1$ or~$g = 2$,
there are  well-known exceptional isomorphisms which allow us to replace~$\SO_{2g+1}$ by the groups~$\GL_2$ and~$\GSp_4$
respectively). Following~\cite{MR3165645}, Gross has made some  precise conjectures 
concerning the level structures of newforms associated to such conjectural automorphic
representations in~\cite{MR3535358}. 

The automorphic representations contributing to the Betti cohomology
groups of locally symmetric spaces have regular infinitesimal
characters, so can only be used for~$g=1$. The automorphic
representations contributing to the coherent cohomology of orthogonal
Shimura varieties are representations of the inner form~$\SO(2g-1,2)$
of~$\SO_{2g+1}$ (which is non-split if~$g>1$), 
 whose infinity components~$\pi_\infty$ are furthermore either discrete series, or
 non-degenerate limits of discrete series.

If~$g=1$, the representations considered by Gross in~\cite{MR3535358}
are discrete series, and if~$g=2$, they are non-degenerate limits of
discrete series, but if~$g\ge 3$, then neither possibility occurs, so the automorphic representations do not contribute to the cohomology (of any kind)
of the corresponding Shimura variety. (Another way of seeing this is
to compute the possible infinitesimal characters of the automorphic
representations corresponding to automorphic vector bundles on the
Shimura variety, or equivalently the Hodge--Tate weights of the
expected $2g$-dimensional symplectic Galois representations; one finds that no Hodge--Tate weight
can occur with multiplicity bigger than~$2$, while the symplectic
Galois representations coming from the \'etale $H^1$
of a curve of genus~$g$ have weights~$0,1$ each occurring with multiplicity~$g$.)
 In particular, the general modularity problem for curves of genus~$g \ge 3$ seems at least as hard as proving non-solvable
cases of the Artin conjecture for totally even representations, and even proving the modularity of a single such curve
with Mumford--Tate group~$\GSp_{2g}$ seems completely
out of reach.

On the other hand, there are some special families in higher genus which may well be amenable to our method. In particular, 
the Tate module of a cyclic trigonal genus three curve
(so-called Picard curves, with affine equations of the form~$y^3 = x^4 + a x^2 + b x+ c$) 
defined
over~$\Q$ splits (over~$\Q(\sqrt{-3})$) into two essentially conjugate self-dual irregular~$3$-dimensional representations of~$G_{\Q(\sqrt{-3})}$.
These Galois representations conjecturally correspond (see the appendix to~\cite{MR2264659}) to automorphic representations~$\pi$ for a form of~$U(2,1)/\Q$ (splitting over~$\Q(\sqrt{-3})$) 
such that~$\pi_{\infty}$ is a non-degenerate limit of discrete series and contributes to the coherent cohomology of the associated Shimura variety.
The methods of this paper should apply (in principle) to these
curves. 

\subsubsection{\emph{K3} surfaces}Our results should also have
applications to the Hasse--Weil conjecture for K3 surfaces over totally real fields with geometric Picard number~$\ge 17$. While we do not
undertake a detailed study of this problem here, we discuss it in~\S\ref{subsec: K3}. 

\subsubsection{A comparison of this paper with~\cite{10author}} \label{section:comparison}
It follows from Theorem~\ref{theorem:main} that any elliptic curve~$E$ over a CM field~$K/F$ is potentially modular (simply
consider the abelian surface given by Weil restriction of scalars of~$E$ from~$K$ to~$F$). This result is also
proved in~\cite{10author}.	 Perhaps surprisingly, there is relatively little overlap between the two proofs.
For example, our argument does not require any of the results of Scholze~\cite{scholze-torsion} on the construction of Galois
representations, nor the derived version of Ihara avoidance required in~\cite{10author}. The only common theme is the use
of the modified Taylor--Wiles method of~\cite{CG}. To further illustrate the difference, it is also proved in~\cite{10author}
that the~$n$th symmetric power of any such~$E$ is potentially automorphic, which is not directly accessible from our approach.
On the other hand,
we also deduce (Theorem~\ref{thm:genusone}) the potential modularity of elliptic curves over fields like~$F = \Q(\sqrt[4]{2})$, which seems out of reach
using the methods of~\cite{10author}.

\subsection{Acknowledgements}We would like to thank all of the
referees for their numerous helpful comments and corrections.
Both David Geraghty and Jacques Tilouine have made important contributions to the problem
of modularity for abelian surfaces; we would especially like to thank them for
many helpful discussions over the years.
We would also like to thank Patrick Allen,
Kevin Buzzard,
Brian Conrad,
Matthew Emerton, Najmuddin Fakhruddin, Dick Gross, Robert Guralnick, Florian Herzig,
Christian Johansson, Keerthi Madapusi Pera, Rutger Noot, Madhav Nori, Ralf Schmidt,
Olivier Ta\"{\i}bi, Richard Taylor, Jack Thorne,
Andrew Wiles, and Liang Xiao for
helpful conversations. We finally want to thank Bruno Klingler and the
University of Paris 7 for hosting us during part of this project.

\section{Background
  material}\label{sec:background}In this
section we recall a variety of more or less
well-known results that we will use in the
body of the paper.
\subsection{Notation and conventions}\label{subsec: notation and conventions}\subsubsection{$\GSp_4$}We define $\GSp_4$ to be the reductive group
over $\Z$
defined as a subgroup of $\GL_4$
by \[\GSp_4(R)=\{g\in\GL_4(R):gJ g^{t} =\nu(g)J\}\]
where $\nu(g)$ is the
similitude factor (which is uniquely determined by $g$, and which we
sometimes call the multiplier factor), and $J$ is
the antisymmetric matrix \[\begin{pmatrix}0& s\\-s&0\end{pmatrix}\]where $s=
\begin{pmatrix}
  0&1\\1&0
\end{pmatrix}
$. Note that the map $\nu:g\mapsto\nu(g)$ is a
homomorphism $\GSp_4\to\Gm$.

We let~$\Sp_4$ be the subgroup with~$\nu=1$, and we let
$B\subset G=\GSp_4$ be the Borel subgroup of upper triangular
matrices, and $T\subset B$ be the diagonal maximal torus. 
 Write $W_G=N_G(T)/T$ for the Weyl group of $(G,T)$. It acts on the
 character group via
$w\cdot \lambda (t)=\lambda(w^{-1}tw)$.
It is generated by $s_1=\left(\begin{array}{cccc}s&0_2\\0_2&s\end{array}\right)$
and
$s_2=\left(\begin{array}{ccc}1&0&0\\0&s'&0\\0&0&1\end{array}\right)$ where $s'=
\begin{pmatrix}
  0&1\\-1&0
\end{pmatrix}
$,
and admits the presentation
$$W_G=\langle s_1,s_2|s_1^2=s_2^2=(s_1 s_2)^4=1 \rangle.$$

Write $X^*(T)$ (resp.\ $X_*(T)$) for the group of characters (resp.\
cocharacters) of $T$. We identify $X^*(T)$ with the lattice in $\Z^3$ of triples
$(a,b;c)\in\Z^3$ such that $c\equiv a+b\pmod 2$ via
\[ \lambda:t=\diag(t_1,t_2,\nu t_2^{-1},\nu t_1^{-1})\mapsto t_1^a t_2^b \nu^{(c-a-b)/ 2}. \] 
In particular, the central character is given by
$\lambda(\diag(z,z,z,z))=z^c$.  The simple roots are  $\alpha_1=(1,-1;0)$ and $\alpha_2=(0,2;0)$;
$\alpha_1$ is the short root. Note that the $\alpha_i$ determine the reflections
$s_i$. The similitude factor
is $(0,0;2)$.

The root datum~$(G,B,T)$ determines the dual root
datum~$(\widehat{G},\widehat{B},\widehat{T})$, where~$\widehat{G}$ is
the dual group~$\mathrm{GSpin}_5$. We always identify~$\mathrm{\GSpin}_5$
with~$\GSp_4$ via the spin isomorphism (see for example~\cite[\S
3.2]{Til4} for a detailed explanation of this). In particular,
the cocharacter in~$X_*(\widehat{T})$ corresponding to the
character~$(a,b;c)\in X^*(T)$ defined above is given
by \[t\mapsto\diag(t^{(a+b+c)/2},t^{(a-b+c)/2},t^{(-a+b+c)/2},t^{(-a-b+c)/2}). \]

 We write~$\fg$ and~$\fb$ for the Lie algebras of~$\GSp_4$
and~$B$, and~$\fg^0$ and~$\fb^0$ for the Lie algebras of~$\Sp_4$
and~$B\cap\Sp_4$. If~$v$ is a finite place of a number field~$F$, with
residue field~$k(v)$, then we have the standard parahoric subgroups of~$\GSp_4(F_v)$:
\begin{itemize}
\item  The hyperspecial subgroup $\mathrm{GSp}_4( \CO_{F_v} )$. 
\item The paramodular subgroup  $\Par(v)$, the stabilizer
  in~$\GSp_4(F_v)$ of~$\cO_{F_v}\oplus \cO_{F_v}\oplus \cO_{F_v}\oplus
  \varpi_{v}\cO_{F_v}$, where $\varpi_v\in\cO_{F_v}$ is a uniformizer.
\item The Siegel parahoric $\Si(v)$, the preimage
  in~$\GSp_4(\cO_{F_v})$ of those matrices in~$\GSp_4(k(v))$ of the
  form \[
    \begin{pmatrix}
      *&*&*&*\\ *&*&*&*\\0&0&*&*\\0&0&*&*
    \end{pmatrix}
.\]
\item the Klingen parahoric $\Kli(v)$, the preimage
  in~$\GSp_4(\cO_{F_v})$ of those matrices in~$\GSp_4(k(v))$ of the
  form \[
    \begin{pmatrix}
      *&*&*&*\\0&*&*&*\\0&*&*&*\\0&0&0&*
    \end{pmatrix}
.\]
\item the Iwahori subgroup $\Iw(v)$, the preimage of~$B(k(v))$ in~$\GSp_4(\cO_{F_v})$.
\end{itemize}

\subsubsection{Algebra}If $R$ is a local ring we write $\mf{m}_{R}$ for the
maximal ideal of~$R$.

If $M$ is a perfect field, we let $\barM$ denote an algebraic closure of $M$
and $G_M$ the absolute Galois group $\Gal(\barM/M)$.  For each
prime~$p$ not equal to the characteristic of~$M$, we let $\varepsilon_p$ denote the $p$-adic cyclotomic
character and $\varepsilonbar_p$ its reduction modulo $p$. We will  usually drop~$p$ from
the notation and simply write~$\varepsilon$,
$\varepsilonbar$.

If $K$ is a finite extension of $\Qp$ for some $p$, we write
$K^\nr$ for its maximal unramified extension; $I_K$ for the inertia
subgroup of $G_K$; $\Frob_K \in G_K/I_K$ for the geometric Frobenius;
and $W_K$ for the Weil group.  If $L/K$ is a Galois extension we will
write $I_{L/K}$ for the inertia subgroup of $\Gal(L/K)$.  We will
write $\Art_K:K^\times \iso W_K^\ab$ for the Artin map normalized to
send uniformizers to geometric Frobenius elements.

If $\rho$ is a continuous representation of~$G_K$ over~$\Qlbar$ for
some~$l\ne p$, valued either in some~$\GL_n$ or in~$\GSp_4$, then we write $\WD(\rho)$ for the corresponding
Weil--Deligne representation. (By definition, a $\GSp_4$-valued
Weil--Deligne representation is just a $\GSp_4$-valued representation of the
Weil--Deligne group, i.e.\ it is considered up
to~$\GSp_4$-conjugacy). If $\rho$ is a de Rham representation of
$G_K$ on a  $\Qpbar$-vector space~$W$, then we will write $\WD(\rho)$ for the
corresponding Weil--Deligne representation of $W_K$, and if
$\tau:K \into \Qpbar$ is a continuous embedding of fields, then we will
write $\HT_\tau(\rho)$ for the multiset of Hodge--Tate numbers of
$\rho$ with respect to $\tau$, which by definition contains $i$ with
multiplicity
$\dim_{\Qpbar} (W \otimes_{\tau,K} \widehat{\barK}(i))^{G_K} $. Thus,
for example, $\HT_\tau(\varepsilon)=\{ -1\}$.

Let $K/\Q$ be a finite extension. 
If $v$ is a finite place of $K$ we  write $k(v)$ for its residue field, $q_v$ for $\# k(v)$, and $\Frob_v$ for $\Frob_{K_v}$. If $v$ is a real place of $K$, then we will let
$[c_v]$ denote the conjugacy class in $G_K$ consisting of complex conjugations associated to $v$. 

We will frequently adopt the following notation: we let $p>2$ be
prime, and we let $E$ be a finite extension of $\Qp$ with ring of
integers $\cO$, uniformizer~$\lambda$ and residue field $k$.

We will sometimes use the following well-known lemma without comment.

\begin{lemma} \label{lemma:allareconjugate}Let~$\Gamma$ be a group and let~$L$ be an algebraically closed
field. Then a semisimple representation $\Gamma\to\GSp_4(L)$ is
determined up to conjugacy by the composite
$\Gamma\to\GSp_4(L)\to\GL_4(L)\times\GL_1(L)$, where the second factor
records the similitude character.
\end{lemma}
\begin{proof}
This follows
(for example) from the proof of Lemma~6.1 of~\cite{gantakeda}.
\end{proof}

\subsubsection{Galois cohomology}If $L/K$ is an extension of fields,
$k$ is a field, and $V$ is a finite-dimensional $k$-vector space with
an action of~$\Gal(L/K)$, then we write $H^i(L/K,V)$ for
$H^i(\Gal(L/K),V)$, and $h^i(L/K,V)$ for $\dim_k H^i(L/K,V)$. We write
$H^i(K,V)$ and~$h^i(K,V)$ for $H^i(\overline{K}/K,V)$ and~$h^i(\overline{K}/K,V)$ respectively.

\subsubsection{Automorphic representations}We will use the
letter~$\pi$ for automorphic representations of~$\GSp_4$, $\Pi$ for
automorphic representations of~$\GL_n$ (usually with~$n=4$), and
$\boldsymbol{\pi}$ for automorphic representations of~$\GL_2$. We
decorate these in various ways, and aim to be consistent in such
decorations. For example,~$\Pi$ will usually denote the transfer
to~$\GL_4$ of~$\pi$ in the sense of~\S\ref{subsec: Arthur
  classification}, so that for example~$\Pi'_2$ will denote the
transfer of~$\pi'_2$.
\subsection{Induction of two-dimensional representations} \label{section:inductions} We will sometimes want to induce
representations from~$\GL_2$ to~$\GSp_4$. Suppose that~$K/F$ is a
quadratic extension of fields, and that~$r:G_K\to\GL_2(L)$ is a
representation, for some
field~$L$. Choose~$\sigma\in G_F\setminus G_K$, and assume 
that~$\det r$ extends to a character~$\chi$ of~$G_F$
Let~$\rho:=\Ind_{G_K}^{G_F}r:G_F\to\GL_4(L)$. The representation~$\wedge^2\rho$  admits
the characters~$\chi$  and~$\chi \otimes \eta_{F/K}$ as constituents, where~$\eta_{F/K}$ denotes
 the quadratic character. In particular, the representation~$\rho$  generally preserves two symplectic forms, and hence gives rise to two
representations~$\rho_1,\rho_2:G_F\to\GSp_4(L)$
with similitude factors~$\chi$ and $\chi \otimes \eta_{F/K}$ respectively.  To describe these more explicitly,
let~$V$ denote a model for~$r$ so that~$W = V \oplus \sigma V$ is a model for~$\rho$.
Then the Galois action of~$W$ preserves (up to scalar) the symplectic form given by choosing
an arbitrary non-degenerate symplectic form on~$V$, letting~$\sigma V$ and~$V$ be orthogonal, 
and then defining~$\sigma v_1 \wedge \sigma v_2$ consistently to be
either~$\chi(\sigma) v_1 \wedge v_2$ or~$- \chi(\sigma) v_1 \wedge
v_2= \chi \otimes \eta_{F/K}(\sigma) v_1 \wedge v_2$. 
The image of~$(A,B) \in \GL_2(E) \times \GL_2(E)$ with~$\det(A) = \det(B)$ inside~$\GSp_4$ relative to our choice of~$J$ can be given by
$$\left(\begin{matrix} * & 0 & 0 & * \\ 0 & * & * & 0 \\ 0 & * & * & 0 \\ * & 0 & 0 & * \end{matrix} \right) \cap \GSp_4(E).$$

In our applications, it will always be the case that~$\det r$ is the
inverse of the cyclotomic character of~$G_K$, and we will write
simply write $\Ind_{G_K}^{G_F}r$ for the corresponding symplectic
representation with similitude factor the inverse of the cyclotomic
character of~$G_F$. For example, if $K/F$ is a quadratic extension of number
fields, $E$ is an elliptic curve over~$K$, and~$r$ is the dual of the
$p$-adic Tate module of~$E$, then $\Ind_{G_K}^{G_F}r$ is the dual of
the $p$-adic Tate module of the abelian surface~$A = \Res_{K/F}E$,
and the corresponding symplectic structure on this representation
coincides with the one coming from the Weil pairing on~$A$. This is because the representation
on the Tate module of~$A$ is the induction of the corresponding representation on the Tate module
of~$E$, and because the similitude character on the Tate module of an abelian variety is always given by
the
cyclotomic character.

\subsection{The non-archimedean local Langlands
  correspondence}\label{subsec: nonarch LL}Let~$K/\Ql$ be a finite extension for some~$l$. We
will let $\rec_K$ be the local Langlands correspondence of~\cite{ht},
so that if $\pi$ is an irreducible complex admissible representation
of $\GL_n(K)$, then $\rec_K(\pi)$ is a Frobenius semi-simple
Weil--Deligne representation of the Weil group $W_K$. We will write
$\rec$ for $\rec_K$ when the choice of $K$ is clear.

If $(r,N)$ is a Weil--Deligne representation of $W_K$ we will write
$(r,N)^{F-\semis}$ for its Frobenius semisimplification. If $\pi_i$ is
an irreducible smooth representation of $\GL_{n_i}(K)$ for $i=1,2$ we
will write $\pi_1 \boxplus \pi_2$ for the irreducible smooth
representation of $\GL_{n_1+n_2}(K)$ with
$\rec(\pi_1 \boxplus \pi_2)=\rec(\pi_1) \oplus \rec(\pi_2)$.  If
$L/K$ is a finite extension and if $\pi$ is an irreducible smooth
representation of $\GL_n(K)$ we will write $\BC_{L/K}(\pi)$ for the
base change of $\pi$ to $L$ which is characterized by
$\rec_{L}(\BC_{L/K}(\pi))= \rec_K(\pi)|_{W_{L}}$.

We denote the local Langlands correspondence of~\cite{gantakeda} by
$\recGT$; this is a surjective finite-to-one map from the set of
equivalence classes of irreducible smooth complex representations
of~$\GSp_4(K)$ to the set of $\GSp_4$-conjugacy classes of
$\GSp_4(\C)$-valued Weil--Deligne representations of~$W_K$, which we
normalize so that
$\recGT(\pi\otimes(\chi\circ\nu))=\recGT(\pi)\otimes\rec(\chi)$, and
$\nu\circ\recGT(\pi)=\rec(\omega_\pi)$, where~$\omega_\pi$ is the
central character of~$\pi$. 

We fix once and for all for each prime~$p$ an isomorphism
$\imath=\imath_p:\C\cong\Qpbar$. We will generally omit these
isomorphisms from our notation, in order to avoid clutter. In
particular, we will frequently use that~$\imath$ determines a square
root of~$p$ in~$\Qpbar$ (corresponding to the positive square root
of~$p$ in~$\C$). We 
write~$\rec_p$ and~$\recGTp$ for the local Langlands correspondences
for~$\Qpbar$-representations given by conjugating by~$\imath$. These
depend on~$\imath$, but in practice this does not cause us any
difficulty; see Remark~\ref{rem: we do not care about imath for rec}.

\begin{defn}\label{defn: L packet for GSp4}
  If $\rho:G_K\to\GSp_4(\Qpbar)$ is a continuous representation for
  some~$p\ne l$, then we write~$L(\rho)$ for the $L$-packet associated 
  to~$\rho$, which by definition is the set of equivalence classes of
  irreducible smooth $\Qpbar$-representations~$\pi$ of~$\GSp_4(K)$
  with the property that
  $\recGTp(\pi\otimes|\nu|^{-3/2})\cong\WD(\rho)^{F-\semis}$.
\end{defn}
(In accordance with the convention explained above, note
that~$|\nu|^{-3/2}$ makes sense because we have a fixed square root of~$p$.)
\begin{rem}\label{rem: we do not care about imath for rec}
  It is presumably possibly to show that the twist of~$\recGT$ in
  Definition~\ref{defn: L packet for GSp4} (which will be present
  whenever we consider~$\recGTp$) gives
  a local Langlands correspondence for~$\Qpbar$-representations which
  is independent of the choice of~$\imath$, but we have not tried to
  establish this, as we do not need it. We make (implicit) use of this
  for unramified representations, and of the statement that the rank
  of the monodromy operator associated to a representation with
  Iwahori-fixed vectors is independent of the choice of~$\imath$, both
  of which are easily verified explicitly.\end{rem}
\begin{rem}
  We will from now on usually regard automorphic representations as
  being defined over~$\Qpbar$, rather than~$\C$, by means of the fixed
  isomorphism $\imath:\C\cong\Qpbar$. We will not in general draw
  attention to this, and no confusion should arise on the few
  occasions (for example, when considering compatible systems) where
  we think of them as being over~$\C$.
\end{rem}

If $L/K$ is a finite solvable Galois extension of number fields and
if $\pi$ is a cuspidal automorphic representation of $\GL_n(\A_K)$, we
will write $\BC_{L/K}(\pi)$ for its base change to $L$ (which exists
by the main results of~\cite{MR1007299}), an
(isobaric) automorphic representation of $\GL_n(\A_{L})$ satisfying
\[ \BC_{L/K}(\pi)_{w}=\BC_{L_w/K_{v}}(\pi_{v}) \] for all places
$w$ of $L$ where~$v = w|_{K}$ is the restriction of~$w$ to~$K$. If $\pi_i$ is an automorphic representation of
$\GL_{n_i}(\A_K)$ for $i=1,2$ we will write $\pi_1 \boxplus \pi_2$ for
the automorphic representation of $\GL_{n_1+n_2}(\A_K)$ satisfying
\[ (\pi_1 \boxplus \pi_2)_v=\pi_{1,v} \boxplus \pi_{2,v} \]
for all places $v$ of $K$.

If~$(r,N)$ is a Weil--Deligne representation, then we write~$n((r,N))$
for the rank of~$N$. If~$\pi$ is an irreducible admissible
representation of~$\GL_n(K)$ (resp.\ $\GSp_4(K)$), then we
write~$n(\pi)$ for~$n(\rec(\pi))$ (resp.\ $n(\recGT(\pi))$).

\subsection{Local representation theory}\label{subsec: local
  representation theory} 

In this section, we recall a number of more or less well-known results
about the representation theory of~$\GSp_4(K)$, where~$K$ is a local
field of characteristic zero. Some of these
results are in~\cite{MR2234862}, but for convenience we have gathered
them all together here, and have usually given proofs. Since our
applications of this material are all global, and some of the definitions we make (such as the
normalizations of Hecke operators at places dividing~$p$) depend on
global information, we have chosen to work in the same global setting
that we consider in the rest of the paper.

Let~$p>2$ be prime, and let~$F$ be a totally real field in which $p$
splits completely. Let~$E/\Qp$ be a finite extension with ring of
integers~$\cO$ and residue field~$k$. Let $v$ be a finite place of $F$, and fix a uniformizer
$\varpi_v\in\cO_{F_v}$. For most of this section, we will allow~$v$ to
divide~$p$, although at the end of the section, we will prove some
results (which follow those of~\cite{KT} for~$\GL_n$) under the
assumption that $q_v\equiv 1\pmod{p}$. We fix once and for all a
square root~$q_v^{1/2}\in E$.
\subsubsection{Generalities} We begin by recalling some results on
Iwahori Hecke algebras. It costs us nothing to recall these in a more
general setting, so we temporarily let $G/\cO_{F_v}$ be a split
reductive group with $T\subset B=T\cdot U$ a maximal torus and Borel
(with unipotent radical~$U$),
and let $N$ be the normalizer of $T$ in~ $G$.  Let $W=N(F_v)/T(F_v)$ be
the Weyl group.  Let $\Delta\subset X^*(T)$ be the simple roots.  We write $\widetilde{W}=N(F_v)/T(\cO_{F_v})$ for the extended affine Weyl group.

Let $\Iw(v)=\ker(G(\cO_{F_v})\to B(k(v)))$ be an Iwahori subgroup, and
let $\Iw_1(v)=\ker(G(\cO_{F_v})\to U(k(v)))$ be a pro-$v$ Iwahori
subgroup.  Let \[\mathcal{H}_1=\cH_1(v)=\cO[\Iw_1(v)\backslash G(F_v)/\Iw_1(v)] = \cO[G(F_v) \doubleslash \Iw_1(v)]\] 
be the pro-$v$ Iwahori Hecke algebra. (Here~$G \doubleslash K$ denotes~$K \backslash G/K$ --- we tend
to prefer the first notation but we also sometimes use the second notation since it is more compact and some of our expressions
are already typographically somewhat complicated.)

We let $T(\cO_{F_v})_1=(\ker T(\cO_{F_v})\to T(k(v)))$.  We also let
\begin{equation*} T(F_v)^+ =\{x\in T(F_v)\mid
\alpha(x)\in\cO_{F_v},\forall\alpha\in\Delta\}.
\end{equation*} For $g\in G(F_v)$, we write
$[\Iw_1(v)g\Iw_1(v)]\in\cH_1$ for the characteristic function of the
double coset $\Iw_1(v)g\Iw_1(v)$.

\begin{prop}\label{prop: invertibility of Iwahori Hecke operators} For $x,y\in T(F_v)^+$, we have
\begin{equation*}
[\Iw_1(v)x\Iw_1(v)]\cdot[\Iw_1(v)y\Iw_1(v)]=[\Iw_1(v)xy\Iw_1(v)]
\end{equation*} and moreover $[\Iw_1(v)x\Iw_1(v)]\in
(\cH_1[1/p])^\times$. If $v\nmid p$, then in fact $[\Iw_1(v)x\Iw_1(v)]\in
\cH_1^\times$.
\end{prop}\begin{proof}
The first statement is a special case of~\cite[Lem.\
4.1.5]{CasselmanNotes}, while the rest is immediate from~\cite[Cor.\ 1]{MR2122539}.
\end{proof}

As a result, there is a homomorphism
\begin{equation*} T(F_v)\to (\cH_1[1/p])^\times
\end{equation*} which is defined as follows: write $x\in T(F_v)$ as
$x=yz^{-1}$ with $y,z\in T(F_v)^+$ and send $x$ to
\begin{equation*}
(\delta^{1/2}_B(y)[\Iw_1(v)y\Iw_1(v)])(\delta^{1/2}_B(z)[\Iw_1(v)z\Iw_1(v)])^{-1}
\end{equation*} 
where~$\delta_{B}$ is the modulus character. The kernel of this homomorphism is
$T(\cO_{F_v})_1$. If $v\nmid p$, then the image of the homomorphism is
in $\cH_1^\times$. 

\begin{prop}\label{prop: Jacquet module iso} Let $\pi$ be a smooth admissible $\barE[G(F_v)]$-module.
Then the map $\pi\to \pi_U$, where $\pi_U$ is the {\em (}normalized{\em)} Jacquet
module, induces an isomorphism of $\barE[T(F_v)]$-modules
\begin{equation*} \pi^{\Iw_1(v)}\to (\pi_U)^{T(\cO_{F_v})_1}.
\end{equation*}
\end{prop}\begin{proof}
  By~\cite[Lem.\ 4.1.1]{CasselmanNotes} (noting that the Jacquet
  module in this reference is \emph{not} the normalized Jacquet
  module), the map $\pi\to\pi_U$ induces an $\barE[T(F_v)]$-module
  homomorphism $\pi^{\Iw_1(v)}\to (\pi_U)^{T(\cO_{F_v})_1}$. It is an
  isomorphism by~\cite[Prop.\ 4.1.4]{CasselmanNotes} and
  Proposition~\ref{prop: invertibility of Iwahori Hecke operators}.
\end{proof}

For a character $\chi:T(F_v)\to\barE^\times$, write
$\pi(\chi)=\nInd_{B(F_v)}^{G(F_v)}\chi$ for the corresponding
principal series representation.  Then we recall
\begin{prop}\label{prop: Jacquet module} For $\chi:T(F_v)\to\overline{E}^\times$ there is an
isomorphism of $\barE[T(F_v)]$-modules
\begin{equation*} (\pi(\chi)_U)^{\mathrm{ss}}\simeq \bigoplus_{w\in W}
\barE(w\cdot\chi).
\end{equation*}
\end{prop}\begin{proof}
  This is a special case of~\cite[Thm.\ 6.3.5]{CasselmanNotes}.
\end{proof}

We say that $\pi(\chi)$ is a \emph{tame principal series} if $\chi$ is
trivial on $T(\cO_{F_v})_1$ and an \emph{unramified principal series} if
$\chi$ is trivial on $T(\cO_{F_v})$.  The results recalled above
immediately imply the well-known facts that if $\pi$ is an irreducible
smooth $\barE[G(F_v)]$-module, then $\pi^{\Iw_1(v)}\not=\{0\}$ if and
only if $\pi$ is a constituent of a tame principal series, and
$\pi^{\Iw(v)}\not=\{0\}$ if and only if $\pi$ is a constituent of an
unramified principal series. 

Write~$\cH:=\cO[\Iw(v)\backslash
G(F_v)/\Iw(v)]$ for the Iwahori Hecke algebra. This enjoys similar
properties to those of~$\cH_1$ recalled above; in particular, the
analogue of Proposition~\ref{prop: invertibility of Iwahori Hecke
  operators} gives an embedding  $E[X_*(T)]\into\cH[1/p]$, and if
$v\nmid p$, then this restricts to an embedding $\cO[X_*(T)]\into\cH$.

\subsubsection{Principal series for $\GSp_4$.}  We now specialize our
discussion to $G=\GSp_4$. We recall some known results on constituents
of unramified principal series representations; many of these results
are originally due to~\cite{MR1212952}, but for convenience we refer
to the tables in~\cite[App.\ A]{MR2344630}. (Note that the
compatibility of the proposed Langlands parameters in~\cite[App.\
A.5]{MR2344630} with the correspondence~$\recGT$ is proved
in~\cite[Prop.\ 13.1]{MR2846304}.)  

If $\chi_1,\chi_2,\sigma$ are characters of~$F_v^\times$, then we
write \[\chi_1\times\chi_2\rtimes\sigma:=\nInd_{B(F_v)}^{\GSp_4(F_v)}\chi_1\otimes\chi_2\otimes\sigma,\]where \[\chi_1\otimes\chi_2\otimes\sigma:
  \begin{pmatrix} a&*&*&*\\&b&*&*\\&&cb^{-1}&*\\&&&ca^{-1}
  \end{pmatrix}\mapsto\chi_1(a)\chi_2(b)\sigma(c).
\]

\begin{prop}\label{prop: rec of principal series} \leavevmode
\begin{enumerate}
\item $\chi_1\times\chi_2\rtimes\sigma$ is irreducible if and only if
none of~$\chi_1,\chi_2,\chi_1\chi_2^{\pm 1}$ is equal
to~$|\cdot|_v^{\pm 1}$.
\item If $\pi$ is an irreducible constituent of
$\chi_1\times\chi_2\rtimes\sigma$,
then \[\recGTp(\pi)^{\mathrm{ss}}=\sigma\circ\Art_{F_v}^{-1}\otimes\left((\chi_1\chi_2)\circ\Art_{F_v}^{-1}\oplus
\chi_1\circ\Art_{F_v}^{-1}\oplus \chi_2\circ\Art_{F_v}^{-1}\oplus 1
\right).\]
\item If $\chi_1\times\chi_2\rtimes\sigma$ is
irreducible, then $\recGTp(\chi_1\times\chi_2\rtimes\sigma)$ is 
semisimple \emph{(}that is, $N = 0$\emph{)}.
\end{enumerate} 
\end{prop}\begin{proof}Part~(1) is~\cite[Lem.\ 3.2]{MR1212952}. Parts~(2) and~(3) follow immediately from rows 
I--VI of~\cite[Table A.7]{MR2344630}.
\end{proof}\subsubsection{Spherical Hecke operators}\label{sssection spherical hecke} Define matrices
\[\beta_{v,0}=\diag(\varpi_v,\varpi_v,\varpi_v,\varpi_v),\]
\[\beta_{v,1}=\diag(\varpi_v,\varpi_v,1,1),\] \[\beta_{v,2}=\diag(\varpi_v^2,\varpi_v,\varpi_v,1).\]

We have the spherical Hecke operators
$T_{v,i}=[\GSp_4(\cO_{F_v})\beta_{v,i}\GSp_4(\cO_{F_v})]$, which are
independent of~$\varpi_v$. It is easy to check (using Proposition~\ref{prop: rec of principal series}~(2)) that
if~$\pi$ is an unramified representation of~$\GSp_4(F_v)$ (that is, if
$\pi^{\GSp_4(\cO_{F_v})}\ne 0$, so that~$\pi$ is a constituent of an
unramified principal series), then the
characteristic polynomial of~$\recGTp(\pi\otimes|\nu|^{-3/2})(\Frob_v)$
is \numequation\label{eqn: char poly for unramified Hecke}Q_v(X):=X^4-t_{v,1}X^3+(q_vt_{v,2}+(q_v^3+q_v)t_{v,0})X^2-q_v^3t_{v,0}t_{v,1}X+q_v^6t_{v,0}^2,\end{equation}where
we are writing~$t_{v,i}$ for the eigenvalue of the operator~$T_{v,i}$
on~$\pi^{\GSp_4(\cO_{F_v})}$. 

\begin{df}
We say that the \emph{Hecke parameters} of~$\pi$ are the roots
of~$Q_v(X)$, ordered in such a way that the pairs of roots~$(1,4)$ and~$(2,3)$ both multiply
to give the value~$\gamma_v$ of the similitude character evaluated on~$\Frob_v$.
We  write these Hecke parameters as
$[\alpha_v,\beta_v,\gamma_v\beta_v^{-1},\gamma_v\alpha_v^{-1}]$, where implicitly we view these terms as labelling the vertices of
a square:
$$
\begin{diagram}
\alpha_v & \rLine & \beta_v \\
\dLine & & \dLine \\
\gamma_v \beta^{-1}_v & \rLine & \gamma_v \alpha^{-1}_v \end{diagram}
$$
and the ordering is unique up to the action of the Weyl group~$D_8 = \Sym(\square)$.
In particular, the data of the quadruple~$[\alpha_v,\beta_v,\gamma_v\beta_v^{-1},\gamma_v\alpha_v^{-1}]$ carries with
it the value of the similitude character.
\end{df}

We will  be concerned with the case that the central
character of~$\pi$ is given by $a\mapsto |a|^{2}$, in which case the Hecke parameters have
the form~$[\alpha_v,\beta_v,q_v\beta_v^{-1},q_v\alpha_v^{-1}]$.

\subsubsection{Iwahori Hecke operators}\begin{defn}\label{defn: general unram PS}
  We say that an unramified principal series $\pi(\chi)$ is \emph{general}
  if the Hecke parameters are pairwise distinct and no ratio of them is
  $q_v$. In particular, $\pi(\chi)$ is irreducible, and $|W\cdot\chi|=8$. \end{defn}

We have Iwahori Hecke operators $\Unaive_{\Iw(v),i}=[\Iw(v)\beta_{v,i}\Iw(v)]$. The notation
``$\Unaive$'' is intended to indicated that we have not yet
appropriately normalized these operators, as we will shortly do in the
case that~$v|p$.
Then we have \begin{prop}\label{prop: general iwahori up eigenvectors} Let $\pi$ be a general unramified principal series with
Hecke parameters
$[\alpha_v,\beta_v,q_v\beta_v^{-1},q_v\alpha_v^{-1}]$.  Then
$\pi^{\Iw(v)}$ is a direct sum of 8 one-dimensional simultaneous eigenspaces
for the $\Unaive_{\Iw(v),i}$.  For a given (ordered) choice of $\alpha_v$ and
$\beta_v$ the corresponding
eigenvalues are $u_{v,0}=q_v^{-2}$, $u_{v,1}=\alpha_v$, and
$u_{v,2}=q_v^{-1}\alpha_v\beta_v$.
\end{prop}\begin{proof}The first part is immediate from Propositions~\ref{prop:
Jacquet module iso} and~\ref{prop:
Jacquet module}.  To compute the eigenvalues, by the definition of the Hecke parameters
and Proposition \ref{prop: rec of principal series} we have
$\alpha_v=q_v^{3/2}(\chi_1\chi_2\sigma)(\varpi_v)$, $\beta_v=q_v^{3/2}(\chi_1\sigma)(\varpi_v)$ and $q_v=q_v^3(\chi_1\chi_2\sigma^2)(\varpi_v)$.
 We then have
  $u_{v,i}=\delta_B(\beta_{v,i})^{-1/2}(\chi_1\otimes\chi_2\otimes\sigma)(\beta_{v,i})$, so that
  $u_{v,0}=(\chi_1\chi_2\sigma)^2(\varpi_v)=q_v^{-2}$,
  $u_{v,1}=q_v^{3/2}(\chi_1\chi_2\sigma)(\varpi_v)=\alpha_v$,
  $u_{v,2}=q_v^{2}(\chi_1^2\chi_2\sigma^2)(\varpi_v)=q_v^{-1}\alpha_v\beta_v$,  as
  required.
\end{proof}Proposition \ref{prop: general iwahori up eigenvectors} has the following converse:
\begin{prop}\label{prop: iwahori up eigenvectors converse}
Let $\pi$ be an irreducible admissible representation of $\GSp_4(F_v)$, and suppose that $\pi^{\Iw(v)}$ contains an eigenvector for the $\Unaive_{\Iw(v),i}$ with eigenvalues $u_{v_i}$ satisfying $u_{v_0}=q_v^{-2}$, $u_{v,1}=\alpha_v$ and $u_{v,2}=q_v^{-1}\alpha_v\beta_v$ such that no ratio of a pair of $[\alpha_v,\beta_v,q_v\beta_v^{-1},q_v\alpha_v^{-1}]$ is $q_v$.  Then $\pi$ is the unramified principal series with Hecke parameters $[\alpha_v,\beta_v,q_v\beta_v^{-1},q_v\alpha_v^{-1}]$.
\end{prop}
\begin{proof}
  Reversing the calculation in the previous proof, we
  let~$\chi=\chi_1\otimes\chi_2\otimes\sigma$ be the unramified
  character with $\chi_1(\varpi_v)=\alpha_v\beta_vq_v^{-1}$,
  $\chi_2(\varpi_v)=\alpha_v\beta_v^{-1}$, and
  $\sigma(\varpi_v)=\alpha_v^{-1}q_v^{-1/2}$.  
     We see that there is an inequality
  $\Hom_{T(F_v)}(\pi^{\Iw(v)},\chi)\not=\{0\}$, and hence
  $\Hom(\pi,\pi(\chi))\not=\{0\}$ by Proposition~\ref{prop: Jacquet
    module iso} and Frobenius reciprocity.  Finally, by
  Proposition~\ref{prop: rec of principal series}, $\pi(\chi)$ is also
  irreducible.
\end{proof}

\subsubsection{Parahoric level Hecke operators for $\mathrm{GL}_2$}We will also need to consider certain  parahoric Hecke algebra and investigate how they relate to the Iwahori Hecke algebra. 

We begin by recalling some standard results for the
group $\mathrm{GL}_2$. We  let $\Iw(v)' \subset \mathrm{GL}_2(\ocal_{F_v})$ be the Iwahori subgroup  of matrices which are upper triangular modulo $\varpi_v$ (we put a prime because $\Iw(v)$ is used to denote the Iwahori subgroup in $\mathrm{GSp}_4(\ocal_{F_v})$). 

We introduce the following operators in the spherical Hecke algebra $\mathcal{H}_{\Sph}[1/p]$: 

\begin{enumerate}
\item  $T_{v,1}^{\GL_2} =  [\mathrm{GL}_2(\ocal_{F_v})   \begin{pmatrix} \varpi_v & 0 \\
      0& 1\\
   \end{pmatrix} \mathrm{GL}_2(\ocal_{F_v})  ] $, 
   
      \item  $T_{v,0}^{\GL_2} =   [\mathrm{GL}_2(\ocal_{F_v})   \begin{pmatrix} \varpi_v & 0 \\
      0& \varpi_v\\
   \end{pmatrix} \mathrm{GL}_2(\ocal_{F_v})  ] $.
   \end{enumerate}
   
   We also define the following operators in the Iwahori Hecke algebra $\mathcal{H}_{\Iw(v)'}[1/p]$: 
   \begin{enumerate}
\item  $U_{v,1}^{\GL_2} =  [\Iw(v)'   \begin{pmatrix} \varpi_v & 0 \\
      0& 1\\
   \end{pmatrix} \Iw(v)'  ] $, 
   \item  $U_{v,0}^{\GL_2} =   [\Iw(v)'   \begin{pmatrix} \varpi_v & 0 \\
      0& \varpi_v\\
   \end{pmatrix} \Iw(v)'  ] $,
   
   \item $e_{\mathrm{Sph}}^{\GL_2} = [\mathrm{GL}_2(\ocal_{F_v})  ] $.
   \end{enumerate}
   
   For any element $f$ of the centre of the Iwahori Hecke algebra, the element $e_{\mathrm{Sph}}^{\GL_2}  f$ defines an element of the spherical Hecke algebra. 
   
\begin{lem}\label{lem-GL2-Iwahori}  The centre  $Z(\mathcal{H}_{\Iw(v)'}[1/p])$ of the Iwahori Hecke algebra is generated  by $U_{v,0}^{\GL_2}$ and $q_vU^{\GL_2}_{v,0}( U^{\GL_2}_{v,1})^{-1} + U^{\GL_2}_{v,1}$, the map $e_{\mathrm{Sph}}^{\GL_2}:  Z(\mathcal{H}_{\Iw(v)'}[1/p]) \rightarrow \mathcal{H}_{\Sph}[1/p] $ is an isomorphism and we have  the following identities: 

\begin{enumerate}
\item $e_{\mathrm{Sph}}^{\GL_2}  U^{\GL_2}_{v,0}  = T^{\GL_2}_{v,0}$,
\item $e_{\mathrm{Sph}}^{\GL_2}  (q_vU^{\GL_2}_{v,0} (U^{\GL_2}_{v,1})^{-1} + U^{\GL_2}_{v,1})  =
  T_{v,1}$. \end{enumerate}

\end{lem}

\begin{proof} This follows from~\cite[\S1, \S2, \S4.6]{MR2642451}.\end{proof}
\subsubsection{Klingen level Hecke operators}

We have Klingen Hecke operators
$\Unaive_{\Kli(v),i}=[\Kli(v)\beta_{v,i}\Kli(v)]$.  

\begin{prop}\label{prop: general Klingen up eigenvectors}
Let $\pi$ be a general unramified principal series with
Hecke parameters
$[\alpha_v,\beta_v,q_v\beta_v^{-1},q_v\alpha_v^{-1}]$.  Then
$\pi^{\Kli(v)}$ is a direct sum of 4 one-dimensional simultaneous eigenspaces
for the $\Unaive_{\Kli(v),i}$.  For a given choice of $\{\alpha_v,\beta_v\}$,
the eigenvalues are $u_{v,0}=q_v^{-2}$, $u_{v,1}=\alpha_v+\beta_v$,
and $u_{v,2}=q_v^{-1}\alpha_v\beta_v$.
\end{prop}
\begin{proof} This follows from a direct computation, see~\cite[Prop.\ 3.2.1,
  Cor.\ 3.2.2]{MR2234862}.
\end{proof}

\begin{rem}\label{remfromrio} We sketch another (related) proof of Proposition~\ref{prop: general Klingen up eigenvectors}.   Let us denote by $\mathcal{H}_{\Iw(v)}[1/p]$ the Iwahori Hecke algebra and by $Z_{\Kli(v)} (\mathcal{H}_{\Iw(v)}[1/p])$ the sub-algebra generated by $\Unaive_{\Iw(v),1} + q_v(\Unaive_{\Iw(v),1})^{-1} \Unaive_{\Iw(v),2},  \Unaive_{\Iw(v),2}, \Unaive_{\Iw(v),0}$. One checks that $Z_{\Kli(v)} (\mathcal{H}_{\Iw(v)}[1/p])$ commutes with $e_{\Kli(v)} = [\Kli(v)]$ by using Bernstein's relation (\cite[\S1.15]{MR2642451}). Therefore we get a map: 
$e_{\Kli(v)}: Z_{\Kli(v)} (\mathcal{H}_{\Iw(v)}[1/p]) \rightarrow \mathcal{H}_{\Kli(v)}[1/p]$ where $\mathcal{H}_{\Kli(v)}[1/p]$ is the Klingen Hecke algebra. 
We claim that:

\begin{itemize}\item $e_{\Kli(v)}  (\Unaive_{\Iw(v),1} + q_v(\Unaive_{\Iw(v),1})^{-1} \Unaive_{\Iw(v),2}) = \Unaive_{\Kli(v), 1}$,
\item  $ e_{\Kli(v)}  \Unaive_{\Iw(v),2} = \Unaive_{\Kli(v), 2}$,
\item $e_{\Kli(v)}  \Unaive_{\Iw(v),0} = \Unaive_{\Kli(v), 0}$.
\end{itemize}

The claim can  be checked after restricting all these functions to the
Levi $\mathrm{GL}_2 \times \mathrm{GL}_1$ of the Klingen parabolic by
\cite[Prop.\ II.5]{Vigneras1998}, so it follows from
Lemma~\ref{lem-GL2-Iwahori}. The result then follows from
Proposition~\ref{prop: general iwahori up eigenvectors}.\end{rem}

 \begin{rem}
  \label{rem: could also use Jacquet modules}
  Proposition~\ref{prop: general Klingen up eigenvectors} could also be proved using Jacquet
  modules (as could analogous results for invariants at other level
  structures which admit parahoric factorizations).
\end{rem}

\begin{prop}\label{prop: Klingen up eigenvectors converse}
Let $\pi$ be an irreducible admissible representation of $\GSp_4(F_v)$, and suppose that $\pi^{\Kli(v)}$ contains an eigenvector for the $\Unaive_{\Kli(v),i}$ with eigenvalues $u_{v_i}$ satisfying $u_{v_0}=q_v^{-2}$, $u_{v,1}=\alpha_v+\beta_v$ and $u_{v,2}=q_v^{-1}\alpha_v\beta_v$ such that no ratio of a pair of~$\{\alpha_v,\beta_v,q_v\beta_v^{-1},q_v\alpha_v^{-1}\}$ is $q_v$.  Then $\pi$ is the unramified principal series with Hecke parameters~$[\alpha_v,\beta_v,q_v\beta_v^{-1},q_v\alpha_v^{-1}]$.
\end{prop}
\begin{proof}   As in the proof of Proposition~\ref{prop: iwahori up eigenvectors
    converse},  we deduce from~$\pi^{\Iw(v)} \ne 0$ that~$\pi$ is a constituent of an unramified principal series representation.  The central character of such a constituent is unramified and so is determined by the value on Frobenius.
   From the equation~$u_{v_0} = q^{-2}_v$, we deduce that the central character of~$\pi$ is~$| \cdot |^2$, and
   hence the central character of~$\pi\otimes|\nu|^{-3/2}$ is~$| \cdot |^{-1}$, and hence that the similitude character 
   of~$\recGTp(\pi\otimes|\nu|^{-3/2})$ is the inverse of the cyclotomic character~$\varepsilon^{-1}$. In particular, the value
   of the similitude character of the Weil--Deligne representation on~$\Frob_v$ is~$q_v$, and thus~$\pi$ is  
a constituent of an
  unramified principal series representation with Hecke parameters
  $[\alpha'_v,\beta'_v,q_v(\beta'_v)^{-1}, q_v(\alpha'_v)^{-1}]$.
  (Note that the ordering of these eigenvalues above is determined up to the action of~$D_8$.)
  Comparing to  Proposition~\ref{prop: general Klingen up
    eigenvectors}, without loss of generality, we may rearrange the Hecke parameters of~$\pi$ so that we deduce the two equations
    $$\alpha'_v + \beta'_v = \alpha_v + \beta_v, \qquad \alpha'_v \beta'_v = \alpha_v \beta_v,$$
    and thus (again up to reordering) $\alpha'_v = \alpha_v$ and~$\beta'_v = \beta_v$.
    By  Proposition~\ref{prop: rec of principal series}, the principal series~$\pi$ is irreducible.
 \end{proof}

\subsubsection{Generic unipotent representations}

We say that a $\GSp_4(\overline{E})$-valued Weil--Deligne
representation~$r$ is \emph{generic} if
$\ad(r)(1)$ has no invariants, and is
\emph{unipotent} if $r^{\mathrm{ss}}$ is unramified. 
\begin{prop}\label{prop: unipotent generic implies generic} Let $r$ be
  unipotent. Then  the $L$-packet corresponding to~$r$ contains a generic
representation if and only if~$r$ is generic.\end{prop}\begin{proof}By the main theorem of~\cite{gantakeda} (part vii), the $L$-packet $L(r)$ contains a generic
representation if and only if the adjoint $L$-factor
$L(s,\ad(r^{F-\semis}))$ is holomorphic at $s=1$, which, by definition, is easily seen to be equivalent to the statement that~$\ad(r^{F-\semis})(1)$ has no invariants.
 Thus we are reduced to
checking that~$r$ is generic if and only if~$r^{F-\semis}$ is
generic. Let~$W$ denote the vector space underlying the representation~$\ad(r)(1)$.
We are reduced to showing that~$\Hom(\overline{E},W) = 0$
if and only if~$\Hom(\overline{E},W^{F-\semis}) = 0$.

One implication is trivial. For the reverse implication, a map from~$\overline{E}$ to~$W^{F-\semis}$ is
the same as giving a vector~$x$ in~$W$ which lies in the kernel of~$N$
and is a generalized eigenvector for the Frobenius~$\phi$
with eigenvalue~$1$. For a suitable choice of~$n \in \mathbf{N}$, the vector~$y = (\phi - 1)^n x$ will be non-zero and a genuine  eigenvector
for~$\phi$ with eigenvalue one. On the other hand, since~$x$ lies in the kernel of~$N$, so does~$\phi x$,  because~$N \phi x = q^{-1}_v \phi Nx  = 0$.
Similarly,  any polynomial in~$\phi$ applied to~$x$ also lies in the kernel of~$N$.
Thus~$y$ also lies in the kernel of~$N$ and gives rise to a nonzero
element of~$\Hom(\overline{E},W)$.
\end{proof}

\subsubsection{Normalized Hecke operators, ordinary representations,
  and ordinary projectors}\label{subsubsec: ordinary Hecke operators}

In this section, we assume that $v|p$.  We fix integers $k\geq l\geq
2$, and $k\equiv l\pmod{2}$ (these will correspond to the weights of
our automorphic forms; see Section~\ref{subsec: archimedean L
packets}). Then we will consider normalized Hecke operators at Iwahori and
Klingen level defined by\begin{align*}
\normU_{\Iw(v),0}&=p^2\Unaive_{\Iw(v),0}\quad& \normU_{\Kli(v),0}&=p^2\Unaive_{\Kli(v),0}\quad\\
\normU_{\Iw(v),1}&=p^{(k+l)/2-2}\Unaive_{\Iw(v),1}\quad& \normU_{\Kli(v),1}&=p^{(k+l)/2-2}\Unaive_{\Kli(v),1}\quad\\
\normU_{\Iw(v),2}&=p^{k-1}\Unaive_{\Iw(v),2}\quad& \normU_{\Kli(v),2}&=p^{k-1}\Unaive_{\Kli(v),2}\quad
\end{align*}We will often write
$U_{v,i}$ for the operators~$\normU_{\Iw(v),i}$ when the context is clear. We will
also keep writing~$U_{v,0}$ for the Hecke operator $p^2[K_v\beta_{v,0} K_v]$ for
any subgroup~$K_v$ of~$\Iw(v)$ (because $U_{v,0}$ lies in the centre
of the Iwahori Hecke algebra and  therefore $p^2[K_v\beta_{v,0} K_v] =
e_{K_v}  U_{v,0}$). We will also often write ~$U_{v,2}$
for the Hecke operator $U_{\Kli(v),2}$ for the same reason (see
Remark~ \ref{remfromrio}). We can and do also normalize the Siegel
Hecke operators in the same way, so that for example $U_{\Si(v),1}=p^{(k+l)/2-2}\Unaive_{\Si(v),1}$.

An irreducible smooth $\overline{E}[\GSp_4(F_v)]$-representation with
central character $|\cdot|^2$ is said to be \emph{ordinary} of weights
$k\geq l\geq 2$ if there exists an eigenvector $v\in\pi^{\Iw(v)}$ for
$\normU_{\Iw(v),i}$ with eigenvalues $u_{v,i}$ with
$v_p(u_{v,i})=0$.  If $\alpha_v$ and $\beta_v$ are defined by
$\alpha_v=u_{v,1}$, $\beta_v=u_{v,2}/u_{v,1}$, then, by Proposition
\ref{prop: Jacquet module}, $\pi$ is a constituent of an unramified
principal series with Hecke parameters \[[\alpha_v p^{2-(k+l)/2},\beta_vp^{-(k-l)/2},\beta_v^{-1}p^{1+(k-l)/2},\alpha_v^{-1}p^{(k+l)/2-1}].\]  We say that $\pi$ is \emph{$p$-distinguished} if these four Satake parameters are pairwise distinct, or in other words if either $l>2$ or $\alpha_v\not=\beta_v$.

If $l>2$, then again by Proposition \ref{prop: Jacquet module},
$v\in\pi^{\Iw(v)}$ is the unique eigenvector (up to scale) with unit
eigenvalues for the $\normU_{\Iw(v),i}$.  In this case, the ordered
pair $(\alpha_v,\beta_v)$ is uniquely determined by $\pi$, and we
call~$(\alpha_v,\beta_v)$ the \emph{ordinary Hecke parameters} of
$\pi$.  If $l=2$ and $\pi$ is $p$-distinguished, then there may also
be an eigenvector $v'\in\pi^{\Iw(v)}$ with unit eigenvalues
$\normU_{\Iw(v),1}v'=\beta_vv'$, $\normU_{\Iw(v),2}v'=\alpha_v\beta_vv'$
(we will see below that in fact such a $v'$ always exists.)  Thus at
least the set $\{\alpha_v,\beta_v\}$ is determined by $\pi$ and we again call them the \emph{ordinary Hecke parameters} of $\pi$.

We let~$e_{\text{reg}}$ be the ordinary projector (in the sense of Section~\ref{subsec: projectors}) associated
to~$\normU_{\Iw(v),1}\normU_{\Iw(v),2}$, and
let~$e_{{\text{irreg}}}$ be the ordinary projector associated
to~$\normU_{\Kli(v),2}$. 
\begin{prop}\label{prop: ordclass}
Let $\pi$ be an ordinary $p$-distinguished representation of weights
$k\geq l\geq 2$, with ordinary Hecke parameters $(\alpha_v,\beta_v)$
\emph{(}or~$\{\alpha_v,\beta_v\}$ if~$l=2$\emph{)}.  Assume that either $k>l>3$ or $l=2$.
\begin{enumerate}
\item If $k>l>3$ or if $l=2$ and $k>2$, then $\pi$ is an irreducible
  principal series. \item If $k=l=2$, then in the sense of the tables of \cite[\S1]{schmidt-tables}, $\pi$ is a representation of type  \emph{Va} if $\{\alpha_v,\beta_v\}=\{1,-1\}$, \emph{IIIa}  if $\alpha_v\beta_v=1$, \emph{IIa} if $\#\{\alpha_v,\beta_v\}\cap\{1,-1\}=1$, or otherwise is an irreducible unramified principal series.
\end{enumerate}

In all cases, $\pi$ is generic and the $L$-packet $L(\pi)$ of $\pi$
contains no other ordinary representations. Moreover: 
\begin{enumerate}
\item $k>l>3$ then\begin{equation*}
\dim e_{\text{reg}}\pi^{\Iw(v)}=1
\end{equation*}
on which $\normU_{\Iw(v),i}$ has eigenvalues $1,\alpha_v,\alpha_v\beta_v$ for $i=0,1,2$.
\item If $l=2$ then
\begin{equation*}
\dim e_{\text{reg}}\pi^{\Iw(v)}=2
\end{equation*}
and there are two eigenspaces for $\normU_{\Iw(v),i}$, with
eigenvalues $1,\alpha_v,\alpha_v\beta_v$ and
$1,\beta_v,\alpha_v\beta_v$ respectively, and moreover
\begin{equation*}
\dim e_{\text{irreg}}\pi^{\Kli(v)}=1
\end{equation*}
with $\normU_{\Kli(v),i}$ eigenvalues $1,\alpha_v+\beta_v,\alpha_v\beta_v$, for $i=0,1,2$.
\end{enumerate}

\end{prop}
\begin{proof}
As remarked above, by Proposition~\ref{prop: Jacquet module}, $\pi$ is a constituent of an unramified principal series with Hecke parameters \[[\alpha_vp^{2-(k+l)/2},\beta_vp^{-(k-l)/2},\beta_v^{-1}p^{1+(k-l)/2},\alpha_v^{-1}p^{(k+l)/2-1}].\]  If either $k>l>3$ or $l=2$ and $k>2$, no ratio of a pair of these parameters can be $p$, and hence $\pi$ is an irreducible principal series by Proposition \ref{prop: rec of principal series}.

In the remaining case, $k=l=2$,  the Satake parameters are
$[\alpha_v,\beta_v,\beta_v^{-1}p,\alpha_v^{-1}p]$, and the
corresponding principal series may be reducible when one of
$\alpha_v^2,\beta_v^2,\alpha_v\beta_v$ is equal to~$1$.  The constituents of these principal series are listed in the tables \cite[\S1]{schmidt-tables}.  The case that either $\alpha_v^2=1$ or $\beta_v^2=1$ but not both corresponds to type II, the case that $\alpha_v^2=\beta_v^2=1$ corresponds to type V, and the case that $\alpha_v\beta_v=1$ corresponds to type III.

For each constituent $\pi$ of such a principal series, the tables give
a computation of the Jacquet module $\pi_U^{\text{ss}}$, which is
equal to $\pi_U$ because $\alpha_v\not=\beta_v$.  This allows us, by Proposition \ref{prop: Jacquet module iso}, to determine the simultaneous eigenvalues of the $U_{\Iw(v),i}$ on $\pi^{\Iw(v)}$.  At this point the result follows from an inspection of the tables.
\end{proof}

We now turn to the global situation. Recall that
we have fixed an isomorphism~$\imath:\Qpbar\cong\C$, so that in
particular in the following definition we can and do identify the
infinite places of~$F$ with the places dividing~$p$. See Section~\ref{subsec: archimedean L
packets} for our conventions regarding the weights of automorphic representations.

\begin{defn}
\label{defn: imath ordinary}Let~$\pi$ be a cuspidal automorphic
representation of~$\GSp_4(\A_F)$ with central character~$|\cdot|^2$
and weight~$(k_v,l_v)_{v|\infty}$, where $k_v\ge l_v\ge 2$ and
$k_v\equiv l_v\pmod{2}$ for all
$v|\infty$. Then we say that~$\pi$ is \emph{ordinary} if for each
place~$v|p$, $\pi_v$ is ordinary of weights~$k_v\ge l_v\ge 2$.
\end{defn}

The following proposition will be useful for going between ordinary
$p$-adic modular forms and ordinary automorphic
representations. For each subset $I\subset S_p$ we set \[K_p(I) = \prod_{v \in I} \Kli(v) \prod_{v \in I^c} \Iw(v).\]
We also let $e(I)=\prod_{v\in I}e_{\text{irreg}}\prod_{v\not\in
  I}e_{\text{reg}}$.

\begin{prop}\label{prop: ordinary eigenform in autrep}Let $\pi$ be a cuspidal automorphic representation of $\GSp_4(\A_F)$
of weight $(k_v,l_v)_{v|\infty}$ with $k_v\geq
l_v\geq 2$ and with central character $|\cdot|^2$, and fix tuples of $p$-adic units
$(\alpha_v,\beta_v)_{v|p}$.  Assume that for each $v\in S_p$, either $k_v>l_v>3$ or $l_v=2$ and $\alpha_v\not=\beta_v$.

Let $I'=\{v\in S_p\mid l_v=2\}$ and let $I\subset I'$ be a subset.
Then $\pi$ is ordinary with ordinary Hecke
parameters $(\alpha_v,\beta_v)_{v|p}$ if and only if
\begin{equation*}
(\otimes_{v\in S_p} \pi_v)^{K_p(I)}\end{equation*}
contains a vector which is:
\begin{itemize}
\item for each $v\in I^c$, an eigenvector for the normalized $\normU_{\Iw(v),0}$, $\normU_{\Iw(v),1}$,
  and $\normU_{\Iw(v),2}$, with respective eigenvalues $1,\alpha_v$, and
  $\alpha_v\beta_v$, and
\item  for each $v\in I$, an eigenvector for
  $\normU_{\Kli(v),0}$, $\normU_{\Kli(v),1}$, and $\normU_{\Kli(v),2}$ with respective
  eigenvalues $1,\alpha_v+\beta_v$, and
  $\alpha_v\beta_v$.
\end{itemize}
Moreover in this case
\begin{equation*}
\dim e(I)(\otimes_{v\in S_p} \pi_v)^{K_p(I)}=2^{\#(I'-I)}.
\end{equation*}\end{prop}

Note that if~$\pi$ is ordinary with ordinary Hecke
parameters $(\alpha_v,\beta_v)_{v|p}$ but~$v \notin I'$, then the~$U_{\Kli(v),1}$ eigenvalue will~\emph{not} be of the form~$\alpha_v + \beta_v$, but rather,
up to some ordering of~$\alpha_v$ and~$\beta_v$, be of the form~$\alpha_v + p^{l_v-2} \beta_v$.

\begin{proof}
This is simply Proposition \ref{prop: ordclass} applied for each $v\in S_p$.
\end{proof}

\subsubsection{An instance of the local Langlands correspondence}

Given a pair of characters
$\chi_{v,1},\chi_{v,2}:k(v)^\times\to\cO^\times$, which we regard as
characters of~$\cO_{F_v}^\times$ by inflation, we define a character of $\chi_v$ of $T(\cO)$ by
\begin{align*}
\chi_v:T(\cO_{F_v})&\to \cO^\times\\
(a,b,cb^{-1},ca^{-1})&\mapsto\chi_{v,1}(ab^{-1})\chi_{v,2}(abc^{-1}).
\end{align*}

Then if $M$ is an $\cH_1$-module, we write
\begin{equation*}
M^{\chi_v}=\{m\in M\mid tm=\chi_v(t)m\ \forall t\in T(k(v))\}
\end{equation*}
and
\begin{equation*}
M_{\chi_v}=M/\langle tm-\chi_v(t)m\mid t\in T(k(v)),m\in M\rangle.
\end{equation*}

Then we record:
\begin{prop}\label{prop: Ihara local global}If $\pi$ is an irreducible
  smooth $\barE[\GSp_4(F_v)]$-module with the property
  that~$(\pi^{\Iw_1(v)})^{\chi_v}\not=\{0\}$, then, for all $\sigma\in
  W_{F_v}$,
$$
\begin{aligned}
\det(X-\recGTp(\pi)(\sigma))= & \ (X-\chi_{v,1}(\Art^{-1}_{F_v}(\sigma)))(X-\chi_{v,1}(\Art^{-1}_{F_v}(\sigma))^{-1}) \\
 & \ (X-\chi_{v,2}(\Art^{-1}_{F_v}(\sigma)))(X-\chi_{v,2}(\Art^{-1}_{F_v}(\sigma))^{-1}). \end{aligned}$$
If, moreover, the characters $\chi_{v,1},\chi_{v,1}^{-1},\chi_{v,2},\chi_{v,2}^{-1}$ are pairwise distinct, then there is an equality~$\dim_{\overline{E}}(\pi^{\Iw_1(v)})^{\chi_v}=1$.
\end{prop}
\begin{proof}
This is an immediate consequence of Propositions \ref{prop: Jacquet module iso}, \ref{prop: Jacquet module} and \ref{prop: rec of principal series}.
\end{proof}
\subsubsection{The case $q_v\equiv 1\pmod{p}$}\label{subsubsec q
  equals 1}
We suppose from now on for the rest of this section that
$q_v\equiv 1\pmod p$. Recall that we have a homomorphism
$T(F_v)/T(\cO_{F_v})_1\to \cH_1^\times$, and thus an (injective) homomorphism
$\cO[T(F_v)/T(\cO_{F_v})_1]\to \cH_1$; we identify
$\cO[T(F_v)/T(\cO_{F_v})_1]$ with its image in~$\cH_1$. Given elements
$\alphabar_1,\alphabar_2\in\Fbar_p^\times$, we let
$\mathfrak{m}_{\alphabar_1,\alphabar_2}$ denote the 
kernel of the homomorphism
$\cO[T(F_v)/T(\cO_{F_v})_1]\to\Fpbar$ induced by the character 
$T(F_v)/T(\cO_{F_v})_1\to \Fpbar^\times$ sending $T(\cO_{F_v})\mapsto 1$,
$\diag(\varpi_v,\varpi_v,\varpi_v,\varpi)\mapsto 1$,
$\diag(\varpi_v,\varpi_v,1,1)\mapsto \alphabar_1$, and
$\diag(\varpi_v^2,\varpi_v,\varpi_v,1)\mapsto \alphabar_1\alphabar_2$.

\begin{prop}\label{prop: TW local global}Let $\pi$ be an irreducible smooth $\barE[\GSp_4(F_v)]$-module with central character $|\cdot|^2$ and with $(\pi^{\Iw_1(v)})_{\m_{\overline{\alpha}_1,\overline{\alpha}_2}}\not=\{0\}$.  Suppose $\overline{\alpha}_1^{\pm1},\overline{\alpha}_2^{\pm 1}$ are pairwise distinct.  Then
\begin{equation*}
\recGTp(\pi)=\gamma_1\oplus\gamma_2\oplus\varepsilon^{-1}\gamma_2^{-1}\oplus\varepsilon^{-1}\gamma_1^{-1}
\end{equation*}
for characters $\gamma_i$ of $G_{F_v}$ with
$\overline{\gamma}_i=\lambda_{\overline{\alpha}_i}$ \emph{(}the unramified
character taking~$\Frob_v$ to~$\overline{\alpha}_i$\emph{)}, and $T(k(v))$ acts on $(\pi^{\Iw_1(v)})_{\m_{\overline{\alpha}_1,\overline{\alpha}_2}}$ via $(\gamma_i\circ\Art_{F_v})|_{\cO_{F_v}^\times}$.
\end{prop}
\begin{proof}
From Proposition~\ref{prop: Ihara local global}, we know the characteristic polynomial of the corresponding
representation, and thus immediately deduce that the semi-simplification of the Galois
representation has the required form.
It thus suffices to show that, under the hypothesis on~$\overline{\alpha}_i$, that
all Galois representations are semi-simple.  Suppose otherwise. Two tamely ramified characters admit an extension if and only if their ratio is unramified and takes
the value~$q_v$ on Frobenius.  Since~$q_v \equiv 1 \mod p$ and~$\varepsilon$ is trivial modulo~$p$,
this implies that ~$\overline{\alpha}_1^{\pm1},\overline{\alpha}_2^{\pm 1}$ are not distinct,
a contradiction.
\end{proof}

\begin{remark} \label{remark:becausewewantto}Let~$Z$ be the centre of~$\GSp_4$, let $\Delta_v$ be the maximal $p$-power quotient of $T(k(v))/Z(k(v))$, and  let $\Delta_v'=\ker(T(k(v))\to\Delta_v)$. 
If the~$\pi$ of Proposition~\ref{prop: TW local global} additionally satisfies the condition that~$(\pi^{\Iw_1(v)})^{\Delta_v'}_{\m_{\overline{\alpha}_1,\overline{\alpha}_2}}\not=\{0\}$, then we immediately deduce that~$\Delta_v$ also acts on $(\pi^{\Iw_1(v)})^{\Delta_v'}_{\m_{\overline{\alpha}_1,\overline{\alpha}_2}}$ via $(\gamma_i\circ\Art_{F_v})|_{\cO_{F_v}^\times}$.
\end{remark}

We now prove some results about the Iwahori Hecke algebra
(under our running assumption that $q_v\equiv
1\pmod{p}$). We follow~\cite[\S
5]{KT} closely, and our proofs are essentially an immediate adaptation
of their arguments from~$\GL_n$ to~$\GSp_4$.
As recalled above, we have an embedding $\cO[X_*(T)]\into\cH$. This
can be refined to give the Bernstein presentation of~$\cH$ (see e.g.\ \cite[\S 1]{MR2642451}), which
is an algebra isomorphism \[\cH\cong \cO[X_*(T)]\widetilde{\otimes}_{\cO}\cO[\Iw(v)\backslash
  \GSp_4(\cO_{F_v})/\Iw(v)],\]where the twisted tensor product
$\widetilde{\otimes}_{\cO}$ is determined by the following relations, where
$s_\alpha\in W$ is simple, corresponding to the simple root~$\alpha$,
and   $\mu\in X_*(T)$:  \numequation\label{eqn: Bernstein presentation
relation}T_{s_\alpha}\theta_\lambda=\theta_{s_\alpha(\lambda)}T_{s_\alpha}+\bigl(q_v-1\bigr)\frac{\theta_{s(\lambda)}-\theta_\lambda}{1-\theta_{-\alpha^\vee}} .\end{equation} Here we are
writing $\theta_\mu$ for the image in $\cH$ of the group element~$e_\mu$ of
$\cO[X_*(T)]$ corresponding to~$\mu$, and for $w\in W$ we
write~$T_w:=[\Iw(v)\dot{w}\Iw(v)]$ where~$\dot{w}\in  \GSp_4(\cO_{F_v})$ is
any representative for~$w$.

\begin{lem}\label{lem: Iwahori Hecke algebra mod p q equals 1}
  There is a natural isomorphism $\cH\otimes_{\cO}k\cong
  k[X_*(T)\rtimes W]$.
\end{lem}
\begin{proof} We claim that the natural $k$-linear map
  $k[W]\to k[\Iw(v)\backslash \GSp_4(\cO_{F_v})/\Iw(v)]$ sending $w\mapsto
  T_w$ is an algebra isomorphism. Admitting this claim, note that since~$q_v\equiv 1\pmod{p}$, the
  relation~(\ref{eqn: Bernstein presentation relation})
  becomes
  \[T_{s_\alpha}\theta_\lambda=\theta_{s_\alpha(\lambda)}T_{s_{\alpha}} \]in
  $\cH\otimes_{\cO}k$, so that there is an isomorphism
  $k[X_*(T)\rtimes W]\to\cH\otimes_{\cO}k$ sending
  $e_\lambda w\mapsto \theta_\lambda T_w$, as required.
  
It remains to prove the claim. The Weyl group $W$ is generated by
$s_1, s_2$ with $s_1^2=s_2^2=(s_1s_2)^4=1$, so it is enough to show
that $k[\Iw(v)\backslash \GSp_4(\cO_{F_v})/\Iw(v)]$ is generated by the
elements~$T_{s_1}, T_{s_2}$, subject to the same relations. This
follows from the assumption that $q_v\equiv 1\pmod{p}$; indeed, we
have the usual relations $T_{s_i}^2=(q_v-1)T_{s_i}+q_v$ ($i=1,2$), and
$T_{s_1}T_{s_2}T_{s_1}T_{s_2}=T_{s_2}T_{s_1}T_{s_2}T_{s_1}$, which
are easily seen to be equivalent to
$T_{s_1}^2=T_{s_2}^2=(T_{s_1}T_{s_2})^4=1$, as required.
\end{proof}
Recall that by definition an $\cO[\GSp_4(F_v)]$-module~$M$ is \emph{smooth} if
every element of~$M$ is fixed by some open compact subgroup of~$\GSp_4(F_v)$,
and it is \emph{admissible} if it is smooth, and if for each open
compact subgroup $U\subset \GSp_4(F_v)$, $M^U$ is a finite $\cO$-module.

\begin{lem}\label{lem: splitting off invariants inside Iwahori
    invariants}If~$M$ is a smooth $\cO[\GSp_4(F_v)]$-module, then the
  natural inclusion $M^{\GSp_4(\cO_{F_v})}\subset M^{\Iw(v)}$ is
  canonically split by the Hecke operator 
  $$\frac{1}{[\GSp_4(\cO_{F_v}):\Iw(v)]}e_{\Sph(v)}.$$
\end{lem}
\begin{proof}
  The Hecke operator $e_{\Sph(v)}\in\cH$ induces the natural
  trace map $ M^{\Iw(v)}\to M^{\GSp_4(\cO_{F_v})}$, so that the composite
  map $M^{\GSp_4(\cO_{F_v})}\to M^{\Iw(v)} \to M^{\GSp_4(\cO_{F_v})}$ is given
  by multiplication by~$[\GSp_4(\cO_{F_v}):\Iw(v)]$. Since
  $[\GSp_4(\cO_{F_v}):\Iw(v)]\equiv |W|=8\pmod{p}$ is a unit in~$\cO$, we
  are done.
\end{proof}

\begin{cor}
  \label{cor: mod p max compact invariants are Weyl group
    invariants}If $M$ is a smooth $k[G]$-module, then $M^{\Iw(v)}$ is
  naturally a $k[W]$-module, and $M^{\GSp_4(\cO_{F_v})}=(M^{\Iw(v)})^W$.
\end{cor}
\begin{proof}
  This is immediate from Lemmas~\ref{lem: Iwahori Hecke algebra mod p
    q equals 1} and~\ref{lem: splitting off invariants inside Iwahori
    invariants}.
\end{proof}

 The centre of~$\cH$ is $\cO[X_*(T)]^W$, and there is an
isomorphism \[\cO[X_*(T)]^W\cong \cO[\GSp_4(\cO_{F_v})\backslash
\GSp_4(F_v)/\GSp_4(\cO_{F_v})]\] given by $x\mapsto e_{\Sph(v)}x$ (where we
are regarding~$x$ as an element of~$\cH$); this
isomorphism agrees with the isomorphism given by the usual Satake isomorphism (see~\cite[\S
4.6]{MR2642451}). The classical description of~$\cO[X_*(T)]$ is as
follows. 
Let~$x_0$, $x_1$, and~$x_2$ denote the following three cocharacters:
$$x_0: t \rightarrow \diag(t,t,1,1),$$
$$x_1: t \rightarrow \diag(1/t,1,1,t),$$
$$x_2: t \rightarrow \diag(1,1/t,t,1).$$
Then~$x^2_0 x_1 x_2$ is the cocharacter~$t \mapsto
\diag(t,t,t,t)$ and
$$\cO[X_*(T)]=\cO[x_0,x_1,x_2,(x^2_0 x_1 x_2)^{-1}] = \cO[x_0,x_1,x_2,(x_0 x_1 x_2)^{-1}].$$
The effect of the involutions ~$s_1$, $s_2$, and~$s_1 s_2 s_1 \in W$ on these
cocharacters is to send~$(x_0,x_1,x_2)$ to
$$(x_0,x_2, x_1), (x_0 x_2,x_1, x^{-1}_2), (x_0 x_1,x^{-1}_1, x_2)$$
respectively. 
All of these involutions preserve~$(x_0,x_0
x_1, x_0 x_2, x_0 x_1 x_2)$ considered as an  unordered quadruple.
Define elements~$e_i(x_0,x_1,x_2) \in \cO[X_*(T)]^{W}$, $0\le i \le 4$, by the following formulae:
$$(X - x_0)(X - x_0 x_1)(X - x_0 x_2)(X - x_0 x_1 x_2) = \sum  e_i(x_0,x_1,x_2)X^i.$$
The relation between the~$e_i$ and the  Hecke
operators $T_{v,i}$ is given by \[\sum  e_i(x_0,x_1,x_2)X^i =X^4 - q_v^{3/2}T_{v,1} X^3 + (q_v^2 T_{v,2} + (1+ q^2_v) T_{v,0}) X^2
- q_v^{3/2}T_{v,0} T_{v,1}X  + T^2_{v,0}.\] Since we are assuming
that~$q_v\equiv 1\pmod{p}$, 
and in our applications of these results in the global setting there
is a twist which makes all of the powers of~$q_v$ integral (as
in~\eqref{eqn: char poly for unramified Hecke}), we will ignore all powers of~$q_v^{1/2}$
from now on.

Given any triple~$\gamma:=(\gamma_0,\gamma_1,\gamma_2)$ and~$w \in W$, 
let~$((w \gamma)_0,(w \gamma)_1,(w \gamma)_2)$ denote the triple obtained by substituting in~$\gamma_i$
for~$x_i$ in the action of~$W$ on~$\cO[X_*(T)]$ described above.

\begin{lem}\label{lem: eigenvalues at level q are in the
    support}Let~$M$ be an $\cH\otimes_{\cO}k$-module which is
  finite-dimensional over~$k$. Suppose that $e_{\Sph(v)}M\ne 0$,
  and that there is a triple~$\gamma_0$, $\gamma_1$, $\gamma_2$ with
~$\gamma^2_0 \gamma_1 \gamma_2=1$
  such that
  $(\gamma_1 - 1)(\gamma_2 - 1)(\gamma_1 - \gamma_2)(\gamma_1 \gamma_2 - 1) \ne 0$;
  equivalently, writing   $\alpha_1 = \gamma_0$, $\alpha_2 = \gamma_0 \gamma_1$,
  suppose that \[\alpha_1,\alpha_2,1/\alpha_2,1/\alpha_1\]are pairwise distinct.
Suppose also 
 that the following operators act by zero on the module~$e_{\Sph(v)}M$:
  $$T_{v,0} - 1 \ T_{v,1} - e_1(\gamma_0,\gamma_1,\gamma_2),
T_{v,2}  + 2T_{v,0} - e_2(\gamma_0,\gamma_1,\gamma_2).$$
   Then, for
  each~$w\in W$, the maximal ideal 
  $$\m_w=(x_0 - (w \gamma)_0,x_1 - (w \gamma)_1,x_2 - (w \gamma)_2) \subset k[X_*(T)]$$
   is in
  the support of~$M$.  
\end{lem}
\begin{proof}Let~$\mathfrak{n} \subset k[X_*(T)]^W$ be the ideal
$$\mathfrak{n}=(e_1(x_0,x_1,x_2) - e_1(\gamma_0,\gamma_1,\gamma_2),
\ldots, e_4(x_0,x_1,x_2) - e_4(\gamma_0,\gamma_1,\gamma_2), x^2_0 x_1 x_2 - \gamma^2_0 \gamma_1 \gamma_2).$$
  Then, by assumption, we have $e_{\Sph(v)}M\subset M[\mathfrak{n}]$,
  so that in
  particular $M_\mathfrak{n}\ne 0$. 
      The assumptions on~$\gamma_i$ imply that 
  all the ideals~$\m_w$ are distinct.
  We may  view~$\mathfrak{n}$ as an ideal in~$k[X_*(T)]$. 
The support of~$\mathfrak{n}$ in~$k[X_*(T)]$ corresponds to 
 triples~$(\gamma_0,\gamma_1,\gamma_2)$ (or equivalently,
pairs~$(\alpha_1,\alpha_2)$) such that~$\alpha_1$, $\alpha_2$, $\alpha^{-1}_2$, and~$\alpha^{-1}_1$
are roots of the polynomial~$\sum e_i(\gamma_0,x_1,x_2) X^i$. Hence the support
of~$\mathfrak{n} \subset k[X_*(T)]$ consists exactly of
the maximal ideals~$\m_w$, and   the product of the~$\m_w$ is precisely
the radical of~$\mathfrak{n}$. 
The ring~$k[X_*(T)]_{\mathfrak{n}}$ is thus a semi-local ring which is isomorphic
  to~$\oplus_{w\in W}k[X_*(T)]_{\m_w}$, and correspondingly we may write
  $M_\mathfrak{n}=\oplus_{w\in W}M_{\m_w}$. It follows that~$M_{\m_w}\ne 0$
  for at least one~$w\in W$. Considering the action of~$W$ on the set
  of maximal ideals of~$k[X_*(T)]$ in the support of~$M$, we see that
  in fact $M_{\m_w}\ne 0$ for all ~$w\in W$, as required.  
\end{proof}

\begin{lem}\label{lem: bijection between localizations at level q and
    level one}Let~$M$ be an $\cH\otimes_{\cO}k$-module which is
  finite-dimensional over~$k$. Suppose that for each maximal ideal
  $\mathfrak{n}\subset k[X_*(T)]^W$ in the support of~$M$,  the degree four polynomial
 $$ \sum  e_i(x_0,x_1,x_2)X^i \in k[X_*(T)]^{W}[X]$$
 has roots~$(\gamma_0,\gamma_0\gamma_1,\gamma_0\gamma_2,\gamma_0\gamma_1\gamma_2)$ modulo~$\mathfrak{n}$
 satisfying~$(\gamma_1 - 1)(\gamma_2 - 1)(\gamma_1 -
 \gamma_2)(\gamma_1 \gamma_2 - 1) \ne 0$ and~$\gamma_0^2\gamma_1\gamma_2=1$.
 Equivalently, writing $\gamma_0=\alpha_1$,
 $\gamma_0\gamma_1=\alpha_2$, assume
 that~$\gamma_0^2\gamma_1\gamma_2=1$ and
 that \[\alpha_1,\alpha_2,1/\alpha_2,1/\alpha_1\]are pairwise distinct. Then~$e_{\Sph(v)}M\ne 0$. If, furthermore, there is a unique maximal
  ideal  $\mathfrak{n}\subset k[X_*(T)]^W$ in the support of~$M$, then
  for each maximal ideal $\m\subset k[X_*(T)]$ in the support of~$M$,
  the maps \[k[W]\otimes_k M_\m\to M,\] \[ w\otimes x\mapsto w\cdot
    x, \]and\[M_\m\to e_{\Sph(v)}M,\]  \[x\mapsto e_{\Sph(v)}\cdot
    x  \]are both isomorphisms.  
\end{lem}
\begin{proof}After possibly enlarging~$k$, we can and do assume that
  the~$\gamma_i$ arising from the roots of the degree four polynomial above lie in~$k$. 
    As in the proof of Lemma~\ref{lem: eigenvalues at level q are in the
    support}, there exist~$|W| = 8$ distinct ideals~$\m_w$ such that
    $M_{\mathfrak{n}} \simeq  \oplus_{w\in W} M_{\m_w}$,
    where
     $\m=(t_0 - \gamma_0,t_1 -  \gamma_1,t_2 - \gamma_2)  \subset
  k[X_*(T)]$.
    Since~$M_{\mathfrak{n}} \ne 0$, we may assume that~$M_{\m_w} \ne 0$ for
    some and hence all~$\m_w$.
    The operator~$e_{\Sph(v)}$ acts by averaging over the action  of the Weyl
    group. It follows (because the~$\m_w$ are distinct)
    that the map~$e_{\Sph(v)}: M_{\m} \rightarrow  \oplus_{w\in W} M_{\m_w}
    = M_{\mathfrak{n}}$ is an injection, and thus~$e_{\Sph(v)} M \ne 0$.

Suppose that~$\mathfrak{n}$ is the only maximal ideal of~$k[X_*(T)]$
in the support of~$M$. Then the maximal ideals of~$k[X_*(T)]$ in the
support of~$M$ are necessarily of the form~$\m_w$, and we have
$M=\oplus_{w\in W}M_{\m_w}=\oplus_{w\in W}w\cdot M_\m$, and the rest
of the lemma follows immediately.
\end{proof}
\begin{rem}
  \label{rem: connection between the x and the U}Note that (using as
  usual that ~$q_v=1$ in~$k$) we have that~$\Unaive_{v,0}=x_0^2x_1x_2$, and
  if this equals~$1$, then $\Unaive_{v,1}=x_0$
  and~$\Unaive_{v,2}=(s_1s_2s_1)x_1$. Consequently we see for example
  that if the hypotheses of Lemma~\ref{lem: eigenvalues at level q are in the
    support} hold then
  $(U_{v,0}-1,U_{v,1}-\alpha_1,U_{v,2}-\alpha_1\alpha_2)$ is in the
  support of~$M$.
\end{rem}

\subsection{Purity}Let~$K$ be a finite extension of~$\Qp$ for some~$p$, with residue
field of order~$q$. Following~\cite[\S1]{ty}, we say that a
Weil--Deligne representation $(W,r,N)$ of~$W_K$ on a  vector space $W$ over an algebraically closed
field~$\Omega$ which is of characteristic~$0$ and of the same cardinality as~$\C$ is  {\em pure of weight $w$} if
there is an exhaustive and separated ascending filtration $\Fil_i$ of $W$
such that
\begin{itemize}
\item each $\Fil_iW$ is invariant under $r$;
\item if $\sigma \in W_K$ maps to $\Frob_K^{v(\sigma)}$, then all
  eigenvalues of $r(\sigma)$ on $\gr_iW$ are Weil
  $q^{iv(\sigma)}$-numbers;\item and for all $j$ we have $N^j:\gr_{w+j}W \iso \gr_{w-j}W$. (Note that necessarily we have $N \Fil_i W \subset \Fil_{i-2}W$.)
\end{itemize}

Recall that for a Weil--Deligne representation~$(r,N)$, we defined in
Section~\ref{subsec: nonarch LL} $n(r,N)$ to be the rank of~$N$.
\begin{lem}
  \label{lem: pure implies N is maximal}If~$(V,r)$ is a semisimple
  representation of~$W_K$, then there is at most one choice of~$N$ for
  which $(V,r,N)$ is a pure Weil--Deligne representation. If such
  an~$N$ exists, then the corresponding Weil--Deligne representation
  is the unique choice which maximizes~$n(r,N)$.  \end{lem}
\begin{proof}
  The uniqueness of~$N$ is~\cite[Lem.\ 1.4(4)]{ty}. The maximality
  follows easily, using that by definition all of the induced maps
  $N^j:\gr_{w+j}W \to \gr_{w-j}W$ are isomorphisms if and only
  if~$(V,r,N)$ is pure. \end{proof}

\subsection{Archimedean
  \texorpdfstring{$L$}{L}-parameters}\label{subsec: archimedean L
packets}We now recall some notation
for archimedean $L$-parameters following~\cite[\S 3.1]{MR3200667}
(although our~$w$ has the opposite sign to this reference).
Recall that $W_\R=\C^\times\cup\C^\times j$, where $jzj^{-1}=\overline{z}$
and $j^2=-1$. Let~$w \in \R$. For an integer~$n \ge 0$,
let~$\phi_{w,n}:W_\R\to \GL_2(\C)$ be the~$L$-parameter given by
\[z\mapsto |z|^{w}\cdot
  \begin{pmatrix}
    (z/\zbar)^{n/2} &\\
& (z/\zbar)^{-n/2}
   \end{pmatrix} = 
    |z|^{w}  \begin{pmatrix}
    z^n |z|^{-n} & \\
    & z^{-n} |z|^{n}
  \end{pmatrix} 
  \]and \[j\mapsto
  \begin{pmatrix}
& 1\\
(-1)^{n}
\end{pmatrix}.\]

 The determinant of~$\phi_{w,n}$ is equal to~$|z|^{2w}$ if~$n$ is odd
 and~$\sgn \cdot |z|^{2w}$ if~$n$ is even, where~$\sgn: W_{\R} \rightarrow \C^{\times}$
 is the degree two character which is~$-1$ on~$j$ (and trivial
 on~$\C^\times$). We also write~$\phi_{w,n}$ for the restriction
 of~$\phi_{w,n}$ to~$W_\C$.
 The~$\GL_2(\R)$ and~$\GL_2(\C)$ representations
 corresponding to the $L$-parameter~$\phi_{w,n}$  
are cohomological
if and only if~$n>0$ and~$w \in \Z$ satisfies~$w + n \equiv 1 \mod
2$.

Let~$m_1 > m_2\ge 0$ be integers, and let~$w \in \R$.
Then we write $\phi_{(w;m_1,m_2)}:W_\R\to\GSp_4(\C)$ for the
$L$-parameter sending 
\[z\mapsto |z|^{w}\cdot
  \begin{pmatrix}
    (z/\zbar)^{(m_1+m_2)/2} & & &\\
& (z/\zbar)^{(m_1-m_2)/2} && \\
&& (z/\zbar)^{-(m_1-m_2)/2}&\\
&&& (z/\zbar)^{-(m_1+m_2)/2}
  \end{pmatrix}
\]
and
 \[j\mapsto
  \begin{pmatrix}
    &&& 1\\
&&1&\\
& (-1)^{m_1 + m_2}&&\\
(-1)^{m_1 + m_2}&&&
\end{pmatrix}.\]
Note that~$\phi_{(w;m_1,m_2)}$ is viewed as having image in~$\GSp_4(\C)$ with respect
to our particular choice of model for~$\GSp_4(\C)$ where~$J$ is anti-diagonal. In particular,
 the image of~$j$ under the composite of~$\phi_{(w;m_1,m_2)}$ with
the similitude character is~$(-1)^{m_1 + m_2}$. 
With respect to the explicit inclusion of
$$\{(A,B) \subset \GL_2(\C) \times \GL_2(\C) \ | \ \det(A) = \det(B)\}
\subset \GSp_4(\C)$$
 given in~\S\ref{section:inductions}, we immediately observe that 
the composite of~$\phi_{(w;m_1,m_2)}$ with the inclusion~$\GSp_4(\C) \rightarrow \GL_4(\C)$
 identifies~$\phi_{(w;m_1,m_2)}$ with~$\phi_{w,m_1 + m_2} \oplus \phi_{w,m_1 - m_2}$ (note that~$(-1)^{m_1 + m_2} = (-1)^{m_1 - m_2}).$
 The $L$-packet of~$\GSp_4(\R)$ corresponding to~$\phi_{(w;m_1,m_2)}$ consists
 of two elements $\pi_{(w;m_1,m_2)}^{H}$ and $\pi_{(w;m_1,m_2)}^{W}$.
When~$m_2=0$, they are (up to twist) non-degenerate limits of discrete
series, and when~$m_2>0$, they are (up to twist) discrete series. The
representations $\pi_{(w;m_1,m_2)}^{H}$ and $\pi_{(w;m_1,m_2)}^{W}$
are respectively holomorphic and generic. Their central character is
given by~$a\mapsto a^{w}$, and they are tempered when~$w=0$.
The minimal $K$-type of
$\pi^H_{(w;m_1,m_2)}$ is the representation
$\det^{m_2+2}\otimes\Sym^{m_1-m_2-1}\C^2$ of $U(2)$. (See for
example~\cite{SchmidtHodge} for these facts and their proofs.)

\begin{lem}[Inductions of real archimedean parameters to~$\GL_4(\C)$] \leavevmode
  \label{lem: induction of real archimedean L parameters from GL2 to GSp4} \label{lemma:inductions}
 \begin{enumerate}
 \item The induction~$\Ind_{W_\C}^{W_\R}\phi_{w,n}: W_{\R} \rightarrow \GL_4(\C)$ is conjugate
 to~$\phi_{w,n} \oplus \phi_{w,n}$.
  \item The composite map
  $$\phi_{(w;n,0)}:W_\R\to\GSp_4(\C) \rightarrow \GL_4(\C)$$
  is conjugate to~$\phi_{w,n}\oplus\phi_{w,n}$.
  \item If~$\varphi:W_\C\to\GL_2(\C)$ is such that $\Ind_{W_\C}^{W_\R}\varphi$ is conjugate
 to~$\phi_{w,n} \oplus \phi_{w,n}$, and~$n \ne 0$, then either $\varphi\cong\phi_{w,n}$, or~$\varphi$
 is one of the scalar $L$-parameters sending $z$ to one of
 $$ |z|^{w} \cdot  \left( \begin{matrix} z^n |z|^{-n} & 0 \\ 0 & z^n |z|^{-n} \end{matrix} \right) \ \text{or} \ 
 |z|^{w} \cdot    \left( \begin{matrix} z^{-n} |z|^{n} & 0 \\ 0 & z^{-n} |z|^{n} \end{matrix} \right).$$
 \item If~$\varphi, \varphi':W_\R\to\GL_2(\C)$ are such that $\varphi\oplus\varphi'$ is conjugate 
to~$\phi_{w,n} \oplus \phi_{w,n}$, and~$n \ne 0$,
 then $\varphi \cong \varphi' \cong\phi_{w,n}$.
  \end{enumerate}
\end{lem}
\begin{proof}
Since~$\phi_{w,n}$ is already a representation
of~$W_{\R}$, the first induction is isomorphic to~$\phi_{w,n} \oplus \phi_{w,n} \otimes \sgn$.
Yet~$\phi_{w,n}$ is itself induced from~$\C^{\times}$,
and so~$\phi_{w,n} \otimes \sgn \simeq \phi_{w,n}$. 
The second claim was already noted above.
Now suppose that~$\varphi: W_{\C} \rightarrow \GL_2(\C)$ is a complex~$L$-parameter.
All such parameters are of the form
$$z^{a_1} |z|^{-a_1} |z|^{w_1} \oplus  z^{a_2} |z|^{-a_2} |z|^{ w_2}$$
for integers~$a_1$ and~$a_2$. The induction of this representation to~$W_{\R}$ is~$\phi_{w_1,a_1} \oplus \phi_{w_2,a_2}$.
Now consider the equality of~$\GL_4(\C)$-representations
$$\phi_{w_1,a_1} \oplus \phi_{w_2,a_2} = \phi_{(w;n,0)}
 = \phi_{w,n} \oplus \phi_{w,n}.$$
 Restricting to~$S^1 \subset \C^{\times} \subset W_{\R}$,  we deduce that~$|a_1| = |a_2| = n$,
 and then restricting to the  action of~$\C^{\times}$ on the eigenspace where~$S^1 \subset\C^\times $ acts by~$z^n$
 (which is distinct from~$z^{-n}$),
 we deduce that~$w_1=w_2=w$, and thus~$\phi_{w_1,a_1} = \phi_{w_2,a_2} = \phi_{w,n}$. If~$a_1$ and~$a_2$ have opposite signs, then~$\varphi= \phi_{w,n}$;
 otherwise we get the possibilities outlined in the statement of the lemma.
 Finally, (4) is immediate from the irreducibility of~$\phi_{w,n}$.
\end{proof}

We note in passing that the~$\GSp_4(\C)$-parameter cannot be
recovered, in general,  from the~$\GL_4(\C)$-parameter.
This is true in particular for~$\phi_{(0;1,0)}$, since one may compute that the~$\GL_4(\C)$
representation preserves two symplectic forms whose similitude characters
differ by~$\sgn$.

If~$K$ is a number field and~$\boldsymbol{\pi}$ is an automorphic representation
of~$\GL_2(\A_K)$, we say that~$\boldsymbol{\pi}$ has weight~$0$ if for each
place~$v|\infty$ of~$K$, $\boldsymbol{\pi}_v$ corresponds to~$\phi_{0,1}$. If~$F$
is a totally real field and~$\pi$ is an automorphic representation
of~$\GSp_4(\A_F)$, then we say that~$\pi$ has
weight~$(k_v,l_v)_{v|\infty}$ if for each place~$v|\infty$ of~$F$, we
have $k_v\ge l_v\ge 2$ and
$k_v\equiv l_v\pmod{2}$, and $\pi_v$ is in the $L$-packet
corresponding to~$\phi_{(2;k_v-1, l_v-2)}$. We say that~$\pi$ has
parallel weight~$2$ if it has
weight~$(2,2)_{v|\infty}$ (we note that the congruence $k_v\equiv l_v\pmod{2}$ is imposed in order to ensure that $\pi$ is algebraic.)

\subsection{Galois representations associated to automorphic
  representations}\label{subsec: Galois associated to automorphic}
We now recall some results from~\cite{MR3200667} on the existence of
Galois representations (adapted to the particular setting of interest
for us), beginning with the
existence of Galois representations for certain cuspidal automorphic
  representation of~$\GSp_4(\A_F)$. The following theorem is
  essentially due to Sorensen~\cite{sorensen}, although at the time
  that~\cite{sorensen} was written, some additional assumptions needed
  to be made, due to the lack of unconditional results on the transfer
  of automorphic representations between~$\GSp_4$ and~$\GL_4$.
\begin{thm}
  \label{thm: existence and properties of Galois representations for
    automorphic representations general Mok version}Suppose that~$F$ is a
  totally real field, and that~$\pi$ is a cuspidal automorphic
  representation of~$\GSp_4(\A_F)$ of weight~$(k_v,l_v)_{v|\infty}$,
  where $k_v\ge l_v>2$ and $k_v\equiv l_v\pmod{2}$ for all $v|\infty$.
  Suppose also that~$\pi$ has central character~$|\cdot|^2$.

Fix a prime~$p$. Then there
is a continuous semisimple representation
$\rho_{\pi,p}:G_F\to\GSp_4(\Qpbar)$ satisfying the following properties.
\begin{enumerate}
\item 
  $\nu\circ\rho_{\pi,p}=\varepsilon^{-1}$. \item For each finite place~$v$, we have \[\WD(\rho_{\pi,p}|_{G_{F_v}})^{\semis}\cong\recGTp(\pi_v\otimes|\nu|^{-3/2})^{\semis}.\]If furthermore~$\rho_{\pi,p}$ is irreducible, then
\[\WD(\rho_{\pi,p}|_{G_{F_v}})^{F-\semis}\cong\recGTp(\pi_v\otimes|\nu|^{-3/2}).\]
\item If $v|p$, then $\rho_{\pi,p}|_{G_{F_v}}$ is de Rham with Hodge--Tate
weights $((k_v+l_v)/2-1,(k_v-l_v)/2+1,-(k_v-l_v)/2,2-(k_v+l_v)/2)$.

\item If~$\rho_{\pi,p}$ is irreducible, then for each finite place~$v$ of~$F$, $\rho_{\pi,p}|_{G_{F_v}}$
  is pure.
\end{enumerate}

\end{thm}
\begin{proof}The existence of a representation~$\rho_{\pi,p}$
  valued in~$\GL_4(\Qpbar)$ and satisfying~(2) and~(3) is part
  of~\cite[Thm.\ 3.5]{MR3200667} (note that the results of~\cite{MR2058604} cited
  in~\cite{MR3200667} hold unconditionally by~\cite{GeeTaibi}). That the representation
  actually takes values in~$\GSp_4(\Qpbar)$ with the claimed
  multiplier follows from~\cite[Cor.\ 1.3]{belchen} (cf.\ \cite[Rem.\
  3.3(3)]{MR3200667}). Finally, for
 part~(4), note that if~$\rho_{\pi,p}$ is irreducible, then~$\pi$
 is of general type in the sense of~\cite{MR2058604} (see Section~\ref{subsec: Arthur classification}), and thus
 corresponds to an essentially self-dual algebraic automorphic
 representation~$\Pi$ of~$\GL_4$. Purity then follows from the main
 results of~\cite{MR2972460,MR3272276}.
\end{proof}
For representations which are ordinary in the sense of
Section~\ref{subsubsec: ordinary Hecke operators}, we have the
following variant on Theorem~\ref{thm: existence and properties of Galois representations for
    automorphic representations general Mok version}.
\begin{thm}\label{thm: existence and properties of Galois representations for
    automorphic representations}Suppose that~$F$ is a
  totally real field, and that~$\pi$ is a cuspidal automorphic
  representation of~$\GSp_4(\A_F)$ of weight~$(k_v,l_v)_{v|\infty}$,
  where $k_v\ge l_v>2$ and $k_v\equiv l_v\pmod{2}$ for all $v|\infty$.
  Suppose also that~$\pi$ has central character~$|\cdot|^2$.

Fix a prime~$p$. Assume
that~$\pi_v$ is unramified at all places $v|p$, and that~$\pi$ is
ordinary, with ordinary Hecke
parameters~$(\alpha_v,\beta_v)_{v|p}$. Then there
is a continuous semisimple representation
$\rho_{\pi,p}:G_F\to\GSp_4(\Qpbar)$ satisfying the following properties.
\begin{enumerate}
\item 
  $\nu\circ\rho_{\pi,p}=\varepsilon^{-1}$. \item For each finite place~$v\nmid p$, we have \[\WD(\rho_{\pi,p}|_{G_{F_v}})^{\semis}\cong\recGTp(\pi_v\otimes|\nu|^{-3/2})^{\semis}.\]If furthermore~$\rho_{\pi,p}$ is irreducible, then
\[\WD(\rho_{\pi,p}|_{G_{F_v}})^{F-\semis}\cong\recGTp(\pi_v\otimes|\nu|^{-3/2}).\]
\item If $v|p$, then
 \[\rho_{\pi,p}|_{G_{F_v}}\cong  \begin{pmatrix}
    \lambda_{\alpha_v}\varepsilon^{(k_{v}+l_{v})/2-2}&*&*&*\\
0 &\lambda_{\beta_v}\varepsilon^{(k_{v}-l_{v})/2}&*&*\\
0&0&\lambda_{\beta_v}^{-1}\varepsilon^{-1-(k_{v}-l_{v})/2}&*\\
0&0&0&\lambda_{\alpha_v}^{-1}\varepsilon^{1-(k_{v}+l_{v})/2}
  \end{pmatrix}.\] 
\item If~$\rho_{\pi,p}$ is irreducible, then for each finite place~$v$ of~$F$, $\rho_{\pi,p}|_{G_{F_v}}$
  is pure.
\end{enumerate}

\end{thm}
\begin{proof}

This follows from Theorem~\ref{thm: existence and properties of Galois representations for
    automorphic representations general Mok version}; part~(3) is a
 standard consequence of $p$-adic Hodge theory, and is in particular immediate from~\cite[Lem. 2.32]{ger} 
  (and
 Proposition~\ref{prop: rec of principal series}). \end{proof}

The following theorem is a variant of the main result
of~\cite{MR3200667}, which proves the existence of Galois
representations associated to certain automorphic representations
of~$\GL_2(K)$, $K$ a CM field. 
\begin{thm}
  \label{thm: Galois representations for GL2 over quadratic over
    totally real}Let $F$ be a totally real field, and let $K/F$ be a
  quadratic extension. Write~$\Gal(K/F)=\{1,\tau\}$. Suppose
  that~$\boldsymbol{\pi}$ is a cuspidal automorphic representation of~$\GL_2(K)$ of weight~$0$
  with trivial central character.

Then there
is a continuous irreducible representation
$\rho_{\boldsymbol{\pi},p}:G_K\to\GL_2(\Qpbar)$ such that for each finite place~$w\nmid p$ of~$K$, we have \[\WD(\rho_{\boldsymbol{\pi},p}|_{G_{K_w}})^{\semis}\cong\rec_p(\boldsymbol{\pi}_w\otimes|\cdot|^{-1/2})^{\semis}.\]
  If~$\boldsymbol{\pi}_w$ is not a twist of a Steinberg representation, then in fact \[\WD(\rho_{\boldsymbol{\pi},p}|_{G_{K_w}})^{F-\semis}\cong\rec_p(\boldsymbol{\pi}_w\otimes|\cdot|^{-1/2}).\]
\item For each place $w|p$ of~$K$, the representation
  $\rho_{\boldsymbol{\pi},p}|_{G_{K_v}}$ is Hodge--Tate, and for
  each~$\tau:K\into\Qpbar$, the $\tau$-labelled Hodge--Tate weights
  of~$\rho_{\boldsymbol{\pi},p}$ are $(0,1)$.

\end{thm}
\begin{proof}[Proof of Theorem~\ref{thm: Galois representations for GL2 over quadratic over
    totally real}]In the case that~$K$ is CM this is a special case of the main theorem
  of~\cite{MR3200667}, and essentially the same proof works in the
  general case. The argument of~\cite[\S 5.1]{MR3200667} goes over
  unchanged to produce a cuspidal automorphic representation~$\pi$
  of~$\GSp_4(\A_F)$ (see Theorem~\ref{thm: arthur classification
    results} below); to see that $\pi_v$ is in the $L$-packet
  corresponding to~$\phi_{(2;1,0)}$ at each place~$v|\infty$
  of~$F$, one uses Lemma~\ref{lem: induction of real archimedean L parameters from GL2 to
    GSp4} at the places which split in~$K$, and~\cite[Prop.\
  5.2]{MR3200667} at the places for which~$K_v$ is complex. One then easily checks that the arguments
  of~\cite[\S 5.2-5.3]{MR3200667} go over without any changes to the case of
  general~$K$, as required.  
\end{proof}

\subsection{Compatible systems of Galois representations,
  \texorpdfstring{$L$}{L}-functions, and Hasse--Weil zeta
  functions}\label{subsec: compatible systems}We now recall
some definitions concerning compatible systems from~\cite[\S
5]{BLGGT} and~\cite[\S 1]{MR3314824}; in fact, our definition of a
``strictly compatible system'' differs
slightly from the definitions in those papers, because we find it
convenient to include local-global compatibility at places dividing~$p$. 
Let $F$ denote a number field. 
By a {\em rank $n$ weakly compatible system of $l$-adic representations} $\CR$
{\em of} $G_F$ {\em defined over} $M$ we mean a $5$-tuple 
\[ (M,S,\{ Q_v(X) \}, \{r_\lambda \}, \{H_\tau\} ) \]
where
\begin{enumerate}
\item $M$ is a number field considered as a subfield of~$\C$;
\item $S$ is a finite set of primes of $F$;
\item for each  prime $v\not\in S$ of $F$, $Q_v(X)$ is a monic degree $n$
polynomial in $M[X]$;
\item for each prime $\lambda$ of $M$ (with residue characteristic $l$, say) 
\[r_\lambda:G_F \lra \GL_n(\barM_\lambda) \]
is a continuous, semi-simple, representation such that 
\begin{itemize}
\item if $v \not\in S$ and $v \ndiv l$ is a prime of $F$, then $r_\lambda$
is unramified at $v$ and $r_\lambda(\Frob_v)$ has characteristic
polynomial $Q_v(X)$,
\item while if $v|l$, then $r_\lambda|_{G_{F_v}}$ is de Rham  and in the case $v \not\in S$ crystalline;
\end{itemize}
\item for $\tau:F \into \barM$, $H_\tau$ is a multiset of $n$ integers such that 
for any $\barM \into \barM_\lambda$ over $M$ we have 
$\HT_\tau(r_\lambda)=H_\tau$.
\end{enumerate}
If $\cR=(M,S,\{ Q_v(X) \}, \{r_\lambda \}, \{H_\tau\} )$ and
$\cR'=(M',S',\{ Q'_v(X) \}, \{r'_\lambda \}, \{H'_\tau\} )$ are two
compatible systems, then we write $\cR\cong\cR'$ if $Q_v(X)=Q'_v(X)$
for a set of places~$v$ of Dirichlet density one. This implies that
$Q_v(X)=Q'_v(X)$ for all~$v\notin S\cup S'$, and that $r_\lambda\cong
r'_\lambda$ for all~$\lambda$, and $H'_\tau=H_\tau$ for all~$\tau$.

We say that~$\CR$ is {\em regular} if for each $\tau:F\into\barM$, the
elements of~$H_\tau$ are pairwise distinct. We will call $\CR$ {\em strictly compatible} if for each finite place
$v$ of $F$ there is a Weil--Deligne representation $\WD_v(\CR)$ of
$W_{F_v}$ over $\barM$ such that for each place $\lambda$ of $M$ and every $M$-linear
embedding $\varsigma:\barM \into \barM_\lambda$ we have
$\varsigma\WD_v(\CR)\cong \WD(r_\lambda|_{G_{F_v}})^\Fsemis$.

We will call a strictly compatible system $\CR$ \emph{pure} of weight $w$ if
for each finite place $v$ of $F$ the
  Weil--Deligne representation $\WD_v(\CR)$ is pure of weight
  $w$.

The following result is well-known (see for example~\cite[Rem.\
2.4.6]{MR1293977}), but as we do not know
of a convenient reference for a proof, we briefly explain how it
follows from results in the literature.
\begin{prop}\label{prop: pure abelian variety compatible system}If~$A$
  is an abelian variety  over a number
  field~$F$, then, for each~$0\le i\le 2\dim X$, the $l$-adic cohomology
  groups~$H^i(A_{\overline{F}},\Qlbar)$ form a strictly  compatible
  system which is pure of
  weight~$i$ and which is defined over~$\Q$.
\end{prop}

\begin{proof} Since $H^i(A,\Qlbar)=\wedge^i H^1(A,\Qlbar)$, it is
  enough to check the case~$i=1$. The compatible system satisfies strict compatibility at the places
   not dividing~$l$ by~\cite[Cor.\
  2.7]{MR3123639}. In the case that~$A$ has semistable reduction, it
  is furthermore strictly compatible by~\cite[Cor.\ 2.2]{MR3573156}. One can
  deduce the general case from this by a base change trick due to
  Saito~\cite{MR1465337}, which was exploited
  in~\cite{kisindefrings,MR2538615,BLGGT11}. Indeed, as in the proof
  of~\cite[Thm.\ 2.1]{BLGGT11}, it suffices for each finite place~$v$ of~$F$ to check that
  whenever~$g\in W_{F_v}$ maps to a positive power of Frobenius in the
  absolute Galois group of the residue field, then the trace of~$g$
  on~$\WD(H^1(A,\Qlbar))$ is independent of~$l$. One can choose an
  extension~$E/F$ (for example, the fixed field of the subgroup
  of~$W_F$ generated by~$g$ and the kernel of the restriction to~$I_F$
  of~$\WD(H^1(A,\Qlbar))$ for some~$l$) and a place~$v|w$ of~$E$ such that~$A_E$ is semistable
  and $g\in W_{E_w}$, and the claim then follows from the independence
  of~$l$ for~$A_E$.

  It remains to check purity. By~\cite[Thm.\ 4.2.2]{MR1293976}, it is enough to
  check purity for the Weil--Deligne representations associated to
  1-motives with potentially good reduction, which is~\cite[Prop.\
  4.6.1, Prop.\ 4.7.4]{MR1293976}.

  The above is of course not a historically accurate account of a
  proof; indeed, the strict compatibility of the compatible system at
  places not dividing~$l$ is
  stated in~\cite[Ex.\ 8.10]{MR0349635}, and given Fontaine's
  definition of the Weil--Deligne representation associated to a
  potentially semistable representation, the entire proposition can be deduced from the results
  of~\cite{MR0354656}. We omit the details, but we would like to thank
  Brian Conrad for explaining them to us. \end{proof}

\begin{defn}\label{defn: Galois repn associated to abelian surface}If~$A/F$ is an abelian surface,
  then we write~$\rho_{A,l}$ for $H^1(A_{\overline{F}},\Qlbar)$, and $\cR_A$
  for the compatible system~$\{\rho_{A,l}\}$. We can think
  of~$\rho_{A,l}$ as a representation
  $\rho_{A,l}:G_F\to\GSp_4(\Qlbar)$ with multiplier~$\varepsilon_l^{-1}$, and will frequently do so without
  comment.
\end{defn}

\begin{rem} \label{rem: compatible systems of GSp4 representations}It will sometimes be
  convenient to say that a set of $\GSp_4$-valued representations form
  a compatible system, by which we simply mean that the corresponding
  $\GL_4$-valued representations form a compatible system. In
  particular, the representations $\rho_{A,l}:G_F\to\GSp_4(\Qlbar)$
  considered in Definition~\ref{defn: Galois repn associated to
    abelian surface} form a compatible system in this sense. (In
  general, one might wish to ask for a compatibility between the
  symplectic structures; such a compatibility always holds in the
  cases that we consider, and in particular we will only consider
  representations whose multiplier character is the inverse cyclotomic
  character, so we ignore this point.)
\end{rem}

We can define the
$L$-function of~$\cR$  as follows: \[ L( \CR,s) = \prod_{v\ndiv \infty} L( \WD_v(\CR), s). \]
Furthermore, if~$\CR$ comes from an abelian variety (or more generally, 
arises in a geometric structure where the Hodge structure is apparent)
then (as in~\cite{Serre1969-1970})
we can
define Gamma factors $L_v(\CR,s)$ for each place~$v|\infty$
of~$F$, and we set \numequation\label{eqn: completed L function
  motivic version} \Lambda( \CR,s)=L( \CR,s)
  \prod_{v|\infty} L_v(\CR,s).  \end{equation}
 In particular, if~$\CR$ arises from an abelian surface over a totally real field~$F$,
  then the corresponding Gamma factor is given by~$L_v(\CR,s) =   \Gamma_{\C}(s)^2$ for all~$v|\infty$ where~$\Gamma_{\C}(s) 
  =  (2 \pi)^{-s} \Gamma(s)$.

   We also have 
   a conductor~$N(\cR)$
which is a product of local factors depending only on the
$\WD_v(\cR)$. 
Conjecturally, if~$\CR$ is a strictly compatible system, then~$\Lambda( \CR,s)$ admits a meromorphic continuation to the
entire complex plane and satisfies a functional equation of the form
\numequation\label{eqn: functional equation in automorphic case} \Lambda( \CR,s)=\varepsilon( \CR)N(\cR)^{-s} \Lambda(
  \CR^\vee, 1-s) \end{equation} 
  for some factor~$\varepsilon( \CR)$. (When~$\CR$ arises geometrically, there
  are natural definitions of the epsilon factor~$\varepsilon(\CR)$, but it is not
  immediately apparent how to read off~$\varepsilon(\CR)$ directly from the compatible system.)

In particular, if~$A/F$ is an abelian variety, then by Proposition~\ref{prop: pure
  abelian variety compatible system}
\[\Lambda_i(A,s):=\Lambda(H^i(A_{\overline{F}},\Qlbar),s)\] is
well-defined, and we
define the completed Hasse--Weil zeta function of~$A$ to
be \[\Lambda(A,s):=\prod_{i=0}^{2\dim
    A}\Lambda_i(A,s)^{(-1)^i}.\]
Note that if~$v$ is a finite place of~$F$ at which $A$ has good
reduction with corresponding reduction~$\overline{A}$, then the local $L$-factor \[L_v(A,s):=\prod_{i=0}^{2\dim
    A}L(\WD(H^i(A_{\overline{F}},\Qlbar)),s)^{(-1)^s}\] can be written
as \[L_v(A,s)=\exp\left(\sum_{m=1}^\infty\frac{\#\overline{A}(k(v)_m)}{m}{\#k(v)}^{-ms}
  \right) \]where~$k(v)$ is the residue field of~$F_v$ and
$k(v)_m/k(v)$ is the extension of degree~$m$.

We have the following conjectures for the~$\Lambda_i(A,s)$, which we
will prove for abelian surfaces over totally real fields by showing
that they are potentially automorphic.
\begin{conj}[\cite{Serre1969-1970}, Conj.\ C9] \label{conj: L function meromorphic continuation}For each~$i$, $\Lambda_i(A,s)$ has a meromorphic
  continuation to the entire complex plane, and satisfies a functional
  equation of the form \[\Lambda_i(A,s)=wN^{\frac{i+1}{2}-s} \Lambda_i(A,i+1-s)\]
  where~$w=\pm 1$ and $N\in\Z_{\ge 1}$.
\end{conj}\begin{cor}
  \label{cor: HW zeta function meromorphic continuation}If
  Conjecture~\ref{conj: L function meromorphic continuation} holds,
  then $\Lambda(A,s)$ has a meromorphic
  continuation to the entire complex plane, and satisfies a functional
  equation of the form $\Lambda(A,s)=\varepsilon N^{-s} \Lambda(A,1+\dim
  A-s)$ where~$\varepsilon\in \R$ and~$N\in\Q_{>0}$.\end{cor}\begin{proof}
  This follows immediately from Conjecture~\ref{conj: L function meromorphic
    continuation} by Poincar\'e duality.\end{proof}

\subsection{Arthur's classification}\label{subsec: Arthur classification}We now recall some consequences of Arthur's
classification~\cite{MR2058604} of discrete
automorphic representations of~$\GSp_4$. The analogous classifications
for~$\Sp_4$ and~$\mathrm{SO}_5$ are special cases of the very general results
proved in~\cite{MR3135650}, and a proof of the classification
announced in~\cite{MR2058604}, making use of the results and
techniques of~\cite{MR3135650} is given
in~\cite{GeeTaibi}. This reference establishes the compatibility of Arthur's
classification with the local Langlands correspondence~$\recGT$, which
we use below without further comment.

We say that an
automorphic representation~$\pi$ of~$\GSp_4(\A_F)$
is discrete if it occurs in the discrete
spectrum of the $L^2$-automorphic forms (with
fixed central character~$\omega=\omega_\pi$). Note in particular that
all cuspidal automorphic representations are discrete. Arthur's classification divides the discrete
spectrum into six families of automorphic
representations. We will not need the full details
of this classification, but rather just some 
consequences that we now recall.

If~$\Pi$ is a cuspidal automorphic representation
    of~$\GL_4(\A_F)$, then we say that~$\Pi$ is of \emph{symplectic
    type} with multiplier~$\chi$ if
  the partial $L$-function
    $L^S(s,\Pi,\bigwedge^2\otimes\chi^{-1})$ has a pole at $s=1$
    (where $S$ is any finite set of places of $F$). Note that this implies in particular that
  $\Pi\cong\Pi^\vee\otimes\chi$.

We say that a discrete automorphic representation~$\pi$ of~$\GSp_4(\A_F)$ is \emph{of
general type}  in the sense
    of~\cite{MR2058604} if there
    is a cuspidal automorphic
    representation~$\Pi$ of $\GL_4(\A_F)$ of symplectic  type with multiplier~$\omega_\pi$
    such that for each place~$v$ of~$F$, the $L$-parameter obtained from
    $\recGT(\pi_v)$ by composing with the usual embedding
    $\GSp_4\into\GL_4$ is $\rec(\Pi_v)$. We say that~$\Pi$ is the \emph{transfer} of~$\pi$. 

In practice, all of the automorphic representations~$\pi$ that we
consider in our main arguments will be of general type. We will often use the
following lemma to guarantee this. (For example, the lemma will be
used to show that when we localize a cohomology group at a non-Eisenstein
maximal ideal, the only automorphic representations that
contribute are of general type.)
    \begin{lem}
      \label{lem: getting to general type}Suppose that~$F$ is totally
      real, and that~$\pi$ is a discrete
  automorphic representation
  of~$\GSp_4(\A_F)$, and
  that at each place~$v|\infty$, $\pi_v$
has the same infinitesimal character as the representations in the
$L$-packet corresponding to
  $\varphi_{(2;k_v-1,l_v-2)}$ with~$k_v\equiv l_v\pmod{2}$ and~$k_v\ge
  l_v\ge 2$. Suppose that $\pi$ is not of general type. 

Then there is a compatible system of \emph{reducible} Galois
  representations~$\rho_{\pi,p}:G_F\to\GSp_4(\Qpbar)$ such that
  for all but finitely many places~$v$ of~$F$, we have  $\WD(\rho_{\pi,p}|_{G_{F_v}})^{\semis}\cong
    \recGTp(\pi_v\otimes|\nu|^{-3/2})^{\semis}$.
    \end{lem}
    \begin{proof}
      We follow the proof of ~\cite[Thm.\ 7.11]{CGGSp4}. Since~$\pi$
      is not of general type, $\pi$ falls
      into one of the five classes  (b)-(f) listed at the end
      of~\cite{MR2058604}. In cases~(e) and~(f), we see that the Hecke
      parameters of~$\pi$ agree with those of a direct sum of 4 idele
      class characters. By the hypothesis on the infinitesimal
      character, these characters are algebraic, so we may take the
      direct sum of the corresponding compatible systems of Galois
      representations.

      In case~(d), the Hecke parameters of~$\pi$ agree with those of
      an isobaric direct sum of the form
      $\lambda|\cdot|^{1/2}\boxplus\lambda|\cdot|^{-1/2}\boxplus\mu$,
      where~$\lambda$ is an idele class character, and~$\mu$ is a
      cuspidal automorphic representation of~$\GL_2(\A_F)$, satisfying
      $\omega_\mu=\lambda^2=\omega_\pi$. Considering infinitesimal
      characters, we see that~$\lambda$ is algebraic, so that
      $\lambda|\cdot|^{1/2}\boxplus\lambda|\cdot|^{-1/2}$ is regular
      algebraic. This implies that~$\mu$ is also regular algebraic,
      and thus has an attached compatible system of Galois
      representations.

      In case~(b), the Hecke parameters of~$\pi$ agree with those of
      an isobaric direct sum of the form $\mu_1\boxplus\mu_2$, where
      $\mu_1\ne\mu_2$ are cuspidal automorphic representations
      of~$\GL_2(\A_F)$ with central character~$\mu_\pi$. Since their
      central characters agree, it follows easily that they both
      correspond to holomorphic Hilbert modular eigenforms of
      paritious weight. Finally in
      case~(c), the Hecke parameters of~$\pi$ agree with those of an
      isobaric direct sum of the form
      $\mu|\cdot|^{1/2}\boxplus\mu|\cdot|^{-1/2}$, where~$\mu$ is a
      cuspidal automorphic representation of~$\GL_2(\A_F)$ of
      orthogonal type; that is,
      it is induced from a quadratic extension of~$F$. Since~$\mu$ is
      certainly algebraic, we again have an attached compatible system
      of reducible Galois representations, as required.
  \end{proof}\begin{rem}
    \label{rem: in general type case we expect Galois representations
      to be irreducible} 
      Suppose that~$\pi$ is of general type but otherwise satisfies
      the conditions of Lemma~\ref{lem: getting to general type}. 
      Then the corresponding Galois
      representations constructed in~\cite{MR3200667} (see also 
      Theorem~\ref{thm: existence and properties of Galois representations for
    automorphic representations general Mok version}) give rise to
      a compatible system of Galois representations 
      which ---  in contrast to those occurring in
       Lemma~\ref{lem: getting to general type} --- are  expected to
       always be irreducible.
\end{rem}
The following theorem summarizes the consequences that we need from
  Arthur's multiplicity formula.
\begin{thm}
  \label{thm: arthur classification
    results}Suppose that~$F$ is a totally
  real field, and that~$\Pi$ is a cuspidal automorphic representation
  of~$\GL_4(\A_F)$ of symplectic type with multiplier~$\chi$. Then
  there exists at least one  discrete
  automorphic representation~$\pi$
  of~$\GSp_4(\A_F)$ with central character
  $\chi$ such that~$\Pi$ is the transfer of~$\pi$. 

More precisely,  for each place $v$ of~$F$, let~$\pi_v$ be an
      element of the $L$-packet corresponding to  $(\rec_p(\Pi_v),\chi_v)$. Then $\pi:=\otimes_v'\pi_v$ is
      automorphic, and occurs with multiplicity one in the discrete spectrum.

If, furthermore, $\Pi$ is algebraic, then~$\pi$ is cuspidal.
\end{thm}
\begin{proof}
The statements of the first two paragraphs are immediate from the
  multiplicity formula of~\cite{MR2058604} as proved in~\cite{GeeTaibi} (note that since~$\pi$ is
  of general type by definition, the group~$\cS_\psi$ considered
  in~\cite{MR2058604} is trivial). Suppose then that~$\Pi$ is algebraic;
  then~$\Pi_\infty$ is essentially tempered by~\cite[Lem.\
  4.9]{MR1044819}, so that~$\pi_\infty$ is also essentially tempered
  (as its $L$-parameter is essentially bounded), so that $\pi$ is
  cuspidal by~\cite[Thm.\ 4.3]{MR733320}.
\end{proof}

\subsection{Balanced modules}\label{subsec: balanced}Let~$S$ be
a Noetherian local ring with residue field~$k$, and let~$M$ be a
finitely generated $S$-module. As in~\cite[\S 2.1]{CG}, we define the
\emph{defect} $d_S(M)$ to be \[d_S(M):=\dim_k M/\m_S
  M-\dim_k\Tor_S^1(M,k).\]

\begin{defn}
  \label{defn: balanced module}We say that~$M$ is \emph{balanced}
  if~$d_S(M)\ge 0$.
\end{defn}

\begin{lem}
  \label{lem: presentations of balanced modules}If~$M$ is balanced,
  then there is a presentation \[S^d\to S^d\to M\to 0\] with~$d=\dim_k
  M/\m_S M$.

Conversely if $M$ admits a presentation \[S^r\to S^r\to M\to 0\] for
some~$r\ge 0$, then~$M$ is balanced.
\end{lem}
\begin{proof}
  Assume firstly that~$M$ is balanced, and choose a (possibly infinite) minimal resolution \[\dots\to P_i\to\dots\to
    P_1\to P_0\to M\to 0\]by finite free $S$-modules $P_i$ of
  rank~$r_i$. (Recall that a minimal resolution is one whose
  differentials vanish modulo~$\m_S$, and that such a resolution
  always exists.) Tensoring this resolution with~$k$ over~$S$, we see
  that $r_i=\dim_k \Tor^i_S(M,k)$, so that in particular by our
  assumptions we have $d=r_0\ge r_1$, so that there is a presentation
  of the form $P_1\oplus S^{\oplus(d-r_1)}\to P_0\to M\to 0$, as
  required.

Conversely, if  $M$ admits a presentation $S^r\to S^r\to M\to 0$, then
let~$K$ be the image of the map $S^r\to S^r$. Then from the exact
sequence \[0\to\Tor_S^1(M,k)\to K/\m_S K\to k^r\to M/\m_S M\to 0\] we
see that \[d_S(M)=r-\dim_k K/\m_S K;\]since $K$ admits a surjection
from~$S^r$, it follows that $d_S(M)\ge 0$, as required.
\end{proof}

\subsection{Projectors}\label{subsec: projectors}

Let $R$ be a complete local Noetherian ring with maximal ideal $\mathfrak{m}_R$ and finite residue field. We let $\mathbf{Mod}^{\comp}(R)$ be the category of $\mathfrak{m}_R$-adically complete and separated $R$-modules. Let $M \in \Ob( \mathbf{Mod}^{\comp}(R))$ and $T \in \mathrm{End}_R(M)$. 
\begin{defn} We say that $T$ is locally finite on $M$ if for all $n \geq 0$, $M/\mathfrak{m}_R^n$ is an inductive limit of finite type $R$-modules which are stable under the action of $T$. 
\end{defn}

\begin{lem}
  \label{lem: commuting locally finite is locally finite}If $T_1$,
  $T_2$ commute and are both locally finite on~$M$, then $T_1T_2$ is
  also locally finite on~$M$. 
\end{lem}
\begin{proof}
  By definition we can assume that $M$ is $\m_R^n$-torsion for
  some~$n$. If $v\in M$ then since $T_1$ is locally finite, the
  $R$-submodule of~$M$ generated by the $T_1^iv$ is finitely
  generated. Since $T_2$ is locally finite, it follows that 
  the $R$-submodule generated by the $T_2^iT_1^jv$ is also finitely
  generated, and since $T_1$, $T_2$ commute, this submodule is stable
  under the action of~$T_1T_2$, as required.
\end{proof}
The following results from~\cite{pilloniHidacomplexes} will be used to
construct the ordinary projectors associated to certain Hecke
operators.

\begin{lem}[{\cite[Lem.\ 2.1.2]{pilloniHidacomplexes}}]\label{lem:
    can check locally finiteness on special fibre}If~$M$ is an object
  of~$\mathbf{Mod}^{\comp}(R)$ and~$T$ is an endomorphism of~$M$, then
  $T$ is locally finite on~$M$ if and only if it is locally finite on~$M/\m_R$.  
\end{lem}
\begin{lem}[{\cite[Lem.\ 2.1.3]{pilloniHidacomplexes}}]\label{lem:
    ordinary projector exists} If $T$ is
  locally finite on $M$, then $\lim_{n \rightarrow \infty} T^{n!}$
  converges pointwise in the $\m_R$-adic topology to a projector
  $e(T)$ on $M$. 

The operators~$T$ and~$e(T)$ commute, and we have a
  $T$-stable decomposition \[M=e(T)M\oplus(1-e(T))M,\]where $T$ is
  bijective on~$e(T)M$ and topologically nilpotent on $(1-e(T))M$.
\end{lem}

We call $e(T)$ the \emph{ordinary projector} attached to $T$.   Let $\mathbf{D}(R)$ be the derived category of $R$-modules, let $\mathbf{D}^{\fflat}(R)$ be the full subcategory of
$\mathbf{D}(R)$ generated by bounded complexes of flat, $\mathfrak{m}_R$-adically complete and  separated 
$R$-modules and let $\mathbf{D}^{\mathrm{perf}}(R)$ be the full subcategory of
$\mathbf{D}(R)$ generated by bounded complexes of finite free
$R$-modules. Let $M \in \Ob( \mathbf{D}^{\fflat}(R))$. We say that an
operator $T \in \mathrm{End}(M)$ is locally finite if there is a
bounded complex of flat modules $N$ representing  $M$ and an operator
$T_0 \in \mathrm{End}(N)$ representing $T$ which is degree-wise locally finite. By \cite[Lem.\ 2.3.1]{pilloniHidacomplexes}, $T$ is locally finite on $M$ if and only if $T$ is locally finite on the cohomology groups $H^i(M \otimes_R^L R/\mathfrak{m}_R)$ and there is a
bounded complex of flat modules $N$ representing  $M$ and an operator
$T_0 \in \mathrm{End}(N)$ representing $T$. Given a choice of representatives $(N,T_0 \in \mathrm{End}(N))$ for a locally finite operator $T$, we get an associated idempotent $e(T_0) \in \mathrm{End}(N)$. 
In general, we do not know whether
two choices of representatives~$(N,T_0 \in \mathrm{End}(N))$ give the same projector in $\mathrm{End}_{\mathbf{D}(R)}(M)$. But by~\cite[Lem.\
2.3.2]{pilloniHidacomplexes}, if we assume that for one choice of representative $e(T_0) M$ is an object of $\mathbf{D}^{perf}(R)$ then, for another  choice of locally finite representative $(N',T_1 \in \mathrm{End}(N'))$, $e(T_1) M$ is an object of $\mathbf{D}^{perf}(R)$ and there is a canonical quasi-isomorphism
 $e(T_0)M \rightarrow e(T_1) M$. In the sequel, these conditions will always be satisfied and we will write $e(T)$ by abuse of notation.

\section{Shimura varieties}\label{sec: integral
  Shimura varieties}In this section, we
discuss the Hilbert--Siegel Shimura varieties that we work with, and
some properties of their integral models.  There are two closely related algebraic groups here: $G_1 = \mathrm{Res}_{F/\Q} \mathrm{GSp}_4$ and its subgroup $G$ of elements with similitude factor in $\mathbb{G}_m \hookrightarrow \mathrm{Res}_{F/\Q} \mathbb{G}_m$. 

The group $G$ admits a standard PEL Shimura variety and there is a
good moduli interpretation, integral models, and a good theory of
integral compactification. Nonetheless, from an automorphic view point
we must work with the group $G_1$ which gives rise to a Shimura
variety of abelian type.  

Going back to the work of Deligne (see in
particular~\cite[\S 2.7]{MR546620}), there is a standard strategy for
handling abelian type Shimura varieties by relating their connected
components to quotients of connected components of Hodge type Shimura
varieties by finite groups. As a particular instance of this strategy,
  the Shimura varieties for $G$ and $G_1$ are closely related: the connected components of  $G_1$-Shimura varieties are quotients of the connected components of  $G$-Shimura varieties  by finite groups. We therefore  study both of them at the same time. 

For convenience, our main references for integral models of  PEL Shimura
varieties and their compactifications are the papers~\cite{MR3186092,lan_2016,Lan2016}, although
some of the results we cite from there were proved in earlier papers,
in particular~\cite{MR1124982}; we refer the reader
to the references in~\cite{MR3186092}  for a more detailed
historical account.
\subsection{Similitude groups}
Let $F$ be a totally real field. Let $V = \CO_F^4$ be a free $\CO_F$-module of rank $4$. We equip $V$ with the symplectic $\CO_F$-linear form $ <, >_1: V \times V \rightarrow \CO_F$ given by the matrix $J$. We let $<,>: = (\mathrm{Tr}_{F/\Q} \circ <,>_1)$ be the associated $\Z$-linear symplectic form. 

Let $G_1 = \mathrm{Res}_{F/\Q} \mathrm{GSp}_4$ be the algebraic group of symplectic $F$-linear automorphisms of $(V_{\Q}, <,>_1)$, up to a similitude factor $\nu$ in $\mathrm{Res}_{F/\Q} \mathbb{G}_m$.  

Let $G \subset G_1$ be the algebraic group of symplectic $F$-linear
automorphisms of $(V_\Q, <,>)$ up to a similitude factor in
$\mathbb{G}_m$; that is, $G = G_1 \times_{\nu, \mathrm{Res}_{F/\Q} \mathbb{G}_m} \mathbb{G}_m$. 
\subsection{Shimura varieties over
  \texorpdfstring{$\C$}{C}} \label{sect-overc}We firstly briefly discuss some Shimura varieties over~$\C$. We
caution the reader that in the bulk of the paper we will work with
Shimura varieties over~$\Z_{(p)}$ which are not quite integral models
of these Shimura varieties, but whose geometrically connected
components are the same as these; see Proposition~\ref{prop: existence
  of integral model for Y} below for a precise statement. We begin by recalling the
definition of a neat compact open subgroup from~\cite[Defn.\ 1.4.1.8]{MR3186092}.
\begin{defn}\label{defn:neat}Write $g=(g_l)_l\in
  G_1(\A^{\infty})$, and for each~$l$, write $\Gamma_{g_l}$ for the
  subgroup of $\Qlbartimes$ generated by the eigenvalues of~$g_l$
  (under any faithful linear representation of~$G_1$). Then we say
  that $g$ is \emph{neat}
  if \[\bigcap_l(\Qbar^\times\cap\Gamma_{g_l})_{\mathrm{tors}}=1. \]
  Similarly, if $g\in G_1(\A^{\infty,p})$, then we say that $g$
  is neat if \[\bigcap_{l\ne p}(\Qbar^\times\cap\Gamma_{g_l})_{\mathrm{tors}}=1. \]
We say that a compact open subgroup $K \subset G_1(\A^{\infty})$
(resp.\ $K^p\subset G_1(\A^{\infty,p})$) is \emph{neat} if all of its elements are neat.
\end{defn}
We consider  
the Shimura variety
associated to the group $G_1$ and a neat compact open subgroup $K \subset G_1(\A^{\infty})$:
$$ S^{G_1}_K (\mathbb{C}) = G_1(\Q) \backslash \big( G_1(\R) \times G_1(\A^{\infty}) \big)/ Z(\R)^0 K_\infty^0 K$$
where $Z(\R)^0 \simeq \R_{>0}^{\HomFR}$ is the connected component of the
centre in $G_1(\R) \simeq \mathrm{GSp}_4(\R)^{{\HomFR}}$ and $K_\infty^0$ is
the connected component of the maximal compact subgroup inside
$G_1(\R)$, so that $K_\infty^0$
 is a product of copies of~$U(2)$. This Shimura variety
carries a natural structure of complex quasi-projective variety, as we have $G_1(\R)/ Z(\R)^0 K_\infty^0 = (\mathcal{H} \cup - \mathcal{H} )^{\HomFR}$, where $\mathcal{H}$ is the Siegel half space of symmetric matrices $M = A + i B \in \mathrm{M}_{2 \times 2}(\C)$ with $B$ positive definite. 

Let $G_1(\Q)^+$ be the subgroup of $G_1(\Q)$ equal to
$\nu^{-1}(F^{\times, +} )$, where $F^{\times, +}$ is the subgroup of
totally positive elements in $F^\times$.  Then by strong
approximation, \[ G_1(\A^{\infty}) =  \coprod_{c}G_1(\Q)^+ c K\]
where $c$ runs through a (finite) set of elements in $G_1(
\A^{\infty})$ such that $\nu(c)$ are representatives of the strict
class group $F^{\times, +} \backslash (\A^{\infty} \otimes_{\Q}
F)^\times / \nu(K)$. 

One can then write $$S^{G_1}_K (\mathbb{C}) = \coprod_{c}  \Gamma_1(c, K) \backslash \mathcal{H}^{\HomFR} $$ where $\Gamma_1(c,K) = G_1(\Q)^+  \cap c K c^{-1}$. 
 
This Shimura variety, although natural from the point of view of
automorphic forms, is not of PEL type. Therefore, it is also necessary
to work with another Shimura variety.  We can consider the double quotient
$$S^G_{K}(\C) =  G(\Q) \backslash \big( G(\R) \times G_1(\A^{\infty})
\big)/ \R_{>0} K_\infty^0 K;$$ this is not strictly speaking a Shimura
variety, and in particular we emphasise that it is not the PEL Shimura variety associated to $G$!
By strong approximation we may write \[ G_1(\A^{\infty}) =  \coprod_{c}G(\Q)^+ c K\] where $c$ runs through a set of elements of $G_1( \A^{\infty})$ such that $\nu(c)$ are representatives of the infinite set  $\Q^{\times, +} \backslash (\A^{\infty} \otimes_{\Q} F)^\times / \nu(K)$. For all $c$, we consider the group $\Gamma(c,K) =   G(\Q)^+ \cap
cKc^{-1}$, so that \[S^G_K(\C)  = \coprod_c  \Gamma(c, K) \backslash \mathcal{H}^{\HomFR}.\]

The inclusion $G(\R)\into G_1(\R)$ induces  a
natural surjective map $S^G_K(\C) \rightarrow S^{G_1}_K(\C)$.  On
connected components, it induces the natural map \[\Q^{\times, +} \backslash (\A^{\infty} \otimes_{\Q} F)^\times / \nu(K) \rightarrow F^{\times, +} \backslash (\A^{\infty} \otimes_{\Q} F)^\times / \nu(K).\]

For any $c \in G_1(\A^{\infty})$ we have an associated surjective map on the
connected components corresponding to $c$, given by 
\[\Gamma(c,K) \backslash \mathcal{H}^{\HomFR} \rightarrow \Gamma_1(c,K) \backslash \mathcal{H}^{\HomFR}.\]

Let $Z(\Gamma_1(c,K)) \subset \Gamma_1(c,K)$ be the centre. Then $Z(\Gamma_1(c,K))$ is a finite index subgroup of $\ocal_F^{\times}$ that we denote by $\ocal_F^{\times}(K)$.  Let \[\Delta(K) =  \Gamma_1(c,K) / \big(\ocal_F^{\times}(K), \Gamma(c,K) \big).\] 

This is a finite  group, independent of $c$ and  isomorphic to $$\nu(\Gamma_1(c,K))/  \nu (\ocal_F^{\times}(K)) = ( \ocal_F^{\times, +} \cap \nu(K))/ \nu(\ocal_F^{\times}(K)),$$
having noted the following:
   \begin{lem} There is an equality~$\nu(\Gamma_1(c,K)) = \ocal_F^{\times, +} \cap \nu(K)$.
\end{lem}

\begin{proof}Recall that by definition $\Gamma_1(c,K) = G_1(\Q)^+  \cap c K
  c^{-1}$, so certainly we have an inclusion~$\nu(\Gamma_1(c,K))\subseteq
  \cO_F^{\times,+}\cap\nu(K)$. Conversely, suppose that $\nu(\gamma)
  = x \in  \OL^{\times,+}_F \cap  \nu(K)$ for some element~$\gamma \in
  cKc^{-1}$. Since $x\in\cO_F^{\times,+}=\nu(G_1(\Q)^+)$, we can
  choose~$g\in G_1(\Q)^+ $ with
  $\nu(g)=x$. Then~$\nu(\gamma^{-1}g)=1$, so by strong approximation, we may find an~$h \in G_1(\Q)^{+}$
with trivial similitude character such that~$h$ is arbitrarily close to~$g \gamma^{-1}$, and in particular
close enough that $h\gamma g^{-1} $ lies in~$cKc^{-1}$. Then~$h \gamma \in  G_1(\Q)^{+} \cap cKc^{-1}$ and has similitude
character~$x$, as required.
\end{proof}

We also have

 \begin{lemma}\label{lem: Delta Pilloni lemma over C} The  map $\Gamma(c,K) \backslash \mathcal{H}^{\HomFR} \rightarrow \Gamma_1(c,K) \backslash \mathcal{H}^{\HomFR}$ is finite \'etale with group $\Delta(K)$. 
\end{lemma}
\begin{proof} The group $\Gamma_1(c,K)$ acts through its quotient
   $\Gamma_1(c,K)/\ocal_F^{\times}(K)$ on $\mathcal{H}^{\HomFR}$, and
   since~$K$ is neat, this action is free. 
\end{proof}

\subsection{Integral models of Shimura varieties}\label{subsec:
  integral models}We now introduce the integral models of Shimura varieties that we will
consider in the rest of the paper.

\subsubsection{Compact open subgroups at $p$} We let $p$ be a prime that is totally split in $F$. Let $v$ be a prime ideal in $\CO_F$ above $p$. Consider the following chain of $\CO_{F_{v}}$-sub modules of $F_v^4$:

$$ V_0  \rightarrow V_1 \rightarrow V_2 \rightarrow V_3 \rightarrow V_4$$
where $V_0 = V\otimes_{\CO_F} \CO_{F_v} = \oplus_{i=1}^4 \CO_{F_v} e_i$ and $V_j = \langle  p^{-1} e_1, \cdots, p^{-1} e_j, e_{j+1}, \cdots, e_4 \rangle$. 
We can identify $V_0$ and $V_4$ through multiplication by $p$ and
sometimes think of the indices as being in $\Z/4\Z$. 

From the perfect pairing $<,>$ on $V_0$  we obtain  perfect pairings on  $V_2 \times V_2$ and on $V_1 \times V_3$.  

We now recall the definitions of the parabolic subgroups that we use
in terms of flags; this description is well suited to the definitions
of our integral models.
\begin{itemize}
\item  $\mathrm{GSp}_4( \CO_{F_v} ) = \mathrm{Aut} (V_0) \cap
  \mathrm{GSp}_4(F_v)$ (the hyperspecial subgroup),  
\item $\Par(v) = \mathrm{Aut} ( V_1 \rightarrow V_3)
  \cap \mathrm{GSp}_4(F_v)$ (the paramodular subgroup),
\item $\Si(v) = \mathrm{Aut} (V_0 \rightarrow V_2)
  \cap \mathrm{GSp}_4(F_v)$ (the Siegel parahoric),
\item $\Kli(v) = \mathrm{Aut} ( V_0 \rightarrow V_1
  \rightarrow V_3 \rightarrow V_0) \cap \mathrm{GSp}_4(F_v)$ (the
  Klingen parahoric),
\item $\Iw(v) = \mathrm{Aut} (V_0 \rightarrow V_1
  \rightarrow V_2 \rightarrow V_3 \rightarrow V_0) \cap
  \mathrm{GSp}_4(F_v)$ (the Iwahori subgroup).
\end{itemize}

\subsubsection{The moduli problem}
Let $ALG/{\Z_{(p)}}$ be the category of Noetherian $\Z_{(p)}$-algebras
and $\AFF/ \Z_{(p)}$ the opposite category. Let $K
\subset G_1( \A^{\infty})$ be a compact open subgroup; we will also
refer to such a compact open subgroup as a \emph{level structure}.
\begin{defn}\label{defn: reasonable} We say that a level structure  $K = K^p K_p$ is
  \emph{reasonable} if~$K^p\subset G(\A^{\infty,p})$ is neat, and if $K_p =  \prod_{v |p} K_v$
  where for each $v|p$ we have $$K_v \in \{\mathrm{GSp}_4( \CO_{F_v}), \Par(v), \Si(v), \Kli(v), \Iw(v)\}.$$ 
\end{defn}

Let $K$ be a reasonable level structure.  We consider the  groupoid
$\mathbb{Y}_K$ over $\AFF/\Z_{(p)}$ whose fibre over $S=\Spec R \in \Ob(
\AFF/\Z_{(p)})$ is the category with objects  $(A, \iota, \lambda,
\eta, \eta_p)$, where:
\begin{enumerate}
\item $A \rightarrow \Spec R$ is an abelian scheme,
\item $\iota: \CO_F \rightarrow \mathrm{End} (G) \otimes \Z_{(p)}$ is an action,
\item $\mathrm{Lie} (A)$ is a locally free $\CO_F \otimes_{\Z} R$-module of rank $2$, 
\item $\lambda: A \rightarrow A^t$ is a prime to $p$, $\cO_F$-linear quasi-polarization such that for all $v |p$, $\mathrm{Ker}(\lambda: A[v^\infty]
  \rightarrow A^t[v^\infty])$ is trivial if $K_v \neq
  \Par(v)$ and is an order $p^2$ group scheme if $K_v = \Par(v)$,
\item $\eta$ is a $K^p$-level structure, 
\item $\eta_p$ is a $K_p$-level structure. 
\end{enumerate}

Here by a prime to $p$ quasi-polarization $\lambda:A\to A^t$ we mean a
$\Z_{(p)}^\times$-polarization in the sense of~\cite[Defn.\
1.3.2.19]{MR3186092}. By a $K_p$-level structure~$\eta_p$, we mean the following list of data: 
\begin{enumerate}
\item For all $v|p $ such that $K_v = \Kli(v)$, $H_v \subset A[v]$ is an order $p$-group scheme, 
\item For all $v |p $ such that $K_v = \Si(v)$, $L_v \subset A[v]$ is an order $p^2$ group scheme that is totally isotropic for the Weil pairing. 
\item For all $v|p $ such that $K_v= \Iw(v)$, $H_v \subset L_v \subset A[v]$ are subgroups such that $H_v$ is of order $p$, $L_v$ is of order $p^2$ and $L_v$ is totally isotropic for the Weil pairing. 
\end{enumerate}

Let us spell out the definition of $K^p$-level structure. We may
assume without loss of generality that $S$ is connected, and we
fix $\overline{s}$ a geometric point of $S$. The adelic Tate module
${\HH}_1( A\vert_{\overline{s}}, \A^{\infty,p})$ carries a
symplectic Weil pairing \[<,>_{\lambda}: {\HH}_1(
  A\vert_{\overline{s}}, \A^{\infty,p}) \times {\HH}_1(
  A\vert_{\overline{s}}, \A^{\infty,p}) \rightarrow {\HH}_1(
  \mathbb{G}_m\vert_{\overline{s}}, \A^{\infty,p})\] or equivalently
an $F$-linear symplectic pairing: \[<,>_{1,\lambda}: {\HH}_1( A\vert_{\overline{s}}, \A^{\infty,p}) \times {\HH}_1( A\vert_{\overline{s}}, \A^{\infty,p}) \rightarrow {\HH}_1( \mathbb{G}_m\vert_{\overline{s}}, \A^{\infty,p}) \otimes F.\] 
The level structure $\eta$ is  a  $K^p$-orbit of pairs of isomorphisms~$(\eta_1,\eta_2)$, where
(with~$V=\cO_F^4$ the standard symplectic space defined above):
\begin{enumerate}
\item An $\ocal_F\otimes_{\ZZ} \A^{\infty,p}$-linear isomorphism of $\Pi_1(S, \overline{s})$-modules $\eta_1:   V \otimes_{\ZZ} \A^{\infty,p}  \simeq {\HH}_1( A\vert_{\overline{s}}, \A^{\infty,p})$.
\item An  $\ocal_F\otimes_{\ZZ} \A^{\infty,p}$-linear isomorphism of
  $\Pi_1(S, \overline{s})$-modules $\eta_2:   F \otimes_{\ZZ}
  \A^{\infty,p}  \simeq F \otimes_{\ZZ} {\HH}_1(
  \mathbb{G}_m\vert_{\overline{s}}, \A^{\infty,p})$. 
\end{enumerate}
We moreover impose that the following diagram is commutative:
\numequation\label{eqn: pairings eta commute one}
\xymatrix{ V \otimes_{\ZZ} \A^{\infty,p}   \times V \otimes_{\ZZ} \A^{\infty,p}   \ar[rr]^{\eta_1 \times \eta_1} \ar[d]^{<,>_1}& & {\HH}_1( A\vert_{\overline{s}}, \A^{\infty,p}) \times {\HH}_1( A\vert_{\overline{s}}, \A^{\infty,p})  \ar[d]^{<,>_{1,\lambda}}\\ 
 F \otimes_{\ZZ} \A^{\infty,p}   \ar[rr]^{\eta_2} & &F \otimes_{\ZZ} {\HH}_1( \mathbb{G}_m\vert_{\overline{s}}, \A^{\infty,p})}
\end{equation}
The action of an element~$k\in K^p$ takes~$(\eta_1,\eta_2)$
to~$(\eta_1k,\nu(k)\eta_2)$.\begin{rem}
  The reader will observe that~$\eta_2$ is uniquely determined
  by~$\eta_1$, but we find it convenient to record it as part of the
  data for the sake of comparison to the PEL setting in
  Proposition~\ref{prop: existence of integral model for Y} below.
\end{rem}

A map between  quintuples $(A, \iota, \lambda, \eta,\eta_p)$ and $(A',
\iota', \lambda', \eta',\eta'_p)$ is an $\CO_F$-linear prime to $p$
quasi-isogeny (in the sense of~\cite[Defn.\ 1.3.1.17]{MR3186092}) $f: A \rightarrow A'$  such that
\begin{itemize}
\item $f^{*} \lambda = r \lambda'$ for a locally constant function
  $r: S \rightarrow \Z_{(p)}^{\times,+}$, 
\item $f (\eta_p) = \eta'_p$, and
\item   $ {\HH}_1(f) \circ \eta= \eta'$.
\end{itemize}
This last condition means that $\eta'$ is defined by ${\HH}_1(f) \circ \eta_1=  \eta'_1$ and $\eta'_2 = r^{-1} \eta_2$. Also, we have denoted $\Z_{(p)}^{\times,+}=\Q^\times_{>0}\cap\Z_{(p)}^\times$.

\begin{rem}
  Note that we are allowing  the similitude factor in the level structure to be
  in $\A^{\infty,p} \otimes_{\Q} F(1)$, but we only allow
  quasi-isogenies with similitude factor in $\A^{\infty,p}(1)$.
\end{rem}

We denote by  $\ZZ_{(p)}^{\times,+} \backslash (\A^{\infty,p}\otimes
F)^\times/ \nu(K^p)(1)$ the set $\ZZ_{(p)}^{\times,+} \backslash
(\A^{\infty,p}\otimes F)^\times/ \nu(K^p)$ equipped with the action
of $\mathrm{Gal}(\overline{\Q}/\Q)$ through the cyclotomic character
$\mathrm{Gal}(\overline{\Q}/\Q) \rightarrow \prod_{\ell \neq p}
\ZZ_\ell^\times\into(\A^{\infty,p})^\times$. This action is unramified at
$p$. It follows easily that $\ZZ_{(p)}^{\times,+} \backslash
(\A^{\infty,p}\otimes F)^\times/ \nu(K^p)(1)$  is represented by an
infinite disjoint union of finite \'etale schemes over
$\Spec \ZZ_{(p)}$. \begin{rem} The group $\ZZ_{(p)}^{\times,+}$ acts freely on $(\A^{\infty,p}\otimes F)^\times/ \nu(K^p)$. 
\end{rem}
\begin{rem} When $\nu(K^p) =  (\ocal_F \otimes_{\Z} \prod_{\ell \neq p} \ZZ_\ell)^{\times}$, then the above Galois action is trivial and $\ZZ_{(p)}^{\times,+} \backslash (\A^{\infty,p}\otimes F)^\times/ \nu(K^p)(1)$ is simply an infinite disjoint union of copies of $\Spec \ZZ_{(p)}$.
\end{rem}

There is a structural map  $\Pi_0: \mathbb{Y}_K \rightarrow
\ZZ_{(p)}^{\times,+} \backslash (\A^{\infty,p}\otimes F)^\times/
\nu(K^p)(1)$ which associates to  an object $(A, \iota, \lambda,
\eta,\eta_p)$ of $\mathbb{Y}_K$ the class of $\eta_2(1)$ (where we are
identifying ${\HH}_1( \mathbb{G}_m\vert_{\overline{s}},
\A^{\infty,p})$ with $\A^{\infty,p}(1)$). 

As we mentioned at the beginning of~\S\ref{sect-overc}, the
complex points of our integral models are not precisely
the double coset spaces considered in~\S\ref{sect-overc}, because our moduli problem only allows
polarizations of degree prime to~$p$. However, the difference amounts
to throwing away some geometrically connected components, as the
following result explains.

\begin{prop}\label{prop: existence of integral model for Y} The groupoid $\mathbb{Y}_K$ is  representable by a
  quasi-projective scheme $\Pi_0: Y_K \rightarrow
  \ZZ_{(p)}^{\times,+} \backslash (\A^{\infty,p}\otimes F)^\times/
  \nu(K^p)(1)$. The morphism $\Pi_0$ has geometrically connected  fibres. Let $c \in \ZZ_{(p)}^{\times,+} \backslash (\A^{\infty,p}\otimes F)^\times/
  \nu(K^p)$ and let $c: \Spec \C \rightarrow \ZZ_{(p)}^{\times,+} \backslash (\A^{\infty,p}\otimes F)^\times/
  \nu(K^p)(1)$ be the associated morphism  \emph{(}for the usual choice of primitive roots of unity in $\C$\emph{)}. Let $Y_{K,c}$ be the fibre of $Y_K$ over $c$. Then there is an isomorphism of analytic spaces 
 $ (Y_{K,c})^{an} = \Gamma(c, K) \backslash \mathcal{H}^{\mathrm{Hom}(F, \R)}$.

\end{prop}
\begin{proof}
  This follows from the usual description of integral models of PEL type
Shimura varieties; in the case of hyperspecial level this goes back to
Kottwitz~\cite{MR1124982}, but for convenience we
follow the notation of~\cite{MR3186092}. To this end, we recall the
description of these integral models for the usual Shimura varieties
for~$G$. We let~$\tK=\tK^p\tK_p$ denote a compact open subgroup
of~$G(\A^{\infty})$, where~$\tK^p$ is a compact open subgroup
of~$G(\A^{\infty,p}$), and~$\tK_p$ is of one of the parahoric subgroups
considered above.

Then we let~$\bY^{G,\Kot}$ be the groupoid over $\AFF/\Z_{(p)}$ whose fibre over $S \in \Ob(
\AFF/\Z_{(p)})$ is the category with objects  $(A, \iota, \lambda,
\etat, \eta_p)$, where~$(A,\iota,\lambda,\eta_p)$ is as in the
definition of~$\bY_K$ above, but now~$\etat$ is given by a~$\tK^p$-orbit
of pairs of isomorphisms~$(\etat_1,\etat_2)$, consisting of:

\begin{enumerate}
\item An $\cO_F\otimes\A^{\infty,p}$-linear isomorphism of $\Pi_1(S, \overline{s})$-modules $\etat_1:   V \otimes_{\ZZ} \A^{\infty,p}  \simeq {\HH}_1( A\vert_{\overline{s}}, \A^{\infty,p})$.
\item An $\A^{\infty,p}$-linear isomorphism of $\Pi_1(S, \overline{s})$-modules $\etat_2:  \A^{\infty,p} \simeq \A^{\infty,p}  \otimes_{\ZZ} {\HH}_1( \mathbb{G}_m\vert_{\overline{s}}, \A^{\infty,p})$.
\end{enumerate}

We moreover impose that the following diagram is commutative:
\numequation\label{eqn: pairings eta commute two}
\xymatrix{ V \otimes _{\ZZ} \A^{\infty,p}   \times V \otimes_{\ZZ} \A^{\infty,p}   \ar[rr]^{\etat_1 \times \etat_1} \ar[d]^{<,>}& & {\HH}_1( A\vert_{\overline{s}}, \A^{\infty,p}) \times {\HH}_1( A\vert_{\overline{s}}, \A^{\infty,p})  \ar[d]^{<,>_{\lambda}}\\ 
 \A^{\infty,p}   \ar[rr]^{\etat_2} & &{\HH}_1( \mathbb{G}_m\vert_{\overline{s}}, \A^{\infty,p})}
\end{equation}

A map between  quintuples $(A, \iota, \lambda, \etat,\eta_p)$ and $(A',
\iota', \lambda', \etat',\eta'_p)$ is an $\cO_F$-linear prime to $p$
quasi-isogeny $f: A \rightarrow A'$  such that \begin{itemize}
\item $f^{*} \lambda = r \lambda'$ for a locally constant function
  $r: S \rightarrow \Z_{(p)}^{\times,+}$, 
\item $f (\eta_p) = \eta'_p$, and
\item   $ {\HH}_1(f) \circ \etat= \etat'$.
\end{itemize}

It follows immediately from the definition that there is a natural
isomorphism\[\mathbb{Y}_K\cong\coprod_{g\in G(\A^{\infty,p})\backslash
    G_1(\A^{\infty,p})/K^p}\bY^{G,\Kot}_{gK^pg^{-1}\cap
    G(\A^{\infty,p})}, \] given by the
maps \[g:\bY^{G,\Kot}_{gK^pg^{-1}\cap G(\A^{\infty,p})}\to\bY_K \]
which are defined
by
\[(A,\iota,\lambda,(\etat_1,\etat_2),\eta_p)\mapsto
  (A,\iota,\lambda,(\etat_1,\etat_2\otimes_\Z \cO_F)g,\eta_p).\]
(Indeed, one easily checks that this already gives a bijection of
tuples before passing to isogeny classes, and that this bijection is
compatible with isogenies.) The result now follows from~\cite[\S 5,
\S 8]{MR1124982}.\end{proof}

We now define an action of $(\ocal_F)_{(p)}^{\times, +}$ (totally
positive elements in $F^\times$ which are prime to $p$) on
$\mathbb{Y}_K$, by scaling the polarization~$\lambda$. Since this
scales the $\lambda$-Weil pairing $\langle,\rangle_{1,\lambda}$, we
see from~(\ref{eqn: pairings eta commute one}) that it also 
scales~$\eta_2$. Explicitly, $x \in (\ocal_F)_{(p)}^\times$ sends $(A,
\iota, \lambda, (\eta_1,\eta_2),\eta_p)$ to  $(A, \iota, x\lambda,
(\eta_1,x\eta_2),\eta_p)$. By definition, the subgroup
$\Z_{(p)}^{\times,+}$ acts trivially on $\mathbb{Y}_K$. 

The group $(\ocal_F)_{(p)}^{\times, +}$ acts on the set of connected
components $\Pi_0({Y}_K)$. Since the cyclotomic character surjects
onto~$\prod_{l\ne p}\Z_l^\times$, the stabilizer of each connected component
is \[\ocal_{F}^{\times,+}( \Pi_0): = \left((\ocal_F)_{(p)}^{\times, +}\cap
\Z_{(p)}^{\times,+}\nu(K^p) \prod_{\ell \neq p} \ZZ_\ell^{\times}\right)/\Z_{(p)}^{\times,+},\]
which we can and do naturally identify with \[ \ocal_F^{\times, +} \cap
\nu(K^p) \prod_{\ell \neq p} \ZZ_\ell^{\times}.\]

\begin{rem}If $\nu(K^p) =  (\ocal_F \otimes_{\Z} \prod_{\ell \neq p} \ZZ_\ell)^{\times}$, then  $\ocal_{F}^{\times,+}( \Pi_0) =  \ocal_F^{\times, +}$. 
\end{rem}

 The subgroup $\ocal_{F}^{\times,+}( \nu(K^p)): =  \ocal_F^{\times, +}
 \cap \nu(K^p)$ acts trivially on each connected component of
 $\Pi_0({Y}_K)$. The quotient stack of connected components is \[[ ((\ocal_F)_{(p)}^{\times, +}/\Z_{(p)}^{\times,+}) \backslash \big( \ZZ_{(p)}^{\times,+} \backslash (\A^{\infty,p}\otimes F)^\times/ \nu(K^p)(1) \big)].\]  It admits a coarse moduli space $ (\ocal_F)_{(p)}^{\times,+} \backslash (\A^{\infty,p}\otimes F)^\times/ \nu(K^p)(1)$ which is a finite \'etale covering of $\Spec \ZZ_{(p)}$.

We  now take the quotient stack \[\mathbb{Y}_K^{G_1}:=
  [{Y}_K/((\ocal_F)_{(p)}^{\times, +}/\Z_{(p)}^{\times,+})].\] This is the ``Shimura stack''
associated to $G_1$ and the level $K$. 

Let us define 
$$\ocal_F^{\times,+}(K^p)  = \{ x^2\mid x \in  \ocal_F^{\times} \cap K^p\},$$
where~$\OL^{\times}_F$ is thought of inside~$G_1(\A^{(p)}_f)$
as a subgroup of the scalar matrices. The multiplier of the scalar matrix given by~$x$
is~$x^2$,
and hence the multiplier of~$\ocal_F^{\times,+}(K^p)$ lands inside~$\nu(K^p)$,
and hence~$\ocal_F^{\times,+}(K^p)$ is a finite index subgroup of 
$\OL^{\times,+}_{F}(\nu(K^p))$ and of~$\ocal_{F}^{\times,+}( \Pi_0)$.

\begin{lem}\label{lem: x2 acts trivially} The restriction of the action of $(\ocal_F)_{(p)}^{\times, +}$ on $Y_{K}$ to  $\ocal_F^{\times,+}(K^p)$ is trivial. More precisely, there is a canonical natural transformation going from the action of   $\ocal_F^{\times,+}(K^p)$ on $\mathbb{Y}_{K}$ to the trivial action of $\ocal_F^{\times,+}(K^p)$ on $\mathbb{Y}_{K}$.
\end{lem}
\begin{proof} Let $x^2 \in \ocal_F^{\times,+}(K^p)$ for a unique $x
  \in \ocal_F^{\times,+} \cap K^p$. The action of $x^2$ sends $(A,
  \iota, \lambda, \eta,\eta_p)$ to $(A, \iota, x^2 \lambda, \eta,\eta_p)$ (note that
  since $x\in K^p$, and $\eta$ is by definition a $K^p$-orbit, the
  action of~$x^2$ on~$\eta$ is trivial). On the other hand
  multiplication by $x^{-1}: A \rightarrow A$ provides a map $(A,
  \iota, \lambda, \eta,\eta_p) \rightarrow (A, \iota, x^2\lambda, \eta,\eta_p)$ in
  the groupoid $\mathbb{Y}_K$. This provides the natural
  transformation from the action of $x^2$ obtained from the action of $(\ocal_F)_{(p)}^{\times, +}$ to the trivial action. 
\end{proof}

\begin{lem}\label{lem-freedelta} For any geometric point $x \in Y_K$, the stabilizer of $x$ for the action of $(\ocal_{F})_{(p)}^{\times, +}$ is $\ocal_F^{\times,+}(K^p)$.
\end{lem}

\begin{proof}By Lemma~\ref{lem: x2 acts trivially}, 
  $\ocal_F^{\times,+}(K^p)$ is contained in the stabilizer of any~$x
  = (A, \iota, \lambda, \eta)$. Let $\epsilon \in
  (\ocal_{F})_{(p)}^{\times, +}$ and assume that there is a morphism \[
  f: (A, \iota, \lambda, (\eta_1,\eta_2),\eta_p)  \rightarrow (A, \iota,
  \epsilon\lambda,(\eta_1,\epsilon\eta_2), \eta_p)\] in the groupoid $\mathbb{Y}_K$. We need
  to show that $f \in \ocal_F^{\times} \cap K^p$. Since~$f$
  respects~$\eta_1$, it follows from~\cite[Lem.\ 1.3.5.2]{MR3186092}
  that $f$ is 
  an automorphism of $A$ (and not just a quasi-isogeny).

The polarization
$\lambda$ induces an involution $ x \mapsto \bar{x}$ on
$F(f)$, and we consider the automorphism
$\alpha = f \kern-0.1em{\bar{f}^{-1}}$ of $A$. It stabilizes the polarization:
$\alpha^{*} \lambda = \lambda \alpha \bar{\alpha} = \lambda$. It
also stabilizes the level structure: $\bar{f}$ acts like the adjoint
of $f$ on ${\HH}_1( A, \A^{\infty,p})$.  Since $K^p$ is neat,
this implies that $\alpha=1$; indeed, all the eigenvalues of~$\alpha$
are roots of unity, because they are algebraic numbers all of whose
conjugates have absolute value~$1$. It follows that $f=\bar{f}$, and $f^2 = f\kern-0.1em{\bar{f}}= \epsilon$.
Since~$f$ is an automorphism, it follows that~$\epsilon \in \OL^{\times}_F$. Hence
it suffices to show that~$f \in F$, since we then have~$\epsilon \in \ocal_F^{\times,+}(K^p)$.

Assume first that $A$ is simple, so that $\mathrm{End}(A)_{\Q}$ is a
division algebra and $F(f) \subset \mathrm{End}(A)_{\Q}$ is a
commutative field on which the Rosati involution $ x \mapsto \bar{x}$ is  complex conjugation. Since $f=\bar{f}$ and $f^2=\epsilon$, $F(f)$ is a
totally real extension of $F$ of degree at most~$2$. If $F(f) = F$,
we are done. Otherwise  $F(f)$ is a quadratic  extension of $F$. The level
structure $\eta$ provides  a $K^p$-orbit of isomorphisms ${\HH}_1(
A, \A^{\infty,p}) \simeq V \otimes \A^{\infty,p}$, and the element $f$ acts via some conjugate of    $$\begin{pmatrix} 0 & \epsilon & 0 & 0 \\
      1 & 0  & 0 & 0\\
      0 & 0 & 0 & 1\\
      0 & 0 & \epsilon & 0
   \end{pmatrix}$$ and has eigenvalues in $\bar{F}$: $\{
   \sqrt{\epsilon}, -\sqrt{\epsilon}\}$ with multiplicity two. By
   neatness, no conjugate of this matrix is in $K^p$, a
   contradiction. 

We now assume that $A$ is not simple. It is easy to see (using the
$\cO_F$-action) that the only possibility is that  $A$ is isogenous to
$A_1 \times A_2$ where $A_1$ and $A_2$ are two abelian schemes of
dimension $[F:\Q]$ with~$F \subset \End(A_i)_{\Q}$. If $A_1$ and $A_2$ are not isogenous, then
$\mathrm{End}(A)_{\Q} = \mathrm{End}(A_1)_{\Q} \times
\mathrm{End}(A_2)_{\Q}$. Moreover, $F(f)$ is a commutative subalgebra
of $\mathrm{End}(A_1)_{\Q} \times \mathrm{End}(A_2)_{\Q}$ and is
therefore included in a product of fields $F_1 \times F_2$ where $F_i$
is either $F$ or a CM extension of $F$. Since $f=\overline{f}$, we see that $f  = (f_1, f_2) \in F \times F$ and that
$f^2 = (f_1^2, f_2^2) = \epsilon$. So either $f_1=f_2$, and we are done, or
$f_1=-f_2$; but this second case is again prohibited by neatness. 

Lastly, we assume that $A$ is isogenous to $A_1^2$. Then $\mathrm{End}(A)_{\Q} \simeq M_2(\mathrm{End}(A_1)_{\Q})$ and $F(f)$ is a commutative subalgebra, therefore included in $M_2(E)$ where $E$ is either $F$ or a CM extension of $F$. Writing $f =    \begin{pmatrix} a & b \\
      c & d \\
   \end{pmatrix} \in \mathrm{GL}_2(E)$, we have $\epsilon=f^2 =    \begin{pmatrix} a^2 +bc & b(a+d) \\
         c(a+d) & bc+d^2 \\
      \end{pmatrix}$. If  $a+d = 0$, the matrix of $f$ has eigenvalues
      $\{ \sqrt{\epsilon}, - \sqrt{\epsilon}\}$ and this is again
      impossible by neatness. We deduce that $a+d \neq 0$, so that
      $b=c= 0$ and $a = d = \sqrt{\epsilon}$ or $a=d = -
      \sqrt{\epsilon}$. Since $\overline{f} = f$ and the Rosati
      involution induces the complex conjugation on $E$, we deduce
      that $\sqrt{\epsilon} \in F$ and that $f \in F$, as required.
\end{proof}

We write \numequation\label{eqn:delta}\Delta = (\ocal_{F})_{(p)}^{\times, +}/
\ocal_F^{\times,+}(K^p), \end{equation}  \numequation\label{eqn:
delta pi zero}\Delta(\Pi_0) = \ocal_{F}^{\times,+}(
\Pi_0)/ \ocal_F^{\times,+}(K^p),\end{equation}
 \numequation\label{eqn:
delta K p}\Delta(K^p) = \ocal_{F}^{\times,+}(
\nu(K^p))/ \ocal_F^{\times,+}(K^p).\end{equation} These last two  groups are finite groups.
Let us set $Y_K^{G_1} = \Delta \backslash Y_K$. This last quotient
exists as a scheme. Indeed, $\Delta$ permutes the connected components of $Y_K$ and the
stabilizer of any connected component is a finite group
$\Delta(\Pi_0)$, while the stabilizer of any geometrically connected
component if $\Delta(K^p)$. Moreover, the action of $\Delta$ can be
lifted to an 
action on an ample line bundle on $Y_K$ (for instance the tensor product of the line bundles $\det(\Omega^1_{(A/C)/Y_K})$ where $C$ runs over all subgroups $C=\prod_{v\mid p}C_v$ where for each $v\mid p$, $C_v$ is either $1$ or whichever of $H_v,L_v$ exist as part of the level structure, see \cite[\S6]{lan_2016}). The group $\Delta(\Pi_0)$
acts without fixed points by Lemma~\ref{lem-freedelta}. The following proposition then follows
immediately from Proposition~\ref{prop: existence of integral model
  for Y} and Lemma~\ref{lem: Delta Pilloni lemma over C}.\begin{prop}\label{prop: integral model for G1} There is a canonical
  map  $\mathbb{Y}_K^{G_1} \rightarrow Y^{G_1} _K$, and $Y^{G_1}_K$ is
  the coarse moduli of $\mathbb{Y}_K^{G_1}$. There is a
  quasi-projective morphism  $\Pi_0: Y^{G_1}_K \rightarrow
  (\ocal_F)_{(p)}^{\times,+} \backslash (\A^{\infty,p}\otimes
  F)^\times/ \nu(K^p)(1)$ with geometrically connected
  fibres. Moreover, the map $Y_K \rightarrow Y^{G_1}_K$ is \'etale and
  surjective. 

Let $c \in \ZZ_{(p)}^{\times,+} \backslash (\A^{\infty,p}\otimes F)^\times/
  \nu(K^p)$ and let 
  $$c: \Spec \C \rightarrow \ZZ_{(p)}^{\times,+} \backslash (\A^{\infty,p}\otimes F)^\times/
  \nu(K^p)(1)$$
   be the associated morphism  \emph{(}for the usual choice of
  primitive roots of unity in $\C$\emph{)}. Let $Y_{K,c}$ be the fibre of
  $Y_K$ over $c$ and let ${Y}_{K,c}^{G_1}$  be the fibre of $Y_{K}^{G_1}$ over $c$. Then there is a commutative diagram of  analytic spaces where the horizontal maps are isomorphisms and the vertical maps are finite \'etale with groups $\Delta(K^p) = \Delta(K)$: 
\begin{eqnarray*}
\xymatrix{ (Y_{K,c})^{an} \ar[r] \ar[d]  & \Gamma(c, K) \backslash \mathcal{H}^{\mathrm{Hom}(F, \R)}\ar[d] \\
(Y^{G_1}_{K,c})^{an} \ar[r] & \Gamma_1(c, K) \backslash \mathcal{H}^{\mathrm{Hom}(F, \R)} }
\end{eqnarray*}

\end{prop}

\subsection{Local models}\label{subsubsec: local models}We now
recall some basic results about local models for~$\GSp_4$; the cases
that we need essentially go back to~\cite{MR1227472}. Continue to
let~$K$ be a reasonable level structure. For each place~$v|p$, we let
$M^{\loc}_{K_v}$ be the moduli space over~$\cO_{F_v}$ of chains of
lattices corresponding to~$K_v$; so for example $M^{\loc}_{\Par(v)}$
is the moduli space of totally isotropic direct factors of
$V_1\otimes_{\cO_F}\cO_{F_v}$ of rank~$2$. We write
$M^{\loc}_{K_p}:=\times_{v|p}M^{\loc}_{K_v}$. Then by the results
of~\cite[\S6]{MR1393439}, each geometric point of the special fibre
of~$Y_K^G$ has an \'etale neighbourhood which is isomorphic to an
\'etale neighbourhood of a geometric point in  $M^{\loc}_{K_p}$. (The description of the
local model in~\cite[\S6]{MR1393439} is in terms of chains of
$\cO_F\otimes\Zp$-lattices, but this description can be immediately
rewritten in terms of products over the places $v|p$ of chains of
$\cO_{F_v}$-lattices.)

\begin{prop}\label{prop: local models results for G} The scheme
  $\mathbb{Y}_K$ is flat over $\Spec \ZZ_{(p)}$, normal,
  and a  local complete intersection \emph{(}so in particular
  Cohen--Macaulay\emph{)} of pure relative dimension~$3[F:\Q]$. If~$K_v=\GSp_4(\cO_{F_v})$ for
  all $v|p$, then it is smooth, while in general it is smooth away
  from codimension~$2$. \end{prop}
\begin{proof}Note that normality follows from being  smooth
 away from codimension~$2$ and Cohen--Macaulay. The properties of being flat and a local complete intersection
  over~$\Spec\Z_{(p)}$, and of being smooth, or smooth away from
  codimension~$2$, can all be checked
  \'etale locally
  (\cite[\href{http://stacks.math.columbia.edu/tag/03E7}{Tag
    03E7},\href{http://stacks.math.columbia.edu/tag/04R3}{Tag
    04R3},\href{http://stacks.math.columbia.edu/tag/06C3}{Tag
    06C3}]{stacks-project}). Furthermore, these properties are all
  preserved by taking products.  It therefore suffices to show that
  they hold for the local models~$M^{\loc}_{K_v}$. This has already been carried out in the
  literature: the case that $K_v=\GSp_4(\cO_{F_v})$ is trivial, and
  the cases that~$K_v=\Kli(v),$ $\Si(v)$ or $\Iw(v)$ are covered
  in~\cite[\S2]{MR2290605}. In the case $K_v=\Par(v)$ see~\cite[Prop.\
  2.5, Thm.\ 2.11]{MR2811273}.\end{proof}

\begin{cor}\label{cor: local models for G1} The scheme $Y^{G_1}_K$  is normal, flat over $\Spec \ZZ_{(p)}$, and
  a local complete intersection. \end{cor} 
\begin{proof}Since $Y_K\to Y^{G_1}_K$ is an \'etale surjection by
  Proposition~\ref{prop: integral model for G1}, this is immediate
  from Proposition~\ref{prop: local models results for G}.  
\end{proof}

\subsection{Compactifications}\label{sec: compactifications}

In this section, we state results on the existence of toroidal
compactifications. Toroidal compactifications depend on some
combinatorial data which we first explain. We will follow closely the
presentation of~\cite{PINK} and \cite{hltt}, see in particular
~\cite[\S5.2]{hltt} (that this presentation is equivalent to Lan's
presentation is explained in ~\cite[App.\ B]{hltt}).     

In this section, we write~$V_F$ for~$V\otimes_{\cO_F} F$. Let $\mathfrak{C}$ be the set of  totally isotropic $F$-subspaces $W
\subset V_F$. For all $W\in \mathfrak{C}$, consider  the  $F \otimes
\R$-module of  $\Q$-bilinear
forms \[ \phi: V_F/W^\bot \times V_F/W^\bot
  \rightarrow  \R\] which satisfy $\phi( \lambda x, y) =
\phi(x, \lambda y)$ for all $\lambda \in F$, $x,y \in V_F/W^\bot$.
Let
$C(V_F/W^\bot)$ be the cone inside this   $\R$-vector space given by those
forms which are positive semidefinite and    whose radical is defined over $F$.  Let $\mathcal{C}$ be the conical complex  which is the quotient of $ \coprod_{W \in \mathfrak{C}} C(V_F/W^\bot) $ by the equivalence relation induced by the inclusions $C(V_F /W^\bot) \subset C(V_F/Z^\bot)$ for $W \subset Z$.

 A \emph{non-degenerate rational polyhedral cone} of  ${\mathcal{C} \times G_1(\A^{\infty})}$ is a subset contained in $\mathcal{C}(V_F/W^\bot) \times \{\gamma\}$  for some $(W, \gamma)$ which is of the form $ \sum_{i= 1}^k\mathbb{R}_{>0}s_i$ for elements $s_i:  V_F/W^\bot \times V_F/W^\bot \rightarrow \Q$.

 A \emph{rational polyhedral cone decomposition} $\Sigma$ of ${\mathcal{C} \times G_1(\A^{\infty})}$ is a partition ${\mathcal{C} \times G_1(\A^{\infty})} = \coprod_{\sigma \in \Sigma} \sigma$ by    non-degenerate rational polyhedral cones $\sigma$  such that the closure of each cone is a union of cones.

 Let $W \in \mathfrak{C}$. We let  $P_W$ be the parabolic subgroup of
 $G_1$ which is the stabilizer of $W$. Let us denote by $M_{W,l}$ the
 group of $F$-linear automorphisms of $V_F/W^\bot$.   We also denote by
 $M_{W,h}$ the group of symplectic similitudes of $W^\bot/W$ (so that
 this group is isomorphic to $\mathrm{Res}_{{F}/\Q} \mathrm{GSp}_{4- 2
   \dim W}$, and in particular is non-trivial even when $\dim W=2$).  The group $M_{W} = M_{W,h} \times M_{W,l}$ is the Levi quotient of $P_W$. We have a surjective map $P_W \rightarrow M_{W,l}$, and we denote by $P_{W, h}$ its kernel. There is a surjective map $P_{W,h} \rightarrow M_{W,h}$. 

The group $G_1( \QQ)^+$ acts on $\mathfrak{C}$ and also on $\mathcal{C}$. Let $W \in \mathfrak{C}$, let $\gamma \in G_1(\Q)^+  \cap P_W$ and $\phi \in C(V_F/W^\bot)$. Let $\gamma_l$ be the projection of $\gamma$ in $M_{W,l}$. Then we set $\gamma \phi (x, y) = \nu(\gamma) \phi(\gamma_l .x , \gamma_l. y)$.

The  set ${\mathcal{C} \times G_1(\A^{\infty})}$ carries a diagonal left action of $G_1( \QQ)$ and  left and  right actions of $G_1(\A^{\infty})$ (by left and right multiplication on the second factor).    For any compact open  subgroup $K \subset  G_1(\A^{\infty})$, a rational polyhedral cone decomposition $\Sigma$ is $K$-equivariant if for all $h \in G_1( \QQ), k \in K$ and $\sigma \in \Sigma$, $h. \sigma . k \in \Sigma$.

 For any compact open subgroup  $K \subset  G_1(\A^{\infty})$ we say that a rational polyhedral cone decomposition $\Sigma$ of ${\mathcal{C} \times G_1(\A^{\infty})}$   is $K$-admissible if: 
 
 \begin{enumerate} \label{assumption}
 \item The decomposition is $K$-equivariant.
 \item For all $ \sigma \subset  C(V_F/W^\bot)  \times \{ \gamma\}$, and all $p \in P_{W,h}( \mathbb{A}^{\infty})$, we have $ p. \sigma  \in \Sigma$.\item  For all cones $\sigma $, let $W \in \mathfrak{C}$ be  such that $\sigma \subset \mathcal{C}( V_F/W^\bot)$ is in the interior of $\mathcal{C}( V_F/W^\bot)$. Then if there are $p \in  P_{W,h}(\mathbb{A}^{\infty})$, $u \in K$ and $h \in G_1(\Q)$  satisfying  $\sigma \cap  h p\sigma u \neq \emptyset$,  then in fact $h \in P_{W,h}(\mathbb{A}^{\infty})$.
 
\item $G_1(\Q) \backslash \Sigma / K$ is finite.
\end{enumerate}

 There exist $K$-admissible rational polyhedral cone decompositions.   Any two $K$-admissible rational polyhedral  cone decompositions can be refined by a third one.

 If~$L_W \subset \mathrm{Hom}_\Q ( \mathrm{Sym}^2_{F} V_F/W^\bot,
 \mathbb{Q}) $ is a lattice, then a cone 
 $$\sigma \subset \mathrm{Hom}_\Q (
 \mathrm{Sym}^2_{\ocal_F} V_F/W^\bot, \mathbb{Q}) $$
 is said to be smooth with respect to $L_{W}$ if the $s_i$ can be taken to be part of a  basis of $ L_{W}$. Assume that for all $(W, \gamma) \in {\mathfrak{C} \times G_1(\A^{\infty}) }$ we have lattices \[L_{W, \gamma} \subset  \mathrm{Hom}_\Q ( \mathrm{Sym}^2_{F} V_F/W^\bot,\mathbb{Q}).\]     We say that a  rational polyhedral  cone decomposition $\Sigma$ is smooth with respect to  these lattices if each cone $\sigma \in \Sigma $ is smooth.

 We now assume that $K = K^pK_p$  is a reasonable compact open subgroup. We choose a lattice $V' \subset V_F$ with the property that $K^p$ stabilizes $V' \otimes_{\Z} \mathbb{A}^{\infty, p}$ and that $V' \otimes \ocal_{F_v} = V \otimes \ocal_{F_v}$ for all places $v |p$ such that $K_v \neq \Par(v)$ and $ V' \otimes \ocal_{F_v} = V_3 $ for all places $v$ such that $K_v =  \Par(v)$. 
 
Then  $(\ocal_F, V', \langle. \rangle)$ defines an integral PEL datum and $K \subset G'_1( \hat{\mathbb{Z} })$ where $G'_1$ is the group scheme over $\Spec \Z$ of symplectic similitudes of $V'$. 
 
 The theory of toroidal compactification associates a lattice $L_{W, K,
   \gamma} \subset  C(V_F/W^\bot)$ to this integral PEL datum, compact open $K$, $W \in
 \mathfrak{C}$ and $\gamma \in  G_1(\A^{\infty})$ 
 (see~\cite[\S5.3]{MR3186092} and~\cite[\S3]{lan_2016}). The  $K$-admissible rational polyhedral  cone decompositions which satisfy the following extra properties form a cofinal subset of the set of all $K$-admissible rational polyhedral  cone decompositions:
\begin{enumerate}
\item  The decomposition is projective (in the sense of~\cite{AMRT}).
\item  The decomposition is smooth with respect to the lattices $L_{W,K, \gamma}$. 
\end{enumerate}

  In the rest of the paper, we will  consider  $K$-admissible rational polyhedral cone decompositions which satisfy these extra properties unless explicitly stated.

\begin{thm}\label{thm: main KWL compactification thm}  \leavevmode
\begin{enumerate}
 \item Let $\Sigma$ be a $K$-admissible  polyhedral cone decomposition which is projective.  There is a toroidal compactification $X_{K,\Sigma}$ of $Y_K$.   It has a stratification indexed by $ (G(\Q)^+ \cap K_p) \backslash \Sigma / K^p =  G(\Q)^+ \backslash \Sigma / K $. The boundary is  the reduced complement of $Y_{K}$ in $X_{K, \Sigma}$.   This is a  relative Cartier divisor denoted by $D_{K, \Sigma}$.  
 
\item\label{item: semi-abelian} The universal abelian scheme $A \rightarrow Y_{K}$ extends to a semi-abelian scheme $A \rightarrow X_{K,\Sigma}$. 
  
\item\label{item: pi sigma sigma prime} If $\Sigma'$ is a refinement of~$\Sigma$, then there are
  projective maps $\pi_{\Sigma', \Sigma}: X_{K, \Sigma'}\rightarrow
  X_{K, \Sigma}$, and $(\mathrm{R}\pi_{\Sigma',\Sigma})_{*}
  \ocal_{X_{K, \Sigma'}} = \ocal_{X_{K, \Sigma}}$. Let
  $\mathcal{I}_{X_{K, \Sigma}}$ and $\mathcal{I}_{X_{K, \Sigma'}}$  be
  the invertible sheaves of the boundary in $X_{K, \Sigma}$ and $X_{K,
    \Sigma'}$. Then $\pi_{\Sigma', \Sigma}^{*} \mathcal{I}_{X_{K,
      \Sigma}} = \mathcal{I}_{X_{K, \Sigma'}}$ and $(\mathrm{R}\pi_{\Sigma',\Sigma})_{*} \mathcal{I}_{X_{K, \Sigma'}}  = \mathcal{I}_{X_{K,
      \Sigma}}$.

\item \label{item: lci X Sigma} Suppose that~$K$ is reasonable \emph{(}in the sense of Definition~\ref{defn: reasonable}\emph{)}. Then the toroidal compactification $X_{K,\Sigma}$ is flat over
  $\Spec \ZZ_{(p)}$, normal, and Cohen-Macaulay.  If $\Sigma$ is smooth, then $X_{K,\Sigma}\to\Spec \ZZ_{(p)}$ is further a local complete intersection.  Finally if $K_v= \mathrm{GSp}_4(\ocal_{F_v})$ for all $v| p$ and $\Sigma$ is smooth then $X_{K,\Sigma}\to\Spec \ZZ_{(p)}$ is smooth.\end{enumerate}
 \end{thm}
 \begin{proof}This follows from~\cite[Thm.\ 6.1]{Lan2016}. We simply need to specify the choices we made to construct the toroidal compactification by normalization (see \cite[\S2]{lan_2016}).  In the first  case that $K_p = \mathrm{G}_1( \Z_p)$ (the nice case: no level at $p$, prime to $p$ polarization), the compactification is constructed  in \cite{MR3186092}. In the  second case that $K_p = \prod_{v |p} K_v$ where $K_v \in \{ \mathrm{GSp}_4(\ocal_{F_v}), \Par(v) \}$, the compactification can be constructed as a closed subscheme of some  toroidal compactification of a Siegel modular variety with a prime to $p$ polarization (Zarhin's trick) (and possibly performing again a blow up or a blow down at the boundary as explained in \cite{Lan2016}).  In the general case where we have a parahoric level structure, we consider all possible degeneration maps $Y_{K} \rightarrow \prod_{K'_p} Y_{K^p K'_p}$ where $K_p  \rightarrow  K'_p$  and $K'_p = \prod_{v |p} K'_v$ with  $K'_v \in \{ \mathrm{GSp}_4(\ocal_{F_v}), \Par (v)\}$ and obtain the toroidal compactification as a closed subscheme  of the product of the toroidal compactifications of the $Y_{K^p K'_p}$ (and possibly performing again a blow up or a blow down at the boundary as explained in \cite{Lan2016}). 
 
  Now, everything apart from~(\ref{item: lci X Sigma}) is
   immediate, while~(\ref{item: lci X Sigma}) follows
   from Proposition~\ref{prop: toroidal compactifications of G1}
   together with the explicit description of the formal completions
   along boundary strata given in~\cite[Thm.\ 6.1~(4)]{Lan2016}.\end{proof}

 We also need to consider the action of the group $\ocal_{F,
   (p)}^{\times, +}$. Recall that we defined a quotient~$\Delta$ of
 this group in~(\ref{eqn:
delta pi zero}).

 \begin{lem} The action of $\ocal_{F, (p)}^{\times, +}$  on $Y_K$ extends to $X_{K, \Sigma}$ and factors through ~$\Delta$.
 \end{lem}
 \begin{proof} It is possible to prove this directly by looking at the
   construction of the toroidal compactification and the boundary
   charts. We will instead give a simpler indirect argument. Since
   $X_{K, \Sigma}$ is normal, it follows that $X_{K, \Sigma}$ is the
   normalization of $Y_{K}$ in $X_{K, \Sigma} \times \Spec \C$. It is
   therefore sufficient to show that the action extends over $\C$. 

We can now use~\cite{AMRT}. Let $c \in G_1(\mathbb{A}_f)$. By
Proposition~\ref{prop: existence of integral model for Y}, the analytification of the 
 component $Y_{K,c} \subset Y_K \times \Spec \C$  corresponding to $c$     is $\Gamma(c, K)
\backslash \mathcal{H}^{ \mathrm{Hom} (F, \R)}$, and we need to show that the group $\Delta(K)$
(which is the subgroup of $\Delta$ acting trivially on the geometrically connected
components) acts on the compactification of $\Gamma(c, K ) \backslash
\mathcal{H}^{ \mathrm{Hom} (F, \R)}$. By the main results of~\cite{AMRT}, our choice of $\Sigma$ provides a partial compactification $\mathcal{H}^{ \mathrm{Hom} (F, \R)}_{\Sigma}$ which carries an action of $\Gamma(K,c)$.  The component of   $(X_{K, \Sigma} \times \Spec \C)^{an}$ corresponding to $c$ is isomorphic to $\Gamma (c, K ) \backslash \mathcal{H}^{ \mathrm{Hom} (F, \R)}_{\Sigma}$. This space still carries an action of $\Gamma_1(c, K ) / \Gamma (c, K) $, which is what we claimed.
 \end{proof}
 
 \begin{lem} The action of $\Delta$ on $X_{K, \Sigma}$ is free. 
 \end{lem} \begin{proof} Over $Y_{K}$, this is the content of Lemma~ \ref{lem-freedelta}.  We claim that the action of~ $\Delta$ is free on the set of non-trivial strata in $X_{K,\Sigma}$. This set is simply $G^+(\Q) \backslash \big( \Sigma \setminus \{0\} \times \mathbb{G}(\mathbb{A}^\infty) \big)/K$.
 Let  $c \in G_1(\mathbb{A}^\infty)$, $\Gamma(c, K) = G(\Q)^+ \cap c K c^{-1}$ and  $\Gamma_1(c, K) = G_1(\Q)^+ \cap c K c^{-1}$. Let $\Sigma_c$ be the restriction of $\Sigma$ to $\mathcal{C} \times \{c\}$. We need to show that the stabilizer of $\Gamma_1(c, K)$ acting on  $\Sigma_c \setminus \{0\}$ is included in  $\Gamma (c,K)$. This will imply that the group $\Delta (K)$ acts freely on $\Gamma(c,K) \backslash (\Sigma_c \setminus \{0\})$. 
 
 Let $W \in \mathfrak{C} \setminus \{0\}$. We denote by
 $\Gamma_W(c,K)$ and $\Gamma_{1,W}(c,K)$  the intersections of $P_W$
 with  $\Gamma(c,K)$ and $\Gamma_1(c,K)$ respectively.  Let $\sigma \subset  C(V_F/W^\bot) \times \{c\}$ in the interior. By our assumption on the cone decomposition, if an element $\gamma \in \Gamma_{1,W}(c,K)$ stabilizes~$\sigma$, then its linear part $\gamma_l$ is trivial. We need to see that $\nu(\gamma)$ is trivial. It is easy to see that we can find an element $\gamma' \in \Gamma_W(c,K)$  and $n \in \mathbb{Z}_{\geq 0}$ such that  $\nu(\gamma)^n \phi = \gamma'. \phi$ for all $\phi \in C(V_F/W^\bot) $ (it follows from the very definition of the action that the image of $\Gamma_W(c,K)$ in the space of automorphisms of $C(V_F/W^\bot)$ contains a finite index subgroup of $\ocal_{F}^{ \times, +}$). We deduce that $\gamma'$ stabilizes $\sigma$ and therefore $\gamma'_l=1$, so that  $\nu(\gamma)^n =1$ and $\nu(\gamma) =1$ since $\ocal_{F}^{\times, +}$ is torsion free.  \end{proof}

 We form the quotient of $X_{K, \Sigma}$ by the action of
 $\O_{F,(p)}^{\times,+}$. This quotient exists because, on a given connected
 component of $X_{K, \Sigma}$, this is the quotient by a finite group,
 and the component is projective because $\Sigma$ is a projective cone decomposition.
 We shall call such a quotient a toroidal compactification $X_{K,\Sigma}^{G_1}$ of $Y_{K}^{G_1}$. We summarize our findings in the following proposition:

 \begin{prop}\label{prop: toroidal compactifications of G1}  The space
   $X_{K,\Sigma}^{G_1}$ has a stratification indexed by $G_1(\Q)^+
   \backslash \Sigma/ K$. The map $X_{K,\Sigma} \rightarrow
   X_{K,\Sigma}^{G_1}$ is \'etale and surjective. If~$K$ is
   reasonable, then $X^{G_1}_{K,\Sigma}$ is a flat local complete
   intersection over~$\Spec \Z_{(p)}$, and is normal.
 \end{prop}

If not necessary, we drop the subscripts $K$ or $\Sigma$ and simply write $X$. We denote the boundary divisor by $D$.

\subsection{Functorialities}\label{sect-functorialities}We now briefly
discuss some functorial maps between Shimura varieties at different
levels, which we will make use of when we discuss Hecke operators in~\S\ref{subsec:
   coherent cohomology and Hecke operators}. All of the
 functorialities that we consider here extend to the toroidal
 compactification for suitable choices of cone decompositions, so we
 confine our discussions to the interior.
\subsubsection{ Change of level away from $p$} \label{subsect-functor-1}

Let $K = K^pK_p$ and $K' =  (K^p)' K_p$ be two compact open subgroups
of $G_1(\A^{\infty})$ such that $K \subset K'$. Then we have finite
\'etale maps $ Y_{K} \rightarrow Y_{K'}$ and $Y_{K}^{G_1} \rightarrow
Y_{K'}^{G_1}$, given by ``forgetting the level structure''; that is,
by replacing the $K^p$-orbit by the corresponding $(K^p)'$-orbit.
\subsubsection{Action of the group $G_1(\A^{\infty,p})$}\label{subsect-functor-2}

Let $g \in G_1(\A^{\infty,p})$. Then we can define an isomorphism 
$$ [g]: Y_K \rightarrow Y_{g^{-1} K g}$$
by sending an object $(A, \iota, \lambda, \eta, \eta_p)$ of $\mathbb{Y}_K$ to $(A, \iota, \lambda, \eta \circ g , \eta_p)$, which is immediately seen to be an object of $\mathbb{Y}_{ g^{-1} K g}$. 

We deduce isomorphisms $ [g]: Y^{G_1} _K \rightarrow Y^{G_1}_{g^{-1} K g}$.

\subsubsection{ Change of level at $p$:  Klingen type correspondences}\label{subsect-functor-3} We now fix $K^p$  and a place $w$ above $p$. We let  $K_p = \prod_{v |p } K_v \subset G_1(\ZZ_p) $ be a reasonable compact open such that  $K_w = \mathrm{GSp}_4( \ocal_{F_w})$.  We let $K'_p = \prod_{ v \neq w} K_v \times \Kli(w)$ be another reasonable level structure at $p$ and let $K''_p =  \prod_{ v \neq w} K_v \times \Par( w)$. Set $K = K^pK_p$, $K' = K^p K'_p$ and $K'' = K^p K''_p$. 

\begin{lem} There are natural proper surjective, generically finite \'etale forgetful maps $p_1: Y_{K'} \rightarrow Y_{K}$ and $p_1:Y_{K'}^{G_1} \rightarrow Y_{K}^{G_1}$. 
\end{lem}

\begin{proof} We simply forget the level structure $H_w$ at $w$.
\end{proof} 

We now choose once and for all an element $x_w \in F^{\times, +}$ which is a uniformizing element in ${F_w}$ and a unit in $F_v$ for all $v\neq w$ above $p$. This element is well defined up to multiplication by an element of $(\ocal_{F})_{(p)}^{\times, +}$. 

\begin{lem} There is a proper, surjective, generically finite \'etale map $p_2:Y_{K'} \rightarrow Y_{K''}$ depending on $x_w$ and sending $A$ to $A/H_w^\bot$. It induces a canonical map $p_2:Y_{K'}^{G_1} \rightarrow Y_{K''}^{G_1}$. 
\end{lem} 

\begin{proof}This map is defined to take an object $(A, \iota,
  \lambda, \eta, \eta_p)$ of $\mathbb{Y}_{K'}$ to the object $(A',
  \iota', \lambda', \eta', \eta_p') \in \mathbb{Y}_{K''}$ defined as follows:

\begin{itemize}
\item $A' = A/H^\bot_w$, where $H^\bot_w \subset A[w]$ is an order $p^3$
  group scheme, the orthogonal complement of $H_w$ for the Weil pairing. Write
  $\pi: A \rightarrow A'$ for the natural isogeny. 
\item $\iota'(x) =  \pi \circ \iota(x) \circ \pi^{-1}$,\item The quasi-polarization $\lambda'$ is obtained by descending the
  quasi-polarization $x_w^2. \lambda$ from $A$ to $A'$. \item $ \eta' = \pi \circ\eta $.\item $\eta'_p$ is the data of level structures at places $v \neq w$ above $p$ deduced from $\eta_p$ by the isomorphisms $\pi: A[ v] \rightarrow A'[v]$.
\end{itemize}

The ambiguity in the choice of $x_w$ disappears when we pass to the quotient stacks by the action of $(\ocal_{F})_{(p)}^{\times, +}$ and pass to the associated coarse moduli. 
\end{proof}

\begin{rem} There is another map $Y_{K'} \rightarrow Y_{K''}$ obtained
  by sending an abelian surface $A$ to $A/H_w$; however, we will  not
  need to make use of this map. 
\end{rem}

\subsubsection{ Change of level at $p$:  Siegel type correspondences}\label{subsect-functor-4} We now fix $K^p$ and  a place~$w$ above $p$. We let  $K_p = \prod_{v |p } K_v \subset G_1(\ZZ_p) $ be a reasonable compact open such that  $K_w = \mathrm{GSp}_4( \ocal_{F_w})$ (resp.\  $\Kli(w)$).  We let $K'_p = \prod_{ v \neq w} K_v \times \Si(w)$ be another reasonable level structure at $p$ (resp.\ $K'_p = \prod_{ v \neq w} K_v \times \Iw(w)$).  Set $K = K^pK_p$, $K' = K^p K'_p$. 

\begin{remark}[Warning] Note that the use of~$K$ and~$K'$ (and~$p_2$) in
this section (\S\ref{subsect-functor-4}) differs from that in
the previous section (\S\ref{subsect-functor-3}). Thus the reader should
be careful when these maps are used to note whether we are in
the Klingen or Siegel  setting (we
 indicate in any ambiguous  context by giving references to the corresponding section). We made
 this choice  since otherwise the number of required
subscripts would become excessively cumbersome.
\end{remark}

\begin{lem} There are natural forgetful maps $p_1:Y_{K'} \rightarrow
  Y_{K}$ and $p_1:Y_{K'}^{G_1} \rightarrow Y_{K}^{G_1}$ which are
  surjective and generically finite.\end{lem}

\begin{proof} We simply forget the level structure $L_w$ at $w$.
\end{proof} 

Recall that we have chosen an element $x_w \in F^{\times, +}$ which is a uniformizing element in ${F_w}$ and a unit in $F_v$ for all $v\neq w$ above $p$.

\begin{lem} There is a map $p_2:Y_{K'} \rightarrow Y_{K}$ depending on $x_w$. It induces a canonical map $p_2:Y_{K'}^{G_1} \rightarrow Y_{K}^{G_1}$. 
\end{lem} 

\begin{proof}  We take an object $(A, \iota, \lambda, \eta, \eta_p)$
  of $\mathbb{Y}_{K'}$. We define $(A', \iota', \lambda', \eta',
  \eta_p') \in \mathbb{Y}_{K}$ as follows:

\begin{itemize}
\item $A' = A/L_w$, call $\pi: A \rightarrow A'$ the isogeny.
\item $\iota'(x) =  \pi \circ \iota(x) \circ \pi^{-1}$.
\item The quasi-polarization $\lambda'$ is obtained by descending the quasi-polarization $x_w \lambda$ from $A$ to $A'$.
\item $ \eta' = \pi \circ\eta$.\item $\eta'_p$ is a data of level structures at places $v \neq w$ above $p$ deduced from $\eta_p$ by the isomorphisms $\pi: A[ v] \rightarrow A'[v]$.
\item In the case $K_w = \Kli(w)$, we define  $H'_w  = H_w^\bot/L_w \subset A'[w]$. 
\end{itemize}

The ambiguity in the choice of $x_w$ disappears when we pass to the quotient stacks by the action of $(\ocal_{F})_{(p)}^{\times, +}$ and pass to the associated coarse moduli. 
\end{proof}

\subsection{Automorphic vector bundles}\label{subsec: auto vector bundles} We now work over $\ZZ_p$, and
assume from now on that $p$ splits completely in $F$.  We let $S_p$ be the set of places of $F$ above $p$. We have a decomposition $\ocal_F \otimes_{\ZZ} \ZZ_{p}= \prod_{v |p} \ZZ_{p}$. We also denote by $v: \ocal_F \rightarrow \ZZ_p$ the projection on the $v$-component. 

\subsubsection{The principal bundle}\label{subsubsec: principal bundle} Over $Y_K$ we have a prime-to-$p$
isogeny class of abelian schemes and therefore we have a canonical
Barsotti--Tate group scheme~$\cG$. We let $\omega_{\cG}$ be its conormal sheaf.  The sheaf $\omega_{\cG}$ carries an action of $\ocal_F$. We have a decomposition $\ocal_F \otimes_{\ZZ} \ZZ_{p}= \prod_{v |p} \ZZ_{p}$ and accordingly, the sheaf $\omega_{\cG}$ decomposes as a product: $\omega_{\cG} = \prod_{v |p} \omega_{\cG,v}$ where each $\omega_{\cG,v}$ is a locally free sheaf of  rank $2$ over $Y_K$. 

\subsubsection{Weights for $G$ and $G_1$}\label{subsubsec: weights for
G and G1} By a dominant algebraic weight~$\kappa$ for~$G$  we mean a tuple $(k_v,l_v)_{v\in
  S_p}$ of integers such that $k_v \geq l_v$ for all $v \in S_p$. By a
classical algebraic weight we mean a dominant algebraic weight which furthermore satisfies~$l_v\ge 2$ for all~$v\in S_p$. We
will frequently write ``weight'' for ``dominant algebraic weight''
where no confusion can result (note though that we will later also
consider $p$-adic weights).  We
associate a locally free sheaf~$\omega^\kappa$ on $Y_K$ to each  weight    $\kappa$ by \[\omega^\kappa =
  \prod_v \mathrm{Sym}^{k_v-l_v} \omega_{\cG,v} \otimes \mathrm{det}^{l_v}
  \omega_{\cG,v}.\]

   By a weight~$\kappa$ for~$G_1$ we mean a tuple $((k_v,l_v)_{v\in
  S_p}, w)$ of integers with the property that $k_v \geq l_v$ and
$k_v-l_v\equiv w\pmod{2}$ for each~$v$; again, we say that~$\kappa$ is
classical algebraic if~$l_v\ge 2$ for all~$v\in S_p$. In fact, we
will insist that~$w$ is even, and we will shortly fix the
choice~$w=2$.  We  claim that given~$w$, there is a  canonical descent datum on $\omega^\kappa$ for the map $Y_{K} \rightarrow Y^{G_1}_K$. 

For clarity, we describe this descent datum on the level of the
groupoid $\mathbb{Y}_K$. For all $x \in  (\ocal_F)^{\times,+}_{(p)}$, we define an
isomorphism \[\omega^\kappa_{(A, \iota, x^{-1} \lambda, \eta,
  \eta_p)} = \omega^\kappa_{(A, \iota, \lambda, \eta, \eta_p)}
\rightarrow \omega^\kappa_{(A, \iota,  \lambda, \eta, \eta_p)}\] by
multiplication by $\prod_v v(x)^{(k_v+l_v-w)/2}$ (here the
first identification is the tautological one, noting that the
definition of~$\omega^\kappa$ does not depend on the
polarization). 

To check that this defines a descent datum, we have to show that it
respects the existing identifications from the action
of~$\ocal_F^{\times,+}(K^p)$.   If  $x \in \ocal_F^{\times,+}(K^p)$, then we may write $x=\epsilon^2$
for some $\epsilon\in \cO_F^\times\cap K^p$, and we have an isomorphism $\epsilon: A \rightarrow A$ which induces an  isomorphism in the groupoid $\mathbb{Y}_K$: 
\[ \epsilon: (A, \iota, \lambda, \eta, \eta_p) \rightarrow (A, \iota,
\epsilon^{-2} \lambda, \eta, \eta_p)\] and  an
isomorphism $$\omega^\kappa_{(A, \iota, \epsilon^{-2} \lambda, \eta,
  \eta_p)} = \omega^\kappa_{(A, \iota, \lambda, \eta, \eta_p)}
\stackrel{\epsilon^{*}} \rightarrow \omega_{(A, \iota,  \lambda,
  \eta, \eta_p)}  $$ which is multiplication by $\kappa(\epsilon)$
(again, the first equality is the tautological one, since $\omega^\kappa$ does not depend on the polarization). Now, $\kappa(\epsilon) = \prod_v v(\epsilon)^{k_v+l_v}=
 \prod_v v(\epsilon^2)^{(k_v+l_v-w)/2}\mathrm{N}_{F/\Q} (\epsilon)^w= \prod_v v(x)^{(k_v+l_v-w)/2}$
since $\mathrm{N}_{F/\Q} (\epsilon)^w =1$ by our assumption that~$w$ is even, so this agrees
with our the isomorphism defined above, as required.

This defines a
descent datum for the \'etale map $Y_K \rightarrow Y_K^{G_1}$.  This
descent datum is effective. Indeed, after first identifying  the sheaf
$\omega^\kappa$ on  various connected components of $Y_K$ we are
reduced to a finite \'etale descent for the group $\Delta(\Pi_0)$.

Although the descent datum depends on~$w$, we will regard~$w$ as fixed
(indeed, in the main arguments of the paper, we always take~$w=2$),
so we omit it from the notation, and simply denote the resulting sheaf
on $Y_K^{G_1}$ by~$\omega^\kappa$.

\begin{rem}\label{rem-descent-Hasse} We assume in this remark that we
  work over $\F_p$ rather than $\Z_p$. We denote by $Y_{K,1}$ and
  $Y_{K,1}^{G_1}$ the fibres of $Y_K$ and $Y_K^{G_1}$ over
  $\Spec\F_p$.  Let $\kappa = (k_v, l_v)_{v |p}$ be a  weight for
  $G$. We further assume that $k_v \equiv l_v \equiv 0 \mod
  (p-1)$. In this case, we claim that we can define a canonical descent datum for the sheaf $\omega^{\kappa}$, from $Y_{K,1}$ to $Y_{K,1}^{G_1}$. This rests on the observation that the character  $\ocal_{F}^\times \rightarrow \F_p^\times$  given by $\epsilon \mapsto \prod_{v |p} [ v( \epsilon)^{k_v +l_v}~\mod p]$ is trivial. Therefore we can define a descent datum for the action of $x \in  (\ocal_F)^{\times,+}_{(p)}$, via the tautological 
isomorphism \[\omega^\kappa_{(A, \iota, x^{-1} \lambda, \eta,
  \eta_p)} = \omega^\kappa_{(A, \iota, \lambda, \eta, \eta_p)}.\]
  This remark will be applied to the various Hasse invariants we will construct later. 
  \end{rem} 

Finally we will need to consider the canonical extensions of these sheaves to toroidal compactifications.  The conormal sheaf $\omega_{\mathcal{G}}/Y_K$ has a canonical
extension to $X_{K, \Sigma}$ given by $e^*\Omega^1_{A/X_{K, \Sigma}}$, where $A$ is the semi-abelian scheme of Theorem \ref{thm: main KWL compactification thm} (\ref{item: semi-abelian}) and $e$ is its identity section.  This gives an extension of the sheaves $\omega^\kappa$
to $X_{K,\Sigma}$ and an extension of the sheaves $\omega^{\kappa_{G_1}}$
to $X_{K, \Sigma}^{G_1}$.  We will denote these extensions by the same symbol.

 \subsection{Coherent cohomology and Hecke operators}\label{subsec:
   coherent cohomology and Hecke operators}
 \subsubsection{Basics}\label{subsec-basics-coherent}
 Let $\kappa = (k_v, l_v)$ be a weight. We will study  the cohomologies $\mathrm{R} \Gamma (X_{K, \Sigma}, \omega^{\kappa})$ and $\mathrm{R}\Gamma (X^{G_1}_{K, \Sigma}, \omega^{\kappa})$ as well as their cuspidal variants $\mathrm{R} \Gamma (X_{K, \Sigma}, \omega^{\kappa}(-D))$ and $\mathrm{R}\Gamma (X^{G_1}_{K, \Sigma}, \omega^{\kappa}(-D))$.
 
 \begin{lem}\label{lem: cohomologies independent of compactification} The cohomologies  $\mathrm{R} \Gamma (X_{K, \Sigma}, \omega^{\kappa})$,  $\mathrm{R}\Gamma (X^{G_1}_{K, \Sigma}, \omega^{\kappa})$, $\mathrm{R} \Gamma (X_{K, \Sigma}, \omega^{\kappa}(-D))$ and $\mathrm{R}\Gamma (X^{G_1}_{K, \Sigma}, \omega^{\kappa}(-D))$ are independent of $\Sigma$. 
 \end{lem}
 \begin{proof}
   This is immediate from Theorem~\ref{thm: main KWL compactification thm}~(\ref{item: pi sigma sigma prime}).
 \end{proof}
Because of this lemma,  we often drop $\Sigma$ from the notation. We now clarify the relationship between $\mathrm{R} \Gamma (X_{K}, \omega^{\kappa})$ and $\mathrm{R}\Gamma (X^{G_1}_{K}, \omega^{\kappa})$.

\begin{prop}\label{prop-splitGG_1} The pull back maps  \[\mathrm{R} \Gamma (X^{G_1}_{K}, \omega^{\kappa}) \rightarrow \mathrm{R}\Gamma (X_{K}, \omega^{\kappa})\] and \[\mathrm{R} \Gamma (X^{G_1}_{K}, \omega^{\kappa}(-D)) \rightarrow \mathrm{R}\Gamma (X_{K}, \omega^{\kappa}(-D))\] split in the derived category of $\Z_p$-modules. 
\end{prop}

\begin{rem} It is often easier to work over $X_{K}$ rather than $X^{G_1}_{K}$ because the former has a clear moduli interpretation. Proposition~\ref{prop-splitGG_1} tells us that we can easily transfer a good property of the cohomology over $X_{K}$ to a property over $X^{G_1}_K$. 
\end{rem} 
 
\begin{proof}[Proof of Proposition~\ref{prop-splitGG_1}] Attached to the weight $\kappa$ is a descent datum (see~\S\ref{subsubsec: weights for
G and G1}) which takes the form of an action of $(\ocal_F)^{\times,+}_{(p)}$ on the sheaf $\omega^\kappa$ over $X_{K}$. Namely, for all $\epsilon \in (\ocal_F)^{\times,+}_{(p)}$, there is an isomorphism $ \epsilon: \epsilon^{*} \omega^\kappa \rightarrow \omega^\kappa$ satisfying the usual cocycle relation. This map induces a map on cohomology: $$\epsilon: \mathrm{R} \Gamma (X_{K}, \omega^{\kappa}) \rightarrow \mathrm{R} \Gamma (X_{K}, \epsilon^{*} \omega^{\kappa}) \rightarrow \mathrm{R} \Gamma (X_{K}, \omega^{\kappa})$$ and defines the group action. 

Recall that there is a commutative diagram: 

\begin{eqnarray*}
\xymatrix{ X_{K} \ar[d]^{\Pi_0} \ar[r] &  X_{K}^{G_1}  \ar[d]^{\Pi_0^{G_1}}\\
 \ZZ_{(p)}^{\times,+} \backslash (\A^{\infty,p}\otimes F)^\times/ \nu(K^p)(1) \ar[r]^{\pi}&  (\ocal_F)_{(p)}^{\times,+} \backslash (\A^{\infty,p}\otimes F)^\times/ \nu(K^p)(1)}
\end{eqnarray*}

Each Galois orbit $c \in  [\ZZ_{(p)}^{\times,+} \backslash (\A^{\infty,p}\otimes F)^\times/ \nu(K^p)(1)]/ \mathrm{Gal}(\overline{\Q}/\Q)$ determines a connected component of $\ZZ_{(p)}^{\times,+} \backslash (\A^{\infty,p}\otimes F)^\times/ \nu(K^p)(1)$, and its fibre is   a connected component $X_{K,c}$ of $X_{K, \Sigma}$ which is a proper scheme over $\Spec \ZZ_p$.  Obviously   $\mathrm{R} \Gamma (X_{K}, \omega^{\kappa}) = \prod_{c} \mathrm{R} \Gamma (X_{K,c}, \omega^{\kappa})$ and for all $\epsilon \in 
 (\ocal_F)^{\times,+}_{(p)}$, we have an isomorphism $ \epsilon: \mathrm{R} \Gamma (X_{K,  \epsilon\cdot c}, \omega^{\kappa}) \rightarrow \mathrm{R} \Gamma (X_{K, c}, \omega^{\kappa})$. 
 
 The subgroup that fixes a component $X_{K, c}$ is denoted by $\ocal_{F}^{\times, +}(\Pi_0)$ and the action of this group on $X_{K,  c}$ and $\mathrm{R} \Gamma (X_{K,  c}, \omega^{\kappa})$ actually factors through the finite group  $\Delta(\Pi_0)$. Let $\pi(c) $ be the image of $c$ in $$[(\ocal_F)_{(p)}^{\times,+} \backslash (\A^{\infty,p}\otimes F)^\times/ \nu(K^p)(1)]/ \mathrm{Gal}(\overline{\Q}/\Q).$$ This determines a connected component $X^{G_1}_{K, \pi(c)}$ of $X_{K}^{G_1}$ and the map $X_{K,  c} \rightarrow X_{K,  \pi(c)}^{G_1}$ is a finite \'etale  cover with group $\Delta(\Pi_0)$. 
 
 It follows from Lemma~\ref{lem: finite group etale splitting} below that $\mathrm{R} \Gamma (X^{G_1}_{K, \pi(c)}, \omega^{\kappa})$ is split in  $\mathrm{R} \Gamma (X_{K, c}, \omega^{\kappa})$, and therefore  the map $$\mathrm{R} \Gamma (X^{G_1}_{K}, \omega^{\kappa}) = \bigoplus_{\pi(c)} \mathrm{R} \Gamma (X^{G_1}_{K,\pi(c)}, \omega^{\kappa}) \rightarrow  \prod_{c} \mathrm{R} \Gamma (X_{K,c}, \omega^{\kappa}) =  \mathrm{R} \Gamma (X_{K}, \omega^{\kappa}) $$ is split. 
 \end{proof}

  \begin{lem}\label{lem: finite group etale splitting} Let $G$ be a finite group. Let $I_G \subset \Z[G]$ be
    the augmentation ideal. Let $f: T \rightarrow S$ be a finite
    \'etale morphism with Galois group $G$. Then $f_{*} \oscr_T = \oscr_S \oplus I_G \otimes_{\Z[G]} f_{*} \oscr_T$. 
 \end{lem}
 
 \begin{proof} There is an obvious map of coherent sheaves $\oscr_S \oplus I_G \otimes_{\Z[G]} f_{*} \oscr_T \rightarrow \ f_{*} \oscr_T$. The sheaf $f_{*} \oscr_T$ is a locally free sheaf (for the \'etale topology) of $\oscr_{X}[G]$-modules.  Therefore, the above map is an isomorphism as this can be checked locally for the \'etale topology.  
 \end{proof}

\subsubsection{Abstract Hecke algebras} Let $ \mathcal{H} = \mathcal{C}^\infty_c ( G_1(\A^{\infty})
    \doubleslash K, \ZZ_p)$ be the convolution algebra of locally
    constant, bi-$K$ invariant, compactly supported functions  on
    $G_1(\A^{\infty})$ with coefficients in $\ZZ_p$. (The Haar measure
    is a product of local Haar measures, normalized by $\mathrm{vol}(K_t)=1$ for all finite places $t$ of $F$.) If $S$ is a finite set of places of $F$, we let $\mathcal{H}^{S}$ be the subalgebra of $\mathcal{H}$ of functions whose restriction to $\mathrm{GSp}_4(F_s)$ is the characteristic function of~$K_s$ for all $s \in S$.  
 For all finite places $s$, we let $\mathcal{H}_s$ be the local Hecke algebra $\mathcal{C}^\infty_c ( \mathrm{GSp}_4(F_s) \doubleslash K_s, \ZZ_p)$, so that $\mathcal{H} = \otimes'_s \mathcal{H}_s$.

\subsubsection{Cohomological correspondences ---
  motivation}\label{subsubsec: cohomological correspondences motivation}We begin by
giving some brief motivation for the way in which we define Hecke
operators on coherent cohomology (following~\cite{pilloniHidacomplexes}).

As usual, the geometric interpretation of Hecke operators is via
correspondences 
  \begin{eqnarray*}     
 \xymatrix{ &C   \ar[rd]^{p_1}\ar[ld]_{p_2}&  \\
  X & & Y}
  \end{eqnarray*}
 (Giving an \emph{integral} definition of the correspondence
 associated to a Hecke operator at a place dividing~$p$ is in general difficult. This question will be addressed later in the paper in some very special cases.)

  Let  ~$\cF,\cG$ be coherent sheaves on~$X, Y$. We assume that we
  have a map of sheaves $p_2^{*} \cF \rightarrow p_1^{*} \cG$. When
  $\cF$ and $\cG$ are automorphic vector bundles (which will typically
  be the case for us), this map is provided by the differential of the universal isogeny over $C$. 
  
  One would like to use the correspondence to define a  corresponding map
 on cohomology $\mathrm{R}\Gamma (X, \cF) \rightarrow \mathrm{R}\Gamma(Y, \cG)$.  This map could be  defined by first taking the pull back via~$p_2: \mathrm{R}\Gamma (X, \cF) \rightarrow \mathrm{R}\Gamma (C, p_2^{*} \cF)$, then  using the map $p_2^{*} \cF \rightarrow p_1^{*} \cG$ to get to $\mathrm{R}\Gamma (C, p_1^{*} \cG)$, and finally applying some trace map to $\mathrm{R}\Gamma(Y, \cG)$. In other words, the action of the correspondence on cohomology should take the form of a map $T:R(p_1)_*p_2^*\cF\to\cG$. 
  There are, however, at least two serious difficulties with making such a definition in our
  context.
  
  The first obvious difficulty is the existence of the trace map, because in general one cannot assume that $p_1$ is finite flat.  Nevertheless, in our cases the existence of the trace map will follow from the machinery of duality in coherent cohomology and the existence of certain fundamental classes, which can be constructed because the schemes $C,X,Y$ will  have reasonable geometric properties over the base. 
  
  The second difficulty
  (which already arises for modular forms for~$\GL_2/\Q$) is
  that the action of the correspondences defining the Hecke operators at places
  dividing~$p$ is typically divisible by a positive power of~$p$, so
  that one has to divide by this power in order to define the correct
  operator mod~$p$. It is hard to check this
  divisibility at the level of the derived category.

The solution  to this introduced in~\cite{pilloniHidacomplexes} (which we also employ here) is as follows. 
By adjunction we can view~$T$ as a map $T:p_2^*\cF\to
p_1^!\cG$, and in favourable circumstances~$p_1^!\cG$ will be a sheaf (and not
merely a complex). Furthermore it will be sufficiently nice that we
can check the condition that~$T$ is divisible by a power of~$p$ after
restricting to the complement of a codimension~$2$ locus, and define
our normalized Hecke operators.

\subsubsection{Duality for coherent complexes}\label{subsec: duality
  for coherent complexes}

We let~$S$ be an affine Noetherian scheme. We say that a morphism
$f:X\to Y$ of $S$-schemes is \emph{embeddable} if there is a smooth
$S$-scheme $P$ such that $f$ can be factored as a
composite \[X\stackrel{i}{\to}P\times_S Y\to Y \] where~$i$ is finite
and the second map is the natural projection. We say that~$f$ is
\emph{projectively embeddable} if~$p$ can be taken to be a projective
space over~$S$. In our applications of this material all of our maps
will be obviously projectively embeddable (essentially because our
Shimura varieties are quasi-projective), and we will not comment
further upon this.

As usual we write $D_{qcoh}(\cO_X)$ for the
derived category of $\cO_X$-modules with quasi-coherent cohomology
sheaves, and $D^+_{qcoh}(\cO_X)$ for the bounded-below version.
Then if $f:X\to Y$ is an embeddable morphism of $S$-schemes, there is an exact
functor of triangulated categories \[f^!: D^+_{qcoh}(\cO_Y)\to
  D^+_{qcoh}(\cO_X). \] 
  
  If  $f$ is projectively embeddable,  the
functor $f^!$ is a right adjoint to~$Rf_*$ and  there is a natural transformation
$Rf_*f^!\implies\Id$ of endofunctors of~$ D^+_{qcoh}(\cO_Y)$, which we
refer to as the trace map.

 If $X\to S$ is a local complete intersection then we write~$K_{X/S}$
for the relative canonical sheaf, which may be defined as the
determinant of the corresponding cotangent complex. The following is~\cite[Cor.\ 4.1.3.1]{pilloniHidacomplexes}.
\begin{lem}
  \label{lem: map between lci means dualizing is a sheaf}Let $f:X\to
  Y$ be an embeddable morphism between two embeddable  $S$-schemes, such that $X\to S$, $Y\to
  S$ are both local complete intersections of pure relative
  dimension~$n$. Then
  $f^!\cO_Y=K_{X/S}\otimes_{\cO_X}f^*K_{Y/S}^{-1}$ is an
  invertible sheaf.
\end{lem}

We will make repeated use of the following lemma.

\begin{lem}\label{lem: base change of correspondences}Suppose that $f:X\to Y$
  is an embeddable morphism of embeddable $S$-schemes, each of which is a local complete
  intersection of pure relative dimension~$n$ over~$S$. Let $h$ be a
  section of a line bundle~$\cL$ over~$Y$, and suppose that neither $h$ nor
  $f^*h$ is a zero-divisor. Write $Y_{h=0}$ for the vanishing locus
  of~$h$, and $X_{h=0}$ for the vanishing locus of~$f^*h$.

Then for any locally free sheaf~$\cF$ on~$Y$, we have an equality of
invertible sheaves \[(f^!\cF)|_{X_{h=0}}=f^!(\cF|_{Y_{h=0}}).\]
\end{lem}
\begin{proof}
  This follows from~\cite[Prop.\ III.8.8]{Hartshorne}. More precisely, note that $\cO_{Y_{h=0}}$  is represented by the
  perfect complex of~$\cO_Y$-modules $\cL^{-1}\stackrel{h}{\to}\cO_{Y}$
  (here we use that~$h$ is not a zero-divisor). In addition, by Lemma~\ref{lem: map between
    lci means dualizing is a sheaf}, $f^!\cF$ is a sheaf, and it
  follows from the
  assumption that neither $h$ nor $f^*h$ is a zero-divisor that the 
   derived tensor products in~\cite[Prop.\ III.8.8]{Hartshorne} are in our case 
  given by the usual tensor product~$\otimes$.
\end{proof}

\subsubsection{Fundamental classes}\label{subsubsec: fundamental class} In two particular situations, we now construct a natural map
$$\Theta: \cO_X= f^{*} \ocal_{Y} \rightarrow f^! \ocal_{Y}$$ which
we call the \emph{fundamental class}.

We firstly consider what we call the \emph{lci situation}, which is
the case that: 
\begin{itemize}
\item $X$ and $Y$ are local complete intersections over $S$ of
  the same relative dimension,
\item $X$ is normal, and 
\item  there is an open
  $V \subset X$ which is smooth over $S$, whose complement is of
  codimension $2$ in $X$, and an open $U \subset Y$ which is smooth and
  such that $f(V) \subset U$.
\end{itemize}
In this situation,
$f^! \ocal_Y$ is an invertible sheaf by Lemma~\ref{lem: map
  between lci means dualizing is a sheaf}, so by the algebraic
Hartogs' lemma,  it is enough to specify the fundamental
class over $V$ (note that $X$ is normal by assumption). Again by
Lemma~\ref{lem: map between lci means dualizing is a sheaf} we have
$f^! \ocal_Y\vert_V = \det \Omega^{1}_{V/S} \otimes f^*(\det
\Omega^{1}_{U/S})^{-1}$, so over $V$, we can define the fundamental
class to be the determinant of the map
\[\mathrm{d} f: f^*\Omega^{1}_{U/S} \rightarrow
\Omega^{1}_{V/S}.\]

The other case we consider is the \emph{finite flat} situation, in
which $f:X\to Y$ is a finite flat map, so that~$f_*$ is exact, and 
\[f_*f^! \ocal_Y = \SheafHom_{\cO_Y} ( f_{*} \ocal_{X}, \ocal_{Y}).\]
We have the usual trace morphism
$\tr_f: f_{*} \ocal_{X} \rightarrow \ocal_Y$, and we define the
fundamental class
$f_*\cO_X\to \SheafHom_{\cO_Y} ( f_{*} \ocal_{X}, \ocal_{Y})$ by
$\Theta(1) = \tr_f$.

Note that if $X\to Y$ is a finite flat
morphism and $X,Y$ are both smooth over~$S$, then the morphism $X\to
Y$ is automatically a local complete intersection. The following compatibility between these definitions is~\cite[Lem.\
4.2.3.1]{pilloniHidacomplexes}.

\begin{lem}Suppose that $X\to Y$ is finite flat, and that $X$, $Y$ are
  both smooth over~$S$. Then \[\mathbb{L}_{X/Y}\isoto
    [\Omega^1_{Y/S}\otimes_{\cO_Y}\cO_X\stackrel{df}{\to}\Omega^1_{X/S}], \]
  and the determinant $\det(df)\in \omega_{X/Y}=f^{!}\cO_Y$ is the
  trace map~$\tr_f$.
  \end{lem}

\subsubsection{Base change for open immersions}\label{subsubsec: base change}Consider a
Cartesian diagram \[\xymatrix{X'\ar[r]^{j}\ar[d]^{f'}& X\ar[d]^f\\
    Y'\ar[r]^{i}&Y} \]If $i$ is an open immersion, and~$f$ is in
either of the finite flat or lci situations, then so is~$f'$. Since
$i^!=i^*$ and~$j^!=j^*$, we have $j^*f^!=(f')^!i^*$, and if $f$ has  fundamental
class~$\Theta$, then $j^*\Theta$ is the fundamental class of~$f'$. 

\subsubsection{Fundamental classes and divisors}\label{subsubsec:
  fundamental class and divisors}We now briefly recall the results of~\cite[\S
4.2.4]{pilloniHidacomplexes}, which show that the
correspondences we define below are suitably well behaved on the boundaries of our
compactified Shimura varieties.

Let $D_X\into X$, $D_Y\into Y$ be two effective reduced Cartier
divisors with respect to~$S$, with the properties that $f:X\to Y$
restricts to a map $f|_{D_X}:D_X\to D_Y$, and the induced map $D_X\to
f^{-1}(D_Y)$ is an isomorphism of topological spaces. Write $X^{sm}$,
$Y^{sm}$ for the smooth loci of $X,Y$. The following is~\cite[Lem.\ 4.2.4.1]{pilloniHidacomplexes}.

\begin{lem}Suppose either that we are in the finite flat situation; or
  that we are in the lci situation and that furthermore $D_X\cap X^{sm}$
  and $D_Y\cap Y^{sm}$ are normal crossings divisors. 

Then the fundamental class $\Theta:\cO_X\to f^!\cO_Y$ restricts to a
morphism $\cO_X(-D_X)\to f^!\cO_Y(-D_Y) $.
\end{lem}

\subsubsection{Traces and restriction}\label{section-traceandres}

In this paper we will have to study how Hecke operators behave with respect to restriction to subschemes of the Shimura variety. This section contains some  preliminary material. Consider the following setup:\begin{itemize}
\item $f:X\to Y$ a finite flat map between smooth varieties over a field $k$.
\item $D\subset Y$ is a smooth Cartier divisor.
\item $f^{-1}(D)=nD'$ for $D'\subset X$ a smooth Cartier divisor.
\end{itemize}
In this setting we have the following:
\begin{itemize}
\item Trace maps on canonical bundles
\begin{equation*}
f_*K_X\to K_Y
\end{equation*}
and
\begin{equation*}
f_*K_{D'}\to K_D.
\end{equation*}
\item Adjunction isomorphisms
\begin{equation*}
K_D\simeq K_Y(D)|_{D}
\end{equation*}
and
\begin{equation*}
K_{D'}\simeq K_X(D')|_{D'}.
\end{equation*}
\end{itemize}

If $\cL$ is a line bundle on $Y$, we can use the projection formula to get a map: $f_*(K_X \otimes_{\cO_X} f^*\mathcal{L}) \to K_Y\otimes_{\cO_Y} \mathcal{L}$. We call such a map a twisted trace map. We use a similar terminology over $D$. The goal of this section is to prove the following compatibility between them.
\begin{prop}\label{prop: trace restrict to cartier}
There is a commutative diagram
\begin{equation*}
\xymatrix{f_*(K_X(-(n-1)D'))\ar[r]\ar[d]&K_Y\ar[d]\\
f_*(K_{D'}\otimes \O_X(-nD')|_{D'})\ar[r]&K_D\otimes\O_Y(-D)|_{D}}
\end{equation*}
Here the vertical maps are restriction followed by adjunction, the top
horizontal map comes from the inclusion of $K_X(-(n-1)D')$ in $K_X$
followed by the trace, and the bottom horizontal arrow is the twisted
trace for $f:D'\to D$ and the line bundle $\O_Y(-D)|_D$ \emph{(}note that $f^*\O_Y(-D)|_D=\O_X(-nD)|_{D'}$\emph{)}.  \end{prop}
\begin{proof}
We write $\cI=\O_Y(-D)$ for the ideal sheaf of $D$ and $\cI'=\O_X(-D')$ for the ideal sheaf of $D'$.  First consider the following commutative diagram:
\begin{equation*}
\xymatrix{f_*{\cI'}^{n-1} \SheafHom_{\O_Y}(\O_X,\O_Y)\ar@{^{(}->}[r]\ar@{..>}[d]&f_*\SheafHom_{\O_Y}(\O_X,\O_Y)\ar[r]\ar[d]&\O_Y\ar[d]\\
f_*\SheafHom_{\O_Y/\cI}(\O_X/{\cI'},\O_Y/\cI)\ar@{^{(}->}[r]&f_*\SheafHom_{\O_Y/\cI}(\O_X/\cI\O_X,\O_Y/\cI)\ar[r]&\O_Y/\cI
}
\end{equation*}
where $\SheafHom_{\O_Y}(\O_X,\O_Y)$ is sheaf of $\O_Y$-homomorphisms from $f_\star\O_X$ to $\O_Y$, which we view as a coherent sheaf of $\O_X$-modules. By definition $\SheafHom_{\O_Y}(\O_X,\O_Y) = f^! \O_Y$. 

Consider first the square on the right: the horizontal maps  are evaluation at 1, while the vertical maps are given by reduction modulo $\cI$, and it is clear that this square commutes.

Now we consider the left hand square: the horizontal maps are the obvious inclusions so we must explain why the dotted arrow exists.  But a local section $s$ of $\cI'^{n-1}\SheafHom_{\O_Y}(\O_X,\O_Y)$ will send $\cI'$ into $\cI$ (using that ${\cI'}^n=\cI\O_X$) and hence the reduction of $s$ mod $\cI$ factors through $\O_X/\cI'$.

Finally we note that the square in the statement of the proposition tensored with $K_Y^{-1}$ may be identified with the outer rectangle of this diagram because we have $K_X\otimes K_Y^{-1}\simeq f^!\O_Y=\SheafHom_{\O_Y}(\O_X,\O_Y)$. \end{proof}

\subsection{Cohomological correspondences --- definitions} Let $S$ be a Noetherian scheme.  Let $X$, $Y$ be two    $S$-schemes. 

\begin{defn}  A \emph{correspondence} $C$ over $X$ and $Y$ is  a diagram of  $S$-morphisms:
\begin{eqnarray*}
\xymatrix{ & C \ar[rd]^{p_1}  \ar[ld]_{p_2} & \\
X & & Y }
\end{eqnarray*}
where $X$, $Y$, $C$ have the same pure relative dimension over $S$ and the morphisms $p_1$ and $p_2$ are projectively embeddable. 
\end{defn}

 Let $\mathcal{F}$ be a coherent sheaf over $X$ and $\mathcal{G}$ a
 coherent sheaf over $Y$.

\begin{defn} A \emph{cohomological correspondence} from $\mathcal{F}$ to $\mathcal{G}$ is the data of a correspondence $C$ over $X$ and $Y$ and a map $ T: \mathrm{R}(p_1)_{*} p_2^{*} \mathcal{F} \rightarrow \mathcal{G}$.
\end{defn}

The map $T$ can be seen, by adjunction, as a map $p_2^{*} \mathcal{F} \rightarrow p_1^! \mathcal{G}$. It gives rise to a map still denoted by $T$ on cohomology:
$$  \mathrm{R} \Gamma( X, \mathcal{F}) \stackrel{p_2^{*}}\rightarrow
\mathrm{R}\Gamma(C, p_2^{*} \mathcal{F}) = \mathrm{R}\Gamma(Y,
\mathrm{R}(p_1)_{*} p_2^{*} \mathcal{F}) \stackrel{T}\rightarrow
\mathrm{R}\Gamma( Y, \mathcal{G}).$$

\subsubsection{Hecke action away from $p$}\label{sect-hecke-away-from-p}  Let $K = K^p
K_p $ be a reasonable
compact open subgroup of $G_1(\A_f)$. Let $\mathcal{H}^p = \mathcal{C}^\infty_c ( G_1(\A^{p,\infty})
    \doubleslash K^p, \ZZ_p)$ be the Hecke algebra away from
    $p$. 

We claim that there is an action
of $\mathcal{H}^p$ on $\mathrm{R}\Gamma ( X_{K, \Sigma},
\omega^{\kappa})$ and $\mathrm{R}\Gamma ( X^{G_1}_{K, \Sigma},
\omega^{\kappa})$. To this end, let $g \in G_1( \A^{\infty,p})$. We will define an endomorphism of $\mathrm{R}\Gamma ( X_{K}, \omega^{\kappa})$ which corresponds to  the action of the double class $[{K^pgK^p}]$.

 We define (for suitable choices of cone decompositions omitted from
 the notation) a correspondence: \begin{eqnarray*}     
 \xymatrix{ &X_{K \cap g K g^{-1} }   \ar[rd]^{p_1}\ar[ld]_{p_2}&  \\
  X_K & & X_K}
  \end{eqnarray*}
  where $p_1$ is the map induced from the inclusion $K \cap g K g^{-1}
  \subset K$ and the functoriality of~\S\ref{subsect-functor-1}.
  
  The map $p_2$ is the composite of the map $[g]: X_{K \cap g K
    g^{-1} } \rightarrow X_{K \cap g^{-1} K g }$ (see
  \ref{subsect-functor-2})  and the natural map $X_{K \cap g^{-1} K g
  } \rightarrow X_{K }$ deduced from the inclusion $K \cap g K g^{-1}
  \subset K$ and functoriality of~\S\ref{subsect-functor-1}. 
  
  We have a canonical isomorphism $p_2^{*} \omega^\kappa \isoto
  p_1^{*} \omega^\kappa$, because the construction of the sheaf
  $\omega^\kappa$ depends only on the $p$-divisible group. Moreover,
  because $X_K$ and $X_{ K \cap g K g^{-1}}$ are lci and smooth outside codimension 2 
  (for a cofinal
  subset of the set of all polyhedral cone decompositions),
  there is a fundamental class $p_1^{*} \ocal_{X_K} \rightarrow
  p_1^! \ocal_{X_K}$, extending the trace for the finite \'etale map $p_1$ on the interior, which we can tensor with $p_1^*\omega^\kappa$ to obtain a map $p_1^*\omega^\kappa\to p_1^!\omega^\kappa=p_1^!\ocal_{X_K}\otimes p_1^*\omega^\kappa$.

 Composing the maps $p_2^*\omega^\kappa\to p_1^*\omega^\kappa$ and $p_1^*\omega^\kappa\to p_1^!\omega^\kappa$ we obtain a cohomological correspondence $\Theta_g: p_2^{*} \omega^\kappa \rightarrow p_1^! \omega^\kappa$ which induces the operator $[{K g K}]$ on cohomology:  
 $$ \mathrm{R}\Gamma(X_K, \omega^\kappa) \rightarrow \mathrm{R}\Gamma( X_{K \cap g K g^{-1}}, p_2^{*} \omega^\kappa) \stackrel{\Theta_g}\rightarrow \mathrm{R}\Gamma( X_{K \cap g K g^{-1}}, p_1^! \omega^\kappa) \stackrel{\mathrm{Tr}}\rightarrow \mathrm{R}\Gamma( X_{K }, \omega^\kappa)$$
 where the last map is induced by the adjunction $\mathrm{Tr}: \mathrm{R}(p_1)_{*} p_1^! \omega^\kappa \rightarrow \omega^\kappa$. 
 
 We have a similar definition on cuspidal cohomology. Moreover, all these definitions commute with the action of $(\ocal_F)_{(p)}^{\times, +}$ and therefore we also get an action on the cohomology $\mathrm{R}\Gamma ( X^{G_1}_{K}, \omega^{\kappa})$ and $\mathrm{R}\Gamma ( X^{G_1}_{K}, \omega^{\kappa}(-D))$.
 
 The characteristic functions of the double classes
 $[{K^pgK^p}]$  generate $\mathcal{H}^p$ as a
 $\ZZ_p$-module. In Proposition~ \ref{prop-cousin-heckeaction} below
 we prove that when $K_p=\prod_{v\mid p}\GSp_4(\ZZ_p)$ is spherical, the actions we just defined  of the
 $[{K^pgK^p}]$ are compatible with products in
 $\mathcal{H}^p$  (the composite action of  $[{K^pg_1K^p}]$ and
 $[{K^pg_2K^p}]$ is equal to the action of
 $[{K^pg_1K^p}] [{K^pg_2K^p}]$ decomposed into sum of
 elementary double classes) so that we get an action of the Hecke
 algebra $\mathcal{H}^p$.

 The difficulties  come from the boundary. Away from the boundary, all the correspondences are finite \'etale and one can follow the discussion of \cite[Chap.\ VII, \S3]{MR1083353},  to show the compatibility. Following that reference, it should be  possible to show in a similar fashion  that the action of the double class is compatible with product in the Hecke algebra on the compactified Shimura variety, but giving all the details  would involve a delicate study of the composition of the correspondences  at the boundary. We instead give a different  ad hoc proof by exhibiting special complexes computing the cohomology. These complexes are Cousin complexes associated with the Ekhedal--Oort stratification on the Shimura variety. The action of all double classes 
 $[{K^pgK^p}]$  on the cohomology is given by a canonical   action on the complex. Moreover, each term of the complex is the global sections  of a certain sheaf and the restriction of the sections of this sheaf  to the interior of the Shimura variety is an embedding. We are therefore  able to prove  that the action of the double classes is compatible with products in the Hecke algebra because we know this holds on the non-compact Shimura variety.  
 
\begin{rem} Over $\Q_p$, the property   that the action of the double class is compatible with product in the Hecke algebra follows from  \cite[Prop.\ 2.6]{MR1064864}. The strategy  of that paper is to define an action of the group $\mathrm{G}_1(\mathbb{A}^p_f)$ after passing to the limit over  the level $K^p$ and then deduce an action of the Hecke algebra at a finite level, but this strategy requires more work over $\ZZ_p$ because at some points one needs to control the cohomology of finite groups (which vanishes in characteristic zero). 
Nevertheless, this is enough to prove that the Hecke algebra $\mathcal{H}^p$ acts on the torsion free part of the cohomology (which embeds in the cohomology with $\Q_p$-coefficients). 
\end{rem}
\subsubsection{Cousin complexes}\label{subsec: cousin complexes} Our main reference  for this section is \cite{KEMPF1978310}. 
Let $X$ be a topological space. Let $Sh_X(Ab)$ be the category of abelian sheaves on $X$. For a subset $Z \subseteq X$ and abelian sheaf $\mathcal{F}$ we denote by $\Gamma_Z(\mathcal{F})$ the subsheaf of $\mathcal{F}$ of sections supported on $Z$.  Let $Z: Z_0 = X \supseteq Z_1 \supseteq \cdots \supseteq Z_{n} \cdots$ be a decreasing sequence of closed subsets of $X$ (called a filtration). For any abelian sheaf $\mathcal{F}$ on $X$, one can build the Cousin complex of $\mathcal{F}$ with respect to  the filtration $Z$, denoted by $\Cous_{Z}(\mathcal{F})$ \cite[p.\ 357]{KEMPF1978310}.

 The Cousin complex $\Cous_Z(\mathcal{F})$ is a complex of abelian sheaves in positive degree. The object in degree $i$ is $ \mathcal{H}^i_{Z_i/Z_{i+1}}(\mathcal{F})$, where
 $\mathcal{H}^k_{Z_i/Z_{i+1}}(\cdot)$ is (by~\cite[Lem.\ 7.3]{KEMPF1978310}) the $k$-th derived functor of
 the functor: \begin{eqnarray*}
 \mathrm{Sh}_X(Ab) & \rightarrow &  \mathrm{Sh}_X(Ab)  \\
  \mathcal{G} &  \mapsto & [ U \mapsto \Gamma_{Z_i\setminus Z_{i+1}}( U \setminus Z_{i+1}, \mathcal{G})] 
  \end{eqnarray*}
  The differential $\mathcal{H}^i_{Z_i/Z_{i+1}}(\mathcal{F}) \rightarrow \mathcal{H}^{i+1}_{Z_{i+1}/Z_{i+2}}(\mathcal{F}) $ is induced by a certain boundary map. The Cousin complex has an augmentation $\mathcal{F} \rightarrow \Cous(\mathcal{F})$. 
  
  We now specialize the discussion: $X$ is a Noetherian  scheme and $\mathcal{F}$ is a quasi-coherent sheaf. Then $\Cous(\mathcal{F})$ is a complex of quasi-coherent sheaves.    
 
 We have the following theorem:
 
 \begin{thm}\label{thmkempf} Let  $X$ be a Noetherian scheme with a
   filtration $Z$ by closed subschemes that satisfies:
  
  \begin{enumerate}
  
  \item $\codim_{X} (Z_i) \ge i $. 
 \item The morphism  $Z_i \setminus Z_{i+1} \rightarrow X$ is affine for all $i$.
 \end{enumerate}
Let $\mathcal{F}$ be a maximal  Cohen--Macaulay coherent sheaf on $X$. Then $\Cous_Z(\mathcal{F})$ is quasi-isomorphic to  $\mathcal{F}$.
 \end{thm}
 
 \begin{proof} This  follows from \cite[Thm.\
   10.9]{KEMPF1978310} (by definition a sheaf $\mathcal{F}$ is locally
   Cohen--Macaulay with respect to a filtration~$Z$ if 
   $\Cous_Z(\mathcal{F})$ is a resolution of~ $\mathcal{F}$, see \cite[p.\ 358]{KEMPF1978310}).\end{proof}
 
 \begin{rem}\label{rem: affine implies acyclic kempf} If we further assume that each  $Z_i \setminus Z_{i+1}$ is
   affine, then  $\Cous_Z(\mathcal{F})$ is a complex of acyclic sheaves
   by  \cite[Thm.\ 9.6]{KEMPF1978310}.
 \end{rem} 
 
 One can sometimes compute  the complex $\Cous_Z(\mathcal{F})$ more
 explicitly. Write $U_{i+1} = X \setminus Z_{i+1}$, and write
 $j_{i+1}: U_{i+1} \hookrightarrow X$ for the inclusion.  Under  the
 assumption that  $Z_i \setminus Z_{i+1} \rightarrow X$ is affine, 
 we have by \cite[Lem.\ 8.5(e)]{KEMPF1978310} (note that the spectral
 sequence there degenerates by \cite[Thm.\
 9.6(c)]{KEMPF1978310}, as in the proof of~\cite[Thm.\
 9.5]{KEMPF1978310}):\numequation\label{eqn: term in Cousin as
   pushforward of local cohomology} \mathcal{H}^k_{Z_i/Z_{i+1}}(\mathcal{F}) = (j_{i+1})_* \mathrm{R}^k\Gamma_{Z_i \setminus Z_{i+1}}(\mathcal{F}\vert_{U_{i+1}}).\end{equation}
In general, for a  Noetherian scheme $X$ and a closed subset $Z$
defined by an ideal sheaf $\mathcal{I}$, we have $$ \mathrm{R}^i\Gamma_{Z} (\mathcal{F}) = \colim_n \mathcal{E}xt^i ( \ocal_{X}/\mathcal{I}^n, \mathcal{F})$$and these  $\mathcal{E}xt$ sheaves can be computed by taking projective  resolutions of $ \ocal_{X}/\mathcal{I}^n$. We also remark that in the previous limit, we can replace the ideals $\mathcal{I}^n$ by any other decreasing sequence of ideals $\{ J_n\}$ with the property that for all $n$, there is $k$ and $k'$   such that $J_{k'} \subset I^n \subset J_k$. 

\begin{example}\label{ex-computation} We are going to compute these
  $\mathcal{E}xt$ sheaves in a special case.  Assume that we have
  effective Cartier divisors $ 
  \mathcal{O}_X \stackrel{s_t}\rightarrow \mathcal{L}_t$ for $1 \leq t \leq i$ and assume  that they intersect
  properly, by which we mean  that  for all $n$, the ``twisted'' Koszul complex: 

$$\Kos(s_1^n, \cdots, s_i^n): 0 \rightarrow \bigotimes_{t} \mathcal{L}^{-n}_t  \rightarrow \oplus_{t'} \bigotimes_{t \neq t'}  \mathcal{L}^{-n}_t \rightarrow \cdots \rightarrow \oplus_t \mathcal{L}^{-n}_t \rightarrow \ocal_X  \rightarrow 0$$ 
is a projective resolution of $\ocal_X/ (\mathcal{L}^{-n}_1, \cdots, \mathcal{L}^{-n}_t).$ 

We let  $Z =  V( \mathcal{L}_1^{-1}, \cdots, \mathcal{L}_i^{-1})$ and let $\mathcal{F}$ be a locally free coherent sheaf. We find that: 
$\mathcal{E}xt^j ( \ocal_{X}/(\mathcal{L}_1^{-n}, \cdots, \mathcal{L}_i^{-n}), \mathcal{F}) =0$ unless $j=i$, and 
$$\mathcal{E}xt^i ( \ocal_{X}/(\mathcal{L}_1^{-n}, \cdots,
\mathcal{L}_i^{-n}), \mathcal{F})  = \mathrm{Coker} ( \oplus_{t}
s_t^{n}: \bigotimes_{t' \neq t} \mathcal{L}_{t'}^{n}  \otimes
\mathcal{F} \rightarrow    (\bigotimes_t \mathcal{L}_t^{n}) \otimes
\mathcal{F}).$$Taking the direct limit over~$n$ gives $\mathrm{R}^i\Gamma_{Z} (\mathcal{F})$.
\end{example} 
\subsubsection{The Cousin complex of the Ekedahl--Oort stratification}\label{sect-cousin-compEO}

We now assume that $K_p=\prod_v\GSp_4(\ZZ_p)$, and let $X = X_{K, \Sigma}$ and denote by $X^\star$ the minimal compactification. We have a morphism $f: X \rightarrow X^\star$. 
We fix an integer $n$ and work over $X_n = X \times \Spec \ZZ/p^n\ZZ$
and $X^\star_n = X^\star \times \Spec \ZZ/p^n\ZZ$, and let~$Y_n$
denote the interior of~$X_n$. 
We consider the filtration $Z^\star$ on $X^\star_n $ given by taking
$Z^\star_i $  to be the closure of all Ekedahl--Oort strata of
codimension $i$. Here are some known facts (see~\cite[Thm.\ 6.2.3]{2015arXiv150705922B}): \begin{enumerate}
\item $Z^\star$ is a filtration.
\item $Z^\star_i \setminus Z^\star_{i+1}$ is affine.
\item $Z^\star_i \setminus Z^\star_{i+1}$ is a set-theoretic local complete
  intersection in $U^\star_{i+1} = X^\star_n \setminus
  Z^\star_{i+1}$. \end{enumerate} 

 We now consider the pull-back  $Z$ of $Z^\star$ on  $X_n$. We deduce that:
\begin{enumerate}
\item $Z$ is a filtration.
\item $Z_i \setminus Z_{i+1} \hookrightarrow X_n$ is affine.
\item $Z_i \setminus Z_{i+1}$ is a set-theoretic  local complete intersection in $U_{i+1} = X_n \setminus  Z_{i+1}$.
\end{enumerate} 
 
 \begin{prop}\label{prop-lastone}  

 The cohomology $\mathrm{R}\Gamma (X_n,  \omega^\kappa(-D))$ is computed by  $$ \Gamma( X_n, \Cous_{Z} (\omega^\kappa(-D))).$$ 
\end{prop} 
 
 \begin{proof}  It follows from  Theorem~ \ref{thmkempf} that $\omega^\kappa(-D) \rightarrow \Cous_{Z} (\omega^\kappa(-D))$ is a quasi-isomorphism. It suffices to prove that $\Cous_{Z} (\omega^\kappa(-D))$ is a complex of acyclic sheaves.  
By~\eqref{eqn: term in Cousin as
   pushforward of local cohomology}, the sheaf in degree $i$ is equal to $(j_{i+1})_* \mathrm{R}^i\Gamma_{Z_i \setminus Z_{i+1}}(\omega^\kappa(-D)\vert_{U_{i+1}})$. This sheaf is supported on $Z_i \setminus Z_{i+1}$. We claim that $Rf_\star (j_{i+1})_* \mathrm{R}^i\Gamma_{Z_i \setminus Z_{i+1}}(\omega^\kappa(-D)\vert_{U_{i+1}})$ is concentrated in degree $0$. Since  $f_\star (j_{i+1})_* \mathrm{R}^i\Gamma_{Z_i \setminus Z_{i+1}}(\omega^\kappa(-D)\vert_{U_{i+1}})$ is an acyclic sheaf because it is supported on $Z^\star_i \setminus Z^\star_{i+1}$, the proposition will follow from our claim.

Let us prove the claim. By construction, $Z_i \setminus Z_{i+1}$ is a
finite disjoint union of Ekedahl--Oort strata and $(j_{i+1})_*
\mathrm{R}^i\Gamma_{Z_i \setminus
  Z_{i+1}}(\omega^\kappa(-D)\vert_{U_{i+1}})$ is  a finite direct sum
indexed by these Ekedahl--Oort strata.  Let $E$ be an  Ekedahl--Oort
stratum appearing in  $Z_i \setminus Z_{i+1}$. It    can be written as
the intersection  of $i$ Cartier divisors in $U_{i+1}$, using the
theory of generalized Hasse invariants (one can also assume that these
Cartier divisors are pulled back from $X^\star_n$; note that a
sufficiently large power of each generalized Hasse invariant can be
lifted to~$X^\star_n$). Let us denote these Cartier divisors by $ \ocal_X \stackrel{s_t}\rightarrow \mathcal{I}_t$.  It follows from Example \ref{ex-computation}  that the direct summand of the sheaf  $\mathrm{R}^i\Gamma_{Z_i \setminus Z_{i+1}}(\omega^\kappa(-D)\vert_{U_{i+1}})$ corresponding to $E$ is the  inductive limit  of the sheaves: 
$$\mathcal{H}^i \mathcal{H}om(\Kos(s_1^n,\cdots, s_i^n), \omega^\kappa(-D)).$$
The complex $\mathcal{H}om(\Kos(s_1^n,\cdots, s_i^n),
\omega^\kappa(-D))$ is a complex of sheaves acyclic relatively to the
minimal compactification, and concentrated in degree $i$. \end{proof}

\begin{lem}\label{lem: injection of Cousin complexes} There is an
  injection of complexes $$ \Gamma( X_n, \Cous_{Z} (\omega^\kappa(-D)))
  \hookrightarrow  \Gamma( Y_n, \Cous_{Z}
  (\omega^\kappa(-D))).$$
\end{lem}

\begin{proof} This follows directly from the description of the objects  of the complex $\Cous_{Z} (\omega^\kappa(-D))$ given in the course of the preceding proof.
\end{proof}

It remains to prove that our Hecke operators act on $ \Gamma( X_n, \Cous_{Z} (\omega^\kappa(-D)))$ and $\Gamma( Y_n, \Cous_{Z} (\omega^\kappa(-D))).$ Let $g \in G( \mathbb{A}_f^p)$. We consider the correspondence: 
 \begin{eqnarray*}     
 \xymatrix{ &X_{K \cap g K g^{-1} }   \ar[rd]^{p_1}\ar[ld]_{p_2}&  \\
  X_K & & X_K}
  \end{eqnarray*}
  and more precisely its reduction modulo $p^n$.    We have a
  cohomological correspondence $p_2^\star \omega^\kappa(-D)
  \rightarrow p_1^! \omega^\kappa(-D)$, as defined in \S\ref{sect-hecke-away-from-p}.
  
  \begin{lem}\label{lem: cohomological correspondence of complexes} This cohomological correspondence  induces  a cohomological correspondence  of complexes compatible with the augmentation: 
  $$ p_2^\star \Cous_Z(\omega^\kappa(-D)) \rightarrow p_1^! \Cous_Z(\omega^\kappa(-D))$$
  \end{lem}
  
  \begin{rem} In the above correspondence, the functors $p_2^\star$ and $p_1^!$ are applied to each object of the complex. Moreover, for each object $\Cous_Z(\omega^\kappa(-D))^i$ of $\Cous_Z(\omega^\kappa(-D))$, $p_1^! \Cous_Z(\omega^\kappa(-D))^i$ is a sheaf (i.e it is concentrated in degree $0$). 
  \end{rem}
  \begin{proof}[Proof of Lemma~\ref{lem: cohomological correspondence of complexes}] For each index $i$, we have  $$\Cous_Z(\omega^\kappa(-D))^i = (j_{i+1})_* \mathrm{R}^i\Gamma_{Z_i \setminus Z_{i+1}}(\omega^\kappa(-D)\vert_{U_{i+1}}).$$
  
  We choose (by considering powers of generalized Hasse invariants) an increasing sequence  $(Z_i \setminus Z_{i+1})_k$ of  subschemes $U_{i+1}$ with support $Z_i \setminus Z_{i+1}$, which are local complete intersections, and are cofinal among all subschemes of $U_{i+1}$ with support $Z_i \setminus Z_{i+1}$.  
  
 We have  $\mathrm{R}^i\Gamma_{Z_i \setminus Z_{i+1}}(\omega^\kappa(-D)\vert_{U_{i+1}}) = \colim_k \mathcal{E}xt^i( \ocal_{(Z_i \setminus Z_{i+1})_k}, \omega^\kappa(-D)\vert_{U_{i+1}})$. 
 Also recall that $\mathcal{E}xt( \ocal_{(Z_i \setminus Z_{i+1})_k}, \omega^\kappa(-D)\vert_{U_{i+1}}) = \mathcal{E}xt^i( \ocal_{(Z_i \setminus Z_{i+1})_k}, \omega^\kappa(-D)\vert_{U_{i+1}})[- i]$. 
 
 We have $$p_1^! (\mathcal{E}xt( \ocal_{(Z_i \setminus Z_{i+1})_k},
 \omega^\kappa(-D)\vert_{U_{i+1}}) = \mathcal{E}xt( p_1^\star
 \ocal_{(Z_i \setminus Z_{i+1})_k}, p_1^!
 \omega^\kappa(-D)\vert_{U_{i+1}})$$ by~\cite[Prop.\
 III.8.8]{Hartshorne}. One checks that $p_1^\star \ocal_{(Z_i
   \setminus Z_{i+1})_k} = ~^{\mathbf{L}}p_1^\star \ocal_{(Z_i
   \setminus Z_{i+1})_k}$ because the pull back of a local regular
 sequence defining $(Z_i \setminus Z_{i+1})_k$ is again a local
 regular sequence; we will not comment on the vanishing of higher pullbacks in the rest of this argument.  We deduce  that 
  $$p_1^! (\mathcal{E}xt^i( \ocal_{(Z_i \setminus Z_{i+1})_k}, \omega^\kappa(-D)\vert_{U_{i+1}}) = \mathcal{E}xt^i( p_1^\star \ocal_{(Z_i \setminus Z_{i+1})_k}, p_1^! \omega^\kappa(-D)\vert_{U_{i+1}}).$$
 On the other hand there is by adjunction a  map: 
 $$ p_2^\star \mathcal{E}xt^i( \ocal_{(Z_i \setminus Z_{i+1})_k}, \omega^\kappa(-D)\vert_{U_{i+1}}) \rightarrow \mathcal{E}xt^i( p_2^\star \ocal_{(Z_i \setminus Z_{i+1})_k}, p_2^\star \omega^\kappa(-D)\vert_{U_{i+1}}).$$
 Since the Ekedahl--Oort stratification is invariant under prime to $p$ isogenies, we deduce that $p_2^{-1}( Z_i\setminus Z_{i+1}) = p_1^{-1}(Z_i \setminus Z_{i+1})$. Therefore, for each $k$, for all large enough $t$, there is   a natural map $p_1^\star  \ocal_{(Z_i\setminus Z_{i+1})_{t}} \rightarrow 
 p_2^\star  \ocal_{(Z_i\setminus Z_{i+1})_{k}}$.
 
  We therefore get a map:
 $$ p_2^\star \mathcal{E}xt^i( \ocal_{(Z_i \setminus Z_{i+1})_k}, \omega^\kappa(-D)\vert_{U_{i+1}}) \rightarrow \mathcal{E}xt^i( p_2^\star \ocal_{(Z_i \setminus Z_{i+1})_k}, p_2^\star \omega^\kappa(-D)\vert_{U_{i+1}}) $$ $$ \rightarrow \mathcal{E}xt^i( p_1^\star \ocal_{(Z_i \setminus Z_{i+1})_t}, p_1^! \omega^\kappa(-D)\vert_{U_{i+1}}) = p_1^! (\mathcal{E}xt^i( \ocal_{(Z_i \setminus Z_{i+1})_t}, \omega^\kappa(-D)\vert_{U_{i+1}}).$$
 
  Passing to the inductive limit over $k$ and $t$ yields the cohomological correspondence: 
  $$p_2^\star  \mathrm{R}^i\Gamma_{Z_i \setminus Z_{i+1}}(\omega^\kappa(-D)\vert_{U_{i+1}}) \rightarrow p_1^! \mathrm{R}^i\Gamma_{Z_i \setminus Z_{i+1}}(\omega^\kappa(-D)\vert_{U_{i+1}}).$$
  
 Moreover this construction is  canonical and  is compatible with all differentials in the Cousin complex and with the augmentation.   \end{proof}
  
  \begin{prop}\label{prop-cousin-heckeaction}The Hecke algebra $\mathcal{H}^p$ acts on $\mathrm{R}\Gamma(X^{G_1}, \omega^\kappa(-D))$ and also on $\mathrm{R}\Gamma(X^{G_1}, \omega^\kappa)$.
  \end{prop}
  \begin{proof} By Serre duality, it suffices to treat the case of $\mathrm{R}\Gamma(X^{G_1}, \omega^\kappa(-D))$. 
  The cohomology $\mathrm{R}\Gamma(X, \omega^\kappa(-D))$ is represented by 
   $$ \varprojlim_n \Gamma( X_n, \Cous_{Z} (\omega^\kappa(-D)))$$  and this complex injects into $$ \varprojlim_n \Gamma( Y_n, \Cous_{Z} (\omega^\kappa(-D))).$$
   The double class $\mathbb{1}_{K^pgK^p}$ acts everywhere. We can
   restrict to the ``$G_1$'' direct factor (the Ekedahl--Oort
   stratification is preserved by the action of
   $(\ocal_F)_{(p)}^{\times, +}$) and  the compatibility with the product in the Hecke algebra is known on 
   $$ \varprojlim_n \Gamma( Y^{G_1}_n, \Cous_{Z} (\omega^\kappa(-D))).$$
    Therefore it holds everywhere. 
   \end{proof}
   \begin{rem}
     \label{rem: mixing the compacts}It follows by an identical
     argument to the proof of
     Proposition~\ref{prop-cousin-heckeaction} that if~$K^p_1$,
     $K^p_2$, $K^p_3$ are three choices of tame level, and~$\Sigma_1$,
     $\Sigma_2$, $\Sigma_3$ are suitable choices of polyhedral cone
     decompositions, then the
     composite of the Hecke operators \[ [{K_1^pg_1K_2^p}]: \mathrm{R}\Gamma ( X^{G_1}_{K^p_2K_p, \Sigma_2},
\omega^{\kappa})\to \mathrm{R}\Gamma ( X^{G_1}_{K^p_1K_p, \Sigma_1},
\omega^{\kappa})\] and
     \[ [{K_2^pg_2K_3^p}]: \mathrm{R}\Gamma ( X^{G_1}_{K^p_3K_p, \Sigma_3},
\omega^{\kappa})\to \mathrm{R}\Gamma ( X^{G_1}_{K^p_2K_p, \Sigma_2},
\omega^{\kappa})\] is the Hecke operator \[  [{K_1^pg_1K_2^p}] [{K_2^pg_2K_3^p}]:\mathrm{R}\Gamma ( X^{G_1}_{K^p_3K_p, \Sigma_3},
\omega^{\kappa})\to \mathrm{R}\Gamma ( X^{G_1}_{K^p_1K_p, \Sigma_1},
\omega^{\kappa}).\]
   \end{rem}

\subsubsection{Hecke operators at $p$: Siegel type
  operator}\label{sect-def-Siegel-Hecke} We assume here that $K = K^p
K_p$ and $K_p = G_1(\ZZ_p)$.  Let us fix a place $w$ above $p$. We are
going to define an action of a Hecke operator $T_{w,1}$ on $\mathrm{R}\Gamma ( X_{K, \Sigma}, \omega^{\kappa})$, $\mathrm{R}\Gamma ( X^{G_1}_{K, \Sigma}, \omega^{\kappa})$, and their cuspidal versions. The action on $\mathrm{R}\Gamma ( X_{K, \Sigma}, \omega^{\kappa})$ and $\mathrm{R}\Gamma ( X_{K, \Sigma}, \omega^{\kappa}(-D))$ is not canonical, and depends on the choice of $x_w$ made in~\S\ref{subsect-functor-3}, but the action on $\mathrm{R}\Gamma ( X^{G_1}_{K, \Sigma}, \omega^{\kappa})$ and $\mathrm{R}\Gamma(X^{G_1}_{K, \Sigma}, \omega^{\kappa}(-D))$ is canonical. 
 
 Set $K' =  K^p K'_p$ where $K'_p = \prod_{v \neq w } \mathrm{GSp}_4(
 \ocal_{F_v}) \times \Si(w)$. In~\S\ref{subsect-functor-4} we defined maps $p_1,p_2:X_{K'}\to
 X_K$ giving a Hecke correspondence:  
  \begin{eqnarray*}     
 \xymatrix{ &X_{K' }   \ar[rd]^{p_1}\ar[ld]_{p_2}&  \\
  X_K & & X_K}
  \end{eqnarray*} 
 The key geometric properties of this correspondence are (see Proposition~\ref{prop: local models results for G}):
 \begin{enumerate}
 \item $X_{K'}$ and $X_K$ are relative complete intersections over
   $\Spec \ZZ_p$, and are pure of the same dimension,
 \item $X_{K}$ is smooth over $\Spec \ZZ_p$,
 \item $X_{K'}$ is smooth over $\Spec \ZZ_p$ up to codimension $2$ and
   normal.\end{enumerate}
 In particular, we are in the lci situation in the sense of~\S\ref{subsubsec: fundamental class}, so  we have an invertible dualizing sheaf $p_1^! \ocal_{X_K}$ and a fundamental class $p_1^{*} \ocal_{X_K} \rightarrow p_1^! \ocal_{X_K}$.  Moreover, for all weights $\kappa = (k_v, l_v)$ with $l_v \geq 0$, we have a natural map $p_2^{*} \omega^\kappa \rightarrow p_1^{*} \omega^\kappa$ provided by the differential of the isogeny  $ p_1^{*} \cG \rightarrow p_2^{*}  \cG $ on $X_{K'}$. 
 
 Composing these maps, we obtain a cohomological correspondence $ \Theta: p_2^{*}
 \omega^\kappa \rightarrow p_1^! \omega^\kappa$.

 \begin{lem}\label{lem-normalization1} When $l_w \geq 2$,  this map is divisible by $p^3$. 
 \end{lem}
 \begin{proof} We need to prove that  $\Theta$ factors through
   $p^3 p_1^! \omega^\kappa$.  As $X_{K'}$ is normal and the source and target are locally free sheaves, it is enough to establish this factorization in codimension one.
   As this factorization is furthermore trivial over the
   generic fibre of~$X_{K'}$, it is enough to prove it over the completed local rings of the generic points
   of the special fibre of $X_{K'}$.

   There are three types of generic points in the special fibre
   classified by the multiplicative rank $r=0,1$ or $2$ of the isogeny
   $p_1^{*} \cG \rightarrow p_2^{*} \cG $. In each case one calculates separately the $p$-divisibility of the map
   $p_2^{*} \omega^\kappa \rightarrow p_1^{*}
   \omega^\kappa$ and of the fundamental class as in the proof of~\cite[Lem.\
   7.1.1]{pilloniHidacomplexes}.  One finds that the fundamental class $p_1^{*} \ocal_{X_{K}} \rightarrow p_1^! \ocal_{X_{K}}$ is divisible by $p^3$ when $r=0$, $p$ when $r=1$, and $p^0$ when $r=2$, and the map $p_2^{*} \omega^\kappa \rightarrow p_1^{*}
   \omega^\kappa$ is divisible by $p^0$ when $r=0$, $p^{l_w}$ when $r=1$, and $p^{k_w+l_w}$ when $r=2$.  The result follows as $3,l_w+1,k_w+l_w\geq 3$.
\end{proof}  
 
 We can thus consider the normalized cohomological correspondence  $T_{w,1}: p^{-3} \Theta: p_2^{*} \omega^\kappa \rightarrow p_1^! \omega^\kappa$, and we obtain a Hecke operator:  $$ T_{w,1}: \mathrm{R}\Gamma(X_K, \omega^\kappa) \rightarrow \mathrm{R}\Gamma( X_{K'}, p_2^{*} \omega^\kappa) \stackrel{p^{-3}\Theta}\rightarrow \mathrm{R}\Gamma( X_{K '}, p_1^! \omega^\kappa) \stackrel{\mathrm{Tr}}\rightarrow \mathrm{R}\Gamma( X_{K }, \omega^\kappa).$$
A similar definition applies to cuspidal cohomology and works over
$X_K^{G_1}$.

\begin{rem}\label{rem-hecke-mea-culpa}
One readily checks that the Hecke correspondence used to define $T_{w,1}$ corresponds to the double coset $[\GSp_4(\O_{F_w})\diag(1,1,p^{-1},p^{-1}) \GSp_4(\O_{F_w})]$ (see \cite[Rem. 5.6]{F-Pilloni} for instance) which differs by an element of the centre from the spherical Hecke operator considered in \S\ref{sssection spherical hecke}.  We justify this discrepancy as follows: when doing geometry and working with the moduli interpretation, we prefer to use this Hecke operator, while when doing local representation theory and considering Galois representations, we prefer to use the Hecke operators considered in \S\ref{subsec: local representation theory}.  In this paper we will systematically work on spaces with fixed central character so that the (normalized) action of these two Hecke operators are the same.  The same remark will apply to all the Hecke operators at $p$ considered in this paper.  We hope this will not cause any confusion.
\end{rem}

\subsubsection{Hecke operators at $p$: Klingen type
  operator}\label{sect-def-Klingen-Hecke} We again assume  that $K = K^p K_p$ and $K_p = G_1(\ZZ_p)$.  Let us fix a place $w$ above $p$. We are going to define an action of a  Hecke operator $T_{w}$ on $\mathrm{R}\Gamma ( X_{K, \Sigma}, \omega^{\kappa})$, $\mathrm{R}\Gamma ( X^{G_1}_{K, \Sigma}, \omega^{\kappa})$ and their cuspidal versions. As before, the action on $\mathrm{R}\Gamma ( X_{K, \Sigma}, \omega^{\kappa})$ and $\mathrm{R}\Gamma ( X_{K, \Sigma}, \omega^{\kappa}(-D))$   is not canonical and depends on the choice of $x_w$ made in~\S\ref{subsect-functor-3}, but the action on $\mathrm{R}\Gamma ( X^{G_1}_{K, \Sigma}, \omega^{\kappa})$ and $\mathrm{R}\Gamma ( X^{G_1}_{K, \Sigma}, \omega^{\kappa}(-D))$ is canonical. 

\begin{rem} The Hecke operator that we define in this section does not correspond to the double coset operator $[\GSp_4(\O_{F_w})\diag(1,p^{-1},p^{-1},p^{-2}) \GSp_4(\O_{F_w})]$ but rather some
variant of it that we call $T_{w}$.  The formula for $T_w$ in terms of double cosets is
\begin{align*}
T_w&=[\GSp_4(\O_{F_w})\diag(1,p^{-1},p^{-1},p^{-1})\Par(w)][\Par(w)\diag(1,1,1,p^{-1})\GSp_4(\O_{F_w})]\\
&=p[\GSp_4(\cO_{F_w})\diag(1,p^{-1},p^{-1},p^{-2})\GSp_4(\cO_{F_w})]\\
&\qquad+(1+p+p^2+p^3)p^{-2}[\GSp_4(\cO_{F_w}) \diag(p^{-1},p^{-1},p^{-1},p^{-1})\GSp_4(\cO_{F_w})].
\end{align*}
\end{rem}

 Set $K' =  K^p K'_p$ where $K'_p = \prod_{v \neq w } \mathrm{GSp}_4(
 \ocal_{F_v}) \times \Kli(w)$ and $K'' = K^p
 K''_p$ where $K''_p = \prod_{v \neq w } \mathrm{GSp}_4( \ocal_{F_v})
 \times \Par ( w)$. In~\S\ref{subsect-functor-3}, we  defined morphisms $p_1:X_{K'}\to X_K$,
 $p_2:X_{K'}\to X_{K''}$ giving a Hecke correspondence: 
 
  \begin{eqnarray*}     
 \xymatrix{ &X_{K' }   \ar[rd]^{p_1}\ar[ld]_{p_2}&  \\
  X_{K''} & & X_K}
  \end{eqnarray*}
 
 The key geometric  properties are again (see Proposition~\ref{prop: local models results for G}): 
 \begin{enumerate}
 \item $X_{K'}$, $X_{K''}$ and  $X_K$ are relative complete
   intersections over $\Spec \ZZ_p$ of the same (pure) dimension,
 \item $X_{K}$ is smooth over $\Spec \ZZ_p$,
 \item $X_{K'}$ and $X_{K''}$  are smooth over $\Spec \ZZ_p$ up to
   codimension $2$ and normal.\end{enumerate} 
   We are again in the lci  situation, so we have invertible dualizing sheaves $p_1^! \ocal_{X_K}$ and $p_2^! \ocal_{X_{K''}}$ and  fundamental classes $p_1^{*} \ocal_{X_K} \rightarrow p_1^! \ocal_{X_K}$ and $p_2^{*} \ocal_{X_{K''}} \rightarrow p_2^! \ocal_{X_{K''}}$.
 
 For all weights $\kappa = (k_v, l_v)$ with $l_v \geq 0$, we have
 natural maps $p_2^{*} \omega^\kappa \rightarrow p_1^{*}
 \omega^\kappa$ provided by the differential of the isogeny  $
 p_1^{*} \cG \rightarrow p_2^{*}  \cG $  on $X_{K'}$, and $p_1^{*} \omega^\kappa \rightarrow p_2^{*} \omega^\kappa$ provided by the differential of the isogeny  $ p_2^{*} \cG \rightarrow p_1^{*}  \cG $. We therefore obtain two cohomological correspondences $ \Theta_1: p_2^{*} \omega^\kappa \rightarrow p_1^! \omega^\kappa$ and $\Theta_2: p_1^{*} \omega^\kappa \rightarrow p_2^! \omega^\kappa$.
 
 \begin{lem}\label{lem-normalization2} When $l_w \geq 2$,  the  map $\Theta_1$ is divisible by $p^{2+l_w}$ and the map $\Theta_2$ is divisible by $p$. 
 \end{lem}
 \begin{proof}This can be proved in exactly the same way as
   Lemma~\ref{lem-normalization1}, by an explicit  check over the
   completed local rings of generic points of the special fibre of
   $X_{K'}$. The details may be found in the proofs of~\cite[Lem.\
   7.1.1, 7.1.2]{pilloniHidacomplexes}.
 \end{proof} 
 
 We can therefore consider the normalized fundamental classes
 $T'_w = p^{-2-l_w} \Theta_1: p_2^{*} \omega^\kappa \rightarrow p_1^!
 \omega^\kappa$ and $T''_w = p^{-1} \Theta_2: p_1^{*} \omega^\kappa
 \rightarrow p_2^! \omega^\kappa$, and we obtain   Hecke operators: 
  
  $$ T'_{w}: \mathrm{R}\Gamma(X_{K''}, \omega^\kappa) \rightarrow \mathrm{R}\Gamma(X_K, \omega^\kappa)$$
  and 
  
  $$ T''_{w}:  \mathrm{R}\Gamma(X_{K}, \omega^\kappa) \rightarrow \mathrm{R}\Gamma(X_{K''}, \omega^\kappa). $$ 
We set $T_{w}:= T'_{w} \circ T''_{w}$. Similar definitions apply to cuspidal cohomology and work over
$X_K^{G_1}$. \begin{rem}
\label{rem: independence of compactification}Just as the complexes
that we are considering are independent of the choice of
compactification by Lemma~\ref{lem: cohomologies independent of
  compactification}, so too are the actions of~$T_{w,1}$ and~$T_w$ on
them. See~\cite[Prop.\
7.2.1]{pilloniHidacomplexes} for the case of~$T_w$; the argument
for~$T_{w,1}$ is similar, but easier, and is left to the interested reader.
\end{rem}

\subsection{Cohomology and automorphic representations}Let $K=\prod_vK_v\subset \GSp_4(\A_F^\infty)=G_1(\A_\Q^\infty)$ be an open compact subgroup,
let~$S\supset S_p$ be a finite set of places such that~$K_v=\GSp_4(\cO_{F_v})$
for $v\notin S$, and let
$$\TTt=\bigotimes_{v\not\in S}\cO[\GSp_4(F_v) \doubleslash \GSp_4(\cO_{F_v})]$$
 be the ring of
spherical Hecke operators away from~$S$.  We say that a maximal ideal
$\m\subset\TTt$ is \emph{non-Eisenstein} if the residue field~$\TTt/\m$ is a
finite extension of~$\Fp$, and for (any) inclusion~$\TTt/\m \rightarrow \Fbar_p$ there exists
 an irreducible representation ~$\rhobar:G_F\to\GSp_4(\Fbar_p)$
with the property that, for each~$v\notin S$, we
have~$\det(X-\rhobar(\Frob_v))$
 \[\equiv X^4-T_{v,1}X^3+(q_vT_{v,2}+(q_v^3+q_v)T_{v,0})X^2-q_v^3T_{v,0}T_{v,1}X+q_v^6T_{v,0}^2\pmod{\m}.\] 
 (cf.\ (\ref{eqn: char poly
  for unramified Hecke})). 

Our main aim in this section is to prove the following result.
\begin{thm}\label{thm: cohomology in terms of automorphic forms}
  Let $\kappa=(k_v,l_v)_{v|\infty}$ with $k_v\geq l_v\geq 2$ and
$k_v\equiv l_v\pmod 2$ be a weight and let $\m$ be  non-Eisenstein.
\begin{enumerate}
\item   For~$i=0, 1$, there is an $\overline{E}[\GSp_4(\A_F^\infty) \doubleslash K]$-equivariant inclusion \numequation\label{eqn: description of cohomology in terms of
    automorphic forms} \bigoplus_\pi
  (\pi^\infty)_\m^K\otimes\barE\subseteq
  H^i(X_K^{G_1},\omega^\kappa(-D))_\m^{|\cdot|^2}\otimes\barE
\end{equation}
where, on the right hand side, the superscript~$|\cdot|^2$ indicates the
space on which the diamond operators at places~$v\not\in S$ act
via~$|\cdot|^2$; and on the left hand side, $\pi$ runs over the cuspidal automorphic representations of
$\GSp_4(\A_F)$ with weight $\kappa$ and central character~$|\cdot|^2$
such that
\begin{itemize}
\item $\pi_v$ is holomorphic for those $v|\infty$ for which $l_v>2$,
  and
\item $\#\{v|\infty\mid \pi_v\text{ is not holomorphic}\}=i$.
\end{itemize}

\item There is an absolute constant~$R$ such that if for each $v|\infty$
\begin{itemize}
\item $k_v-l_v>R$, and
\item either $l_v=2$ or $l_v>R$,
\end{itemize}
 then the inclusion \emph{(\ref{eqn:
    description of cohomology in terms of automorphic forms})} is an
equality.
\item If~$i=0$, then~\emph{(\ref{eqn: description of cohomology in terms of
    automorphic forms})} is an
equality. In fact, a version of this statement holds without having to localize
at a non-Eisenstein maximal ideal; there is an $\overline{E}[\GSp_4(\A_F^\infty) \doubleslash K]$-equivariant 
isomorphism \numequation\label{eqn: description of cohomology in terms
  of automorphic forms i equals 0 incl Eisenstein}
H^0(X_K^{G_1},\omega^\kappa(-D))^{|\cdot|^2}\otimes\barE=\bigoplus_\pi
(\pi^\infty)^K\end{equation}where $\pi$ runs over the cuspidal automorphic
representations of $\GSp_4(\A_F)$ with weight $\kappa$ and central
character~$|\cdot|^2$ which are holomorphic at all infinite places.\end{enumerate}
  \end{thm}\begin{rem}Theorem~\ref{thm: cohomology in terms of automorphic forms}
  is by no means optimal; the same results should hold for any cohomological degree $i$, and with a much weaker regularity assumption on $\kappa$ in part (2). However,
  it seems difficult to deduce results in this
  generality from the literature, so we have restricted ourselves to
  this result, for which we only need to consider the cohomology of
  the boundary in degree~$0$. We explain the proof below, after
  proving a corollary and a preparatory lemma.
  \end{rem}
Theorem~\ref{thm: cohomology in terms of automorphic forms} has the
following useful corollary.
\begin{cor}
  \label{cor: dimensions of cohomology in various degrees}Suppose
  that we are in the setting of Theorem~\ref{thm: cohomology in terms
    of automorphic forms} and the hypothesis on $\kappa$ in (2) holds. Let~$l_0$
  denote the number of infinite places~$v$ with~$l_v=2$. Then
\begin{equation*}
\dim_{\barE}H^1(X_K^{G_1},\omega^\kappa(-D))_\m^{|\cdot|^2}\otimes\barE=l_0 \dim_{\barE}H^0(X_K^{G_1},\omega^\kappa(-D))_\m^{|\cdot|^2}\otimes\barE.
\end{equation*}

\end{cor}
\begin{proof}
  Since~$\m$ is non-Eisenstein, the automorphic representations~$\pi$
  which contribute to~(\ref{eqn: description of cohomology in terms of
  automorphic forms}) are all of general type in the sense
of~\cite{MR2058604} by Lemma~\ref{lem: getting to general type}. There are~$l_0$ ways to choose an
  infinite place~$v$ with~$l_v=2$, and we let~$\pi_v$ be generic for
  this place and holomorphic at the other infinite places. The
  result then follows from Theorem~\ref{thm: arthur classification
    results}.
\end{proof}
\begin{rem}
  \label{rem: harris zucker did it already}The following lemma is
  essentially a special case of the much more general results proved
  in~\cite{HZIII}, and can presumably be proved using the techniques
  of that paper, but since our Shimura varieties do not satisfy the
  precise assumptions needed to cite the results of~\cite{HZIII}, we
  have chosen to give a direct proof.
\end{rem}
\begin{lemma}
  \label{lem:harris-zucker} 
Let $\kappa=(k_v,l_v)_{v|\infty}$ be a weight, with $k_v\geq l_v\geq 2$ and
$k_v\equiv l_v\pmod 2$, and let~$\m$ be a
non-Eisenstein maximal ideal. Let $D$ denote the boundary of
$X_K^{G_1}$.  Then  $\HH^0(D, \omega^\kappa)_{\m} \otimes
\Ebar=0$. \end{lemma}
\begin{proof}In the case that~$F=\Q$ this follows from~\cite[IV, Satz
  4.4]{MR871067}, as in the proof of~\cite[Cor.\
  15.2.3.1]{pilloniHidacomplexes}, so we can and do assume
  that~$F\ne \Q$ in what follows. We let
  $ \pi: X_{K}^{G_1} \rightarrow X_K^{\star, G_1}$ be the map between
  toroidal and minimal compactifications. We let $\partial X^\star\subset X_K^{\star,G_1}$ be the (reduced) boundary of the minimal compactification, which we can write
  as~$\partial_0 X^\star\coprod \partial_1 X^\star$,
  where~$\partial_1X^\star$ is a union of Hilbert modular varieties
  for the group $\mathrm{Res}_{F/\Q} \mathrm{GL}_{2}$, and the
  complement $\partial_0 X^\star$ is a finite union of points.

  Suppose firstly that we are in the case that~$k_v=l_v=k$ for
  some~$k$ independent of~$v$. Then~$\omega^\kappa$ is pulled back
  from the minimal compactification, and
  since~$\pi_*(\omega^\kappa|_{D})=\omega^\kappa|_{\partial X^*}$, we
  have
  $\HH^0(D, \omega^\kappa) = \HH^0(\partial X^\star,
  \omega^\kappa)$. (To see
  that we have an identification~$\pi_*(\omega^\kappa|_D)=\omega^\kappa|_{\partial X^*}$, it
  suffices by the projection formula to show that
  $\pi_*(\cO_{D})=\cO_{\partial X^*}$. This follows from the
  facts that $\pi_*\cO_{X_K^{G_1}}=\cO_{X_{K}^{*,G_1}}$,
  $\pi_*\cI_{D}=\cI_{\partial X^*}$, and
  $R^1\pi_*\cI_{D}=0$.)

Suppose now that we are not in the case that~$k_v$ and~$l_v$ are equal and
independent of~$v$. Then it follows from the results
of~\cite{MR3186092}, see~\cite[Prop.\ 1.5.8]{BrascaRosso}, that any element of~$\HH^0(D,
\omega^\kappa)$ vanishes on~$\pi^{-1}(\partial_0 X^*)$; so  the map $\HH^0(D,
\omega^\kappa)\to \HH^0(\pi^{-1}(\partial_1 X^*),
\omega^\kappa)$ is injective, and it suffices to show that~$\HH^0(\pi^{-1}(\partial_1 X^*),
\omega^\kappa)$ is Eisenstein. Again by~\cite[Prop.\ 1.5.8]{BrascaRosso} it
follows that   $\pi_\star (\omega^\kappa
\vert_{\pi^{-1}(\partial_1 X^\star)})$ is zero if the~$l_v$ are not
all equal, and otherwise is equal to the sheaf $\omega^{(k_v)_{v|p}}$
(where we are using the usual labelling of weights for sheaves on
Hilbert--Blumenthal modular schemes).

In either case, we have seen that the space that we are considering
either vanishes, or injects into~$H^0(\partial_1X^*,\omega^{(k_v)})$.
Now it is convenient to work adelically. Let us fix $W \in \mathfrak{C}$ with $\dim_F W =1$.
 Then $\partial_1 X^\star$ is as follows (where~$P_W$ and~$M_W$
 are defined in~\S\ref{sec: compactifications}):
$$ P_{W}^+(\Q) \backslash \mathcal{H}^{[F:\Q]}_1 \times G_1(\mathbb{A}^\infty)/ K  $$ $$ =P_{W}^+(\Q) \backslash \mathcal{H}^{[F:\Q]}_1 \times P_W(\mathbb{A}^\infty) \times^{P_{W}(\mathbb{A}^\infty)} G_1(\mathbb{A}^\infty)/ K $$ $$ =M_{W}^+(\Q) \backslash \mathcal{H}^{[F:\Q]}_1 \times M_W(\mathbb{A}^\infty) \times^{P_{W}(\mathbb{A}^\infty)} G_1(\mathbb{A}^\infty)/ K. $$

We therefore find that 

$$\HH^0(\partial_1 X^\star, \pi_\star \omega^{(k_v)}) = \big( \mathrm{Ind}^{G_1(\mathbb{A}^\infty)}_{P_W(\mathbb{A}^\infty)} ( M) \big)^K$$where $M = \lim_{K \subset M_W(\mathbb{A}^\infty)} \HH^0( M_{W}^+(\Q) \backslash \mathcal{H}^{[F:\Q]}_1 \times M_W(\mathbb{A}^\infty)/K, \omega^{(k_v)})$ 
from which it follows that the eigensystems arising from
$\HH^0(\partial_1 X^\star, \omega^{(k_v)})$ are Eisenstein. Indeed,
since we are assuming that~$F\ne\Q$, it follows from Koecher's
principle that the cohomology groups $\HH^0( M_{W}^+(\Q) \backslash
\mathcal{H}^{[F:\Q]}_1 \times M_W(\mathbb{A}^\infty)/K,
\omega^{(k_v)})$ are spaces of Hilbert modular forms, and thus
have associated two-dimensional Galois representations. More precisely, we have Satake transforms between spherical algebras (say at some unramified place $v$): 
$$  \mathcal{H}_{G_1} \rightarrow \mathcal{H}_{M_W} \rightarrow \C[X_\star(T)]$$
for which the element $[1,1,\varpi_v,\varpi_v]$ is mapped to $$[1,
\varpi_v, 1, \varpi_v] + [\varpi_v, 1, \varpi_v, 1] + [\varpi_v,
\varpi_v,1,1] + [1, \varpi_v, 1, \varpi_v]$$ $$= ([1, \varpi_v, 1,
\varpi_v] +    [1, \varpi_v, 1, \varpi_v])  + [\varpi_v, 1, 1,
\varpi_v^{-1}]( [1, \varpi_v, 1, \varpi_v] +    [1, \varpi_v, 1,
\varpi_v]).$$  This expresses the relation between the Hecke operators
on $G_1$ and $M_W$, so that if $\chi \times \pi$ is an automorphic
representation contributing to $M$, it will contribute to the
cohomology of $G_1$ via  the compatible system of representations~$
\rho_{\pi} \oplus  (\rho_{ \pi} \otimes \chi)$.
\end{proof}

Before proving Theorem~\ref{thm: cohomology in terms of automorphic
  forms}, we introduce some notation. Let 
\[ h: \Res_{\C/\R}(\mathbf G_m)(\R)= \C^{\times} \to \GSp_4(\R) \]
be the homomorphism sending $x+iy$ to the matrix
\[
\begin{pmatrix}
  xI_2 & yS \\
  -yS & xI_2
\end{pmatrix}.
\]
Let $K^h$ denote the
centralizer of $h$ in $\GSp_4(\R)$ (acting by conjugation). Then since
$h(i) = J$, we see that we may identify $K^h = \R^\times U(2)$, so that $U(2)$ is a
maximal compact subgroup of the identity component
of~$\GSp_4(\R)$. Let~$\fg_{\C} = \fg^{0,0} \oplus \fg^{-1,1} \oplus \fg^{1,-1}$ denote
the Hodge structure on~$\fg$, where~$\fg^{0,0} = \gk_{h,\C}$ is the complex Lie algebra of~$K^h$.
Let~$\p^{+} = \fg^{-1,1}$, $\p^{-} = \fg^{1,-1}$, and~$\gP_h = \gk_{h,\C} \oplus \p^{-}$.
We now define~$P^{-}$ to be the parabolic with Lie algebra~$\gP_h$ with~$P^{-} \cap \GSp_4(\R) = K^h$.  We warn 
the reader that this parabolic is denoted by~$\mathcal{P}_h$ in~\cite{harris-ann-arb}, by~$Q^-$ in~\cite{blasius-harris-ramak},
and by~$Q$ in~\cite{CG}. (Note, however, that the fundamental object is really the Lie algebra~$\gP_h$ because~$P^{-}$
only intervenes below via its Lie algebra.)

 For each place~$v|\infty$, we
write~$P^{-}_v$ and~$ K^h_v$ for the corresponding groups
for~$\GSp_4(F_v)\cong\GSp_4(\R)$. We write~$V_{\kappa_v}$ for the
representation of~$ K^h_v\cong\R^\times U(2)$ such that
the  automorphic vector bundle 
corresponding to~$V_\kappa:=\otimes_{v|\infty}V_{\kappa_v}$ via Definition~1.3.2 of~\cite{blasius-harris-ramak}
is identified with~$\omega^{\kappa}$. 
We set
$P^-_\infty:=\prod_{v|\infty}P^-_v$ and~$K^h_\infty:=\prod_{v|\infty}K^h_v$.

\begin{proof}[Proof of Theorem~\ref{thm: cohomology in terms of automorphic forms}]

We begin by proving~(\ref{eqn: description of cohomology in terms of
  automorphic forms}). By Lemma~\ref{lem:harris-zucker} we can and do
replace $H^i(X_K^{G_1},\omega^\kappa(-D))$ by the interior
cohomology \[\overline{H}^i(X_K^{G_1},\omega^\kappa):=\im\left(H^i(X_K^{G_1},\omega^\kappa(-D))\to
  H^i(X_K^{G_1},\omega^\kappa)\right).\]
(This is the only place that we use our assumption that~$i\le 1$, or the non-Eisenstein localization.) Let $\CA_{\cusp}(G_1)\subset\CA_{(2)}(G_1)$ be respectively the
space of cuspidal automorphic forms on~$G_1$ with central
character~$|\cdot|^2$, and the space of square integrable forms with
this central character.  By~\cite[Thm.\ 2.7]{harris-ann-arb}, we have  inclusions 

\[ \bigoplus_{\pi
  \in\CA_{\cusp}(G_1)} \left(\left(\pi^{\infty}\right)^{K} \otimes
  H^i(\Lie P^{-}_\infty, K^h_\infty; \pi_{\infty}\otimes V_{\kappa}) \right)^{\oplus
  m_{\cusp}(\pi)}\subseteq \overline{H}^i(X_K^{G_1}, \omega^\kappa_{\C})^{|\cdot|^2}  \]
and 
\[ \overline{H}^i(X_K^{G_1}, \omega^\kappa_{\C})^{|\cdot|^2}
  \subseteq \bigoplus_{\pi \in \CA_{(2)}(G_1)}
  \left(\left(\pi^{\infty}\right)^{K} \otimes H^i(\Lie P^{-}_\infty, K^h_\infty;
    \pi_{\infty}\otimes V_{\kappa}) \right)^{\oplus m_{(2)}(\pi)}\]
where $m_{*}(\pi)$ denotes the multiplicity of $\pi$ in
$\CA_{*}(G_1)$. By Arthur's multiplicity formula
for~$\GSp_4$~\cite{MR2058604,GeeTaibi}, we in fact have
$m_{\cusp}(\pi)=m_{(2)}(\pi)=1$ for all~$\pi$.

In the proof of~\cite[Thm.\ 4.2.3]{blasius-harris-ramak}, it is shown that if $\pi\in \CA_{(2)}(G_1)$ with $H^i(\Lie P^{-}_\infty, K^h_\infty; \pi_{\infty}\otimes V_{\kappa})\ne 0$ and if the infinitesimal character of $\pi_\infty$ is sufficiently far away from all the root hyperplanes that it does not lie on,
then~$\pi_\infty$ is essentially tempered.  In view of the relation between the infinitesimal character of $\pi_\infty$ and $\kappa$ arising from the Casselman--Osborne theorem (see~\cite[Prop.\ 2.4.5]{blasius-harris-ramak}), the regularity condition on the infinitesimal character is exactly what we have assumed on $\kappa$ in (2).  Then by~\cite[Thm.\
4.3]{MR733320}, we in fact have $\pi\in\CA_{\cusp}(G_1)$, and so the inclusions above are equalities.  In addition, by~\cite[Thm.\ 3.5]{harris-ann-arb} (a
theorem of Mirkovi\'c) and ~\cite[Thm.\ 3.2.1]{blasius-harris-ramak},
for each~$v|\infty$ we have that $H^j(\Lie P^{-}_v, K^h_v;
\pi_{v}\otimes V_{\kappa_v})=0$ unless either:
\begin{itemize}
\item $l_v>2$, $j=0$, and~$\pi_v$ is the holomorphic discrete series of weight~$(k_v,l_v)$, or;
\item $l_v=2$, $j=0$, and~$\pi_v$ is the holomorphic limit of discrete series of
  weight~$(k_v,l_v)$, or;
\item $l_v=2$, $j=1$, and~$\pi_v$ is the generic limit of discrete series of weight~$(k_v,l_v)$.
\end{itemize}
Moreover, in each of these cases that $H^j(\Lie P^{-}_v, K^h_v;
\pi_{v}\otimes V_{\kappa_v})$ is nonzero, it is one-dimensional. The first two parts
of the theorem then follow from the K\"unneth formula.

  We now prove~(\ref{eqn: description of cohomology in terms of
  automorphic forms i equals 0 incl Eisenstein}). In this case the map
from $H^0(X_K^{G_1},\omega^\kappa(-D))$ to the interior cohomology is
an isomorphism by definition. Furthermore, by~\cite[Prop.\
  2.7.2]{harris-ann-arb}, the only~$\pi$ that contribute are
  automatically cuspidal (without needing to assume any regularity conditions). It follows from the theory of
  lowest weight representations, see for example~\cite[\S
  2.3]{MR2495302}, that if
  $H^0(\Lie P^{-}_v, K^h_v; \pi_{v}\otimes V_{\kappa_v})\ne 0$, 
  then~$\pi_v$ is the holomorphic (limit of) discrete series of
  weight~$(k_v,l_v)$, as required. \end{proof}

\section{Hida complexes}\label{sec: Hida complexes}In this section, we
construct  (higher) Hida theories for~$\GSp_4(\A_F)$. 
The classical Hida theory is developed in \cite{MR2055355} and takes
the form of a projective module over the total  weight space (which is
$2[F:\Q]$-dimensional). The construction of higher Hida theory was
carried out when $F = \Q$  in~\cite{pilloniHidacomplexes},  and takes
the shape of a perfect complex of amplitude $[0,1]$ over a one
dimensional hyperplane of the weight space.

We assume that~$p$ splits completely in~$F$ and we
construct all possible Hida theories, allowing the weight space at
each place above~$p$ to be either $1$- or $2$-dimensional. Many of our arguments are
simply the ``product over the places~$v|p$'' of the arguments
of~ \cite{MR2055355} and \cite{pilloniHidacomplexes}. To keep this paper at a reasonable
length, we will often refer to~\cite{pilloniHidacomplexes} for the
details of arguments which go over directly to our case.

The bookkeeping needed to deal with having multiple places above~$p$
is considerable, and in the hope of orienting the reader, we begin
this section with an overview of the arguments we will make. The main
theorem of this section (and the only theorem that we will need later
in the paper) is Theorem~\ref{theorem-p-adic-complex}, which proves
the existence of integral Hida complexes, and gives a control theorem
for them in sufficiently high weight. Say that a classical
weight~$\kappa=(k_v,l_v)_{v|p}$ with $k_v\geq l_v\geq2$ is ``singular'' at~$v$ if~$l_v=2$, and
``regular'' at~$v$ if~$l_v>2$. Fix
 any set~$I$ of places above~$p$; these will be the places at which we
 interpolate automorphic forms of singular weight, while at the places
 in~$I^c$, we interpolate forms of regular weight. (Thus traditional
 Hida theory considers the case~$I=\emptyset$, while the higher Hida
 theory of~\cite{pilloniHidacomplexes} is the case~$F=\Q$
 and~$I=\{p\}$.)

There is a Hecke operator~$U^I$ (an analogue of the~$U_p$
 operator for elliptic modular forms), which acts locally finitely on a complex
 of $p$-adic automorphic forms. The $U^I$-ordinary part~$M_I$ of this
 complex is a perfect complex over a weight space~$\Lambda_I$,
 concentrated in degrees~$[0,\#I]$. Furthermore, there is a
 constant~$C$ such that if~$k_v-l_v\ge C$ and~$l_v=2$ for~$v\in I$, and $l_v\ge C$
 for~$v\notin I$, then the~$H^0$ of the specialization of~$M_I$ in
 weight~$\kappa$ agrees with the ordinary part of the degree~$0$
 cohomology of~$X_K^{G_1}$. (We expect that in fact this
 specialization should be quasi-isomorphic to the ordinary part of the
 classical cohomology, but we do not prove this. We do prove that
 there is also an injection of~$\HH^1$s from the classical cohomology
 into that of~$M_I$, which we will make use of in~\S\ref{sec:
   higher Coleman theory}.)

 The definition of~$M_I$ is motivated by  the traditional case~$I=\emptyset$ considered in \cite{MR2055355}. In that case
 one considers the cohomology at infinite Iwahori level over the ordinary locus, with coefficients
 in a certain interpolation sheaf which can be thought of as an
 interpolation of the highest weight vectors in the finite dimensional representations
 of the group $\mathrm{GL}_2/F$. Since the ordinary locus is affine in the minimal compactification, one can prove that
 there is only cohomology in degree~$0$.  Then one cuts out the ordinary part using a projector  attached to $U^{I}$ and proves that this defines a finite projective module over the Iwasawa algebra.

  For general~$I$,  we instead consider the cohomology at infinite Klingen level of the locus which
 has $p$-rank at least~$1$ at places $w\in I$, and infinite Iwahori level over the ordinary (that
 is, $p$-rank $2$) locus at places $w\in I^c$. This locus is no longer affine and it has cohomology in higher
 degrees. In fact by relating the cohomology of the toroidal and minimal compactifications, one can show that the cohomology is supported in degrees~$[0,\#I]$.

 One of the major difficulties in  the proof of Theorem~\ref{theorem-p-adic-complex} is to show that the operator $U^I$ acts locally finitely (in order to be able to associate an ordinary projector) and that the ordinary projection defines a perfect complex. By
 Nakayama's lemma for complexes, one reduces to showing  that $U^I$ has these properties for  the cohomology modulo $p$, in some fixed weight. In particular, it suffices to consider the case of
 sufficiently large weight, i.e.\ the case that ~$k_v-l_v\ge C$
 and~$l_v=2$ for~$v\in I$, and $l_v\ge C$ for~$v\notin I$, for some
 constant~$C$. 
 The first part of the argument is to relate  this  cohomology  with the   cohomology of the automorphic vector bundle  of the corresponding weight over  the locus ~$X_{K,1}^{G_1,I}$ of the
 special fibre of the Shimura variety which has $p$-rank at least~$1$ at the places $w\in I$,
 and is ordinary at places~$w\in I^c$. This  boils down to a computation  at the level of the sheaf~$\omega^\kappa$ itself,  and to a computation in the Hecke algebra to show that the $U^I$-operator decreases the Klingen and Iwahori level.

 In the case of Hida theory for $0$-dimensional ``Shimura varieties''
 (e.g.\ $p$-adic families of automorphic forms on definite unitary groups, as considered
 in~\cite{MR2075765,ger}) these arguments at the level of the sheaf and
 the Hecke algebra are all that is needed. In the geometric setting,
 more work is needed to establish the required finiteness of the
 ordinary part of the cohomology in characteristic~$p$; recall that we
 are considering the cohomology on the locus ~$X_{K,1}^{G_1,I}$, so the cohomology groups are
 infinite dimensional before taking the ordinary parts. One has to
 show that (in sufficiently large weight, in characteristic~$p$)
 ordinary cohomology classes on this locus extend to the whole Shimura
 variety (which is proper and has finite cohomology). 

 In order to do this, one shows that (again, in sufficiently regular
 weight, in characteristic~$p$) the Hecke operator $U^I$ acts by zero
 on the complement of~$X_{K,1}^{G_1,I}$, so that after passing to ordinary
 parts, the cohomology agrees with that of the full Shimura variety,
 and is in particular finite-dimensional. 

 The vanishing of the Hecke operators on the part of the Shimura
 variety which is either of $p$-rank 0 (if~$w\in I$) or is
 non-ordinary is accomplished by local calculations, using the
 definitions of the Hecke operators as cohomological
 correspondences. The case of~$w\in I^c$ is relatively
 straightforward, as we are able to use the Hecke operator~$T_{w,1}$
 to prove this vanishing. (Note though that in this case we need to
 use the operator~$U_{w,2}$, which is the operator at Klingen level
 corresponding to~$T_w$, in the part of the argument explained above
 which takes place at the level of the sheaf.) The case~$w\in I$ is
 much more delicate, as we need to use the operator~$T_w$, which is
significantly harder to control. (In this case, though, we use the
same Hecke operators in the argument at the level of the sheaf as we
do for the geometric part of the argument.)

The arguments below are in fact written in roughly the reverse order
of the explanation above. We begin in~\S\ref{subsec: Hasse
  invariants and stratifications} by recalling some standard results
on Hasse invariants and the $p$-rank stratification, before proving
the vanishing of the Hecke operators in small $p$-rank at spherical
level in~\S\ref{subsec:
 vanishing for ordinary cohomology}. In~\S\ref{subsec: formal
 geometry} and~\ref{subsec: p adic forms over G1} we introduce the Igusa tower over the Shimura variety at
Klingen and Iwahori level, and define the interpolation sheaves whose cohomology we use to
define~$M_I$. We then define the Hecke operator~$U^I$ in~\S\ref{subsec: Hecke actions on p adic modular forms}, and 
in~\S\ref{subsec: perfect Hida complex} we prove
Theorem~\ref{theorem-p-adic-complex}, by relating the ordinary parts
of the cohomology at spherical and Klingen level, and then carrying
out the argument sketched above.
 
\subsection{Mod \texorpdfstring{$p$}{p}-geometry:  Hasse invariants and
  stratifications}\label{subsec: Hasse invariants and stratifications}In this section, we introduce the $p$-rank stratification on our Siegel variety and the definition of several Hasse invariants attached to this stratification.  The discussion follows
~\cite[\S 6.3, 6.4]{pilloniHidacomplexes}.

\subsubsection{Over $X$} We assume that $K =
K^pK_p$, $K_p = \prod_{v |p} K_v$ with $$K_v \in  \{
\mathrm{GSp}_4( \ocal_{F_v}), \Par(v)\}.$$

We fix a
polyhedral cone decomposition~$\Sigma$, and write $X=X_{K,\Sigma}$ if the context is clear. We let $\cG = A[p^\infty]$ be the $p$-divisible group corresponding to the semi-abelian scheme $A$ defined over $X$ (well defined up to prime-to-$p$ quasi-isogeny). This $p$-divisible group decomposes as $\cG = \prod_{v |p} \cG_v$.
If $K_v = \mathrm{GSp}_4( \ocal_{F_v})$, the $p$-divisible group $\cG_v$ defined over $X$ carries a principal quasi-polarization.   If $K_v = \Par(v)$ then the $p$-divisible group $\cG_v$ carries a  quasi-polarization of degree $p^2$: $\mathcal{G}_v \rightarrow \mathcal{G}_v^D$.  
Let $X_1$ be the reduction of $X$ modulo $p$. Then we let \[{\Ha} (\cG_v) \in  {\HH}^0(X_1  , \det
\omega_{\cG,v}^{ {p-1}})\] be the Hasse invariant corresponding to~
$\cG_v$; it is compatible with \'etale isogenies (by construction) and also with duality. 

For any place $v |p $, we let $X_1^{ \geq_v 2} = X_1^{ =_v 2}$ be
the open subscheme defined by $ {\Ha}(\cG_v) \neq 0$.  This is the
ordinary locus at $v$.  We let $X_1^{ \leq_v 1}$ be its complement
defined by ${\Ha}(\cG_v) = 0$. This is the non-ordinary locus at
$v$. It carries the reduced schematic structure by the proof of~\cite[Lem.\
6.4.1]{pilloniHidacomplexes}. (Whenever we use notation of the form~$X_1^{ \geq_v 2}$,
$X_1^{ \leq_v 1}$ etc., the superscript is referring to the
multiplicative rank of the  group scheme ~$\cG_v$.)

As a very special case of the general constructions
of~\cite{2015arXiv150705922B,2015arXiv150705032G}, there is a secondary Hasse invariant \[{\Ha}'(\cG_v) \in
{\HH}^0( X_1^{ \leq_v 1}, \det \omega_{\cG,v}^{p^2-1})\] (see also
\cite[\S 6.3.2]{pilloniHidacomplexes} when $K =\Par(v)$).  Its non-vanishing locus is
$X_1^{ =_v 1}$, the rank $1$ locus at $v$. We define its schematic  complement  $X_1^{
  =_v 0} $, the supersingular locus at $v$,  by the equation ${\Ha}'(\cG_v) = 0$. It carries a non-reduced
schematic structure, see~\cite[Rem.\ 6.4.1]{pilloniHidacomplexes}.

We can intersect the locally closed subschemes we have
defined. Consider disjoint subsets  $ I_1, \dots, I_r \subset \{v
|p\} $, symbols ${*}(i) \in \{ \leq, \geq , = \}$ for $1 \leq i \leq r$ and numbers
$a_i \in \{ 0, 1, 2\}$ for $1\le i\le r$.  Then we  define $X_1^{ {*}(i)_{I_i} a_i,~i = 1, \cdots, r}$ as the intersection of the spaces $X_1^{{*}(i)_{v} a_i}$ for all $1 \leq i \leq r$ and $v \in I_i$.  
 It will be convenient  to denote by $X_1^{ \geq 2} = X_1^{ \geq_{ \{v
     |p\}} 2}$ the ordinary locus and by $X_1^{ \geq 1} = X_1^{
   \geq_{ \{v |p \}} 1}$ the rank~$1$ locus.
 
Note that for any disjoint sets $I,J,K$, the scheme $X_1^{\le_I 1,\ge_J
  1,\ge_K 2}$ is Cohen--Macaulay, and indeed is a local complete
intersection over~$\Spec\F_p$. To see this, note that $X_1^{\le_I 1,\ge_J
  1,\ge_K 2}$ is open in~$X_1^{\le_I 1}$, and $X_1^{\le_I 1}$ is a complete intersection in~$X_1$,
because it is given by the vanishing of the Hasse
invariants~$\Ha(\cG_v)$ for~$v\in I$. Since $X_1$ itself is local
complete intersection by Proposition~\ref{prop: toroidal compactifications of G1}, the result follows. We will in particular repeatedly use
this fact in order to apply Lemma~\ref{lem: base change of
  correspondences}.

We will also frequently use some well-known results on the density of
the ordinary locus, and on the density of the $p$-rank 1 locus in the
$p$-rank less than or equal to~$1$ locus. We will need these results
in slightly greater generality than has been considered above. To this
end, consider disjoint subsets $ I_1, \cdots, I_r $ as above, and
let~$v|p$ be a place not contained in $ I_1\cup \dots\cup I_r $. 

We assume that $K = K^p K_p$, $K_p = \prod_{w |p, w \neq v} K_w \times K_v$, and that $$K_v\in\{\GSp_4(\cO_{F_v}),\Par(v),\Kli(v),\Si(v),\Iw(v)\},$$ while $K_w \in \{
\mathrm{GSp}_4( \ocal_{F_w}), \Par(w)\}$ for~$w\ne v$.  We can 
define topological spaces
$$|X_{K,1}^{ {*}(i)_{I_i} a_i} |, \  |X_{K,1}^{ {*}(i)_{I_i} a_i,=_v 2} |, \  |X_{K,1}^{ {*}(i)_{I_i}
  a_i,\le_v 1} |, \ \text{and} \  |X_{K,1}^{ {*}(i)_{I_i} a_i,=_v 1} |$$
   using the $p$-rank stratification as before. The point is that the $p$-rank is invariant under isogeny so we can consider the $p$-rank of any of the Barsotti--Tate groups of the chain. Note
that one could give these spaces a schematic structure by using the Hasse invariants, but  this structure will in general depend on which Barsotti--Tate group of the chain we use to define the Hasse invariants.  

Our claims about density are then the following:  ~$|X_{K,1}^{ {*}(i)_{I_i} a_i,=_v 2} |$
is dense in~$|X_{K,1}^{ {*}(i)_{I_i} a_i} |$, while
if we further assume that ~$K_v\in\{\GSp_4(\cO_{F_v}),\Par(v),\Kli(v)\}$,
$|X_{K,1}^{ {*}(i)_{I_i} a_i,=_v 1} |$ is dense in
$|X_{K,1}^{ {*}(i)_{I_i} a_i,\le_v 1} |$. To see this, it suffices
to prove the first statement in the case~$K_v=\Iw(v)$, and the second
statement in the case $K_v=\Kli(v)$. It then suffices to prove the
corresponding statements for the corresponding local models, which follows easily
from an explicit calculation. Indeed, the first statement is already
proved in~\cite{MR1227472}, while the second follows from an analysis of the
Kottwitz--Rapoport stratification at Iwahori level, and its image at
Klingen level; see~\cite[Thm.\ 4.2]{MR2811273} for a precise statement.

\subsubsection{Over $X^{G_1}$} The $p$-rank stratification is
independent of the polarization and therefore all of the spaces we have defined in this section carry an induced action of $(\ocal_{F})^{\times, +}_{(p)}$. It follows that the stratification descends to a stratification on $X_{1}^{G_1}$. We  moreover observe that the sheaf $\det \omega_{\cG,v}^{
  {p-1}}$ can be canonically descended  to a sheaf $\det
\omega_{\cG,v}^{p-1}$ on
$X_1^{G_1}$ (see Remark~\ref{rem-descent-Hasse}).  It follows that  the Hasse invariants
${\Ha}(\cG_v)$ and ${\Ha}'(\cG_v)$ (whose definition is independent of the polarization) also  descend to  sections of this sheaf over $X_1^{G_1}$ and $X_1^{G_1, \leq_v 1}$ respectively.

Therefore, if we consider disjoint subsets  $ I_1, \cdots, I_r \subset \{v
|p\} $, symbols ${*}(i) \in \{ \leq, \geq , = \}$ for $1 \leq i \leq r$ and numbers
$a_i \in \{ 0, 1, 2\}$ for $1\le i\le r$,  there is a unique locally closed subscheme  $(X^{G_1}_{K,1})^{ {*}(i)_{I_i} a_i,~i = 1, \cdots, r}$  of $X_{K,1}^{G_1}$ whose inverse image in $X_{K,1}$ is $(X_{K,1})^{ {*}(i)_{I_i} a_i,~i = 1, \cdots, r}$.

\begin{rem}
  In~\S\ref{subsubsec: Igusa towers} below we will define some
  other locally closed subschemes of the special fibres of the
  spaces~$X_{K,1}$ at Klingen and Iwahori level, which will be
  important in the rest of the paper. We caution the reader that these
  will \emph{not} be defined in terms of the $p$-ranks of the
  ~$\cG_v[p]$, but will rather depend on subschemes of~$\cG_v[p]$
  given by the Klingen and Iwahori level structures.
\end{rem}

\subsection{Vanishing theorem for ordinary cohomology} \label{subsec:
 vanishing for ordinary cohomology} 
 We assume that $K_p = G_1( \ZZ_p)$. Let
$\kappa = (k_v, l_v)_{v |p}$ be a weight (recall from~\S\ref{subsubsec: weights for
G and G1} that we are assuming that it satisfies the parity condition $k_v- l_v  =  0 \mod 2$). Let $S_p:=\{v|p\} =
I \coprod I^c$ be a partition. We write $ X_1^I:=X_1^{ \geq_{I} 1, \geq_{I^c} 2} \hookrightarrow X_1,$  an open subscheme, and similarly $X_1^{G_1,I} \hookrightarrow X_{1}^{G_1}$.  The main theorem of this subsection is: 

\begin{thm}\label{thm-vanishing}Let $T^I = \prod_{w \in I} T_w
  \prod_{w \in I^c} T_{w,1}$. There is a universal  constant $C$
  depending only on $p$ and $F$ but not on the tame level $K^p$ such
  that if~$l_w\ge 2$ for all~$w$, $k_w-l_w \geq C$ for all $w \in I$,
  and $l_w \geq C$ for all $w \in I^c$, then  ~$\mathrm{R}\Gamma(X_1^{G_1,I},
  \omega^{\kappa}(-D)) $carries a locally finite action of $T^I$. Furthermore, under this assumption on~$\kappa$,
\begin{enumerate}
\item $e(T^I)\mathrm{R}\Gamma(X_1^{G_1, I}, \omega^{\kappa}(-D))$ is a perfect complex of amplitude $[0, \# I ]$. 
\item The map $ e(T^I){\HH}^0(X^{G_1}_1, \omega^{\kappa}(-D)) \rightarrow e(T^I){\HH}^0(X_1^{ G_1, I}, \omega^{\kappa}(-D))$ is an isomorphism. 
\item The map $ e(T^I){\HH}^1(X^{G_1}_1, \omega^{\kappa}(-D)) \rightarrow e(T^I){\HH}^1(X_1^{ G_1, I}, \omega^{\kappa}(-D))$ is injective. 
\item If furthermore $l_w \geq 3$ for all $w \in I$, then $$ e(T^I)\mathrm{R}\Gamma(X^{G_1}_1, \omega^{\kappa}(-D)) \rightarrow e(T^I)\mathrm{R}\Gamma(X^{G_1,I}_1, \omega^{\kappa}(-D))$$ is a quasi-isomorphism. 
\end{enumerate}
\end{thm}

 Here $e(T^I)$ is the  ordinary projector associated to the operator $T^I$ (see~\S\ref{subsec: projectors}). We remark that $(2), (3), (4)$ of the theorem hold true for non-cuspidal cohomology as well. 
 \begin{rem}\label{rem: multiple style vanishing theorem could be
     improved in various ways}Various improvements on
   Theorem~\ref{thm-vanishing} should be possible. For example, the reader will see from the proof below that it is possible to prove
   that the Hecke operators at each place act locally finitely (rather
   than just proving it for their product), provided they satisfy
   explicit mild bounds on the weights (rather than depending on the
   indeterminate constant~$C$); see Remark~\ref{rem: could use this
     complex to prove local finiteness of individual Hecke operators}
   for one approach to this. It may also be possible to give explicit
   values of~$C$. For   the purposes of this paper the
   statement of Theorem~\ref{thm-vanishing} suffices, and is
   well-adapted to a (somewhat involved) inductive proof working one
   place at a time.
 \end{rem}

 We now briefly explain the main idea of the proof.  We will often  work at
 the level of $X_1$ rather than $X_1^{G_1}$. It is easier to work on
 $X_1$ because of the moduli interpretation. One can always deduce
 results for the cohomology on $X_1^{G_1}$ from results on the
 cohomology for $X_1$  by Proposition \ref{prop-splitGG_1}. We
 nevertheless warn the reader that $X_1$ has infinitely many connected
 components and therefore one cannot expect any finiteness results
 for the cohomology over $X_1$; accordingly, we work over~$X_1^{G_1}$
 when we want to show that a Hecke operator acts locally finitely.

The basic principle underlying these arguments is that  the ordinary projectors  $e(T_{w,1})$, $e(T_{w})$  can be used to kill many cohomology
 classes. This idea is already used in~\cite[\S
 7, \S 8]{pilloniHidacomplexes} (this is what we call Klingen
 vanishing below, because the Hecke operator $T_{w}$ is associated
 with the Klingen parabolic) and of course also in \cite{MR2055355}
 (this is what we call Siegel vanishing, because the Hecke operator
 $T_{w,1}$ is associated with the Siegel parabolic).

We will typically not comment on the commutativity of the actions of Hecke
 operators at one place with multiplication by Hasse invariants at
 other places, which is easily checked.

\subsubsection{Vanishing theorems: Siegel vanishing} \label{subsubsec:
Siegel vanishing}

Let $K = K^pK_p$ be a reasonable level
at $p$. We assume that $K_p = G_1( \ZZ_p)$. We let $X = X_{K, \Sigma}$ and $X^{G_1} = X_{K, \Sigma}^{G_1}$. Let
$\kappa = (k_v, l_v)_{v |p}$ be a weight. We begin with the
following theorem.

\begin{thm}\label{thm: Siegel locally finite quasi iso}There is a universal
  constant $C$ depending only on ~$p$ and ~$F$ but not on the tame level
  $K^p$ such that if ~$J\subseteq S_p$, and for each~$w\in J$, we have
  $l_w \geq C$, then $\mathrm{R}\Gamma(X_1^{G_1, \geq_J 2},
  \omega^{\kappa})$ has a locally finite action of
  $T^J:=\prod_{w\in J}T_{w,1}$, and \[ e(T^J)\mathrm{R}\Gamma(X_1^{G_1},
  \omega^{\kappa})\to e(T^J)\mathrm{R}\Gamma(X_1^{G_1, \geq_J 2},
  \omega^{\kappa})\]
is a quasi-isomorphism. In particular~$e(T^J)\mathrm{R}\Gamma(X_1^{G_1, \geq_J 2},
  \omega^{\kappa})$ is a perfect complex. The analogous statements
  also hold for cuspidal cohomology.
\end{thm} 
\begin{rem}
  In fact Theorem~\ref{thm: Siegel locally finite quasi iso} is not
  quite strong enough for our purposes; we will later replace it with
  Theorem~\ref{thm: Siegel locally finite quasi iso with Ia thrown
    in}, which is proved in exactly the same way. Our justification
  for presenting the material in this way is that the proof of
  Theorem~\ref{thm: Siegel locally finite quasi iso} is a good warmup
  for the arguments that we will later make to prove
  ``Klingen vanishing'', and it seems simplest to make these arguments
  before considering the Klingen level Hecke operators and the much more
  complicated statements and arguments that we make in that context.
\end{rem}

We only give the proof in the non-cuspidal case. The arguments go
through unchanged in the cuspidal setting. 
The proof of Theorem~\ref{thm: Siegel locally finite quasi iso} is by
induction on~$\# J$ (the case~$J=\emptyset$ being vacuous), and
depends on several lemmas. In our inductive argument we will feel free
to increase the constant~$C$ in a manner depending only on~$J$ without
comment. Write~$J=J'\cup\{w\}$, and assume that
Theorem~\ref{thm: Siegel locally finite quasi iso} holds for~$J'$.

Recall from \S
\ref{sect-def-Siegel-Hecke} that the correspondence underlying the
operator $T_{w,1}$ is $X_{K'}$ with $K' = K^p K'_p$ and $K'_p =
\prod_{v\neq w } K_v \times \Si( w)$. We let
$(X_{K'})_1$ denote the special fibre of this correspondence.  Let
$\kappa = (k_v, l_v)_{v |p}$ be a weight such that  $l_w \geq 2$, so that we have a
cohomological correspondence $T_{w,1}: p_2^{*} \omega^{\kappa}
\rightarrow p_1^! \omega^\kappa$.   By reduction modulo~ $p$, it follows
from Lemma~\ref{lem: base change of correspondences} (and the
flatness of $X_{K'}$ and $X_K$ over~$\Zp$) that we get a
cohomological correspondence still denoted $T_{w,1}: p_2^{*}
(\omega^{\kappa}\vert_{X_1}) \rightarrow p_1^! (
\omega^\kappa\vert_{X_1})$. This cohomological correspondence is a map of locally free sheaves over $(X_{K'})_1$. As in~\S\ref{subsubsec: base
  change}, this correspondence pulls back to the open subscheme~$X_1^{G_1, \geq_{J'} 2}$.

Adopting the notation of~\S\ref{subsec: Hasse invariants and
  stratifications} we consider the dense open subscheme $X_{K',1}^{\ge_{J}2}= (X_{K',1}^{\ge_{J'}2})^{=_w
  2}$ of $X_{K',1}^{\ge_{J'}2}$, which is by definition the ordinary locus at
$w$ (that is, the locus for which~$\cG_w$ is ordinary). This scheme is the union of
several types of connected components. Let $p_1^{*} \cG_w \rightarrow
p_2^{*} \cG_w$ be the universal map on the $p$-divisible group. We let
$(X_{K',1}^{\ge_{J'}2})^{=_w 2, et}$ be the \'etale components (that is, those for which the kernel of this
isogeny doesn't contain a multiplicative group), and we let $(X_{K',1}^{\ge_{J'}2})^{=_w 2, net}$ be the other components. We can therefore decompose the cohomological correspondence $T_{w,1}$ over $(X_{K',1}^{\ge_{J'}2})^{=_w 2}$ into $T_{w,1} = T_{w,1}^{et} + T_{w,1}^{net}$ where $T_{w,1}^{et}$ is the projection of $T_{w,1}$ on the \'etale components and $T_{w,1}^{net}$ is the projection on the other components.

\begin{lem}\label{lem: Siegel sheaves commute over ordinary} The map $T_{w,1}^{net}$ is zero as soon as $l_w \geq 3$. For all $l_w \geq 3$ we have a commutative diagram of maps of sheaves over   $(X_{K',1}^{\ge_{J'}2})^{=_w 2}$: 
\begin{eqnarray*}
\xymatrix{ p_2^{*} \omega^\kappa \ar[r]^{T_{w,1}^{et}} \ar[d]^{ p_2^{*} {\Ha}(\cG_w)}& p_1^! \omega^\kappa  \ar[d]^{p_1^{*} {\Ha}(\cG_w)}\\
p_2^{*}( \omega^\kappa \otimes \det \omega_{\cG_w}^{p-1}) \ar[r]^{T_{w,1}^{et}} & p_1^! (\omega^\kappa \otimes \det \omega_{\cG_w}^{p-1})}
\end{eqnarray*}
\end{lem}
\begin{proof} The first point follows from an inspection of the proof
  of Lemma~\ref{lem-normalization1} (since $l_w+1,k_w+l_w>3$). The second point follows from the fact that the Hasse invariant commutes with \'etale isogenies. 
\end{proof}

\begin{lem}\label{lem: Siegel sheaves commute} The following diagram of locally free sheaves on $X_{K',1}^{\ge_{J'}2}$ is commutative for $l_w \geq 3$: 
\begin{eqnarray*}
\xymatrix{ p_2^{*} \omega^\kappa \ar[r]^{T_{w,1}} \ar[d]^{ p_2^{*} {\Ha}(\cG_w)}& p_1^! \omega^\kappa  \ar[d]^{p_1^{*} {\Ha}(\cG_w)}\\
p_2^{*} (\omega^\kappa \otimes \det \omega_{\cG_w}^{p-1}) \ar[r]^{T_{w,1}} & p_1^! (\omega^\kappa \otimes \det \omega_{\cG_w}^{p-1})}
\end{eqnarray*}
\end{lem} 
\begin{proof} Since $X_{K',1}^{\ge_{J'}2}$ is Cohen--Macaulay and all of the
  sheaves are locally free, it suffices to prove the commutativity
  over a dense open subscheme. We may therefore prove it over
  $(X_{K',1}^{\ge_{J'}2})^{=_w 2}$, so we are done by Lemma~\ref{lem: Siegel sheaves commute over ordinary}.
\end{proof}

\begin{lem}
  \label{lem: Hasse not a zero divisor}${\Ha}(\cG_w)$ is not a
  zero divisor on~$X_1^{\ge_{J'}2}$, and $p_2^{*}{\Ha}(\cG_w)$ is not a
  zero divisor on~$X_{K',1}^{\ge_{J'}2}$.
\end{lem}
\begin{proof}
Since $X_1^{\ge_{J'}2}$ and $X_{K',1}^{\ge_{J'}2}$ are Cohen--Macaulay, this follows
from the fact that the non-ordinary loci have codimension~$1$.\end{proof}

In what follows, we warn the reader that while the schemes $p_i^{-1}(X_1^{\geq_{J'}2,\leq_w1})=X_{K',1}^{\ge_{J'}2} \times_{p_i, X_1^{\ge_{J'}2}} X_1^{\geq_{J'}2,\leq_w 1}$ for $i=1,2$ have the same underlying topological spaces, they have different (non reduced) scheme structures.  In particular sheaves like $p_i^*\omega^\kappa|_{X^{\geq_{J'}2,\leq_w1}}$ for $i=1,2$ have different (scheme theoretic) support.

\begin{lem}\label{lem: restricting Siegel correspondence}For all~$l_w\ge
  p+2$ the cohomological correspondence~$T_{w,1}$ restricts to give a
  cohomological correspondence \[T_{w,1}: p_2^{*}
    (\omega^\kappa\vert_{X_1^{\geq_{J'}2,\leq_w 1}}) \rightarrow p_1^!
    (\omega^\kappa\vert_{X_1^{\geq_{J'}2,\leq_w 1}}).\]\end{lem}

\begin{proof}Since the cokernel of
  \[p_2^{*}\omega^\kappa\stackrel{p_2^{*}{\Ha}(\cG_w)}{\longrightarrow}
  p_2^{*} (\omega^\kappa \otimes \det \omega_{\cG_w}^{p-1})\] is $p_2^{*}
    (\omega^\kappa\otimes \det \omega_{\cG_w}^{p-1}\vert_{
      X_1^{\geq_{J'}2,\leq_w 1}})$, and (by Lemmas~\ref{lem: base change of
      correspondences} and~\ref{lem: Hasse not a zero divisor}) the cokernel of \[p_1^! \omega^\kappa
    \stackrel{p_1^{*}{\Ha}(\cG_w)}{\longrightarrow}  p_1^! (\omega^\kappa
    \otimes \det \omega_{\cG_w}^{p-1})\] is~$p_1^!
    (\omega^\kappa\otimes \det \omega_{\cG_w}^{p-1}\vert_{
      X_1^{\geq_{J'}2,\leq_w 1}})$, it follows from Lemma~\ref{lem: Siegel
      sheaves commute} that provided that~$l_w\ge 3$, $T_{w,1}$
    restricts to give a cohomological correspondence \[p_2^{*}
    (\omega^\kappa\otimes \det \omega_{\cG_w}^{p-1}\vert_{
      X_1^{\geq_{J'}2,\leq_w 1}})\to p_1^!
    (\omega^\kappa\otimes \det \omega_{\cG_w}^{p-1}\vert_{
      X_1^{\geq_{J'}2,\leq_w 1}}), \]as required.
\end{proof}

\begin{lem}\label{lem-vanishing-Siegel} There is a universal constant
  $C$ which depends only on $F$ and $p$ \emph{(}but not on the tame level~$K^p$\emph{)} such
  that the map of Lemma \ref{lem: restricting Siegel correspondence} \[T_{w,1}: p_2^{*}
    (\omega^\kappa\vert_{ X_1^{\geq_{J'}2,\leq_w 1}}) \rightarrow p_1^!
    (\omega^\kappa\vert_{ X_1^{\geq_{J'}2,\leq_w 1}})\] is zero for all $l_w \geq C$. 
\end{lem} 

\begin{proof}We may and do assume that $J'=\emptyset$, as the general
  case follows immediately from this by restriction to an open.  We moreover note that it suffices to find such a
  constant~$C$ for a single tame level~$K^p$. Indeed,
  if~$K_1^p\subseteq K_2^p$ are two choices of tame level, the natural
  forgetful map $X_{K'_pK_1^p}\to X_{K'_pK_2^p}$ commutes
  with~$p_1$ and~$p_2$, and is faithfully flat, from which it follows
  that~$C$ works for~$K_1^p$ if and only if it works for~$K_2^p$.

 Let $\mathcal{J}$ be the ideal defining $X^{\leq_w 1}_1$
  in $X$ and let $\mathcal{I} = p_1^{*} \mathcal{J}$. We need to
  prove that for $l_w$ sufficiently large, the cohomological correspondence over $X_{K'}$,  $T_{w,1}: p_2^{*} \omega^{\kappa} \rightarrow p_1^! \omega^\kappa$ factors through $T_{w,1}: p_2^{*} \omega^{\kappa} \rightarrow \mathcal{I} p_1^! \omega^\kappa$. 

By definition, we have $T_{w,1} = p^{-3} \Theta(\kappa)$, where
$\Theta(\kappa)$  is the composite of a map $\Theta_1(\kappa):
p_2^{*} \omega^\kappa \rightarrow p_1^{*} \omega^\kappa$ and a
fundamental class $\Theta_2(\kappa):  p_1^{*} \omega_\kappa
\rightarrow p_1^! \omega^\kappa$, so in turn we need to show that
$\Theta(\kappa)$ factors through $p^3\mathcal{I} p_1^! \omega^\kappa$
(of course, we have already shown that it factors through~$p^3
p_1^{!} \omega^\kappa$). 

Let $x$ be  a generic point of $V(\mathcal{I}) \subset X_{K'}$. It
corresponds to a Barsotti--Tate group in characteristic $p$ whose
$p$-rank at $w$ is exactly one.  The map $p_2^{*} \det
\omega_{\cG_w} \rightarrow p_1^{*} \det \omega_{\cG_w}$ is zero over
$k(x)$ because the isogeny $p_1^{*} \cG_w \rightarrow p_2^{*}
\cG_w$ is not \'etale at $x$. Let $I(x)$ be the ideal defining
the Zariski closure $\overline{x}$ in $X_{K'}$. We deduce that the map $\Theta_1(\kappa)
: p_2^{*} \omega^\kappa \rightarrow p_1^{*} \omega^\kappa$ factors
through $I(x)^{l_w} p_1^{*} \omega^{\kappa}$. 

It follows that the map
$\Theta(\kappa)$ factors through $p^3 p_1^{*} \omega^\kappa \bigcap
\cap_{x} I(x)^{l_w} p_1^{*} \omega^\kappa$. By the Artin--Rees
lemma, it factors through $p^3 \mathcal{I} p_1^{*} \omega^\kappa$
for $l_w$ larger than a constant $C$, as required.  (We note that strictly speaking $V(I)$ has infinitely many connected components, however there are only finitely many orbits for the action of $(\O_F)^{\times,+}_{(p)}$, so there is some constant $C$ which works for all of them.)\end{proof} 

\begin{proof}[Proof of Theorem~\ref{thm: Siegel locally finite quasi
    iso}]
    
Take~$C$ as in Lemma~\ref{lem-vanishing-Siegel}. Recall that we write~$J=J'\cup\{w\}$, and we are assuming that the
theorem holds for~$J'$. We begin by showing that the action of~$T^J$ on $\mathrm{R}\Gamma(X^{G_1,\ge_{J'}2}_1,
  \omega^{\kappa})$ is locally finite. By the inductive hypothesis,
  the action of~$T^{J'}$ on this complex is locally finite, and $e(T^{J'})\mathrm{R}\Gamma(X^{G_1,\ge_{J'}2}_1,
  \omega^{\kappa})$ is perfect, so that in particular the action
  of~$T^J$ on $e(T^{J'})\mathrm{R}\Gamma(X^{G_1,\ge_{J'}2}_1,
  \omega^{\kappa})$ is locally finite. It is therefore enough to show
  that~$T^J$ acts locally finitely on $(1-e(T^{J'}))\mathrm{R}\Gamma(X^{G_1,\ge_{J'}2}_1,
  \omega^{\kappa})$.   Since~$T^{J'}$ acts
  locally nilpotently on the complex $(1-e(T^{J'}))\mathrm{R}\Gamma(X^{G_1,\ge_{J'}2}_1,
  \omega^{\kappa})$ by definition, so does~$T^J$, so in particular it acts locally
  finitely, as required.

  Now we
  consider the exact triangle
$$ \mathrm{R}\Gamma(X^{G_1,\ge_{J'}2}_1, \omega^{\kappa}) \rightarrow
\mathrm{R} \Gamma (X^{G_1,\ge_{J'}2}_1, \omega^\kappa \otimes \det
\omega_{\cG_w}^{p-1} ) \rightarrow \mathrm{R} \Gamma (X^{G_1,\ge_{J'}2,\leq_w 1}_1, \omega^\kappa \otimes \det \omega_{\cG_w}^{p-1} ) \stackrel{+1}\rightarrow $$
The operator $T^J$ acts everywhere and is zero on
$\mathrm{R} \Gamma (X^{G_1,\ge_{J'}2,\leq_w 1}_1, \omega^\kappa \otimes \det
\omega_{\cG_w}^{p-1} )$ by Lemma~\ref{lem-vanishing-Siegel}. We therefore deduce that
$$e(T^J) \mathrm{R}\Gamma(X^{G_1,\ge_{J'}2}_1, \omega^{\kappa}) = e(T^J)
\mathrm{R}\Gamma(X^{G_1,\ge_{J'}2}_1, \omega^{\kappa} \otimes \det
\omega_{\cG_w}^{p-1} ).$$
Since
$\mathrm{R}\Gamma(X^{G_1,\geq_J 2}_1, \omega^{\kappa}) = 
\colim_n\mathrm{R}\Gamma(X^{G_1,\ge_{J'}2}_1, \omega^{\kappa} \otimes \det
\omega_{\cG_w}^{n(p-1)} )$, the theorem follows.
\end{proof}

\subsubsection{Vanishing theorem: Siegel and  Klingen vanishing}\label{subsubsec:
  Klingen vanishing} We now turn to the more general situation, which involves the study of 
 the Hecke operator~$T_w$. Our analysis is similar to that of~\S\ref{subsubsec:
Siegel vanishing}, but it is rather more involved because~$T_w$ is
defined as the composite of two correspondences and because we need to study  the supersingular locus at $w$ rather than the non-ordinary locus at $w$. This subsection is
devoted to the proof of the following theorem, which implies most of Theorem~\ref{thm-vanishing}.

Let~$I_a,I_b,J$ be pairwise
  disjoint subsets of~$S_p$. Then we will write \[X_1^{I_{a,b},J}:=X_1^{\le_{I_a}1,\ge_{I_b}1,\ge_J
  2}~\textrm{and}~X_1^{G_1,I_{a,b},J}:=X_1^{G_1,\le_{I_a}1,\ge_{I_b}1,\ge_J
  2}. \]
\begin{thm}\label{thm-Klingen-vanish}Let~$I_a,I_b,J$  be pairwise
  disjoint subsets of~$S_p$.  Set
  $T^{I_b,J} = \prod_{w \in I_b} T_w \prod_{w \in J} T_{w,1}$. Then there
  is a universal constant $C$ depending only on $p$ and $F$ but not on
  the tame level $K^p$ such that if $k_w-l_w \geq C$ and~$l_w\ge 2$
  for all $w \in I_b$, and $l_w \geq C$ for all $w \in J$, then
  $\mathrm{R}\Gamma(X_1^{G_1,I_{a,b},J}, \omega^{\kappa})
  $ carries a locally finite action of $T^{I_b,J}$.  Furthermore:

\begin{enumerate}\item $e(T^{I_b,J})\mathrm{R}\Gamma(X_1^{G_1,I_{a,b},J}, \omega^{\kappa})$ is a perfect complex. 
\item The map \[ e(T^{I_b,J}){\HH}^0(X^{G_1,\le_{I_a}1}_1, \omega^{\kappa}) \rightarrow e(T^{I_b,J}){\HH}^0(X_1^{G_1,I_{a,b},J}, \omega^{\kappa})\] is an isomorphism. 
\item The map \[ e(T^{I_b,J}){\HH}^1(X^{G_1,\le_{I_a}1}_1, \omega^{\kappa}) \rightarrow e(T^{I_b,J}){\HH}^1(X_1^{G_1,I_{a,b},J}, \omega^{\kappa})\] is injective. 
\item If furthermore $l_w \geq 3$ for all $w \in I_b$, then $$
  e(T^{I_b,J})\mathrm{R}\Gamma(X^{G_1,\le_{I_a}1}_1, \omega^{\kappa}) \rightarrow
  e(T^{I_b,J})\mathrm{R}\Gamma(X_1^{G_1,I_{a,b},J}, \omega^{\kappa})$$ is a
  quasi-isomorphism.  
\end{enumerate}
Moreover, the same results hold for cuspidal cohomology.
\end{thm}

We only give the proof in the non-cuspidal setting. The same arguments work in the  cuspidal case. The proof of this result again depends on several lemmas. We will
firstly prove the result in the case~$I_b=\emptyset$, by induction on~$\#J$. We will then prove the general case by induction on~$\# I_b$.

   Recall from~\S\ref{sect-def-Siegel-Hecke}, that   if we set   $K' = K^p K'_p$ with $K'_p = \prod_{v\neq w } K_v \times \Si( w)$,  there  is  a  cohomological correspondence  of Siegel type: 
 $$T_{w,1}: p_2^{*}
   (\omega^{\kappa}\vert_{X_{K}}) \rightarrow p_1^! (
   \omega^\kappa\vert_{X_{K}}).$$
    By reduction
   modulo $p$  and Lemma~\ref{lem: base change of correspondences}, we
   again get a
   cohomological correspondence:
$T_{w,1}: p_2^{*}
   (\omega^{\kappa}\vert_{X_{K,1}}) \rightarrow p_1^! (
   \omega^\kappa\vert_{X_{K,1}})$.
We let $X^{I_{a,b}, J}_{K',1}$ be the pre-image of $X_
   {1}^{I_{a,b},J}$ in $X_{K',1}$ (via any of the projections, it doesn't matter). 
   These correspondences may be
   restricted to $X^{I_{a,b}, J}_{K',1}$ whenever $w\notin I_a$ by another application of Lemma~\ref{lem: base change of correspondences}, because this correspondence obviously commutes with the Hasse invariants at places in $I_a$.

\begin{thm}\label{thm: Siegel locally finite quasi iso with Ia thrown in}There is a universal
  constant $C$ depending only on ~$p$ and ~$F$ but not on the tame level
  $K^p$ such that if $I_a$, $J$ are disjoint, and if $l_w \geq C$ for all $w \in J$,  then $\mathrm{R}\Gamma(X_1^{G_1,\le_{I_a}1, \geq_J 2},
  \omega^{\kappa})$ has a locally finite action of
  $T^{\emptyset,J}:=\prod_{w\in J}T_{w,1}$, and \[ e(T^{\emptyset,J})\mathrm{R}\Gamma(X_1^{G_1,\le_{I_a}1},
  \omega^{\kappa})\to e(T^{\emptyset,J})\mathrm{R}\Gamma(X_1^{G_1,\le_{I_a}1, \geq_J 2},
  \omega^{\kappa})\]
is a quasi-isomorphism. In particular~$e(T^{\emptyset,J})\mathrm{R}\Gamma(X_1^{G_1,\le_{I_a} 1, \geq_J 2},
  \omega^{\kappa})$ is a perfect complex.
\end{thm} 
\begin{proof}
  The
  case~$I_a=\emptyset$ is Theorem~\ref{thm: Siegel locally finite
    quasi iso}, and the theorem at hand may be proved by an identical
  inductive argument on $\# J$, once we have proved the base
  case~$J=\emptyset$. But in this case $X_1^{G_1,\le_{I_a}1}$ is
  proper, so $\mathrm{R}\Gamma(X_1^{G_1,\le_{I_a} 1},
  \omega^{\kappa})$ is a perfect complex, and we are done.
\end{proof}

We now reintroduce Klingen type correspondences. Let $w |p$. By~\S\ref{sect-def-Klingen-Hecke}     if we set   $K' = K^p K'_p$ with $K'_p = \prod_{v\neq w } K_v \times {\Kli}( w)$, and $K'' = K^p K''_p$ with $K''_p = \prod_{v\neq w } K_v \times {\Par}( w)$,
there  are cohomological correspondences  of Klingen type:  $T'_w: p_2^{*}
   (\omega^{\kappa}\vert_{X_{K''}}) \rightarrow p_1^! (
   \omega^\kappa\vert_{X_{K}})$ and $T''_w: p_1^{*}
   (\omega^{\kappa}\vert_{X_{K}}) \rightarrow p_2^! (
   \omega^\kappa\vert_{X_{K'}})$ for all weights $\kappa = (k_v, l_v)$ with $k_w \geq l_w \geq 2$. 
    By reduction
   modulo $p$  and Lemma~\ref{lem: base change of correspondences}, we again get
   cohomological correspondences: $T'_w: p_2^{*}
   (\omega^{\kappa}\vert_{X_{K'',1}}) \rightarrow p_1^! (
   \omega^\kappa\vert_{X_{K,1}})$ and $T''_w: p_1^{*}
   (\omega^{\kappa}\vert_{X_{K,1}}) \rightarrow p_2^! (
   \omega^\kappa\vert_{X_{K'',1}})$. 
We let $X^{I_{a,b}, J}_{K',1}$ be the pre-image of $X_
   {1}^{I_{a,b},J}$ in $X_{K',1}$. 
   These correspondences may be
   restricted to $X^{I_{a,b}, J}_{K',1}$ whenever $w\notin I_a$  by another application of Lemma~\ref{lem: base change of correspondences}, because they obviously commute with the Hasse invariants at places in $I_a$.

   In the rest of this section we prove Theorem~\ref{thm-Klingen-vanish}
by induction on~$\# I_b$. To this end, choose~$w\in I_b$, write~$I=I'\coprod\{w\}$, and
write~$I'_a=I'\cap I_a=I_a$, and $I'_b=I'\cap I_b$. We assume that
Theorem~\ref{thm-Klingen-vanish} holds (for some value of~$C$, which
we fix) for all smaller values of~$\#
I_b$ (as we may, having proved the case~$I_b=\emptyset$ in Theorem~\ref{thm: Siegel locally finite quasi iso}).

   We now consider the scheme $(X_{K',1}^{I'_{a,b},J})^{=_w2}$ (which is again by
   definition the subscheme where~$\cG_w$ is ordinary), which decomposes
   into several components. Let $p_1^{*} \cG \rightarrow p_2^{*}
   \cG$ be the universal isogeny of degree $p^3$. We denote by
   $(X_{K',1}^{I'_{a,b},J})^{=_w2, et }$ the ``\'etale'' components, namely those where the kernel of the universal isogeny has multiplicative rank $1$ (so it is as \'etale as possible),  and by $(X_{K',1}^{I'_{a,b},J})^{=_w2, net }$ the other components where the kernel has multiplicative rank $2$. 
   
 This provides a  decomposition of the correspondence $T'_{w} = T'_{w,et} + T'_{w,net}$ where $T'_{w,et}$ stands for the projection on the ``\'etale'' components and $T'_{w,net}$ for the projection on the ``non-\'etale'' components.

 We also have an isogeny $p_2^{*} \cG \rightarrow p_1^{*} \cG$ of degree $p$ and over $(X_{K',1}^{I'_{a,b},J})^{=_w2, et }$ this isogeny has multiplicative kernel, while it is \'etale over $(X_{K',1}^{I'_{a,b},J})^{=_w2, net }$ (observe that the \'etale and non-\'etale components are interchanged when we pass from the isogeny $p_1^{*} \cG \rightarrow p_2^{*}
   \cG$ to the isogeny $p_2^{*} \cG \rightarrow p_1^{*}
   \cG$). This provides a second decomposition  
 $T''_{w} = T''_{w,et} + T''_{w,net}$ where $T''_{w,et}$  stands for the projection on  $(X_{K',1}^{I'_{a,b},J})^{=_w2, net }$ and $T''_{w,net}$ for the projection on $(X_{K',1}^{I'_{a,b},J})^{=_w2, et }$. 
 
 \begin{lem}\label{lem-commutTW ordinary} If $l_w \geq 2$, and $k_w \geq 3$ then over  $(X_{K',1}^{I'_{a,b},J})^{=_w2}$ we have $T'_{w,net} = T''_{w,net} = 0$. Moreover, the following diagrams are commutative: 

\begin{eqnarray*}
\xymatrix{ p_2^{*} \omega^\kappa \ar[r]^{T'_{w,et}} \ar[d]^{ p_2^{*} {\Ha}(\cG_w)}& p_1^! \omega^\kappa  \ar[d]^{p_1^{*} {\Ha}(\cG_w)}\\
p_2^{*} (\omega^\kappa \otimes \det \omega_{\cG_w}^{p-1}) \ar[r]^{T'_{w,et}} & p_1^! (\omega^\kappa \otimes \det \omega_{\cG_w}^{p-1})}
\end{eqnarray*}

\begin{eqnarray*}
\xymatrix{ p_1^{*} \omega^\kappa \ar[r]^{T''_{w,et}} \ar[d]^{ p_1^{*} {\Ha}(\cG_w)}& p_2^! \omega^\kappa  \ar[d]^{p_2^{*} {\Ha}(\cG_w)}\\
p_1^{*} (\omega^\kappa \otimes \det \omega_{\cG_w}^{p-1}) \ar[r]^{T''_{w,et}}& p_2^! (\omega^\kappa \otimes \det \omega_{\cG_w}^{p-1})}
\end{eqnarray*}

\end{lem} 
\begin{proof} That $T'_{w,net} = T''_{w, net} = 0$ follows from an
  inspection of the proof of Lemma \ref{lem-normalization2}; more
  precisely, by the proofs of~\cite[Lem.\ 7.1.1,
  7.1.2]{pilloniHidacomplexes}, $T_{w,net}'$ is divisible
  by~$p^{k_w-2}$, and~$T''_{w, net}$ is divisible
  by~$p^{l_w-1}$. The commutativity of the second diagram follows
  immediately from the fact that the Hasse invariant commutes with
  \'etale isogenies. The commutativity of the first diagram is
  slightly more delicate; see the proof of~
  \cite[Prop.\ 7.4.1.1]{pilloniHidacomplexes}, which explains how it
  reduces to~\cite[Lem.\ 6.3.4.1]{pilloniHidacomplexes}.
\end{proof}

\begin{lem}\label{lem-commutTW} The following diagrams of locally free
  sheaves on  $X_{K',1}^{I'_{a,b},J}$ are commutative for $l_w \geq 2$  and $k_w \geq 3$: 

\begin{eqnarray*}
\xymatrix{ p_2^{*} \omega^\kappa \ar[r]^{T'_{w}} \ar[d]^{ p_2^{*} {\Ha}(\cG_w)}& p_1^! \omega^\kappa  \ar[d]^{p_1^{*} {\Ha}(\cG_w)}\\
p_2^{*} (\omega^\kappa \otimes \det \omega_{\cG_w}^{p-1}) \ar[r]^{T'_w} & p_1^! (\omega^\kappa \otimes \det \omega_{\cG_w}^{p-1})}
\end{eqnarray*}

\begin{eqnarray*}
\xymatrix{ p_1^{*} \omega^\kappa \ar[r]^{T''_{w}} \ar[d]^{ p_1^{*} {\Ha}(\cG_w)}& p_2^! \omega^\kappa  \ar[d]^{p_2^{*} {\Ha}(\cG_w)}\\
p_1^{*} (\omega^\kappa \otimes \det \omega_{\cG_w}^{p-1}) \ar[r]^{T''_w} & p_2^! (\omega^\kappa \otimes \det \omega_{\cG_w}^{p-1})}
\end{eqnarray*}

\end{lem} 
\begin{proof} Since all sheaves are locally free and $X_{K',1}^{I'_{a,b},J}$ is
  Cohen--Macaulay, it is enough to check the commutativity over the
  dense open subscheme $(X_{K',1}^{I'_{a,b},J})^{=_w2}$, which is
  Lemma~\ref{lem-commutTW ordinary}.
\end{proof} 

\begin{cor}
  \label{cor: Tw locally finite on ordinary}Assume that for all~$v\in
  I'$ we have~$l_v\ge 2$ and~$k_v-l_v\ge C$, and that for all $v\in J$
  we have~$l_v\ge C$.   Then the action of $T^{I_b,J}$ on
  ${\RGamma}(X_{K,1}^{G_1, I'_{a,b},J,=_w2},\omega^\kappa)$ is locally finite if $l_w \geq 2$  and $k_w \geq 3$.
\end{cor} 
\begin{proof}We have
  ${\HH}^{*}(X_{K,1}^{G_1,I'_{a,b},J,=_w2},\omega^\kappa)=\varinjlim_n{\HH}^{{*}}(X^{G_1,I'_{a,b},J}_{K,1},\omega^\kappa\otimes\det\omega_{\cG_w}^{(p-1)n})$
  where the transition maps (given by multiplication by~$\Ha(\cG_w)$)
  are $T^{I_b,J}$-equivariant by Lemma~\ref{lem-commutTW}.  By the
  inductive hypothesis, each $e(T^{I'_b,J}){\HH}^{{*}}(X^{G_1,I'_{a,b},J}_{K,1},\omega^\kappa\otimes\det\omega_{\cG_w}^{(p-1)n})$ is
finite-dimensional and $T^{I_b,J}$-stable, while $T^{I_b,J}$ acts
locally nilpotently on $(1-e(T^{I'_b,J})){\HH}^{{*}}(X^{G_1,I'_{a,b},J}_{K,1},\omega^\kappa\otimes\det\omega_{\cG_w}^{(p-1)n})$ (because
$T^{I'_b,J}$ does). The result follows.
\end{proof}

Exactly as in Lemma~\ref{lem: restricting Siegel correspondence}, for
all $l_w \geq p+1 = p-1 + 2$ and $k_w \geq p+2 = 3 + p-1$ we obtain cohomological correspondences:

$$ T'_{w}: p_2^{*} (\omega^\kappa\vert_{X_{K'',1}^{I'_{a,b},J,\leq_w 1}}) \rightarrow p_1^! (\omega^\kappa\vert_{X_{K,1}^{I'_{a,b},J,\leq_w 1}})$$ and $$ T''_{w}: p_1^{*} (\omega^\kappa\vert_{X_{K,1}^{I'_{a,b},J,\leq_w 1}}) \rightarrow p_2^! (\omega^\kappa\vert_{X_{K'',1}^{I'_{a,b},J,\leq_w 1}}).$$

We now consider the space $(X_{K'',1}^{I'_{a,b},J})^{=_w1}$
where~$\cG_w[p]$ has $p$-rank $1$. We have a universal
quasi-polarization $\cG_w \rightarrow \cG_w^D$ over $X_{K''}$. Over the  interior of the
moduli space,  the  kernel  of the quasi-polarization is a self dual
rank $p^2$ group scheme which is either connected or an extension of a
multiplicative by an \'etale group scheme.  The space
$(X_{K'',1}^{I'_{a,b},J})^{=_w1}$   decomposes as the union of connected components
$(X_{K'',1}^{I'_{a,b},J})^{=_w1,00}$  and  $(X_{K'',1}^{I'_{a,b},J})^{=_w1, \met}$  for which
the kernel of the quasi-polarization doesn't contain (respectively
contains) a multiplicative group (see~\cite[Lem.\ 7.4.2.3]{pilloniHidacomplexes}).

We now consider the space $(X_{K',1}^{I'_{a,b},J})^{=_w1}$, which we view here only as a topological space (it has multiple natural non-reduced scheme structures defined by the vanishing of either $p_1^*\Ha(\G_w)$ or $p_2^*\Ha(\G_w)$).  We have the chain of
isogenies $p_1^{*} \cG \rightarrow p_2^{*} \cG \rightarrow
p_1^{*} \cG$ where the composite is multiplication by $p$. We have a
decomposition of $(X_{K',1}^{I'_{a,b},J})^{=_w1}$ as a union of connected
components: $(X_{K',1}^{I'_{a,b},J})^{=_w1,m}$, $(X_{K',1}^{I'_{a,b},J})^{=_w1, et}$ and
$(X_{K',1}^{I'_{a,b},J})^{=_w1,00}$. Here the open and closed subspace
$(X_{K',1}^{I'_{a,b},J})^{=_w1,m}$ is the locus where the kernel of $p_2^{*} \cG
\rightarrow p_1^{*} \cG$ is an \'etale group scheme; the open and
closed subspace $(X_{K',1}^{I'_{a,b},J})^{=_w1,et}$ is the locus where the kernel
of $p_2^{*} \cG \rightarrow p_1^{*} \cG$ is a multiplicative group
scheme; and the open and closed subspace $(X_{K',1}^{I'_{a,b},J})^{=_w1, 00}$ is the locus where the kernel of $p_2^{*} \cG \rightarrow p_1^{*} \cG$ is a bi-connected group scheme.

It follows from the definitions (see~\cite[Lem.\
7.4.2.4]{pilloniHidacomplexes}) that \numequation\label{eqn: 00 Hecke}p_2( (X_{K',1}^{I'_{a,b},J})^{=_w1, 00}) \subseteq (X_{K'',1}^{I'_{a,b},J})^{=_w1,00}\end{equation} and that at the level of topological spaces,
\numequation\label{eqn: met Hecke}p_2( (X_{K',1}^{I'_{a,b},J})^{=_w1, m} \cup (X_{K',1}^{I'_{a,b},J})^{=_w1,et}  ) \subseteq (X_{K'',1}^{I'_{a,b},J})^{=_w1,\met}.\end{equation}

Over $(X_{K,1}^{I'_{a,b},J})^{=_w1}$ and $(X_{K',1}^{I'_{a,b},J})^{=_w1}$ we can decompose the cohomological correspondences $T'_{w}$ and
$T''_{w}$ by projecting on the various components (in other words,
composing with the various idempotents associated to each of these
connected components). This gives us decompositions $T'_{w} = T'_{w,m} + T'_{w,et} + T'_{w,00}$
and $T''_{w} = T''_{w,m} + T''_{w,et} + T''_{w,00}$ obtained by
projecting on the multiplicative, \'etale and bi-connected components respectively.

\begin{lem}\label{lem: commutativity of etale Klingen correspondences}
  The following diagrams of  sheaves on $X_{K',1}^{I'_{a,b},J}$ are commutative for $l_w \geq p+1$  and $k_w \geq 2p+3$:
  \begin{eqnarray*}
\xymatrix{ p_2^{*} \omega^\kappa \vert_{ (X^{I'_{a,b},J}_{K'',1})^{=_w 1}} \ar[r]^{T'_{w,et}} \ar[d]^{ p_2^{*} {\Ha}'(\cG_w)}& p_1^! \omega^\kappa \vert_{(X^{I'_{a,b},J}_{K,1})^{=_w 1}}  \ar[d]^{p_1^{*} {\Ha}'(\cG_w)}\\
p_2^{*} (\omega^\kappa \otimes \det \omega_{\cG_w}^{p^2-1})\vert_{ (X^{I'_{a,b},J}_{K'',1})^{=_w 1}}  \ar[r]^{T'_{w,et}} & p_1^! (\omega^\kappa \otimes \det \omega_{\cG_w}^{p^2-1})  \vert_{(X^{I'_{a,b},J}_{K,1})^{=_w 1}}}
\end{eqnarray*}
\begin{eqnarray*}
\xymatrix{ p_1^{*} \omega^\kappa \vert_{ (X^{I'_{a,b},J}_{K,1})^{=_w 1}} \ar[r]^{T''_{w,et}} \ar[d]^{ p_1^{*} {\Ha}'(\cG_w)}& p_2^! \omega^\kappa \vert_{(X^{I'_{a,b},J}_{K'',1})^{=_w 1}}  \ar[d]^{p_2^{*} {\Ha}'(\cG_w)}\\
p_1^{*} (\omega^\kappa \otimes \det \omega_{\cG_w}^{p^2-1})\vert_{ (X^{I'_{a,b},J}_{K,1})^{=_w 1}}  \ar[r]^{T''_{w,et}} & p_2^! (\omega^\kappa \otimes \det \omega_{\cG_w}^{p^2-1})  \vert_{(X^{I'_{a,b},J}_{K'',1})^{=_w 1}}}
\end{eqnarray*}
 Moreover, $T'_{w,m} = T'_{w,00} = 0$ and $T''_{w,m} = 0$. If $l_w
 \geq p+2$, then $T''_{w,00} =0$.

\end{lem}

\begin{proof} See~\cite[Prop.\ 7.4.2.1]{pilloniHidacomplexes}.
\end{proof} 

We recall that by definition we have $T_w = T'_{w} \circ T''_{w} $ as operators on the cohomology over $X^{I'_{a,b},J}_K$. It
will also be useful to consider the composition 
$\tT_w:=T''_w\circ T'_w$ defining an operator on the cohomology over $X_{K''}^{I'_{a,b},J}$.

\begin{lem}\label{lem-commutativityTwrank1}If $l_w\ge p+1$ and $k_w\ge
  2p+3$, then we have $T_w = T'_{w, et } \circ T''_{w,et}  \in
  \mathrm{End} ( \mathrm{R} \Gamma (X^{I'_{a,b},J,=_w1}_{K,1}, \omega^\kappa))$, and \[T_w {\Ha}'(\cG_w) = {\Ha}'(\cG_w) T_w \in
  \mathrm{Hom} (\mathrm{R} \Gamma (X^{I'_{a,b},J,=_w1}_{K,1}, \omega^\kappa), \mathrm{R} \Gamma (X^{I'_{a,b},J,=_w1}_{K,1}, \omega^\kappa
  \otimes \det \omega_{\cG_w}^{(p^2-1)})).\]

Similarly, if  $k_w\ge
  2p+3$, then we have  $\tT_w = T''_{w, et } \circ T'_{w,et}  \in
  \mathrm{End} ( \mathrm{R} \Gamma (X^{I'_{a,b},J,=_w1}_{K'',1}, \omega^\kappa))$, and
   \[\tT_w {\Ha}'(\cG_w) = {\Ha}'(\cG_w) \tT_w \in
  \mathrm{Hom} (\mathrm{R} \Gamma (X^{I'_{a,b},J,=_w1}_{K'',1}, \omega^\kappa), \mathrm{R} \Gamma (X^{I'_{a,b},J,=_w1}_{K'',1}, \omega^\kappa
  \otimes \det \omega_{\cG_w}^{(p^2-1)})).\]\end{lem}

\begin{proof}We give the argument for~$T_w$; the argument for~$\tT_w$
  is essentially the same, and is left to the reader.
  From~(\ref{eqn: 00 Hecke}) and~(\ref{eqn: met Hecke}) we see that
  we can write $T_w$ as the sum of the two operators $T'_{w,00} \circ
  T''_{w,00}$ and~$(T'_{w,et}+T'_{w,m}) \circ (T''_{w,et}+T''_{w,m})$.

By Lemma~\ref{lem: commutativity of etale Klingen correspondences}, we
have $T'_{w,m} = T'_{w,00} = 0$ and
$T''_{w,m} = 0$, so that $T_w = T'_{w,et} \circ T''_{w,et}$. The
commutativity with~${\Ha}'(\cG_w)$ then follows from the
commutative diagrams in Lemma~\ref{lem: commutativity of etale Klingen correspondences}.
\end{proof} 

\begin{cor}\label{cor: T locally finite on w1}Assume that
  for all places~$v\in I'_b$, we have $l_v\geq 2$ and $k_v \geq C$,
  and that for all places~$v\in J$, we have $l_v\ge C$. If $l_w\ge p+1$ and $k_w\ge
  2p+3$, then the action of $T^{I_b,J}$ on $\mathrm{R}\Gamma (
  X^{G_1,I'_{a,b},J,=_w1}_{K,1} , \omega^\kappa)$ is locally finite. Similarly,
  if $k_w\ge
  2p+3$, then the action of  $T^{I'_b,J}\tT_w$  on $\mathrm{R}\Gamma (
  X^{G_1,I'_{a,b},J,=_w1}_{K'',1} , \omega^\kappa)$  is also locally finite.\end{cor}

\begin{proof} Again, we give the proof for~$T^{I_b,J}$, the argument
  for~$\tT_wT^{I',J}$ being essentially identical. We have ${\HH}^{*} (
  X^{G_1,I'_{a,b},J,=_w1}_{K,1} , \omega^\kappa) = \colim_{n} {\HH}^{*}
  (X^{G_1,I'_{a,b},J,\leq_w1}_{K,1} , \omega^\kappa \otimes \det
  \omega_{\cG_w}^{n(p^2-1)})$.  More precisely, for all $n \geq 0$, the map ${\HH}^{*}
  (X^{G_1,I'_{a,b},J,\leq_w1}_{K,1} , \omega^\kappa \otimes \det
  \omega_{\cG_w}^{n(p^2-1)}) \rightarrow {\HH}^{*} (
  X^{G_1,I'_{a,b},J,=_w1}_{K,1} , \omega^\kappa) $   is defined by the composition:
\[{\HH}^{*} (X^{G_1,I'_{a,b},J,\leq_w1}_{K,1} , \omega^\kappa \otimes \det \omega_{\cG_w}^{n(p^2-1)}) \rightarrow {\HH}^{*} (X^{G_1,I'_{a,b},J,=_w1}_{K,1} , \omega^\kappa \otimes \det \omega_{\cG_w}^{n(p^2-1)})\]\[ \stackrel{{\Ha}'(\cG_w)^{-n}} \rightarrow {\HH}^{*} (X^{G_1,I'_{a,b},J,=_w1}_{K,1} , \omega^\kappa).\]

The vector space  ${\HH}^{*} ( X^{G_1,I'_{a,b},J,=_w1}_{K,1} , \omega^\kappa)$ is therefore an inductive limit of  the images of the spaces  $${\HH}^{*} (X^{G_1,I'_{a,b},J,\leq_w1}_{K,1} , \omega^\kappa \otimes \det \omega_{\cG_w}^{n(p^2-1)})$$ under these maps.
The first map is obviously $T^{I_b,J}$-equivariant, and  by
Lemma~\ref{lem-commutativityTwrank1}, the second map is also
$T^{I_b,J}$-equivariant. The result follows from our inductive
hypothesis that Theorem~\ref{thm-Klingen-vanish} holds for all smaller
values of~$\# I_b$, because by definition we have
$X^{G_1,I'_{a,b},J,\leq_w1}_{K,1}=X^{G_1,\le_{I_a\cup\{w\}}1,\ge_{I'_b}1,\ge_J2}_{K,1}$. 
\end{proof}

\begin{lem} \label{lem-commutativityTw}The following diagram of
  sheaves on $X_{K',1}^{I'_{a,b},J}$ is commutative for $l_w \geq p+1$  and $k_w \geq 2p+3$:

\begin{eqnarray*}
\xymatrix{ p_2^{*} \omega^\kappa \vert_{ (X^{I'_{a,b},J}_{K'',1})^{\leq_w 1}} \ar[r]^{T'_{w}} \ar[d]^{ p_2^{*} {\Ha}'(\cG_w)}& p_1^! \omega^\kappa \vert_{(X^{I'_{a,b},J}_{K,1})^{\leq_w 1}}  \ar[d]^{p_1^{*} {\Ha}'(\cG_w)}\\
p_2^{*} (\omega^\kappa \otimes \det \omega_{\cG_w}^{p^2-1})\vert_{ (X^{I'_{a,b},J}_{K'',1})^{\leq_w 1}}  \ar[r]^{T'_w} & p_1^! (\omega^\kappa \otimes \det \omega_{\cG_w}^{p^2-1})  \vert_{(X^{I'_{a,b},J}_{K,1})^{\leq_w 1}}}
\end{eqnarray*}
If $l_w \geq p+2$ and $k_w \geq 2p+3$, the following diagram is commutative:  

\begin{eqnarray*}
\xymatrix{ p_1^{*} \omega^\kappa \vert_{ (X^{I'_{a,b},J}_{K',1})^{\leq_w 1}} \ar[r]^{T''_{w}} \ar[d]^{ p_1^{*} {\Ha}'(\cG_w)}& p_2^! \omega^\kappa \vert_{(X^{I'_{a,b},J}_{K'',1})^{\leq_w 1}}  \ar[d]^{p_2^{*} {\Ha}'(\cG_w)}\\
p_1^{*} (\omega^\kappa \otimes \det \omega_{\cG_w}^{p^2-1})\vert_{ (X^{I'_{a,b},J}_{K',1})^{\leq_w 1}}  \ar[r]^{T''_w} & p_2^! (\omega^\kappa \otimes \det \omega_{\cG_w}^{p^2-1})  \vert_{(X^{I'_{a,b},J}_{K'',1})^{\leq_w 1}}}
\end{eqnarray*}
\end{lem}

\begin{proof} It is enough to check the commutativity over a dense
  open subscheme of the support of these Cohen--Macaulay sheaves, and this follows from Lemma~\ref{lem: commutativity of etale Klingen correspondences}.
\end{proof} 

Since $p_1^{*} {\Ha}'(\cG_w)$ is not a zero divisor on
$X_{K',1}\vert_{(X^{I'_{a,b},J}_{K,1})^{\leq_w 1}}$, it follows exactly as in the
proof of
Lemma~\ref{lem: restricting Siegel correspondence} that there is for all $l_w \geq p^2+p$ and $k_w \geq p^2+2p+2$ a cohomological correspondence: 
$$T'_{w}: p_2^{*} (\omega^\kappa\vert_{(X^{I'_{a,b},J}_{K'',1})^{=_w 0}}) \rightarrow p_1^! (\omega^\kappa\vert_{(X^{I'_{a,b},J}_{K,1})^{=_w 0}})$$
 and similarly for all $l_w \geq p^2+p+1$ and $k_w \geq p^2+2p+2$ a cohomological correspondence: 
 $$T''_{w}: p_1^{*} (\omega^\kappa\vert_{(X^{I'_{a,b},J}_{K,1})^{=_w 0}}) \rightarrow p_2^! (\omega^\kappa\vert_{(X^{I'_{a,b},J}_{K'',1})^{=_w 0}}).$$

\begin{lem}\label{lem: constant C so that T prime vanishes on w equals
    0} There is a universal constant $C'$ which depends only on
  $F$ and $p$ but not on the tame level such that \[T'_{w}: p_2^{*} (\omega^\kappa\vert_{(X^{I'_{a,b},J}_{K'',1})^{=_w 0}}) \rightarrow p_1^! (\omega^\kappa\vert_{(X^{I'_{a,b},J}_{K,1})^{=_w 0}})\]  is zero for all $k_w-l_w \geq C'$ and all $l_w \geq p^2+p$. 
\end{lem} 

\begin{proof} See~\cite[Prop.\ 7.4.2.2]{pilloniHidacomplexes}.\end{proof}
We now increase our constant~$C$ if necessary, so that $C\ge C'$,
where~$C'$ is as in Lemma~\ref{lem: constant C so that T prime vanishes on w equals
    0}.
\begin{lem}\label{lem: quasi isomorphism to equals w 1} Assume that
  for all places~$v\in I'_b$, we have $l_v \geq 2$ and $k_v-l_v \geq C$, that $l_v \geq C$ for
  all $v \in J$, and that $l_w \geq p+2 $ and $k_w-l_w \geq C$. Then the map $ e(T^{I_b,J})\mathrm{R}\Gamma(X^{G_1,I'_{a,b},J,\leq_w 1}_{K,1}, \omega^{\kappa}) \rightarrow e(T^{I_b,J}) \mathrm{R}\Gamma(X^{G_1,I'_{a,b},J,=_w 1}_{K,1}, \omega^{\kappa})$ is a quasi-isomorphism. In particular, $e(T^{I_b,J}) \mathrm{R}\Gamma(X^{G_1,I'_{a,b},J,=_w 1}_{K,1}, \omega^{\kappa})$ is a perfect complex. 
\end{lem}

\begin{proof}
Consider the following diagram of exact triangles:
\begin{eqnarray*}
\xymatrix{  
\mathrm{R}\Gamma(X^{G_1,I'_{a,b},J,\leq_w 1}_{K'',1}, \omega^{\kappa}) 
\ar[r]^{T'_{w}} \ar[d]^{{\Ha}'(\cG_w)} &  \mathrm{R}\Gamma(X^{G_1,I'_{a,b},J, \leq_w 1}_{K,1}, \omega^{\kappa})
 \ar[d]^{{\Ha}'(\cG_w)} \\
  \mathrm{R} \Gamma (X^{G_1,I'_{a,b},J, \leq_w 1} _{K'',1}, \omega^\kappa \otimes \det \omega_{\cG_w}^{p^2-1} ) 
  \ar[r]^{T'_{w}} \ar[d] &   \mathrm{R} \Gamma (X^{G_1,I'_{a,b},J, \leq_w 1} _{K,1}, \omega^\kappa \otimes \det \omega_{\cG_w}^{p^2-1} ) 
  \ar[d] \\
  \mathrm{R} \Gamma (X^{G_1,I'_{a,b},J,=_w 0}_{K'',1}, \omega^\kappa \otimes \det \omega_{\cG_w}^{p^2-1} )
\ar[r]^{T'_{w}} & 
  \mathrm{R} \Gamma (X^{G_1,I'_{a,b},J,=_w 0}_{K,1}, \omega^\kappa \otimes \det \omega_{\cG_w}^{p^2-1} ) }
  \end{eqnarray*}

By Lemma~\ref{lem: constant C so that T prime vanishes on w equals
    0}, the  rightmost operator $T'_{w}$ acts by zero. We have the
  ordinary projectors~$e(T^{I_b,J})$ on $ \mathrm{R}\Gamma(X^{G_1,I'_{a,b},J,\leq_w 1}_{K,1}, \omega^{\kappa})
  $ and $\mathrm{R} \Gamma (X^{G_1,I'_{a,b},J,\leq_w 1} _{K,1}, \omega^\kappa \otimes
  \det \omega_{\cG_w}^{p^2-1} )$, and the
  ordinary projectors~$e(\tT_wT^{I'_b,J})$ on $ \mathrm{R}\Gamma(X^{G_1,I'_{a,b},J,\leq_w 1}_{K'',1}, \omega^{\kappa})
  $ and $\mathrm{R} \Gamma (X^{G_1,I'_{a,b},J,\leq_w 1} _{K'',1}, \omega^\kappa
  \otimes \det \omega_{\cG_w}^{p^2-1} )$. It follows from the defining
  properties of the ordinary projectors that after applying them, the
  left two vertical arrows~$T'_w$ are quasi-isomorphisms. 

By Lemma~\ref{lem-commutativityTw} the projectors commute with
multiplication by ${\Ha}'(\cG_w)$. It follows from a short
diagram chase that the map $$e(T^{I_b,J})\mathrm{R}\Gamma(X^{G_1,I'_{a,b},J, \leq_w
  1}_{K,1}, \omega^{\kappa}) \stackrel{{\Ha}'(\cG_w)}\rightarrow
e(T^{I_b,J}) \mathrm{R} \Gamma (X^{G_1,I'_{a,b},J, \leq_w 1} _{K,1}, \omega^\kappa \otimes
\det \omega_{\cG_w}^{p^2-1} ) $$ is a quasi-isomorphism. The claimed quasi-isomorphism
now follows by taking an inductive limit as in the proof of
Corollary~\ref{cor: T locally finite on w1}. By our inductive hypothesis, $e(T^{I_b,J})\mathrm{R}\Gamma(X^{G_1,I'_{a,b},J,\leq_w 1}_{K,1},
\omega^{\kappa})$ is a perfect complex, so we are done.\end{proof}

\begin{lem}\label{lem-isomTw}Assume that
  for all places~$v\in I$, we have $l_v \geq 2$ and $k_v-l_v \geq C$, that $l_v \geq C$ for
  all $v \in J$, and that $l_w=p+1$.  Then the ordinary cohomology
  $e(T^{I_b,J}) \mathrm{R}\Gamma(X^{G_1,I'_{a,b},J,=_w 1}_{K,1}, \omega^{\kappa})$  is a
  perfect complex, and the map $$ e(T^{I_b,J}){\HH}^0(X^{G_1,I'_{a,b},J,\leq_w
    1}_{K,1}, \omega^{\kappa}) \rightarrow e(T^{I_b,J}) {\HH}^0(X^{G_1,I'_{a,b},J,=_w
    1}_{K,1}, \omega^{\kappa})$$ is an isomorphism. 
\end{lem}

\begin{proof} The map  $$\mathrm{R}\Gamma(X^{G_1,I'_{a,b},J,=_w 1}_{K,1},
  \omega^{\kappa}) \stackrel{{\Ha}'(\cG_w)}\rightarrow
  \mathrm{R}\Gamma(X^{G_1,I'_{a,b},J,=_w 1}_{K,1}, \omega^{\kappa} \otimes \det
  \omega_{\cG_w}^{p^2-1})$$ is a quasi-isomorphism, which commutes
  with the projector $e(T^{I_b,J})$ by Lemma~\ref{lem-commutativityTwrank1}.
  It
  follows from Lemma~\ref{lem: quasi isomorphism to equals w 1} that
 $e(T^{I_b,J}) \mathrm{R}\Gamma(X^{G_1,I'_{a,b},J,=_w 1}_{K,1}, \omega^{\kappa})$ is perfect. 

 We now prove the claimed isomorphism on degree $0$ cohomology. Since $X^{G_1,I'_{a,b},J,\leq_w 1}_{K,1}$
 is Cohen--Macaulay, and $X^{G_1,I'_{a,b},J,=_w 1}_{K,1}$ is an open
 dense subscheme,  we have
 injections \[{\HH}^0(X^{G_1,I'_{a,b},J,\leq_w 1}_{K,1},
 \mathcal{F} )\hookrightarrow {\HH}^0(X^{G_1,I'_{a,b},J,=_w 1}_{K,1},
 \mathcal{F})\] for any locally free sheaf $\mathcal{F}$,  so it is enough to prove surjectivity. 
In order to do this, it is enough to prove that for all $n\geq 0$, the
map $${\HH}^0(X^{G_1,I'_{a,b},J,\leq_w 1}_{K,1}, \omega^{\kappa})
\stackrel{({\Ha}'(\cG_w))^n}\rightarrow {\HH}^0(X^{G_1,I'_{a,b},J, \leq_w
  1} _{K,1}, \omega^\kappa \otimes \det \omega_{\cG_w}^{n(p^2-1)} )$$
(which commutes with~$e(T_w)$ by Lemma~\ref{lem-commutativityTwrank1}
and the injectivity of the restrictions ${\HH}^0(X^{G_1,I'_{a,b},J, \leq_w 1}_{K,1},
 \mathcal{F} )\hookrightarrow {\HH}^0(X^{G_1,I'_{a,b},J,=_w 1}_{K,1},
 \mathcal{F})$ discussed above) induces a surjection: 
$$e(T^{I_b,J}){\HH}^0(X^{G_1,I'_{a,b},J,\leq_w 1}_{K,1}, \omega^{\kappa}) \stackrel{({\Ha}'(\cG_w))^n}\rightarrow e(T^{I_b,J}){\HH}^0(X^{G_1,I'_{a,b},J,\leq_w 1} _{K,1}, \omega^\kappa \otimes \det \omega_{\cG_w}^{n(p^2-1)} ).$$

In fact, by Lemma~\ref{lem: quasi isomorphism to equals w 1}, it
suffices to prove the surjectivity for $n=1$. We consider the
following diagram:

\begin{eqnarray*}
\xymatrix{  
{\HH}^0(X^{G_1,I'_{a,b},J, \leq_w 1}_{K'',1}, \omega^{\kappa})  \ar[r]^{T'_{w}} \ar[d]^{{\Ha}'(\cG_w)} &
{\HH}^0(X^{G_1,I'_{a,b},J,\leq_w 1}_{K,1}, \omega^{\kappa})   \ar[d]^{{\Ha}'(\cG_w)}   \\
 {\HH}^0 (X^{G_1,I'_{a,b},J, \leq_w 1} _{K'',1}, \omega^\kappa \otimes \det \omega_{\cG_w}^{p^2-1} )   \ar[r]^{T'_{w}} \ar[d]  &
{\HH}^0(X^{G_1,I'_{a,b},J, \leq_w 1} _{K,1}, \omega^\kappa \otimes \det \omega_{\cG_w}^{p^2-1} )   \ar[d]  \\
{\HH}^0 (X^{G_1,I'_{a,b},J, =_w 0}_{K'',1}, \omega^\kappa \otimes \det \omega_{\cG_w}^{p^2-1} )    \ar[r]^{T'_{w}}  &
{\HH}^0 (X^{G_1,I'_{a,b},J, =_w 0}_{K,1}, \omega^\kappa \otimes \det \omega_{\cG_w}^{p^2-1} )} 
\end{eqnarray*}

Let $f \in e(T^{I_b,J}){\HH}^0(X^{G_1,I'_{a,b},J, \leq_w 1} _{K,1}, \omega^\kappa
\otimes \det \omega_{\cG_w}^{p^2-1} )$. As in the proof of
Lemma~\ref{lem: quasi isomorphism to equals w 1}, $T_w'$ induces a
bijection
 \[e(\tT_wT^{I'_b,J}){\HH}^0(X^{G_1,I'_{a,b},J, \leq_w 1} _{K'',1}, \omega^\kappa
\otimes \det \omega_{\cG_w}^{p^2-1} )  \to  
 e(T^{I_b,J}){\HH}^0(X^{G_1,I'_{a,b},J, \leq_w 1} _{K,1}, \omega^\kappa
\otimes \det \omega_{\cG_w}^{p^2-1} ). \]

 In particular, $f = T'_w g$ for some $g \in {\HH}^0(X^{G_1,I'_{a,b},J, \leq_w 1}
 _{K,1}, \omega^\kappa \otimes \det \omega_{\cG_w}^{p^2-1} )$ and
 therefore,  since the rightmost operator $T_w'$ acts by zero  by Lemma~\ref{lem: constant C so that T prime vanishes on w equals
    0}, $f$ has trivial image in  ${\HH}^0 (X^{G_1,I'_{a,b},J, =_w 0}_{K,1},
 \omega^\kappa \otimes \det \omega_{\cG_w}^{p^2-1} ) $. It follows
 that~$f$ comes from a class $\tilde{f} \in {\HH}^0(X^{G_1,I'_{a,b},J, \leq_w
   1}_{K,1}, \omega^{\kappa})$. Replacing $\tilde{f}$ by $e(T^{I_b,J})
 \tilde{f}$ we deduce the required surjectivity. \end{proof}

\begin{cor}\label{cor: Klingen quasi iso l at least 3} Assume that
  for all places~$v\in I_b$, we have $l_v \geq 2$ and $k_v-l_v\ge C$,
  and that for all places~$v\in J$, we have $l_v\ge C$. 
Then~$T^{I_b,J}$ acts locally finitely on~$\mathrm{R}\Gamma(X_{1}^{ G_1,I_{a,b},J }, \omega^{\kappa})$  and $e(T^{I_b,J})\mathrm{R}\Gamma(X_{1}^{ G_1, I_{a,b},J}, \omega^{\kappa})$ is a perfect complex.

  We have an isomorphism: $$ e(T^{I_b,J})
 {\HH}^0(X^{G_1,\le_{I_a}1}_1, \omega^\kappa) \isoto e(T^{I_b,J})
 {\HH}^0(X_1^{G_1,I_{a,b},J}, \omega^\kappa)$$ and an injection: 
$$ e(T^{I_b,J}) {\HH}^1(X^{G_1,\le_{I_a}1}_1, \omega^\kappa) \into e(T^{I_b,J}) {\HH}^1(X_1^{G_1,I_{a,b},J}, \omega^\kappa).$$

If furthermore~$l_w\ge 3$, then the map $$
  e(T^{I_b,J})\mathrm{R}\Gamma(X^{G_1,I'_{a,b},J}_{1}, \omega^{\kappa}) \rightarrow
  e(T^{I_b,J})\mathrm{R}\Gamma(X_{1}^{ G_1, I_{a,b},J}, \omega^{\kappa})$$is a quasi-isomorphism. 

    \end{cor} 

\begin{proof} We begin by showing that~$T^{I_b,J}$ acts locally finitely on both ~$\mathrm{R}\Gamma(X_1^{G_1,I'_{a,b},J, \geq_w
    2}, \omega^{\kappa})$ and~$\mathrm{R}\Gamma(X_1^{ G_1,I_{a,b},J}, \omega^{\kappa})=\mathrm{R}\Gamma(X_1^{ G_1,I'_{a,b},J, \geq_w
    1}, \omega^{\kappa})$. For $\mathrm{R}\Gamma(X_1^{G_1,I'_{a,b},J, \geq_w
    2}, \omega^{\kappa})$, this  is Corollary  \ref{cor: Tw locally finite on ordinary}.

Our argument for $\mathrm{R} \Gamma ( X^{G_1,I'_{a,b},J,\geq_w 1} _1, \omega^\kappa)
$ is slightly more involved. We have an exact triangle
$$\begin{diagram}
 \mathrm{R} \Gamma ( X^{G_1,I'_{a,b},J, \geq_w 1} _1, \omega^\kappa) & \rTo &
  \mathrm{R} \Gamma ( X^{G_1,I'_{a,b},J, \geq_w 2}_1, \omega^\kappa) \\
  & & \dTo \\
  & &    \mathrm{R}
  \Gamma_{X^{G_1,I'_{a,b},J,=_w 1}} ( X^{G_1,I'_{a,b},J, \geq_w1}_1, \omega^\kappa)  [+1], \end{diagram}$$

  so it is
enough to prove that the action of~$T^{I_b,J}$ on $\mathrm{R}
  \Gamma_{X_1^{G_1,I'_{a,b},J,=_w 1}} ( X^{G_1,I'_{a,b},J, \geq_w1}_1, \omega^\kappa)  [+1]$ is locally
  finite. We have 
  \numequation\label{eqn: relative cohomology le 1 as colimit}
  \begin{diagram}
\mathrm{R} \Gamma_{X_1^{G_1,I'_{a,b},J, =_w 1}} ( X_1^{G_1,I'_{a,b},J, \geq_w 1}, \omega^\kappa)[+1] \\
\dEquals \\ 
\mathrm{R} \Gamma( X_1^{G_1,I'_{a,b},J, \geq_w 1},  \colim_n \omega^{\kappa} \otimes
\det \omega_{\cG_w}^{n(p-1)}|_{V({\Ha}(\cG_w)^{n})}) 
\end{diagram}
\end{equation} so
it is enough to prove that the action of~$T^{I_b,J}$ on each $\mathrm{R} \Gamma( X_1^{G_1,I'_{a,b},J, \ge_w 1}, \omega^{\kappa} \otimes
\det \omega_{\cG_w}^{n(p-1)}|_{V({\Ha}(\cG_w)^{n})})$ is locally
finite. In the case $n=1$, this is Corollary~\ref{cor: T locally
  finite on w1}, and the general case follows by induction by taking
the cohomology of the
short exact sequence of sheaves
\numequation\label{eqn: restricting sheaves to vanishing of powers of Hasse} 
0 \rightarrow  \omega^{\kappa} \otimes \det
\omega_{\cG_w}^{(n-1)(p-1)}|_{V({\Ha}(\cG_w)^{n-1})}
\stackrel{ \times {\Ha}(\cG_w)}\rightarrow\omega^{\kappa} \otimes \det
\omega_{\cG_w}^{n(p-1)}|_{V({\Ha}(\cG_w)^{n})} \end{equation} \[\rightarrow
\omega^{\kappa} \otimes \det
\omega_{\cG_w}^{n(p-1)}|_{V({\Ha}(\cG_w))} \rightarrow 0 \]

Consider now the following diagram of exact triangles:

\begin{eqnarray*}
\xymatrix{ 
 \mathrm{R} \Gamma ( X^{G_1,I'_{a,b},J}_1, \omega^\kappa) \ar[r] \ar[d]  & \mathrm{R} \Gamma ( X^{G_1,I'_{a,b},J,\geq_w 1} _1, \omega^\kappa) 
 \ar[d] \\
  \mathrm{R} \Gamma ( X^{G_1,I'_{a,b},J, \geq_w 2}_1, \omega^\kappa) \ar[r] \ar[d] &  
  \mathrm{R} \Gamma ( X^{G_1,I'_{a,b},J, \geq_w 2}_1, \omega^\kappa) \ar[d] \\
\mathrm{R} \Gamma_{X_1^{G_1,I'_{a,b},J, \leq_w 1}} ( X^{G_1,I'_{a,b},J}_1, \omega^\kappa)  [+1] \ar[r] & \mathrm{R} \Gamma_{X_1^{G_1,I'_{a,b},J, =_w 1}} ( X^{G_1,I'_{a,b},J, \geq_w1}_1, \omega^\kappa)  [+1]
   }
   \end{eqnarray*}

We have already seen that~$T^{I_b,J}$
acts locally finitely on all but the last term of the first row, so it
acts locally finitely on every term in the diagram. The middle vertical arrow is the identity map.  In order to show that
the left vertical arrow is a quasi-isomorphism after applying
$e(T^{I_b,J})$, it is therefore enough to prove it for the right vertical
arrow.  By~(\ref{eqn: relative cohomology le 1 as colimit}) and the
similar expression
 $$\mathrm{R} \Gamma_{X_1^{G_1,I'_{a,b},J, \leq_w 1}} ( X^{G_1,I'_{a,b},J}_1, \omega^\kappa)[+1]=
\mathrm{R} \Gamma( X_1^{G_1,I'_{a,b},J},  \colim_n \omega^{\kappa} \otimes
\det \omega_{\cG_w}^{n(p-1)}|_{V({\Ha}(\cG_w)^{n})}),$$it suffices to show that for each~$n\ge 1$, the map 

 $$  \begin{diagram}
 e(T^{I_b,J})\mathrm{R}
\Gamma(X^{G_1}_1, \omega^{\kappa} \otimes \det
\omega_{\cG_w}^{n(p-1)}|_{V({\Ha}(\cG_w)^{n})}) \\ \dTo \\  e(T^{I_b,J})
\mathrm{R} \Gamma( X_1^{G_1,I'_{a,b},J, \geq_w 1},  \omega^{\kappa} \otimes
\det \omega_{\cG_w}^{n(p-1)}|_{V({\Ha}(\cG_w)^{n})}) \end{diagram}$$
is a
quasi-isomorphism. To see this, note that the case $n=1$ is Lemma~\ref{lem: quasi isomorphism to equals w 1}, and the general case follows by induction on~$n$, by taking the
cohomology of the exact sequence of sheaves~(\ref{eqn: restricting
  sheaves to vanishing of powers of Hasse}). The remaining claims  follow
from  the quasi-isomorphism just proved and the inductive hypothesis.
\end{proof}

\begin{cor}\label{coro-classicality-Tw-rw3}Assume that
  for all places~$v\in I$, we have  $k_v-l_v \geq C$ and $l_v \geq 2$,
  that for all places~$v\in J$, we have $l_v\ge C$.
Then the ordinary cohomology  $$ e(T^{I_b,J})\mathrm{R}\Gamma(X_1^{G_1,I_{a,b},J}, \omega^{\kappa})$$ is represented by a perfect complex.
Moreover we have an isomorphism: $$ e(T^{I_b,J}) {\HH}^0(X^{G_1,\le_{I_a}1}_1,
\omega^\kappa) \isoto e(T^{I_b,J}) {\HH}^0(X_1^{G_1,I_{a,b},J},
\omega^\kappa)$$ and an injection: 
$$ e(T^{I_b,J}) {\HH}^1(X^{G_1,\le_{I_a}1}_1, \omega^\kappa) \into e(T^{I_b,J}) {\HH}^1(X_1^{G_1,I_{a,b},J}, \omega^\kappa).$$
\end{cor}

\begin{proof}To see that  $$ e(T^{I_b,J})\mathrm{R}\Gamma(X_1^{G_1,I_{a,b},J}, \omega^{\kappa})$$ is represented by a perfect complex, we
consider the exact triangle: $$\mathrm{R} \Gamma(X^{G_1,I'_{a,b},J, \geq_w 1} _1, \omega^{\kappa}) \stackrel{{\Ha}(\cG_w)} \rightarrow \mathrm{R} \Gamma(X^{G_1,I'_{a,b},J, \geq_w 1} _1, \omega^{\kappa} \otimes \det
\omega_{\cG_w}^{p-1})$$ $$ \rightarrow \mathrm{R} \Gamma(X^{G_1,I'_{a,b},J,=_w1}_1, \omega^{\kappa} \otimes \det \omega_{\cG_w}^{p-1}) \stackrel{+1}\rightarrow $$
Applying the projector $e(T^{I_b,J})$ everywhere (which commutes with the
various maps by Lemma~\ref{lem-commutTW}) we deduce this from
Corollary~\ref{cor: Klingen quasi iso l at least 3} and
Lemma~\ref{lem-isomTw}. By our inductive hypothesis, in order to prove the claims about the
morphisms \[e(T^{I_b,J}) {\HH}^i(X^{G_1,\le_{I_a}1}_1,
\omega^\kappa) \isoto e(T^{I_b,J}) {\HH}^i(X_1^{G_1,I_{a,b},J},
\omega^\kappa),\] it is enough to prove the corresponding statements
for the morphisms 
$$e(T^{I_b,J}) {\HH}^i(X^{G_1,I'_{a,b},J}_1,
\omega^\kappa) \to e(T^{I_b,J}) {\HH}^i(X_1^{G_1,I_{a,b},J},
\omega^\kappa).$$
Firstly, the natural restriction map
\[ {\HH}^0(X^{G_1,I'_{a,b},J}_1, \omega^\kappa) \rightarrow
  {\HH}^0(X_1^{G_1,I_{a,b},J}, \omega^\kappa)\]is an isomorphism,
because $X_1^{G_1,I'_{a,b},J}$ is Cohen--Macaulay, and the complement of $X_1^{G_1,I_{a,b},J}$ is of codimension at least
$2$.  It remains to prove the injectivity of the
map of~$H^1$s. 

We have a commutative diagram of exact triangles:

\begin{eqnarray*}
\xymatrix{ 
 \mathrm{R} \Gamma ( X^{G_1,I'_{a,b},J}_1, \omega^\kappa) \ar[r] \ar[d] &
 \mathrm{R} \Gamma ( X^{G_1,I_{a,b},J}_1, \omega^\kappa) \ar[d] \\
 \mathrm{R} \Gamma ( X^{G_1,I'_{a,b},J}_1, \omega^\kappa \otimes \det \omega_{\cG_w}^{(p-1)}) \ar[r] \ar[d] &
\mathrm{R} \Gamma ( X^{G_1,I_{a,b},J}_1, \omega^\kappa \otimes \det \omega_{\cG_w}^{(p-1)}) \ar[d] \\
  \mathrm{R} \Gamma ( X^{G_1,I'_{a,b},J, \leq_w 1}_1, \omega^\kappa \otimes \det \omega_{\cG_w}^{(p-1)}) \ar[r]   & 
\mathrm{R} \Gamma ( X^{G_1,I'_{a,b},J, =_w 1}_1, \omega^\kappa \otimes \det \omega_{\cG_w}^{(p-1)})  }
\end{eqnarray*}

The injectivity of $ e(T^{I_b,J}) {\HH}^1(X^{G_1,I'_{a,b},J}_1, \omega^\kappa)
\rightarrow e(T^{I_b,J}) {\HH}^1(X_1^{G_1,I_{a,b},J}, \omega^\kappa)$ 
therefore follows from a short diagram chase, using the
quasi-isomorphisms provided by Corollary~\ref{cor: Klingen
  quasi iso l at least 3} and the isomorphism on~$\HH^0$s of Lemma~\ref{lem-isomTw}.\end{proof}

\begin{proof}[Proof of Theorem~\ref{thm-Klingen-vanish}]This is immediate from Corollaries~\ref{cor: Klingen quasi iso l at least 3} and~\ref{coro-classicality-Tw-rw3}.\end{proof}

 \subsubsection{A Cousin complex computing $\mathrm{R}\Gamma(X^{G_1, I}_1,
   \omega^{\kappa}(-D))$ }\label{section-computingcohoexplicitely} Our
 goal in this section is to  provide an explicit Hecke stable complex
 computing $\mathrm{R}\Gamma(X^{G_1, I}_1, \omega^{\kappa}(-D))$. This
 complex will be used to complete the proof of
 Theorem~\ref{thm-vanishing}, and will also be used in~\S\ref{subsec: perfect Hida
   complex} to compare the cohomology at spherical and
 Klingen levels (by considering the corresponding complex at Klingen
 level, and the natural map between these complexes).   
 This complex is the Cousin complex associated to $X^{G_1, I}_1$ and the stratification given by the $p$-rank (see 
 \S\ref{subsec: cousin complexes}).    This section is very similar to \S\ref{sect-cousin-compEO} where we introduced the Cousin complex over the full Shimura variety associated with the Ekedahl--Oort stratification. The case we consider here is, however, much simpler because we have  canonical global equations provided by the partial Hasse invariants for our stratification. We has thus decided, despite redundancy,  to give a complete  and explicit construction of the Cousin complex in this case.  

Let $S$ be  a smooth scheme over a field $k$ and let $\cL$ be an
invertible sheaf on $S$. We assume that $\cL = \otimes_{i=1}^d \cL_i$
and that we have non-vanishing sections $s_i \in \HH^0(S, \cL_i)$. We
let $D_i = V(s_i)$, an effective Cartier divisor on $S$.  Set $ s =
\prod_{i=1}^d s_i$. Set $D = V(s) = \cup_i D_i$. We assume that $D =
\cup_i D_i$ is a strict normal crossing divisor on $S$. 

 For all $n$, consider the following exact complex   of coherent sheaves on $S$:
$$ 0 \rightarrow \oscr_{S} \stackrel{ s^n} \rightarrow \cL^n \rightarrow \bigoplus_{i = 1}^d \cL^n /s_i^n \rightarrow \bigoplus_{1 \leq i < j \leq d} \cL^n/(s_i^n, s_j^n)  \rightarrow \cdots \rightarrow \cL^n/(s_1^n, \cdots, s_d^n) \rightarrow 0$$
This is a complex of length $d + 2$. For all $0 \leq k \leq d$, the object placed in degree $k+1$ is $\bigoplus_{ 1 \leq i_1 < \cdots < i_k \leq d} \cL^n/(s_{i_1}^n, \cdots, s_{i_k}^n)$ (when we write $\cL^n/(s_{i_1}^n, \cdots, s_{i_k}^n)$, we  mean $\cL^n/ \cL^n(s_{i_1}^n\cL_{i_1}^{-n}, \cdots, s_{i_k}^n\cL_{i_k}^{-n})$).  The differential 
\[\bigoplus_{ 1 \leq i_1 < \cdots < i_k \leq d} \cL^n/(s_{i_1}^n, \cdots, s_{i_k}^n)\to \bigoplus_{ 1 \leq i_1 < \cdots < i_{k+1} \leq d} \cL^n/(s_{i_1}^n, \cdots, s_{i_k}^n)\]  takes a section $ (f_{i_1, \cdots, i_k})_{1 \leq i_1 < \cdots < i_k \leq d}$ to the section $$(\sum (-1)^{i_j} \overline{f_{i_1, \cdots, \widehat{i}_j, \cdots, i_{k+1}}})_{1 \leq i_1 < \cdots < i_{k+1} \leq d}$$
where  $\overline{f_{i_1, \cdots, \widehat{i}_j, \cdots, i_{k+1}}}$ is the class modulo $s_j^n$ of $f_{i_1, \cdots, \widehat{i}_j, \cdots, i_{k+1}}$.

The following diagram is commutative:
\begin{eqnarray*}
\xymatrix{ 0 \ar[r] & \oscr_{S} \ar[r]^{s^{n+1}} & \cL^{n+1} \ar[r] & \bigoplus_{i = 1}^d \cL^{n+1} /s_i^{n+1}  \ar[r] & \cdots \\
 0 \ar[r] & \oscr_{S}\ar[u]^{\id} \ar[r]^{s^n} & \cL^{n} \ar[r]\ar[u]^{s} & \bigoplus_{i = 1}^d \cL^n /s_i^n \ar[r]\ar[u]^{s} & \cdots}
 \end{eqnarray*}
 
Passing to the limit over $n$, we get the following exact complex: 
 $$ 0 \rightarrow \oscr_{S}  \rightarrow \colim_n \cL^n \rightarrow  \colim_n \bigoplus_{i = 1}^d \cL^n /s_i^n \rightarrow $$
 $$ \colim_n \bigoplus_{1 \leq i < j \leq d} \cL^n/(s_i^n, s_j^n)  \rightarrow \cdots \rightarrow  \colim_n \cL^n/(s_1^n, \cdots, s_d^n) \rightarrow 0$$
where in all the direct limits, the transition maps  are given by multiplication by powers of  $s$. 

\begin{lem}\label{lem: L versus L prime}  Let $1 \leq i_1 < \cdots < i_k \leq d$.  Set $\cL' = \otimes_{j= 1}^k \cL_{i_j}$, $s' = \prod_{j=1}^k s_{i_j}$, $D' = V(s(s')^{-1})$.  There is a canonical isomorphism 

$$ \colim_{ \times s^n}   \cL^n/(s_{i_1}^n, \cdots, s_{i_k}^n) \simeq  \colim_{ \times (s')^n} ((\cL')^n/ (s_{i_1}^n, \cdots, s_{i_k}^n))\vert_{S \setminus D'}$$
\end{lem}
\begin{proof} Easy and left to the reader. \end{proof} 

\begin{rem} The complex $$0 \rightarrow  \colim_n \cL^n \rightarrow  \colim_n \bigoplus_{i = 1}^d \cL^n /s_i^n \rightarrow $$
 $$ \colim_n \bigoplus_{1 \leq i < j \leq d} \cL^n/(s_i^n, s_j^n)  \rightarrow \cdots \rightarrow  \colim_n \cL^n/(s_1^n, \cdots, s_d^n) \rightarrow 0$$ is just the Cousin complex of $\ocal_S$ associated with the stratification given by the divisors $D_i$.
 \end{rem}

We now work over $S= X^{G_1, I}_1$. We take $\cL = \otimes_{w \in I}
(\det \cG_w)^{p-1}$, $\cL_w = (\det \cG_w)^{p-1}$ and  $s_w =
\Ha(\cG_w)$,  and we consider the complex 
$ K^0 \rightarrow K^1 \cdots \rightarrow K^d$ obtained by applying  $\HH^0$ to  the complex $$\colim_n \cL^n \rightarrow  \colim_n \bigoplus_{i = 1}^d \cL^n /s_i^n \rightarrow 
  \colim_n \bigoplus_{1 \leq i < j \leq d} \cL^n/(s_i^n, s_j^n)
  \rightarrow \cdots \rightarrow  \colim_n \cL^n/(s_1^n, \cdots,
  s_d^n)$$ tensored with $\omega^\kappa(-D)$. (So in the above notation, the indices~$i$ will correspond to the different
places~$v \in I$,  the~$s_i$ will correspond to Hasse invariants, and the
assumption that the divisor~$V(s)$ has strict normal crossings is an
easy consequence of the Serre--Tate theorem and the product structure on the
$p$-divisible group~$\cG$.)

 It follows from Lemma~\ref{lem: L versus L prime} that $K^k$ equals $$  
 \textstyle{\bigoplus_{ J \subset I, \# J = k}} \displaystyle{ \colim_{\times \prod_{w \in J} \Ha(\cG_w)}} \HH^0\big(X^{ G_1,I, \geq_{ J^c} 2}_1, \omega^{\kappa} (-D)\otimes \bigotimes_{w \in J} (\det \cG_w)^{n(p-1)}/(\sum_{w\in J} (\Ha(\cG_w)^n))\big).$$

\begin{prop}\label{prop: the complex K} The complex $K^\bullet$ computes $\mathrm{R}\Gamma(X^{G_1, I}_1, \omega^{\kappa}(-D))$. 
\end{prop} 
\begin{proof} The argument is the same as  in the proof of  Proposition~ \ref{prop-lastone}. It suffices to show that 
 each of the sheaves  $$ \omega^{\kappa} (-D)\otimes \bigotimes_{w \in J} (\det \cG_w)^{n(p-1)}/(\sum_{w\in J} (\Ha(\cG_w)^n))$$  
 when restricted to~$X^{ G_1,I, \geq_{ J^c} 2}_1$ and then pushed
 forward to~$X^{G_1, I}_1$ is acyclic on $X^{G_1, I}_1$. Since the
 inclusion $X^{ G_1,I, \geq_{ J^c} 2}_1\into X^{G_1, I}_1$ is affine, it suffices to show that the restriction of this sheaf to~$X^{ G_1,I, \geq_{ J^c} 2}_1$ is acyclic.
 By \cite[Thm.\ 8.6]{Lan2016}, this sheaf is acyclic relative to the minimal compactification and its support in the minimal compactification is the locally closed subscheme given by the set of equations: 
\begin{itemize}
\item ${\Ha}(\cG_w)^n = 0$ for $w \in J$,
\item ${\Ha}'(\cG_w) \neq 0$ for $w \in J$,
\item ${\Ha}(\cG_w) \neq 0$ for $w \in J^c$,
\end{itemize}
 which is affine. 
\end{proof}

\begin{rem}
  \label{rem: could use this complex to prove local finiteness of
    individual Hecke operators}Using Proposition~\ref{prop: the
    complex K}, one can show that if we have $l_w \geq 2$ and
  $k_w \geq 2p+3$ for $w \in I$, and $l_w \geq 3$ for $w \in I^c$,
  then the individual Hecke operators~$T_w$ for~$w\in I$ and~$T_{w,1}$
  for~$w\in I^c$ act locally finitely
  on~$R\Gamma(X_1^{G_1,I},\omega^\kappa(-D))$ (by showing that they act
  locally finitely on each term of the complex~$K^\bullet$). We leave
  the details to the interested reader.
\end{rem}

We can finally complete the proof of Theorem \ref{thm-vanishing}.

\begin{proof}[Proof of Theorem~ \ref{thm-vanishing}]Everything is immediate from
  Theorem~\ref{thm-Klingen-vanish} (taking~$I_a=\emptyset$
  and~$J=I^c$), except
  for the claim that
  $e(T^{I})\mathrm{R}\Gamma(X_1^{G_1, I}, \omega^{\kappa}(-D))$ has
  amplitude $[0, \# I ]$,  which follows from Proposition \ref{prop: the complex K}.
  \end{proof}

 \subsubsection{Commutativity over the ordinary
   locus}\label{subsubsec: commutativity over ordinary locus}While we do not
 prove the commutativity of the correspondences~$T_{w,1}$ and~$T_w$,
 we do prove it over the ordinary locus at~$w$, where all of the
 correspondences are finite flat over the interior. We will need this result at the
 places~$w\in I^c$, because we need to make use of both of these Hecke
 operators in this case (because the Hecke operator~$U_{w,2}$ at
 Klingen level which corresponds to~$T_w$ is needed for those parts of the control theorem which
 take place at the level of the sheaf, but we need to use~$T_{w,1}$ to
 prove the finiteness of cohomology).   
 \begin{lem}\label{lem: Hecke operators commute on the ordinary locus}Suppose that $w\in I^c$, and that~$l_w\ge 2$.  Then on $\mathrm{R}\Gamma(X_1^{G_1,I},\omega^\kappa(-D))$ we have $T_{w,1}\circ T_w=T_w\circ T_{w,1}$.
  \end{lem}\begin{proof}We can easily compose cohomological correspondences when the projections are finite flat.  In particular we may form the compositions $T_w\circ T_{w,1}$ and $T_{w,1}\circ T_w$ over the interior, and it is easy to see that the compositions give the same cohomological correspondence.
    
    In order to check that they commute on
    $\mathrm{R}\Gamma(X_1^{G_1,I},\omega^\kappa(-D))$ we use a similar
    trick to the one that we used to prove Proposition~\ref{prop-cousin-heckeaction}: recall the complex $K^\bullet$ of Proposition
    \ref{prop: the complex K} which computes
    $\mathrm{R}\Gamma(X_1^{G_1,I},\omega^\kappa(-D))$.  We may form
    another complex ${K'}^\bullet$ by applying the same construction
    to the interior $Y^{G_1,I}_1\subset X_1^{G_1,I}$.  As we have
    explained above, the Hecke operators $T_w$ and $T_{w,1}$ commute on
    each term \[\HH^0\big(Y^{ G_1,I, \geq_{ J^c} 2}_1, \omega^{\kappa} (-D)\otimes \bigotimes_{w \in J} (\det \cG_w)^{n(p-1)}/(\sum_{w\in J} (\Ha(\cG_w)^n))\big) \]
  in the definition of~  ${K'}^\bullet$, and hence on the
  subcomplex~$K^\bullet$ of~${K'}^\bullet$ and thus on $\mathrm{R}\Gamma(X_1^{G_1,I},\omega^\kappa(-D))$.
  \end{proof}

  We end this section by proving the following technical result, whose
  formulation relies on Lemma~\ref{lem:
    Hecke operators commute on the ordinary locus}. We will make use
  of it in~\S\ref{subsec: perfect Hida
    complex}, in order to compare the complex of
  Proposition~\ref{prop: the complex K} to the analogous complex at
  Klingen level.  Fix a subset
  $J \subset I$; we now consider the space $X^{I, =_J 1, =_{J^c}
    2}_{K,1}=X_{K,1}^{=_J1,=_{J^c}2}$.
  \begin{lem}\label{lem: control of Ttilde I} There is a universal  constant $C$
  depending only on $p$ and $F$ but not on the tame level $K^p$ such
  that if~$l_v\ge 2$, $k_v-l_v \geq C$ for all $v \in I$,
  $l_v \geq C$ for all $v \in I^c$ and $l_v \geq p+1$ for all $v \in J$, then  ~$\mathrm{R}\Gamma(X^{I, =_J 1, =_{J^c} 2}_{K,1},
  \omega^{\kappa}(-D)) $carries a locally finite action of~$\widetilde{T}^I = \prod_{w |p} T_{w} \prod_{w \in I^c} T_{w,1}$.  
\end{lem}
\begin{proof}
  Note that by Lemma~\ref{lem: Hecke operators commute on the ordinary
    locus}, all of the Hecke operators in the definition
  of~$\widetilde{T}^I$ commute. We begin by showing that the action
  of~$\widetilde{T}^I$
  on $\mathrm{R}\Gamma(X^{I}_{K,1}, \omega^{\kappa}(-D)) $ is locally finite. To this end,
  note that by Theorem~\ref{thm-Klingen-vanish}, the action of~$T^I$
  on
  $\mathrm{R}\Gamma(X^{I}_{K,1},
  \omega^{\kappa}(-D)) $ is locally finite, and
  $e(T^I)\mathrm{R}\Gamma(X^{I}_{K,1}, \omega^{\kappa}(-D)) $ is a perfect
  complex if  $l_v\ge 2$, $k_v-l_v \geq C$ for all $v \in I$, and
  $l_v \geq C$ for all $v \in I^c$. Since~$\widetilde{T}^I=T^I\prod_{w\in I^c}T_w$, it follows
  that the action of ~$\widetilde{T}^I$ is also locally finite (as it
  acts locally  nilpotently on  $(1-e(T^I))\mathrm{R}\Gamma(X^{I}_{K,1}, \omega^{\kappa}(-D)) $).

  Taking the exact triangles induced by 
  $$\omega^\kappa \stackrel{\Ha(\cG_w)}\rightarrow \omega^\kappa \otimes \det\omega_{\cG_w}^{p-1}
   \rightarrow \omega^\kappa \otimes \det\omega_{\cG_w}^{p-1}/ \Ha(\cG_w)$$
    for all $w \in J$, and using Lemma \ref{lem-commutTW ordinary},  we deduce that $\widetilde{T}^I$ is locally finite on 
  $ \mathrm{R}\Gamma(X^{I, =_J1}_{K,1}, \omega^{\kappa}(-D))$ for all
  weights $\kappa = (k_v, l_v)$ with  $l_v\ge 2$, $k_v-l_v \geq C$ for
  all $v \in I$, $l_v \geq C$ for all $v \in I^c$ and $l_v \geq p+1$
  if $v \in J$. Passing to the limit over multiplication by
  $\Ha(\cG_w)$ for  $w \in I\setminus J$, we deduce that $\widetilde{T}^I$ is
  locally finite on  $ \mathrm{R}\Gamma(X^{I, =_J1, =_{J^c} 2}_{K,1},
  \omega^{\kappa}(-D))$, as required.\end{proof}

\subsection{Formal geometry}\label{subsec: formal geometry} In this section we continue to assume that $K =
K^pK_p$, $K_p = \prod_{v |p} K_v$ with $K_v \in  \{
\mathrm{GSp}_4( \ocal_{F_v}), \Par( v)\}$. Our goal
in this section is to define the Igusa tower at Klingen level, and the
$p$-adic sheaves whose cohomology defines our spaces of $p$-adic
automorphic forms.

\subsubsection{Completion of $X$}We adopt the convention that if~$Z$
is a scheme over~$\Spec\Z_{(p)}$, then we write~$Z_n$ for
$Z\otimes_{\Z_{(p)}}\Z/p^n\Z$, and~$\mathfrak{Z}:=\varinjlim_nZ_n$ for
the formal $p$-adic completion of~$Z$, which is by definition a $p$-adic formal scheme.
In particular, we let $\mathfrak{X}_{K}$ be the formal $p$-adic completion
of $X_{K}$, and we write $ \mathfrak{X}^{\geq 2}_{K} \hookrightarrow
\mathfrak{X}^{\geq 1}_{K}  \hookrightarrow \mathfrak{X}_{K}$ for
the open formal subschemes corresponding to $X^{\geq 1}_{K,1}$ and
$X^{ \geq 2}_{K,1}$. We write~$\mathfrak{Y}_{K}$ for the complement of the boundary
of~$\mathfrak{X}_{K}$, with special fibre~$Y_{K,1}$. We
write~$\mathfrak{Y}^{\ge 2}_{K}$ for the ordinary locus on the
interior, and so on.

\subsubsection{Deep Klingen level structure} 
For all $m \geq 1$ we consider the formal scheme $\mathfrak{X}^{ \geq
  1}_{K,\Kli}(p^m) \rightarrow \mathfrak{X}_{K}^{ \geq 1}$ which
parametrizes a subgroup $H_m \subset \cG[p^m]$ which is locally for
the \'etale topology isomorphic to $\mu_{p^m} \otimes \ocal_F$;
equivalently, $H_m=\prod_{w|p}H_{m,w}$  where for each~$w|p$,
$H_{m,w}\subset\cG_w[p^m]$ is isomorphic to~$\mu_{p^m}$.

\begin{prop}\label{prop: deep Klingen is affine etale} The morphism $\mathfrak{X}^{\geq 1}_{K,\Kli}(p^m) \rightarrow \mathfrak{X}_{K}^{ \geq 1}$ is affine and \'etale. Its fibre $\mathfrak{X}^{\geq 2}_{K,\Kli}(p^m)$ over $\mathfrak{X}_{K}^{\geq 2}$ is finite \'etale. 
\end{prop}
\begin{proof}This can be proved in 
  exactly the same way as~\cite[Lem.\ 9.1.1.1]{pilloniHidacomplexes}.
\end{proof}

We denote by $\mathfrak{X}^{\geq 1}_{K,\Kli}(p^\infty) = \varprojlim_m
\mathfrak{X}^{\geq 1}_{K,\Kli}(p^m)$ the $p$-adic formal scheme obtained
by taking the inverse limit (in the category of $p$-adic formal schemes). It exists because the transition morphisms are
affine (see for example~\cite[Prop.\ D.4.1]{MR2441312}, or~\cite[\href{http://stacks.math.columbia.edu/tag/01YT}{Tag
    01YT}]{stacks-project} for the corresponding statement for
  schemes, from which this follows easily).  Over  $\mathfrak{X}^{\geq 1}_{K,\Kli}(p^\infty)$ we have for all places $v |p$ a Barsotti--Tate group of height one and dimension one  $H_{\infty,v} \hookrightarrow \cG_v$. 

\subsubsection{Igusa towers}\label{subsubsec: Igusa towers}
 We  fix a partition $\{ v |p \} = I \coprod I^c$, and we let $X_{K,1}^{ \geq_I 1,\Parmet, \geq_{I^c}
   2}$ be the open subscheme of $X_{K,1}^{ \geq_I 1, \geq_{I^c}
   2}$ where for each place $v\in I$ with $K_v =\Par(\cO_{F_v})$,  the kernel
 of the quasi-polarization $\lambda:\cG_v\to\cG_v^D$ contains a
 multiplicative group (so away from the boundary, this kernel is an
 extension of an \'etale group of rank~$p$ by a multiplicative group
 of rank~$p$). We then let
 $\mathfrak{X}_{K}^{ \geq_I 1,\Parmet, \geq_{I^c} 2}$  be the corresponding open of
 $\mathfrak{X}_{K}$, and in order to save some notation, we will for the moment set 
$\mathfrak{X}_{K}^{ I}:=\mathfrak{X}_{K}^{ \geq_I 1,\Parmet,
  \geq_{I^c} 2}$. The fibre of $\mathfrak{X}^{ \geq 1}_{K,\Kli}(p^m)$  over
$\mathfrak{X}_{K}^{ I}$ is denoted by $\mathfrak{X}^{I}_{K,\Kli}(p^m)$.  Over
$\mathfrak{X}^{I}_{K,\Kli}(p^\infty)$ we have for all places $v \in I$ a
Barsotti--Tate group of height one and dimension one  $H_{\infty,v}
\hookrightarrow \cG_v$. Observe that  for all $v \in I^c$, we have a rank $2$ multiplicative Barsotti--Tate group $\cG_v^m \hookrightarrow \cG_v$ and  that 
$H_{\infty,v} \hookrightarrow \cG_v^{m}$ is a rank one
sub-Barsotti--Tate group.  

If~$K_v=\GSp_4(\cO_{F_v})$ for all~$v|p$, this gives us a convenient alternative description
of~$\mathfrak{X}^{I}_{K,\Kli}(p)$. Set
$$K_p(I) = \prod_{v \in I} \Kli(v) \prod_{v \in I^c} \Iw(v).$$ Then we
have~$\mathfrak{X}^{I}_{K,\Kli}(p)=\mathfrak{X}^I_{K_p(I) K^p}$, where
the superscript~$I$ refers to the fact that for each $v\in I$, the
Klingen level structure~$H_v$ is multiplicative, and at each~$v\in
I^c$, we have extended the given multiplicative Klingen level
structure to the canonical (ordinary) Iwahori level structure.

We denote by $\mathfrak{IG}^{I} \rightarrow  \mathfrak{X}^{I}_{K,\Kli}(p^\infty)$ the profinite-\'etale torsor of trivializations:
$$ \psi_v: \ZZ_p \simeq T_p(H_{\infty,v}^D),~v |p; ~~~~~~\phi_v: \ZZ_p \simeq  T_p((\cG_v^{m}/H_{\infty,v})^D), ~v \in I^c.$$
The upper script $D$ stands for the  dual of these Barsotti--Tate group schemes and $T_p$ stands for the Tate module which here is a pro-\'etale sheaf. For all $v \in I$, there is an action of $\lambda \in \ZZ_p^\times$ on $\psi_v$, mapping $\psi_v$ to $\psi_v \circ \lambda$. For all $v \in I^c$, there is an action of $(\lambda, \mu) \in( \ZZ_p^\times)^2$ on $(\psi_v, \phi_v)$, mapping $(\psi_v, \phi_v)$ to $(\psi_v \circ \lambda, \phi_v \circ \mu)$.

  The Galois  group of the torsor $\mathfrak{IG}^{I} \rightarrow  \mathfrak{X}^{I}_{K,\Kli}(p^\infty)$ is 
  $$T_I:=\prod_{v \in I} \OL^\times_{F_v} \prod_{v \in I^c}( \OL^\times_{F_v})^2 \simeq \prod_{v \in I} \ZZ_p^\times \prod_{v \in I^c}( \ZZ_p^{\times})^2.$$

\subsubsection{Sheaves of $p$-adic modular forms}

Let~$\tLambda_{1,v} = \ZZ_p\llbracket \OL^\times_{F_v} \rrbracket \simeq \ZZ_p\llbracket \ZZ_p^\times\rrbracket$ and 
$\tLambda_{2,v} = \ZZ_p\llbracket(\OL^\times_{F_v})^2\rrbracket  \simeq \ZZ_p\llbracket(\ZZ_p^\times)^2\rrbracket$ be the one and two variable Iwasawa algebras. Let \[\tLambda_{I} = \hat{\otimes}_{v \in I} \tLambda_{1,v} \hat{\otimes}_{v \in I^c} \tLambda_{2,v} = \ZZ_p\llbracket T_I\rrbracket\] and let $\tilde{\kappa}_I:  T_I \rightarrow \tLambda_I^\times$ be the universal character. 

We define a sheaf $\Omega_0^{\tilde{\kappa}_I}$ over $\mathfrak{X}^I_{K,\Kli}(p^\infty)$ by the formula: 
$$ \Omega^{\tilde{\kappa}_I}_0 = ((\pi_{*} \ocal_{ \mathfrak{IG}^I }) \hat{\otimes}_{\Zp} \tLambda_I)^{T_I}$$
where $\pi: \mathfrak{IG}^I  \rightarrow
\mathfrak{X}^I_{K,\Kli}(p^\infty)$ is the affine projection, and the
group $T_I$ acts diagonally (via its natural action on~$\pi_{*}
\ocal_{ \mathfrak{IG}^I }$, and via~$\tilde{\kappa}_I$ on~$\tLambda_I$). 

We set $\Omega^{\tilde{\kappa}_I} = \Omega_0^{\tilde{\kappa}_I} \bigotimes_{v |p }
\det^2 \omega_{\cG_v}$. The explanation for the twist by this
invertible sheaf is given below in~\S\ref{section-expl} (see in
particular Lemma~\ref{lem: HT description of Omega}). This is an invertible sheaf of $\ocal_{\mathfrak{X}^I_{K,\Kli}(p^\infty)} \hat{\otimes}_{\ZZ_p} \tLambda_I$-modules.

\subsubsection{Comparison with classical sheaves}\label{section-expl}

Over $\mathfrak{X}^I_{K,\Kli}(p^\infty)$ we have for all $v |p$ a surjective map  $  \omega_{\cG_v} \rightarrow  \omega_{H_{\infty,v}}$ arising from the differential of the inclusion $H_{\infty,v} \hookrightarrow \cG_v$.

 Let $\kappa = ( (k_v, 2)_{v \in I}, (k_v,l_v)_{v \in I^c}) \in
(\ZZ^2)^{S_p}$, $k_v \geq 2$ if $v \in I$, $k_v \geq l_v$ if $v \in I^c$,  be an algebraic  weight.  By construction
 there is a surjective map\numequation\label{eqn: highest weight of sheaf} \omega^{\kappa}\vert_{ \mathfrak{X}^I_{K_p,\Kli}(p^\infty)}\to (\otimes_{v \in I}  \omega_{H_{\infty,v}}^{k_v-2} ) \bigotimes ( \otimes_{v \in I^c}  \omega_{\cG_v/H_{\infty,v}}^{l_v-2} \otimes \omega_{H_{\infty,v}}^{k_v-2}) \bigotimes (\otimes_{v |p} \det \omega_{\cG_v}^{2})\end{equation}

This map should be interpreted as the projection to the highest weight
vector. Moreover over $\mathfrak{IG}^I $ the Hodge--Tate map (see
\cite[p.\ 117]{MR0347836}, as well as Section
\ref{recollection-Hodgetate} below) provides
maps: $$ \ZZ_p \stackrel{\psi_v}\rightarrow T_p( H_{ \infty,v}^D) \stackrel{\HT}\rightarrow  \omega_{H_{\infty,v}}$$
for all $v\in S_p$, and $$\ZZ_p \stackrel{\phi_v} \rightarrow   T_p( (\cG_v^m/
H_{\infty,v})^D)  \stackrel{\HT} \rightarrow \omega_{\cG_v/H_{\infty,v}}$$
for all $v\in I^c$.

These maps induce isomorphisms after tensoring with $\oscr_{\mathfrak{IG}^I}$ on the left. Therefore the Hodge--Tate map provides a $\ZZ_p^\times$-reduction of the $\mathrm{GL}_1$-torsors $\omega_{H_{\infty,v}}$ and $\omega_{\cG_v/H_{\infty,v}}$.

 Let $\kappa = ( (k_v, 2)_{v \in I}, (k_v,l_v)_{v \in I^c})$ be a
 classical algebraic weight. Then we can  naturally identify $\kappa$  with a
 $p$-adic weight (that is, an element of $\mathrm{Hom} ( \tLambda_I, \ZZ_p)$) via the character:  $$ ( (x_v)_{v \in I}, (x_v, y_v)_{v \in I^c}) \in T_I \mapsto \prod_{v \in I}  x_v^{k_v-2} \prod_{v \in I^c} x_v^{k_v-2} y_v^{l_v-2}$$
   Let us define $\Omega^{\kappa} = \Omega^{\tilde{\kappa}_I} \otimes_{\tLambda_I, \kappa} \ZZ_p$.
 
 \begin{lem}\label{lem: HT description of Omega} For all $\kappa = ( (k_v, 2)_{v \in I}, (k_v,l_v)_{v \in I^c}) \in
(\ZZ^2)^{S_p}$ with $k_v \geq 2$ if $v \in I$, $k_v \geq l_v$ if $v \in I^c$,  there is a canonical isomorphism $$\mathrm{HT}^{*}:   (\otimes_{v \in I}  \omega_{H_{\infty,v}}^{k_v-2} ) \bigotimes ( \otimes_{v \in I^c}  \omega_{\cG_v/H_{\infty,v}}^{l_v-2} \otimes \omega_{H_{\infty,v}}^{k_v-2}) \bigotimes (\otimes_{v |p} \det \omega_{\cG_v}^{2}) \simeq \Omega^{\kappa}.$$
 \end{lem}
 
 \begin{proof}  By definition, it suffices to construct an isomorphism $$\mathrm{HT}^{*}:   (\otimes_{v \in I}  \omega_{H_{\infty,v}}^{k_v-2} ) \bigotimes ( \otimes_{v \in I^c}  \omega_{\cG_v/H_{\infty,v}}^{l_v-2} \otimes \omega_{H_{\infty,v}}^{k_v-2}) \simeq \Omega_0^{\kappa}$$
 where $\Omega^{\kappa}_0 = \Omega_0^{\tilde{\kappa}_I} \otimes_{\tLambda_I, \kappa} \ZZ_p$. Sections of the sheaf $\Omega_0^\kappa$ are rules $f$ associating to
 $(x, (\phi_v)_{v\in I^c}, (\psi_v)_{v\in S_p}) \in
 \mathfrak{IG}^I(R)$,    an element $$f(x, (\phi_v), (\psi_v)) \in R$$
 such that
 \begin{align*}
f( x, (\phi_v\circ \lambda_v^{-1}), (\psi_v\circ
 \beta_v^{-1}) ) &= \kappa( (\lambda_v, \beta_v)) f( x, (\phi_v),
 (\psi_v) )\\ &= \prod_{v \in S_p}  \lambda_v^{k_v-2} \prod_{v \in I^c} \beta_v^{l_v-2} f( x, (\phi_v),
 (\psi_v) )
 \end{align*}
for $((\lambda_v)_{v\in S_p}, (\beta_v)_{v\in I^c}) \in T_I$.

Sections of the sheaf $$ (\otimes_{v \in I}  \omega_{H_{\infty,v}}^{k_v-2} ) \bigotimes ( \otimes_{v \in I^c}  \omega_{\cG_v/H_{\infty,v}}^{l_v-2} \otimes \omega_{H_{\infty,v}}^{k_v-2})$$    are rules $g$ associating to triples $$(x, (a_v)_{v |p}, (b_v)_{v \in I^c})$$ for  $R$ a $p$-adically complete $\ZZ_p$-algebra, $x \in \mathfrak{X}^I_{K,\Kli}(p^\infty)(R)$, $a_v:  R \simeq x^{*} \omega_{H_{v, \infty}}$, $b_v:  R \simeq x^{*} ~\omega_{\cG_{v}/H_{v , \infty}}$   an element $$f(x, a_v, b_v) \in R$$  such that $$f( x, a_v\circ \lambda_v^{-1},  b_v \circ \delta_{v}^{-1}  ) = \prod_{v|p} \lambda_v^{k_v-2}   \prod_{v\in I^c} \delta_v^{l_w-2} f( x, a_v,  b_v )$$ for all $(\lambda_v) \in (R^\times)^{S_p},   (\delta_v) \in (R^\times)^{I^c}$.

To a rule $g$ as above, we associate a rule $$\mathrm{HT}^{*}(g) (x,
(\phi_v)_{v \in I^c}, (\psi_v)_{v \in S_p})  = g (x,
(\mathrm{HT}(\phi_v(1)))_{v \in I^c}, (\mathrm{HT}(\psi_v(1)))_{v \in S_p} ).$$
It is easy to check that the map $\mathrm{HT}^{*}$ is  an
isomorphism. \end{proof}

We can now summarize the interpolation property of the sheaf $\Omega^{\tilde{\kappa}_I}$.

\begin{cor}\label{cor: surjection for Komega}  For all $\kappa = ( (k_v, 2)_{v \in I}, (k_v,l_v)_{v \in I^c}) \in
(\ZZ^2)^{S_p}$ with $k_v \geq 2$ if $v \in I$, $k_v \geq l_v$ if $v \in I^c$, there is a canonical surjective map:
$$ \omega^\kappa \vert_{ \mathfrak{X}^I_{K,\Kli}(p^\infty)} \rightarrow \Omega^\kappa=\Omega^{\tilde{\kappa}_I} \otimes_{\tilde{\Lambda}_I, \kappa} \ZZ_p.$$
\end{cor}

\begin{proof} In view of Lemma~\ref{lem: HT description of Omega},
  this is just the map~(\ref{eqn:
  highest weight of sheaf}).\end{proof}

\subsection{Sheaves of \texorpdfstring{$p$}{p}-adic modular forms for
  \texorpdfstring{$G_1$}{G1}}\label{subsec: p adic forms over G1}

  In this section we explain how we can descend our construction to
  the Shimura variety for $G_1$. This section is the analogue   for $p$-adic sheaves of~\S\ref{subsubsec: weights for
G and G1}.

\subsubsection{Weight space for
  \texorpdfstring{$G_1$}{G1}}\label{subsec: weight space for G1} 

We now assume that $p \neq 2$. We let $T_I^0$ be the pro-$p$ sub-group of $T_I$, so that $T_I = T_I^f \times T_I^0$ is the product of a finite group~$T_I^f$  and $T_I^0$.  We let $\Lambda_I = \ZZ_p \llbracket T_I^0 \rrbracket$.  There is a canonical projection $T_I \rightarrow T_I^0$ and a canonical character $\kappa_1: T_I \rightarrow \Lambda_I^\times$ which identifies $\Lambda_I$ with the deformation space of the trivial character of $T_I$. This canonical projection makes $\Lambda_I$ a quotient of $\tilde{\Lambda}_I$. We let $\kappa_I = \kappa_1 \otimes ((2,2)_{v |p})$. The pair $( \kappa_I, \Lambda_I)$ is the universal deformation space of the character $((2,2)_{v |p})$ mod $p$. We let $\Omega^{\kappa_I} = \Omega^{\tilde{\kappa}_I}\otimes_{\tilde{\Lambda}_I} \Lambda_I$. 

If $\kappa = ((k_v, 2)_{v \in I}, (k_v, l_v)_{v \in I^c})$ is a
classical algebraic weight such that $k_v\equiv l_v \equiv 2~\pmod{p-1}$, then $\kappa$ defines a $\ZZ_p$-point of $\Spf~\Lambda_I$ in the following way: we associate to $\kappa$ the character 
$$ ( (x_v)_{v \in I}, (x_v, y_v)_{v \in I^c}) \in T_I \mapsto \prod_{v \in I}  x_v^{k_v-2} \prod_{v \in I^c} x_v^{k_v-2} y_v^{l_v-2}$$
which factors through a character of $T_I^0$ and therefore defines a morphism $f_\kappa: \Lambda_I \rightarrow \ZZ_p$.  The specialization of $\kappa_I$ along the map $f_\kappa$ recovers the character $\kappa$.

\subsubsection{Descent}\label{subsub-descent}

The group $  (\ocal_{F})_{(p)}^{\times, +} $  can be embedded
``diagonally'' in $T_I$ by sending $x\in  (\ocal_{F})_{(p)}^{\times,
  +} $ to  $((x_v)_v \in I, (x_v,x_v)_{v \in I^c})$ where for all
places $v |p$ we denote by $x_v \in \ocal_{F_v}^\times = \ZZ_p^\times$ the image of $x$
in~$F_v$.   For an element $x \in T_I$, we denote by $x_0$ the
projection of $x$ to $T_I^0$. For an element $x \in T_I^0$, we denote by $\tilde{x}$ the corresponding group element in $\Lambda_I$. 

Since $T_I^0$ is a pro-$p$ group, and $p>2$,  the map $ x \mapsto x^2$
is bijective on $T_I^0$. Accordingly if $x \in T_I^0$, then we define
$\sqrt{x}  \in T_I^0$ by the equation $(\sqrt{x})^2 = x$. We then
define a character $ d: (\ocal_{F})_{(p)}^{\times, +} \rightarrow
(\Lambda_I)^\times$ (where ``$d$'' stands for ``descent'') by the formula: $$ d(x) =  \sqrt{x_0}.$$

The group $(\ocal_F)^{\times,+}_{(p)}$ acts on
$\mathfrak{X}_{K}^I$; in the notation of~\S\ref{subsec:
  integral models}, the
element $x \in (\ocal_F)^{\times,+}_{(p)}$ sends $(A, \iota, \lambda,
\eta = (\eta_1,\eta_2), \eta_p)$ to $(A, \iota, x\lambda,
(\eta_1,x\eta_2), \eta_p)$. We can lift this action to $\Omega^{\kappa_I}$ by setting
$$ x: x^{*} \Omega^{\kappa}_I \rightarrow \Omega^{\kappa}_I$$ to be
the composition of the tautological isomorphism (the construction of
$\Omega^{\kappa_I}$ doesn't depend on the polarization) and
multiplication by $d(x)$. The reader can easily check that this defines an action and is compatible with the construction of~\S \ref{subsubsec: weights for
G and G1}; as always, we are making the choice $w=2$.

For all $n \in \Z_{\geq 0} \cup \{ \infty\}$, we can form the quotient of $\mathfrak{X}^I_{K,\Kli}(p^n)$ by the action of $(\ocal_F)^{\times,+}_{(p)}$ (which factors through a finite group acting freely) and we denote by $\mathfrak{X}^{G_1,I}_{K,\Kli}(p^n)$ the corresponding quotient. The maps $\mathfrak{X}^{I}_{K, \Kli}(p^n) \rightarrow \mathfrak{X}^{G_1,I}_{K, \Kli}(p^n)$ are \'etale. We can also descend the sheaf $\Omega^{\kappa_I}$ to a sheaf $\Omega^{\kappa_I}$ over  $\mathfrak{X}^{G_1,I}_{K,\Kli}(p^\infty)$ using the descent datum provided by~$d$. 

We let $M_I^{\pad, \kappa _I} = \mathrm{R}\Gamma( \mathfrak{X}^{G_1,I}
_{K,\Kli}(p^\infty), \Omega^{\kappa_I}(-D))$ be the cohomology of the $p$-adic cuspidal modular forms of weight $\kappa_I$.

\begin{prop} The canonical map $ M_I^{\pad, \kappa _I}  \rightarrow \mathrm{R}\Gamma( \mathfrak{X}^{I} _{K,\Kli}(p^\infty), \Omega^{\kappa_I}(-D))$ is split in the derived category of $\ZZ_p$-modules.
\end{prop}
\begin{proof} See the proof of Proposition~ \ref{prop-splitGG_1}.
\end{proof}

\subsection{Hecke operators at \texorpdfstring{$p$}{p} on the cohomology of \texorpdfstring{$p$}{p}-adic
  modular forms}\label{subsec: Hecke actions on p adic modular forms}

In this section we define Hecke operators at $p$ acting on the cohomology of $p$-adic modular forms. Recall that  we have fixed a partition $S_p=\{v |p\} = I \coprod I^c$.

\subsubsection{Hecke operators of Siegel type}\label{subsubsec:
  Siegel type padic Hecke operators}

Let $w \in I^c$ be a place above $p$. Let $K = K^pK_p$ be a reasonable
compact open subgroup with $K_p = G_1( \ZZ_p)$. In~\S
\ref{sect-def-Siegel-Hecke} we  defined a Hecke operator attached
to the correspondence (for suitable choices of polyhedral cone decompositions omitted from the notation):
 \begin{eqnarray*}     
 \xymatrix{ &X_{K' }   \ar[rd]^{p_1}\ar[ld]_{p_2}&  \\
  X_K & & X_K}
  \end{eqnarray*}
where $K' = K^p K'_p$ and $K'_p = \prod_{v |p, v \neq w}
\mathrm{GSp}_4( \ocal_{F_v}) \times \Si( w)$. The
map $p_2$ depends on the choice of an element $x_w \in F^{\times, +}$.
We are now going to pull back this correspondence to a deep Klingen
level structure and isolate the ``essential part''. As in
Remark~\ref{rem: independence of compactification}, the resulting
Hecke operators are easily seen to be independent of the choice of
polyhedral cone decomposition.
  
Taking formal $p$-adic completions, we obtain a correspondence:
     \numequation  \label{eqn: completion of Siegel type Hecke operator}
 \xymatrix{ &\mathfrak{X}_{K' }   \ar[rd]^{p_1}\ar[ld]_{p_2}&  \\
  \mathfrak{X}_K & & \mathfrak{X}_K}
  \end{equation}

  We consider the fibre product
  $\mathfrak{X}_{K'} \times_{p_1, \mathfrak{X}_K}
  \mathfrak{X}^{I}_{K,\Kli}(p^m)$.  Recall that by our assumption that
  $w\in I^c$, $\cG_w$ is an ordinary Barsotti--Tate group. We denote by
  $\mathfrak{C}_{w,1}(p^m)$ the open and closed formal subscheme of
  this fibre product where the kernel of the canonical isogeny
  $p_1^{*} \cG \rightarrow p_2^{*} \cG$ has trivial multiplicative
  part (it is open and closed by the rigidity of multiplicative groups).

There is an obvious map $u_1: \mathfrak{C}_{w,1}(p^m) \rightarrow
\mathfrak{X}^I_{K,\Kli}(p^m)$, given by projection onto the second
factor of the fibre product. We claim that the projection 
$\mathfrak{C}_{w,1}(p^m) \rightarrow \mathfrak{X}_K $ induced
by~$p_2$ can be lifted to
a map $u_2: \mathfrak{C}_{w,1}(p^m) \rightarrow
\mathfrak{X}^I_{K,\Kli}(p^m)$. Indeed, since~$H_m$ is multiplicative,
the isogeny $p_1^{*} \cG \rightarrow p_2^{*} \cG$ induces an
isomorphism from $p_1^{*} H_m$ to
its image in $p_2^{*} \cG$. We call this image $p_2^{*} H_m$.  We therefore have a correspondence
\begin{eqnarray*}     
 \xymatrix{ &\mathfrak{C}_{w,1}(p^m)   \ar[rd]^{u_1}\ar[ld]_{u_2}&  \\
 \mathfrak{X}^I_{K,\Kli}(p^m) & & \mathfrak{X}^I_{K,\Kli}(p^m)}
  \end{eqnarray*}
We now associate to this correspondence a Hecke operator $U_{w,1}$. 

\begin{rem} The Hecke operator $U_{w,1}$ is  the standard ``$U_p$'' operator (at the place~$w$) that is considered in the usual theory of $p$-adic modular forms. 
\end{rem}

\begin{lem}\label{lem-existencetracemapu_1}There is a normalized trace
  map ``$ \frac{1}{p^3} \mathrm{Tr}_{u_1}$'' $:  \mathrm{R} (u_1)_*  \cO_{\mathfrak{C}_{w,1}(p^m)} \to \cO_{\mathfrak{X}^I_{K,\Kli}(p^m)}$.
\end{lem}

\begin{proof} The formal schemes $\mathfrak{C}_{w,1}(p^m)$ and $\mathfrak{X}^I_{K,\Kli}(p^m)$ are smooth over $\Z_p$. Consider the  map induced by $u_1$ on top-differentials:
$$ \mathrm{d} u_1:  \det \Omega^1_{\mathfrak{X}^I_{K, \Kli}(p^m)/\ZZ_p} \rightarrow \det \Omega^1_{\mathfrak{C}_{w,1}(p^m)/\ZZ_p}.$$

This map is divisible by $p^3$  by the same arguments as in the proof
of Lemma~ \ref{lem-normalization1}. Namely, the map $u_1$ is totally
inseparable and hence a homeomorphism.  
For any closed point $x \in \mathfrak{X}^I_{K,\Kli}(p^m)$ in the interior, one sees by Serre--Tate theory that  the map of completed local rings $\widehat{\cO_{\mathfrak{X}^I_{K,\Kli}(p^m),x}} \rightarrow \widehat{\cO_{\mathfrak{C}_{w,1}(p^m),x}}$ is given by $\otimes_{v| p} W(k(x)) [[ T_{1,v}, T_{2,v}, T_{3,v}]] \rightarrow \otimes_{v| p} W(k(x)) [[ T_{1,v}, T_{2,v}, T_{3,v}]]$ where $T_{i,v} \mapsto T_{i,v}$ if $v \neq w$,  and $T_{i,w} \mapsto (1+T_{i,w})^p-1$. 

 By reduction modulo $p^n$ of
$u_1$, we get a proper map  $u_1: {C}_{w,1}(p^m)_n \rightarrow
X^I_{K,\Kli}(p^m)_n$ of smooth schemes over $\Spec \ZZ/p^n\ZZ$. The
above map $\frac{1}{p^3}\mathrm{d}u_1$ induces a map $
\oscr_{{C}_{w,1}(p^m)_n} \rightarrow u_1^!
\oscr_{X^I_{K,\Kli}(p^m)_n}$ or by adjunction a map  ``$\frac{1}{p^3}
\mathrm{Tr}_{u_1}$'' $:  \mathrm{R} (u_1)_*  (\cO_{\mathfrak{C}_{w,1}(p^m)}/p^n) \to (\cO_{\mathfrak{X}^I_{K,\Kli}(p^m)})/p^n$ (see~\S\ref{subsubsec: fundamental class}). Passing to the limit over $n$ yields the map of the lemma. 
\end{proof} 

\begin{rem}\label{rem: alternative proof of normalization away from boundary} We sketch another argument for the proof of the lemma. Write $\mathfrak{C}_{w,1}(p^m)_{\mathfrak{Y}}$ for the
restriction of~$\mathfrak{C}_{w,1}(p^m)$ to~$\mathfrak{Y}^I_{K,\Kli}(p^m)$. It is easy to show that the map $u_1:\mathfrak{C}_{w,1}(p^m)_{\mathfrak{Y}}\to
  \mathfrak{Y}^I_{K,\Kli}(p^m)$ is finite flat of degree
  $p^3$.  Therefore, the restriction of the map of the lemma to
  $\mathfrak{Y}^I_{K,\Kli}(p^m)$ is  the usual trace map, normalized
  by a factor $p^{-3}$. Let $\Sigma$ be the $K$-admissible polyhedral cone decomposition  such that  $\mathfrak{X}_{K} = \mathfrak{X}_{K, \Sigma}$ and $ \mathfrak{X}^I_{K,\Kli}(p^m) =  \mathfrak{X}^I_{K,\Kli}(p^m)_{\Sigma}$ . We can use the same $\Sigma$ to get the (non smooth) toroidal compactification of $\mathfrak{X}_{K', \Sigma}$ and then of $\mathfrak{C}_{w,1}(p^m)_{\Sigma}$. Now we observe that the map $\mathfrak{C}_{w,1}(p^m)_{\Sigma} \rightarrow \mathfrak{X}^I_{K,\Kli}(p^m)_{\Sigma}$ is finite flat and therefore has a trace map. It remains to recall that for any refinement $\Sigma'$ of $\Sigma$, the map $\pi: \mathfrak{C}_{w,1}(p^m)_{\Sigma'}  \rightarrow \mathfrak{C}_{w,1}(p^m)_{\Sigma} $ induces a quasi-isomorphism: $\mathrm{R} \pi_\star  \oscr_{\mathfrak{C}_{w,1}(p^m)_{\Sigma'} } = \oscr_{\mathfrak{C}_{w,1}(p^m)_{\Sigma}}$.

    \end{rem}

To define the Hecke operator $U_{w,1}$ on
$\mathrm{R}\Gamma( \mathfrak{X}^{I}_{K,\Kli}(p^\infty),
\Omega^{\kappa_I})$,
$\mathrm{R}\Gamma( \mathfrak{X}^{I}_{K,\Kli}(p^m), \omega^{\kappa})$,
and so on, we argue as follows. By the usual formalism, if~$\cF$ is
one of~$\Omega^{\kappa_I}$ or~$\omega^{\kappa}$, it is enough to define
morphisms $u_2^*\cF\to u_1^*\cF$; we can then compose with the trace
map of Lemma~\ref{lem-existencetracemapu_1}. 

To this end, note that over $\mathfrak{C}_{w,1}(p^m)$ we have the
canonical \'etale isogeny $u_1^*\cG\to u_2^*\cG$, which  determines an isomorphism on
differentials. We thus have a canonically determined isomorphism
$u_2^*\omega^\kappa\to u_1^*\omega^\kappa$ (with no need to
normalize). Similarly, since the canonical isogeny induces an
isomorphism $u_1^*H_m\to u_2^*H_m$, we have a canonical isomorphism
$u_2^*\Omega^{\kappa_I} \to u_1^*\Omega^{\kappa_I}$ (again with no
need to normalize).

\subsubsection{The operator $U_{\Kli(w),1}$} \label{subsec: U Kli 1
  first time around}
We now introduce another  Hecke operator of Siegel type for~$w \in I$,
which we denote~ $U_{\Kli(w),1}$.
This operator will not be used 
until~\S\ref{sec:doubling} and~\S\ref{sec:CG}. We decided to introduce it here because its definition is similar to the other operators of Siegel type introduced in \S\ref{subsubsec:
  Siegel type padic Hecke operators}, and because it is convenient to
discuss the commutativity of all of our Hecke operators at~$p$ in one
go (see Lemma~\ref{lem: commutativity of Hecke ops on infinite Klingen
    big sheaf} below). We defer the details of the normalization of
  this Hecke operator to~\S\ref{sect-UKLI}, where we will consider the
  operator~$U_{\Kli(w),1}$ in a more general context.

  We again consider the correspondence~(\ref{eqn:
  completion of Siegel type Hecke operator}), and the product
$\mathfrak{X}_{K'} \times_{p_1, \mathfrak{X}_K}
\mathfrak{X}^{I}_{K,\Kli}(p^m)$. We denote by
$\mathfrak{C}_{\Kli(w),1}(p^m)$ the open and closed formal subscheme
of this fibre product where the kernel of the canonical isogeny
$p_1^{*} \cG \rightarrow p_2^{*} \cG$ has trivial intersection with
the group~$p_1^*H_m$. Exactly as above, we obtain a
correspondence \begin{eqnarray*}
                 \xymatrix{ &\mathfrak{C}_{\Kli(w),1}(p^m)   \ar[rd]^{v_1}\ar[ld]_{v_2}&  \\
                 \mathfrak{X}^I_{K,\Kli}(p^m) & &
                                                  \mathfrak{X}^I_{K,\Kli}(p^m)}
  \end{eqnarray*}

  We show in Lemma~ \ref{lem-onemorenormalization} below that there is a
  Hecke operator
  \numequation\label{eqn: U Kli 1 weight 2} U_{\Kli(w),1}: \mathrm{R} (v_1)_* v_2^*\omega^2\to\omega^2,\end{equation}
  defined using a trace map normalized by a factor of~$1/p^3$.
On the other hand, we have natural isomorphisms $v_2^\star
\Omega^{\tilde{\kappa}_I}_0 \rightarrow v_1^\star
\Omega^{\tilde{\kappa}_I}_0 $, and tensoring this map with~\eqref{eqn: U Kli 1 weight 2} produces the desired cohomological correspondence (and associated Hecke operator):   $$ U_{\Kli(w),1}:  \mathrm{R} (v_1)_* v_2^*\Omega^{\kappa_I}\to\Omega^{\kappa_I}.$$

\subsubsection{Hecke operators of Klingen type}\label{subsubsec:
Klingen type padic Hecke operators} 
Let $w |p$ be a
place.  Let $K = K^pK_p$ be a reasonable compact open subgroup with
$K_p = G_1( \ZZ_p)$. In~\S\ref{sect-def-Klingen-Hecke} we have
defined a Hecke operator attached to the correspondence (again, for suitable choices of polyhedral cone decompositions omitted from the notation):
 \begin{eqnarray*}     
 \xymatrix{ &X_{K' }   \ar[rd]^{p_1}\ar[ld]_{p_2}&  \\
  X_{K''} & & X_K}
  \end{eqnarray*}
where $K' = K^p K'_p$ with $K'_p = \prod_{v |p, v \neq w}
\mathrm{GSp}_4( \ocal_{F_v}) \times \Kli(w)$, and
$K'' = K^p K''_p$ with $K''_p = \prod_{v |p, v \neq w}
\mathrm{GSp}_4( \ocal_{F_v}) \times \Par(w)$. The
map $p_2$ depends on the choice of an element $x_w \in F^{\times, +}$,
and over $X_{K'}$ we have natural isogenies \[p_1^{*}\cG\to p_2^{*}\cG\to
  p_1^{*}\cG \] whose composite is multiplication by~$x_w$.
As in~\S\ref{subsubsec:
  Siegel type padic Hecke operators}, we are going to pull back this
correspondence to a deep Klingen level structure and isolate the
``essential part'' of the correspondence. 

  Taking formal $p$-adic completions, we obtain:
     \begin{eqnarray*}     
 \xymatrix{ &\mathfrak{X}_{K' }   \ar[rd]^{p_1}\ar[ld]_{p_2}&  \\
  \mathfrak{X}_{K''} & & \mathfrak{X}_K}
  \end{eqnarray*}

We consider the fibre product  $\mathfrak{X}_{K'} \times_{p_1,
  \mathfrak{X}_K} \mathfrak{X}^{I}_{K,\Kli}(p^m)$.  We denote by
$\mathfrak{C}_{w,2,1}(p^m)$ the  formal
subscheme where the kernel of the isogeny $p_1^{*} \cG \rightarrow p_2^{*} \cG$ has trivial intersection with the  group $p_1^{*} H_m$. 

\begin{lem} The formal subscheme $\mathfrak{C}_{w,2,1}(p^m)$ is open and closed in $\mathfrak{X}_{K'} \times_{p_1,
  \mathfrak{X}_K} \mathfrak{X}^{I}_{K,\Kli}(p^m)$.
  \end{lem}
  
 \begin{proof} Let $L = p_1^*H_m \cap  \mathrm{Ker} (p_1^{*} \cG \rightarrow p_2^{*} \cG)$. Since $L$ is a closed subscheme of $p_1^*H_m$,
 it is finite over $\mathfrak{X}_{K'} \times_{p_1, \mathfrak{X}_K}
  \mathfrak{X}^{I}_{K,\Kli}(p^m)$ and the condition that $L = \{0\}$ is therefore open. It is also closed because if at some point $x$ there is a non-trivial map $p_1^\star H_m\vert_x \rightarrow \mathrm{Ker} (p_1^{*} \cG \rightarrow p_2^{*} \cG)\vert_x$, this map will extend on the completed local ring at $x$ by the rigidity of multiplicative groups. 
  \end{proof}

There is an obvious map $r_1: \mathfrak{C}_{w,2,1}(p^m) \rightarrow
\mathfrak{X}^I_{K,\Kli}(p^m)$ induced by the projection~$p_1$. We
claim that the second projection~$p_2$, which induces a map
$\mathfrak{C}_{w,2,1}(p^m) \rightarrow \mathfrak{X}_{K''} $, can be
lifted to a map $r_2: \mathfrak{C}_{w,2,1}(p^m) \rightarrow
\mathfrak{X}^I_{K'',\Kli}(p^m)$. Indeed, over $
\mathfrak{C}_{w,2,1}(p^m)$ the isogeny $p_1^{*} \cG \rightarrow
p_2^{*} \cG$ induces an isomorphism from $p_1^{*} H_m$ to its
image in $p_2^{*} \cG$  (which we call $p_2^{*} H_m$). We therefore have a correspondence
\begin{eqnarray*}     
 \xymatrix{ &\mathfrak{C}_{w,2,1}(p^m)   \ar[rd]^{r_1}\ar[ld]_{r_2}&  \\
 \mathfrak{X}^I_{K'',\Kli}(p^m) & & \mathfrak{X}^I_{K,\Kli}(p^m)}
  \end{eqnarray*}

We now associate to this correspondence a Hecke operator $$U'_{w} \in
\mathrm{Hom}(\mathrm{R}\Gamma( \mathfrak{X}^{I}_{K'',\Kli}(p^\infty),
\Omega^{\kappa_I}), \mathrm{R}\Gamma(
\mathfrak{X}^{I}_{K,\Kli}(p^\infty), \Omega^{\kappa_I})).$$

To do so, we have the following lemma.

\begin{lem}\label{lem-tracemapr1} There is a normalized trace map  ``$ \frac{1}{p^2} \mathrm{Tr}_{r_1}$''$:(r_1)_*\cO_{\mathfrak{C}_{w,2,1}(p^m)} \to \cO_{\mathfrak{X}^I_{K,\Kli}(p^m)}$
\end{lem}

\begin{proof} The formal schemes $\mathfrak{C}_{w,2,1}(p^m)$ and $\mathfrak{X}^I_{K,\Kli}(p^m)$ are smooth over $\Z_p$. Consider the induced map on top differentials:
$$ \mathrm{d} r_1:  \det \Omega^1_{\mathfrak{X}^I_{K, \Kli}(p^m)/\ZZ_p} \rightarrow \det \Omega^1_{\mathfrak{C}_{w,2,1}(p^m)/\ZZ_p}$$

This map is divisible by $p^2$  for the same reason as in the proof of
Lemma~ \ref{lem-normalization2}. Namely, let us fix a closed point $x \in \mathfrak{X}^I_{K,\Kli}(p^m)$ which is in the interior and  ordinary  at $w$. The fibre of $r_1$ at $x$ parametrizes the subgroup $L = \mathrm{Ker} (p_1^{*} \cG \rightarrow p_2^{*} \cG)$ of $\mathcal{G}_w[p]$ of \'etale rank $2$, multiplicative rank $1$ and trivial intersection with $H_m$. The total degree of $r_1$ is $p^3$.  The fibre of $r_1$ over $x$ has $p$ points (corresponding to the choice of the multiplicative part $L^m$ of $L$). The inseparability degree is $p^2$ (corresponding to finding sections of  $\mathcal{G}_w[p]/L^m \rightarrow \mathcal{G}_w[p]^{et})$. For any $x' \in \mathfrak{C}_{w,2,1}(p^m)$ lying above $x$, Serre--Tate theory shows that the map on completed local rings $\widehat{\cO_{\mathfrak{X}^I_{K,\Kli}(p^m),x}} \rightarrow \widehat{\cO_{\mathfrak{C}_{w,2,1}(p^m),x'}}$ is isomorphic  to  $\otimes_{v| p} W(k(x)) [[ T_{1,v}, T_{2,v}, T_{3,v}]] \rightarrow \otimes_{v| p} W(k(x')) [[ T_{1,v}, T_{2,v}, T_{3,v}]]$ where $T_{i,v} \mapsto T_{i,v}$ if $v \neq w$ or $i=1$,  and $T_{i,w} \mapsto (1+T_{i,w})^p-1$ for $i=2,3$.

  By reduction modulo $p^n$ of $r_1$,
we get a proper map  $r_1: {C}_{w,2,1}(p^m)_n \rightarrow
X^I_{K,\Kli}(p^m)_n$ of smooth schemes over $\Spec \ZZ/p^n\ZZ$. The
above map $\frac{1}{p^2}\mathrm{d}r_1$ induces a map $
\oscr_{{C}_{w,2,1}(p^m)_n} \rightarrow r_1^!
\oscr_{X^I_{K,\Kli}(p^m)_n}$ or by adjunction a map  ``$\frac{1}{p^2}
\mathrm{Tr}_{r_1}$'' $:  \mathrm{R} (r_1)_*  (\cO_{\mathfrak{C}_{w,2,1}(p^m)}/p^n) \to (\cO_{\mathfrak{X}^I_{K,\Kli}(p^m)})/p^n$ (see~\S\ref{subsubsec: fundamental class}). Passing to the limit over $n$ yields the map of the lemma. 
\end{proof} 

\begin{rem}  We could give an alternative proof of
  Lemma~\ref{lem-tracemapr1} as in Remark~\ref{rem: alternative proof of normalization away from boundary}. 
  \end{rem}

As usual over $\mathfrak{C}_{w,2,1}(p^m)$ we have the
canonical isogeny $r_1^*\cG\to r_2^*\cG$, whose differential  determines  a morphism
$r_2^*\omega^\kappa\to r_1^*\omega^\kappa$.  We have a
commutative diagram\[\xymatrix{r_2^*\omega_{\cG_w} \ar[r]\ar[d] &
    r_2^*\omega_{ H_{m,w}}\ar[d]\ar[r] & 0\\ r_1^* \omega_{\cG_w}\ar[r]&
   r_1^* \omega_{H_{m,w}}\ar[r]&0 } \] which Zariski locally on affine
opens~$\Spf R $ is isomorphic to \numequation\label{eqn: local description of U prime on differentials}\xymatrix{R^2 \ar[r]^{(0\ 1)}\ar[d]_{
       \begin{pmatrix}
         p&0\\0&1
       \end{pmatrix}
} &
   R\ar[r]\ar[d]^{1_R} & 0\\ R^2\ar^{(0\ 1)}[r]&
   R\ar[r]&0 } \end{equation}
It follows that we can and do normalize the morphism
$r_2^*\omega^\kappa\to r_1^*\omega^\kappa$ by dividing by~ $p^{l_w}$.  
When $m = \infty$, the isogeny induces an isomorphism $r_1^*
H_{\infty,w} \rightarrow r_2^* H_{\infty, w}$, and we  therefore obtain an isomorphism  $r_2^*\Omega^{\kappa_I} \to r_1^*\Omega^{\kappa_I}$.  Combining this with Lemma~ \ref{lem-tracemapr1} gives the desired operator $U'_w$.

We now exchange the roles of $p_1$ and $p_2$, and consider the fibre
product  $\mathfrak{X}_{K'} \times_{p_2, \mathfrak{X}_K''}
\mathfrak{X}^{I}_{K'',\Kli}(p^m)$.  We denote by
$\mathfrak{C}_{w,2,2}(p^m)$ the open and closed formal
subscheme where the kernel of the isogeny $p_2^{*} \cG \rightarrow p_1^{*}
\cG$ has trivial connected component (so that away from the
boundary,  the kernel of this isogeny is \'etale). Note that by
definition this kernel is contained in the kernel of the quasi-polarization $p_2^{*} \cG \rightarrow p_2^{*}\cG^D$,  so the kernel of $p_2^{*} \cG \rightarrow
p_1^{*}\cG$ has multiplicative rank at least~$1$.

The projection~$p_2$ induces a map $s_2: \mathfrak{C}_{w,2,2}(p^m) \rightarrow
\mathfrak{X}^I_{K'',\Kli}(p^m)$. We claim that the first projection
$p_1:\mathfrak{C}_{w,2,2}(p^m) \rightarrow \mathfrak{X}_{K} $ can be lifted
to a map $s_1: \mathfrak{C}_{w,1,2}(p^m) \rightarrow
\mathfrak{X}^I_{K,\Kli}(p^m)$. Indeed, since~$H_m$ is connected, we see
that over $
\mathfrak{C}_{w,2,2}(p^m)$ the isogeny $p_2^{*} \cG \rightarrow
p_1^{*} \cG$ induces an isomorphism from $p_2^{*} H_m$ to its
image in $p_1^{*} \cG$ (which we call $p_1^{*} H_m$); and the map
$\mathfrak{C}_{w,2,2}(p^m) \rightarrow \mathfrak{X}_{K}$ factors
through $\mathfrak{X}^{I}_{K}$. Accordingly, we
have a correspondence 
\begin{eqnarray*}     
 \xymatrix{ &\mathfrak{C}_{w,2,2}(p^m)   \ar[rd]^{s_2}\ar[ld]_{s_1}&  \\
 \mathfrak{X}^I_{K,\Kli}(p^m) & & \mathfrak{X}^I_{K'',\Kli}(p^m)}
  \end{eqnarray*}
We can associate to this correspondence a Hecke operator $$U''_{w} \in
\mathrm{Hom}(\mathrm{R}\Gamma( \mathfrak{X}^{I}_{K,\Kli}(p^\infty),
\Omega^{\kappa_I}), \mathrm{R}\Gamma( \mathfrak{X}^{I}_{K'',\Kli}(p^\infty), \Omega^{\kappa_I})),$$ which
again depends on the construction of a trace map:

\begin{lem}
  \label{cor: U2 double prime normalized trace map exists away from boundary}There is a normalized
  trace map $p^{-1}\Tr_{s_2}:
 \mathrm{R}(s_2)_*\cO_{\mathfrak{C}_{w,2,2}(p^m)}\to \cO_{\mathfrak{X}^I_{K'',\Kli}(p^m)}$.
\end{lem}
\begin{proof} This is a calculation in Serre--Tate theory which is
  similar to the proof of Lemma~ \ref{lem-tracemapr1}.   Namely, let us fix a closed point $x \in \mathfrak{X}^I_{K'',\Kli}(p^m)$ which is in the interior. The fibre of $s_2$ at $x$ parametrizes  rank $p$ \'etale subgroups in the kernel of the quasi-polarization $\mathcal{G}_w \rightarrow \mathcal{G}_w^D$ (which is a rank $p^2$ finite flat group scheme, extension of an \'etale by a multiplicative subgroup). We deduce that the map $s_2$ is totally inseparable at $x$ of degree $p$. 
   Serre--Tate theory shows that the map on completed local rings $\widehat{\cO_{\mathfrak{X}^I_{K'',\Kli}(p^m),x}} \rightarrow \widehat{\cO_{\mathfrak{C}_{w,2,2}(p^m),x}}$ is isomorphic  to  $\otimes_{v| p} W(k(x)) [[ T_{1,v}, T_{2,v}, T_{3,v}]] \rightarrow \otimes_{v| p} W(k(x)) [[ T_{1,v}, T_{2,v}, T_{3,v}]]$ where $T_{i,v} \mapsto T_{i,v}$ if $v \neq w$ or $i=1,2$,  and $T_{3,w} \mapsto (1+T_{3,w})^p-1$. 
 \end{proof}

Over $\mathfrak{C}_{w,2,2}(p^m)$ we have the
canonical isogeny $s_2^*\cG\to s_1^*\cG$, which is \'etale, and therefore determines isomorphisms
$s_1^*\omega^\kappa\to s_2^*\omega^\kappa$ and
$s_1^*\Omega^{\kappa_I} \to s_2^*\Omega^{\kappa_I}$. Combining with Lemma~ \ref{cor: U2 double prime normalized trace map exists away from boundary}, we get the desired Hecke operator $U''_{w}$. 

We set $U_{w,2}:= U'_{w} \circ U''_{w}$. 

\subsubsection{Commutativity of the Hecke operators} 

We remind the reader that whenever we write~$U_{v,1}$ below, we mean~$U_{\Iw(v),1}$.

\begin{lem}\label{lem: commutativity of Hecke ops on infinite Klingen
    big sheaf} The operators $\{U_{\Kli(v),1}, U_{v,2} \}_{ v \in I}$ and $\{U_{v,1}, U_{v,2} \}_{v \in I^c}$ commute with each other on  $\RGamma(
\mathfrak{X}^{G_1,I}_{K, \Kli(p^\infty)}, \Omega^{\kappa_I}(-D))$.
\end{lem}

\begin{proof}We prove this in the same way as Lemma~\ref{lem: Hecke operators commute on the ordinary locus}.  We  first introduce  a similar complex as in \S\ref{section-computingcohoexplicitely} to compute the cohomology. Namely, there is a complex $L^\bullet$ computing $\RGamma( \mathfrak{X}^{G_1,I}_{K, \Kli(p^\infty)}, \Omega^{\kappa_I}(-D))$, such that $L^k = \varprojlim_t L^k_t$ where~$L^k_t$ is
  $$  \textstyle{\bigoplus_{ J \subset I, \# J = k}} 
  \displaystyle{ \colim_{\times \prod_{w \in J}
   \Ha(\cG_w)}} \HH^0\big(\mathfrak{X}^{ G_1,I, \geq_{ J^c} 2}_{K,\Kli(p^\infty)}, \Z/p^t \Z \otimes
 \Omega^{\kappa_I} (-D)\otimes \bigotimes_{w \in J} (\det
 \cG_w)^{n(p-1)}/(\sum_{w\in J} (\Ha(\cG_w)^n) \big) $$
Note that in that formula we use the fact that $\Ha(\cG_w)^n$ has a canonical lift to $\Z/p^t\Z$ for all $n$'s which are multiples of $p^{t-1}$ and these are cofinal among all natural numbers. 

Each term \numequation\label{eqn: thing in complex term}\displaystyle{ \colim_{\times \prod_{w \in J}
   \Ha(\cG_w)}} \HH^0\big(\mathfrak{X}^{ G_1,I, \geq_{ J^c} 2}_{K,\Kli(p^\infty)}, \Z/p^t \Z \otimes
 \Omega^{\kappa_I} (-D)\otimes \bigotimes_{w \in J} (\det
 \cG_w)^{n(p-1)}/(\sum_{w\in J} (\Ha(\cG_w)^n) \big)\end{equation}
appearing in this complex is stable under the Hecke action and it is
therefore enough to prove the commutativity for each of these
terms. Each of the terms~\eqref{eqn: thing in complex term} can be
embedded into the corresponding direct limit of
cohomology groups taken over the interior of the moduli space. Over the
interior of the moduli space all our correspondences are finite flat
and the commutativity follows from standard properties of the Iwahori and Klingen Hecke algebras. 
\end{proof}

\begin{lem}
  \label{lem: quadratic relation for U Kli 1 big sheaf}If $w\in I$
  then we have an equality of Hecke operators \[U_{\Iw(w),1} (U_{\Kli(w),1}-U_{\Iw(w),1}) =
  U_{w,2}\] on $\RGamma(
\mathfrak{X}^{G_1,I,=_w2}_{K, \Kli(p^\infty)}, \Omega^{\kappa_I}(-D)).$ 
\end{lem}
  \begin{proof}Using a Cousin complex computing the cohomology as in the proof of
    Lemma~\ref{lem: commutativity of Hecke ops on infinite Klingen
    big sheaf}, we reduce to proving  that the underlying (cohomological) correspondences agree away from the boundary. By definition, the correspondence
    associated to $U_{\Iw(w),1}$ parameterizes triples $(\cG,H_m,L)$
    where~$L\subset \cG_w[p]$ is \'etale, totally isotropic of
    degree~$p^2$, and has $L\cap H_m=\{0\}$. Similarly, the
    correspondence associated to~$Z_w$
    parameterizes triples $(\cG,H_m,M)$ where~$M\subset \cG_w[p]$ has
    multiplicative rank~$1$, is totally isotropic of degree~$p^2$, and
    has $M\cap H_m=\{0\}$. Finally, as in~\cite[Prop.\ 10.2.1]{pilloniHidacomplexes}, the correspondence associated
    to~$U_{w,2}$ parameterizes triples $(\cG,H_m,N)$
    where~$N\subset \cG_w[p^2]$ is totally isotropic of degree~$p^4$,
    $N[p]$ has degree~$p^3$, and $N\cap H_m=\{0\}$. 

Comparing these definitions, we see that on the level of underlying
correspondences, we have  $pU_{w,2}=U_{\Iw(w),1}Z_w$. Since the
normalization factors involved in the Hecke operators $U_{w,2},
U_{\Iw(w),1}, Z_w$ are respectively $p^{-5},p^{-3},p^{-3}$, and since
$5+1=3+3$, the result follows.
  \end{proof}

\subsection{Perfect complexes of \texorpdfstring{$p$}{p}-adic modular
  forms} \label{subsec: perfect Hida complex}

Let $K = K^pK_p$ be a reasonable
compact open subgroup with $K_p = G_1( \ZZ_p)$. Set \[U^I = \prod_{v \in I} U_{v,2} \prod_{v \in I^c} U_{v,1}U_{v,2}.\]
This is an endomorphism of  $M_I^{\pad, \kappa_I} = \RGamma(
\mathfrak{X}^{G_1,I}_{K, \Kli(p^\infty)}, \Omega^{\kappa_I}(-D))$, an
object of the bounded derived category of $\Lambda_I$-modules. In this
section we prove the following theorem.

 \begin{thm}\label{theorem-p-adic-complex} \leavevmode \begin{enumerate}
\item The operator $U^I$ is locally finite on $M_I^{\pad,
    \kappa_I}$. 
\item Let $e(U^I)$ be the ordinary projector attached  to $U^I$ and
  let  $M_I:= e (U^I) M_I^{\pad, \kappa_I}$ be the associated direct
  summand. Then the complex $M_I$ is a perfect complex of $\Lambda_I$-modules concentrated in the interval $[0, \# I]$.
\item For all  classical algebraic weights $\kappa = (( k_v, l_v)_{v |p})$ with $l_v = 2$ when $v \in I$ and $k_v\equiv l_v\equiv 2\pmod{p-1}$ for all $v |p$, there is a canonical quasi-isomorphism: 
$$ e (U^I) \RGamma(\mathfrak{X}^{ G_1, I}_{K_p(I) K^p},
\omega^\kappa(-D)) \rightarrow  M_I \otimes^{\mathbf{L}}_{\Lambda_I, \kappa}
\ZZ_p.$$ \item There is a universal constant $C$ depending only on $p$ but not
  on the tame level $K^p$ such that for all  classical algebraic
  weights $\kappa = (( k_v, l_v)_{v |p})$ with $l_v = 2$ when $v
  \in I$, $k_v\equiv l_v\equiv 2\pmod{p-1}$ for all $v |p$, $k_v-l_v \geq
  C$ when $v | p$, and $l_v \geq  C$ when $v \in I^c$, the map: 
$$ e(  \prod_{v |p} T_{v} \prod_{v \in I^c} T_{v,1}) {\HH}^i ({X}_K^{G_1}, \omega^\kappa(-D)) \rightarrow
{\HH}^i(M_I \otimes^{\mathbf{L}}_{\Lambda_I, \kappa} \ZZ_p)$$ is an
isomorphism for $i =0$ and injective for $i=1$. \end{enumerate}
\end{thm}

We will deduce the theorem from a number of intermediate results.  In
particular, we need to analyze the Hecke operators $U_{w,1}$ and~$U_{w,2}$ and relate them to $T_{w,1}$ and~$T_w$ in order to be able to use the results of~\S\ref{subsec: vanishing for ordinary cohomology}.

\subsubsection{Reduction of the correspondence modulo $p$} Let $w \in
I^c$.   We begin by considering the special fibre of the
correspondence over~$\mathfrak{X}_{K,\Kli(p)}^I$ underlying the
operator $U_{w,1}$; we write
$C_{w,1}(p)_1$ for this
special fibre.  By reduction modulo~$p$, it follows
from Lemma~\ref{lem: base change of correspondences} that for each classical algebraic weight~$\kappa$ we obtain a
cohomological correspondence  which we continue to denote by $U_{w,1}: u_2^{*}
(\omega^{\kappa}\vert_{X^I_{K,\Kli}(p)_1}) \rightarrow u_1^! (
\omega^\kappa\vert_{X^I_{K,\Kli}(p)_1})$.

\begin{lem}\label{lem: U1 commutes with Hasse}For any
  place~$w\in I^c$, we
  have a commutative diagram

\begin{eqnarray*}
\xymatrix{ u_2^{*} \omega^\kappa \ar[r]^{U_{w,1}} \ar[d]^{ u_2^{*} {\Ha}(\cG_w)}& u_1^! \omega^\kappa  \ar[d]^{u_1^{*} {\Ha}(\cG_w)}\\
u_2^{*} (\omega^\kappa \otimes \det \omega_{\cG_w}^{p-1}) \ar[r]^{U_{w,1}} & u_1^! (\omega^\kappa \otimes \det \omega_{\cG_w}^{p-1})}
\end{eqnarray*}

\end{lem} 
\begin{proof}Since the kernel of $u_1^*\cG\to u_2^*\cG$ is \'etale,
  and the formation of~$\Ha(\cG_w)$ commutes with \'etale isogenies,
  this is immediate. \end{proof}

We now consider the operator~$U_{w,2}$ on ~$X^I_{K,\Kli}(p)_1$,
where~$w$ is any place lying over~$p$. Taking the special fibres of the
correspondences of~\S\ref{subsubsec:
Klingen type padic Hecke operators} with~$m=1$, we have
a correspondence \begin{eqnarray*}     
 \xymatrix{ &C_{w,2,1}(p)_1  \ar[rd]^{r_1}\ar[ld]_{r_2}&  \\
 X^I_{K'',\Kli}(p)_1 & & X^I_{K,\Kli}(p)_1}
  \end{eqnarray*} and by Lemma~\ref{lem: base change of
  correspondences}, a cohomological correspondence
$r_2^*\omega^\kappa|_{X^I_{K'',\Kli}(p)_1} \to
r_1^!\omega^\kappa|_{X^I_{K,\Kli}(p)_1}$; and a correspondence
\begin{eqnarray*}     
 \xymatrix{ &C_{w,2,2}(p)_1   \ar[rd]^{s_2}\ar[ld]_{s_1}&  \\
 X^I_{K,\Kli}(p)_1 & & X^I_{K'',\Kli}(p)_1}
  \end{eqnarray*}
and again by Lemma~\ref{lem: base change of
  correspondences}, a cohomological correspondence  $s_2^*\omega^\kappa|_{X^I_{K,\Kli}(p)_1} \to
s_1^!\omega^\kappa|_{X^I_{K'',\Kli}(p)_1}$.

We can associate to these cohomological  correspondences Hecke operators which we
denote as before as $$U'_{w} \in
\mathrm{Hom}(\mathrm{R}\Gamma(X^I_{K'',\Kli}(p)_1 ,
\omega^\kappa), \mathrm{R}\Gamma(X^I_{K,\Kli}(p)_1
, \omega^\kappa)),$$
 $$U''_{w} \in
 \mathrm{Hom}(\mathrm{R}\Gamma(X^I_{K,\Kli}(p)_1,
 \omega^\kappa), \mathrm{R}\Gamma(
 X^I_{K'',\Kli}(p)_1, \omega^\kappa)).$$ We
 continue to write $U_{w,2} = U'_{w} \circ U''_{w}$.

 \begin{lem}
   \label{lem: U operators commute with Hasse}For any $w |p$, we have commutative
   diagrams 
\begin{eqnarray*}
\xymatrix{ r_2^{*} \omega^\kappa \ar[r]^{U'_{w}} \ar[d]^{r_2^{*} {\Ha}(\cG_w)}& r_1^! \omega^\kappa  \ar[d]^{r_1^{*} {\Ha}(\cG_w)}\\
r_2^{*} (\omega^\kappa \otimes \det \omega_{\cG_w}^{p-1}) \ar[r]^{U'_w} & r_1^! (\omega^\kappa \otimes \det \omega_{\cG_w}^{p-1})}
\end{eqnarray*}

\begin{eqnarray*}
\xymatrix{ s_1^{*} \omega^\kappa \ar[r]^{U''_{w}} \ar[d]^{s_1^{*} {\Ha}(\cG_w)}& s_2^! \omega^\kappa  \ar[d]^{s_2^{*} {\Ha}(\cG_w)}\\
s_1^{*} (\omega^\kappa \otimes \det \omega_{\cG_w}^{p-1}) \ar[r]^{U''_w} & s_2^! (\omega^\kappa \otimes \det \omega_{\cG_w}^{p-1})}
\end{eqnarray*}
 \end{lem}
 \begin{proof} See  \cite[Lem.\ 10.5.2.1]{pilloniHidacomplexes}.
 \end{proof}

\subsubsection{Reduction of the correspondences to the non-ordinary locus}
 
 \begin{lem}
   \label{lem: Hasse invariant not zero divisor on Klingen
     correspondences}For any $w |p$, the Hasse invariant~$\Ha(\cG_w)$ is not a zero divisor
   on each of~$C_{w,2,1}(p)_1 $ and~$C_{w,2,2}(p)_1$.
 \end{lem}
 \begin{proof} See \cite[Lem.\ 10.5.2.2]{pilloniHidacomplexes}.
    \end{proof}

We now assume $w \in I$ (otherwise the schemes we consider would be empty) and  consider the rank one locus at $w$,~$X^{I,=_w1}_{K,\Kli}(p)_1$, which by definition is the vanishing locus of~$\Ha(\cG_w)$ in
~$X^{I}_{K,\Kli}(p)_1$.  Taking the zero locus of~$\Ha(\cG_w)$ at all entries of
the correspondences~$C_{w,1}(p)_1$, $C_{w,2,2}(p)_1$ and
$C_{w,2,1}(p)_1$ (and taking into account  Lemmas \ref{lem: U1
  commutes with Hasse} and~\ref{lem: U operators commute with Hasse}), we obtain 
correspondences 
\begin{eqnarray*}     
 \xymatrix{ &C^{=_w1}_{w,1}(p)_1  \ar[rd]^{u_1}\ar[ld]_{u_2}&  \\
 X^{I,=_w1}_{K,\Kli}(p)_1 & & X^{I,=_w1}_{K,\Kli}(p)_1}
  \end{eqnarray*} 
\begin{eqnarray*}     
 \xymatrix{ &C^{=_w1}_{w,2,1}(p)_1  \ar[rd]^{r_1}\ar[ld]_{r_2}&  \\
 X^{I,=_w1}_{K'',\Kli}(p)_1 & & X^{I,=_w1}_{K,\Kli}(p)_1}
  \end{eqnarray*} 
\begin{eqnarray*}     
 \xymatrix{ &C^{=_w1}_{w,2,2}(p)_1   \ar[rd]^{s_2}\ar[ld]_{s_1}&  \\
 X^{I,=_w1}_{K,\Kli}(p)_1 & & X^{I,=_w1}_{K'',\Kli}(p)_1}
  \end{eqnarray*}

By Lemmas~\ref{lem: base change of correspondences} and~\ref{lem: Hasse invariant not zero divisor on Klingen
     correspondences}, we also obtain
cohomological correspondences 
$$u_2^*\omega^\kappa|_{X^{I,=_w1}_{K,\Kli,1}(p)} \to
u_1^!\omega^\kappa|_{X^{I,=_w1}_{K,\Kli}(p)_1}, \quad
r_2^*\omega^\kappa|_{X^{I,=_w1}_{K'',\Kli}(p)_1} \to
r_1^!\omega^\kappa|_{X^{I,=_w1}_{K,\Kli}(p)_1},$$
 and  
 $$s_2^*\omega^\kappa|_{X^{I,=_w1}_{K,\Kli}(p)_1} \to
s_1^!\omega^\kappa|_{X^{I,=_w1}_{K'',\Kli}(p)_1}.$$

We can associate to these cohomological  correspondences Hecke
operators which we again
write as $$U_{w,1} \in
\mathrm{Hom}(\mathrm{R}\Gamma(X^{I,=_w1}_{K,\Kli}(p)_1 ,
\omega^\kappa), \mathrm{R}\Gamma(X^{I,=_w1}_{K,\Kli}(p)_1
, \omega^\kappa)),$$
$$U'_{w} \in
\mathrm{Hom}(\mathrm{R}\Gamma(X^{I,=_w1}_{K'',\Kli}(p)_1 ,
\omega^\kappa), \mathrm{R}\Gamma(X^{I,=_w1}_{K,\Kli}(p)_1
, \omega^\kappa)),$$
 $$U''_{w} \in
 \mathrm{Hom}(\mathrm{R}\Gamma(X^{I,=_w1}_{K,\Kli}(p)_1,
 \omega^\kappa), \mathrm{R}\Gamma(
 X^{I,=_w1}_{K'',\Kli}(p)_1, \omega^\kappa)).$$ We of course
 continue to write $U_{w,2} = U'_{w} \circ U''_{w}$. 

By Lemmas~\ref{lem: U1 commutes with Hasse} and~\ref{lem: U operators commute with Hasse}, the long exact
sequence
\[ H^*(X^{I}_{K,\Kli}(p)_1 ,
  \omega^\kappa)\stackrel{\times\Ha(\cG_w)}{\to}H^*(X^{I}_{K,\Kli}(p)_1
  , \omega^\kappa \otimes \det \omega_{\cG_w}^{p-1} )\to H^*(X^{I,=_w1}_{K'',\Kli}(p)_1
  , \omega^\kappa \otimes \det \omega_{\cG_w}^{p-1} )\]
is $U_{w,2}$- and $U_{w,1}$-equivariant.

\begin{lem}\label{lem: U commutes with secondary Hasse}We have commutative
  diagrams 
\begin{eqnarray*}
\xymatrix{ u_2^{*} \omega^\kappa|_{ X^{I,=_w1}_{K,\Kli}(p)_1} \ar[r]^{U_{w,1}} \ar[d]^{u_2^{*} {\Ha}'(\cG_w)}& u_1^! \omega^\kappa|_{ X^{I,=_w1}_{K,\Kli}(p)_1}  \ar[d]^{u_1^{*} {\Ha}'(\cG_w)}\\
u_2^{*} (\omega^\kappa \otimes \det \omega_{\cG_w}^{p^2-1}|_{ X^{I,=_w1}_{K,\Kli}(p)_1}) \ar[r]^{U_{w,1}} & r_1^! (\omega^\kappa \otimes \det \omega_{\cG_w}^{p^2-1}|_{ X^{I,=_w1}_{K,\Kli}(p)_1})}
\end{eqnarray*}

\begin{eqnarray*}
\xymatrix{ r_2^{*} \omega^\kappa|_{ X^{I,=_w1}_{K'',\Kli}(p)_1} \ar[r]^{U'_{w}} \ar[d]^{r_2^{*} {\Ha}'(\cG_w)}& r_1^! \omega^\kappa|_{ X^{I,=_w1}_{K,\Kli}(p)_1}  \ar[d]^{r_1^{*} {\Ha}'(\cG_w)}\\
r_2^{*} (\omega^\kappa \otimes \det \omega_{\cG_w}^{p^2-1}|_{ X^{I,=_w1}_{K'',\Kli}(p)_1}) \ar[r]^{U'_w} & r_1^! (\omega^\kappa \otimes \det \omega_{\cG_w}^{p^2-1}|_{ X^{I,=_w1}_{K,\Kli}(p)_1})}
\end{eqnarray*}

\begin{eqnarray*}
\xymatrix{ s_1^{*} \omega^\kappa|_{ X^{I,=_w1}_{K,\Kli}(p)_1} \ar[r]^{U''_{w}} \ar[d]^{s_1^{*} {\Ha}'(\cG_w)}& s_2^! \omega^\kappa|_{ X^{I,=_w1}_{K'',\Kli}(p)_1}  \ar[d]^{s_2^{*} {\Ha}'(\cG_w)}\\
s_1^{*} (\omega^\kappa \otimes \det \omega_{\cG_w}^{p^2-1}|_{ X^{I,=_w1}_{K,\Kli}(p)_1}) \ar[r]^{U''_w} & s_2^! (\omega^\kappa \otimes \det \omega_{\cG_w}^{p^2-1}|_{ X^{I,=_w1}_{K'',\Kli}(p)_1})}
\end{eqnarray*}
\end{lem}
\begin{proof} See \cite[Lem.\ 10.5.3.1]{pilloniHidacomplexes} .
  \end{proof}

\subsubsection{Comparison of $U_{w,2}, U_{w,1}$ and $T_w, T_{w,1}$ in a special case}  

We fix  $ J \subset I$.    The space $X^{I, =_J 1, =_{J^c} 2}_{K_p,
  \Kli}(p)_1$ carries a finite   \'etale map  to the space  $X^{I, =_J 1,
  =_{J^c} 2}_{K,1}$ studied in~\S\ref{subsubsec: commutativity over ordinary locus}. This map is given by forgetting the
multiplicative groups $H_v$ of order $p$ at the places $v \in J^c$. Therefore, it has degree $(p+1)^{ \# J^c}$. Let $\kappa = (k_v, l_v)$ be a classical algebraic weight. We have an injective map 
$$ \HH^0( X^{I, =_J 1, =_{J^c} 2}_{K,1}, \omega^\kappa(-D)) \rightarrow  \HH^0( X^{I, =_J 1, =_{J^c} 2}_{K,\Kli}(p)_1, \omega^\kappa(-D)).$$

We assume that $l_v \geq 3$ if $v \in I^c$, that $k_v \geq 3$, $l_v
\geq 2$ if $v \in I$, and moreover that $l_v \geq p+1$ and $k_v \geq 2p+3$ if $v \in J$. On the left hand side, we have an action of $T_{w}$ for $w |p$  and $T_{w,1}$ for $w \in I^c$. On the right hand side, we have an action of $U_{w,2}$  for $w |p$ and $U_{w,1}$ for $w \in I^c$.  This follows from the fact that all these Hecke operators have been proved to commute with the Hasse invariants  (by Lemmas \ref{lem: Siegel sheaves commute},   \ref{lem-commutTW}, \ref{lem: commutativity of etale Klingen correspondences}, \ref{lem: U1 commutes with Hasse}, \ref{lem: U operators commute with Hasse}, \ref{lem: U commutes with secondary Hasse}).

The main result of this subsection is:

\begin{prop}\label{prop-comparisonTUprimitive}There is a universal  constant $C$
  depending only on $p$ and $F$ but not on the tame level $K^p$ such
  that if~$l_w\ge 2$ for all~$w$, $k_w-l_w \geq C$ for all $w | p$,      and $l_w \geq C$ for all $w \in I^c$, then:  \begin{enumerate}
\item The operator $U^I = \prod_{w |p} U_{w,2} \prod_{w \in I^c} U_{w,1}$ is locally finite on $$\HH^0( X^{G_1, I, =_J 1, =_{J^c} 2}_{K, \Kli}(p)_1, \omega^\kappa(-D)).$$
\item\label{item: T to U map is an iso} Let $\widetilde{T}^I = \prod_{w |p} T_{w} \prod_{w \in I^c} T_{w,1}$.   The map  $$ e(\widetilde{T}^I) \HH^0( X^{G_1, I, =_J 1, =_{J^c} 2}_{K,1}, \omega^\kappa(-D)) \rightarrow  e( U^I ) \HH^0( X^{G_1, I, =_J 1, =_{J^c} 2}_{K,\Kli}(p)_1, \omega^\kappa(-D))$$ is an isomorphism.
\item This isomorphism is equivariant for the action of $T_{w,1}$ on the left and $U_{w,1}$ on the right for all $w \in I^c$ and of $T_{w}$ and $U_{w,2}$ for all $w \in J$. 
\end{enumerate}
\end{prop}

This result establishes a first relation between the cohomology  at Klingen level and spherical level and will  allow us to reduce a big proportion of  the proof of Theorem~ \ref{theorem-p-adic-complex} to Theorem~ \ref{thm-vanishing}.

\begin{rem} One can interpret this result as saying that the ordinarity condition prevents the existence of ``newforms'' of Klingen level. 
\end{rem}

We have a
 finite  \'etale map: \[X^{I, =_J 1, =_{J^c} 2}_{K,
    \Kli}(p)_1 \rightarrow  
    (X^{I, =_J 1, =_{J^c} 2}_{K})_1\]  which parametrizes  multiplicative subgroups of order $p$,  $H_w \subset
\cG_w[p]$ for all $w \in J^c$. We introduce various Hecke operators that decrease the level at places $w \in J^c$ and compare them with our existing Hecke operators.

We first define a  correspondence for each $w \in J^c$ (with~$X_w$ defined below), 
\begin{eqnarray*}     
 \xymatrix{ &X_w \ar[rd]^{x_1}\ar[ld]_{x_2}&  \\
 (X^{I, =_J 1, =_{J^c} 2}_{K})_1 & & (X^{I, =_J 1, =_{J^c} 2}_{K})_1}
  \end{eqnarray*}as follows. We let $ x_1: X_w \rightarrow (X^{I, =_J 1,
    =_{J^c} 2}_{K})_1$ be the natural forgetful map,
  where~$X_w$ parametrizes subgroups $  L_w \subset \cG_w[p^2]$, where  $L_w$ is  totally
  isotropic of \'etale rank $p^3$ and multiplicative rank $p$.  A
  standard computation shows that $x_1$ is finite flat. For a suitable
  choice of polyhedral decomposition, there is a map $x_2: X_w
  \rightarrow (X^{I, =_J 1, =_{J^c} 2}_{K})_1$, which on the
  $p$-divisible group is given by $\cG \mapsto \cG/L_w$ (since we are
  only dealing with $\HH^0$ cohomology groups, we will for the most part
 suppress the discussion of the boundary in this section).

We can define an operator, using the usual procedure,  associated to
$X_w$,  $\widetilde{T}_w  \in  \mathrm{End} (\HH^0 ( (X^{I, =_J 1,
  =_{J^c} 2}_{K})_1, \omega^\kappa(-D)))$.

 If we denote by $X^{I, =_J 1, =_{J^c} 2}_{K,
  \Kli_w}(p)_1 \rightarrow (X^{I, =_J 1,
  =_{J^c} 2}_{K})_1$ the finite \'etale cover that parametrizes subgroups $H_w \subset \cG_w[p]$ of order $p$, then we observe that the projection
$x_2$ lifts to a map $X_w \rightarrow X^{I, =_J 1, =_{J^c} 2}_{K,
  \Kli}(p)_1$ by sending $(\cG, L_w)$ to $(\cG/L_w, \cG_w[p]/L_w)$ and
therefore one can promote $\widetilde{T}_w$ to  maps
$\widetilde{T}'_w$ and $\widetilde{T}''_w$ fitting in  a commutative
diagram (where vertical maps are the injections given by the obvious pull back maps): 

\numequation\label{eqn: a label we have yet to use}
\xymatrix{\HH^0( X^{I, =_J 1, =_{J^c} 2}_{K, \Kli_w}(p)_1, \omega^\kappa(-D)) \ar[r]^{\widetilde{T}''_w} \ar[rd]^{\widetilde{T}'_w} & \HH^0 ( X^{I, =_J 1, =_{J^c} 2}_{K, \Kli_w}(p)_1, \omega^\kappa(-D)) \\
\HH^0 ( (X^{I, =_J 1, =_{J^c} 2}_{K})_1, \omega^\kappa(-D))\ar[u] \ar[r]^{\widetilde{T}_w} & \HH^0 ( (X^{I, =_J 1, =_{J^c} 2}_{K})_1, \omega^\kappa(-D)) \ar[u]}
\end{equation}

On the other hand we have already defined a Hecke operator $T_w$ on $\HH^0 ( (X^{I, =_J 1, =_{J^c} 2}_{K})_1, \omega^\kappa(-D))$.
We also have a
chain of finite  \'etale maps: \[X^{I, =_J 1, =_{J^c} 2}_{K,
    \Kli}(p)_1 \rightarrow  X^{I, =_J 1, =_{J^c} 2}_{K, \Kli_w}(p)_1
  \rightarrow  (X^{I, =_J 1, =_{J^c} 2}_{K})_1\] where the first map
forgets the multiplicative subgroup of order $p$,  $H_{w'} \subset
\cG_{w'}[p]$ for $w \in J^c \setminus \{w\}$.  We have defined an operator $U_{w,2}$ on $\HH^0( X^{I, =_J 1, =_{J^c} 2}_{K,
    \Kli}(p)_1, \omega^{\kappa}(-D))$, but clearly it descends to an operator on $\HH^0( X^{I, =_J 1, =_{J^c} 2}_{K,
    \Kli_w}(p)_1, \omega^{\kappa}(-D))$ because only the Klingen level structure at $w$ matters in the definition of $U_{w,2}$.

\begin{lem}\label{lem: UT=UU} Assume that  $k_w-l_w \geq 1$. \begin{enumerate}
\item We have $\widetilde{T}_w = T_w$.
\item We have $U_{w,2} \circ \widetilde{T}''_{w} = U_{w,2} \circ U_{w,2}$. 
\end{enumerate}
\end{lem}

\begin{proof} See \cite[Lem.\ 11.1.1.1, Lem.\
  11.1.1.3]{pilloniHidacomplexes}.\end{proof}

\begin{lem}
  \label{lem: U2 equals T on w equals 1}Let $w \in J$. The canonical map $ \HH^0( X^{G_1, I, =_J 1, =_{J^c} 2}_{K,1}, \omega^\kappa(-D)) \rightarrow  \HH^0( X^{G_1, I, =_J 1, =_{J^c} 2}_{K,\Kli}(p)_1, \omega^\kappa(-D))$ intertwines the actions of $T_w$ and $U_{w,2}$.
\end{lem}
\begin{proof}
  By Lemma~\ref{lem: commutativity of etale Klingen correspondences}
  we have $T_w=T'_{w,et}\circ T''_{w,et}$, which corresponds to~$U_{w,2}$ by definition of the right hand side.
\end{proof}

\begin{lem}\label{lem: U1 equals T1 on w equals 1} Let $w \in I^c$. The canonical map $ \HH^0( X^{G_1, I, =_J 1, =_{J^c} 2}_{K,1}, \omega^\kappa(-D)) \rightarrow  \HH^0( X^{G_1, I, =_J 1, =_{J^c} 2}_{K,\Kli}(p)_1, \omega^\kappa(-D))$ intertwines the actions of $T_{w,1}$ and $U_{w,1}$.
\end{lem}
\begin{proof} By Lemma~ \ref{lem: Siegel sheaves commute over ordinary}, we have $T_{w,1}=T^{et}_{w,1}$, which corresponds to~$U_{w,1}$ by definition on the right hand side.
\end{proof}

\begin{cor}\label{cor: UI is locally finite and surjective} Suppose
   that
  we have  ~$l_w\ge 2$ and $k_w-l_w \geq 1$ for all~$w\in  J^c$. Then the action
  of ~$\prod_{w \in J^c} U_{w,2}$ on  $\HH^0( X^{G_1, I, =_J 1, =_{J^c}
    2}_{K,\Kli}(p)_1, \omega^\kappa(-D))$ is locally finite,   the action of $\prod_{w \in J^c} \tilde{T}_w$ is locally finite on $\HH^0( X^{G_1, I, =_J 1,
    =_{J^c} 2}_{K,1}, \omega^\kappa(-D))$, and the map $$ e( \prod_{w \in J^c} \tilde{T}_w) \HH^0( X^{G_1, I, =_J 1,
    =_{J^c} 2}_{K,1}, \omega^\kappa(-D)) \rightarrow  e( \prod_{w \in J^c} U_{w,2} ) \HH^0(
  X^{G_1, I, =_J 1, =_{J^c} 2}_{K,\Kli}(p)_1, \omega^\kappa(-D))$$ is
  surjective. \end{cor}
\begin{proof}Combining the diagrams~\eqref{eqn: a label we have yet to use} for
all~$w\in J^c$,  we see that there is  a
  commutative diagram: \begin{eqnarray*}
\xymatrix{\HH^0( X^{G_1, I, =_J 1, =_{J^c} 2}_{K, \Kli}(p)_1, \omega^\kappa(-D)) \ar[rr]^{\prod_{w \in J^c} \tilde{T}_w'' } \ar[rrd]^{\prod_{w \in J^c} \tilde{T}_w'} && \HH^0 ( X^{G_1, I, =_J 1, =_{J^c} 2}_{K, \Kli}(p)_1, \omega^\kappa(-D)) \\
\HH^0 ( X^{G_1, I, =_J 1, =_{J^c} 2}_{K,1}, \omega^\kappa(-D))\ar[u]^{\iota} \ar[rr]^{\prod_{w \in J^c} \tilde{T}_w} && \HH^0 ( X^{G_1, I, =_J 1, =_{J^c} 2}_{K,1}, \omega^\kappa(-D)) \ar[u]_{\iota}}
\end{eqnarray*}where the vertical maps~$\iota$ are the natural
injections. We deduce from Lemma~\ref{lem: UT=UU} that $ \prod_{w \in J^c} U_{w,2}\circ(\prod_{w \in J^c} \tilde{T}_w)''=  \prod_{w \in J^c} U_{w,2}\circ  \prod_{w \in J^c} U_{w,2}$. The argument now follows the proofs of~ \cite[Cor.\ 11.1.1.1, Cor.\
11.1.1.2]{pilloniHidacomplexes}. We deduce
that~$\prod_{w \in J^c} \tilde{T}_w$ acts locally finitely on $\HH^0 ( X^{G_1, I,
  =_J 1, =_{J^c} 2}_{K,1}, \omega^\kappa(-D))$ by Lemma~\ref{lem-commutTW}. It follows that for any~$f\in \HH^0( X^{G_1, I, =_J 1, =_{J^c}
  2}_{K, \Kli}(p)_1, \omega^\kappa(-D)) $, there is a
$\prod_{w \in J^c} \tilde{T}_w$-stable finite-dimensional vector space~$V$
containing~$(\prod_{w \in J^c} \tilde{T}_w)'f$, and then the subspace of $\HH^0(
X^{G_1, I, =_J 1, =_{J^c} 2}_{K, \Kli}(p)_1, \omega^\kappa(-D))$
spanned by $\iota(V)$, $ \prod_{w \in J^c} U_{w,2}\iota(V)$, $f$, and~$ \prod_{w \in J^c} U_{w,2}f$ is
finite-dimensional and $ \prod_{w \in J^c} U_{w,2}$-stable, so $ \prod_{w \in J^c} U_{w,2}$ acts locally finitely,
as claimed.

  To prove the claimed surjectivity, if~$f\in e( \prod_{w \in J^c} U_{w,2})\HH^0( X^{G_1, I, =_J 1, =_{J^c}
  2}_{K, \Kli}(p)_1, \omega^\kappa(-D))$, then one checks from the
definitions that we have \[f=e( \prod_{w \in J^c} U_{w,2})\iota(e(\prod_{w \in J^c} \tilde{T}_w)(\prod_{w \in J^c} \tilde{T}_w')( \prod_{w \in J^c} U_{w,2})^{-1}f).\qedhere\]
  \end{proof}

It remains to prove the injectivity of the map considered in
Proposition~\ref{prop-comparisonTUprimitive}~\eqref{item: T to U map
  is an iso}. This will be done by exhibiting an inverse up to a
certain power of $p$. For this reason, it is necessary to lift the
situation to a $\ZZ_p$-flat base. This is done by considering certain formal schemes.  
Let us denote by   $\mathfrak{X}^{G_1, I, =_J 1, =_{J^c} 2}_{K}$ the formal completion of $\mathfrak{X}^{G_1, I, =_{J^{c}} 2}_K$ along $X^{G_1, I, =_J 1, =_{J^c} 2}_{K,1}$. We also denote by   $\mathfrak{X}^{G_1, I, =_J 1, =_{J^c} 2}_{K,\Kli}(p)$ the formal completion of $\mathfrak{X}^{G_1, I, =_{J^c} 2}_{K, \Kli}(p)$ along $X^{G_1, I, =_J 1, =_{J^c} 2}_{K,\Kli}(p)_1$. We denote by $I$  the ideal of definition of these formal schemes. Observe that $p \in I$ and that 
$ I/p = ( \mathrm{Ha}(\cG_w) \det \omega_{\cG_w}^{(1-p)},~w \in J)$.

We consider the modules $$\HH^0( \mathfrak{X}^{G_1, I, =_J 1, =_{J^c} 2}_{K}, \omega^\kappa(-D))~~\textrm{and}~~ \HH^0( \mathfrak{X}^{G_1, I, =_J 1, =_{J^c} 2}_{K, \Kli}(p), \omega^\kappa(-D)),$$ which are $I$-adically complete and separated, and also  $\ZZ_p$-flat. Moreover, the natural map $\HH^0( \mathfrak{X}^{G_1, I, =_J 1, =_{J^c} 2}_{K}, \omega^\kappa(-D)) \rightarrow  \HH^0( \mathfrak{X}^{G_1, I, =_J 1, =_{J^c} 2}_{K, \Kli}(p), \omega^\kappa(-D))$ reduces modulo $I$ to the map  $\HH^0( ({X}^{G_1, I, =_J 1, =_{J^c} 2}_{K})_1, \omega^\kappa(-D)) \rightarrow  \HH^0( {X}^{G_1, I, =_J 1, =_{J^c} 2}_{K, \Kli}(p)_1, \omega^\kappa(-D))$ (we are using here that $({X}^{G_1, I, =_J 1, =_{J^c} 2}_{K})_1$ and ${X}^{G_1, I, =_J 1, =_{J^c} 2}_{K, \Kli}(p)_1$ have affine image in the minimal compactification and thus that higher  cuspidal cohomology over these spaces vanishes). 

We can lift the map $$\tilde{T}_w : \HH^0( ({X}^{G_1, I, =_J 1, =_{J^c} 2}_{K})_1, \omega^\kappa(-D)) \rightarrow \HH^0( ({X}^{G_1, I, =_J 1, =_{J^c} 2}_{K})_1, \omega^\kappa(-D))$$ to a map  $ \tilde{T}_w : \HH^0( \mathfrak{X}^{G_1, I, =_J 1, =_{J^c} 2}_{K}, \omega^\kappa(-D)) \rightarrow \HH^0( \mathfrak{X}^{G_1, I, =_J 1, =_{J^c} 2}_{K}, \omega^\kappa(-D))$ (the correspondences $X_w$ lift to correspondences on the formal schemes).

There is a trace map 
$$\mathrm{Tr} : \HH^0( \mathfrak{X}^{G_1, I, =_J 1, =_{J^c} 2}_{K, \Kli}(p), \omega^\kappa(-D)) \rightarrow \HH^0( \mathfrak{X}^{G_1, I, =_J 1, =_{J^c} 2}_{K}, 
\omega^\kappa(-D))$$
 associated to the finite \'etale map $\mathfrak{X}^{G_1, I, =_J 1, =_{J^c} 2}_{K, \Kli}(p) \rightarrow \mathfrak{X}^{G_1, I, =_J 1, =_{J^c} 2}_{K}$.

\begin{lem}\label{lemma-new-lemma} For any $n \in \ZZ_{\geq 1}$,  we have the congruence $$\mathrm{Tr} \circ  (\prod_{w \in J^c} U_{w,2})^n \circ  \iota (f) \equiv p^{\# J^c} (\prod_{w \in J^c} \tilde{T}_w)^n( f)\pmod{ p^{\inf_{w \in J^c} k_w-l_w}}$$ for any $f \in \HH^0 ( \mathfrak{X}^{G_1, I, =_J 1, =_{J^c} 2}_{K}, \omega^\kappa(-D))$. 
\end{lem}

\begin{proof}  We have  $$ \mathrm{Tr} \circ  (\prod_{w \in J^c}   U_{w,2})^n \circ  \iota ( f( \cG, \omega ) ) = \frac{1}{p^{n\sum_{w \in J^c} l_w + 3}}\sum_{w \in J^c} \sum_{L_{w,n}} (\sum_{H_w} f (\cG/ (\oplus_w L_{w,n}) , \omega' )$$ where $H_w$ runs over all multiplicative subgroups of rank $p$ of $\cG_w[p]$ and $L_{w,n}$ runs over all totally isotropic subgroups of order $p^{3n}$ of $\cG_w[p^{2n}]$, with trivial intersection with $H_w$ (this implies that $L_{w,n}$ is locally in the \'etale topology an extension of $\Z/p^n\Z \oplus \Z/p^{n^2}\Z$ by $\mu_{p^n}$), and where $\omega$ is a trivialization of $\omega_{\cG}$ and $\omega'$ is a rational trivialization of $\omega_{\cG/(\oplus_w L_{w,n})}$, defined by the condition that  $\pi^\star \omega' = \omega$ for the isogeny $\cG \rightarrow \cG/(\oplus L_{w,n})$. Given a group $L_{w,n}$, we can find $p$ subgroups $H_w$ of order $p$ and of multiplicative type such that $L_{w,n} \cap H_w = \{0\}$. This means that the groups $L_{w,n}$ in  the formula defining $ \mathrm{Tr} \circ  (\prod_{w \in J^c}  U_{w,2}) \circ  \iota (f) $ occur with multiplicity $p^{\# J^c}$. 
On the other hand, $$\prod_{w \in J^c} \tilde{T}_w f(\cG, \omega) = \frac{1}{p^{\sum_{w \in J^c} l_w + 3}} \sum_{w \in
  J^c} \sum_{L_w} f(\cG/(\oplus_wL_w), \omega')$$ where $L_w$ runs over all totally
isotropic subgroups of order $p^3$  of $\cG_w[p^2]$ with
multiplicative rank $1$.  Now we observe that  $$(\prod_{w \in J^c} \tilde{T}_w)^n f(\cG, \omega) \equiv \frac{1}{p^{n\sum_{w \in J^c} l_w + 3}} \sum_{w \in
  J^c} \sum_{L_{w,n}} f(\cG/(\oplus_w L_{w,n}), \omega')\pmod {p^{\inf_{w \in J^c} k_w-l_w} }$$
  where $L_{w,n}$ runs over all totally isotropic subgroups of order $p^{3n}$ of $\cG_w[p^{2n}]$, which are  locally in the \'etale topology an extension of $\Z/p^n\Z \oplus \Z/p^{n^2}\Z$ by $\mu_{p^n}$.  Indeed, if we write $$(\prod_{w \in J^c} \tilde{T}_w)^n f(\cG, \omega) = \frac{1}{p^{n\sum_{w \in J^c} l_w + 3}} \sum_{w \in
  J^c} \sum_{L'_{w,n}} f(\cG/(\oplus_w L'_{w,n}), \omega')$$ using the definition of $\prod_{w \in J^c} \tilde{T}_w$, we find  that all the groups $L_{w,n}$ appear exactly one time among the groups $L'_{w,n}$,  and that all the remaining groups precisely contain the multiplicative subgroup of $\cG_w[p]$ (and these give a contribution divisible by $p^{(k_w-l_w)}$). 
\end{proof}

\begin{lem}
  \label{lem: ordinary projectors quasi isomorphic on ordinary locus
    XI XKliI} Assume that for all $w \in J^c$, we have $k_w-l_w >  p^{\#
    J^c}$. Then the natural map $$ e(\prod_{w \in J^c} \tilde{T}_w ) \HH^0( X^{G_1, I, =_J 1, =_{J^c} 2}_{K,1}, \omega^\kappa(-D)) \rightarrow   e(\prod_{w \in J^c} U_{w,2} ) \HH^0( X^{G_1, I, =_J 1, =_{J^c} 2}_{K,\Kli}(p)_1, \omega^\kappa(-D))$$ is bijective.
   \end{lem}

\begin{proof} We will show that the map: $$ \iota :  e(\prod_{w \in J^c} \tilde{T}_w  ) \HH^0( \mathfrak{X}^{G_1, I, =_J 1, =_{J^c} 2}_{K}, \omega^\kappa(-D)) \rightarrow   e( \prod_{w \in J^c} U_{w,2} ) \HH^0( \mathfrak{X}^{G_1, I, =_J 1, =_{J^c} 2}_{K,\Kli}(p), \omega^\kappa(-D))$$
is bijective (note that it is legitimate to apply the ordinary
projectors on these spaces, because they can be written
 as projective limits (modding out by  $I^n$) of spaces carrying a
locally finite action). The result will then follow by taking
reduction modulo $I$. The map is surjective by
Corollary~\ref{cor: UI is locally finite and surjective} (and using
$I$-adic approximation). It remains to prove injectivity. Let us take
$$f \in e(\prod_{w \in J^c} \tilde{T}_w  ) \HH^0( \mathfrak{X}^{G_1, I,
  =_J 1, =_{J^c} 2}_{K}, \omega^\kappa(-D)),$$
   with $f \neq 0$ and
$\iota(f) =0$. Without loss of generality, we can suppose that $$f
\notin pe(\prod_{w \in J^c} \tilde{T}_w  ) \HH^0( \mathfrak{X}^{G_1,
  I, =_J 1, =_{J^c} 2}_{K}, \omega^\kappa(-D)).$$  It follows from
Lemma~ \ref{lemma-new-lemma} that  $\mathrm{Tr} ( \iota f) \equiv p^{\#
  J^c} f \pmod {p^{\inf_w k_w -l_w}}$. Therefore, $f \in p^{ (\inf_w k_w
  -l_w)- \# J^c}e(\prod_{w \in J^c} \tilde{T}_w  ) \HH^0(
\mathfrak{X}^{G_1, I, =_J 1, =_{J^c} 2}_{K}, \omega^\kappa(-D))$. This
is a contradiction. 
\end{proof}

\begin{proof}[Proof of
  Proposition~\ref{prop-comparisonTUprimitive}]This is immediate from
  Corollary~\ref{cor: UI is locally finite and surjective},
  Lemma~\ref{lem: ordinary projectors quasi isomorphic on ordinary
    locus XI XKliI}, and Lemmas~\ref{lem: UT=UU},~\ref{lem: U2 equals T on w equals 1}
  and~\ref{lem: U1 equals T1 on w equals 1}. 
\end{proof}

\subsubsection{ Comparison of the cohomology on ${X}^{I, G_1}_{K, \Kli}(p)_1$ and ${X}^{I, G_1}_{K,1}$ }

We now deduce the following proposition.

\begin{prop}\label{prop: T to U in medium weight}There is a universal  constant $C$
  depending only on $p$ and~ $F$ but not on the tame level $K^p$ such
  that if~$l_w\ge 2$ for all~$w$, $k_w-l_w \geq C$ for all $w | p$,
  and $l_w \geq C$ for all $w \in I^c$, then the operator $U^I$ is locally finite on  $\RGamma({X}^{I, G_1}_{K, \Kli}(p)_1, \omega^\kappa(-D))$ and  there is a  canonical quasi-isomorphism: 
$$ e (\widetilde{T}^I) \RGamma({X}^{I, G_1}_{K,1},
\omega^\kappa(-D)) \rightarrow  e (U^I) \RGamma({X}^{I, G_1}_{K, \Kli}(p)_1,
\omega^\kappa(-D)).$$
\end{prop}

\begin{proof}
In~\S\ref{section-computingcohoexplicitely}, we constructed a complex $K^\bullet$ computing explicitly the cohomology $\RGamma({X}^{I, G_1}_{K,1}, \omega^\kappa(-D))$. We recall that $K^k =$ $$  \oplus_{ J \subset I, \# J = k} \colim_{\times \prod_{w \in J} \Ha(\cG_w)}\HH^0\big(X^{ G_1,I, \geq_{ J^c} 2}_{K,1}, \omega^{\kappa} (-D)\otimes \bigotimes_{w \in J} (\det \cG_w)^{n(p-1)}/(\sum_{w\in J} (\Ha(\cG_w)^n))\big).$$

In exactly the same way, there is a complex $L^\bullet$ computing $\RGamma({X}^{I, G_1}_{K, \Kli}(p)_1, \omega^\kappa(-D))$, such that $L^k = $
 $$  \oplus_{ J \subset I, \# J = k} \colim_{\times \prod_{w \in J}
   \Ha(\cG_w)}\HH^0\big(X^{ G_1,I, \geq_{ J^c} 2}_{K,\Kli}(p)_1,
 \omega^{\kappa} (-D)\otimes \bigotimes_{w \in J} (\det
 \cG_w)^{n(p-1)}/(\sum_{w\in J} (\Ha(\cG_w)^n))\big).$$

It therefore suffices to prove that for each~$J \subset I$ and each~$n\ge 1$, $U^I$ is
locally finite on \[\HH^0\big(X^{ G_1,I, \geq_{ J^c} 2}_{K,\Kli}(p)_1,
 \omega^{\kappa} (-D)\otimes \bigotimes_{w \in J} (\det
 \cG_w)^{n(p-1)}/(\sum_{w\in J} (\Ha(\cG_w)^n))\big),\] and the
map \[e(\widetilde{T}^I) \HH^0\big(X^{ G_1,I, \geq_{ J^c} 2}_{K,\Kli}(p)_1,
 \omega^{\kappa} (-D)\otimes \bigotimes_{w \in J} (\det
 \cG_w)^{n(p-1)}/(\sum_{w\in J} (\Ha(\cG_w)^n))\big) \rightarrow\]\[ e(U^I)\HH^0\big(X^{ G_1,I, \geq_{ J^c} 2}_{K,\Kli}(p)_1,
 \omega^{\kappa} (-D)\otimes \bigotimes_{w \in J} (\det
 \cG_w)^{n(p-1)}/(\sum_{w\in J} (\Ha(\cG_w)^n))\big) \]is an
isomorphism. In the case~$n=1$, this is
Proposition~\ref{prop-comparisonTUprimitive}, and the general case
follows by induction on~$n$, using the short exact sequence \[0\to \omega^{\kappa} (-D)\otimes \bigotimes_{w \in J} (\det
 \cG_w)^{(n-1)(p-1)}/(\sum_{w\in J} (\Ha(\cG_w)^{n-1}))\to\]\[ \omega^{\kappa} (-D)\otimes \bigotimes_{w \in J} (\det
 \cG_w)^{n(p-1)}/(\sum_{w\in J} (\Ha(\cG_w)^n))\to\]\[ \omega^{\kappa} (-D)\otimes \bigotimes_{w \in J} (\det
 \cG_w)^{n(p-1)}/(\sum_{w\in J} (\Ha(\cG_w)))\to 0 \]
and the acyclicity of these sheaves (for which see the proof of Proposition~\ref{prop: the complex K}).\end{proof}

\subsubsection{The proof of Theorem~\ref{theorem-p-adic-complex}} 

\begin{lem}\label{lem: control from Kli p infinity to Kli p} If
  $\kappa = (( k_v, l_v)_{v |p})$ is a classical algebraic weight with $l_v = 2$
  when $v \in I$ and $k_v\equiv l_v\equiv 2\pmod{p-1}$ for all $v|p$, then for each~$n\ge 2$ there is a diagonal map making a commutative diagram: 

\begin{eqnarray*}
\xymatrix{ \mathrm{R}\Gamma( \mathfrak{X}^{G_1, I}_{\Kli}(p^n), \omega^{\kappa}(-D)) \ar[r]^{U^I} \ar[rd]^{U^I} &  \mathrm{R}\Gamma( \mathfrak{X}^{G_1, I}_{\Kli}(p^n), \omega^{\kappa}(-D)) \\
 \mathrm{R}\Gamma( \mathfrak{X}^{G_1, I}_{\Kli}(p^{n-1}), \omega^{\kappa}(-D)) \ar[r]^{U^I} \ar[u] &  \mathrm{R}\Gamma( \mathfrak{X}^{G_1, I}_{\Kli}(p^{n-1}), \omega^{\kappa}(-D))  \ar[u] }
 \end{eqnarray*}
\end{lem}

\begin{proof} This is an easy computation in the Hecke
  algebra, see the proof of~ \cite[Thm.\ 11.3.1]{pilloniHidacomplexes}.\end{proof} 

We now make repeated use of Nakayama's lemma for complexes, in the
form of~\cite[Prop.\ 2.2.1, Prop.\ 2.2.2]{pilloniHidacomplexes}. In
fact, we need the following slight strengthening of~\cite[Prop.\
2.2.1]{pilloniHidacomplexes}, which is proved in the same way; for
ease of reference we explain how it follows from results in the
literature.

\begin{lem}
  \label{lem: finite cohomology implies perfect}Let $R$ be a complete
  local Noetherian ring with maximal ideal~$\m$, and let~$M^\bullet$
  be a bounded complex of $\m$-adically complete and separated, flat
  $R$-modules, with the property that the cohomology groups of
  $M^\bullet\otimes_RR/\m$ are finite-dimensional and concentrated
  in degrees~$[a,b]$. Then~$M^\bullet$ is a perfect complex,
  concentrated in degrees~$[a,b]$.
\end{lem}
\begin{proof}
  It follows from~\cite[Prop.\ 2.2.1]{pilloniHidacomplexes}
  that~$M^\bullet$ is a perfect complex, and it then follows
  from~\cite[Lem.\ 2.3, Cor.\ 2.7]{KT} that it is concentrated in
  degrees~$[a,b]$.
\end{proof}
All the complexes we consider below can be represented by
bounded complexes of flat, complete and separated $\Z_p$-modules
(resp.\ $\Lambda_I$-modules), as can be seen by considering a \v{C}ech
complex for any finite affine cover, so the hypotheses of
Lemma~\ref{lem: finite cohomology implies perfect} apply in our
situation.
\begin{lem}\label{lem: classicity at level of sheaf}   For all classical algebraic weights $\kappa = (( k_v,
  l_v)_{v |p})$  with $l_w\ge 2$ for all~$w$, $k_w-l_w \geq C$ for all $w | p$,
  and $l_w \geq C$ for all $w \in I^c$, the operator $U^I$ is locally finite on  $\RGamma(\mathfrak{X}^{G_1,I}_{K, \Kli}(p^\infty), \omega^\kappa(-D))$ and  there is a  canonical quasi-isomorphism: 
$$ e (\widetilde{T}^I) \RGamma(\mathfrak{X}^{I, G_1}_{K},
\omega^\kappa(-D)) \rightarrow  e (U^I) \RGamma(\mathfrak{X}^{I, G_1}_{K, \Kli}(p^\infty),
\omega^\kappa(-D)).$$
\end{lem}

\begin{proof} By Proposition~ \ref{prop: T to U in medium weight},
  together with~\cite[Prop.\ 2.2.2, Prop.\ 2.3.1]{pilloniHidacomplexes}, the action of $U^I$ is locally finite on  $\RGamma(\mathfrak{X}^{I, G_1}_{K, \Kli}(p),
\omega^\kappa(-D))$, and  the map $$ e (\widetilde{T}^I) \RGamma(\mathfrak{X}^{I, G_1}_{K},
\omega^\kappa(-D)) \rightarrow  e (U^I) \RGamma(\mathfrak{X}^{I, G_1}_{K, \Kli}(p),
\omega^\kappa(-D))$$  is a quasi-isomorphism. It follows easily from Lemma~ \ref{lem: control from Kli p infinity to Kli p} that $U^I$ is locally finite on $\RGamma(\mathfrak{X}^{I, G_1}_{K, \Kli}(p^\infty),
\omega^\kappa(-D))$  and that the map $e (U^I) \RGamma(\mathfrak{X}^{I, G_1}_{K, \Kli}(p),
\omega^\kappa(-D)) \rightarrow e (U^I) \RGamma(\mathfrak{X}^{I, G_1}_{K, \Kli}(p^\infty),
\omega^\kappa(-D))$ is a quasi-isomorphism, as required.
\end{proof}

 Let $\kappa = (k_v, l_v)$ be a classical algebraic weight. Let $K\omega^\kappa$
denote the kernel of the surjection of Corollary~\ref{cor: surjection for Komega}, so that over $\mathfrak{X}^I_{K,\Kli}(p^\infty)$, so we
have a short exact sequence of sheaves \[0\to
  K\omega^\kappa\to\omega^\kappa\to \Omega^\kappa\to 0. \] A key step
in the comparison between the ordinary forms of these weights is the
following basic lemma.

\begin{lem}
  \label{lem: comparison of p adic classical weights U operators}For
  any $w|p$, we
  have $U_{w,2}\in
  p\mathrm{End}( \mathrm{R}\Gamma(\mathfrak{X}^I_{K,\Kli}(p^\infty),K\omega^\kappa))$.\end{lem}
\begin{proof}
  This follows immediately from an examination
  of~(\ref{eqn: local description of U prime on differentials}).\end{proof}

\begin{lem}\label{lem-classicity-last} For all classical algebraic weights
  $\kappa = (( k_v, l_v)_{v |p})$ with $l_v = 2$ when $v \in I$
  and $k_v\equiv l_v\equiv 2\pmod{p-1}$ for all $v |p$, the
  operator $U^I$ is locally finite on
  $\mathrm{R}\Gamma(\mathfrak{X}^{G_1,I}_{K,\Kli}(p^\infty),
  \omega^\kappa(-D))$ and
  $\mathrm{R}\Gamma(\mathfrak{X}^{G_1,I}_{K,\Kli}(p^\infty),
  \Omega^\kappa(-D))$,  and the map \[e(U^I) \mathrm{R}\Gamma(
    \mathfrak{X}^{G_1, I}_{K, \Kli}(p), \omega^{\kappa}(-D))
    \rightarrow e(U^I) \mathrm{R}\Gamma( \mathfrak{X}^{G_1,
      I}_{K,\Kli}(p^\infty), \Omega^{\kappa}(-D))\] is a quasi-isomorphism. \end{lem}

\begin{proof} We consider the exact triangle   $$ \mathrm{R}\Gamma(\mathfrak{X}^{G_1,I}_{K,\Kli}(p^\infty), K\omega^\kappa(-D)) \rightarrow \mathrm{R}\Gamma(\mathfrak{X}^{G_1,I}_{K,\Kli}(p^\infty), \omega^\kappa(-D)) $$ $$ \rightarrow \mathrm{R}\Gamma(\mathfrak{X}^{G_1,I}_{K,\Kli}(p^\infty), \Omega^\kappa(-D)).$$
By Lemma ~\ref{lem:
  comparison of p adic classical weights U operators}, the operator
$U^I$ is topologically nilpotent on
$$ \mathrm{R}\Gamma(\mathfrak{X}^{G_1,I}_{K,\Kli}(p^\infty),
K\omega^\kappa(-D)),$$
so in particular it acts locally finitely with
$e(U^I) = 0$.

If we further
assume that $l_w\ge 2$ for all~$w$, $k_w-l_w \geq C$ for all $w | p$,
  and $l_w \geq C$ for all $w \in I^c$, then it follows from Lemma \ref{lem: classicity at
  level of sheaf} that $U^I$ is locally finite on
$\mathrm{R}\Gamma(\mathfrak{X}^{G_1,I}_{K,\Kli}(p^\infty),
\omega^\kappa(-D))$.  Therefore in this case, it follows from the
above exact triangle that~$U^I$ is locally finite on
$\mathrm{R}\Gamma(\mathfrak{X}^{G_1,I}_{K,\Kli}(p^\infty),
\Omega^\kappa(-D))$.  

Again using~\cite[Prop.\ 2.3.1]{pilloniHidacomplexes}, we deduce that $U^I$ is locally finite on
the complex
$\mathrm{R}\Gamma(\mathfrak{X}^{G_1,I}_{K,\Kli}(p^\infty),
\Omega^\kappa(-D))$ for any weight $\kappa$, and therefore (again
using the above exact triangle) it is also
locally finite on
$\mathrm{R}\Gamma(\mathfrak{X}^{G_1,I}_{K,\Kli}(p^\infty),
\omega^\kappa(-D))$ for any weight $\kappa$, as required. \end{proof}

 \begin{lem}\label{lem: this had the same label as another lemma} For
   all classical algebraic weights $\kappa = (( k_v, l_v)_{v |p})$ with $l_v
   = 2$ when $v \in I$, $k_v\equiv l_v\equiv 2\pmod{p-1}$, for all $v
   |p$, $l_w\ge 2$ for all~$w$, $k_w-l_w \geq C$ for all $w \in I$,
  and $l_w \geq C$ for all $w \in I^c$, the complex  $e(U^I) \mathrm{R}\Gamma( \mathfrak{X}^{I, G_1}_{K, \Kli}(p^\infty), \Omega^{\kappa}(-D))  $ is a perfect complex  of $\ZZ_p$-modules concentrated in degrees $[0, \# I ]$. 
\end{lem}
\begin{proof}  This follows from Lemma \ref{lem-classicity-last},
  Lemma  \ref{lem: classicity at level of sheaf} and Theorem
  \ref{thm-vanishing} (noting that~$T^I$ divides~$\widetilde{T}^I$, so
  that $e(\widetilde{T}^I)\mathrm{R}\Gamma(X_1^{G_1, I}, \omega^{\kappa}(-D))$ is
  a direct summand of $e(T^I)\mathrm{R}\Gamma(X_1^{G_1, I},
  \omega^{\kappa}(-D))$).\end{proof}

\begin{lem}\label{lem: p-adic complex is perfect}The operator  $U^I$ is locally finite on  $\mathrm{R}\Gamma( \mathfrak{X}^{I, G_1}_{K,\Kli}(p^\infty), \Omega^{\kappa_I}(-D)) $, and $e(U^I) \mathrm{R}\Gamma( \mathfrak{X}^{I, G_1}_{K,\Kli}(p^\infty), \Omega^{\kappa_I}(-D))  $ is a perfect complex  of $\Lambda_I$-modules concentrated in degree $[0, \# I ]$. 
\end{lem} 
\begin{proof} This follows from  Lemma~\ref{lem: this had the same label as another lemma} by
   Nakayama's lemma, in the form of Lemma~\ref{lem: finite cohomology
     implies perfect} and~\cite[Prop.\ 2.3.1]
  {pilloniHidacomplexes}. \end{proof}

\begin{proof}[Proof of Theorem~\ref{theorem-p-adic-complex}] Parts~(1)
 and ~(2) are Lemma~\ref{lem: p-adic
    complex is perfect}. Part~(3) is
  Lemma~\ref{lem-classicity-last}, together with Lemma~\ref{lem:
    control from Kli p infinity to Kli p}, which shows that the natural map 
    $$e (U^I) \RGamma(\mathfrak{X}^{I, G_1}_{K, \Kli}(p),
\omega^\kappa(-D)) \rightarrow e (U^I) \RGamma(\mathfrak{X}^{I, G_1}_{K, \Kli}(p^\infty),
\omega^\kappa(-D))$$
 is a quasi-isomorphism. Part~(4)  follows from
  Theorem~\ref{thm-vanishing}, together with Proposition~\ref{prop: T
    to U in medium weight} and Lemma~\ref{lem-classicity-last}.
\end{proof}

\section{Doubling}\label{sec:doubling}

In this section, 
we prove a doubling result (see Theorems~\ref{thm: the doubling map is injective}
and~\ref{thm: main doubling thm}) which is the key ingredient for proving local--global
compatibility  in~\S\ref{subsec: TW systems2}.
The general ideal of doubling  is that certain spaces of ordinary low weight modular forms
admit (at least) two degeneracy maps to spaces of ordinary modular forms of either higher
weight or higher level.  For example, the space of weight one elliptic modular
forms modulo~$p$ of level~$\Gamma_1(N)$, $p\nmid N$,  admits degeneracy maps~$f \mapsto \Ha \cdot f$ 
and~$f \mapsto f^p$ (where~$\Ha$ is the Hasse invariant) to spaces of forms of weight~$p$ and level~$\Gamma_1(N)$. (Alternatively, after dividing by~$\Ha$, these degeneracy maps can also be thought of
as maps from classical forms of level~$\Gamma_1(N)$  and weight one 
to ordinary $p$-adic modular forms
of level~$\Gamma_1(N)$ and weight one.)
If one can show that the direct sum of two copies of the original
space  embeds under the direct sum of these degeneracy maps, then,
following ideas going back to Gross~\cite{MR1074305} and
 isolated and expanded by~\cite{MR3247800} 
 (see also~\cite{CG} for further
 exploitation of these ideas), one can
make deductions about the local properties of the Galois
representations of interest. 

Let us explicate this in the example of weight one forms mentioned above (the following is implicit in the
first few lines of~\cite[p.499]{MR1074305}
and explicit in~\cite[Prop.~2.7]{edix:weights}).
If~$f$ is  a weight one elliptic modular cuspidal eigenform with
Nebentypus character~$\chi$ and $T_p$-eigenvalue $a_p$ satisfying $a^2_p\ne 4 \chi(p)$,
one can show that the associated Galois representation is unramified
at~$p$ in the following way. Since the polynomial $X^2-a_pX+\chi(p)$ has
distinct roots, one can show using the degeneracy maps above
(and having established doubling) that
there are two weight~$p$ ordinary forms congruent to~$f$ with
level~$\Gamma_1(N)$ and $T_p$-eigenvalues given by the roots~$\alpha$ and~$\beta$ of
$X^2-a_pX+\chi(p)$. Using the known properties of the corresponding Galois
representations, one shows that the restriction to~$p$ is an extension
of distinct unramified characters, and thus that the extension is split
(because the representations corresponding to the two weight~$p$ forms
are extensions in the opposite orders). 

The above argument for local--global compatibility at~$p$ works equally well in the ordinary symplectic case once
we have established a doubling theorem, and we will use this in
\S\ref{subsec: TW systems2} below. Before proceeding, 
we begin by recalling the doubling argument in more detail in the case
of~$\GL_2/\Q$. 

\subsection{The case of \texorpdfstring{$\GL_2/\Q$}{GL2Q}} For the moment, let~$X$ denote the
special fibre of a classical modular curve of level~$\Gamma_1(N)$
with~$N \ge 5$  with~$p\nmid N$, and let~$\omega$ denote the usual invertible
 line bundle on~$X$ (as in~\cite[\S2]{MR1074305}).
The doubling strategy of~\cite{CG,C2} may be reduced to ruling out the existence of simultaneous eigenforms~$f \in H^0(X,\omega)$ for the operators~$T_p$ and~$U_p$.
This is easily seen: indeed if~$f$ is a simultaneous eigenform for~$T_p$ and~$U_p$,
 then it is also an eigenform for~$V_p  = T_p - U_p$, which is
 immediately seen to be impossible by examining the action on
 $q$-expansions. This argument does not directly generalize to the symplectic case (even over~$\Q$), and instead,
 the paper~\cite{CGGSp4} employs a rather labyrinthian argument
 involving $q$-expansions to prove an analogous result for~$\GSp_4/\Q$. In this paper, we give a different argument which is based on analyzing the behavior of the $U_p$ operator at the non ordinary locus. This argument in this form appears to be new even for modular forms
 of weight one (although there are certainly some echos of this argument in papers such as~\cite{MR642341,MR0404145,Cais}),
 and so we present it first as a warm up for the general symplectic case.

 If~$f  \in H^0(X,\omega)$, we may think of~$U_p f$ as a section of~$H^0(X \setminus SS,\omega)$ for the finite set~$SS$ of supersingular points of~$X$.  We claim that there is a commutative diagram
\numequation\label{eqn: commutative diagram restriction to SS}
\xymatrix{
  H^0(X,\omega)\ar[r]^{\Ha  \cdot U_p}\ar[d]&H^0(X,\omega^p)\ar[d]\\
H^0(SS,\omega)\ar[r]&H^0(SS,\omega^p)\
}
\end{equation}
where the vertical maps are the natural restriction maps, and the
lower horizontal map is
an isomorphism. The existence of such a diagram can be proved in several
ways; for example it can be checked in the same way as the
corresponding statements for~$\GSp_4/F$ later in this section, by using the
Kodaira--Spencer isomorphism to describe the $U_p$ operator as a trace map on differentials.

Suppose that~$f$ is a ~$U_p$-eigenform in~$H^0(X,\omega)$
 with non-zero eigenvalue. Considering the commutative diagram~\eqref{eqn: commutative diagram
  restriction to SS}, we see that since~$\Ha \cdot U_pf$ maps
to zero in~$H^0(SS,\omega^p)$, the restriction of~$f$ to~$SS$
must vanish. Thus~$f=\Ha\cdot g$ for some~$g\in H^0(X,\omega^{2-p})$,
and this cohomology group vanishes if~$p>2$, so~$f=0$ in this case. If~$p = 2$, the only non-zero sections of~$H^0(X,\OL_X)$ are constants,
 and we deduce that~$f$ is a multiple of the Hasse invariant.

In the rest of this section we prove a generalization of
this  to the Hilbert--Siegel case. The analogue of
the commutative diagram~\eqref{eqn: commutative diagram restriction to SS}  in the Siegel case (with
$F=\Q$) is the following commutative diagram (where we write~$Y^{\geq 1}$ for the locus in the interior of the special fibre of the Shimura variety with Klingen level $H$ which is multiplicative, and we write $Y^{=1}$ for the divisor where the abelian variety is non ordinary.)

\begin{equation*}
\xymatrix@+2pc{
H^0(Y^{\geq 1},\omega^2)\ar[r]^{\Ha \cdot U_{\Iw(p),1}}\ar[d]&H^0(Y^{\geq 1},\omega^{p+1})\ar[d]\\
H^0(Y^{= 1},\omega^2)\ar[r]&H^0(Y^{=
  1},\omega^{p+1})
}
\end{equation*}
However, in contrast to the modular curve case, the map on the bottom line of this diagram is probably not
injective, so we cannot conclude as before. Instead, we construct a
larger commutative diagram

\begin{equation*}
\xymatrix@+2pc{
H^0(Y^{\geq 1},\omega^2)\ar[r]^{\Ha \cdot U_{\Iw(p),1}}\ar@/^2pc/[rr]^{U_{p,2}}\ar[d]&H^0(Y^{\geq 1},\omega^{p+1})\ar[d]&H^0(Y^{\geq 1},\omega^2)\ar[d]\\
H^0(Y^{= 1},\omega^2)\ar[r]&H^0(Y^{=
  1},\omega^{p+1})\ar[r]&H^0(Y^{= 1},\omega^2).
}
\end{equation*}
If we assume that~$f\in H^0(Y^{\ge 1},\omega^2)$ is also a~$U_{p,2}$-eigenform with nonzero
eigenvalue (which suffices for our purposes), we can use this diagram
to make a similar argument to the above, considering the composite
morphisms from the top left to the lower right hand corner. It may help the reader to
note that there is an analogous diagram for~$\GL_2$:
\begin{equation*}
\xymatrix{
  H^0(X,\omega)\ar[r]^{\Ha \cdot U_p}\ar@/^2pc/[rr]^{\langle p\rangle}\ar[d]&H^0(X,\omega^p)\ar[d]&H^0(X,\omega)\ar[d]\\
H^0(SS,\omega)\ar[r]&H^0(SS,\omega^p)\ar[r]&H^0(SS,\omega)
}
\end{equation*}(again, the existence of this
diagram can be checked in the same way as our calculations
below). We see that if~$U_p f$ has no poles, then the image
of~$f$ in the bottom right hand copy of~$H^0(SS,\omega)$ vanishes;
since the diamond operator~$\langle p\rangle$ is an isomorphism, it
follows that the restriction of~$f$ to~$SS$ vanishes, and we conclude
as before. 

There is an additional complication in the Hilbert--Siegel case, which
is that rather than considering the entire Shimura variety, we are
only working on an open
subspace~$X_{K_p(I) K^p,1}^{ \geq_{v\in I} 1, \geq_{v \in I^c}
  2}$. This means that the vanishing of the space of (partial)
negative weight modular forms is not
obvious. We sketch a proof for this vanishing in~\S\ref{section:vanishinginnegativeweight} below, using  Fourier--Jacobi expansions (which ultimately reduces
to the vanishing of spaces of Hilbert modular forms of partial
negative weight), but we do not rely on this result. Instead, we give
a complete proof of a slightly weaker result which is nonetheless
sufficient for our purposes; this argument does not use the boundary, but rather
considers the behavior of another Hecke operator~$U_{\Kli(w),1}-U_{\Iw(w),1}$ (called $Z_w$ below) along the $w$-non-ordinary locus.

\subsection{Conventions}Throughout this section, we fix a set $I\subset S_p$ and a prime $w\in
I$. Recall that, as in~\S\ref{subsubsec: Igusa towers}, for each subset $I\subset S_p$ we set \[K_p(I) = \prod_{v
    \in I} \Kli(v) \prod_{v \in I^c} \Iw(v),\] and we write
$X_1^I:=X_{K_p(I) K^p,1}^{ \geq_{v\in I} 1, \geq_{v \in I^c}
  2}$. We will use the following simplified notation:
\begin{itemize}
\item We write $X_1^I$ for the space $X_{K_p(I) K^p,1}^{\geq_{v\in I}1,\geq_{v\in I^c}2}$.
\item We write $X_1^{I,=_w2}$ for the open subspace where $A[w^\infty]$ is ordinary. 
\item We write $X_1^{I,=_w1}$ for the (reduced) complement of this open subspace, which is a divisor in $X_1^I$.
\end{itemize}We also write $Y_1^I$, $Y_1^{I,=_w2}$ and $Y_1^{I,=_w1}$ for their
interiors. We use the analogous notation $\fX^I$, $\fY^I$ etc. for the corresponding formal schemes. We denote $\det\omega_\G$ by $\omega$ and $\det\omega_{\G_w}$ by $\omega_w$.
 We will  finally denote the partial Hasse invariant $\Ha(\G_w)\in H^0(X_1^I,\omega_w^{p-1})$ by $\Ha_w$.

\subsection{The operator \texorpdfstring{$U_{\Kli(w),1}$}{UKli(1)}}\label{sect-UKLI}

We now define a Hecke operator~$U_{\Kli(w),1}$ (see also~\S~\ref{subsec: U Kli 1
  first time around}). We again consider the $p$-adic completion of  the correspondence considered in \S
\ref{sect-def-Siegel-Hecke} 
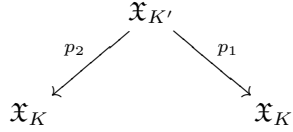
\begin{eqnarray*}
                 \xymatrix{ &\mathfrak{X}_{K'}  \ar[rd]^{p_1}\ar[ld]_{p_2}&  \\
                 \mathfrak{X}_K & &
                                                  \mathfrak{X}_K}
  \end{eqnarray*}where we recall that $K = K^pK_p$ with $K_p = \prod_v \GSp_{4}(\ocal_{F_v})$  and $K' = K^pK'_p$ with $K_p' = \Si(w) \times \prod_{v \neq w} \GSp_{4}(\ocal_{F_v})$.

We can form  the fibre product $\mathfrak{X}_{K'} \times_{p_1, \mathfrak{X}_K}\mathfrak{X}^{I}$.  As $\mathfrak{X}^I\to\mathfrak{X}_K$ is \'etale by Proposition \ref{prop: deep Klingen is affine etale}, this inherits the properties of $\mathfrak{X}_{K'}$ deduced from the theories of local models and toroidal compactifications.  In particular it is flat over $\Z_p$, normal, Cohen--Macaulay, and the ordinary locus is dense in the special fibre, see \S\ref{subsubsec: local models}, Theorem \ref{thm: main KWL compactification thm}, and \S\ref{subsec: Hasse invariants and stratifications}.
We denote by $\mathfrak{C}_{\Kli(w),1}$ the open and closed formal subscheme
of this fibre product where the kernel of the canonical isogeny
$p_1^{*} \cG \rightarrow p_2^{*} \cG$ has trivial intersection with
the multiplicative group~$p_1^*H_w$.

We obtain a
correspondence \begin{eqnarray*}
                 \xymatrix{ &\mathfrak{C}_{\Kli(w),1}   \ar[rd]^{v_1}\ar[ld]_{v_2}&  \\
                 \mathfrak{X}^I & &
                                                  \mathfrak{X}^I}
  \end{eqnarray*}
where  $v_1: \mathfrak{C}_{\Kli(w),1} \rightarrow \mathfrak{X}^I$ is induced by the projection $\mathfrak{X}_{K'} \times_{p_1, \mathfrak{X}_K}\mathfrak{X}^{I}\to\mathfrak{X}^I$ and $v_2$ is defined as follows: the projection $\mathfrak{X}_{K'} \times_{p_1, \mathfrak{X}_K}\mathfrak{X}^{I}\to\mathfrak{X}_{K'}$ composed with $p_2:\mathfrak{X}_{K'}\to\mathfrak{X}_K$ induces a map $\mathfrak{C}_{\Kli(w),1}\to\mathfrak{X}_K$ which we would like to lift to a map $v_2:\mathfrak{C}_{\Kli(w),1}\to\mathfrak{X}^I$.  In other words, given a point of $\mathfrak{C}_{\Kli(w),1}$, we need to give multiplicative subgroups of order $p$, $H'_{w'}\subseteq p_2^*\cG_{w'}$ for all $w'\in S_p$.  But for all $w'\in S_p$, the kernel of the isogeny $p_1^*\cG_{w'}\to p_2^*\cG_{w'}$ has trivial intersection with $p_1^*H_{w'}$, (for $w'\not=w$ it is an isomorphism, and for $w'=w$ this was assumed in the definition of $\mathfrak{C}_{\Kli(w),1}$) and we take $H'_{w'}$ to be the image of $p_1^*H_{w'}$ under this isogeny.

\begin{lem}\label{lem-existencetracemap}  We have $R (v_1)_* \oscr_{\mathfrak{C}_{\Kli(w),1}} =( v_1)_* \oscr_{\mathfrak{C}_{\Kli(w),1}}$ and there is a trace map $R (v_1)_* \oscr_{\mathfrak{C}_{\Kli(w),1}} \rightarrow \oscr_{\mathfrak{X}^I}$. 
\end{lem}

\begin{proof} We have $\mathfrak{X}_{K} = \mathfrak{X}_{K, \Sigma}$ for a smooth polyhedral cone decomposition $\Sigma$ and we have $\mathfrak{X}_{K'} = \mathfrak{X}_{K', \Sigma'}$. We can now assume that $\Sigma' = \Sigma$ because we have $R \pi_* \oscr_{\mathfrak{X}_{K', \Sigma'}} = \oscr_{\mathfrak{X}_{K', \Sigma}}$ for $\pi:\fX_{K',\Sigma'}\to\fX_{K',\Sigma}$ the projection (we note that the cone decomposition at level $K'$ may not be smooth but we will not need this).
Since~$\Sigma=\Sigma'$, the map $v_1$ is quasi-finite, and (since it is proper) is
therefore finite. Hence we have~$R (v_1)_* \oscr_{\mathfrak{C}_{\Kli(w),1}}
=( v_1)_* \oscr_{\mathfrak{C}_{\Kli(w),1}}$.
Moreover, as $\mathfrak{C}_{\Kli(w),1}$ is Cohen--Macaulay and $\mathfrak{X}^I$ is regular, we deduce that the map $v_1$ is also flat, and so it has an associated trace map.
\end{proof}

We let $\mathfrak{C}_{\Kli(w),1}^{=_w2}$ be the open formal subscheme where $v_1^*\cG_w$ (or equivalently $v_2^*\cG_w$) is ordinary.  It restricts to a correspondence over $\mathfrak{X}^{I,=_w2}$.  Over
$\mathfrak{C}_{\Kli(w),1}^{=_w2}$, the multiplicative rank of $\ker(v_1^*\cG_w\to v_2^*\cG_w)$ is either 0 or 1 (it cannot be 2 because $H_w$ has trivial intersection with $\ker(v_1^*\cG\to v_2^*\cG)$), and hence we have a decomposition $$\mathfrak{C}_{\Kli(w),1}^{=_w2}  = \mathfrak{C}_{\Kli(w),1}^{=_w2,et}  \coprod \mathfrak{C}_{\Kli(w),1}^{=_w2, \met}$$
where $\mathfrak{C}_{\Kli(w),1}^{=_w2,et}$ is the locus where the isogeny $v_1^*\cG\to v_2^*\cG$ is \'etale, while $\mathfrak{C}_{\Kli(w),1}^{=_w2, \met}$ is the locus where $\ker(v_1^*\cG\to v_2^*\cG)$ has multiplicative rank 1.

For any weight $\kappa = (k_v, l_v)$, we have a map 
$$ v_2^\star \omega^\kappa \rightarrow v_1^\star \omega^\kappa [1/p]$$
 induced from the universal isogeny (we note that we are not assuming $l_w\geq 0$.)  Tensoring this map with the trace map of Lemma \ref{lem-existencetracemap}, we obtain a map of sheaves over $\mathfrak{X}^I$: 
  $$ \Theta_\kappa:  (R v_1)_\star v_2^\star \omega^\kappa \rightarrow  \omega^\kappa [1/p].$$ We now define $U_{\Kli(w),1} =  p^{-l_w -1}  \Theta_\kappa $ if $l_w \leq 2$ and $U_{\Kli(w),1} =  p^{-3}  \Theta_\kappa $   if $l_w \geq 2$. 
 
 \begin{lem}\label{lem-onemorenormalization} We have $U_{\Kli(w),1}: (R v_1)_\star v_2^\star \omega^\kappa \rightarrow  \omega^\kappa $. 
 \end{lem}
 
 \begin{proof} We follow the same strategy as the proof of Lemma \ref{lem-normalization1}.  Both the source and target of this map are locally free sheaves over the smooth formal scheme $\mathfrak{X}^I$. To prove that the map is indeed $p$-integral, it is enough to prove it over the ordinary locus, and we check it separately on each type of component.
 
On the component of the map corresponding to $\mathfrak{C}_{\Kli(w),1}^{=_w2,et}$, the isogeny is \'etale over the ordinary locus and therefore the map $ v_2^\star \omega^\kappa \rightarrow v_1^\star \omega^\kappa [1/p]$ is actually  an isomorphism $ v_2^\star \omega^\kappa \rightarrow v_1^\star \omega^\kappa $, while the trace map is divisible by $p^3$ (see the proof of Lemma \ref{lem-normalization1}).
 
On the component of the map corresponding to $\mathfrak{C}_{\Kli(w),1}^{=_w2,\met}$, the isogeny has multiplicative rank one over the ordinary locus and therefore the map $$ v_2^\star \omega^\kappa \rightarrow v_1^\star \omega^\kappa [1/p]$$ is actually  a map  $$ v_2^\star \omega^\kappa \rightarrow p^{l_w}v_1^\star \omega^\kappa,$$ while the trace map is divisible by $p$ (again see the proof of Lemma \ref{lem-normalization1}).

Therefore, on the \'etale component, the map is divisible by $p^{3}$ and on the multiplicative-\'{e}tale component it is divisible by $p^{l_w +1}$. 
\end{proof}

\subsection{The operators \texorpdfstring{$U_{\Iw(w),1}$}{UIw(1)} and \texorpdfstring{$Z_w$}{Zw}}\label{subsec: UIw1 and Zw}
Now we consider some Hecke operators on $\mathfrak{X}^{I,=_w2}$.  The restriction of $\Theta_\kappa$ to $\mathfrak{X}^{I,=_w2}$ decomposes as a sum $\Theta_\kappa=\Theta_\kappa^{et}+\Theta_\kappa^{\met}$, according to the decomposition of $\mathfrak{C}_{\Kli(w),1}^{=_w2}$.  We define normalized cohomological correspondences
$U_{\Iw(w),1} =  p^{-3}\Theta_\kappa^{et} $ and $Z_w =  p^{-l_w-1} \Theta_\kappa^{\met} $.

\begin{lem} The cohomological correspondences $U_{\Iw(w),1}$ and $Z_w$ are $p$-integral. 
 \end{lem}
\begin{proof} This follows from the proof of Lemma \ref{lem-onemorenormalization}. 
\end{proof} 

We have the following identities of Hecke operators over the ordinary locus at $w$: 

\begin{enumerate}
\item $U_{\Kli(w), 1} =  U_{\Iw(w),1} + p^{ l_w-2} Z_w$ if $l_w \geq 2$,
\item  $U_{\Kli(w), 1} =  p^{2-l_w}U_{\Iw(w),1} + Z_w$ if $l_w \leq 2$.
\end{enumerate}

It follows in particular that: 

\begin{enumerate}
\item $U_{\Kli(w), 1} =  U_{\Iw(w),1} ~\mod p$ if $l_w > 2$,
\item  $U_{\Kli(w), 1} = Z_w~\mod p$ if $l_w < 2$.
\end{enumerate}

Another important property is the following: 

\begin{prop}\label{prop-commuteagain} For any weight $\kappa$, we have the following identities of cohomological correspondences over $X^{I, =_w 2}_1$: 

\begin{enumerate}
\item  $Z_w \Ha_w = \Ha_w Z_w$,
\item  $U_{\Iw(w), 1}  \Ha_w = \Ha_w U_{\Iw(w), 1}.$
\end{enumerate}

\end{prop}

\begin{proof}  The correspondence $Z_w$ is the tensor product of the fundamental class (deduced from the trace map normalized by $p^{-1}$) and a map $v_2^\star \omega^{\kappa} \rightarrow v_1^\star \omega^\kappa$ which is obtained by normalizing the natural map by a factor $p^{-l_w}$. It suffices to check that for $\omega^\kappa = \omega_w^{p-1}$ this normalized map matches the Hasse invariants $v_2^\star \mathrm{Ha}( \cG_w)$ and $v_1^\star \mathrm{Ha}(\cG_w)$. This is the content of \cite[Lem.\ 6.2.4.1]{pilloniHidacomplexes}. The case of $U_{\Iw(w),1}$ is clear because the universal isogeny is \'etale.
\end{proof}

Finally, we will need the following property:

\begin{prop}\label{prop-noZw} If $l_w \leq 0$,  the cohomological correspondence  over $X_1^{I}$:
$$ U_{\Kli(w), 1}: (v_1)_\star v_2^\star \omega^\kappa \rightarrow \omega^\kappa$$ factors through
$$  U_{\Kli(w), 1}: (v_1)_\star v_2^\star \omega^\kappa \rightarrow \omega^\kappa (- X_1^{I, =_w1}).$$
\end{prop}

Before giving the proof we need some preparations. Let $BT/\F_p$ be the smooth algebraic stack of quasi-polarized  1-truncated  Barsotti--Tate groups of height $2$ and dimension $1$ over $\Spec\F_p$. 
Let $Y/\F_p$ be a modular curve of level prime to~$p$. The map $Y
\rightarrow BT$ is a presentation of $BT$ (that is, it is a smooth surjection). We denote by $E$ the universal object on $BT$. 
We have a Cartier divisor $\omega_E^{1-p} \stackrel{\mathrm{Ha}(E)} \rightarrow \oscr_{BT}$ whose support is the non-ordinary locus of $BT$. 
Let $\pi: BT_{\Iw} \rightarrow BT$ be the representable finite flat map which parametrizes a subgroup $H \subset E$ of order $p$. Let $Y_0(p)$ be a modular curve of Iwahori level at $p$. The map $Y_0(p) \rightarrow BT_{\Iw}$ is a  presentation of $BT_{\Iw}$. 
Over $BT_{\Iw}$ we have a universal morphism $ g: E/H \rightarrow E$ with kernel $H^D$ (using the polarization to identify $E$ and $E^D$, $E/H$ and $(E/H)^D$).  By differentiating, we get a map  of line bundles $d g: \omega_{E} \otimes \omega_{E/H}^{-1} \rightarrow \oscr_{BT_\Iw}$. 

\begin{lem}\label{lemma-stacky} We have a canonical  factorization $ (d g)^{\otimes 2}: (\omega_{E} \otimes \omega_{E/H}^{-1})^{\otimes 2} \rightarrow \pi^\star \omega_E^{1-p} \rightarrow \oscr_{BT_\Iw}$.
\end{lem}

\begin{proof} It suffices to prove the claim over any  presentation of $BT_{\Iw}$. We therefore reduce to proving the statement over the modular curve $Y_0(p)$. The vanishing locus of $\pi^\star (\mathrm{Ha}(E))$ is a product of  Artinian local rings of length $p+1$ (the degree of $\pi$) indexed by the supersingular points. The vanishing locus of $d g $ is the entire  irreducible component of $Y_0(p)$ which is degree $p$ over $Y$ via $\pi$ (this is the component where $H$ is generically \'etale). Therefore, for any supersingular point $x \in Y$, the image of $dg$ in $\oscr_{Y_0(p)} \otimes_{\oscr_{Y}} k(x)$ defines a closed subscheme of length $p$, and hence the ideal generated by the image of $dg$ in $\oscr_{Y_0(p)} \otimes_{\oscr_{Y}} k(x)$ is both nilpotent and length 1. It follows that $(dg)^2$ maps to zero in $\oscr_{Y_0(p)} \otimes_{\oscr_{Y}} k(x)$. 
\end{proof}

\begin{proof}[Proof of Proposition \ref{prop-noZw}]  Let $\kappa$ be a weight with $l_w \leq 2$. Let  $\kappa'$ be another weight with $(k_v, l_v) = (k'_v, l'_v)$ for $v\not=w$, $k_w -l_w = k'_w-l'_w$ and $l'_w =2$. 
Let us denote by $U_{\Kli(w), 1}(2): (v_1)_\star v_1^\star \omega^{\kappa'} \rightarrow \omega^{\kappa'}$  the cohomological correspondence in weight $\kappa'$. Let $U_{\Kli(w),1}: (v_1)_\star v_1^\star \omega^{\kappa} \rightarrow \omega^{\kappa}$ be the cohomological correspondence in weight $\kappa$.  The proof of Lemma \ref{lem-onemorenormalization} shows that the map $v_2^*\det\omega_w^{l_w-2}\to v_1^*\det\omega_w^{l_w-2}[1/p]$ induces a regular map $p^{l_w-2}v_2^*\det\omega_w^{l_w-2}\to v_1^*\det\omega_w^{l_w-2}$, and that moreover $U_{\Kli(w),1}$ is obtained from $U_{\Kli(w),1}(2)$ by twisting by this map.  It thus suffices to show that on the special fibre, $p^{l_w-2}v_2^*\det\omega_w^{l_w-2}\to v_1^*\det\omega_w^{l_w-2}$ factors through $v_1^*\det\omega_w^{l_w-2}(-X_1^{I,=1})$ when $l_w\leq 0$.

This statement is local in a neighbourhood of $X^{I, =1}_1$ and we can
therefore replace $\mathfrak{X}^{I}$ by its completion along this
closed subscheme. We may also work on the interior of the moduli
space, as the interior of the divisor $X^{I,=1}_1$ is dense.  Therefore, we may suppose that $\cG_w$ comes equipped with a
multiplicative sub-Barsotti--Tate subgroup~$\cG_w^m$ of rank $1$,  and
we denote by $\cG_w^{oo} =(\cG_w^m)^\bot/ \cG_w^m$, which is a Barsotti--Tate group scheme of height $2$ and dimension $1$. 
The isogeny $v_1^\star \cG_w \rightarrow v_2^\star \cG_w$ induces an isomorphism 
$v_1^\star  \cG_w^m \rightarrow v_2^\star \cG_w^m$ and a degree $p$ map 
$v_1^\star \cG_w^{oo} \rightarrow v_2^\star \cG_w^{oo}.$

The normalized map $p^{-1} v_2^\star \omega_w^{-1} \rightarrow v_1^\star \omega_w^{-1}$ is the tensor product of the isomorphism 
$ v_2^\star \omega_{\cG_w^m}^{-1} \rightarrow v_1^\star \omega_{\cG_w^m}^{-1}$
and the map: $p^{-1} v_2^\star \omega_{\cG_w^{oo}}^{-1} \rightarrow v_1^\star \omega_{\cG^{oo}_w}^{-1}$ which is the transpose of the map 
$ v_1^\star \omega_{\cG_w^{oo}} \rightarrow v_2^\star \omega_{\cG^{oo}_w}$ obtained by differentiating the dual isogeny: 
$v_2^\star (\cG_w^{oo}) \rightarrow v_1^\star \cG_w^{oo}.$ The result
follows from Lemma~ \ref{lemma-stacky}. 
\end{proof}

\begin{cor}\label{cor: the S conjecture} Let $\kappa $ be a weight with $l_w \leq 0$. 
  Let  $f\in H^0(X_1^{I},\omega^{\kappa})$  be such that $Z_w f = \beta_w f$ for some $\beta_w \neq 0$. Then $f=0$.
 \end{cor}
  
  \begin{proof}  Since~$l_w\le 0$, we have $Z_w =
    U_{\Kli(w),1}$ on $H^0(X_1^{I},\omega^{\kappa})$. Assume that $f \neq 0$, and let  $n$ be the order
    of vanishing of $f$ along $X_1^{=1}$.   By considering
    $\mathrm{Ha}(\cG_w)^{-n} f$ and using Proposition~
    \ref{prop-commuteagain}  we can suppose that $n=0$.  This contradicts Proposition \ref{prop-noZw}.
  \end{proof}

\subsection{Preliminaries on Kodaira--Spencer}\label{section:KS}

In this section, we recall the Kodaira--Spencer map and its
compatibility with certain functorialities. A convenient reference for
what we need is~\cite{MR3186092}.

Let $S$ be a $\mathbb{Z}_{(p)}$-scheme and let $X$ be a smooth
$S$-scheme of relative dimension $3[F:\Q]$.  Suppose that we have a tuple $(A,\iota,\lambda)$ with
\begin{itemize}
\item $A/X$ an abelian scheme of dimension $2[F:\Q]$.
\item $\iota:\O_F\to\End(A)\otimes\mathbb{Z}_{(p)}$ making $\Lie(A)$ into a locally free $\O_F\otimes_{\mathbb{Z}}\O_X$-module of rank 2.
\item $\lambda:A\to A^t$ a prime to $p$, $\O_F$-linear quasi-polarization such that $\lambda[p^\infty]:A[p^\infty]\to A^t[p^\infty]$ is an isomorphism.
\end{itemize}

Then we have the first de Rham cohomology of $A/X$ together with its Hodge filtration
\begin{equation*}
0\to \omega_A\to H^1_{dR}(A/X)\to\omega_{A^t}^\vee\to 0
\end{equation*}
as well as the Gauss--Manin connection
\begin{equation*}
H^1_{dR}(A/X)\to H^1_{dR}(A/X)\otimes\Omega^1_{X/S}.
\end{equation*}
Passing to subquotients for the Hodge filtration we obtain the Kodaira--Spencer map for $A$
\begin{equation*}
\omega_A\to \omega_{A^t}^\vee\otimes\Omega^1_{X/S}.
\end{equation*}
The polarization $\lambda$ induces an isomorphism $\lambda^*:\omega_{A^t}\to\omega_A$.  Using this we may obtain a Kodaira--Spencer map for $(A,\lambda)$
\begin{equation*}
\omega_A\otimes\omega_A\to\Omega^1_{X/S}.
\end{equation*}
Then one checks (see~\cite[Prop.\ 6.2.5.18]{MR3186092}) that this map factors through the quotient $\Sym^2_{\O_X\otimes\O_F}\omega_A$ of $\omega_A\otimes\omega_A$, so that we obtain a map
\begin{equation*}
\Sym^2_{\O_X\otimes\O_F}\omega_A\to\Omega^1_{X/S}.
\end{equation*}
As usual, if $Y/S$ is a smooth scheme of relative dimension $d$, we write
$K_{Y/S}$ for the relative canonical bundle
$\wedge^d\Omega^1_{Y/S}$. We will be especially interested in the induced map on top exterior powers
\begin{equation*}
\wedge^{3[F:\Q]}(\Sym_{\O_X\otimes\O_F}^2\omega_A)=\det(\omega_A)^3\to
K_{X/S}.\end{equation*}

\begin{prop}\label{prop: KS isogeny}
Suppose that $(A,\iota,\lambda)$ and $(A',\iota',\lambda')$ are tuples as above and that we have a prime to $p$ quasi-isogeny $\phi:A\to A'$ satisfying $\phi\iota=\iota'\phi$ and $\phi^t\lambda'\phi=x\lambda$ for some $x\in \O_F\otimes\mathbb{Z}_{(p)}$.  Then we have a commutative diagram
\begin{equation*}
\xymatrix{
\det(\omega_{A'})^3\ar[r]\ar[d]^{\phi^*}& K_{X/S}\ar[d]^{\cdot N_{F/\Q}(x)^3}\\
\det(\omega_A)^3\ar[r]& K_{X/S}}
\end{equation*}
\end{prop}
\begin{proof}
  It follows from the definitions that under the Kodaira--Spencer
  maps, $\phi^*:\Sym^2_{\O_X\otimes\O_F}\omega_{A'}\to
  \Sym^2_{\O_X\otimes\O_F}\omega_A$ induces the endomorphism
  of~$\Omega^1_{X/S}$ given by multiplication by~$x$. The result
  follows on passing to top exterior powers.
\end{proof}

\begin{prop}\label{prop: KS traces}
Let $f:X\to Y$ be a finite flat map of smooth $S$-schemes of relative dimension $3[F:\Q]$ and let $(A,\iota,\lambda)/Y$ be a tuple as above.  Then the Kodaira--Spencer map is compatible with base change in the sense that there is a commutative diagram
\begin{equation*}
\xymatrix{
f^*\det(\omega_A)^3\ar[r]\ar[d]_{\wr}& f^*K_{Y/S}\ar[d]\\
\det(\omega_{A_X})^3\ar[r]& K_{X/S}
}
\end{equation*}
where the horizontal maps are the Kodaira--Spencer maps for $A$ and
$A_X$, the right vertical map is pullback on differentials, and
the left vertical map is the natural isomorphism.

Moreover it is compatible with traces in the sense that there is a commutative diagram
\begin{equation*}
\xymatrix{
f_*\det(\omega_{A_X})^3\ar[r]\ar[d]&f_*K_{X/S}\ar[d]\\
\det(\omega_A)^3\ar[r]&K_{Y/S}
}
\end{equation*}
where again the horizontal arrows are the Kodaira--Spencer maps for $A_X$ and~$A$ while the vertical map on the left comes from the \emph{(}unnormalized\emph{)} trace map on functions $f_*\O_X\to\O_Y$ and the isomorphism $\omega_{A_X}\simeq f^*\omega_A$, and the right vertical map is the trace map on dualizing sheaves.
\end{prop}
\begin{proof}The commutativity of the first diagram follows from the
  compatibility of the formation of de Rham cohomology with flat base
  change, and the compatibility of the Gauss--Manin connection with
  flat base change (which in turn follows from the compatibility of
  the Hodge to de Rham spectral sequence with flat base change). 

To see that the second diagram commutes, it is  by
adjunction equivalent to show that the lower square in the following
diagram commutes. 
\[\xymatrix{
f^*\det(\omega_A)^3\ar[r]\ar[d]_{\wr}& f^*K_{Y/S}\ar[d]\\
\det(\omega_{A_X})^3\ar[r]\ar[d]& K_{X/S}\ar[d]^{\wr}\\
f^!\det(\omega_A)^3\ar[r]& f^!K_{Y/S}
}\]
Since we have already seen that the upper square commutes, and since
the indicated vertical arrows are isomorphisms, the commutativity of
the lower square is equivalent to the commutativity of the outer
square. This commutativity follows from unwinding the definitions;
indeed, this outer square is the natural one obtained from the
Kodaira--Spencer morphisms and the natural transformation from~$f^*$ to~$f^!$ (which is given by the trace of the morphism $f$).
\end{proof}

Finally, we recall the Kodaira--Spencer isomorphism for our Shimura varieties.

\begin{prop}\label{prop: KS on divisor}
The Kodaira--Spencer map
\begin{equation*}
\omega^3\to K_{\fY^I/\Z_p}
\end{equation*}
is an isomorphism.
\end{prop}
\begin{proof}
This follows from the usual Kodaira--Spencer isomorphism~\cite[Thm.\ 6.4.1.1]{MR3186092}
 and the compatibility with \'{e}tale base change
proved in Proposition~\ref{prop: KS traces} (noting that the formation
of the canonical sheaf is compatible with \'etale base change).
\end{proof}

\subsection{The Hecke operator \texorpdfstring{$U_{\Iw(w),1}$}{UIw(1)} and traces for partial Frobenius}\label{section:U1equalsF}

We recall the construction of the Hecke operator $U_{\Iw(w),1}$.  We
have a correspondence (see~\S\ref{subsec: UIw1 and Zw}, where this correspondence was denoted $\fC^{=_w2,et}_{\Kli(w),1}$ but we adopt here a simplified notation $\fC^I_w$)\begin{equation*}
\xymatrix{
&\fC^I_w\ar[dl]_{p_2}\ar[dr]^{p_1}&\\
\fY^{I,=_w2}&&\fY^{I,=_w2}
}
\end{equation*}
where 
 $\fC^I_w$ parameterizes a point $(A,\iota,\lambda,\{H_v\}_{v\in
   S_p},\eta)$ of $\fY^{I,=_w2}$ along with an \'{e}tale maximal
 isotropic subgroup $L_w\subset A[w]$.  The map $p_1$ simply forgets
 $L_w$.  To describe $p_2$, consider the \'{e}tale isogeny $\pi:A\to
 A/L_w$.  Then $p_2$ sends $(A,\iota,\lambda,\{H_v\}_{v\in S_p})$ to
 $A/L_w$ with the induced action of $\O_F\otimes\mathbb{Z}_{(p)}$, the
 prime to $p$ quasi-polarization obtained by descending $x_w\lambda$,
 and the level structures $\pi(H_v)$ and $\pi(\eta)$. Since the subgroup $A[x_w]/L_w$ of $A/L_w$ is the canonical multiplicative subgroup of $A[x_w]$, we see that $p_2$ is an isomorphism.  

For any weight $\kappa$ for $G$, pullback by the universal \'{e}tale isogeny over $\fC_w^I$ induces an isomorphism of sheaves $p_2^*\omega^\kappa\to p_1^*\omega^\kappa$, and the Hecke operator $U_{\Iw(w),1}$ is obtained from the composition of maps of sheaves over $\fY^{I,=_w 2}$
\begin{equation*}
p_{1,*}p_2^*\omega^\kappa\to p_{1,*}p_1^*\omega^\kappa\overset{\text{``$\frac{1}{p^3}\Tr_{p_1}$\kern-0.2em{''}}}{\to} \omega^\kappa.
\end{equation*}
Now we turn to the Kodaira--Spencer isomorphism $\omega^3\simeq K_{\fY^{I,=_w2}}$.
\begin{prop}\label{prop: KS that we can reduce mod p}
There is a commutative diagram of sheaves on $\fY^{I,=_w2}$
\begin{equation*}
\xymatrix{
p_{1,*}p_2^*\omega^3\ar[r]\ar[d]&\omega^3\ar[d]\\
p_{1,*}p_2^*K_{\fY^{I,=_w2}/\Z_p}\ar[r]^{\frac{N_{F/\Q}(x_w)^3}{p^3}\text{tr}}&K_{\fY^{I,=_w2}/\Z_p}
}
\end{equation*}
where the vertical arrows are the Kodaira--Spencer isomorphism, and the
top horizontal arrow is $U_{\Iw(w),1}$. The bottom horizontal arrow is
defined as follows: since~$p_2$ is an isomorphism, we may identify
$p_2^*K_{\fY^{I,=_w2}/\Z_p}$ with $K_{\fC_w^I/\Z_p}$, and the morphism then comes from
 the trace map for $p_1$ on dualizing sheaves, multiplied by a factor of $\frac{N_{F/\Q}(x_w)^3}{p^3}\in\mathbb{Z}_{(p)}^\times$.
\end{prop}
\begin{proof}
This follows from Propositions \ref{prop: KS isogeny} and \ref{prop: KS traces}.
\end{proof}

We note that although we are primarily interested in using Proposition~\ref{prop: KS that we can reduce mod p} on the special fibre, we cannot apply Propositions \ref{prop: KS isogeny} and \ref{prop: KS traces} directly on the special fibre because some of the maps in the commutative square reduce to $0$ modulo $p$. 

We may also describe $U_{\Iw(w),1}$ in weights other than parallel weight 3 using traces on differentials.  For any weight $\kappa=(k_v,l_v)_{v\in S_p}$ we let $\kappa-3=(k_v-3,l_v-3)$.  Then tensoring the Kodaira--Spencer isomorphism with $\omega^{\kappa-3}$ we have an isomorphism $\omega^\kappa\simeq K_{\fY^{I,=_w2}/\Z_p}\otimes\omega^{\kappa-3}$.  Then we have a commutative diagram of sheaves on $\fY^{I,=_w2}$
\begin{equation*}
\xymatrix@C=0.5cm{
p_{1,*}p_2^*\omega^\kappa\ar[rr]\ar[d]&&\omega^\kappa\ar[d]\\
p_{1,*}p_2^*(K_{\fY^{I,=_w2}/\Z_p}\otimes\omega^{\kappa-3})\ar[r]&(p_{1,*}p_2^*K_{\fY^{I,=_w2}/\Z_p})\otimes\omega^{\kappa-3}\ar[r]&K_{\fY^{I,=_w2}/\Z_p}\otimes\omega^{\kappa-3}
}
\end{equation*}
where on the bottom row, the first map is an isomorphism coming from the projection formula and the isomorphism $p_2^*\omega^{\kappa-3}\simeq p_1^*\omega^{\kappa-3}$ and the second map is the tensor product of the map of Proposition \ref{prop: KS that we can reduce mod p} and the identity.

We would now like to understand the behavior of $U_{\Iw(w),1}$ beyond the $w$-ordinary locus on the special fibre.  In order to do this we make the following definition.

\begin{defn}
We define a ``partial Frobenius'' map 
\begin{equation*}
F_w:Y_1^I\to Y_1^I
\end{equation*}
as follows: given a point $(A,\iota,\lambda,\{H_v\}_{v\in S_p},\eta)$ of $Y_1^I$, we may consider the maximal isotropic subgroup $L_w\subset A[w]$ defined by
\begin{equation*}
L_w=\ker(F:A[w^\infty]\to A[w^\infty]^{(p)})
\end{equation*}
and form the subgroup of degree $p^{4[F:\Q]-2}$
\begin{equation*}
\tilde{L}_w=L_w\times\prod_{v\not=w}A[v]\subset A
\end{equation*}
and the isogeny $\pi:A\to \tilde{A}=A/\tilde{L}_w$.  $\tilde{L}_w$ is isotropic for the polarization $\frac{p^2}{x_w}\lambda$ which thus descends to a principal polarization $\tilde{\lambda}$ on $A/\tilde{L}_w$.  Then the map $F_w$ is defined by 
\begin{equation*}
F_w(A,\iota,\lambda,\{H_v\}_{v\in S_p},\eta)=(\tilde{A},\tilde{i},\tilde{\lambda},\{\tilde{H}_v\}_{v\in S_p},\tilde{\eta})
\end{equation*}
where $\tilde{A}$, and $\tilde{\lambda}$ are as described above, $\tilde{i}$ is the induced action of $\O_F\otimes\Z_{(p)}$, $\tilde{H}_v=\pi H_v$ for $v\not=w$, $\tilde{H}_w=H_w^{(p)}\subset A[w^\infty]^{(p)}\simeq (A/L_w)[w^\infty]$, and $\tilde{\eta}=\frac{1}{p}\pi\eta$.

Note that this definition depends on the choice of $x_w$.  \end{defn}

To explain why we call the map $F_w$ a partial Frobenius, observe that according to the product decomposition
\begin{equation*}
A[p^\infty]=\prod_{v|p} A[v^\infty]
\end{equation*}
we have $\tilde{A}[w^\infty]=\tilde{A}[w^\infty]/L_w\simeq A[w^\infty]^{(p)}$ while $\tilde{A}[v^\infty]\simeq A[v^\infty]$ for all $v\not=w$.  In particular, according to the local product structure of $Y_1^I$ coming from the Serre--Tate theorem and the product decomposition of the $p$-divisible group $A[p^\infty]$, $F_w$ looks like Frobenius on the factor corresponding to $w$.

As a consequence of this we may record
\begin{prop}
$F_w$ is finite flat of degree $p^3$.  It restricts to a map $F_w:Y_1^{I,=_w1}\to Y_1^{I,=_w1}$ which is finite flat of degree $p^2$.
\end{prop}
\begin{proof}
Using the Serre--Tate theorem and the description of $F_w$ on the $p$-divisible group above, this follows from the fact that Frobenius on a smooth variety of dimension $n$ is finite flat of degree $p^n$.
\end{proof}

The identification $\tilde{A}[w^\infty]\simeq A[w^\infty]^{(p)}$ induces a canonical isomorphism $F_w^*\omega_w\simeq\omega_w^p$ while the isomorphisms $\tilde{A}[v^\infty]\simeq A[v^\infty]$ for $v\not=w$ induce canonical isomorphisms $F_w^*\omega_v\simeq\omega_v$.

The point of this definition is that if we identify $\fC^I_w$ with
$\fY^{I,=_w2}$ via $p_2$, the map $p_1$ on the special fibre is simply
the partial Frobenius $F_w$ restricted to $Y_1^{I,=_w2}$.  Moreover
making these identifications, the isogeny $p_1^*\cG_w\to p_2^*\cG_w$ becomes $V:\cG_w^{(p)}\to\cG_w$ (as its dual is Frobenius) and so the pullback map $p_2^*\omega_w\to
p_1^*\omega_w$ becomes $\Ha_w:\omega_w\to\omega_w^p$.

As in \S\ref{section-traceandres}, we may consider trace maps for $F_w$ on differentials
$F_{w,*}K_{Y_1^I}\to K_{Y_1^I}$.  Tensoring with any line bundle $\L$
on $Y_1^I$ and using the projection formula
$F_{w,*}K_{Y_1^I}\otimes\L\simeq F_{w,*}(K_{Y_1^I}\otimes F_w^*\L)$, we
obtain a twisted trace map
\begin{equation*}
F_{w,*}(K_{Y_1^I}\otimes F_w^*\L)\to K_{Y_1^I}\otimes \L.
\end{equation*} 
We may similarly consider twisted trace maps for line bundles on the divisor $Y_1^{I,=_w1}$.

Now we restrict to parallel weight 2 and work on the special fibre.
With the identifications we have made, our discussion above
shows that we have a commutative diagram of sheaves on $Y_1^{I,=_2w}$
\begin{equation*}
\xymatrix{
F_{w,*}\omega^2\ar[rr]^{U_{\Iw(w),1}}\ar[d]&&\omega^2\ar[d]\\
F_{w,*}(K_{Y^{I,=_w2}_1}\otimes\omega^{-1})\ar[r]^{\frac{N_{F/\Q}(x_w)^3}{p^3}\Ha_w^{-1}}&F_{w,*}(K_{Y_1^{I,=_w2}}\otimes\omega^{-1}\otimes\omega_w^{1-p})\ar[r]&K_{Y_1^{I,=_w2}}\otimes\omega^{-1}
}
\end{equation*}

Now we want to extend this description to all of $Y_1^I$.  Here is the first main result of this section.
\begin{prop}\label{prop:HaUw1}
The map $\Ha_w \cdot U_{\Iw(w),1}:F_{w,*}\omega^2\to\omega^2\otimes\omega_w^{p-1}$ of sheaves on $Y_1^{I,=_w2}$  extends to $Y^I_1$ and fits in to a commutative diagram of sheaves on $Y^I_1$
\begin{equation*}
\xymatrix{
F_{w,*}\omega^2\ar[rr]^{\Ha_w \cdot U_{\Iw(w),1}}\ar[d]&&\omega^2 \otimes \omega_w^{p-1}\ar[d]\\
F_{w,*}(K_{Y^{I,=_w2}_1}\otimes\omega^{-1})\ar[r]^{\frac{N_{F/\Q}(x_w)^3}{p^3}\Ha_w^{p-1}}&F_{w,*}(K_{Y_1^{I,=_w2}}\otimes\omega^{-1}\otimes\omega_w^{(p-1)^2})\ar[r]&K_{Y_1^{I,=_w2}}\otimes\omega_w^{p-1}\otimes \omega^{-1}
}
\end{equation*}
\end{prop}
\begin{proof}
To prove the proposition, it suffices to establish the commutativity of the diagram over
$Y^{I,=_w2}_1$, as the vertical maps are isomorphisms, and the maps on the bottom are already defined over $Y_1^I$.  This commutativity follows from the discussion above and the fact that $F_w^*\Ha_w=\Ha_w^p$.
\end{proof}

Now we are going to restrict to the divisor $Y_1^{I,=_w1}$.  The
Kodaira--Spencer isomorphism $K_{Y_1^I}\simeq \omega^3$ of
Proposition~\ref{prop: KS on divisor} induces by the adjunction formula an isomorphism
\begin{equation*}
K_{Y_1^{I,=_w1}}\simeq\omega^3\otimes\omega_w^{p-1}|_{Y_1^{I,=_w1}}.
\end{equation*}

\begin{prop}\label{prop:trrestrict}
There is a commutative diagram of sheaves on $Y_1^I$
\begin{equation*}
\xymatrix@C=0.5cm{
F_{w,*}(K_{Y^I_1}\otimes\omega^{-1})\ar[r]^{\Ha_w^{p-1}}\ar[d]&F_{w,*}(K_{Y^I_1}\otimes\omega^{-1}\otimes\omega_w^{(p-1)^2})\ar[r]&K_{Y^I_1}\otimes\omega^{-1}\otimes\omega_w^{p-1}\ar[d]\\
F_{w,*}(K_{Y_1^{I,=_w1}}\otimes\omega^{-1}\otimes\omega_w^{1-p})\ar[rr]&&K_{Y_1^{I,=_w1}}\otimes\omega^{-1}
}
\end{equation*}
where the vertical maps are obtained from restriction and the adjunction formula as recalled above, and the bottom horizontal map is a twisted trace for $F_w$ on the divisor $Y_1^{I,=_w1}$ and the line bundle $\omega^{-1}$.
\end{prop}
\begin{proof}
The commutativity of this diagram follows from Proposition \ref{prop:
  trace restrict to cartier}, where in the notation of that proposition
we take $X=Y=Y_1^I$, $D'=D=Y_1^{I,=_w1}$, $f=F_w$ and $n=p$, and
identify $\O_{Y_1^I}(Y_1^{I,=_w1})$ with $\omega_w^{p-1}$ via $\Ha_w$
(tensor the commutative diagram of Proposition \ref{prop: trace
  restrict to cartier} with
$\omega^{-1}\otimes\omega_w^{p-1}$).\end{proof}

We may then define a map $\gamma_1:F_{w,*}(\omega^2|_{Y_1^{I,=_w1}})\to\omega^2\otimes\omega_w^{p-1}|_{Y_1^{I,=_w1}}$ by the diagram
\numequation\label{KS*}
\xymatrix{
F_{w,*}(\omega^2|_{Y_1^{I,=_w1}})\ar[r]^{\gamma_1}\ar[d]_{\wr}&\omega^2\otimes\omega_w^{p-1}|_{Y_1^{I,=_w1}}\ar[d]^{\wr}\\
F_{w,*}(K_{Y_1^{I,=_w1}}\otimes\omega^{-1}\otimes\omega_w^{1-p})\ar[r]&K_{Y_1^{I,=_w1}}\otimes\omega^{-1}
}
\end{equation}
where the vertical maps are Kodaira--Spencer and the bottom horizontal map is $\frac{N_{F/\Q}(x_w)^3}{p^3} $ times the twisted trace for $\omega^{-1}$.

Now combining Proposition \ref{prop:HaUw1} with Proposition \ref{prop:trrestrict} with the definition of $\gamma_1$ by~\eqref{KS*} we have proved the following.

\begin{prop}\label{prop: main U1 diagram}
There is a commutative diagram
\begin{equation*}
\xymatrix{
H^0(Y_1^I,\omega^2)\ar[r]^-{\Ha_w \cdot U_{\Iw(w),1}}\ar[d]& H^0(Y_1^I,\omega^2\otimes\omega_w^{p-1})\ar[d]\\
H^0(Y_1^{I,=_w1},\omega^2|_{Y_1^{I,=_w1}})\ar[r]^-{\gamma_1}&H^0(Y_1^{I,=_w1},\omega^2\otimes\omega_w^{p-1}|_{Y_1^{I,=_w1}})
}
\end{equation*}
where the vertical maps are restrictions and the horizontal maps are as explained above.
\end{prop}

\subsection{The Hecke operator \texorpdfstring{$U_{w,2}$}{U2} on the \texorpdfstring{$w$}{w}-non ordinary locus.}\label{section:U2equalsF2}

In this section we consider the Hecke operator $U_{w,2}$ that was first introduced in \S\ref{subsubsec:
Klingen type padic Hecke operators}.   We consider the  correspondence
\begin{equation*}
\xymatrix{
&\fC_{w,2}^I\ar[dl]_{p_2}\ar[dr]^{p_1}&\\
\fY^I&&\fY^I
}
\end{equation*}which is the composition of the correspondences $\mathfrak{C}_{w,2,1}(p)$ and $\mathfrak{C}_{w,2,2}(p)$ considered in \S\ref{subsubsec:
Klingen type padic Hecke operators} (or more precisely their
restrictions to the interior of the moduli space). The correspondence
$\fC_{w,2}$  admits the following direct description: it parametrizes
isogenies $p_1^\star \cG \rightarrow p_2^\star \cG$ whose kernel $K_w$
is a totally isotropic subgroup of $\cG_w[p^2]$ which has trivial
intersection with the group $p_1^\star H_w$. To see this, note that  $K_w$ fits into an exact sequence $0 \rightarrow K_w[p] \rightarrow K_w \rightarrow K_w/K_w[p] \rightarrow 0$ where $K_w[p]$ is a finite flat group scheme of rank $p^3$ and \'etale rank $p$, and $K_w/K_w[p]$ is a finite \'etale group scheme of rank $p$.

There is yet another description of $\mathfrak{C}^I_{w,2}$ that will be important for us.
To any point   $(\cG, \iota, \lambda, \{H_v\}_{v \in S_p}, \eta) \in
\fY^I$ we can associate a subgroup $\tilde{L}_w \subset \cG[p^2]$ as
follows: the finite flat group scheme $x_w^{-1}H_w/H_w^\perp \subset
\cG_w/H_w^\perp$ of degree $p^2$ contains a canonical multiplicative subgroup
$L'_w$ of degree $p$ (as $x^{-1}H_w/\cG_w[p]\simeq H_w$ is multiplicative and $\cG_w[p]/H_w^\perp\simeq H_w^D$ is \'etale, we see that $x_w^{-1}H_w/H_w^\perp$ is isomorphic at geometric points to $\mu_p\times\Z/p\Z$, and
hence over the entire (reduced) special fibre, the kernel of Frobenius
on $x_w^{-1}H_w/H_w^\perp$ is a multiplicative group of order $p$  which lifts uniquely over $\fY^I$).  Then we may define the group $L_w$ to be the preimage of $L'_w$ under the isogeny $A\to A/H_w^\perp$. Observe that  $L_w\subset A[w^2]$ is a totally isotropic subgroup of degree $p^4$. Then we take
\begin{equation*}
\tilde{L}_w=L_w\times \prod_{v\not=w}A[v^2].
\end{equation*}

We temporarily write $\fY^I_{w-\mathrm{sph}}$ for the formal completion of $Y_{K_pK^p}$ with $K_p=\GSp_4(\O_{F_w})\prod_{v\in I,v\not=w}\Kli(v)\prod_{v\in I^c}\Iw(v)$, along the open subvariety of the special fibre where $A[w^\infty]$ has $p$-rank $\geq 1$, $H_v$ is multiplicative for $v\in I$, $v\not=w$, and $L_v$ is multiplicative for $v\in I^c$.  Then there is a natural map $ f: \fY^I\to \fY^I_{w-\mathrm{sph}}$ which forgets the Klingen level structure $H_w$ at $w$.  It is \'{e}tale and affine (see Proposition~ \ref{prop: deep Klingen is affine etale}).

We define a  map $\psi_w:\fY^I\to \fY^I_{w-\mathrm{sph}}$ by  sending a point $(\cG,\iota,\lambda,\{H_v\}_{v\in S_p}, \eta)$ to $\cG/\tilde{L}_w$ with the polarization descended from $\frac{p^4}{x_w^2}\lambda$, the induced action of $\O_F$ and the level structures $p^{-2}\pi\eta$ and $\pi H_v$ for $v\not=w$ where $\pi: \cG \rightarrow \cG/\tilde{L}_w$ is the isogeny.

\begin{lem}\label{lem-com-diag-Uw2} The correspondence $\fC_{w,2}^I$ fits in the following Cartesian diagram
\begin{equation*}
\xymatrix{\fC_{w,2}^I\ar[d]^{p_2}\ar[r]^{p_1}&\fY^I\ar[d]^f\\
\fY^I\ar[r]^{\psi_w}&\fY^I_{w-\mathrm{sph}}
}
\end{equation*}

\end{lem}

\begin{proof} Let $ p_1^\star \cG \rightarrow p_2^\star \cG$ be the universal isogeny over $\fC_{w,2}^I$. Then the composite of this  isogeny with the isogeny $p_2^\star \cG \rightarrow p_2^\star \cG/ \tilde{L}_w$ identifies with multiplication by $p^2$ on $p_1^\star \cG$. 
\end{proof} 

\begin{lem} The map  $p_1: \fC^I_{w,2} \rightarrow \fY^I$ is finite flat of degree $p^4$.
\end{lem}

\begin{proof} The correspondence $\fC^I_{w,2}$ is smooth, and the map $p_1$ is generically \'etale of degree $p^4$ and finite. It follows from miracle flatness that it is finite flat. 
\end{proof} 

We can now  deduce the following important relation between $\fC^I_{w,2}$ and $F_w^2$ over $Y_1^{I,=_w1}$.
\begin{prop}\label{prop: Uw2 is F^2}
The restriction of $p_2$ to the scheme theoretic preimage $p_2^{-1}(Y_1^{I,=_w1})$ is an isomorphism to $Y_1^{I,=_w1}$.  Making this identification, $p_1$ becomes $F_w^2:Y_1^{I,=_w1}\to Y_1^{I,=_w1}$.
\end{prop}
\begin{proof} We observe that the restriction of $f$ to $Y_1^{I,=_w1}$
  is an isomorphism. This implies that the restriction of $p_2$ is an
  isomorphism. Now we examine the definition of $\tilde{L}_w$. The key
  observation is that $L_w$ coincides with the kernel of $\Frob^2:
  \cG_w \rightarrow \cG_w^{(p^2)}$.  It suffices to check this on
  geometric points.  Over a geometric point of $Y_1^{I,=_w1}$, the
  $p$-divisible group $\cG_w$ has a  decomposition into a product of multiplicative, slope $\frac{1}{2}$ and \'etale group $p$-divisible groups: $\cG_w = \cG_w^{et} \times \cG_w^{oo} \times \cG_w^{m}$. Then the kernel of $\Frob^2$ is simply $\cG_w^{oo}[p] \times \cG_w^{m}[p^2]$, and this group equals $L_w$. 
\end{proof} 

Pullback by the universal isogeny over $\fC^I_{w,2}$ induces a morphism
$\delta_0:p_2^*\omega_w\to p_1^*\omega_w$, as well as isomorphisms
$p_2^*\omega_v\to p_1^*\omega_v$ for $v\not=w$. The following proposition is implicitly   contained in Lemma \ref{lem: U operators commute with Hasse}, but we briefly recall the argument. 
\begin{prop}\label{prop: delta hasse}\leavevmode
\begin{enumerate}
\item The map $\delta_0$ is divisible by $p$ and the resulting map $\delta=\frac{1}{p}\delta_0:p_2^*\omega_w\to p_1^*\omega_w$ is an isomorphism.
\item Under the isomorphism $\delta^{p-1}:p_2^*\omega_w^{p-1}\simeq p_1^*\omega_w^{p-1}$ we have $p_1^*\Ha_w=p_2^*\Ha_w$.\end{enumerate}
\end{prop}
\begin{proof}  
Because $\fC_{w,2}^I$ is smooth we are free to  check the first claim on the ordinary locus where it simply follows from the fact that the isogeny $p_1^\star \cG \rightarrow p_2^\star \cG$ has kernel of multiplicative rank one. The second claim follows from \cite[Lem.\ 10.5.2.1]{pilloniHidacomplexes}. 
\end{proof}

Making the identifications of the Proposition~\ref{prop: Uw2 is F^2}, we may view the restriction of $\delta$ to $Y_1^{I,=_w1}$ as an isomorphism
\begin{equation*}
\delta|_{Y_1^{I,=_w1}}^{-1}:\omega_w\to (F_w^2)^*\omega_w\simeq\omega_w^{p^2}
\end{equation*}
or equivalently as a non vanishing section
$\delta|_{Y_1^{I,=_w1}}^{-1}\in H^0(Y_1^{I,=_w1},\omega_w^{p^2-1})$.
We also denote by  $\delta': \prod_{v \neq w}
p_2^*\omega_w\simeq \prod_{v \neq w} p_1^*\omega_w$ the isomorphism
coming from the pullback of differentials.

In weights $\kappa=(k_v,l_v)$ with $l_w\geq 0$, the Hecke operator
\begin{equation*}
U_{w,2}:p_{1,*}p_2^*\omega^\kappa\to\omega^\kappa
\end{equation*}
is defined by tensoring the unnormalized trace map $p_{1,*}p_2^*\O_{Y_1^I}\to\O_{Y_1^I}$ with the unnormalized pullback map $p_2^*\omega^\kappa\to p_1^*\omega^\kappa$, and normalizing by a factor of $\frac{1}{p^{3+l_w}}$ (see \S\ref{subsubsec:
Klingen type padic Hecke operators}; equivalently, the normalized map $p_2^*\omega^\kappa\to p_1^*\omega^\kappa$ is constructed with the help of the operator $\delta$).

First we may use the Kodaira--Spencer isomorphism to describe $U_{w,2}$ in weight $3$ in terms of traces on differentials.

\begin{prop}\label{prop: Uw2 KS} 
There is a commutative diagram of sheaves on $\fY^I$
\begin{equation*}
\xymatrix{
p_{1,*}p_2^*\omega^3\ar[r]^{U_{w,2}}\ar[d]&\omega^3\ar[d]\\
p_{1,*}p_2^*K_{\fY^I/\Z_p}\ar[r]^{\frac{N_{F/\Q}(x_w)^6}{p^6}\mathrm{tr}}&K_{\fY^I/\Z_p}
}
\end{equation*}
where the vertical arrows are the Kodaira--Spencer isomorphism, and the bottom horizontal arrow is
defined as follows: since~$p_2$ is \'{e}tale, we may identify
$p_2^*K_{\fY^I/\Z_p}$ with $K_{\fC_{w,2}^I/\Z_p}$, and the morphism then comes from
 the trace map for $p_1$ on dualizing sheaves, multiplied by a factor of $\frac{N_{F/\Q}(x_w)^6}{p^6}\in\mathbb{Z}_{(p)}^\times$.
\end{prop}
\begin{proof}
This follows from Propositions \ref{prop: KS isogeny} and \ref{prop: KS traces}.
\end{proof}

In parallel weight $2$ we can still express the cohomological correspondence $U_{w,2}$ by using  a similar commutative diagram of sheaves on $\fY^I$: 
\begin{equation*}
\xymatrix{
p_{1,*}p_2^*\omega^2\ar[rrr]^{U_{w,2}}\ar[d]&&&\omega^2\ar[d]\\
p_{1,*}p_2^*(K_{\fY^I/\Z_p} \otimes \omega^{-1}) \ar[rrr]^{\frac{N_{F/\Q}(x_w)^6}{p^6}\text{tr} \otimes (\delta' \delta)^{-1} }&&&K_{\fY^I/\Z_p} \otimes \omega^{-1}
}
\end{equation*}

We can restrict   to $Y^{I,=_w1}$ and obtain the following: 
\begin{prop}\label{prop: big U2 diagram}
There is a commutative diagram
\begin{equation*}
\xymatrix@C=0.5cm{
p_{1,*}p_2^*(K_{Y^I_1}\otimes\omega^{-1})\ar[r]\ar[d]&K_{Y^I_1}\otimes\omega^{-1}\ar[d]\\
p_{1,*}p_2^*(K_{Y^{I,=_w1}_1}\otimes\O_{Y^I_1}(-Y^{I,=_w1}_1)|_{Y^{I,=_w1}_1}\otimes\omega^{-1})\ar[r]\ar[d] &K_{Y^{I,=_w1}_1}\otimes\O_{Y^I_1}(-Y^{I,=_w1}_1)|_{Y^{I,=_w1}_1}\otimes\omega^{-1}\ar[d]\\
(F_w^2)_*(K_{Y_1^{I,=_w1}}\otimes\omega^{-1}\otimes\omega_w^{1-p})\ar[r]&K_{Y_1^{I,=_w1}}\otimes\omega^{-1}\otimes\omega_w^{1-p}
}
\end{equation*}
where:
\begin{itemize}
\item The upper vertical maps are obtained by restriction to $Y^{I,=_w1}_1$ and the adjunction isomorphism $K_{Y^I_1}|_{Y^{I,=_w1}}\simeq K_{Y_1^{I,=_w1}}\otimes\O_{Y^I_1}$.
\item The lower vertical maps are obtained from making the identification of $p_2$ with $\id$ and $p_1$ with $F_w^2$ of Proposition \ref{prop: Uw2 is F^2} as well as using the isomorphism $\Ha_w:\O(Y^{I,=_w1}_1)\to\omega_w^{p-1}$.
\item The top horizontal arrow is the composition of $(\delta
  \delta')^{-1}$ and the twisted trace for $\omega^{-1}$ on
  $K_{Y_1}^I$, as on the bottom row of the diagram immediately
  preceding this proposition.
\item The middle horizontal arrow is multiplication by $(\delta \delta')\vert_{Y_1^{I,=_w1}}^{-1}$ followed by the twisted trace for $p_1$ on the sheaf $\O_{Y_1^I}(-Y_1^{I,=_w1})\otimes\omega^{-1}$.
\item The bottom horizontal arrow is  multiplication by $\delta|_{Y_1^{I,=_w1}}^{-p} \delta'\vert_{Y_1^{I,=_w1}}^{-1} $ followed by the twisted trace for $F_w^2$ on the line bundle $\omega^{-1}\otimes\omega_w^{1-p}$ \emph{(}normalized by the $p$-adic unit $\frac{N_{F/\Q}(x_w)^6}{p^6}$\emph{)}.
\end{itemize}
\end{prop}\begin{proof}
The commutativity of the top square follows from Proposition \ref{prop: trace restrict to cartier} while the commutativity of the bottom square follows from Proposition \ref{prop: delta hasse}. The reason we get multiplication by $\delta|_{Y_1^{I,=_w1}}^{-p} \delta'\vert_{Y_1^{I,=_w1}}^{-1} $ in the bottom horizontal arrow is because we multiply the original $(\delta \delta')\vert_{Y_1^{I,=_w1}}^{-1}$ with $\delta|_{Y_1^{I,=_w1}}^{-1-p}$ which arises when relating the isomorphisms $p_i^\star (\Ha_w:\O(Y^{I,=_w1}_1)\to\omega_w^{p-1})$ for $i=1,2$. 
\end{proof}

Our goal from now on is to interpret  the bottom horizontal line of this diagram  in terms of the map $\gamma_1$ of ~\eqref{KS*}. 
We introduce a map $\gamma_2:F_{w,*}(\omega^2\otimes\omega_w^{p-1}|_{Y_1^{I,=_w1}})\to \omega^2|_{Y_1^{I,=_w1}}$ defined by the diagram
\numequation\label{KS**}
\xymatrix@C=0.5cm{
F_{w,*}(\omega^2\otimes\omega_w^{p-1}|_{Y_1^{I,=_w1}})\ar[rr]^{\gamma_2}\ar[d]_{\wr}&&\omega^2|_{Y_1^{I,=_w1}}\ar[d]^{\wr}\\
F_{w,*}(K_{Y_1^{I,=_w1}}\otimes\omega^{-1})
 \ar@/_1.5pc/[r]_{\delta|_{Y_1^{I,=_w1}}^{-1}}
&F_{w,*}(K_{Y_1^{I,=_w1}}\otimes\omega^{-1}\otimes\omega_w^{1-p^2})\ar[r]&K_{Y_1^{I,=_w1}}\otimes\omega^{-1}\otimes\omega_w^{1-p}
}
\end{equation}
where the vertical maps are induced by Kodaira--Spencer, and on the bottom the first horizontal map is multiplication by $\delta|_{Y_1^{I,=_w1}}^{-1}$, while the second is the twisted trace for $\L=\omega^{-1}\otimes\omega_w^{1-p}$ (normalized by the $p$-adic unit $\frac{N_{F/\Q}(x_w)^3}{p^3}$).

We now consider  the composition $\gamma_2\circ\gamma_1: (F_w^2)_*(\omega^2|_{Y_1^{I,=_w1}}) \rightarrow  \omega^2|_{Y_1^{I,=_w1}} $.
\begin{prop}\label{prop:F2compare}
There is a commutative diagram
\begin{equation*}
\xymatrix@C=0.2cm{
(F_w^2)_*(\omega^2|_{Y_1^{I,=_w1}})\ar[r]^-{\gamma_1}\ar[d]_{\wr} &F_{w,*}(\omega^2\otimes\omega_w^{p-1}|_{Y_1^{I,=_w1}})\ar[r]^-{\gamma_2}&\omega^2|_{Y_1^{I,=_w1}}\ar[d]^{\wr}\\
(F_w^2)_*(K_{Y_1^{I,=_w1}}\otimes\omega^{-1}\otimes\omega_w^{1-p})
 \ar@/_1.5pc/[r]_{\delta^{-p}\delta'^{-1}|_{Y_1^{I,=_w1}}}
&(F_w^2)_*(K_{Y_1^{I,=_w1}}\otimes\omega^{-1}\otimes\omega_w^{1-p^3})\ar[r]&K_{Y_1^{I,=_w1}}\otimes\omega^{-1}\otimes\omega_w^{1-p}}
\end{equation*}
where the vertical arrows are given by the Kodaira--Spencer isomorphism and on the bottom row we first multiply by $\delta|_{Y_1^{I,=_w1}}^{-p}$ and then take a twisted trace for $F_w^2$ \emph{(}normalized by the $p$-adic unit $\frac{N_{F/\Q}(x_w)^6}{p^6}$\emph{)} and the line bundle $\L=\omega^{-1}\otimes\omega_w^{1-p}|_{Y_1^{I,=_w1}}$.
\end{prop}
\begin{proof}
This follows from the fact that $F_w^*\delta|_{Y_1^{I,=_w1}}=\delta|_{Y_1^{I,=_w1}}^p$.
\end{proof}

Combining Proposition \ref{prop: big U2 diagram} with Proposition \ref{prop:F2compare} we have proved the following:

\begin{prop}\label{prop: main U2 diagram}
There is a commutative diagram
\begin{equation*}
\xymatrix{
H^0(Y_1^I,\omega^2)\ar[r]^-{U_{w,2}}\ar[d]& H^0(Y_1^I,\omega^2)\ar[d]\\
H^0(Y_1^{I,=_w1},\omega^2|_{Y_1^{I,=_w1}})\ar[r]^-{\gamma_2\circ \gamma_1}&H^0(Y_1^{I,=_w1},\omega^2|_{Y_1^{I,=_w1}})
}
\end{equation*}
where the vertical maps are restrictions and the horizontal maps are as explained above.
\end{prop}

\subsection{Main doubling results} \label{subsec: doubling preliminaries}
 There  is an obvious   injective restriction map:
$$ H^0 ( X^{I}_{1}, \omega^2(-D)) \rightarrow H^0 ( X^{I,=_w2}_{1}, \omega^2(-D)) $$
which is equivariant for the action of the Hecke algebra away from
$w$, and for the actions of $U_{w,2}$ and $U_{w,0}$. We now compare
the action of $U_{\Kli(w), 1}$, which acts on both
the left hand and right hand modules, and $U_{\Iw(w),1}$, which acts
on the right hand module.

We have defined a Hecke operator~$Z_w$ on $H^0 ( X^{I,=_w2}_{1},
\omega^2(-D)) $ with~$U_{\Kli(w),1}=U_{\Iw(w),1}+Z_w$ (see \S\ref{subsec: UIw1 and Zw}). 

    \begin{lem}\label{lem: doubling map quadratic relation}On $H^0 (
  X^{I,=_w2}_{1}, \omega^2(-D))$ we have the identity of
  operators $U_{\Iw(w),1} Z_w =
  U_{w,2}$.
    \end{lem}
  \begin{proof} This is immediate from Lemma~ \ref{lem: quadratic relation for U Kli 1 big sheaf}.

\end{proof}

We introduce the doubling map: \begin{eqnarray*}
 H^0 ( X^{I}_{1}, \omega^2(-D))  \oplus H^0 ( X^{I}_{1}, \omega^2(-D))& \rightarrow &H^0 ( X^{I,=_w2}_{1}, \omega^2(-D))\\
 (f,g) & \mapsto & f + Z_w g 
 \end{eqnarray*}
 In this formula $f$ and $g$ on the right hand side are viewed as sections of $H^0 ( X^{I,=_w2}_{1}, \omega^2(-D))$ via the above restriction map. 
 
We can define an operator that we formally denote by  $U_{\Iw(w), 1}$  on the left hand side by the following matrix:  
 $$ U_{\Iw(w), 1} =    \begin{pmatrix} U_{\Kli(w),1} & U_{w,2} \\
       -1 & 0 \\
    \end{pmatrix}$$
    \begin{lemma}\label{lem: doubling U1 equivariant} The  doubling map  is equivariant for the action of $U_{\Iw(w),1}$. The operator $U_{\Iw(w),1}$ on  $H^0 ( X^{I}_{1}, \omega^2(-D)) \oplus H^0 ( X^{I}_{1}, \omega^2(-D))$ commutes with the action of $U_{w,2}$.  
    \end{lemma} 

\begin{proof} The equivariance follows from Lemma~\ref{lem: doubling map quadratic
    relation}, and the commutativity follows from Lemma~\ref{lem: commutativity of Hecke ops on infinite Klingen
    big sheaf}. \end{proof}

We now consider the $U_{w,2}$-ordinary part: $$ e(U_{w,2}) H^0 ( X^{G_1,I}_{1}, \omega^2(-D))  \oplus e(U_{w,2})H^0 ( X^{G_1,I}_{1}, \omega^2(-D)).$$  We  have restricted to the direct factor $H^0 ( X^{G_1,I}_{1}, \omega^2(-D))$ of $H^0 ( X^{I}_{1}, \omega^2(-D))$ in order to be able to use local finiteness and apply ordinary projectors.

\begin{lem}\label{lem: doubling preserves ordinarity}\leavevmode \begin{enumerate}
\item  The image of $$ e(U_{w,2}) H^0 ( X^{G_1,I}_{1}, \omega^2(-D))  \oplus e(U_{w,2})H^0 ( X^{G_1,I}_{1}, \omega^2(-D))$$ via the doubling map lands in  $e(U_{\Iw(w),1}U_{w,2}) H^0 ( X^{G_1,I,=_w2}_{1}, \omega^2(-D))$. 
\item The operator  $U_{\Iw(w),1}$ is \emph{(}left and right\emph{)} invertible on 
$$e(U_{w,2}) H^0 ( X^{G_1,I}_{1}, \omega^2(-D))  \oplus e(U_{w,2})H^0 ( X^{G_1,I}_{1}, \omega^2(-D)).$$\item  On $e(U_{w,2}) H^0 ( X^{G_1,I}_{1}, \omega^2(-D))  \oplus e(U_{w,2})H^0 ( X^{G_1,I}_{1}, \omega^2(-D))$ we have the identity $U_{\Kli(w),1} = U_{\Iw(v),1} + U_{w,2} U_{\Iw(v),1}^{-1}$. In particular $U_{\Kli(w),1}$ and $U_{\Iw(w),1}$ commute with each other. 
\end{enumerate}
\end{lem}

\begin{proof} The doubling map  \[  H^0 ( X^{G_1,I}_{1}, \omega^2(-D))  \oplus H^0 ( X^{G_1,I}_{1}, \omega^2(-D)) \rightarrow  H^0 ( X^{G_1,I,=_w2}_{1}, \omega^2(-D))\]  can be written as an inductive limit of maps between finite dimensional vector spaces stable under the $U_{\Iw(w),1}$ and $U_{w,2}$ operators, so we will freely use the usual properties of  linear endomorphisms on finite  dimensional vector spaces.  We first observe that the operator $U_{w,2}$ is invertible on $ e(U_{w,2}) H^0 ( X^{G_1,I}_{1}, \omega^2(-D))  $. Concretely, for any $f \in e(U_{w,2})H^0 ( X^{G_1,I}_{1}, \omega^2(-D))$ we have $e(U_{w,2}) f = f = U_{w,2}^{p^{N !}}f$ for $N$ large enough, and $U_{w,2}^{-1} f = U_{w,2}^{p^{N !}-1}f$.  Therefore we may consider the operator $$    U_{w,2}^{-1}\begin{pmatrix} 0  & -U_{w,2} \\
       1 & U_{\Kli(w),1} \\
    \end{pmatrix}$$ on $ e(U_{w,2}) H^0 ( X^{G_1,I}_{1},
    \omega^2(-D))  \oplus e(U_{w,2})H^0 ( X^{G_1,I}_{1},
    \omega^2(-D))$, and it is straightforward to check (using
    that~$U_{\Kli(w),1}$ and~$U_{w,2}$ commute, as we noted in the
    proof of Lemma~\ref{lem: doubling U1 equivariant}) that this is a 2-sided
    inverse of $U_{\Iw(w),1}$. This proves the first  and second
    points. The third point is obvious from the formulae
    defining~$U_{\Iw(w),1}$ and~$U_{\Iw(w),1}^{-1}$. 
    \end{proof} 

We now prove our doubling theorems, combining ingredients from the previous sections.
\begin{thm}\label{thm: main doubling thm}Suppose that $w\in I$ and that $f\in H^0(X_{1}^{G_1,I},\omega^2(-D))$
satisfies~$U_{\Kli(w),1}f=(\alpha_w+\beta_w)f$, $U_{\Iw(w),1}f=\alpha_wf$, and $U_{w,2}f=\alpha_w\beta_w f$, where $\alpha_w,\beta_w\ne 0$. Then $f=0$.
\end{thm}
\begin{proof}First suppose that the restriction of $f$ to $X_1^{G_1,I,=_w1}$ is zero.  Then we may write $f=\Ha_wg$ for some $g\in H^0(X_1^{G_1,I},\omega^2\otimes\omega_w^{1-p}(-D))$.  Moreover because we have $Z_wf=\beta_wf$ by hypothesis, we would then have $Z_wg=\beta_wg$ by Proposition \ref{prop-commuteagain}.  But then by Corollary \ref{cor: the S conjecture}, $g=0$ and hence $f=0$.

Now we may suppose that the restriction of $f$ to $X_1^{G_1,I,=_w1}$ is nonzero.  Combining Proposition \ref{prop: main U1 diagram} with Proposition \ref{prop: main U2 diagram} there is a commutative diagram
\begin{equation*}
\xymatrix{
H^0(Y_1^{G_1,I},\omega^2)\ar[d]^{\Ha_wU_{\Iw(w),1}} \ar[r] \ar@/_5pc/[dd]_{U_{w,2}}
& H^0(Y_1^{G_1,I,=_w1},\omega^2|_{Y_1^{G_1,I,=_w1}})  \ar[d]^-{\gamma_1} \\
H^0(Y_1^{G_1,I},\omega^2\otimes\omega_w^{p-1})\ar[r] &  H^0(Y_1^{G_1,I,=_w1},\omega^2\otimes\omega_w^{p-1}|_{Y_1^{G_1,I,=_w1}})\ar[d]^-{\gamma_2} \\
H^0(Y_1^{G_1,I},\omega^2)\ar[r] & H^0(Y_1^{G_1,I,=_w1},\omega^2|_{Y_1^{G_1,I,=_w1}})
}
\end{equation*}
where the horizontal maps are the natural restriction maps, and the vertical
maps on the right column are as in diagrams~\eqref{KS*} and~\eqref{KS**}.

If we start with $f$ in the top left of the diagram, we obtain
something nonzero on the bottom right because
$U_{w,2}f=\alpha_w\beta_wf$ and the restriction of $f$ to
$Y_1^{G_1,I,=_w1}$ is nonzero.  The commutativity of the diagram  implies that $\Ha_wU_{\Iw(w),1}f$ has nonzero restriction to $Y_1^{G_1,I,=_w1}$.  On the other hand, because $U_{\Iw(w),1}f=\alpha_wf$, $\Ha_wU_{\Iw(w),1}f=\alpha_w\Ha_wf$ which vanishes along $Y_1^{G_1,I,=_w1}$.  This is a contradiction.\end{proof}

\begin{rem}\label{rem: FJ argument}In fact, something stronger than
  Theorem~\ref{thm: main doubling thm} is true:
if $w\in I$, and $f\in H^0(X^{G_1,I}_{1},\omega^2(-D))$ satisfies
$U_{w,2}f\not=0$, then $U_{\Iw(w),1}f\not\in
H^0(X^{G_1,I}_{1},\omega^2(-D))$. This can be proved in exactly the same way
as Theorem~\ref{thm: main doubling thm}, given the following
strengthening of Corollary~\ref{cor: the S conjecture}:
  if~ $w\in I$, then  
  $$H^0(X_{K_p(I) K^p,1}^{G_1,I},\omega^2\otimes\omega_w^{1-p}(-D))=0.$$
In the case~$p>3$, we will sketch a proof of this result in~\S\ref{section:vanishinginnegativeweight} using Fourier--Jacobi expansions, but since a complete
argument in the case~$p=3$ would involve developing considerably more
of the details of toroidal compactifications than we need in the rest
of the paper, we have decided not to give the details.

When $p>3$, this vanishing result holds
even for non cusp forms, so the same is true of Theorem \ref{thm: main doubling thm}. \end{rem}

We can now prove the injectivity of the doubling map.
\begin{thm}[Doubling]\label{thm: the doubling map is injective} The doubling map $$ e(U_{w,2}) H^0 ( X^{G_1,I}_{1}, \omega^2(-D))  \oplus e(U_{w,2})H^0 ( X^{G_1,I}_{1}, \omega^2(-D))$$ $$ \rightarrow e(U_{\Iw(w),1}U_{w,2}) H^0 ( X^{G_1,I,=_w2}_{1}, \omega^2(-D))$$
is injective. 
\end{thm}

 \begin{proof} Assume the map is not injective. By Lemma~\ref{lem:
     doubling U1 equivariant}, the kernel  is an
   inductive limit of finite dimensional vector spaces stable under
   the commuting operators $U_{w,2}$ and $U_{\Iw(w),1}$. We may
   therefore take a nonzero simultaneous eigenvector $(f,g)$ for  $U_{w,2}$
   and $U_{\Iw(w),1}$ in this kernel, with respective
   eigenvalues~$\alpha_w\beta_w$ and~$\beta_w$ for some~$\alpha_w$, $\beta_w\ne
   0$ (the eigenvalues are nonzero because we are by assumption in the
   ordinary space for both~$U_{\Iw(w),1}$ and~$U_{w,2}$). It follows from the definition
   of the action of ~$U_{\Iw(w),1}$ that $f=-\beta_w g$ and
   $U_{\Kli(w),1}f=(\alpha_w+\beta_w)f$. Since we are also assuming that~$f + Z_w g = 0$, we see that  the image of $f$  in $H^0 ( X^{G_1,I,=_w2}_{1}, \omega^2(-D))$
   satisfies  $U_{\Iw(w),1}f=\alpha_w f$. The result follows from
   Theorem~\ref{thm: main doubling thm} (note that the
   eigenvalues for~$U_{\Iw(w),1}$ and~$U_{w,2}$ are nonzero because we
   are in the ordinary space for these operators by hypothesis).
  \end{proof}

  \begin{rem}\label{rem: doubling as mapping to generalized
      eigenspaces} We now put the Theorem~\ref{thm: the doubling map is injective} in a form that is used 
      in~\S\ref{subsec: TW systems2}. Assume that $M \subset e(U_{w,2}) H^0 ( X^{G_1,I}_{1}, \omega^2(-D))\otimes \bar{\FF}_p$ is a finite dimensional vector space, stable under $U_{\Kli(w),1}$ and $U_{w,2}$. Assume that there are distinct  elements $\alpha_w, \beta_w \in \bar{\FF}_p^\times$  such that $U_{\Kli(w),1}- (\alpha_w+\beta_w)$ and $U_{w,2} - \alpha_w\beta_w$ are nilpotent on $M$.  The sub-algebra $\mathcal{E}$ of $\End( M)$ generated by $U_{\Kli(w),1}$ and $U_{w,2}$ is therefore a local Artinian algebra and  there  are elements $\tilde{\alpha}_w, \tilde{\beta}_w \in \mathcal{E}$  satisfying $\tilde{\alpha}_w = \alpha_w ~\mod \m_{\mathcal{E}}$ and $\tilde{\beta}_w = \beta_w ~\mod \m_{\mathcal{E}}$ and  such that on $M \oplus M$ we have $(U_{\Iw(w),1}- \tilde{\alpha}_w)( U_{\Iw(w),1}- \tilde{\beta}_w) = 0$. 
  
  We can define  maps $\iota_{\xi_w}: M \rightarrow M \oplus M$ by $f \mapsto ( f, - \tilde{\xi}_w^{-1}f)$    for $\xi_w \in \{ \alpha_w,\beta_w\}$. 
  Then one checks easily that the map $\iota_{\alpha_w} \oplus \iota_{\beta_w}: M \oplus M  \rightarrow M \oplus M$ is an isomorphism and that the composite with the doubling map takes  the form $(f_1, f_2) \mapsto ((1- \tilde{\beta}_w U_{\Iw(w),1}^{-1})  f_1 +   (1- \tilde{\alpha}_w U_{\Iw(w),1}^{-1})  f_2 )$.  The first and second components of this map  therefore define injective maps  \[M \hookrightarrow e(U_{\Iw(w),1}U_{w,2}) H^0 ( X^{G_1,I,=_w2}_{1}, \omega^2(-D))_{ U_{\Iw,1}- \xi_{w}}.\]  for $\xi_w$ respectively equal to $ \alpha_w$ and $\beta_w$. 
    \end{rem}

\subsection{Vanishing in partial negative
  weight: Fourier--Jacobi expansions} \label{section:vanishinginnegativeweight}We end this section by giving a  proof of the
following vanishing result in ``partial negative weight'', which partially strengthens Corollary \ref{cor: the S conjecture} but is not needed in this paper (see also Remark \ref{rem: FJ argument}).
\begin{prop}
  \label{prop: low weight vanishing} Assume $w\in I$. If~$p > 3$ and $[F: \Q]>1$, then  
  $$H^0(X_{K_p(I) K^p,1}^I,\omega^2\otimes\omega_w^{1-p})=0.$$
   \end{prop}

\begin{rem} We have a sketch of an argument to show that  if~$p = 3$, then    $H^0(X_{K_p(I) K^p,1}^I,\omega^2\otimes\omega_w^{1-p}(-D))=0$. We also have a sketch of an argument for $F = \Q$. But to give  complete proofs would require us
to justify certain properties of Fourier--Jacobi expansion for which we could not find  references (for example we would need  to have good geometric theory of cuspidal Fourier--Jacobi forms).
\end{rem}

We will prove Proposition~\ref{prop: low weight vanishing} by restriction to a boundary stratum, and
ultimately reducing to the vanishing of spaces of Hilbert modular forms
of partial negative weight.

 We let $K_p = \prod_{v |p} \mathrm{GSp}_4(\ocal_{F_v})$, and by possibly shrinking $K^p$ we may assume that it is a principal level structure in the sense of \cite[\S1.3.6]{MR3186092}.  We let $c \in  \ZZ_{(p)}^{\times,+} \backslash (\A^{\infty,p}\otimes F)^\times/
  \nu(K^p)$,
 and we may work with the connected component $X_{K,1,c}$ of $X_{K,1}$. We now choose a boundary stratum $Z \hookrightarrow X_{K,1,c}$ corresponding to a one dimensional totally isotropic factor $W \in \mathfrak{C}$  (see \S\ref{sec: compactifications}). It means that the restriction of the semi-abelian scheme along $Z$ is an extension of an abelian scheme $A$ of dimension $[F: \Q]$ with $\ocal_F$-action  by a torus $T$   of dimension $[F: \Q]$ with $\ocal_F$-action. 
 
Let $H:=\ker(\Res_{F/\Q}\GL_2\to(\Res_{F/\Q}\Gm)/\Gm)$. The abelian scheme of dimension  $[F: \Q]$ is parametrized by a (connected) Shimura variety for the group $H$ (this is a Hilbert--Blumenthal modular variety) that we denote by $Y_{H,1}$ and is a moduli space of isomorphism classes of  triples $(A, \iota, \lambda, \eta)$: 

\begin{enumerate}
\item $A \rightarrow \Spec R$ is an abelian scheme,
\item $\iota: \CO_F \rightarrow \mathrm{End} (G) \otimes \Z_{(p)}$ is an action,
\item $\mathrm{Lie} (A)$ is a locally free $\CO_F \otimes_{\Z} R$-module of rank $1$, 
\item $\lambda: A \rightarrow A^t$ is a prime to $p$, $\cO_F$-linear quasi-polarization such that for all $v |p$, $\mathrm{Ker}(\lambda: A[v^\infty]
  \rightarrow A^t[v^\infty])$ is trivial, \item $\eta$ is a prime to $p$ level structure. \end{enumerate}

We denote by $X_{H,1}$ a toroidal compactification of $Y_{H,1}$. We have partial Hasse invariants $A_v$ for all $v |p$. Let $Y_{H,1}^I  \subset X_{H,1}^I$ be the Zariski opens where
$A_v$ is invertible  for all~
$v\in I^c$. We have a map $Z \rightarrow Y_{H,1}$ and we let $Z^I = Z \times_{Y_{H,1}} Y_{H,1}^I$.

The \'etale map $X_{K_p(I) K^p,1} \rightarrow X_{K,1}$ has a section along $Z^I \hookrightarrow X_{K,1}$ which is provided by the rank one multiplicative groups $T[v]$ for all $v \in I$. Therefore the map $X_{K_p(I) K^p,1} \rightarrow X_{K,1}$ has a section restricted to  the completion of $X_{K,1}$ along~ $Z^I$.

\begin{prop}[Fourier--Jacobi expansion principle]\label{prop: FJ expansion}There is a natural injective Fourier--Jacobi expansion map  \begin{equation*}
H^0(X_{K_p(I) K^p,1,c}^I,\omega^2\otimes\omega_v^{1-p})\to\prod_{\xi \in \mathfrak{a}^+}H^0(A_1^I,\omega^2\otimes\omega_v^{1-p}\otimes\cL_\xi)
\end{equation*}
where $\mathfrak{a}$ is a fractional ideal of $\ocal_F$ and $\mathfrak{a}^+$ are the positive elements, $A_1^I \rightarrow Y_{H,1}^I$ is an abelian scheme isogenous to the universal abelian scheme $A$ and $\cL_\xi$ is an invertible sheaf over $A_1^I$, rigidified along the identity section.

\end{prop}
\begin{proof}
The existence of such a map  follows from the description of the
toroidal boundary charts, as in~\cite[\S V]{MR1083353} or~\cite[\S
6.2.3,\S7.1]{MR3186092}.  It is obtained by restricting sections to the completion along $Z^I$. 
The sheaves $\cL_\xi$ are obtained by pullback from a Poincar\'e bundle which is rigidified along the identity section. 

The injectivity result is clear as long as  we can show that
$X_{K_p(I) K^p,1,c}^I$ is connected.  This follows directly from the connectedness of $X_{K,1,c}$ and the
irreducibility of the Igusa tower, for which see~\cite[Cor.\ 8.17]{MR2055355}
or~\cite[Thm.\ 0.1]{MR2465518}.\end{proof}

We will now prove the vanishing of the groups $H^0(A_1^I,\omega^2\otimes\omega_w^{1-p}\otimes\cL_\xi)$. We first need the following preliminary lemma.

\begin{lem}\label{lem-torsionsheaf} Let $S$ be a scheme and let $A \rightarrow S$ be an abelian scheme. Let $\cL$ be an invertible sheaf on $A$, rigidified along the unit section. Then for all $n \in \Z_{\geq 1}$, $\cL^{n^3}\vert_{A[n]}$ is trivial. 
\end{lem}
\begin{proof} It is well-known (\cite[Chap.\ II, \S6 and \S8]{mumford}) that $n^\star \cL \simeq \cL^{n^2} \otimes \cL_0$ where $\cL_0$ is a sheaf algebraically equivalent to zero. Moreover $n^\star \cL_0 \simeq  \cL_0^{n}$. Therefore $\cL^{n^3} \simeq (n^\star \cL)^n \otimes n^\star \cL_0^{-1}$ is trivial on $A[n]$. 
\end{proof} 

Now we may prove the following sequence of vanishing results for
negative weight forms. (Note that the first part is a very special
case of the main theorem of~\cite{MR3668413}, although the argument
there is different.)  

\begin{prop} \label{prop:impliesvanishinginnegativeweight}
Assume that $[F:\Q]>1$.  Let $\kappa=(k_v)_{v\in S_p}$ be a weight for $H$ and suppose that there is a $w\in I$ such that $k_w < 0$.
Then:
\begin{enumerate}
\item \label{part:previouspart} $H^0(Y_{H,1},\omega^\kappa)=0$.
\item \label{part:more creative way of saying previouspart} $H^0(Y_{H,1}^I,\omega^\kappa)=0$.
\item  \label{prop:part3Zag} For any~$\xi\in\mathfrak{a}^+$, $H^0(A_1^I,\omega^\kappa\otimes\cL_\xi)=0$.
\end{enumerate}
\end{prop}
\begin{proof} We derive each claim in turn from the previous one:
\begin{enumerate}
\item Let $C_w\subset Y_{H,1}$ be the simultaneous vanishing locus of the Hasse invariants $A_v$ for $v\not= w$; it is a (union of) smooth curves (since~$p$ is split completely,
this is an easy local calculation). 
Furthermore, because~$[F:\Q] > 1$, it is also proper  (note that if~$[F:\Q] = 1$, then~$C_w = Y_{H,1}$ is not proper).

 By the existence of the secondary Hasse invariants, $\omega_v|_{C_w}$
 is a torsion line bundle for $v\not=w$, while $\omega_w|_{C_w}$ has
 positive degree on each component. Let $\cI$ be the ideal sheaf of $C_w$ in $Y_{H,1}$.  It follows from
the Kodaira--Spencer isomorphism that we have an isomorphism\begin{equation*}
\cI/\cI^2=\bigoplus_{v\not=w}\omega_v^2.
\end{equation*}
Thus for all $m\geq0$, $\cI^m/\cI^{m+1}=\Sym^m\cI/\cI^2$ is a direct sum of torsion line bundles.  Because $k_w<0$, it follows that~$\cI^m/\cI^{m+1}\otimes\omega^\kappa$ is a direct sum of line bundles of negative degree, and hence has no sections.  
 The result follows from this and the fact that every irreducible component 
of~$X_{H,1}$ contains a component of~$C_w$ (by considering the formal expansion of any form along~$C_w$).

\item If $f\in H^0(Y_{H,1}^I,\omega^\kappa)$, then, for $c_v\gg 0$ for
  all~$v\in I^c$,
\begin{equation*}
f\prod_{v\in I^c}A_v^{c_v}\in H^0(Y_{H,1},\omega^\kappa\otimes\bigotimes_{v\in I^c}\omega_v^{c_v(p-1)})
\end{equation*}
and hence vanishes by part~\eqref{part:previouspart}. Thus the same conclusion holds for~$f$.
\item  \label{magnificent} 
 For all $n \in \ZZ_{\geq 1}$  with $(n,p)=1$ we will show that any section of $f\in
 H^0(A^I_1,\omega^\kappa\otimes\cL_\xi)$ vanishes on the $n$-torsion subgroup
 $A_1^I[n]$, and hence  vanishes identically.  After replacing $f$ by
 $f^{n^3}$, we can assume that $\cL_\xi $ is the trivial sheaf (see
 Lemma \ref{lem-torsionsheaf}).   We  then consider the norm of $f$
 for each irreducible component of the finite \'{e}tale map $A_1^I[n]
 \rightarrow Y^I_{H,1}$ to reduce to part \eqref{part:more creative
   way of saying previouspart}. \qedhere
\end{enumerate}
 \end{proof}

\begin{proof}[Proof of Prop.~\ref{prop: low weight vanishing}]
This is an immediate consequence of
Proposition~\ref{prop: FJ expansion} and Proposition~\ref{prop:impliesvanishinginnegativeweight}~(\ref{prop:part3Zag}), because all the terms in the Fourier--Jacobi expansion will be zero.
\end{proof}

\section{Higher Coleman theory}\label{sec: higher Coleman
  theory}In this section, we construct (higher) Coleman theories
for~$\GSp_4(\A_F)$.  As in~\S\ref{sec: Hida complexes}, we
assume that~$p$ splits completely in~$F$ and we construct all possible
Coleman theories, allowing the weight space at each place above~$p$ to
be either one or two-dimensional. In the case that $F = \Q$ this was
carried out in~\cite{MR3275848} and~\cite{pilloniHidacomplexes}. Many
of our arguments are simply the ``product over the places~$v|p$'' of
the arguments of~\cite{MR3275848} and \cite{pilloniHidacomplexes}. To
keep this paper at a reasonable length, we will often refer to these
papers for the details of arguments which go over directly to our
case.

The main results of this section are Theorem~\ref{thm-classicality} (a
classicality result for overconvergent cohomology classes of small
slope), and Theorem~\ref{thm: ordinary is overconvergent I equal 1}
(which shows that in the case that~$I$ has size at most one, the
cohomology of the Hida complex~$M_I$ constructed in~\S\ref{sec:
  Hida complexes} is overconvergent, once~$p$ is inverted). These
results together improve (at the expense of inverting~$p$) on the classicality results of~\S\ref{sec: Hida complexes}, in that they do not require the
weight to be sufficiently large; this is crucial for our applications
to abelian surfaces, which correspond to modular forms of parallel
weight~$2$.

We begin in~\S\ref{subsec: formal overconvergent sheaves} with
the construction at the level of formal schemes of a version of the
analytic sheaves of overconvergent forms that we will use later in
this section. The purpose of these sheaves is to allow us in~\S\ref{subsec: vanishing theorem} to show that the cohomology of
our analytic complexes is concentrated in degrees~$[0,\#I]$; as usual,
this involves a comparison of the toroidal and minimal
compactifications, and we do not know how to carry out this argument
purely in the analytic setting. In~\S\ref{subsec: analytic
  overconvergent sheaves} we construct the corresponding structures in
the analytic world, and we show that an appropriate Hecke operator (a
product of ``$U_p$'' operators at the places dividing~$p$) acts
compactly.

We then recall in~\S\ref{subsec: comparing overconvergent
  analytic and algebraic classes} the analytic BGG resolution
comparing the cohomology with locally analytic coefficients to that
with algebraic coefficients, which is one of the ingredients in our
small slope classicality theorem, which is proved in~\S\ref{subsec: small slope implies classical}, the other
ingredient being a version of the analytic continuation argument
of~\cite{MR2219265}. Finally, in~\S\ref{subsec: ordinary
  cohomology is overconvergent} we apply our results to the complexes
constructed in~\S\ref{subsec: perfect Hida complex}. We are only
able to show that the ordinary cohomology is overconvergent
if~$\#I\le 1$; fortunately, this suffices for the arguments that we
make in~\S\ref{sec:CG}.

\subsection{Sheaves of overconvergent and locally analytic
  modular forms: the formal construction}\label{subsec: formal
  overconvergent sheaves} In this section the base is
$\whalingship_p$, the $p$-adic completion of an algebraic closure of
$\QQ_p$. We
will construct overconvergent versions of our interpolation sheaves
$\Omega^{\kappa_I}$ and develop a finite slope theory. It is necessary
to connect the ordinary theory and the slope $0$ overconvergent
theory, because we are only able to prove a strong classicality
theorem in the overconvergent setting. In the first part of this
section,  we begin by working at a formal level. The reason is that we
need to prove a vanishing theorem (Theorem~ \ref{theo-van}) for the
overconvergent cohomology and we don't know how to prove it without
using formal models.

\subsubsection{Slope decompositions}We very briefly recall the basics
of the theory of slope decompositions for compact operators, which was
introduced in~\cite{ashstevensslopes} and further developed
in~\cite{MR2846490}. Given a vector space~$M$ over~$\whalingship_p$ with a linear
endomorphism~$U$, and a rational number~$h$, an $h$-slope
decomposition of~$M$ with respect to~$U$ is a
decomposition~$M=M^{\le h}\oplus M^{>h}$ into $U$-stable subspaces,
where
\begin{itemize}
\item $M^{\le h}$ is finite-dimensional,
\item all of the eigenvalues~$a$ of~$U$ on~$M^{\le h}$ have $v(a)\le
  h$, and
\item if~$Q$ is a monic polynomial whose roots all have valuation less
  than~$h$, then~$Q(U)$ acts invertibly on~$M^{>h}$.\end{itemize}
If slope decompositions exist, they are unique. If they exist for
all~$h$, then we say that the finite slope part is the union of
the~$M^{\le h}$ for all~$h\in\Q$.

The notion of a slope decomposition can be generalized to the case of
modules over a $\whalingship_p$-Banach algebra~$A$. In particular, it is known
that compact operators on projective $A$-Banach modules admit slope
decompositions locally on~$\Max A$. It is explained in~\cite[\S 2]{MR2846490} and~\cite[\S 13]{pilloniHidacomplexes}
how to generalize this notion to perfect complexes of modules over Banach
algebras. In brief, an endomorphism~$U$ of a perfect complex is said to be compact if it admits
a representative~$\widetilde{U}$ as an endomorphism of a bounded complex~$M^\bullet$ of projective
Banach modules, which is compact in each degree. Then one may consider the product of characteristic power series
of~$\widetilde{U}$ on the individual~$M^i$, and the corresponding
spectral variety for~$\widetilde{U}$ as in~\cite{1997InMat.127..417C}. The
complex~$M^\bullet$ determines a complex of coherent
sheaves~$\cM^\bullet$ over this spectral variety, and one defines the
spectral variety of~$U$ to be the support of the cohomology
sheaves~$H^\bullet(\cM^\bullet)$. One checks that this is independent
of the choice of~$M^\bullet$ and~$\widetilde{U}$. The sheaves
~$H^\bullet(\cM^\bullet)$ over the spectral variety for~$U$ admit
slope decompositions.

\subsubsection{Recollections about formal Banach sheaves}  An
\emph{admissible} $\ocal_{\whalingship_p}$-algebra  is a flat $\ocal_{\whalingship_p}$-algebra which is
a quotient of a converging power series ring $\ocal_{\whalingship_p} \langle X_1,
\cdots, X_n \rangle$ by a finitely generated ideal. In this section we
work with quasi-compact and separated $p$-adic formal schemes over
$\Spf~\ocal_{\whalingship_p}$ which admit an open covering by  formal spectra of
admissible algebras. We call these formal schemes admissible. (In some parts of the literature, an admissible
affine formal scheme~$\Spf A$ is one for which~$A$ is admissible, in
the sense that it is a complete and
separated topological ring, which is linearly topologized and has an
ideal of definition, i.e.\ an open ideal~$I$ such that every
neighbourhood of~$0$ contains some power of~$I$. Our admissible
algebras are a special case of this definition, and we hope that our
terminology will not cause any confusion.)

We recall  some definitions taken from \cite[Defn.\ A.1.1.1]{MR3275848}.  We let $\mathfrak{S}$ be an admissible formal scheme. A formal Banach sheaf over $\mathfrak{S}$ is   a family   $ (\mathfrak{F}_n)_{n \geq 0}$ of   quasi-coherent sheaves such that: 

\begin{enumerate}
\item $\mathfrak{F}_n$ is a sheaf of  $\ocal_{\mathfrak{S}}/p^n$-modules,
\item   $\mathfrak{F}_n$ is flat over $\ocal_{\whalingship_p}/p^n$,
\item  For all  $0 \leq m \leq n$, we have isomorphisms  $\mathfrak{F}_n \otimes_{\ocal_{\whalingship_p}} \ocal_{\whalingship_p}/p^m \simeq  \mathfrak{F}_m$. 
\end{enumerate}

We can associate to $(\mathfrak{F}_n)_n$ a sheaf $\mathfrak{F}$ over $\mathfrak{S}$ equal to the inverse limit $\varprojlim_n \mathfrak{F}_n$ (the maps in the inverse limit are those provided by $(3)$ above). Since $\mathfrak{F}_n = \mathfrak{F} \otimes_{\ocal_{\whalingship_p}} \ocal_{\whalingship_p}/p^n$, the sheaf $\mathfrak{F}$ clearly determines the $(\mathfrak{F}_n)$ and we identify $\mathfrak{F}$ and the family $(\mathfrak{F}_n)$ in the sequel.   
We say that a Banach sheaf $\mathfrak{F}$ is flat if $\mathfrak{F}_n$ is a flat $\ocal_\mathfrak{S}/p^n$-module for all $n$. 

We say that a Banach sheaf $\mathfrak{F}$ is small if there exists a coherent $\ocal_{\mathfrak{S}}/p$-module  $\mathcal{F}$ such that $\mathfrak{F}_1$ is an inductive limit of coherent sheaves $\colim_{j \in \mathbb{N}} \mathfrak{F}_{1,j}$  and the quotients $\mathfrak{F}_{1,j+1}/\mathfrak{F}_{1, j}$ are direct summands of $\mathcal{F}$. 

 We now recall a vanishing result from~\cite{MR3275848}.
 
 \begin{thm}\label{thm-AIP} Let $\mathfrak{S}$ be an admissible formal scheme. Assume that $\mathfrak{S}$ admits a projective map $\mathfrak{S} \rightarrow \mathfrak{S}'$ to an affine admissible formal scheme which induces an isomorphism of the associated analytic adic  spaces over $\Spa ( {\whalingship_p}, \ocal_{\whalingship_p})$. Let $\mathfrak{F}$ be a small Banach sheaf over $\mathfrak{S}$. Let $\mathfrak{U}$ be an affine cover of $\mathfrak{S}$. Then the \v{C}ech complex $$ Cech ( \mathfrak{U}, \mathfrak{F}) \otimes_{\ocal_{\whalingship_p}} \whalingship_p$$ is acyclic in positive degree. 
 \end{thm}
 \begin{proof}
   This is a special case of~\cite[Thm.\ A.1.2.2]{MR3275848}. Indeed,
   the proof of~\cite[Thm.\ A.1.2.2]{MR3275848} is by reducing to this
   case, which is case~(1) of that proof.
 \end{proof}

\subsubsection{Recollections about the Hodge--Tate period map}\label{recollection-Hodgetate}

If $H \rightarrow  \Spec S$ is a finite flat group scheme, we denote
by $H^D$ its Cartier dual and by $\omega_{H^D}$ the conormal sheaf of
$H^D$ along its unit section. This is a coherent $\ocal_S$-module.  We
can view $\omega_{H^D}$  as an \emph{fppf}-sheaf of abelian groups. If
$q: T \rightarrow S$ is an $S$-scheme, we let $\omega_{H^D}(T) =
{\HH}^0(T, q^{*} \omega_{H^D})$.  There is a well-known Hodge--Tate map
$\mathrm{HT}_H: H \rightarrow \omega_{H^D}$ of \emph{fppf}-sheaves of
abelian groups which associates to any $S$-scheme $T$ and point  $x
\in H(T)$ the differential $x^{*} \frac{\mathrm{d} t}{t}$, where we
are (thanks to Cartier duality) viewing~$x$ as a morphism $x: H^D_{T} \rightarrow \mathbb{G}_m\vert_T$ of $T$-group schemes.

Let $ K = K_p K^p$ be a neat compact open subgroup with $K_p = \prod_{v |p}  \mathrm{GSp}_4( \ocal_{F_v})$. Consider the non-compactified Shimura variety $Y_{K}  \rightarrow \Spec \ocal_{\whalingship_p}$.   We denote by $\mathfrak{Y}_{K} \rightarrow \Spf~\ocal_{\whalingship_p}$ the associated $p$-adic formal scheme. We fix a toroidal compactification $Y_{K} \hookrightarrow X_{K}$ and denote by $\mathfrak{X}_K$ the $p$-adic formal scheme associated to $X_K$. 
Let $\mathcal{Y}_{K} \hookrightarrow \mathcal{X}_{K}$ be the associated analytic adic spaces over $\Spa~( {\whalingship_p}, \ocal_{\whalingship_p})$. 

Let $n = (n_v)_{v \in S_p} \in \mathbb{Z}_{\geq 0}^{S_p}$. We let $K(p^n)$ be the compact open subgroup defined by $K(p^n) = K_p(p^n) K^p$ where $K_p(p^n) = \prod_v \mathrm{Ker} \big( \mathrm{GSp}_4( \ocal_{F_v}) \rightarrow \mathrm{GSp}_4( \ocal_{F_v}/p^{n_v})\big)$ is the principal congruence subgroup of level $n$. 

We let $Y_{K(p^n), {\whalingship_p}}  \rightarrow Y_{K} \times_{\Spec \ocal_{\whalingship_p}}
\Spec {\whalingship_p}$ be the  Shimura variety   with  level  $K(p^n)$ structure
over $\Spec {\whalingship_p}$. This map is finite \'etale with Galois group equal
to $\prod_{v |p} \mathrm{GSp}_4 ( \ocal_{F_v}/p^{n_v})$.
Associated to our choice of polyhedral cone decomposition we have a
toroidal compactification  $Y_{K(p^n), {\whalingship_p}} \rightarrow X_{K(p^n),
  {\whalingship_p}}$. We denote by $\mathcal{Y}_{K(p^n)} \hookrightarrow
\mathcal{X}_{K(p^n)}$  the associated analytic spaces over $\Spa (
{\whalingship_p}, \ocal_{{\whalingship_p}})$. The map $\mathcal{X}_{K(p^n)} \rightarrow
\mathcal{X}_{K}$ is finite flat. We denote by $\mathfrak{X}_{K(p^n)}
\rightarrow \mathfrak{X}_{K}$ the normalization of $\mathfrak{X}_{K}$
in $\mathcal{X}_{K(p^n)}$ and by $\mathfrak{Y}_{K(p^n)}$ the
normalization of $\mathfrak{Y}_{K}$ in $\mathcal{Y}_{K(p^n)}$. These
are  admissible formal schemes (see \cite[\S 1.1]{MR3512528}). There is a universal, $\ocal_F$-linear  map $\prod_{v |p} (\ocal_{F_v}/p^{n_v} \ocal_F)^{4} \rightarrow \prod_{v |p} \mathcal{G}_v[p^{n_v}]$  over $\mathfrak{Y}_{K(p^n)}$, which is symplectic up to a similitude factor and is an isomorphism on the associated analytic adic spaces. 

There is a Hodge--Tate period map $ \mathrm{HT}:  \prod_{v |p} \mathcal{G}_v[p^{n_v}] \rightarrow \prod_{v |p} \omega_{\mathcal{G}_v}/p^{n_v} \omega_{\mathcal{G}_v}$ (we are using the quasi-polarization of $\mathcal{G}_v$ to identify $\mathcal{G}_v$ and $\mathcal{G}_v^D$) which we can compose with $\prod_{v |p} (\ocal_{F_v}/p^{n_v} \ocal_{F_v})^{4} \rightarrow \prod_{v |p}  \mathcal{G}_v[p^{n_v}]$ to obtain an $\ocal_F$-linear map of sheaves over $\mathfrak{Y}_{K(p^n)}$ $$ \mathrm{HT}: \prod_{v |p} (\ocal_{F_v}/p^{n_v} \ocal_{F_v})^{4} \rightarrow \prod_{v |p}  \omega_{\mathcal{G}_v}/p^{n_v} \omega_{\mathcal{G}_v}.$$

We claim that this map admits an extension  $$\mathrm{HT}: \prod_{v |p}
(\ocal_{F_v}/p^{n_v} \ocal_{F_v})^{4} \rightarrow \prod_{v |p}
\omega_{\mathcal{G}_v}/p^{n_v} \omega_{\mathcal{G}_v}$$ over
$\mathfrak{X}_{K(p^n)}$. When $F = \Q$, this is the content of
\cite[Prop.\ 1.2]{MR3512528}. For a general $F$ we can use the Koecher
principle of \cite[Thm.\ 8.7]{Lan2016}.

According to a result of Fargues (\cite[Thm.\ 7]{MR2673421}, see
also~\cite[Thm.\ 1.5]{MR3512528}), the cokernel  of the linearization
of $\mathrm{HT}$  is annihilated by $p^{\frac{1}{p-1}}$.  By \cite[\S 1.4]{MR3512528} when $F= \Q$ (and an immediate generalization for general $F$),   there exists an admissible formal scheme $\mathfrak{X}_{K(p^n)}^{mod} \rightarrow \mathfrak{X}_{K(p^n)}$, which is the normalization of a blow-up (the ideal of the blow-up is finitely generated and contains a power of $p$), and a modification $\omega_{\mathcal{G}}^{mod} \subset \omega_{\mathcal{G}}$ such that: 

\begin{enumerate}

\item $\omega_{\mathcal{G}}^{mod}$ is a locally free $\ocal_F \otimes \ocal_{\mathfrak{X}_{K(p^n)}^{mod}}$-module of rank $2$, 
\item $p^{\frac{1}{p-1}} \omega_{\mathcal{G}} \subset \omega_{\mathcal{G}}^{mod} \subset \omega_{\mathcal{G}}$,
\item The Hodge--Tate map $ \mathrm{HT}$ factorizes into a map  $$ \mathrm{HT}: \prod_{v |p}  (\ocal_{F_v}/p^{n_v} \ocal_{F_v})^{4} \rightarrow \prod_{v |p}  \omega_{\mathcal{G}_v}^{mod}/p^{n_v- \frac{1}{p-1}} \omega_{\mathcal{G}_v}^{mod}$$ and the linearized map $$\mathrm{HT} \otimes 1: \prod_{v |p} (\ocal_{F_v}/p^{n_v} \ocal_{F_v})^{4} \otimes \ocal_{\mathfrak{X}_{K(p^n)}^{mod}}  \rightarrow \prod_{v |p} \omega_{\mathcal{G}_v}^{mod}/p^{n_v- \frac{1}{p-1}} \omega_{\mathcal{G}_v}^{mod}$$ is surjective. 
\end{enumerate}
We say a few words about the construction of this formal model. We
first introduce the subsheaf $\omega_{\mathcal{G}}^{mod}$ of
$\omega_{\mathcal{G}}$ generated over  $\mathfrak{X}_{K(p^n)}$ by
$p^{\frac{1}{p-1}} \omega_{\mathcal{G}}$ and local lifts of
$\mathrm{HT}(\prod_{v |p} (\ocal_{F_v}/p^{n_v} \ocal_{F_v})^{4})$ in
$\omega_{\mathcal{G}}$. The sheaf $\omega_{\mathcal{G}}^{mod}$
constructed in this way is not locally free, but becomes locally free
after pulling back to
$\mathfrak{X}_{K(p^n)}^{mod}$ (and we continue to denote this pulled
back sheaf by $\omega_{\mathcal{G}}^{mod}$). We now describe the procedure used to construct $\mathfrak{X}_{K(p^n)}^{mod}$. Zariski locally over $\mathfrak{X}_{K(p^n)}$  we can find a  map 
$$\prod_{v |p} (\ocal_{F_v} \otimes_{\ZZ_p} \ocal_{\mathfrak{X}_{K(p^n)}})^4 \rightarrow \omega_{\mathcal{G}}$$
by considering local lifts of the Hodge-Tate classes in $\omega_{\mathcal{G}}$, and the image of this map is $\omega_{\mathcal{G}}^{mod}$. Zariski locally, we can trivialize $\omega_{\mathcal{G}}$ and we can represent the above map by a $2\times 4$ matrix at each place $v| p$. 
 The formal scheme is obtained by taking the normalization of the blow-up of the ideal which is the product at all places $v$ dividing $p$ of the ideal  locally generated by the $2\times 2$-minors of the   matrix at $v$. 

We denote by $e_{v, 1}, \dots, e_{v, 4}$ the canonical basis of
$\ocal_{F_v}^4$. We let $ \epsilon = (\epsilon_v)_{v |p} \in
\prod_{v |p} ([0, n_v- \frac{1}{p-1}] \cap \Q)$. We define an
admissible formal scheme $\mathfrak{X}_{K(p^n)}( \epsilon) \rightarrow
\mathfrak{X}_{K(p^n)}^{mod} $ (an open subscheme of an admissible
blowup of $\mathfrak{X}_{K(p^n)}^{mod}$) by the conditions that: 

\begin{itemize}
\item  $\mathrm{HT} (e_{v,1}) \in p^{\epsilon_v} \omega_{\mathcal{G}_v}^{mod}/ p^{n_v- \frac{1}{p-1}} \omega_{\mathcal{G}_v}^{mod}$ for all $v |p$,
\item  $\mathrm{HT} (e_{v,2}) \in p^{\epsilon_v} \omega_{\mathcal{G}_v}^{mod}/ p^{n_v- \frac{1}{p-1}} \omega_{\mathcal{G}_v}^{mod}$ for all $v \in I^{c}$. 
\end{itemize}

For all $v \in I^c$, the Hodge--Tate map factorizes into an isomorphism $\mathrm{HT} \otimes 1: \ocal_{\mathfrak{X}_{K(p^n)}}(\epsilon)/p^{\epsilon_v} e_{v,3} \oplus \ocal_{\mathfrak{X}_{K(p^n)}}(\epsilon)/p^{\epsilon_v} e_{v,4} \rightarrow \omega_{\mathcal{G}_v}^{mod}/p^{\epsilon_v} \omega_{\mathcal{G}_v}^{mod}$.

For all $v |p$, we let $\mathrm{Fil}_v^{\can} \subset \omega_{\mathcal{G}_v}^{mod}/p^{\epsilon_v}$ be the sub-module generated by $\mathrm{HT}(e_{v,2})$ and $\mathrm{HT}(e_{v,3})$. 

\begin{lem} $\mathrm{Fil}_v^{\can}$ is a locally free
  $\ocal_{\mathfrak{X}_{K(p^n)}}(\epsilon)/p^{\epsilon_v} $- module of
  rank one, and is locally a direct factor in $\omega_{\mathcal{G}_v}^{mod}/p^{\epsilon_v} \omega_{\mathcal{G}_v}^{mod}$.
\end{lem}

\begin{proof} See  \cite[Lem.\ 12.2.2.1]{pilloniHidacomplexes}.
\end{proof}

We let $\mathrm{Gr}_v^{\can} = \omega_{\mathcal{G}_v}^{mod}/(p^{\epsilon_v}\omega_{\mathcal{G}_v}^{mod} + \mathrm{Fil}_v^{\can})$. Then for all $v|p$,  the Hodge--Tate map induces an isomorphism: 
$$\mathrm{HT} \otimes 1:  (\ocal_{\mathfrak{X}_{K(p^n)}}(\epsilon)/p^{\epsilon_v})e_{v,4} \rightarrow \mathrm{Gr}^{\can}_v. $$

If $v \in I^c$, the Hodge--Tate map also induces an isomorphism $$\mathrm{HT} \otimes 1:  (\ocal_{\mathfrak{X}_{K(p^n)}}(\epsilon)/p^{\epsilon_v})e_{v,3} \rightarrow \mathrm{Fil}^{\can}_v. $$

\subsubsection{Flag varieties}

We denote by $\mathfrak{FL}_n  \rightarrow \mathfrak{X}_{K(p^n)}^{mod}$  the flag  formal scheme which parametrizes locally free direct summands of rank one  (as $\ocal_F \otimes \ocal_{\mathfrak{X}_{K(p^n)}^{mod}}$-modules)  $\mathrm{Fil} \omega_{\mathcal{G}}^{mod}$ in $ \omega_{\mathcal{G}}^{mod}$.  This space decomposes into a product  $\mathfrak{FL}_n = \prod_{v|p } \mathfrak{FL}_{v,n}$ over all places $v$ above $p$.  

Let $ w = (w_v) \in \prod_{v |p}  [0, \epsilon_v] \cap \QQ$.  We let $\mathfrak{FL}_{n, \epsilon , w} \rightarrow \mathfrak{X}_{K(p^n)}(\epsilon)$ be the moduli space of locally free direct summands of rank one  $\mathrm{Fil} \omega_{\mathcal{G}_v}^{mod} \subset \omega_{\mathcal{G}_v}^{mod}$ such that $\mathrm{Fil} \omega_{\mathcal{G}_v}^{mod} = \mathrm{Fil}^{\can}_v \mod p^{w_v}$. 

We let $w' = (w'_v) \in \prod_{v |p} [0, w_v] \cap \Q$. We let  $\mathfrak{FL}^+_{n, \epsilon , w,w'} \rightarrow \mathfrak{FL}_{n, \epsilon , w}$ be the moduli space parametrizing: 

\begin{enumerate}

\item For all $v |p$ a basis $\rho_v: \ocal_{\mathfrak{FL}^+_{n, \epsilon , w,w'}} \rightarrow \omega_{\mathcal{G}^{mod}_v} / \mathrm{Fil} \omega_{\mathcal{G}^{mod}_v}$  such that $\rho_v(1) = \mathrm{HT}(e_{v,4}) \mod p^{w'_v}$, 

\item For all $v \in I^c$, a basis $\nu_v: \ocal_{\mathfrak{FL}^+_{n, \epsilon , w,w'}} \rightarrow \mathrm{Fil} \omega_{\mathcal{G}^{mod}_v}$ such that $\nu_v(1) = \mathrm{HT}(e_{v,3}) \mod p^{w'_v}$.
\end{enumerate}

\subsubsection{Some groups}

The group $\prod_{v |p} \mathrm{GSp}_4( \ocal_{F_v}/p^{n_v})$ acts
on $\mathfrak{X}_{K(p^n)}$ and $\mathfrak{X}_{K(p^n)}^{mod}$. The
parabolic subgroup $\prod_{v\in I^c} \mathrm{B}(\ocal_{F_v}/p^{n_v})
\prod_{v \in I} \Kli(\ocal_{F_v}/p^{n_v})$ acts on
$\mathfrak{X}_{K(p^n)}( \epsilon)$. 

Let us denote by $\mathfrak{X}_{K(I,p^n)}(\epsilon)$ the quotient of $\mathfrak{X}_{K(p^n)}(\epsilon)$ by the action of this finite group. This is an admissible  formal scheme.

We have maps $\mathrm{B}(\ocal_{F_v}/p^{n_v}) \rightarrow ((\ocal_{F_v}/p^{n_v})^\times)^2$ provided by the last two diagonal entries and 
$\Kli(\ocal_{F_v}/p^{n_v}) \rightarrow (\ocal_{F_v}/p^{n_v})^\times$ provided by the last diagonal entry. 

We denote by $\mathfrak{T}_{w'}^0$ the formal group defined by 
$$\mathfrak{T}_{w'}^0(R) = \prod_{v\in I} (1 + p^{w'_v} R) \prod_{v \in I^c} (1 + p^{w'_v} R)^2$$
  for any admissible $\ocal_{\whalingship_p}$-algebra $R$. 

We denote by  $\mathfrak{T}_{w'}$ the group 
$$\mathfrak{T}_{w'}(R) = \prod_{v\in I}  \ocal_{F_v}^\times(1 + p^{w'_v} R) \prod_{v \in I^c} (\ocal_{F_v}^\times(1 + p^{w'_v} R))^2$$  for any admissible $\ocal_{\whalingship_p}$-algebra $R$. 

Finally, we denote by $\mathfrak{T}_{n,w'}$ the fibre product  $$ \mathfrak{T}_{w'} \times_{ \mathfrak{T}_{w'}/\mathfrak{T}_{w'}^0} \prod_{v\in I^c} \mathrm{B}(\ocal_{F_v}/p^{n_v}) \prod_{v \in I} \Kli(\ocal_{F_v}/p^{n_v}).$$

\subsubsection{Torsors}

The map $\mathfrak{FL}^+_{n ,\epsilon, w, w'} \rightarrow \mathfrak{FL}_{n, \epsilon, w}$ is a $\mathfrak{T}_{w'}^0$-torsor. The group $\mathfrak{T}_{w'}^0$ acts on $\rho_v$ and $\nu_v$. This action extends to an action of $\mathfrak{T}_{n,w'}$ on $\mathfrak{FL}^+_{n ,\epsilon, w, w'} \rightarrow \mathfrak{X}_{K(p^n)}(\epsilon)$, compatible with the action of $\prod_{v\in I^c} \mathrm{B}(\ocal_{F_v}/p^{n_v}) \prod_{v \in I} \Kli(\ocal_{F_v}/p^{n_v})$ on $\mathfrak{X}_{K(p^n)}(\epsilon)$. 

\subsubsection{Formal Banach sheaves}

Let $A$ be a normal admissible $\ocal_{\whalingship_p}$-algebra. Let  $\kappa_A:
\prod_{v\in I} \ocal_{F_v}^\times \prod_{v \in I^c}
(\ocal_{F_v}^\times)^2 \rightarrow A^\times$ be a character, which we
assume is $w'$-analytic, in the sense  that it extends to a pairing $\mathfrak{T}_{w'} \times \Spf~A \rightarrow \mathbb{G}_m$. 

We denote by $\pi_1: \mathfrak{FL}^+_{n ,\epsilon, w, w'} \rightarrow \mathfrak{FL}_{n, \epsilon, w}$ the projection. We can define an invertible sheaf of $\ocal_{ \mathfrak{FL}_{n, \epsilon,w}} \hat{\otimes} A$-modules, $$\mathfrak{L}^{\kappa_A} = ((\pi_1)_{{*}} \ocal_{\mathfrak{FL}^+_{n ,\epsilon, w, w'}} \hat{\otimes} A)^{\mathfrak{T}_{w'}^0}$$ where the invariants are taken for the diagonal action.

We let $\pi_2: \mathfrak{FL}_{n, \epsilon, w} \rightarrow
\mathfrak{X}_{K(p^n)}(\epsilon)$. This is an affine map.  We define a
formal Banach sheaf $\mathfrak{G}^{\kappa_A,w} = (\pi_2)_{*}
\mathfrak{L}^{\kappa_A}$ over $\mathfrak{X}_{K(p^n)}(\epsilon)$; this
is independent of the choice of~$w'$, as is easily seen from the construction.

Finally, we let $\pi_3: \mathfrak{X}_{K(p^n)}(\epsilon)\to \mathfrak{X}_{K(I,p^n)}(\epsilon)$, and we   define $\mathfrak{F}^{\kappa_A,w} = ((\pi_3)_{*} \mathfrak{G}^{\kappa_A, w})^{\mathfrak{T}_{n,w}}$. This is a formal Banach sheaf over $\mathfrak{X}_{K(I,p^n)}(\epsilon)$. 

\subsubsection{ Some properties}

For each $v\in I$ we choose an element $i_{v} \in \{2,3\}$. Let $\mathfrak{X}_{K(p^n)}( \epsilon, (i_v))$ be the open subset of $\mathfrak{X}_{K(p^n)}(\epsilon)$ where $\mathrm{Fil}_{v}^{\can}$ is generated by $\mathrm{HT}(e_{i_v})$ for all $v \in I$. 

\begin{lem}\label{lem-some-prop1} The  quasi-coherent sheaf $\mathfrak{G}^{\kappa_A,w}/p^{\inf_v w_v}$ restricted to $\mathfrak{X}_{K(p^n)}( \epsilon, (i_v))$ is an inductive limit of coherent sheaves which are extensions of the  sheaf $\ocal_{\mathfrak{X}_{K(p^n)}(\epsilon, (i_v))}/p^{\inf_{v} w_v}$. 
\end{lem} 
\begin{proof} This can be proved in the same way as~\cite[Lem.\ 8.1.6.2]{MR3275848}.
\end{proof} 

\begin{lem}\label{lem-some-prop2} The  quasi-coherent sheaf $\mathfrak{G}^{\kappa_A,w}/p$ is a flat sheaf of $\ocal_{\mathfrak{X}_{K(p^n)}(\epsilon)}/p$-modules.
\end{lem}

\begin{proof} See~ \cite[Lem.\
  12.6.2.1]{pilloniHidacomplexes}.\end{proof}

\subsection{Vanishing theorem}\label{subsec: vanishing theorem}
\subsubsection{The minimal compactification}The main result of this
subsection is Theorem~\ref{theo-van}. As in the proof of
Theorem~\ref{thm-vanishing}, we will use the minimal compactification,
and in particular the facts that the pushforward of our sheaves to the
minimal compactification are supported on open subsets that admit an
explicit affine cover, and that  the higher derived
pushforwards from the toroidal to minimal compactifications of the
cuspidal cohomology vanish.

We denote by $\mathfrak{X}^{*}_{K}$ the minimal compactification of $\mathfrak{Y}_{K}$. There is a natural map $\mathfrak{X}_K \rightarrow \mathfrak{X}^{*}_K$. The invertible sheaf $\det \omega_{\mathcal{G}}$ over  $\mathfrak{X}_K$ descends to an invertible sheaf still denoted by $\det \omega_{\mathcal{G}}$ over $\mathfrak{X}^{*}_K$.  Let 
$n = (n_v) \in \ZZ_{\geq 0}^{S_p}$. In this subsection we consider only
the case that $n_v$ is independent of $v$, and accordingly we will
write~$n$ for~$n_v$. We let
$\mathfrak{X}^{*}_{K(p^n)}$ be the Stein factorization of the
morphism: $\mathfrak{X}_{K(p^n)} \rightarrow
\mathfrak{X}_K^{{*}}$. This is a normal admissible formal scheme. In
\cite[Cor.\ 1.4]{MR3512528} it is proved that the determinant of the Hodge--Tate map on $\mathfrak{X}_{K(p^n)}$: 
$$\Lambda^2 \mathrm{HT}:  \bigotimes_{v|p}  \Lambda^2( \ocal_{F_v}/p^{n} \ocal_{F_v}) ^4 \rightarrow \bigotimes_{v |p}  \det \omega_{\mathcal{G}_v}/p^{n}$$is the pull back of a map  denoted the same way:
 $$\Lambda^2 \mathrm{HT}: \bigotimes_{v |p} \Lambda^2 (\ocal_{F_v}/p^{n} \ocal_{F_v}) ^4 \rightarrow  \bigotimes_{v |p} \det \omega_{\mathcal{G}_v}/p^{n}$$
 which is defined over $\mathfrak{X}^{*}_{K(p^n)}$.

\begin{rem} Literally, the determinant of the Hodge--Tate map is a map: $$\Lambda^{2[F:\QQ]} (\ocal_{F}/p^n \ocal_F)^4 \rightarrow \det \omega_{\mathcal{G}}/p^n.$$ But using the action of $\ocal_F$, it is easy to see that it factors through the direct factor $\bigotimes_{v|p}  \Lambda^2( \ocal_{F_v}/p^{n} \ocal_{F_v}) ^4$. 
\end{rem} 

By \cite[\S 1.4]{MR3512528} (for $F= \Q$, and the same construction for general $F$), there is a normal admissible  formal
scheme $\mathfrak{X}_{K(p^n)}^{{*}-mod}  \rightarrow
\mathfrak{X}_{K(p^n)}^{*}$ which is the normalization of a blow up and
carries a locally free modification $\det \omega_{\mathcal{G}}^{mod} \subset \det
\omega_{\mathcal{G}}$  such that:

\begin{enumerate}
\item $p^{ \frac{2[F:\Q]}{p-1}} \det \omega_{\mathcal{G}} \subset \det \omega_{\mathcal{G}}^{mod}  \subset \det \omega_{\mathcal{G}}$.
\item The Hodge--Tate map factorizes into a surjective map: 
$$ \bigotimes_{v |p}  \Lambda^2 ( \ocal_{F_v}/p^{n} \ocal_{F_v} )^4 \otimes \ocal_{ \mathfrak{X}_{K(p^n)}^{{*}-mod}}  \rightarrow   \det \omega_{\mathcal{G}}^{mod}/p^{n-\frac{2[F:\Q]}{p-1}}.$$
\end{enumerate}

The construction  $\mathfrak{X}_{K(p^n)}^{{*}-mod}$ follows a similar procedure as the construction of $\mathfrak{X}_{K(p^n)}^{mod}$ explained in section \ref{recollection-Hodgetate}: one can lift locally the map $\Lambda^2 \mathrm{HT}$ to a map   $\bigotimes_{v|p}  \Lambda^2( \ocal_{F_v}^4) \otimes \ocal_{\mathfrak{X}_{K(p^n)}^{*}} \rightarrow \bigotimes_{v |p}  \det \omega_{\mathcal{G}_v}$ and consider the normalization of the blow up of the ideal which is locally  the product at all places $v$ of the ideals generated by the coefficients of the above map at the place $v$.  By the universal properties of blow-ups and normalizations, there is a map $\mathfrak{X}_{K(p^n)}^{mod} \rightarrow \mathfrak{X}_{K(p^n)}^{{*}-mod}$.  
Let $\epsilon = (\epsilon_v) \in \prod_{v |p} ([0, n-
\frac{2[F:\QQ]}{p-1}] \cap \QQ) $. We denote by
$\mathfrak{X}_{K(p^n)}^{*} (\epsilon) \rightarrow
\mathfrak{X}_{K(p^n)}^{{*}-mod}$ the formal scheme defined by the
condition:

\begin{itemize}
\item $ \mathrm{HT}( e_{v,1}) \wedge \mathrm{HT}( e_{v, i}) \otimes_{v' |p, v' \neq v} \mathrm{HT}( e_{v',j_{v'}}) \wedge \mathrm{HT}( e_{v', k_{v'}}) \in p^{\epsilon_v} \det \omega_{\mathcal{G}}^{mod}/p^{n- \frac{2[F:\Q]}{p-1}} $ for all $v |p$ and $1 \leq i, j_{v'}, k_{v'} \leq 4$,
\item $ \mathrm{HT}( e_{v,2}) \wedge \mathrm{HT}( e_{v, i}) \otimes_{v' |p, v' \neq v} \mathrm{HT}( e_{v',j_{v'}}) \wedge \mathrm{HT}( e_{v', k_{v'}}) \in p^{\epsilon_v} \det \omega_{\mathcal{G}}^{mod}/p^{n- \frac{2[F:\Q]}{p-1}} $ for all $v \in I^c$ and $1 \leq i,j_{v'},k_{v'} \leq 4$.
\end{itemize}

There is a Cartesian diagram (see the proof of~\cite[Lem.\ 12.9.1.1]{pilloniHidacomplexes}): \begin{eqnarray*} 
\xymatrix{ \mathfrak{X}_{K(p^n)}( \epsilon) \ar[r] \ar[d] & \mathfrak{X}_{K(p^n)}^{mod} \ar[d] \\
\mathfrak{X}_{K(p^n)}^{{*}} ( \epsilon) \ar[r] & \mathfrak{X}_{K(p^n)}^{{*}-mod}}
\end{eqnarray*}
By the proof of~ \cite[Thm.\ 4.3.1, pp 1029-30]{scholze-torsion} (see
also \cite[Thm.\ 1.16]{MR3512528}), there is an integer $N$ such that for all $n \geq N$, there is  a normal  admissible formal scheme $\mathfrak{X}_{K(p^n)}^{{*}-\mathrm{HT}}$ and a projective map $\mathfrak{X}_{K(p^n)}^{{*}-mod} \rightarrow \mathfrak{X}_{K(p^n)}^{{*}-\mathrm{HT}}$   which is an isomorphism on the associated analytic spaces and satisfies: 
\begin{enumerate}
\item The invertible sheaf $\det \omega_{\mathcal{G}}^{mod} $ descends
  to an ample invertible sheaf on
  $\mathfrak{X}_{K(p^n)}^{{*}-\mathrm{HT}}$.\item For all  rational numbers $\epsilon >0 $, there
  is $n( \epsilon) \geq N$ such that if $n \geq n(\epsilon)$, then
  there are sections $s_{(i_v,j_v)_{v|p}} \in {\HH}^0( \mathfrak{X}_{K(p^n)}^{{*}-\mathrm{HT}}, \det \omega_{\mathcal{G}}^{mod})$ satisfying $s_{(i_v,j_v)} = \otimes_{v |p} \mathrm{HT}( e_{v,i_v}) \wedge \mathrm{HT}( e_{v,j_v})$ in $\det \omega_{\mathcal{G}}^{mod}/p^{\epsilon}$ for all $1 \leq i_v , j_v \leq 4$. 
\end{enumerate} 

Let $\epsilon = ( \epsilon_v) \in ( \QQ_{>0})^{S_p}$.  Let $n \geq \sup_v n( \epsilon_v)$.  We define a formal scheme $ \mathfrak{X}_{K(p^n)}^{{*}-\mathrm{HT}}( \epsilon)  \rightarrow \mathfrak{X}^{{*}-\mathrm{HT}}_{K(p^n)}$ by the condition: 

\begin{itemize}
\item  for all $v |p$, for all $(i_{v'}, j_{v'})_{v'|p} \in (\{ 1,2,3,4\} \times \{1,2,3,4\})^{S_p}$,  such that  $i_v = 1$, we have   $s_{(i_{v'}, j_{v'})} \in p^{\epsilon_{v}} \det \omega_{\mathcal{G}}^{mod} $,
\item for all $v \in I^c$, for all $(i_{v'}, j_{v'})_{v'|p}  \in (\{ 1,2,3,4\} \times \{1,2,3,4\})^{S_p}$,  such that  $i_v = 2$, we have   $s_{(i_{v'}, j_{v'})} \in p^{\epsilon_v} \det \omega_{\mathcal{G}}^{mod}$. 
\end{itemize}

We have a Cartesian diagram 

\begin{eqnarray*} 
\xymatrix{ \mathfrak{X}_{K(p^n)}^{{*}}( \epsilon) \ar[r] \ar[d] & \mathfrak{X}_{K(p^n)}^{{*}-mod} \ar[d] \\
\mathfrak{X}_{K(p^n)}^{{*}- \mathrm{HT}} ( \epsilon) \ar[r] & \mathfrak{X}_{K(p^n)}^{{*}-\mathrm{HT}}}
\end{eqnarray*}
where both vertical maps are projective maps and induce isomorphisms on the associated analytic generic fibres.

For all $v \in I$, let $i_v \in \{2,3\}$.  We define an open subspace $\mathfrak{X}^{{*}-\mathrm{HT}}_{K(p^n)}( \epsilon, (i_v))$ of  $\mathfrak{X}^{{*}-\mathrm{HT}}_{K(p^n)}( \epsilon)$ by the condition that $s_{(i_v,4)_{v \in I}, (3,4)_{v \in I^c}} \neq 0$. 

We similarly define an open subspace $\mathfrak{X}^{{*}}_{K(p^n)}( \epsilon, (i_v))$ of  $\mathfrak{X}^{{*}}_{K(p^n)}( \epsilon)$ by the condition that $\otimes_{v \in I}  \mathrm{HT}(e_{v, i_v}) \wedge \mathrm{HT}(e_{v,4}) \otimes_{v \in I^c}  \mathrm{HT}(e_{v, 3}) \wedge \mathrm{HT}(e_{v,4})   \neq 0$ for all $v \in I$.

\begin{lem}\label{lem-usingAIP} We have a projective map  $\mathfrak{X}^{{*}}_{K(p^n)}( \epsilon, (i_v)) \rightarrow \mathfrak{X}^{{*}-\mathrm{HT}}_{K(p^n)}( \epsilon, (i_v))$ which is an isomorphism on the associated analytic adic spaces. Moreover, $\mathfrak{X}^{{*}-\mathrm{HT}}_{K(p^n)}( \epsilon, (i_v))$ is an affine formal scheme. 
\end{lem}

\begin{proof} The first point is clear. The second point follows from
  the ampleness of $\det \omega_{\mathcal{G}}^{mod}$ on
  $\mathfrak{X}^{{*}-\mathrm{HT}}_{K(p^n)}( \epsilon)$ and the fact
  that $\mathfrak{X}^{{*}-\mathrm{HT}}_{K(p^n)}( \epsilon, (i_v))$ is
  the open subscheme defined by the non-vanishing of a section of an ample sheaf. 
\end{proof}

\subsubsection{Vanishing}

 We have a map $\pi: \mathfrak{X}_{K(p^n)}(\epsilon) \rightarrow \mathfrak{X}_{K(p^n)}^{*}(\epsilon)$. We denote as usual by $D$ the boundary divisor. 
 
 \begin{prop}\label{prop-van} We have $\mathrm{R}^i\pi_{*} \ocal_{{\mathfrak{X}_{K(p^n)}(\epsilon)}}(-D) = 0$ for all $i >0$. 
 \end{prop}
 \begin{proof} This can be proved in exactly the same way as~
   \cite[Prop. 12.9.2.1]{pilloniHidacomplexes} (which is the case~$F=\Q$).\end{proof}
 
 \begin{thm}\label{theo-van} Let $\epsilon = (\epsilon_v) \in \QQ_{>0}^{S_p}$. Let $n = (n_v)$ with $\inf_{v} n_v \geq \sup_v n(\epsilon_v)$.  
 The complex $$\mathrm{R}\Gamma( \mathfrak{X}_{K(p^n)}(\epsilon), \mathfrak{G}^{\kappa_A,w} \otimes (\det \omega_{\mathcal{G}}^{mod})^2(-D))[1/p]$$ has cohomology concentrated in degrees $[0,\#I]$.
\end{thm}

\begin{proof} Consider the hypercube $[2,3]^{I}$. We can associate to it a category  denoted by $\mathcal{C}$. Its objects are the faces $\sigma$ of the hypercube. By definition, a face is a product $\prod_{v \in I} \lambda_v$  where for each $v \in I$, $\lambda_v \in \{ 2,3, [2,3]\}$ (our convention is that faces are closed). There is a map $\sigma' \rightarrow \sigma$ between faces if $\sigma'$ is included in~$\sigma$. 

We now define a functor $\mathcal{C}^{op} \rightarrow Op(
\mathfrak{X}_{K(p^n)}^{{*}-mod}(\epsilon))$, where the target is the category of  open subsets of
$\mathfrak{X}_{K(p^n)}^{{*}-mod}(\epsilon)$ (whose morphisms are 
open immersions). It sends a face $\sigma$ to $U_\sigma$, the intersection of  all the formal schemes $\mathfrak{X}^{{*}-mod}_{K(p^n)}( \epsilon, (i_v))$ for  $(i_v)_{v \in I} \in \sigma$ (we recall that $i_v \in \{2,3\}$).

Write $f:\mathfrak{X}_{K(p^n)}(\epsilon) \rightarrow \mathfrak{X}^{{*}-mod}_{K(p^n)}(\epsilon)$  for the  map defined above. For all $\sigma \in \mathcal{C}$, consider the following Cartesian diagram: 
\begin{eqnarray*}
\xymatrix{  \mathfrak{X}_{K(p^n)}(\epsilon)_{\sigma} \ar[d]^{f_\sigma} \ar[r]^{i}  & \mathfrak{X}_{K(p^n)}(\epsilon) \ar[d]^{f} \\
U_{\sigma} \ar[r]^{j} & \mathfrak{X}^{{*}-mod}_{K(p^n)}(\epsilon)}
\end{eqnarray*}
where~$j$ is the natural open immersion.

Let $Sh (  \mathfrak{X}_{K(p^n)}(\epsilon))$ be the category of sheaves on $\mathfrak{X}_{K(p^n)}(\epsilon)$. We define a functor $\mathcal{C} \rightarrow Sh (  \mathfrak{X}_{K(p^n)}(\epsilon))$ which sends $\sigma$ to the sheaf \[ \mathfrak{G}_\sigma = i_* i^* (\mathfrak{G}^{\kappa_A,w}\otimes (\det \omega_{\mathcal{G}}^{mod})^2(-D)[1/p]).\]

We deduce from Lemmas \ref{lem-some-prop1} and~ \ref{lem-some-prop2},
together with Proposition~ \ref{prop-van}, that the sheaf $f_* \mathfrak{G}_\sigma$ is a small formal Banach sheaf and that $\mathrm{R}^i f_\star \mathfrak{G}_\sigma = 0$ for all $i >0$.  It follows from Theorem \ref{thm-AIP} (which applies because of Lemma \ref{lem-usingAIP}) that $f_* \mathfrak{G}_\sigma$ is acyclic. 

We deduce that the cohomology $\mathrm{R}\Gamma( \mathfrak{X}_{K(p^n)}(\epsilon), \mathfrak{G}^{\kappa_A,w} \otimes (\det \omega_{\mathcal{G}}^{mod})^2(-D)[1/p])$ is represented by the  complex  $ C^\bullet$ concentrated in degree $0$ to $\# I$, whose $i$th term is $\oplus_{ \sigma, \dim \sigma = i} \HH^0(  \mathfrak{X}_{K(p^n)}(\epsilon), \mathfrak{G}_\sigma)$ and whose differentials are alternating sums of the restriction maps $\HH^0(  \mathfrak{X}_{K(p^n)}(\epsilon), \mathfrak{G}_{\sigma'}) \rightarrow \HH^0(  \mathfrak{X}_{K(p^n)}(\epsilon), \mathfrak{G}_{\sigma})$ for $\sigma' \subset \sigma$, $\dim \sigma' = \dim \sigma -1$. 
 \end{proof}

\subsection{Sheaves of overconvergent and locally analytic modular
  forms: the analytic construction}\label{subsec: analytic
  overconvergent sheaves} We now translate our previous formal constructions to the  analytic setting, which is well adapted for the spectral theory. 

\subsubsection{Analytic Hilbert--Siegel varieties} This section is
parallel to \cite[\S 12.7]{pilloniHidacomplexes}.  We let $\mathcal{X}_{K(p^n)}$ be the generic fibre of $\mathfrak{X}_{K(p^n)}$. We write $\mathcal{X}_K$ for the generic fibre of $\mathfrak{X}_K$.  
 We let $\mathcal{X}_{K(p^n)} ( \epsilon) \subset
 \mathcal{X}_{K(p^n)}$ be the generic fibre of
 $\mathfrak{X}_{K(p^n)}(\epsilon)$. We let
 $\mathcal{X}_{K(I,p^n)}(\epsilon)$ be the  generic fibre of
 $\mathfrak{X}_{K(I,p^n)}(\epsilon)$. We now give a modular interpretation
 of this last space. Let $A$ be the universal semi-abelian scheme and $\cG$ be its $p$-divisible group. Let $\omega_{\mathcal{G}}$ be the conormal sheaf of $A$ at the origin and let $\omega_{\mathcal{G}}^+ \subset \omega_{\mathcal{G}}$ be the subsheaf of integral differentials (we use the slight abuse of notation to write $\omega_{\cG}$ instead of $\omega_A$). These are  sheaves over $\mathcal{X}_K$ on the analytic site. 

We let $\omega_{\mathcal{G}}^{mod,+}$ be the subsheaf of $\omega_{\mathcal{G}}^{+}$ generated by the image of the Hodge--Tate map. This is an \'etale sheaf over $\mathcal{X}_K$. 

The fibres  of the map   $\mathcal{X}_{K(I,p^n)}(\epsilon) \rightarrow
\mathcal{X}_{K}$ parametrize:\begin{itemize}
\item For all $v \in I$, a subgroup $H_{v,n_v} \subset  \mathcal{G}_v[p^{n_v}]$ which is locally for the \'etale topology isomorphic to $\ZZ/p^{n_v}\ZZ$ and is locally for the \'etale topology generated by an element $e_{v,1}$ which satisfies $\mathrm{HT}(e_{v,1})  =0$ in $\omega_{\mathcal{G}_v}^{mod,+} / p^{\epsilon_v}$. 
\item For all $v \in I^c$, totally isotropic subgroups $H_{v,n_v} \subset L_{v,n_v} \subset  \mathcal{G}_v[p^{n_v}]$ such that $H_{v,n_v}$ is locally for the \'etale topology isomorphic to $\ZZ/p^{n_v}\ZZ$, $ L_{v,n_v}$ is locally for the \'etale topology isomorphic to $(\ZZ/p^{n_v}\ZZ)^2$, and is locally for the \'etale topology generated by  elements $e_{v,1}$ and $e_{v,2}$ which satisfy $\mathrm{HT}(e_{v,1}) = \mathrm{HT}(e_{v,2}) =0$ in $\omega_{\mathcal{G}_v}^{mod,+} / p^{\epsilon_v}$.  
\end{itemize}

We can define for all $v |p$ an \'etale sheaf $\mathrm{Fil}_v^{\can} = \mathrm{Im} ( \mathrm{HT}: H_{v,n_v}^\bot \otimes \ocal_{\mathcal{X}_{K(I,p^n)}(\epsilon)}^+ \rightarrow \omega_{\mathcal{G}_v}^{mod,+}/p^{\epsilon_v})$. This is a \'{e}tale locally free sheaf of $\ocal_{\mathcal{X}_{K(I,p^n)}(\epsilon)}^+/p^{\epsilon_v}$-modules of rank $1$.  We let $\mathrm{Gr}_v^{\can} = \omega_{\mathcal{G}_v}^{mod,+}/(p^{\epsilon_v} + \mathrm{Fil}_v^{\can})$. 

We have isomorphisms deduced from the Hodge--Tate map: 
\begin{itemize}
\item $ \mathrm{HT}: H_{v,n_v}^D \otimes \ocal_{\mathcal{X}_{K(I,p^n)}(\epsilon)}^+/p^{\epsilon_v} \rightarrow \mathrm{Gr}_v^{\can} $ for all $v |p$,
\item $ \mathrm{HT}: (L_{v,n_v}/H_{v,n_v})^D \otimes \ocal_{\mathcal{X}_{K(I,p^n)}(\epsilon)}^+/p^{\epsilon_v} \rightarrow \mathrm{Fil}_v^{\can}$ for all $v \in I^c$.
\end{itemize}

We let $\mathcal{FL}_{K(I,p^n),\epsilon, w} \rightarrow \mathcal{X}_{K(I,p^n)}(\epsilon)$ be the moduli space of flags $\mathrm{Fil} \omega_{\mathcal{G}} \subset \omega_{\mathcal{G}}$ satisfying $\mathrm{Fil} \omega_{\mathcal{G}_v}  \cap \omega_{\mathcal{G}_v}^{mod,+} / p^{w_v} = \mathrm{Fil}_v^{\can}/p^{w_v}$. 

We let $\mathcal{X}_{K(I,p^n)^+}(\epsilon) \rightarrow \mathcal{X}_{K(I,p^n)}(\epsilon)$ be the \'etale cover parametrizing trivializations:

\begin{itemize}
\item $\ZZ/p^{n_v} \ZZ e_{v,4} \stackrel{\sim}\rightarrow H_{v,n_v}^D$ for all $v |p$,
\item $\ZZ/p^{n_v}\ZZ e_{v,3} \stackrel{\sim}\rightarrow (L_{v, n_v}/H_{v,n_v})^{D}$ for all $v \in I^c$. 
\end{itemize}

We let $\mathcal{FL}_{K(I,p^n),\epsilon, w,w'}^+ \rightarrow \mathcal{FL}_{K(I,p^n),\epsilon, w}\times_{\mathcal{X}_{K(I,p^n)}(\epsilon)}\mathcal{X}_{K(I,p^n)^+}(\epsilon)$ be the moduli space of trivializations of:
\begin{itemize}\item for all~$v|p$, $\rho_v: \ocal_{\mathcal{FL}_{K(I,p^n),\epsilon, w,w'}^+}
  \rightarrow \mathrm{\mathrm{Gr}} \omega_{\mathcal{G}_v}=\omega_{\cG_v}/\Fil\omega_{\cG_v}$ such that $\rho_v(1) = \mathrm{HT}(e_{v,4})$ modulo $p^{w'_v}$. 
\item  for all~$v\in I^c$, $\nu_v: \ocal_{\mathcal{FL}_{K(I,p^n),\epsilon, w,w'}^+} \rightarrow \mathrm{Fil} \omega_{\mathcal{G}_v}$ such that $\nu_v(1) = \mathrm{HT}(e_{v,3})$ modulo $p^{w'_v}$.  

\end{itemize}

We can connect these definitions with the constructions of the previous sections. Let $\mathcal{FL}_{n, \epsilon, w} \rightarrow \mathcal{X}_{K(p^n)}(\epsilon)$ be the analytic space associated to $\mathfrak{FL}_{n, \epsilon, w}$. Let $\mathcal{FL}_{n, \epsilon, w, w'}^+$ be the analytic space associated to $\mathfrak{FL}_{n, \epsilon, w,w'}^+$. 

\begin{lem} We have $$\mathcal{FL}_{n, \epsilon, w} = \mathcal{FL}_{K(I,p^n), \epsilon, w} \times_{\mathcal{X}_{K(I,p^n)}( \epsilon)} \mathcal{X}_{K(p^n)}(\epsilon)$$ and $$\mathcal{FL}_{n, \epsilon, w, w'}^+ = \mathcal{FL}_{K(I,p^n), \epsilon, w, w'}^+ \times_{\mathcal{X}_{K(I,p^n)^+}( \epsilon)} \mathcal{X}_{K(p^n)}(\epsilon).$$
\end{lem}

\begin{proof} This follows from the definitions.
\end{proof}

\subsubsection{Banach sheaves}  We let $\mathcal{L}^{\kappa_A}$ be the
invertible sheaf over $\mathcal{FL}_{n, \epsilon,w} \times \Spa
(A[1/p], A)$ associated to $\mathfrak{L}^{\kappa_A}$.  We let
$\mathcal{G}^{\kappa_A,w}$ be the Banach sheaf over
$\mathcal{X}_{K(p^n)}(\epsilon)$ associated to
$\mathfrak{G}^{\kappa_A,w}$. We let $\mathcal{F}^{\kappa_A,w}$ be the
Banach sheaf over $\mathcal{X}_{K(I,p^n)}(\epsilon)$ attached to
$\mathfrak{F}^{\kappa_A,w}$. A direct definition of
$\mathcal{F}^{\kappa_A,w}$ is the following. Let $\pi:
\mathcal{FL}^+_{K(I,p^n), \epsilon, w,w'} \rightarrow
\mathcal{X}_{K(I,p^n)}(\epsilon)$ be the affine projection. Let
$\mathcal{T}_{w'}$ be the generic fibre of $\mathfrak{T}_{w'}$. This
group acts naturally on  $\mathcal{FL}^+_{K(I,p^n), \epsilon, w,w'}$,
trivially on $ \mathcal{X}_{K(I,p^n)}(\epsilon)$, and  the morphism
$\pi$ is equivariant for the action. It follows from the definitions
that $$\mathcal{F}^{\kappa_A, w} = (\pi_{*}
\ocal_{\mathcal{FL}^+_{K(I,p^n), \epsilon, w,w'}} \hat{\otimes}
A)^{\mathcal{T}_{w'}}$$ where the invariants are for the diagonal
action (with the action on the second factor being via~$\kappa_A$).

\subsubsection{Locally analytic overconvergent cohomology}

We define the $n, \epsilon$-convergent, cuspidal $w$-analytic cohomology of
weight parametrized by $A$ to be: $$C_{cusp}(n, \epsilon, w, \kappa_A\otimes (2,2)_{v |p}):= \mathrm{R} \Gamma ( \mathcal{X}_{K(I,p^n)}(\epsilon), \mathcal{F}^{\kappa_A,w} \otimes (\det \omega_{\mathcal{G}})^2(-D)).$$
For $\epsilon' \geq \epsilon$, $n' \geq n$, $w' \geq w$, we have maps: 
$C_{cusp}(n, \epsilon, w,  \kappa_A\otimes (2,2)_{v |p}) \rightarrow C_{cusp}(n', \epsilon', w',  \kappa_A\otimes (2,2)_{v |p})$.

Passing to the limit over $n, \epsilon, w$, we  define the $i$th
cohomology groups of cuspidal, overconvergent, locally analytic
cohomology of weight parametrized by $A$:
$$ {\HH}_{cusp} ^i( \dag, \kappa_A\otimes (2,2)_{v |p}) = \colim {\HH}^i ( C_{cusp}(n, \epsilon, w, \kappa_A\otimes (2,2)_{v |p})).$$

\subsubsection{Properties of locally analytic overconvergent cohomology}

\begin{prop} The complex $C_{cusp}(n, \epsilon, w,  \kappa_A\otimes (2,2)_{v |p})$ is represented by a bounded complex of projective  Banach $A[1/p]$-modules.
\end{prop}
\begin{proof}This follows easily by considering a \v{C}ech complex; see~ \cite[Prop.\
  12.8.2.1]{pilloniHidacomplexes}.\end{proof}

\begin{prop}\label{prop: cuspidal cohomology in degrees 0 I}  The cohomology ${\HH}_{cusp}^i (\dag, \kappa_A\otimes (2,2)_{v |p})$ vanishes for $i\notin [0,\# I]$. 
\end{prop} 

\begin{proof} This follows from Theorem \ref{theo-van}. 
\end{proof}

\subsubsection{Descent}\label{subsec-local-anal-ov-coh}
 We now assume that the character 
 $$\kappa_A: \prod_{v \in I} \ocal_{F_v}^\times \prod_{v \in I^c} (\ocal_{F_v}^\times)^2 \rightarrow A^\times$$
  is trivial on the torsion subgroup  of $\prod_{v \in I} \ocal_{F_v}^\times \prod_{v \in I^c} (\ocal_{F_v}^\times)^2 $  (of order prime to $p$ since $p>2$). 

The group $(\ocal_{F})_{(p)}^{\times, +}$ acts on
$\mathcal{X}_{K(I,p^n)}$, and the action factors through a finite
group. We let $\mathcal{X}^{G_1}_{K(I,p^n)}$ be the quotient. The
action of $(\ocal_{F})_{(p)}^{\times, +}$ can be lifted to the sheaf
$\mathcal{F}^{\kappa_A, w}$ by setting \[x: x^{*} \mathcal{F}^{\kappa_A, w}
\rightarrow \mathcal{F}^{\kappa_A, w}\] for all $x \in
(\ocal_{F})_{(p)}^{\times, +}$,  to be the composition of the
tautological isomorphism (the polarization is not used in the
construction of the sheaf) and multiplication by the character $d: (\ocal_{F})_{(p)}^{\times, +} \rightarrow \Lambda_I^\times \rightarrow A^\times$ of~\S\ref{subsub-descent}. 

We denote by $\mathcal{F}^{\kappa_A^{G_1}, w}$ the descended sheaf on $\mathcal{X}_{K(I,p^n)}^{G_1}$. We let \[C_{cusp}(G_1, n, \epsilon, w, \kappa_A\otimes (2,2)_{v |p})\] be the cohomology of the sheaf \[\mathcal{F}^{\kappa_A^{G_1},w} \otimes( \det
\omega_{\mathcal{G}})^2(-D)\] over
$\mathcal{X}^{G_1}_{K(I,p^n)}(\epsilon)$.  This is a direct factor of $C_{cusp}(n, \epsilon, w, \kappa_A\otimes (2,2)_{v |p})$. We also let $ {\HH}_{cusp} ^i( G_1, \dag, \kappa_A\otimes (2,2)_{v |p}) = \colim {\HH}^i ( C_{cusp}(G_1, n, \epsilon, w, \kappa_A\otimes (2,2)_{v |p})).$

\subsubsection{Spectral theory: construction of  the operator $U_{v,2}$} 

Firstly let $v \in I$.  We define an analytic adic space  $s_1: \mathcal{C}_{v,2} \rightarrow \mathcal{X}_{K(I,p^n)}(\epsilon)$ which parametrizes isogenies $A \rightarrow A'$ with associated Barsotti--Tate group $\mathcal{G} \rightarrow \mathcal{G}' $  whose kernel is a group $M_v \subset \mathcal{G}_v[p^2]$ which:

\begin{itemize}
\item is totally isotropic and locally isomorphic to $(\ocal_{F_v}/p\ocal_{F_v})^2 \oplus   \ocal_{F_v}/p^2 \ocal_{F_v}$,
\item has trivial intersection with $H_{v,n_v}$. 
\end{itemize}

There is a second projection $s_2: \mathcal{C}_{v,2} \rightarrow \mathcal{X}_{K(I,p^{n'})}(\epsilon')$ where: 

\begin{itemize}
\item $n' = (n'_{v'})_{v' |p}$ where $n'_{v} = n_{v}+1$, and $n'_{v'} = n_{v'}$ if $v' \neq v$.
\item $\epsilon' = ( \epsilon_{v'})_{v'|p}$ where $\epsilon'_{v} =
  \epsilon_v + 1$, and $\epsilon'_{v'} = \epsilon_{v'}$ if $v' \neq v$.
\end{itemize}
This map is provided by sending $(A, A')$ to $A'$, equipped with the subgroups:

\begin{itemize}
\item $ H'_{v', n_{v'}} = \mathrm{Im} (H_{v', n_{v'}})$ for all $v' \neq v$,
\item $ L'_{v', n_{v'}} = \mathrm{Im} (L_{v', n_{v'}})$ for all $v' \in I^c$,
\item $H'_{v, n_v+1} = \mathrm{Im} (p^{-1} H_{v, n_v})$ where $p^{-1} H_{v,n_v}$ is the pre-image in $\mathcal{G}_v[p^{n_v+1}]$ of $H_{v,n_v}$. 
\end{itemize}

One checks as in \cite[Lem.\  13.2.1.1]{pilloniHidacomplexes} (see
also Lemma~\ref{lem: u2 lands in epsilon prime} below) that the image
of $s_2$ lands in $\mathcal{X}_{K(I,p^{n'})}(\epsilon')$. The natural
map~$\omega_{\cG'}\to\omega_{\cG}$ induces a natural map $ s_2^{*}
\mathcal{F}^{\kappa_A,w'} \rightarrow  s_1^{*} \mathcal{F}^{\kappa_A,
  w}$,  where $w' = (w_{v'})$ with $w'_{v'}=w_{v'}$ if $v' \neq v$, and $w'_{v}   = w_{v} +1$. (See \cite[Lem.\  13.2.2.1]{pilloniHidacomplexes}.)
We deduce that there is  a normalized  map $ s_2^{*} \mathcal{F}^{\kappa_A,w'} \otimes (\det \omega_{\mathcal{G}})^2(-D) \rightarrow  s_1^{*} \mathcal{F}^{\kappa_A, w} \otimes (\det \omega_{\mathcal{G}})^2(-D)$  obtained by taking the tensor product of the above map and the normalized map (by $p^{-2}$) $s_2^{*} (\det \omega_{\mathcal{G}})^2(-D) \rightarrow s_1^{*} (\det \omega_{\mathcal{G}})^2(-D)$. 

We can therefore construct a Hecke operator $U_{v,2}: \mathrm{R}
\Gamma ( \mathcal{X}_{K(I,p^n)}( \epsilon) ,
\mathcal{F}^{\kappa_A,w}\otimes (\det \omega_{\mathcal{G}})^2(-D))
\rightarrow \mathrm{R}\Gamma ( \mathcal{X}_{K(I,p^n)},
\mathcal{F}^{\kappa_A,w} \otimes (\det \omega_{\mathcal{G}})^2(-D))$ by
the following composition:$$ \mathrm{R} \Gamma ( \mathcal{X}_{K(I,p^n)}( \epsilon) , \mathcal{F}^{\kappa_A,w}\otimes (\det \omega_{\mathcal{G}})^2(-D)) \rightarrow \mathrm{R} \Gamma ( \mathcal{X}_{K(I,p^{n'})}( \epsilon') , \mathcal{F}^{\kappa_A,w'}\otimes (\det \omega_{\mathcal{G}})^2(-D)) $$ $$ \rightarrow \mathrm{R} \Gamma ( \mathcal{C}_{v,2} , s_2^{*} \mathcal{F}^{\kappa_A,w'}\otimes (\det \omega_{\mathcal{G}})^2(-D)) \rightarrow $$
$$\mathrm{R} \Gamma ( \mathcal{C}_v , s_1^{*} \mathcal{F}^{\kappa_A,w}\otimes (\det \omega_{\mathcal{G}})^2(-D)) \rightarrow \mathrm{R} \Gamma ( \mathcal{X}_{K(I,p^n)}(\epsilon) ,(s_1)_{*}  s_1^{*} \mathcal{F}^{\kappa_A,w}\otimes (\det \omega_{\mathcal{G}})^2(-D)) $$ $$\stackrel{p^{-3}\mathrm{Tr}} \rightarrow  \mathrm{R} \Gamma ( \mathcal{X}_{K(I,p^n)}( \epsilon) , \mathcal{F}^{\kappa_A,w}\otimes (\det \omega_{\mathcal{G}})^2(-D)).$$

Now let $v \in I^c$.  We define an analytic adic space  $s_1: \mathcal{C}_{v,2} \rightarrow \mathcal{X}_{K(I,p^n)}(\epsilon)$ which parametrizes isogenies $\mathcal{G} \rightarrow \mathcal{G}' $  whose kernel is a group $M_v \subset \mathcal{G}_v[p^2]$ which:

\begin{itemize}
\item is totally isotropic and locally isomorphic to
  $(\ocal_{F_v}/p\ocal_{F_v})^2\oplus \ocal_{F_v}/p^2\ocal_{F_v}$, and
\item locally in the \'etale topology there is a symplectic isomorphism  $(\mathcal{G}_v)[p^\infty] \simeq (F_v/\ocal_{F_v})^4$ such that
  \begin{itemize}
  \item $M_v$ is generated by
    $p^{-1} e_{v,2}, p^{-1} e_{v,3}, p^{-2} e_{v,4}$, 
  \item $H_{v,n_v}$ is
    generated by $p^{-n_v} e_{v,1}$, and
  \item  $L_{v,n_v}$ is generated by
    $p^{-n_v} e_{v,1}$ and $p^{-n_v} e_{v,2}$.
  \end{itemize}

\end{itemize}

There is a second projection $s_2: \mathcal{C}_{v,2} \rightarrow \mathcal{X}_{K(I,p^{n})}(\epsilon)$, sending $(\mathcal{G}, \mathcal{G}')$ to $\mathcal{G}'$, equipped with the subgroups:

\begin{itemize}
\item $ H'_{v', n_{v'}} = \mathrm{Im} (H_{v', n_{v'}})$ for all $v' \neq v$.
\item $ L'_{v', n_{v'}} = \mathrm{Im} (L_{v', n_{v'}})$ for all $v' \in I^c$, $v' \neq v$.
\item In the notation above,  $H'_{v, n_v}$ is the group generated by the
image in $\mathcal{G}'$  of $p^{-n_{v}} e_{v,1}$ and $L'_{v, n_v}$ is
the group generated by the image of $p^{-n_{v}} e_{v,1}$ and
$p^{-n_{v}-1} e_{v,2}$.  One checks easily that these groups only
depend on $M_v$, $H_{v,n_v}$ and $L_{v, n_v}$ (and not on the choice
of symplectic basis).
\end{itemize}

Again, there is a natural map $ s_2^{*} \mathcal{F}^{\kappa_A,w}
\rightarrow  s_1^{*} \mathcal{F}^{\kappa_A, w}$. (See  \cite[Lem.\
13.2.2.1]{pilloniHidacomplexes} and \cite[\S 6.2]{MR3275848}.)
We deduce that there is  a normalized  map $ s_2^{*} \mathcal{F}^{\kappa_A,w} \otimes (\det \omega_{\mathcal{G}})^2(-D) \rightarrow  s_1^{*} \mathcal{F}^{\kappa_A, w} \otimes (\det \omega_{\mathcal{G}})^2(-D)$  obtained by taking the tensor product of the above map and the normalized map (by $p^{-2}$) $s_2^{*} (\det \omega_{\mathcal{G}})^2(-D) \rightarrow s_1^{*} (\det \omega_{\mathcal{G}})^2(-D)$. 

We can therefore construct a Hecke operator $U_{v,2}: \mathrm{R} \Gamma ( \mathcal{X}_{K(I,p^n)}( \epsilon) , \mathcal{F}^{\kappa_A,w}\otimes (\det \omega_{\mathcal{G}})^2(-D)) \rightarrow \mathrm{R}\Gamma ( \mathcal{X}_{K(I,p^n)}, \mathcal{F}^{\kappa_A,w} \otimes (\det \omega_{\mathcal{G}})^2(-D))$ by the following composition: 
$$ \mathrm{R} \Gamma ( \mathcal{X}_{K(I,p^n)}( \epsilon) , \mathcal{F}^{\kappa_A,w}\otimes (\det \omega_{\mathcal{G}})^2(-D)) \rightarrow \mathrm{R} \Gamma ( \mathcal{C}_{v,2} , s_2^{*} \mathcal{F}^{\kappa_A,w}\otimes (\det \omega_{\mathcal{G}})^2(-D)) \rightarrow $$
$$\mathrm{R} \Gamma ( \mathcal{C}_{v,2} , s_1^{*} \mathcal{F}^{\kappa_A,w}\otimes (\det \omega_{\mathcal{G}})^2(-D)) \rightarrow \mathrm{R} \Gamma ( \mathcal{X}_{K(I,p^n)}(\epsilon) ,(s_1)_{*}  s_1^{*} \mathcal{F}^{\kappa_A,w}\otimes (\det \omega_{\mathcal{G}})^2(-D)) $$ $$\stackrel{p^{-3}\mathrm{Tr}} \rightarrow  \mathrm{R} \Gamma ( \mathcal{X}_{K(I,p^n)}( \epsilon) , \mathcal{F}^{\kappa_A,w}\otimes (\det \omega_{\mathcal{G}})^2(-D)).$$

\begin{rem}When $v\in I^c$ we observe that $U_{v,2}$ itself is not
  improving analyticity and convergence in the $v$ direction (while it
  visibly does so in the case~$v\in I$). We next define an
  operator $U_{v,1}$ when $v \in I^c$. We will then show that the
  composite operator $U_{v,1} U_{v,2}$ improves analyticity and
  convergence. (This is related to our needing to use both the
  operators~$T_v$ and~$T_{v,1}$ at places $v\in I^c$ in~\S\ref{sec: Hida complexes}.)\end{rem}

\subsubsection{Spectral theory: construction of  the operator $U_{v,1}$} 

We let $v \in I^c$. We define an analytic adic space  $t_1: \mathcal{C}_{v,1} \rightarrow \mathcal{X}_{K(I,p^n)}(\epsilon)$ which parametrizes isogenies $A \rightarrow A'$ with associated Barsotti--Tate groups  $\mathcal{G} \rightarrow \mathcal{G}' $  whose kernel is a group $M_v \subset \mathcal{G}_v[p]$ which:

\begin{itemize}
\item is totally isotropic and locally isomorphic to $(\ocal_{F_v}/p\ocal_{F_v})^2$,
\item has trivial intersection with $L_{v,n_v}$. 
\end{itemize}

There is a second projection $t_2: \mathcal{C}_{v,1} \rightarrow \mathcal{X}_{K(I,p^{n})}(\epsilon)$, 
given by sending $(A, A')$ to $A'$, equipped with the subgroups:

\begin{itemize}
\item $ H'_{v', n_{v'}} = \mathrm{Im} (H_{v', n_{v'}})$ for all $v' $,
\item $ L'_{v', n_{v'}} = \mathrm{Im} (L_{v', n_{v'}})$ for all $v' \in I^c$.

\end{itemize}

There is a natural map $ t_2^{*} \mathcal{F}^{\kappa_A,w} \rightarrow
t_1^{*} \mathcal{F}^{\kappa_A, w}$ (again see \cite[Lem.\
13.2.2.1]{pilloniHidacomplexes} and \cite[\S 6.2]{MR3275848}). We deduce that there is  a   map $ t_2^{*} \mathcal{F}^{\kappa_A,w'} \otimes (\det \omega_{\mathcal{G}})^2(-D) \rightarrow  t_1^{*} \mathcal{F}^{\kappa_A, w} \otimes (\det \omega_{\mathcal{G}})^2(-D)$  obtained by taking the tensor product of the above map and the  map $t_2^{*} (\det \omega_{\mathcal{G}})^2(-D) \rightarrow t_1^{*} (\det \omega_{\mathcal{G}})^2(-D)$. 

We can therefore construct a Hecke operator $ U_{v,1}: \mathrm{R} \Gamma ( \mathcal{X}_{K(I,p^n)}( \epsilon) , \mathcal{F}^{\kappa_A,w}\otimes (\det \omega_{\mathcal{G}})^2(-D)) \rightarrow \mathrm{R}\Gamma ( \mathcal{X}_{K(I,p^n)}, \mathcal{F}^{\kappa_A,w} \otimes (\det \omega_{\mathcal{G}})^2(-D))$ by the following composition: 
$$ \mathrm{R} \Gamma ( \mathcal{X}_{K(I,p^n)}( \epsilon) , \mathcal{F}^{\kappa_A,w}\otimes (\det \omega_{\mathcal{G}})^2(-D)) \rightarrow \mathrm{R} \Gamma ( \mathcal{C}_{v,1} , t_2^{*} \mathcal{F}^{\kappa_A,w}\otimes (\det \omega_{\mathcal{G}})^2(-D)) \rightarrow $$
$$\mathrm{R} \Gamma ( \mathcal{C}_{v,1} , t_1^{*} \mathcal{F}^{\kappa_A,w}\otimes (\det \omega_{\mathcal{G}})^2(-D)) \rightarrow \mathrm{R} \Gamma ( \mathcal{X}_{K(I,p^n)}(\epsilon) ,(t_1)_{*}  t_1^{*} \mathcal{F}^{\kappa_A,w}\otimes (\det \omega_{\mathcal{G}})^2(-D)) $$ $$\stackrel{p^{-3}\mathrm{Tr}} \rightarrow  \mathrm{R} \Gamma ( \mathcal{X}_{K(I,p^n)}( \epsilon) , \mathcal{F}^{\kappa_A,w}\otimes (\det \omega_{\mathcal{G}})^2(-D)).$$

\subsubsection{Spectral theory: construction of the operator $U_{v,1} U_{v,2}$}

Let $v \in I^c$. We now consider the composite
operator $U_{v,1} U_{v,2}$. Our main task it to show that this
operator improves convergence and analyticity in the
$v$-direction. We begin by giving the correspondence corresponding to
this composite.

We define an analytic adic space  $u_1: \mathcal{C}_v \rightarrow \mathcal{X}_{K(I,p^n)}(\epsilon)$ which parametrizes isogenies $A \rightarrow A'$ with associated Barsotti--Tate groups $\mathcal{G} \rightarrow \mathcal{G}' $  whose kernel is a group $M_v \subset \mathcal{G}_v[p^3]$ which:

\begin{itemize}
\item is totally isotropic and locally isomorphic to $\ocal_{F_v}/p\ocal_{F_v} \oplus \ocal_{F_v}/p^2\ocal_{F_v} \oplus   \ocal_{F_v}/p^3 \ocal_{F_v}$,
\item locally in the \'etale topology there is a symplectic isomorphism  $(\mathcal{G}_v)[p^\infty] \simeq (F_v/\ocal_{F_v})^4$ such that
  \begin{itemize}
  \item $M_v$ is generated by
    $p^{-1} e_{v,2}, p^{-2} e_{v,3}, p^{-3} e_{v,4}$, 
  \item $H_{v,n_v}$ is generated by $p^{-n_v} e_{v,1}$, and
  \item  $L_{v,n_v}$
    is generated by $p^{-n_v} e_{v,1}$ and $p^{-n_v} e_{v,2}$.
  \end{itemize}

\end{itemize}

There is a second projection $u_2: \mathcal{C}_v \rightarrow \mathcal{X}_{K(I,p^{n'})}(\epsilon')$, given by sending $(A, A')$ to $A'$, equipped with the subgroups:

\begin{itemize}
\item $ H'_{v', n_{v'}} = \mathrm{Im} (H_{v', n_{v'}})$ for all $v' \neq v$,
\item $ L'_{v', n_{v'}} = \mathrm{Im} (L_{v', n_{v'}})$ for all $v' \in I^c$, $v' \neq v$,
\item In the notation above,  $H'_{v, n_v+1}$ is the group generated by the image in $\mathcal{G}'$  of $p^{-n_{v}-1} e_{v,1}$ and $L'_{v, n_v+1}$ is the group generated by the image of $p^{-n_{v}-1} e_{v,1}$ and $p^{-n_{v}-2} e_{v,2}$.  One checks easily that these groups only depend on $M_v$, $H_{v,n_v}$ and $L_{v, n_v}$. 

\end{itemize}

\begin{lem}\label{lem: u2 lands in epsilon prime} The image of $u_2$ lands in $\mathcal{X}_{K(I,p^{n'})}(\epsilon')$. 
\end{lem}
\begin{proof}We argue in the same way as in the proof of~\cite[Lem.\ 13.2.1.1]{pilloniHidacomplexes}. We fix symplectic bases $(e_{v,i})_{1 \leq i \leq 4}$ of $T_p(\mathcal{G})$,  $(e'_{v,i})_{1 \leq i \leq 4}$ of $T_p(\mathcal{G}')$, $(f_i)_{1 \leq i \leq 2}$ of  $\omega_{\mathcal{G}_v}^{mod}$, and $(f'_i)_{1 \leq i \leq 2}$ of $\omega_{\mathcal{G}'_v}^{mod}$ (compatible with the canonical filtration) such that there is a commutative diagram: 
\begin{eqnarray*}
\xymatrix{ T_p(\mathcal{G}_v) \ar[rr]^{\mathrm{diag}(1,p,p^2,p^3)}\ar[d]^{\mathrm{HT}} && T_p(\mathcal{G}_v') \ar[d]^{\mathrm{HT} }\\
\omega_{\mathcal{G}_v}^{mod}  \ar[rr]^{ \mathrm{diag}(p^2,p^3)} &&  \omega_{\mathcal{G}'}^{mod}}
\end{eqnarray*}

By definition we have that $\mathrm{HT} (e_{v,1}),
\mathrm{HT}(e_{v,2}) \in
p^{\epsilon_v}\omega_{\mathcal{G}_v}^{mod}$. On the other hand,
$\mathrm{HT}(e_{v,3}), \mathrm{HT}(e_{v,4})$ generate
$\omega_{\mathcal{G}_v}^{mod}$ and  $\mathrm{HT}(e'_{v,3}),
\mathrm{HT}(e'_{v,4})$ generate $\omega_{\mathcal{G}'_v}^{mod}$. The
group $L_{v,n_{v}+1}$ is generated by $\mathrm{diag}(1,p,p^2, p^3)\cdot
e_{v,1} = e_{v,1}'$ and $\mathrm{diag}(1,p,p^2,p^3)\cdot p^{-1} e_{v,2} = e_{v,2}'$. Therefore we deduce that 
$\mathrm{HT} (e_{v,1}'), \mathrm{HT}(e_{v,2}') \in p^{\epsilon_v+1}\omega_{\mathcal{G}}^{mod}$. 
\end{proof}

There is again a natural map $ u_2^{*} \mathcal{F}^{\kappa_A,w'}
\rightarrow  u_1^{*} \mathcal{F}^{\kappa_A, w}$  where $w' = (w_{v'})$
with $w'_{v'} = w_{v'}$ is $v' \neq v$ and $w'_{v}   = w_{v} +1$; see \cite[Lem.\
13.2.2.1]{pilloniHidacomplexes} and \cite[\S 6.2]{MR3275848}. 
We deduce that there is  a normalized  map $ u_2^{*} \mathcal{F}^{\kappa_A,w'} \otimes (\det \omega_{\mathcal{G}})^2(-D) \rightarrow  u_1^{*} \mathcal{F}^{\kappa_A, w} \otimes (\det \omega_{\mathcal{G}})^2(-D)$  obtained by taking the tensor product of the above map and the normalized map (by $p^{-2}$) $u_2^{*} (\det \omega_{\mathcal{G}})^2(-D) \rightarrow u_1^{*} (\det \omega_{\mathcal{G}})^2(-D)$. 

We can therefore construct a Hecke operator $U_{v,2} U_{v,1}: \mathrm{R} \Gamma ( \mathcal{X}_{K(I,p^n)}( \epsilon) , \mathcal{F}^{\kappa_A,w}\otimes (\det \omega_{\mathcal{G}})^2(-D)) \rightarrow \mathrm{R}\Gamma ( \mathcal{X}_{K(I,p^n)}, \mathcal{F}^{\kappa_A,w} \otimes (\det \omega_{\mathcal{G}})^2(-D))$ by the following composition: 
$$ \mathrm{R} \Gamma ( \mathcal{X}_{K(I,p^n)}( \epsilon) , \mathcal{F}^{\kappa_A,w}\otimes (\det \omega_{\mathcal{G}})^2(-D)) \stackrel{res}\rightarrow \mathrm{R} \Gamma ( \mathcal{X}_{K(I,p^{n'})}( \epsilon') , \mathcal{F}^{\kappa_A,w'}\otimes (\det \omega_{\mathcal{G}})^2(-D)) $$ $$ \rightarrow \mathrm{R} \Gamma ( \mathcal{C}_v , u_2^{*} \mathcal{F}^{\kappa_A,w'}\otimes (\det \omega_{\mathcal{G}})^2(-D)) \rightarrow $$
$$\mathrm{R} \Gamma ( \mathcal{C}_v , u_1^{*} \mathcal{F}^{\kappa_A,w}\otimes (\det \omega_{\mathcal{G}})^2(-D)) \rightarrow \mathrm{R} \Gamma ( \mathcal{X}_{K(I,p^n)}(\epsilon) ,(u_1)_{*}  u_1^{*} \mathcal{F}^{\kappa_A,w}\otimes (\det \omega_{\mathcal{G}})^2(-D)) $$ $$\stackrel{p^{-6}\mathrm{Tr}} \rightarrow  \mathrm{R} \Gamma ( \mathcal{X}_{K(I,p^n)}( \epsilon) , \mathcal{F}^{\kappa_A,w}\otimes (\det \omega_{\mathcal{G}})^2(-D))$$

We now set $U^I = \prod_{v \in I} U_{v,2} \prod_{v \in I^c} U_{v,1} U_{v,2}$. 

\begin{lem}\label{lem: U operator is compact} The operator $U^I$ acting on $C_{cusp} (G_1, n, \epsilon, w,\kappa_A \otimes (2,2)_{v |p})$ is compact. Moreover, for $n+1 = (n_v+1)_{v |p}$, $\epsilon +1 = (
\epsilon_v+1)_v$ and $w +1 = (w_v +1)_v$ we have a factorization (where the vertical maps are the natural restriction maps):
\begin{eqnarray*} 
\xymatrix{C_{cusp} (G_1, n, \epsilon, w,\kappa_A \otimes (2,2)_{v |p}) \ar[r]^{U^I} \ar[d] & C_{cusp} (G_1, n, \epsilon, w,\kappa_A \otimes (2,2)_{v |p}) \ar[d] \\
C_{cusp} (G_1, n+1, \epsilon+1, w+1,\kappa_A \otimes (2,2)_{v |p}) \ar[r]^{U^I} \ar[ru] & C_{cusp} (G_1, n+1, \epsilon+1, w+1,\kappa_A \otimes (2,2)_{v |p}) }
\end{eqnarray*}
\end{lem}  
\begin{proof} By construction, the action of $U^I$ can be factored into  $$\mathrm{R} \Gamma ( \mathcal{X}^{G_1}_{K(I,p^n)}( \epsilon) , \mathcal{F}^{\kappa^{G_1}_A,w}\otimes (\det \omega_{\mathcal{G}})^2(-D)) \rightarrow  \mathrm{R} \Gamma ( \mathcal{X}^{G_1}_{K(I,p^{n+1})}( \epsilon+1) , \mathcal{F}^{\kappa^{G_1}_A,w+1}\otimes (\det \omega_{\mathcal{G}})^2(-D))$$$$ \stackrel{\tilde{U}^I} \rightarrow \mathrm{R} \Gamma ( \mathcal{X}^{G_1}_{K(I,p^n)}( \epsilon) , \mathcal{F}^{\kappa^{G_1}_A,w}\otimes (\det \omega_{\mathcal{G}})^2(-D))$$
where $n+1 = (n_v+1)_{v |p}$ and $\epsilon +1 = (
\epsilon_v+1)_v$, $w +1 = (w_v +1)_v$. It is enough to show that the
map  $\mathrm{R} \Gamma ( \mathcal{X}^{G_1}_{K(I,p^n)}( \epsilon) ,
\mathcal{F}^{\kappa^{G_1}_A,w}\otimes (\det \omega_{\mathcal{G}})^2(-D))
\rightarrow  \mathrm{R} \Gamma ( \mathcal{X}^{G_1}_{K(I,p^{n+1})}(
\epsilon+1) , \mathcal{F}^{\kappa^{G_1}_A,w+1}\otimes (\det
\omega_{\mathcal{G}})^2(-D))$ is compact. This follows by
consideration of an appropriate \v{C}ech complex, as in \cite[Lem.\ 13.2.4.1]{pilloniHidacomplexes}.\end{proof}

\subsubsection{Spectral theory: local constancy of the Euler Characteristics}
Let $ W_{cl}$ be the set of weights  $\kappa = ((k_v, l_v))_{v |p}
\in \ZZ^{S_p}$,  with $l_v = 2$ if $v \in I$, $k_v \equiv l_v \equiv 2~
\mathrm{mod} (p-1)$ for all $v |p$. It is equipped with the
$p$-adic topology. 

For all $\kappa \in W_{cl}$,  we let $C_{cusp} (G_1, n, \epsilon, w, \kappa)$ be $n, \epsilon$-convergent, $w$-analytic cohomology of weight $\kappa$ and we set ${\HH}^i_{cusp} (G_1,  \dag , \kappa) = \colim \HH^i(C_{cusp} (G_1, n, \epsilon, w, \kappa))$. In other words, following the notation of~\S\ref{subsec-local-anal-ov-coh}, we have $A = \whalingship_p$ and $\kappa_A = \kappa \otimes (-2,-2)_{v |p}$. 
It follows from Lemma~\ref{lem: U operator is
  compact} that the cohomology groups $e(U^I) {\HH}^i_{cusp} (G_1,
\dag, \kappa)$ are finite-dimensional. The following standard
consequence of our constructions will be crucial in our comparison in~\S\ref{subsec: ordinary
cohomology is overconvergent} of the complexes constructed in~\S\ref{sec: Hida complexes} and the overconvergent cohomology we
are considering in this section.
\begin{thm}\label{thm: local constancy of Euler char}  The map
\begin{eqnarray*}
W_{cl} & \rightarrow & \ZZ \\
\kappa &\mapsto &\sum_{i} (-1)^i \dim e(U^I) {\HH}^i_{cusp} (G_1, \dag, \kappa)
\end{eqnarray*}
 
 is locally constant.
\end{thm}

\begin{proof} This follows from Coleman's theory
  \cite[\S A5]{1997InMat.127..417C}, as in~\cite[\S
  13.4]{pilloniHidacomplexes}. Indeed, there is a perfect complex~$C^\bullet$
  interpolating $C_{cusp} (G_1, n, \epsilon, w, \kappa)$ over the
  spectral variety, and the dimensions of the slope zero parts of
  the~$C^i$ are locally constant.
  \end{proof}

\begin{rem} In particular if $\# I = 1$, we deduce that 
$$ \kappa \mapsto \dim e(U^I) {\HH}^0_{cusp} (G_1, \dag,  \kappa) - \dim e(U^I) {\HH}^1_{cusp} (G_1, \dag, \kappa) $$
is locally constant.  We will use this in~\S\ref{subsec:
  ordinary cohomology is overconvergent} to reduce the comparison of
ordinary and overconvergent cohomology to the case of high weight,
where the control theorems proved in~\S\ref{sec: Hida complexes}
apply.
\end{rem}

\subsection{Locally analytic overconvergent classes and algebraic
  overconvergent classes}\label{subsec: comparing overconvergent
  analytic and algebraic classes}
Let $\kappa = ((k_v, l_v))_{v |p}$ with $l_v = 2$ if $v \in I$,
$k_v \equiv l_v \equiv 2 ~\mathrm{mod} ~(p-1)$ be a dominant algebraic weight. 
\begin{prop}\label{prop: slopes bounded below} On ${\HH}^i_{cusp}(G_1, \dag, \kappa)$, the slopes of
  $(U_{v,1})_{v \in I^c}$ and $(U_{v,2})_{v |p}$ are  $\geq -3$. On
  ${\HH}^0_{cusp}(G_1, \dag, \kappa)$ they are  $\ge 0$. 
\end{prop}

\begin{proof} The proof of~ \cite[Prop.\
  13.3.1.1]{pilloniHidacomplexes} goes through essentially without change.\end{proof}

Below, we denote by $\mathcal{F}^{\kappa,w^-} = \colim_{w' < w} \mathcal{F}^{\kappa,w'}$. 
\begin{prop}\label{prop: BGG}  Let $\kappa = ((k_v, l_v))_{v |p}$ with
  $l_v = 2$ if $v \in I$, $k_v = l_v = 2 ~\mathrm{mod}~(p-1)$ be a
  dominant algebraic weight.  There is a relative analytic BGG resolution:
$$ 0 \rightarrow \omega^\kappa(-D) \rightarrow \mathcal{F}^{\kappa,w^-} \otimes (\det \omega_{\mathcal{G}})^2(-D) \rightarrow  $$ $$\bigoplus_{s \in W^{(1)}}  \mathcal{F}^{s \bullet \kappa,w^-} \otimes (\det \omega_{\mathcal{G}})^2(-D) \rightarrow  \ldots \rightarrow \bigoplus_{s \in W^{(d)}}  \mathcal{F}^{s \bullet \kappa,w^-} \otimes (\det \omega_{\mathcal{G}})^2(-D) \rightarrow 0$$
where $W$ is the  Weyl group of  $\mathrm{GL}_2( F \otimes_{\Q} \Q_p)$, $W^{(i)}$ stands for the elements of length $i$ in $W$, and  $\bullet$ is the twisted Weyl action. 
\end{prop}\begin{proof}  This is a relative version of the main result of
  \cite{JONES20111616}, and is proved in~\cite[\S
  7.2]{MR3275848}. (Note though that there is a minor error there; one
  needs to replace~$\mathcal{F}^{\kappa,w}$
  with~$\mathcal{F}^{\kappa,w^-}$ as defined above, but having made
  this change, the arguments go through unchanged.)\end{proof} 

The actions of $U_{v,1}$ and $U_{v,2}$ by cohomological correspondences on the sheaf $\mathcal{F}^{\kappa,w^-} \otimes (\det \omega_{\mathcal{G}})^2(-D)$ restrict to actions on the subsheaf $\omega^\kappa(-D)$ (and the action of $U^I$ is compact on the cohomology).

\begin{cor}\label{cor-boundslope} Let $\kappa = ((k_v, l_v))_{v |p}$
  with $l_v = 2$ if $v \in I$, $k_v \equiv l_v \equiv 2 ~\mathrm{mod}~
  (p-1)$ be a dominant algebraic  weight. Then the map $$e(U^I){\HH}^i( \mathcal{X}_{K(I,p^n)}(\epsilon), \omega^\kappa(-D)) \rightarrow e(U^I) {\HH}^i( \mathcal{X}_{K(I,p^n)}(\epsilon), \mathcal{F}^{\kappa, w} \otimes( \det \omega_{\mathcal{G}})^2(-D))$$ is an isomorphism for $i =0$ and injective if $i = 1$. It is an isomorphism for $i =1$ if we further assume that $k_v-l_v \geq 3$ for all $v|p$.
\end{cor}

\begin{proof} Proposition \ref{prop: BGG} gives a spectral sequence  $E_1^{p,q} = \oplus_{s \in W^{(p)}} {\HH}^q( \mathcal{X}_{K(I,p^n)}(\epsilon),  \mathcal{F}^{s \bullet\kappa, w} \otimes( \det \omega_{\mathcal{G}})^2(-D))$ converging to  ${\HH}^{p+q}( \mathcal{X}_{K(I,p^n)}(\epsilon), \omega^\kappa(-D))$. We shall see that the ordinary projector kills  the terms $E_1^{p,q}$ of the spectral sequence for $p > 1$ under a suitable normalization of the action of the Hecke operators and suitable assumptions on the weight $\kappa$.  We analyze the differentials  of proposition \ref{prop: BGG}:  $\bigoplus_{s \in W^{(i)}}  \mathcal{F}^{s \bullet \kappa,w^-} \otimes (\det \omega_{\mathcal{G}})^2(-D)  \rightarrow \bigoplus_{s \in W^{(i+1)}}  \mathcal{F}^{s \bullet \kappa,w^-} \otimes (\det \omega_{\mathcal{G}})^2(-D)$. We let $W = \prod_{v | p} \{ 1_v, w_v\}$ with $\ell(w_v)=1$. For any subset $J$ of places diving $p$, we let $w_J= \prod_{v \in J} w_v$.  The above map is given by the product of the maps $\theta_{s,s'}: \mathcal{F}^{s \bullet \kappa,w^-} \otimes (\det \omega_{\mathcal{G}})^2(-D)  \rightarrow  \mathcal{F}^{s' \bullet \kappa,w^-} \otimes (\det \omega_{\mathcal{G}})^2(-D)$ for $s = w_J$ (for  a subset $J$   of cardinality $i$) and $s' = w_{J  \cup \{ v\}}$ for $v \notin J$. 
By \cite[\S
    7.3]{MR3275848}, this map induces on cohomology an  equivariant
    map for the operators $U_{w,i}$ for $w \neq v$ and $U_{v,1}$; and on the other hand, we have  $U_{v,2} \circ \theta_{s,s'} = p^{(k_v-l_v)+1}\theta_{s,s'} \circ U_{v,2}$. 
    A way to interpret this relation is to say that the spectral sequence is equivariant for the action of Hecke operators, if the standard action of $U_{v,2}$ on $\HH^i(\mathcal{X}_{K(I,p^n)}(\epsilon), \mathcal{F}^{w_{J} \bullet \kappa, w} \otimes( \det \omega_{\mathcal{G}})^2(-D))$ is  twisted by multiplication by $p^{k_v-l_v+1}$ if $v \in J$. 
The corollary therefore follows from the slope bounds of Propositions~\ref{prop: slopes bounded
    below}.\end{proof}
\subsection{Small slope forms are classical}\label{subsec: small slope
implies classical}

\subsubsection{Fargues' degree function}\label{subsec: degree function}We
now recall some results on the degree of quasi-finite flat group
schemes, following
the papers~\cite{MR2673421,MR2919687}. Let~$K$ be a complete
valued extension of~$\Qp$ with corresponding valuation
$v:K\to\R\cup\{\infty\}$, which we assume to be normalized so
that~$v(p)=1$. We also write $v:\cO_K/p\cO_K\to [0,1]$ for the induced
map. If~$M$ is a finitely presented torsion $\cO_K$-module, then we
can write $M\cong\oplus_{i=1}^r\cO_K/x_i$ for some $x_i\in\cO_K$, and
we set $\deg M:=\sum_{i=1}^rv(x_i\mod p)$. 

If~$H$ is a group scheme over~$\cO_K$, we let~$\omega_H$ denote the
conormal sheaf to the identity section. If $H$ is finite flat, then
$\omega_H$ is finitely presented and torsion over~$\cO_K$, and
following Fargues we define the degree of~$H$ to be \[\deg H:=\deg\omega_H.\]

More generally, let $A\rightarrow A'$ be an isogeny of semi-abelian schemes with associated $p$-divisible groups $\cG \rightarrow \cG'$ over some analytic adic space $S$.  We denote by $\omega_{\cG}$ and $\omega_{\cG'}$ the conormal sheaves of $A$ and $A'$ along their unit sections and by $\omega_{\cG}^+$ and $\omega_{\cG'}^+$ the subsheaf of integral differentials (which means that locally on $S$ they arise from differentials on a formal model of  $A$ or $A'$).  Let $H$ be the kernel of $\cG \rightarrow \cG'$.  This is a quasi-finite group scheme. To this isogeny we may attach a section  $\delta_{H}$ of the locally free sheaf of rank one $\det \omega_{\cG'} \otimes \det \omega_{\cG}^{-1}$. Moreover, this section lies in the subsheaf $(\det \omega_{\cG'})^+ \otimes (\det \omega_{\cG}^{-1})^+$ of integral differential forms. For each point $x \in S$, we may compute the associated norm $\vert \delta_{H} \vert_x$, by choosing a trivialization of $(\det \omega_{\cG'})^+ \otimes (\det \omega_{\cG}^{-1})^+$ in a neighbourhood of $x$ and viewing $\delta_H$ as a function (the norm $\vert \delta_{H} \vert_x$ is independent of the trivialization). If $x \in S$ is a rank one point with associated valuation normalized by $v_x(p) = 1$,  and if $H_x \subset \cG_x$ extends to a finite flat group scheme on a formal model $\mathfrak{G}_x$ of $\cG_x$ over $\Spec k(x)^+$, then $v_x( \delta_H) = \deg H_x$. 

 \subsubsection{Neighbourhoods of the ordinary locus}  \label{section:neighbourhoods}

 Recall from~\S\ref{subsec: formal geometry} that we define  $$K_p(I) = \prod_{v \in I} \Kli(v) \prod_{v \in I^c} \Iw(v).$$

 We can consider ${X}_{K_p(I) K^p}$. Let $\mathcal{X}_{K_p
   K_p(I)}$ be the associated analytic space. For each $v |p$,  we have an isogeny $\mathcal{G} \rightarrow \mathcal{G}'$ whose kernel is a quasi-finite group scheme $H_v$ which is of order $p$ away from the boundary.   For each $v \in I^c$,  we have an isogeny $\mathcal{G} \rightarrow \mathcal{G}'$ whose kernel is a quasi-finite group scheme $L_v$ which is of order $p^2$ away from the boundary.

    Let $\mathcal{X}_{K^p
   K_p(I)}^{rk 1}$ be the subset of rank one points. To each rank one point  $x$ is associated a rank one valuation $v_x: \ocal_{\mathcal{X}_{K_p(I) K^p,x}} \rightarrow \mathbb{R} \cup \{ \infty \}$ which we normalize by $v_x(p) = 1$. 
   If~$v\in I$, we define   $\mathrm{deg}_v: \mathcal{X}_{K^p
   K_p(I)}^{rk 1} \rightarrow [0,1]$  by $\deg_v( x) = v_x(\delta_{H_v})$. Similarly, for all $v \in I^c$, we define   $\mathrm{deg}_v: \mathcal{X}_{K^p
   K_p(I)}^{rk 1} \rightarrow [0,2]$  by $\deg_v( x) = v_x( \delta_{L_{v}})$. 
   
   We can put all these degree functions together into a function 
   $$ \mathrm{deg}: \mathcal{X}_{K^p
   K_p(I)}^{rk 1} \rightarrow [0,1]^{I} \times [0,2]^{I^c}.$$

   For each rational  interval $J \subset [0,1]^{I} \times [0,2]^{I^c}$, there is a unique quasi-compact open subset $\mathcal{X}_{K^p
   K_p(I)}(J) \subset \mathcal{X}_{K^p
   K_p(I)}$ such that   $\mathcal{X}_{K^p
   K_p(I)}(J)^{rk 1} = \mathrm{deg}^{-1}(J)$. 
   
   Of particular interest is the multiplicative locus: $$\mathcal{X}_{K^p
   K_p(I)}^{\mult} = \mathcal{X}_{K^p
   K_p(I)}(\{1\}^{I} \times \{2\}^{I^c}).$$ Let  $(\epsilon_v) \in   ([0,1]^{I} \times [0,2]^{I^c}) \cap  \mathbb{Q}^{S_p}$ and set $$\mathcal{X}_{K^p
   K_p(I)}((\epsilon_v)_{v \in S_p}) = \mathcal{X}_{K^p
   K_p(I)}( \prod_{v \in I} [1- \epsilon_v, 1]^{I} \times  \prod_{v \in I^c} [2- \epsilon_v, 2]^{I^c}).$$
   
 Observe that $\mathcal{X}_{K^p
   K_p(I)}^{\mult} = \mathcal{X}_{K^p
   K_p(I)}((0)_{v \in S_p})$ while $\{\mathcal{X}_{K^p
   K_p(I)}((\epsilon_v)_{v \in S_p})\}_{\epsilon_v \rightarrow 0^+, \forall v \in S_p}$ is a fundamental system of strict neighbourhoods of $\mathcal{X}_{K^p  K_p(I)}^{\mult}$.

   All these spaces are stable under the action of ${\ocal_{F}}_{(p)}^{\times, +}$ on the polarization, and descend to open subspaces of $\mathcal{X}^{G_1}_{K^p
   K_p(I)}$. We can therefore add a superscript $G_1$ to any of these spaces with the obvious meaning.

\subsubsection{Comparison between spaces of overconvergent cohomology} In this section we make the connection between  the spaces  $\mathcal{X}_{K^p
   K_p(I)}((\epsilon_v)_{v \in S_p})$ (with $\epsilon_v \in ([0,1]^{I} \times [0,2]^{I^c}) \cap  \mathbb{Q}^{S_p}$) that we just introduced and the spaces $\mathcal{X}_{K(I,p^n)}((\epsilon_v)_{v \in S_p})$ (say for parallel $n \in \Z_{\geq 1}$ and with $\epsilon_v \in ([0,n-\frac{1}{p-1}] \cap  \mathbb{Q})^{S_p}$) introduced in \S\ref{recollection-Hodgetate}.   Both types of spaces are neighbourhoods of the multiplicative locus in an appropriate sense. The previous spaces are well adapted to the construction of interpolation sheaves and eigenvarieties while these new spaces appear naturally when one wants  to prove classicity theorems. 
   
   There is a natural forgetful  map $\mathcal{X}_{K(I,p^n)}((\epsilon_v)_{v \in S_p}) \rightarrow \mathcal{X}_{K^p
   K_p(I)}$. By \cite[Lem.\ 14.1.1]{pilloniHidacomplexes} (for the
 places $v \in I$, and a trivial extension for the places $v \in
 I^c$), this  map factors into a map
 $\mathcal{X}_{K(I,p^n)}((\epsilon_v)_{v \in S_p}) \rightarrow
 \mathcal{X}_{K_p(I) K^p}((1- \frac{2}{n}(n-\epsilon_v +
 \frac{1}{p-1}))_{v \in I}\times (2- \frac{2}{n}(n-\epsilon_v +
 \frac{1}{p-1}))_{v \in I^c})$. Observe that when $\epsilon_v = n-\frac{1}{p-1}$ and $n \rightarrow \infty$, $1- \frac{2}{n}(n-\epsilon_v +
 \frac{1}{p-1}) \rightarrow 1$ and $2- \frac{2}{n}(n-\epsilon_v +
 \frac{1}{p-1}) \rightarrow 2$. Conversely, by \cite[Lem.\ 14.1.2]{pilloniHidacomplexes}, there is a natural inclusion: $\mathcal{X}_{K_p(I) K^p}((\epsilon_v)_{v \in S_p}) \hookrightarrow \mathcal{X}_{K(I,p)}((1- \frac{1}{p-1})_{v \in S_p} )$ for all $\epsilon_v \geq 1- \frac{1}{p}$ if $v \in I$ and $\epsilon_v \geq 2- \frac{1}{p}$ if $v \in I^c$.

\begin{lem}\label{lem: U is compact on classical sheaves} Let $(\epsilon_v) \in ([  1- \frac{1}{p},1)^I \times [2- \frac{1}{p},2)^{I^c})\cap \Q^{S_p}$. Let $\kappa$ be a classical algebraic weight. 
\begin{enumerate}
\item The cohomology   $\mathrm{R}\Gamma(\mathcal{X}_{K^p
   K_p(I)}((\epsilon_v)_{v \in S_p}), \omega^\kappa)$ carries an action of the operators $U_{v,1}$ and $U_{v,2}$. 
   \item The operator $ \prod_{v \in I} U_{v,2} \prod_{v \in I^c} U_{v,1}$ is compact $\mathrm{R}\Gamma(\mathcal{X}^{G_1}_{K^p
   K_p(I)}((\epsilon_v)_{v \in S_p}), \omega^\kappa)$.
   \item The canonical map    $\mathrm{R}\Gamma(\mathcal{X}^{G_1}_{K^p
   K_p(I)}((\epsilon_v)_{v \in S_p}), \omega^\kappa) \rightarrow \mathrm{R}\Gamma(\mathcal{X}^{G_1}_{K(I,p)}((1- \frac{1}{p-1})_{v \in S_p} ), \omega^\kappa)$ induces a quasi-isomorphism on the finite slope part for $\prod_{v \in I} U_{v,2} \prod_{v \in I^c} U_{v,1}$.
   \item The same holds for cuspidal cohomology.
   \   \end{enumerate}
   \end{lem} 
   
  \begin{proof} The definition of the operators is a routine
    computation. To prove compactness, we need to show that the
    operators improve convergence. This is entirely parallel to Lemma~
    \ref{lem: U operator is compact}. For the degree functions
    considered here this follows from  of  \cite[Prop.\
    2.3.6]{MR2783930}. The quasi-isomorphism follows from  an easy analytic continuation argument (see Lemma~ \ref{lem-extension-first} below, for example). 
  \end{proof} 
   
   It is sometimes convenient to  consider the dagger
space \numequation\label{eqn: X mult}\mathcal{X}_{K^p
   K_p(I)}^{\mult, \dag}:= \lim_{\epsilon_v \rightarrow 0^+} \mathcal{X}_{K^p
   K_p(I)}((\epsilon_v)_{v \in S_p}).\end{equation}
and its $G_1$-variant. 
In view of the previous lemma we can define the complex $e(U^I) \mathrm{R} \Gamma ( \mathcal{X}^{G_1, \mult, \dag}_{K_p(I) K^p},
\omega^{\kappa}(-D))$ as being equal to $e(U^I)\mathrm{R}\Gamma(\mathcal{X}^{G_1}_{K^p
   K_p(I)}((\epsilon_v)_{v \in S_p}), \omega^\kappa)$. 
   
   \begin{lem}\label{lem: ordinary mult cohomology is perfect} 
   The complex $e(U^I) \mathrm{R} \Gamma ( \mathcal{X}^{G_1, \mult, \dag}_{K_p(I) K^p},
\omega^{\kappa}(-D))$ is a perfect complex supported in degrees~$[0,\#I]$.
\end{lem}
   \begin{proof} That the cohomology vanishes outside of degrees~$[0,\#I]$
  follows as usual by pushing forward to the minimal
  compactification. The finiteness of the cohomology follows from the
  compactness of~$U^I$.
  \end{proof}
\subsubsection{Main classicity theorem} 
We now state our main classicity result for overconvergent cohomology,
which we will prove using a generalization of the analytic continuation method
of~\cite{MR2219265} to higher degree cohomology, which was proved
in~\cite[\S 3]{pilloniHidacomplexes}.  Let $\kappa = (k_v, l_v)_{v
  |p}$ be a dominant algebraic weight. There is a  canonical restriction map $$ \mathrm{R} \Gamma ( \mathcal{X}_{K_p(I) K^p}, \omega^{\kappa}(-D)) \rightarrow \mathrm{R} \Gamma( \mathcal{X}_{K_p(I) K^p}^{\mult, \dag}, \omega^\kappa(-D))$$ 
which is equivariant for the Hecke operators $U_{v,1}$ and $U_{v,2}$. 

\begin{thm}\label{thm-classicality} The canonical map 
$$ \mathrm{R}\Gamma ( \mathcal{X}^{G_1}_{K_p(I) K^p}, \omega^{\kappa}(-D))[ U_{v,2} < k_v+l_v-3~v\in I,~ U_{v,1} < l_v-3~v \in I^c] $$
$$  \rightarrow \mathrm{R}\Gamma ( \mathcal{X}_{K_p(I) K^p}^{G_1,\mult, \dag}, \omega^\kappa(-D))[ U_{v,2} < k_v+l_v-3~v\in I,~U_{v,1} < l_v-3~v \in I^c]$$
is a quasi-isomorphism.\end{thm}
\begin{rem}The meaning of
$[ U_{v,2} < k_v+l_v-3~v\in I,~ U_{v,1} < l_v-3~v \in I^c]$ in
Theorem~\ref{thm-classicality} is the
obvious one: it means the part of slope less than $k_v+l_v-3$ for
$U_{v,2}$ at $v \in I$ and less than $l_v-3$ for $U_{v,1}$ at
$v \in I^c$. (Note that while the individual operators~$U_{v,1},
U_{v,2}$ do not act compactly on the complex on the right hand side,
their product~$U$ does by Lemma~\ref{lem: U is compact on classical
  sheaves} It follows the individual operators~$U_{v,1},
 		U_{v,2}$ act compactly on the part of the
 		complex with bounded slope for~$U$, and so this small slope part is
well-defined by the procedure explained at the start of this
section.)
\end{rem}

\begin{rem} When $I= \emptyset$, Theorem~\ref{thm-classicality} (for ${\HH}^0$) is
  proved in~\cite{MR3488741}. It may be possible to improve on the
  bound~$l_v-3$ at the places $v\in I^c$, but this does not matter for
  our purposes.\end{rem}

\subsubsection{Hecke correspondences again} Let $w \in I$.  We
consider the following correspondence, whose corresponding Hecke operator is $U_{w,2}^n$ (the $n$th
iterate of $U_{w,2}$):
$t_{w,n,1}, t_{w,n,2}: C^{(n)}_{w} \rightarrow \mathcal{X}_{K_p(I) K^p}$, which parametrizes $(\mathcal{G}, \{H_v\}_{v \in I} , \{ H_v
\subset L_v\}_{v \in I^c}, \mathcal{G} \rightarrow \mathcal{G}_n)$
where the isogeny $\mathcal{G} \rightarrow \mathcal{G}_{n}$ has kernel
$M_{n,w} \subset \mathcal{G}_w[p^{2n}]$ which is totally isotropic and
locally isomorphic to $(\ocal_{F_w}/p^n)^2 \oplus \ocal_{F_w}/p^{2n}$,
and satisfies $M_{n,w} \cap H_w = \{0\}$. 
The first projection is  $$t_{w,n,1}\big((\mathcal{G}, \{H_v\}_{v \in I} , \{ H_w \subset L_v\}_{v \in I^c}, \mathcal{G} \rightarrow \mathcal{G}_n)\big)= (\mathcal{G}, \{H_v\}_{v \in I} , \{ H_v \subset L_v\}_{v \in I^c})$$
and the second projection is $$t_{w, n,2}\big((\mathcal{G}, \{H_v\}_{v \in I} , \{ H_v \subset L_v\}_{v \in I^c}, \mathcal{G} \rightarrow \mathcal{G}_n)\big)= (\mathcal{G}_n, \{H'_v\}_{v \in I} , \{ H'_v \subset L'_v\}_{v \in I^c})$$
where $\{H'_v\}_{v \in I}$ and $\{H'_v \subset L'_v\}_{v \in I^c}$ are the images of $\{H_v\}_{v \in I}$ and $\{H_v \subset L_v\}_{v \in I^c}$ in $\mathcal{G}_n$. 

 There are cohomological correspondences 
$$ (t_{w,n,1})_{*} t_{w,n,2}^{*} \omega^\kappa \rightarrow
\omega^\kappa,~  (t_{w,n,1})_{*} t_{w,n,2}^{*} \omega^\kappa(-D)
\rightarrow \omega^\kappa(-D),$$which give~$U_{w,2}^n$.
Moreover, these cohomological correspondences restrict to 
$$ (t_{w,n,1})_{*} t_{w,n,2}^{*} (\omega^\kappa)^{++} \rightarrow
p^{-3n}(\omega^\kappa)^{++},~(t_{w,n,1})_{*} t_{w,n,2}^{*}
(\omega^\kappa(-D))^{++} \rightarrow
p^{-3n}(\omega^\kappa(-D))^{++},$$and they induce maps on cohomology in the usual way.

 Let $w \in I^c$. We consider the correspondence:  $t_{w,n,1},
 t_{w,n,2}: C^{(n)}_{w} \rightarrow \mathcal{X}_{K_p(I) K^p}$ which
 parametrizes $(\mathcal{G}, \{H_v\}_{v \in I} , \{ H_v \subset
 L_v\}_{v \in I^c}, \mathcal{G} \rightarrow \mathcal{G}_n)$ where the
 isogeny $\mathcal{G} \rightarrow \mathcal{G}_n$ has kernel $M_{n,w}
 \subset \mathcal{G}_w[p^n]$ which is totally isotropic, locally
 isomorphic to $(\ocal_{F_w}/p^n)^2 $, and satisfies $M_{n,w} \cap L_w = \{0\}$. 
The first projection is  $$t_{w,n,1}\big((\mathcal{G}, \{H_v\}_{v \in
  I} , \{ H_w \subset L_v\}_{v \in I^c}, \mathcal{G} \rightarrow
\mathcal{G}_n)\big)= (\mathcal{G}, \{H_v\}_{v \in I} , \{ H_v \subset
L_v\}_{v \in I^c})$$ and the second projection is $$t_{w, n,2}\big((\mathcal{G}, \{H_v\}_{v \in I} , \{ H_v \subset L_v\}_{v \in I^c}, \mathcal{G} \rightarrow \mathcal{G}_n)\big)= (\mathcal{G}_n, \{H'_v\}_{v \in I} , \{ H'_v \subset L'_v\}_{v \in I^c})$$
where $H'_v$ and $H'_v \subset L'_v$ are the images of $H_v$ and $H_v \subset L_v$ in $\mathcal{G}_n$. 

The Hecke operator attached to this correspondence is $U_{w,1}^n$ (the $n$th iterate of $U_{w,1}$). More precisely, there are cohomological correspondences 
$$ (t_{w,n,1})_{*} t_{w,n,2}^{*} \omega^\kappa \rightarrow \omega^\kappa,~  (t_{w,n,1})_{*} t_{w,n,2}^{*} \omega^\kappa(-D) \rightarrow \omega^\kappa(-D).$$
Moreover, these cohomological correspondences restrict to 
$$ (t_{w,n,1})_{*} t_{w,n,2}^{*} (\omega^\kappa)^{++} \rightarrow p^{-3n}(\omega^\kappa)^{++},~(t_{w,n,1})_{*} t_{w,n,2}^{*} (\omega^\kappa(-D))^{++} \rightarrow p^{-3n}(\omega^\kappa(-D))^{++},$$
and they induce maps on cohomology.

\begin{lem}\label{lem-dynamic1}Let $w \in I$. Let  $x = (\mathcal{G}, \{H_v\}_{v \in I} , \{ H_v \subset L_v\}_{v \in I^c}, \mathcal{G} \rightarrow \mathcal{G}_1) \in C^{(1)}_{w}( \Spa (K, \ocal_K))$. \begin{enumerate}
\item If $v\in I $ and~$v\ne w$, we have $\deg H_v  = \deg H'_v$. 
\item If $v \in I^c $, we have $\deg L_v= \deg L'_v$.
\item We have $\deg H'_w \geq \deg H_w$, and in case of equality, $\deg H_w \in \{0,1\}$. 
\item $\deg H'_w = 1- \deg M_{1,w}/M_{1,w}[p]$.
\item $\deg M_{1,w}[p]/pM_{1,w} =1$, and $\deg pM_{1,w} \geq \deg M_{1,w}/M_{1,w}[p]$.
\item Let $\epsilon \geq 0$.  If $\deg M_w \leq 3-2 \epsilon$, then $\deg H'_w \geq \epsilon$. 
\end{enumerate}
\end{lem}
\begin{proof}Parts~(1) and~(2) follow because the maps $H_v\to H'_v$,
  $L_v\to L'_v$ are isomorphisms. The remaining parts are~ \cite[Lem.\
  14.3.1, Cor.\ 14.3.1]{pilloniHidacomplexes}.
\end{proof}

\begin{lem}\label{lem-dynamic2} Let $w \in I^c$. Let  $x = (\mathcal{G}, \{H_v\}_{v \in I} , \{ H_v \subset L_v\}_{v \in I^c}, \mathcal{G} \rightarrow \mathcal{G}_1) \in C^{(1)}_{w}( \Spa (K, \ocal_K))$. \begin{enumerate}
\item If $v \in I $, we have $\deg H_v  = \deg H'_v$. 
\item If $v \in I^c $ and $v \neq w $, we have $\deg L_v= \deg L'_v$,
\item We have $\deg L'_w \geq \deg L_w$, and in case of equality, $\deg L_w \in \{0,1, 2\}$. 
\item $\deg L'_w = 2- \deg M_{1,w}$. 
\end{enumerate}
\end{lem}
\begin{proof}Parts~(1) and~(2) follow as in
  Lemma~\ref{lem-dynamic1}. Parts~(3) and~(4) follow from~
  \cite[Prop.\ 2.3.1, 2.3.2, Lem.\ 2.3.4]{MR2783930} (and their proofs).
\end{proof}

\begin{cor}\label{coro-contraction} Let $w \in I^c$. Let $1 > \epsilon ' \geq \epsilon >0$. There exists $n \in \mathbb{Z}_{\geq 0}$ such that for all intervals $\prod_{v \neq w} J_v \subset [0,1]^{I} \times [0,2]^{I^c \setminus \{w\}}$, 

$$ U_{w,1}^n ( \mathcal{X}_{K_p(I) K^p}( \prod_{v \neq w} J_v \times [1 + \epsilon,2])) \subset  \mathcal{X}_{K_p(I) K^p}( \prod_{v \neq w} J_v \times [1 + \epsilon',2])$$
\end{cor} 
\begin{proof} This follows from Lemma~\ref{lem-dynamic2}~(3) and the
  maximum principle; see \cite[Prop.\ 2.3.6]{MR2783930}.
\end{proof}
\begin{cor} Let $w \in I$. Let $1 > \epsilon ' \geq \epsilon >0$. There exists $n \in \mathbb{Z}_{\geq 0}$ such that for all intervals $\prod_{v \neq w} J_v \subset [0,1]^{I\setminus \{w\}} \times [0,2]^{I^c }$, 

$$ U_{w,1}^n ( \mathcal{X}_{K_p(I) K^p}( \prod_{v \neq w} J_v \times [ \epsilon,1])) \subset  \mathcal{X}_{K_p(I) K^p}( \prod_{v \neq w} J_v \times [ \epsilon',1])$$
\end{cor} 
\begin{proof} This follows in the same way as
  Corollary~\ref{coro-contraction}, using Lemma~\ref{lem-dynamic1}~(3).
\end{proof}

\subsubsection{First analytic continuation result}

Let $J = \prod_{v |p} J_v \subset [0,1]^{I} \times [0,2]^{I^c}$ be a product of intervals.

\begin{lem} Let $w \in I^c$. Assume that $J_w = [ 2- \epsilon,2]$. The operator $U_{w,1}$ acts on ${\HH}^i(\mathcal{X}_{K^p
   K_p(I)}(J), \omega^\kappa)$. 
 \end{lem}  
   \begin{proof} In view of Lemma~ \ref{lem-dynamic2}~(3), the correspondence $C^{(1)}_{w}$ restricts to \[t_{w,1,2}: C^{(1)}_{w} \times_{t_{w,1,1},  \mathcal{X}_{K^p
   K_p(I)}} \mathcal{X}_{K^p
   K_p(I)}(J) \rightarrow \mathcal{X}_{K^p
   K_p(I)}(J). \qedhere\]
\end{proof} 
   
   We denote by ${\HH}^i(\mathcal{X}_{K_p(I) K^p}(J), \omega^\kappa)^{fs-U_{w,1}}$ the finite slope subspace for $U_{w,1}$. This is the subspace generated by classes which are annihilated by a polynomial in $U_{w,1}$ with non-zero constant term. 
   
   \begin{lem}\label{lem-extension-first} For all $1 > \epsilon' \geq  \epsilon >0$, the restriction map $${\HH}^i(\mathcal{X}_{K^p
   K_p(I)}(\prod_{v \neq w} J_v \times [ 2-\epsilon', 2]), \omega^\kappa)^{fs-U_{w,1}} \rightarrow {\HH}^i(\mathcal{X}_{K^p
   K_p(I)}(\prod_{v \neq w} J_v \times [ 2-\epsilon, 2]), \omega^\kappa)^{fs-U_{w,1}}  $$ is an isomorphism. 
   \end{lem}
   
   \begin{proof} Take $n$ as in Corollary~ \ref{coro-contraction}. Let $f \in {\HH}^i(\mathcal{X}_{K^p
   K_p(I)}(\prod_{v \neq w} J_v \times [ 2-\epsilon, 2]),
 \omega^\kappa)^{fs-U_{w,1}}  $ be a cohomology class. Let $P(X) = X^m
 + a_{m-1} X^{m-1} + \cdots + a_0$ be a polynomial with $a_0 \neq 0$
 such that $P(U_{w,1}) f = 0$. Therefore, if we set $Q(X) = -a_0^{-1}(
 P(X)-a_0 )$, we obtain that  $Q(U_{w,1}) f = f$. By iteration we get
 that $Q(U_{w,1})^n f = f$. 
The operator $$Q(U_{w,1})^n: {\HH}^i(\mathcal{X}_{K^p
   K_p(I)}(\prod_{v \neq w} J_v \times [ 2-\epsilon, 2]), \omega^\kappa) \rightarrow {\HH}^i(\mathcal{X}_{K^p
   K_p(I)}(\prod_{v \neq w} J_v \times [ 2-\epsilon, 2]), \omega^\kappa)$$ can be factored into:
   $${\HH}^i(\mathcal{X}_{K^p
   K_p(I)}(\prod_{v \neq w} J_v \times [ 2-\epsilon, 2]), \omega^\kappa) \stackrel{\widetilde{Q(U_{w,1})^n}}\rightarrow {\HH}^i(\mathcal{X}_{K^p
   K_p(I)}(\prod_{v \neq w} J_v \times [ 2-\epsilon', 2]), \omega^\kappa)$$ $$ \stackrel{res} \rightarrow {\HH}^i(\mathcal{X}_{K^p
   K_p(I)}(\prod_{v \neq w} J_v \times [ 2-\epsilon, 2]), \omega^\kappa),$$
   where the map~$\widetilde{Q(U_{w,1})^n}$ is the one coming from Corollary~  \ref{coro-contraction}.
   We therefore get an extension $\tilde{f}$ of $f$ to $${\HH}^i(\mathcal{X}_{K^p
   K_p(I)}(\prod_{v \neq w} J_v \times [ 2-\epsilon', 2]), \omega^\kappa)$$ by setting $\tilde{f} = \widetilde{Q(U_{w,1})^n} f$. 
   This proves the surjectivity of the map of the corollary. 

We now prove injectivity. Let $f, g \in {\HH}^i(\mathcal{X}_{K^p
   K_p(I)}(\prod_{v \neq w} J_v \times [ 2-\epsilon', 2]), \omega^\kappa)^{fs-U_{w,1}}$ be two classes having the same restriction to ${\HH}^i(\mathcal{X}_{K^p
   K_p(I)}(\prod_{v \neq w} J_v \times [ 2-\epsilon, 2]), \omega^\kappa)^{fs-U_{w,1}}$. We can find a polynomial $P$ as before such that $P(U_{w,1}) f = P(U_{w,1}) g = 0$. Therefore, using the same notation as before, we get that $Q(U_{w,1}) f = f$ and $Q(U_{w,1}) g = g$. 
   We can factor the operator $Q(U_{w,1})^n$ into:

   $${\HH}^i(\mathcal{X}_{K^p
   K_p(I)}(\prod_{v \neq w} J_v \times [ 2-\epsilon', 2]), \omega^\kappa) \stackrel{res}\rightarrow {\HH}^i(\mathcal{X}_{K^p
   K_p(I)}(\prod_{v \neq w} J_v \times [ 2-\epsilon, 2]), \omega^\kappa) $$ $$ \stackrel{\widetilde{Q(U_{w,1})^n}} \rightarrow {\HH}^i(\mathcal{X}_{K^p
   K_p(I)}(\prod_{v \neq w} J_v \times [ 2-\epsilon', 2]), \omega^\kappa)$$
   
Since  $res (f) = res( g )$, we deduce that $f=g$. 
   \end{proof}

The following two lemmas are the analogue of Lemma~\ref{lem-extension-first} for a place $w \in I$. The proofs are identical and left to the reader. 

\begin{lem} Let $w \in I$.  Assume that $J_w = [ 1- \epsilon,1]$. The operator $U_{w,2}$ acts on ${\HH}^i(\mathcal{X}_{K^p
   K_p(I)}(J), \omega^\kappa)$. 
 \end{lem}

   We denote by ${\HH}^i(\mathcal{X}_{K_p(I) K^p}(J), \omega^\kappa)^{fs-U_{w,2}}$ the finite slope subspace. This is the subspace generated by classes which are annihilated by a polynomial in $U_{w,2}$ with non-zero constant term. 
   
   \begin{lem} For all $1 > \epsilon \geq  \epsilon' >0$, the restriction map $${\HH}^i(\mathcal{X}_{K^p
   K_p(I)}(\prod_{v \neq w} J_v \times [ 1-\epsilon, 1]), \omega^\kappa)^{fs-U_{w,2}} \rightarrow {\HH}^i(\mathcal{X}_{K^p
   K_p(I)}(\prod_{v \neq w} J_v \times [ 1-\epsilon', 1]), \omega^\kappa)^{fs-U_{w,2}}  $$ is an isomorphism. 
   \end{lem}
   
\subsubsection{More analytic continuation results}
Let $w \in I^c$. Let $0 < \epsilon\le 1$. The cohomological correspondences: 

$$ (t_{w,n,1})_{*} (t_{w,n,2})^{*}(( \omega^\kappa)\vert_{\mathcal{X}_{K_p(I) K^p} ( \prod_{v \neq w} J_v \times [ 1+ \epsilon,2])} ) \rightarrow   \omega^\kappa\vert_{\mathcal{X}_{K_p(I) K^p} ( \prod_{v \neq w} J_v \times [ 1+ \epsilon,2])} $$
 and $$ (t_{w,n,1})_{*} (t_{w,n,2})^{*}(( \omega^\kappa)\vert_{\mathcal{X}_{K_p(I) K^p} ( \prod_{v \neq w} J_v \times [ 0,2])} ) \rightarrow   \omega^\kappa\vert_{\mathcal{X}_{K_p(I) K^p} ( \prod_{v \neq w} J_v \times [0,2])} $$ 
 
 can be related if we work with torsion coefficients.

\begin{prop} \label{prop-hecke-corresp} Let $0 < \epsilon < \epsilon'$. There is a factorization of the Hecke correspondence~$U_{w,1}^n$:
\[
\begin{tikzpicture}[
node distance = 10mm and 1mm,
every edge/.style = {draw, -{Straight Barb[scale=0.8]}, semithick}
                    ]
                    \node (b) {
$ ( \omega^\kappa/ p^{n(l_w(1-\epsilon')-3)}(\omega^\kappa)^{++})\vert_{\mathcal{X}_{K_p(I) K^p} ( \prod_{v \neq w} J_v \times [ 1+ \epsilon,2])} $
};
\node (a)  [below left=of b.south] {
$(t_{w,n,1})_{*} (t_{w,n,2})^{*}(( \omega^\kappa)^{++}\vert_{\mathcal{X}_{K_p(I) K^p} ( \prod_{v \neq w} J_v \times [ 1+ \epsilon,2])} )$
};
\node (f) [below of=a]{};
\node (g) [below of=b]{};
\node (hh) [below of=g]{};
\node (h) [below of=hh]{};
\node (ee) [below of=f]{};
\node (e) [below of=ee]{};
\node (c) [below of=e]{
$(t_{w,n,1})_{*} (t_{w,n,2})^{*}(( \omega^\kappa)^{++}\vert_{\mathcal{X}_{K_p(I) K^p} ( \prod_{v \neq w} J_v \times [ 0,2])} ) $
};
\node (d) [below of=h] {
$ ( \omega^\kappa/ p^{n(l_w(1-\epsilon')-3)}(\omega^\kappa)^{++})\vert_{\mathcal{X}_{K_p(I) K^p} ( \prod_{v \neq w} J_v \times [0,2])} $
};
     	\draw[-to] (a) -- (b);
           \draw[-to] (c) -- (d);
           \draw[-to] (c) -- (a);
           \draw[-to] (d) -- (b);
           \draw[-to] (a) --  (d);
 \end{tikzpicture}
\]
\end{prop}
\begin{proof} Let $x\in \mathcal{X}_{K_p(I) K^p} ( \prod_{v \neq w} J_v \times [ 0,2]) )$. We have to find a neighbourhood $U$ of $x$ and to construct a canonical map 
$$ (t_{w,n,2})^{*}(( \omega^\kappa)^{++}\vert_{\mathcal{X}_{K_p(I) K^p} ( \prod_{v \neq w} J_v \times [ 1+ \epsilon,2])} )(t_{w,n,1}^{-1}(U)) $$ $$\rightarrow ( \omega^\kappa/ p^{n(l_w(1-\epsilon')-3)}(\omega^\kappa)^{++})\vert_{\mathcal{X}_{K_p(I) K^p} ( \prod_{v \neq w} J_v \times [0,2])}(U).$$

Pick $\epsilon'' \in ( \epsilon, \epsilon')$ such that for all $y \in t_{w,n,1}^{-1}(x)$ corresponding to a subgroup $M_{w,n} \subset \mathcal{G}$, we have $\vert \delta_{M_{w,n}} \vert_y \neq \vert p^{n(1-\epsilon'')} \vert_y$. 
It follows that there exists an open neighbourhood $U$ of $x$ and a
disjoint decomposition $t_{w,n,1}^{-1} (U) = V \coprod W$, such that
for all $y \in W$, we have $\vert \delta_{M_{w,n}} \vert_y  > \vert
p^{n(1-\epsilon'')} \vert_y$, and for all $y \in V$, we have $\vert \delta_{M_{w,n}} \vert_y  < \vert p^{n(1-\epsilon'')} \vert_y$.

The cohomological correspondence is a map: 
$$ t_{w,n,2}^{*} (\omega^\kappa)^{++} ( V) \oplus  t_{w,n,2}^{*} (\omega^\kappa)^{++} ( W) \rightarrow \omega^\kappa(U).$$

The image of $t_{w,n,2}^{*} (\omega^\kappa)^{++} ( V)$ lands in
$p^{n(l_v(1- \epsilon'')-3)}(\omega^\kappa)^{++}(U)$ (see~\cite[Lem.\ 14.6.1]{pilloniHidacomplexes}). Therefore, we have a factorization: $$ t_{w,n,2}^{*} (\omega^\kappa)^{++} ( t_{w,n,1}^{-1}(U) ) \rightarrow   t_{w,n,2}^{*} (\omega^\kappa)^{++} ( W) \rightarrow \omega^\kappa/ p^{n(l_v(1- \epsilon'')-3)}(\omega^\kappa)^{++}(U).$$
On the other hand, we claim that $t_{w,n,2} (W) \subset \mathcal{X}_{K_p(I) K^p} (
\prod_{v \neq w} J_v \times [ 1+ \epsilon,2]) ) $.

Indeed, let $x' \in U$  and let $y' \in t_{w,n,1}^{-1}(\{x\})$. Without loss of generality, we may assume that $x'$ and $y'$ are rank one points. Let us define $M_{w,i} = M_{w,n}\vert_{y'}[p^i]$ for all $1 \leq i \leq n$. Then we have a sequence of isogenies:
$$ \cG_w\vert_{x'} \rightarrow \cG_w\vert_{x'}/M_{w,1} \rightarrow \cdots \rightarrow \cG_w\vert_{x'}/M_{w,n}$$
We let $L_{w,i}$ be the image of $L_w\vert_{x'}$ in  $\cG_w/M_{w,i}$
(so that $ L_{w,n} = L_{w}\vert_{t_{w,n,2}(y')}$). Then, by Lemma
\ref{lem-dynamic2}~(4), we have that $\deg L_{w,n} = 2- \deg
M_{w,n}/M_{w,n-1}$. On the other hand, for all $1 \leq i \leq n-1$,
the map $p: M_{w,i+1}/M_{w,i} \rightarrow M_{w,i}/M_{w,i-1}$ is a
generic isomorphism. It follows that $\deg  M_{w,i+1}/M_{w,i} \leq
\deg  M_{w,i}/M_{w,i-1}$. Since $\deg M_{w,n} = \sum_{i=1}^{n} \deg
M_{w,i}/M_{w,i-1}$ and $\deg M_{w,n} \leq n(1- \epsilon'')$, we deduce
that $ \deg M_{w,n}/M_{w,n-1} \leq 1- \epsilon''$ and therefore $\deg
L_{w,n} \geq 1 + \epsilon''>1+\epsilon$, as required. 

We can therefore produce the expected map as the composition: 
 $$ (t_{w,n,2})^{*}((
 \omega^\kappa)^{++}\vert_{\mathcal{X}_{K_p(I) K^p} ( \prod_{v \neq w}
   J_v \times [ 1+ \epsilon,2])}
 )(t_{w,n,1}^{-1}(U))  $$ $$\rightarrow (t_{w,n,2})^{*}((
 \omega^\kappa)^{++}\vert_{\mathcal{X}_{K_p(I) K^p} ( \prod_{v \neq w}
   J_v \times [ 1+ \epsilon,2])} )(W ) =  t_{w,n,2}^{*}
 (\omega^\kappa)^{++} ( W) $$ \[ \rightarrow \omega^\kappa/
   p^{n(l_v(1- \epsilon')-3)}(\omega^\kappa)^{++}(U). \qedhere \]
\end{proof} 
\begin{cor}\label{coro-extension-proj-lim}Let $P = X^m + a_{m-1}
  X^{m-1} + \cdots + a_0$ be a polynomial, with the property that all
  the roots $a$ of $P$ satisfy $v(a) < l_w-3$. Then there is a map  \[\ext:
  {\HH}^i(\mathcal{X}_{K_p(I) K^p} ( \prod_{v \neq w} J_v \times [ 1+
  \epsilon,2]), \omega^{\kappa})[P(U_{w,1})=0] \to\]\[ \left(\lim_n
    {\HH}^i ( \mathcal{X}_{K_p(I) K^p} ( \prod_{v \neq w} J_v \times [
    0,2]),  \omega^\kappa/ p^n
    (\omega^{\kappa})^{++})\right)[P(U_{w,1})=0]\] such that the
composite of $\ext$ followed by restriction to~$\mathcal{X}_{K_p(I)Kp} (
\prod_{v \neq w} J_v \times [ 1+\epsilon,2])$ is the natural map
induced by  $ \omega^\kappa\to\omega^\kappa/ p^n
    (\omega^{\kappa})^{++}$.

Furthermore, the composite of the restriction map \[{\HH}^i(\mathcal{X}_{K_p(I) K^p} ( \prod_{v \neq w} J_v \times [0,2]), \omega^{\kappa})[P(U_{w,1})=0]\to \]\[{\HH}^i(\mathcal{X}_{K_p(I) K^p} ( \prod_{v \neq w} J_v \times [ 1+
  \epsilon,2]), \omega^{\kappa})[P(U_{w,1})=0]\] followed by~$\ext$
is the natural map
induced by  $ \omega^\kappa\to\omega^\kappa/ p^n
    (\omega^{\kappa})^{++}$.
\end{cor}

\begin{proof} Let $\epsilon'>0$ be such that  for all roots $a$ of $P$,
  we have $l_w(1- \epsilon') - 3 > v(a)$. Let $\alpha = \inf_a
  \{l_w(1-\epsilon')-3-v(a)\}$ (so that in particular $\alpha>0$). By
  Lemma~ \ref{lem-extension-first}, we can assume that $0 < \epsilon <
  \epsilon'$. Suppose that $f\in {\HH}^i(\mathcal{X}_{K_p(I) K^p} ( \prod_{v \neq w} J_v \times [ 1+
  \epsilon,2]), \omega^{\kappa})$ satisfies
  $P(U_{w,1})f =0$. By rescaling $f$, we can and do also assume that  $f \in
  {\HH}^i( \mathcal{X}_{K_p(I) K^p}( \prod_{v \neq w} J_v \times [
  1+\epsilon,2]), (\omega^\kappa)^{++})$. Let $Q(X) = -a_0^{-1} (P(X)-a_0)$ so that $Q(U_{w,1})f = f$. 

 Since~$Q(U_{w,1})$ can be written as a sum of products of
 the~$\frac{1}{a}U_{w,1}$, where~$a$ runs over the roots of~$P$, it
 follows from  Proposition~ \ref{prop-hecke-corresp} that the map 
$$Q(U_{w,1})^n:  {\HH}^i( \mathcal{X}_{K_p(I) K^p}( \prod_{v \neq w} J_v \times [ 1+\epsilon,2]), (\omega^\kappa)^{++}) \rightarrow  {\HH}^i( \mathcal{X}_{K_p(I) K^p}( \prod_{v \neq w} J_v \times [ 1+\epsilon,2]), \omega^\kappa) $$ $$\rightarrow {\HH}^i( \mathcal{X}_{K_p(I) K^p}( \prod_{v \neq w} J_v \times [ 1+\epsilon,2]), (\omega^\kappa)/p^{n\alpha} (\omega^\kappa)^{++})$$ can actually be factored into:
\[
\begin{tikzpicture}[
node distance = 14mm and 1mm,
every edge/.style = {draw, -{Straight Barb[scale=0.8]}, semithick}
                    ]
\node (a) {
$ {\HH}^i( \mathcal{X}_{K_p(I) K^p}( \prod_{v \neq w} J_v \times [ 1+\epsilon,2]), (\omega^\kappa)^{++})$
};
\node (b) [below right=of a.south] {
$  {\HH}^i( \mathcal{X}_{K_p(I) K^p}( \prod_{v \neq w} J_v \times [ 0,2]), (\omega^\kappa)/p^{n\alpha} (\omega^\kappa)^{++})   $
};
\node (f) [below of=a]{};
\node (g) [below of=b]{};
\node (h) [below of=g]{};
\node (e) [below of=f]{};
\node (c) [below of=e]{
$   {\HH}^i( \mathcal{X}_{K_p(I) K^p}( \prod_{v \neq w} J_v \times [ 1+\epsilon,2]), (\omega^\kappa)/p^{n\alpha} (\omega^\kappa)^{++}) $
};
\draw[-to] (a) -- node[above] {\footnotesize{$\quad \qquad \widetilde{Q(U_{w,1})^n}$}} (b);
\draw[-to] (a) --  (c);
\draw[-to] (b) -- (c);
 \end{tikzpicture}
\]
We define sections $f_n \in {\HH}^i ( \mathcal{X}_{K_p(I) K^p} (
\prod_{v \neq w} J_v \times [ 0,2]),  \omega^\kappa/ p^{n\alpha}
(\omega^{\kappa})^{++})$ by $f_n = \widetilde{Q(U_{w,1})^n}(f)$. It follows
from the definitions   that $f_n = f_{n-1}$ in 
$${\HH}^i ( \mathcal{X}_{K_p(I) K^p} ( \prod_{v \neq w} J_v \times [
0,2]),  \omega^\kappa/ p^{(n-1)\alpha- m(l_w(1-\epsilon')-\alpha)}
(\omega^{\kappa})^{++})$$
 and that $Q(U_{w,1}) f_n = f_n $ in 
 $${\HH}^i
( \mathcal{X}_{K_p(I) K^p} ( \prod_{v \neq w} J_v \times [ 0,2]),
\omega^\kappa/ p^{(n\alpha- m(l_w(1-\epsilon')-\alpha)} (\omega^{\kappa})^{++})$$
 (see
the proof of ~\cite[Cor.\ 14.6.1]{pilloniHidacomplexes} for a similar
verification). We let~$\ext(f)$ be the projective system given by
the~$f_n$.

It remains to check that if~$f$ is the  restriction of a class in
$${\HH}^i( \mathcal{X}_{K_p(I) K^p}( \prod_{v \neq w} J_v \times [
0,2]), \omega^\kappa),$$
 then the~$f_n$ are obtained from the natural
map $ \omega^\kappa\to\omega^\kappa/ p^n
    (\omega^{\kappa})^{++}$. This follows easily from the factorization
    \[
\begin{tikzpicture}[
node distance = 15mm and 1mm,
every edge/.style = {draw, -{Straight Barb[scale=0.8]}, semithick}
                    ]
\node (a) {
$ {\HH}^i( \mathcal{X}_{K_p(I) K^p}( \prod_{v \neq w} J_v \times [0,2]), (\omega^\kappa)^{++})$
};
\node (b) [below right=of a.south] {
$  {\HH}^i( \mathcal{X}_{K_p(I) K^p}( \prod_{v \neq w} J_v \times [0,2]), (\omega^\kappa)/p^{n\alpha} (\omega^\kappa)^{++})  $
};
\node (f) [below of=a]{};
\node (g) [below of=b]{};
\node (h) [below of=g]{};
\node (e) [below of=f]{};
\node (c) [below of=e]{
$   {\HH}^i( \mathcal{X}_{K_p(I) K^p}( \prod_{v \neq w} J_v \times [ 1+\epsilon,2]), (\omega^\kappa)/p^{n\alpha} (\omega^\kappa)^{++}) $
};
\draw[-to] (a) -- node[above] {\footnotesize{$\quad \qquad {Q(U_{w,1})^n}$}} (b);
\draw[-to] (a) --  (c);
\draw[-to] (c) -- node[below] {\footnotesize{$\quad \qquad
    \widetilde{Q(U_{w,1})^n}$}}   (b); 
 \end{tikzpicture}
 \]
\end{proof}

The next proposition and corollary are the analogue of
the above results for a place $w \in I$. The proofs are virtually
identical to the above, and are left to the reader (or look at
\cite[Prop.\ 14.6.1, Cor.\ 14.6.1]{pilloniHidacomplexes}). 

\begin{prop} Let~$w\in I$, and let $0 < \epsilon < \epsilon'$. There is a factorization of the Hecke correspondence~$U_{w,2}^n$:
\[
\begin{tikzpicture}[
node distance = 10mm and 1mm,
every edge/.style = {draw, -{Straight Barb[scale=0.8]}, semithick}
                    ]
                    \node (b) {
$\qquad ( \omega^\kappa/ p^{n(l_w+k_w-3-2\epsilon' k_v)}(\omega^\kappa)^{++})\vert_{\mathcal{X}_{K_p(I) K^p} ( \prod_{v \neq w} J_v \times [\epsilon, 1])}  $
};
\node (a)  [below left=of b.south] {
$(t_{w,n,1})_{*} (t_{w,n,2})^{*}(( \omega^\kappa)^{++}\vert_{\mathcal{X}_{K_p(I) K^p} ( \prod_{v \neq w} J_v \times [\epsilon,1])} )$
};
\node (f) [below of=a]{};
\node (g) [below of=b]{};
\node (hh) [below of=g]{};
\node (h) [below of=hh]{};
\node (ee) [below of=f]{};
\node (e) [below of=ee]{};
\node (c) [below of=e]{
$(t_{w,n,1})_{*} (t_{w,n,2})^{*}(( \omega^\kappa)^{++}\vert_{\mathcal{X}_{K_p(I) K^p} ( \prod_{v \neq w} J_v \times [1, 0])} ) $
};
\node (d) [below of=h] {
$\qquad ( \omega^\kappa/ p^{n(l_w+k_w-3-2\epsilon' k_w)}(\omega^\kappa)^{++})\vert_{\mathcal{X}_{K_p(I) K^p} ( \prod_{v \neq w} J_v \times [0,1])} $
};
     	\draw[-to] (a) -- (b);
           \draw[-to] (c) -- (d);
           \draw[-to] (c) -- (a);
           \draw[-to] (d) -- (b);
           \draw[-to] (a) --  (d);
 \end{tikzpicture}
\]
\end{prop}

\begin{cor}\label{coro-extension-proj-lim2} Let $w\in I$, and let $1> \epsilon \geq 0$. Let $P = X^m + a_{m-1}
  X^{m-1} + \cdots + a_0$ be a polynomial, with the property that all
  the roots $a$ of $P$ satisfy $v(a) <k_w+ l_w-3$. 
 Then there is a map  \[\ext:
  {\HH}^i(\mathcal{X}_{K_p(I) K^p} ( \prod_{v \neq w} J_v \times [ \epsilon,1]), \omega^{\kappa})[P(U_{w,2})=0] \to\]\[ \left(\lim_n
    {\HH}^i ( \mathcal{X}_{K_p(I) K^p} ( \prod_{v \neq w} J_v \times [
    0,1]),  \omega^\kappa/ p^n
    (\omega^{\kappa})^{++})\right)[P(U_{w,2})=0]\] such that the
composite of $\ext$ followed by restriction to~$\mathcal{X}_{K_p(I)Kp} (
\prod_{v \neq w} J_v \times [\epsilon,1])$ is the natural map
induced by  $ \omega^\kappa\to\omega^\kappa/ p^n
(\omega^{\kappa})^{++}$.

Furthermore, the composite of the restriction map \[{\HH}^i(\mathcal{X}_{K_p(I) K^p} ( \prod_{v \neq w} J_v \times [0,1]), \omega^{\kappa})[P(U_{w,2})=0]\to \]\[{\HH}^i(\mathcal{X}_{K_p(I) K^p} ( \prod_{v \neq w} J_v \times [
  \epsilon,1]), \omega^{\kappa})[P(U_{w,2})=0]\] followed by~$\ext$
is the natural map
induced by  $ \omega^\kappa\to\omega^\kappa/ p^n
    (\omega^{\kappa})^{++}$.
\end{cor}

\subsubsection{Proof of the main classicality theorem}

Let $S \subset S_p$ be a subset. Let $J(S, \epsilon) = \prod_{v \in S
  \cap I} [0,1] \times \prod_{v \in S \cap I^c} [0,2] \times \prod_{
  S^c\cap I} [ \epsilon, 1]   \times  \prod_{ S^c\cap I^c} [ 1+
\epsilon, 2] $. We say that a cohomology class  $f \in {\HH}^i(
\mathcal{X}_{K_p(I) K^p}( J(S, \epsilon),\omega^\kappa)$ is of finite
slope if for all $v |p$, there is a polynomial $P_v$ all of whose
roots are nonzero, such that:
\begin{itemize} 
\item if $v \in I^c$, $P_v(U_{v,1}) f  = 0$,
\item if $v \in I$, $P_v(U_{v,2}) f = 0$. 
\end{itemize}

\begin{lem}\label{lem: isomorphic to projective limit on finite slope part} The canonical map 
$$  {\HH}^i( \mathcal{X}^{G_1}_{K_p(I) K^p}( J(S, \epsilon)), \omega^\kappa) \rightarrow   \lim {\HH}^i( \mathcal{X}^{G_1}_{K_p(I) K^p}( J(S, \epsilon)), \omega^\kappa/p^n(\omega^{\kappa})^{++})$$
is surjective and induces an isomorphism on the finite slope part. \end{lem} 

  \begin{proof} The surjectivity follows from~ \cite[Prop.\
    3.2.1]{pilloniHidacomplexes}. The injectivity can be proved in
    exactly the same way as
    \cite[Lem.\ 14.7.1]{pilloniHidacomplexes}.  We have put the superscript $G_1$ because we need some finiteness property to deduce the injectivity. 
 \end{proof}

  \begin{lem}\label{lem: restriction is an isomorphism on finite
      slope}Choose polynomials~$P_v$ such that
    \begin{itemize}
    \item if $v \in I^c$, all the roots $a$ of $P_v$ satisfy
      $v(a) < l_v-3$, and
    \item if $v \in I$, all the roots $a$ of $P_v$ satisfy
      $v(a) < k_v +l_v-3$.
    \end{itemize}
Write~$U_v=U_{v,1}$ if $v\in I^c$, and $U_v=U_{v,2}$ if
    $v\in I$.
    If~$S\subset T$, then the natural restriction map
    $$
    \begin{diagram}
      {\HH}^i(
      \mathcal{X}^{G_1}_{K_p(I) K^p}( J(T, \epsilon)),
      \omega^\kappa)[P_v(U_v)=0]_{v\in S_p} \\
      \dTo \\
       {\HH}^i(
      \mathcal{X}^{G_1}_{K_p(I) K^p}( J(S, \epsilon)),
      \omega^\kappa)[P_v(U_v)=0]_{v\in S_p}
      \end{diagram}
      $$
       is an isomorphism.
  \end{lem}
  \begin{proof}By induction, it is enough to treat the
    case~$T=S\cup\{w\}$ for some~$w$. The result then follows from
    Lemma~\ref{lem: isomorphic to projective limit on finite slope
      part} (applied to both~$S$ and~$T$), together with Corollary~ \ref{coro-extension-proj-lim} and Corollary~ \ref{coro-extension-proj-lim2}.    
  \end{proof}

  \begin{proof}[Proof of Theorem~\ref{thm-classicality}]This follows immediately
    from Lemma~\ref{lem: restriction is an isomorphism on finite
      slope}, applied with the choices~$S=\emptyset$ and~$T=S_p$.    
  \end{proof}
\subsection{Application to ordinary cohomology}\label{subsec: ordinary
cohomology is overconvergent} In this section we
study the case $\# I = 1$, where we are able to relate the Hida
complexes constructed in~\S\ref{sec: Hida complexes} to the
overconvergent cohomology considered in this section. Our first result
is the following, which shows in particular that in this case the
ordinary classes in~$\HH^1$ are overconvergent. The proof can be
viewed as a generalization of the familiar argument for~$\GL_2$ which
shows that ordinary $p$-adic modular forms are overconvergent
(see~\cite[Lem.\ 1]{MR1709306}), by using the continuity of the
ordinary projector to the finite-dimensional space of ordinary forms.

Recall that we defined the complex~$M_I$ in
Theorem~\ref{theorem-p-adic-complex}. By
Theorem~\ref{theorem-p-adic-complex}~(3), for all classical algebraic
  weights $\kappa$ with $l_v=2$ for~$v\in I$ and $k_v\equiv l_v\equiv
  2 \pmod{p-1}$ for all~$v|p$ we have
\[ M_I \otimes^{\mathbf{L}}_{\Lambda_I, \kappa}
\whalingship_p = e(U^I)\mathrm{R} \Gamma(\mathcal{X}^{G_1,\mult}_{K_p(I) K^p},
\omega^{\kappa}(-D)). \]

 \begin{prop}\label{prop-ho-h1}Suppose that~$\# I=1$. For all classical algebraic
  weights $\kappa$ with $l_v=2$ for~$v\in I$ and $k_v\equiv l_v\equiv
  2 \pmod{p-1}$ for all~$v|p$, the restriction map $$ e(U^I) \mathrm{R} \Gamma ( \mathcal{X}^{G_1, \mult, \dag}_{K_p(I) K^p},
\omega^{\kappa}(-D)) \rightarrow M_I \otimes^{\mathbf{L}}_{\Lambda_I, \kappa}
\whalingship_p$$ induces an injective map on ${\HH}^0$ and a surjective map on
${\HH}^1$. \end{prop} 
  
\begin{proof}The injectivity of the map on~$\HH^0$s is clear. In the
  case $F=\Q$, the surjectivity of the map on~$\HH^1$s is proved in ~\cite[Lem.\
    14.8.2]{pilloniHidacomplexes}; we now recall this argument in our
    setting. Let
    $\pi:\cX_{K_p(I) K^p}\to\cX^*_{K_p(I) K^p}$ be the projection to the
    minimal compactification; as usual, we have
    $\mathrm{R}^i\pi_{*} \omega^\kappa(-D)=0$ for~$i>0$. Let
    $\cX^{*,\mult}_{K_p(I) K^p}$ be the image
    of~$\cX_{K_p(I) K^p}^{\mult}$ (the rigid analytic generic fibre
    of~$\mathfrak{X}^{I}_{K_p(I) K^p}$) in the
    minimal compactification; it admits an affinoid
    cover~$\cX^{*,\mult}_{K_p(I) K^p}=U_1\cup U_2$. 

Then the complex $\mathrm{R} \Gamma ( \mathcal{X}^{\mult}_{K_p(I) K^p},
\omega^{\kappa}(-D))$  is represented by the
complex \[\HH^0(U_1,\omega^\kappa(-D))\oplus
  \HH^0(U_2,\omega^\kappa(-D))\to \HH^0(U_1\cap
  U_2,\omega^\kappa(-D)). \]The terms of this complex are Banach
spaces; a norm giving their topology is provided by taking an
appropriate formal model. The topology on $\HH^1( \mathcal{X}^{\mult}_{K_p(I) K^p},
\omega^{\kappa}(-D))$ is the induced quotient topology, which
coincides with the topology obtained by declaring that $\HH^1( \mathfrak{X}^{I}_{K_p(I) K^p},
\omega^{\kappa}(-D))$ is open and bounded.

The complex $\mathrm{R} \Gamma ( \mathcal{X}^{\mult, \dag}_{K_p(I) K^p},
\omega^{\kappa}(-D))$ is represented by the subcomplex of
overconvergent sections \[\HH^0(U_1,\omega^{\kappa,\dagger}(-D))\oplus
  \HH^0(U_2,\omega^{\kappa,\dagger}(-D))\to \HH^0(U_1\cap
  U_2,\omega^{\kappa,\dagger}(-D)), \]where by definition for an
open~$U$, \[\omega^{\kappa,\dagger}(U):=\colim_{V\supset\overline{U}}\omega^\kappa(V) \]where~$V$
runs over the strict neighbourhoods of~$U$.

It follows in particular that the map \[\HH^1( \mathcal{X}^{\mult, \dag}_{K_p(I) K^p},
\omega^{\kappa}(-D))\to \HH^1(\mathcal{X}^{\mult}_{K_p(I) K^p},
\omega^{\kappa}(-D)) \]has dense image. The operator~$U^I$ is
continuous on~$\HH^1(\mathcal{X}^{\mult}_{K_p(I) K^p},
\omega^{\kappa}(-D))$   (consider the action
of~$U^I$ on $\HH^1( \mathfrak{X}^{I}_{K_p(I) K^p},
\omega^{\kappa}(-D))$), so the projection \[\HH^1(\mathcal{X}^{G_1, \mult}_{K_p(I) K^p},
\omega^{\kappa}(-D))\to e(U^I)\HH^1(\mathcal{X}^{G_1, \mult}_{K_p(I) K^p},
\omega^{\kappa}(-D))\] is also continuous (we have introduced the superscript $G_1$ to make sure that the projector $e(U^I)$ is well defined, the passage from the cohomology of $\mathcal{X}^{ \mult}_{K_p(I) K^p}$ to the cohomology of $\mathcal{X}^{G_1, \mult}_{K_p(I) K^p}$ is given by a projector so all density  statements are preserved). It follows that the induced
map \[e(U^I)\HH^1( \mathcal{X}^{G_1, \mult, \dag}_{K_p(I) K^p},
\omega^{\kappa}(-D))\to e(U^I)\HH^1(\mathcal{X}^{G_1, \mult}_{K_p(I) K^p},
\omega^{\kappa}(-D))\] has dense image. But the target is a
finite-dimensional Banach space over~$\whalingship_p$, so its topology is the
unique one extending that on~$\whalingship_p$, and in particular it contains no
proper dense subspaces, so we are done.
\end{proof}
\begin{prop}\label{prop: I equal 1 Euler char}Suppose that~$\# I\le 1$.  For all classical algebraic weights
  $\kappa$ with $l_v=2$ for~$v\in I$, we have the equality of  Euler characteristics: $$ EC( e(U^I) \mathrm{R} \Gamma ( \mathcal{X}^{G_1, \mult, \dag}_{K_p(I) K^p}, \omega^{\kappa}(-D)) ) = EC(M_I \otimes^{\mathbf{L}}_{\Lambda_I, \kappa} \whalingship_p).$$
\end{prop}

\begin{proof}By Theorem~\ref{theorem-p-adic-complex} and Lemma~\ref{lem: ordinary mult cohomology is perfect}, both complexes are perfect complexes in
  degrees~$[0,1]$.  By Corollary~ \ref{cor-boundslope} and
  Proposition~\ref{prop: cuspidal cohomology in degrees 0 I}, we have that  $$  EC(e(U^I) {\HH}^i_{cusp}(G_1,
  \dag, \kappa)) \leq  EC( e(U^I) \mathrm{R} \Gamma ( \mathcal{X}^{G_1, \mult, \dag}_{K_p(I) K^p}, \omega^{\kappa}(-D)) )$$ and the inequality is an equality if $k_v-l_v \geq 3$ for all $v| p$.  By  Proposition \ref{prop-ho-h1}, we have that $$EC( e(U^I) \mathrm{R} \Gamma ( \mathcal{X}^{G_1, \mult, \dag}_{K_p(I) K^p}, \omega^{\kappa}(-D)) ) \leq EC(M_I \otimes^{\mathbf{L}}_{\Lambda_I, \kappa} \whalingship_p).$$Consequently, it suffices to prove that $$EC(e(U^I) {\HH}^i_{cusp}(G_1,
  \dag, \kappa)) \geq EC(M_I \otimes^{\mathbf{L}}_{\Lambda_I, \kappa} \whalingship_p).$$   By Theorem~\ref{thm: local constancy of Euler char}, and Theorem~\ref{theorem-p-adic-complex},  both Euler characteristics under
  consideration are locally constant functions of $\kappa$. It therefore suffices to
  prove the statement when $l_v \geq C$ for all $v \in I^c$,  and $k_v-l_v \geq C$ for all $v | p$.  In this range of weights we can compare these cohomology to classical cohomology. 

It follows from Theorem~ \ref{thm-classicality}
that  $$ e(U^I) \mathrm{R} \Gamma(
\mathcal{X}_{K_p(I) K^p}^{G_1, \mult, \dag}, \omega^{\kappa}(-D)) )=e(U^I) \mathrm{R} \Gamma ( \mathcal{X}^{G_1}_{K_p(I) K^p},
\omega^{\kappa}(-D)) ).$$ We claim that the natural map \[ \mathrm{R} \Gamma ( \mathcal{X}^{G_1}_{K^pK_p},
\omega^{\kappa}(-D)) )\to \mathrm{R} \Gamma ( \mathcal{X}^{G_1}_{K_p(I) K^p},
\omega^{\kappa}(-D)) )\] induces a quasi-isomorphism \[ e(T^I)\mathrm{R} \Gamma ( \mathcal{X}^{G_1}_{K^pK_p},
\omega^{\kappa}(-D)) )\to e(U^I)\mathrm{R} \Gamma ( \mathcal{X}^{G_1}_{K_p(I) K^p},
\omega^{\kappa}(-D)) ).\] Indeed, by Theorem~\ref{thm: cohomology in terms of automorphic forms}~(2), the cohomology groups on each side can be
  computed in terms of automorphic representations, and the claim
  follows  from Proposition~\ref{prop: ordinary eigenform in autrep}
  as explained in Remark \ref{rem-notsoeasy} below. 

Now, it follows from Theorem~ \ref{theorem-p-adic-complex} that the map
$e(T^I) \mathrm{R} \Gamma ( \mathcal{X}^{G_1}_{K^pK_p},
\omega^{\kappa}(-D)) \rightarrow  M_I \otimes^{\mathbf{L}}_{\Lambda_I, \kappa}
\whalingship_p$ is an isomorphism on ${\HH}^0$ and is injective on
${\HH}^1$. Putting this all together, 
the proposition follows.  \end{proof} 

\begin{rem}\label{rem-notsoeasy} Let us point out a subtle point in the proof of Proposition~\ref{prop: I equal 1 Euler char}. In
  order to use Proposition~\ref{prop: ordinary eigenform in autrep}
  one needs to check that for any $v | p$, any representation
  $\pi_v$ of $\mathrm{GSp}_4(\ocal_{F_v})$ contributing to either
  $e(T^I)\mathrm{R} \Gamma ( \mathcal{X}^{G_1}_{K^pK_p},
  \omega^{\kappa}(-D)) )$ or
  $e(U^I)\mathrm{R} \Gamma ( \mathcal{X}^{G_1}_{K_p(I) K^p},
  \omega^{\kappa}(-D)) )$ is ordinary. For all places $v \in I^c$,
  this is true essentially by definition since the two Hecke operators
  at $v|p$ occur in the projector. For places $v \in I$, this is
  a bit more subtle since only one operator $T_v$ or $U_{v,2}$ is
  involved in the definition of the projector. The
  $U_{v,2}$-ordinarity of a local representation $\pi_v$ with Hecke
  parameters
  \[[\alpha_vp^{1-k_v/2},\beta_vp^{1-k_v/2},\beta_v^{-1}p^{k_v/2},\alpha_v^{-1}p^{k_v/2}]\]
  implies that $\alpha_v \beta_v$ is a $p$-adic unit.  Ordinarity
  means that $\alpha_v$ and $\beta_v$ are both $p$-adic units.  This
  is implied by $U_{v,2}$-ordinarity if we assume the Katz--Mazur
  inequality which says the the Newton polygon is above the Hodge
  polygon with the same initial and terminal point.  Indeed, in our
  case, the Katz--Mazur inequality translates into the condition that
  $\alpha_v$ and $\beta_v$ are $p$-adic integers. 

However, this inequality is
  subtle at non-cohomological weights.  For $F = \Q$  the
  Katz--Mazur inequality for $H^0$ and $H^1$ classes is proved in
  \cite[Prop.\ 14.9.1]{pilloniHidacomplexes}, and the argument
  generalizes without difficulty to our case. We also remark that for classes in the
  $H^0$, we can use eigenvarieties to deduce that the Katz--Mazur
  inequality holds in non-cohomological weights because it holds at
  cohomological weights. A similar argument will apply for classes in
  the $H^1$ once eigenvarieties are constructed for $H^1$ cohomology
  classes. Note that alternatively we could
  force the Katz--Mazur inequality by localizing further at certain
 $p$-adically integral eigenvalues of the operators $T_{v,1}$ and $U_{\Kli(v),1}$,
  and in fact such a localization will be in force in the rest of the
  paper. We could also directly deduce the Katz--Mazur inequality for
  classes in the $H^1$ from the corresponding inequality for classes
  in the $H^0$ after making a non-Eisenstein localization (because
  after making such a localization, the Euler characteristic vanishes). Such a localization will
  also be in force in the rest of the paper. In view of this, we do
  not spell out the details of the  generalization of~\cite[Prop.\
  14.9.1]{pilloniHidacomplexes} to our setting.
\end{rem} 

\begin{thm}\label{thm: ordinary is overconvergent I equal 1}Suppose
  that~$\# I\le 1$. For all   classical algebraic weights $\kappa$
  with~$l_v=2$ for~$v\in I$, we have a canonical  isomorphism: 
    $$e(U^I) \mathrm{R} \Gamma ( \mathcal{X}^{G_1, \mult, \dag}_{K_p(I) K^p}, \omega^{\kappa}(-D)) = M_I \otimes^{\mathbf{L}}_{\Lambda_I, \kappa} \whalingship_p$$
  \end{thm}
  \begin{proof}
    For all classical algebraic weights, let us denote by $d_i(\kappa)$
    the dimension of ${\HH}^i(M_I \otimes^{\mathbf{L}}_{\Lambda_I, \kappa} \whalingship_p)$
    and by $d_i ^\dag (\kappa)$ the dimension of
    $e(U^I) {\HH}^i_{cusp} (G_1, \dag, \kappa)$. We deduce from Proposition~\ref{prop: I equal 1 Euler char} that
    $d_1(\kappa)-d_0(\kappa) = d_1^\dag(\kappa) - d_0^\dag (\kappa)$
    for all $\kappa$. Therefore
    $d_1(\kappa) - d_1^\dag (\kappa) = d_0(\kappa) - d_0^\dag(\kappa)$
    for all $\kappa$. By Proposition~ \ref{prop-ho-h1}, the first difference is non-positive and the second
    difference is non-negative,  so both
    are equal to zero. We  deduce in particular that the map
    $e(U^I) {\HH}^i _{cusp}(G_1, \dag, \kappa) \rightarrow H^i(M_I
    \otimes^{\mathbf{L}}_{\Lambda_I, \kappa} \whalingship_p)$ is an isomorphism.

    We now consider the composite
    $$e(U^I) {\HH}^i( \mathcal{X}^{G_1, \mult, \dag}_{K_p(I) K^p}, \omega^{\kappa}(-D))
    \to e(U^I){\HH}^i_{cusp}(G_1, \dag, \kappa)
    \rightarrow {\HH}^i( M_I \otimes^{\mathbf{L}}_{\Lambda_I, \kappa} \whalingship_p).$$  
By Proposition~\ref{prop-ho-h1} this composite map is an isomorphism
for~$i=0$, and is surjective for~$i=1$. On the other hand, the first
map is injective for~$i=1$ by Corollary~\ref{cor-boundslope}, and we
have just seen that the second map is an isomorphism. It follows that
the composite is injective for~$i=1$, and is thus an isomorphism, as required.
\end{proof}
Finally, we deduce the following classicity theorem.
  \begin{thm}\label{thm: ordinary classicity for I at most 1}Suppose that~$\#I\le 1$, and that~$\kappa$ is a
classical algebraic weight with~$l_v=2$ if~$v\in I$, and~$l_v\ge 4$
    if~$v\notin I$. Then the canonical
    map \[e(U^I)\mathrm{R}\Gamma(X^{G_1}_{K_p(I) K^p},\omega^\kappa(-D))[1/p]\to M_I \otimes^{\mathbf{L}}_{\Lambda_I, \kappa}\Qp \]is a quasi-isomorphism.    
  \end{thm}
  \begin{proof} This follows from
    Theorem~\ref{thm: ordinary is overconvergent I equal 1} and
    Theorem~\ref{thm-classicality}.  \end{proof}
\section{The Taylor--Wiles/Calegari--Geraghty method}\label{sec:CG}In this section, we implement the Taylor--Wiles patching method to
patch the complexes~$M_I$
constructed in~\S\ref{sec: Hida complexes}. More precisely, we carry out the
analogue of the patching argument using ``balanced modules'' which was
introduced in~\cite{CG}, and used there to study weight one modular
forms for~$\GL_2/\Q$. This argument works in situations where the
cohomology appears in at most two degrees, which for us means that~$\#
I\le 1$; we are restricted to working in this case due to the
limitations of our understanding of the cohomology of our complexes in
higher degree, as was the case in \S\ref{sec: Hida complexes}
and~\S\ref{sec: higher Coleman theory}. For our modularity result, it is
crucial to be able to work with~$I=S_p$; we will do this
in~\S\ref{sec:gluing} by considering the  spaces of modular forms coming from the various complexes
with~$\#I \le 1$. 

The papers~\cite{MR2234862,MR2920881,CGGSp4} apply the
Taylor--Wiles method to~$\GSp_4$ over~$\Q$, but a number of changes
are needed in order to apply it over general totally real fields. We do not attempt to
prove results in maximal generality, but instead develop the minimal
amount of material that we need. The reader familiar with the
literature on modularity lifting theorems will not find many
surprises, but we highlight a few things that may be less standard:
\begin{itemize}
\item In~\S\ref{subsec:Ordinary
  lifts}, we study the ordinary deformation rings at places
dividing~$p$. We show that their generic fibres are irreducible under
a rather mild $p$-distinguishedness assumption; in particular, this
assumption is not sufficient to guarantee that the deformation rings
are formally smooth, and it takes us some effort to prove the
irreducibility. Working in this generality is important for our
applications to modularity of abelian surfaces in~\S\ref{sec:modularityapplications}.
 For the potential modularity
results of~\S\ref{section:
  potential modularity of abelian surfaces}, however, it would be enough to
work with a stronger $p$-distinguishedness assumption which would
guarantee the formal smoothness.
\item In~\S\ref{subsec: l not
  p deformations}, we prove the statements about local deformation
rings needed for Taylor's ``Ihara avoidance argument''; the proofs are similar to those for~$\GL_n$, although there are some
complications which arise because the
relationship between conjugacy classes and characteristic polynomials
is more complicated. We also need to do some additional work to handle
the case $p=3$; again, this is crucial for~\S\ref{sec:modularityapplications}, although it is not needed for~\S\ref{section:
  potential modularity of abelian surfaces}.
\item In~\S\ref{subsec: big image conditions}, we study the
  ``big image'' conditions needed in the Taylor--Wiles method. Here
  our approach is slightly different from that of~\cite{cht} and the
  papers that followed it; again, this is with the applications of
 ~\S\ref{sec:modularityapplications} in mind, where it is
  important to be able to consider representations with
  image~$\GSp_4(\F_3)$. For the same reason, when we impose a
  condition at an auxiliary prime which will allow us to assume that
  our Shimura varieties are at neat level, we make the weakest
  hypothesis that we can, at the expense of slightly complicating the
  corresponding local representation theory.
\item We make repeated use of the doubling results of
 ~\S\ref{sec:doubling}; they are needed in order to prove
  local-global compatibility for the Galois representations we
  consider, and also to compare the spaces of $p$-adic modular forms
  for different~$I$. 
\item Our implementation of the ``Ihara avoidance'' argument
  of~\cite{tay} uses the framework of~\cite{emertongeerefinedBM,shottonGLn}, and compares the underlying
  cycles of various patched modules. In particular, we use the patched
  modules with~$I=~\emptyset$ to prove a local result, which we then
  apply to the patched modules with $\# I=~1$.  In order to apply Ihara avoidance,  we repeatedly use the fact that the Galois
representation associated to our abelian surface is pure, to deduce
that the corresponding points on the generic fibres of the local deformation rings are
smooth; we use this smoothness to be able to compute the dimensions of
various spaces of $p$-adic automorphic forms, using a
characteristic~$0$ version of the freeness arguments of Diamond and Fujiwara~\cite{MR1440309}. While we do not use the full strength of purity, since we make
arguments with base change we would otherwise need to impose a
hypothesis of being ``stably generic'' on our local Galois
representations, and we do not know of any natural examples where this
condition is known, but purity is not.
\end{itemize}
Having carried out the patching argument, we know from the results of~\S\ref{sec: higher Coleman theory} that for each~$I$
with~$\#I\le 1$ there  is a nonzero space of ordinary $p$-adic
modular forms corresponding to our given Galois representation, which
are ``overconvergent in the direction of~$I$''. We will combine these
spaces in~\S\ref{sec:gluing}, using as an input that by a version
of Diamond's multiplicity one argument~\cite{MR1440309}, we know the
dimensions of these spaces when~$\#I \le 1$ (they are given by the expected product of
local terms). (Here we are again using our assumption that the local Galois
representations are pure, in order to know that the corresponding
points of the generic fibres of the local deformation rings are
smooth. A similar characteristic zero version of Diamond's argument
first appeared in~\cite{MR3546966}.) 

\subsection{Galois deformation rings} \label{subsec: Galois
  deformation rings}
We let $E$ be a finite extension of $\Qp$ with ring of integers $\cO$,
uniformizer~$\lambda$ and residue field $k$. We will always assume
that $E$ is chosen to be large enough such that all irreducible components of all
deformation rings that we consider, and all irreducible components of
their special fibres, are geometrically irreducible. (We are always
free to enlarge~$E$ in all of the arguments that we make, so this is
not a serious assumption.) Given a complete
Noetherian local $\cO$-algebra $\Lambda$ with residue field $k$, we
let $\CNL_{\Lambda}$ denote the category of complete Noetherian local
$\Lambda$-algebras with residue field $k$.  We refer to an object in
$\CNL_{\Lambda}$ as a $\CNL_{\Lambda}$-algebra.

We fix a totally real field $F$, and let $S_p$ be the set of places of
$F$ above $p$.  We assume that $E$ contains all embeddings of $F$ into
an algebraic closure of $E$.  We also fix a continuous absolutely
irreducible homomorphism $\rhobar: G_F \to \GSp_4(k)$.  We assume
throughout that $p>2$.

Let $S$ be a finite set of finite places of $F$ containing $S_p$ and
all places at which $\rhobar$ is ramified. We write $F_S$ for the
maximal subextension of $\overline{F} / F$ which is unramified outside
$S$, and write~$G_{F,S}$ for~$\Gal(F_S/F)$.  For each $v\in S$, we fix
$\Lambda_v \in \CNL_{\cO}$, and set
$\Lambda = \widehat{\otimes}_{v \in S} \Lambda_v$, where the completed
tensor product is taken over $\cO$.  Then $\CNL_\Lambda$ is a
subcategory of $\CNL_{\Lambda_v}$ for each $v\in S$, via the canonical
map $\Lambda_v \to \Lambda$.

\begin{rem}
  \label{rem: choice of Lambda}In our applications, we will
  take~$\Lambda_v=\cO$ if~$v\nmid p$. If~$v|p$, then we will take ~$\Lambda_v$ to be an Iwasawa algebra.\end{rem}
Fix a character~$\psi:G_{F,S}\to\Lambda^\times$ lifting~$\nu\circ\rhobar$.
\begin{defn}\label{defn: local lift}
  A \emph{lift}, also called a \emph{lifting}, of $\rhobar|_{G_{F_v}}$
  is a continuous homomorphism
  $\rho: G_{F_v} \to \GSp_4(A)$ to a
  $\CNL_{\Lambda_v}$-algebra $A$, such that
  $\rho \bmod \frakm_A = \rhobar|_{G_{F_v}}$ and $\nu\circ\rho=\psi|_{G_{F_v}}$.
\end{defn}

We let $\cD_v^\square$ denote the set-valued functor on
$\CNL_{\Lambda_v}$ that sends $A$ to the set of lifts of
$\rhobar|_{G_{F_v}}$ to $A$.  This functor is representable (see for
example~\cite[Thm.\ 1.2.2]{MR3152673}), and we
denote the representing object by $R_v^\square$.

Let $x\in\Spec R_v^\square[1/p]$ be a closed point. By~\cite[Lem.\ 1.6]{tay} the residue field of~$x$ is a finite
extension $E'/E$. Let~$\rho_x:G_{F_v}\to\GL_n(E')$ be the
corresponding specialization of the universal lifting. By an argument of Kisin,  $(R_v^\square[1/p])^{\wedge}_x$ is the universal lifting ring
for $\rho_x$, i.e.\ if $A$ is an Artinian local $E'$-algebra with
residue field $E'$ and if $\rho:\Gamma \to \GSp_4(A)$ is a
continuous representation lifting $\rho_x$, then there is a unique
continuous map of $E'$-algebras
$(R_v^\square[1/p])^{\wedge}_x\to A$ so that the universal lift
pushes forward to $\rho$. (See~\cite[Thm.\ 1.2.1]{MR3546966}
for the analogous result for~$\GL_n$; the result for~$\GSp_4$ 
can be proved by an identical argument.) We say that~$x$ is \emph{smooth} if
$(R_v^\square[1/p])^{\wedge}_x$ is regular. Let~$\ad^0\rho_x$ denote
the Lie algebra~$\mathfrak{g}^0(E')$ with the adjoint action of~$G_F$
via~$\rho_x$; then we have the following convenient criterion for~$x$
to be smooth.

\begin{lem}\label{lem: smooth points of Ihara1}Suppose that~$v\nmid p$. Then the point $x$ is smooth if and only if
  $(\ad^0\rho_x)(1)^{G_{F_v}}=0$. In particular, if~$\rho_x$
  is pure, then~$x$ is smooth.
\end{lem}\begin{proof}The first claim is a special case of~\cite[Cor.\ 3.3.4,
  Rem.\ 3.3.6]{Bellovin}. If~$\rho_x$ is pure, then~$\Hom_{E'[G_{F_v}]}(\rho_x,\rho_x(1))=0$
  (because the definition of purity is easily seen to preclude the
  existence of a morphism between the corresponding Weil--Deligne
  representations), as required.
\end{proof}

\begin{defn}\label{def:locdefprob}
 A \emph{local deformation problem} for $\rhobar|_{G_{F_v}}$ is a subfunctor $\cD_v$ of $\cD_v^\square$ satisfying the following:
  \begin{itemize}
   \item $\cD_v$ is represented by a quotient $R_v$ of $R_v^\square$.
   \item For all $A \in \CNL_{\Lambda_v}$, $\rho \in \cD_v(A)$, and $a \in \ker(\GSp_4(A) \to \GSp_4(k))$, 
   we have $a\rho a^{-1} \in \cD_v(A)$. 
  \end{itemize}
\end{defn}

\begin{defn}\label{def:globdefprob}
 A \emph{global deformation problem} is a tuple
  \[
   \cS = (\rhobar, S, \{\Lambda_v\}_{v\in S},\psi, \{\cD_v\}_{v\in S}),
  \]
 where:
 \begin{itemize}
  \item $\rhobar$, $S$, $\{\Lambda_v\}_{v\in S}$ and~$\psi$ are as above.
  \item For each $v\in S$, $\cD_v$ is a local deformation problem for $\rhobar|_{G_{F_v}}$.
 \end{itemize}
\end{defn}

As in the local case, a \emph{lift} (or \emph{lifting}) of $\rhobar$
is a continuous homomorphism $\rho: G_{F,S} \to \GSp_4(A)$
to a $\CNL_\Lambda$-algebra $A$, such that
$\rho \bmod \frakm_A = \rhobar$ and $\rho\circ\nu=\psi$.  We say that
two lifts $\rho_1,\rho_2: G_{F,S} \to \GSp_4(A)$ are
\emph{strictly equivalent} if there is an
$a\in \ker(\GSp_4(A) \to \GSp_4(k))$ such that
$\rho_2 = a\rho_1 a^{-1}$.  A \emph{deformation of $\rhobar$} is a
strict equivalence class of lifts of $\rhobar$.

For a global deformation problem 
  \[
   \cS = (\rhobar, S, \{\Lambda_v\}_{v\in S},\psi, \{\cD_v\}_{v\in S}),
  \]
  we say that a lift $\rho: G_F \to \GSp_4(A)$ is of \emph{type
    $\cS$} if $\rho|_{G_{F_v}} \in \cD_v(A)$ for each $v\in S$.  Note
  that if $\rho_1$ and $\rho_2$ are strictly equivalent lifts of
  $\rhobar$, and $\rho_1$ is of type $\cS$, then so is $\rho_2$.  A
  \emph{deformation of type $\cS$} is a strict equivalence class of
  lifts of type $\cS$, and we denote by $\cD_{\cS}$ the set-valued
  functor that takes a $\CNL_\Lambda$-algebra $A$ to the set of lifts
  $\rho: G_F \to \GSp_4(A)$ of type $\cS$.

  Given a subset $T\subseteq S$, a \emph{$T$-framed lift of type
    $\cS$} is a tuple $(\rho,\{\gamma_v\}_{v\in T})$, where $\rho$ is
  a lift of type $\cS$, and
  $\gamma_v \in \ker(\GSp_4(A) \to \GSp_4(k))$ for each $v\in T$.  We
  say that two $T$-framed lifts $(\rho_1,\{\gamma_v\}_{v\in T})$ and
  $(\rho_2,\{\gamma'_v\}_{v\in T})$ to a $\CNL_\Lambda$-algebra $A$ are
  strictly equivalent if there is an $a\in \ker(\GSp_4(A) \to \GSp_4(k))$
  such that $\rho_2 = a \rho_1 a^{-1}$, and $\gamma'_v = a\gamma_v$ for
  each $v\in T$.  A strict equivalence class of $T$-framed lifts of
  type $\cS$ is called a \emph{$T$-framed deformation of type $\cS$}.
  We denote by $\cD_{\cS}^T$ the set valued functor that sends a
  $\CNL_\Lambda$-algebra $A$ to the set of $T$-framed deformations to
  $A$ of type $\cS$.

  The functors $\cD_{\cS}$ and $\cD_{\cS}^T$ are representable (as we
  are assuming that~$\rhobar$ is absolutely irreducible), and we
  denote their representing objects by $R_{\cS}$ and $R_{\cS}^T$
  respectively. If~$T$ is empty, then $R_{\cS}=R_{\cS}^T$, and
  otherwise the natural map~$R_{\cS}\to R_{\cS}^T$ is formally smooth
  of relative dimension~$11\#T-1$.  Indeed $\cD_{\cS}^T \rightarrow \cD_{\cS}$ is a torsor under $(\prod_{v \in T}  \widehat{\GSp}_4)/\widehat{\mathbb{G}_m}$. Define $\cT $ to be the coordinate ring  of $(\prod_{v \in T}  \widehat{\GSp}_4)/\widehat{\mathbb{G}_m}$ over $\Lambda$. This is a power series algebra over $\Lambda$ in $11\#T-1$ variables.

  \begin{lem}\label{lem: framed deformation ring over deformation ring}
    The choice of a representative
    $\rho_{\cS} \colon G_F \rightarrow \GSp_4(R_{\cS})$ for the
    universal type $\cS$ deformation determines  a splitting of the torsor $\cD_{\cS}^T \rightarrow \cD_{\cS}$  and   a canonical
    isomorphism $R_{\cS}^T \cong R_{\cS} \widehat{\otimes}_\Lambda \cT$.
  \end{lem}
  \begin{proof} This is obvious. 
  \end{proof}

\subsection{Galois cohomology and presentations}\label{subsec: Galois
  cohomology and presentations}
Fix a global deformation problem
  \[
   \cS = (\rhobar, S, \{\Lambda_v\}_{v\in S},\psi, \{\cD_v\}_{v\in S}),
  \]
  and for each $v\in S$, let $R_v$ denote the object representing
  $\cD_v$.  Let $T$ be a subset of~ $S$ containing~$S_p$, with the property that
  $\Lambda_v = \cO$ and~$\cD_v=\cD_v^\square$ for all $v\in S
  \smallsetminus T$. Define
  $R_{\cS}^{T,\loc} = \widehat{\otimes}_{v\in T} R_v$, with the
  completed tensor product being taken over $\cO$.  It is canonically
  a $\Lambda$-algebra, via the canonical isomorphism
  $\widehat{\otimes}_{v\in T} \Lambda_v \cong \widehat{\otimes}_{v\in
    S} \Lambda_v$.  For each $v\in T$, the morphism
  $\cD_{\cS}^T \to \cD_v$ given by
  $(\rho,\{\gamma_v\}_{v\in T}) \mapsto \gamma_v^{-1}\rho|_{G_{F_v}}
  \gamma_v$ induces a local $\Lambda_v$-algebra morphism
  $R_v \to R_{\cS}^T$.  We thus have a local $\Lambda$-algebra
  morphism $R_{\cS}^{T,\loc} \to R_{\cS}^T$.  

  The relative tangent space of this map is computed by a standard
  calculation in Galois cohomology, which we now recall. We let $\ad\rhobar$ (resp.\
  $\ad^0\rhobar$) denote~$\fg(k)$ (resp.\ $\fg^0(k)$), with the
  adjoint $G_F$-action via $\rhobar$.

The trace pairing $(X,Y) \mapsto \tr(XY)$ on $\ad^0\rhobar$ is perfect and $G_F$-equivariant, so $\ad^0\rhobar(1)$ is isomorphic to the 
Tate dual of $\ad^0\rhobar$. We define
\[
   H_{\cS,T}^1(\ad^0\rhobar):= \ker \left( H^1(F_S/F,\ad^0\rhobar) \rightarrow 
   \prod_{v\in T} H^1(F_v,\ad^0\rhobar)\right),
  \] \[
   H_{\cS^\perp,T}^1(\ad^0\rhobar(1)):= \ker \left( H^1(F_S/F,\ad^0\rhobar(1)) \rightarrow 
   \prod_{v\in S \smallsetminus T} H^1(F_v,\ad^0\rhobar(1))\right).
  \]

\begin{prop}\label{prop: tangent space of local over global}
  Continue to assume that $T$ contains~$S_p$, and that for all
  $v\in S \smallsetminus T$ we have
  $\Lambda_v = \cO$ and~$\cD_v=\cD_v^\square$ .  Then there is a local $\Lambda$-algebra
  surjection $R_{\cS}^{T,\loc}\llbracket X_1,\ldots,X_g \rrbracket \rightarrow R_{\cS}^T$, with
    \begin{multline*}
    g  = h_{\cS^\perp,T}^1(\ad^0\rhobar(1)) - h^0(F_S/F,\ad^0\rhobar(1))  - \sum_{v| \infty} h^0(F_v,\ad^0\rhobar) \\
     + \sum_{v \in S \setminus T} h^0(F_v,\ad^0\rhobar(1))+\#T-1.
    \end{multline*}
\end{prop}
\begin{proof}
  We follow~\cite[\S 3.2]{kis04}. 
   By~\cite[Lem.\ 3.2.2]{kis04} (or
  rather the same statement for~$\GSp_4$, which has an identical
  proof), the claim of the proposition holds with \[g=
    h_{\cS,T}^1(\ad^0\rhobar)-h^0(F_S/F,\ad\rhobar)+\sum_{v\in
      T}h^0(F_v,\ad\rhobar).\]
  By~\cite[Thm.\ 2.19]{MR1605752} (and the assumption that~$\rhobar$
  is absolutely irreducible, which implies that
  $h^0(F_S/F,\ad\rhobar^0)=0$), we have 
  \begin{multline*}
    h_{\cS,T}^1(\ad^0\rhobar)=h_{\cS^\perp,T}^1(\ad^0\rhobar(1))-
    h^0(F_S/F,\ad^0\rhobar(1))- \sum_{v| \infty}
    h^0(F_v,\ad^0\rhobar) \\+\sum_{v\in S\setminus
      T}h^1(F_v,\ad^0\rhobar)-\sum_{v\in S}h^0(F_v,\ad^0\rhobar).
  \end{multline*}
The result follows from the local Euler characteristic formula and
Tate local duality.
\end{proof}

\subsection{Local deformation problems, \texorpdfstring{$l=p$}{l=p}}\label{subsec:Ordinary
  lifts}

Assume from now on that~$p$ splits completely in~$F$. Let~$v$
be a place of~$F$ lying over~$p$. If~$x\in k^\times$, then we
write~$\lambda_x:G_{F_v}\to k^\times$ for the unramified character
with~$\lambda_x(\Frob_v)=x$.
\begin{defn}\label{defn: generic flat ordinary}
  We say that~$\rhobar|_{G_{F_v}}$ is \emph{$p$-distinguished weight~$2$
    ordinary} if it is conjugate to a representation of
  the form \[ \begin{pmatrix}
      \lambda_{\alphabar_v}&0&*&*\\
      0 &\lambda_{\betabar_v}&*&*\\
      0&0&\varepsilonbar^{-1}\lambda_{\betabar_v}^{-1}&0\\
      0&0&0&\varepsilonbar^{-1}\lambda_{\alphabar_v}^{-1}
    \end{pmatrix},\]where $\alphabar_v\ne \betabar_v$.

If~$\rhobar|_{G_{F_v}}$ is $p$-distinguished weight~$2$ ordinary, then we say that a
lift~$\rho:G_{F_v}\to\GSp_4(\cO)$ of~$\rhobar|_{G_{F_v}}$ is  $p$-distinguished weight~$2$ ordinary if
 $\rho$ itself is conjugate to a representation of the form  \[ \begin{pmatrix}
      \lambda_{\alpha_v}&0&*&*\\
      0 &\lambda_{\beta_v}&*&*\\
      0&0&\varepsilon^{-1}\lambda_{\beta_v}^{-1}&0\\
      0&0&0&\varepsilon^{-1}\lambda_{\alpha_v}^{-1}
    \end{pmatrix}\]where $\alpha_v$, $\beta_v$ lift $\alphabar_v$,
  $\betabar_v$ respectively. Note that~$\rho$ is then automatically
  semistable, although not necessarily crystalline.
\end{defn}
\begin{rem}\label{rem: apologize for name of weight 2 ordinary}The terminology
  ``weight $2$ ordinary'' is not ideal, but we were unable to find a
  better alternative. Possibilities include ``$P$-ordinary''
  (referring to the Siegel parabolic subgroup), which clashes with
  ``$p$-distinguished'', or ``semistable ordinary''. We could of
  course restrict to the crystalline case and use ``flat ordinary''
  representations, but as it costs us little to allow semistable
  representations, and it may prove to be useful in future applications, we
  have not done this.
\end{rem}
\begin{rem}
  \label{rem: not very p distinguished}For the purposes of proving the
  potential modularity of abelian surfaces, it would suffice to work
  with a stronger $p$-distinguishedness hypothesis, as
  in~\cite{CGGSp4}. In particular, by assuming that none of
  $\alphabar_v^2,\betabar_v^2,\alphabar_v\betabar_v$ are equal to~$1$, we could arrange that the various
  deformation rings considered in this section are formally
  smooth. However, such a hypothesis is very restrictive in the
  case~$p=3$, and in particular would seriously restrict the
  applicability of our modularity lifting theorems to proving the
  modularity (as opposed to potential modularity) of particular
  abelian surfaces.
\end{rem}

We assume from now on that~$\rhobar|_{G_{F_v}}$ is $p$-distinguished weight~$2$
ordinary for all~$v|p$; the roles of~$\alphabar_v,\betabar_v$ in the
definition of $p$-distinguished weight~$2$ ordinary are symmetric, and we fix a
labelling of~$\alphabar_v,\betabar_v$ for each~$v|p$. 

Set $\Lambda_{v,1}=\cO\llbracket\cO_{F_v}^\times(p)\rrbracket$,
~$\Lambda_{v,2}=\cO\llbracket(\cO_{F_v}^\times(p))^2\rrbracket$, where
$\cO_{F_v}^\times(p) = 1+ p \ocal_{F_v}$ denotes the
pro-$p$ completion of $\cO_{F_v}^\times$. Both $\Lambda_{v,1}$ and $\Lambda_{v,2}$ are formally smooth
over~$\cO$ (because we are assuming that~$F_v=\Qp$). Let~$\Lambda_v$
be either~$\Lambda_{v,1}$ or~$\Lambda_{v,2}$. There is a canonical character
$I_{F_v}\to\cO_{F_v}^\times(p)$ given by~$\Art_{F_v}^{-1}$, and we
define a pair of characters~$\theta_{v,i}:I_{F_v}\to
\Lambda_{v}$, $i=1,2$ as follows: if~$\Lambda_v=\Lambda_{v,1}$, then we let
$\theta_{v,1}=\theta_{v,2}=\theta_v$ be the natural character 
and
if~$\Lambda_v=\Lambda_{v,2}$, then we let~$\theta_{v,i}$
correspond to the embedding  ~$\cO_{F_v}^\times(p)$ to
$(\cO_{F_v}^\times(p))^2$ given by the $i$th copy.

Let~$\alphabetabar_v$ denote a choice of
either~$\alphabar_v$ or~$\betabar_v$, and write $\alphabetabar_v':=\alphabar_v\betabar_v/\alphabetabar_v$
for the other choice.
Recall that we have the Borel subgroup~$B$ of~$\GSp_4$ consisting of
matrices of the form \[ \begin{pmatrix}
     *&*&*&*\\
      0 &*&*&*\\
      0&0&*&*\\
      0&0&0&*
    \end{pmatrix}.\]We let~$P$ denote the subgroup of~$B$ consisting
of matrices of the form \[ \begin{pmatrix}
     *&0&*&*\\
      0 &*&*&*\\
      0&0&*&0\\
      0&0&0&*
    \end{pmatrix}.\]

If~$A\in\CNL_{\Lambda_v}$, then we say that a lift
$\rho_A:G_{F_v}\to\GSp_4(A)$ of~$\rhobar|_{G_{F_v}}$ is $(B,\alphabetabar_v)$-ordinary
if there is an increasing filtration of free~$A$-submodules  \[
     0 = \Fil^0 \subset \Fil^1 \subset \cdots \subset \Fil^4 = A^4
    \] of~$A^4$ by
$A[G_{F_v}]$-submodules such that the action of~$G_{F_v}$ on $\Fil^i/\Fil^{i-1}$ is via a character~$\chi_i$
with $\chibar_1=\lambda_{\alphabetabar_v}$,
$\chibar_2=\lambda_{\alphabeta'_v}$, $\chi_1|_{I_{F_v}}=\theta_{v,1}$,
$\chi_2|_{I_{F_v}}=\theta_{v,2}$, and $\chi_3=\varepsilon^{-1}\chi_2^{-1}$, $\chi_4=\varepsilon^{-1}\chi_1^{-1}$. 

By~\cite[Lem.\ 2.4.6]{cht} such a filtration is unique;
since~$\{(\Fil^{4-i})^{\perp}\}$ gives another filtration satisfying
the same conditions, we see that~$\rho_A$ is
$\ker(\GSp_4(A)\to\GSp(k))$-conjugate to a representation of the form \[ \begin{pmatrix}
    \chi_1&*&*&*\\
0 &\chi_2&*&*\\
0&0&\varepsilon^{-1}\chi_2^{-1}&*\\
0&0&0&\varepsilon^{-1}\chi_1^{-1}
  \end{pmatrix}\]where $\chi_1$, $\chi_2$ are as above. 

If $\Lambda_v=\Lambda_{v,1}$ (so that $\theta_{v,1}=\theta_{2,v}$), then we say that~$\rho_A$  is
$P$-ordinary if it is both $(B,\alphabar_v)$-ordinary
and $(B,\betabar_v)$-ordinary; equivalently, if~$\rho_A$ is $\ker(\GSp_4(A)\to\GSp(k))$-conjugate to a representation of the form \[ \begin{pmatrix}
    \chi_1&0&*&*\\
0 &\chi_2&*&*\\
0&0&\varepsilon^{-1}\chi_2^{-1}&0\\
0&0&0&\varepsilon^{-1}\chi_1^{-1}
  \end{pmatrix}.\]
If~$\Lambda_v=\Lambda_{v,2}$ (resp.\ $\Lambda_v=\Lambda_{v,1}$) then we let~$\cD_v^{B,\alphabetabar_v}$
(resp.\ $\cD_v^{P}$) be the subfunctor of
$(B,\alphabetabar_v)$-ordinary lifts (resp.\ of
$P$-ordinary lifts). By~\cite[Lem.\ 2.4.6]{cht}, we see
that $\cD_v^{B,\alphabetabar_v}$ and $\cD_v^{P}$ are local
deformation problems in the sense of Definition~\ref{def:locdefprob},
so they are represented by $\CNL_{\Lambda_v}$-algebras
$R_v^{B,\alphabetabar_v}$, $R_v^{P}$ respectively.

Most of the rest of this section is devoted to the proof of the following
result. 
\begin{prop}\label{prop: local ordinary deformation ring
    dimension}The generic fibres $R_v^{B,\alphabetabar_v}[1/p]$,
  $R_v^{P}[1/p]$ are irreducible,   and are of relative dimensions~$16$ and~$14$
  respectively over~$\Qp$.
  \end{prop}

Our arguments are rather ad hoc, and will require a
number of preliminary lemmas. 

\subsubsection{Ordinary deformation rings for~$\GL_2$}
We begin by studying some ordinary deformation
rings for~$\GL_2$. As well as being a warmup for our main arguments, we
will often be able to show that our deformation rings for~$\GSp_4$ are formally
smooth over a completed tensor product of deformation rings
for~$\GL_2$, thus reducing to this case. 

Let~$\rbar:G_{\Qp}\to\GL_2(k)$ be of the form \[\begin{pmatrix}
      \lambda_{\alphabar}&*\\
      0
      &\varepsilonbar^{-1}\lambda_{\alphabar}^{-1}    \end{pmatrix}.\]Set~$\Lambda=\cO[[1+p\Zp]]$,
  and write $\theta:I_{\Qp}\to\Lambda$ for the canonical character
  defined above. If~$A\in\CNL_\Lambda$, then we say that a lift
  of~$\rbar$ to $r:G_{\Qp}\to\GL_2(A)$ is ordinary if it
  is $\ker(\GL_2(A)\to\GL_2(k))$-conjugate to a representation of the
  form
\[\begin{pmatrix}
      \chi&*\\
      0
      &\varepsilon^{-1}\chi^{-1}    \end{pmatrix}\]where~$\chibar=\lambda_{\alphabar}$
  and~$\chi|_{I_{\Qp}}=\theta$. As above, this is a local deformation
  problem, and is represented by a $\CNL_\Lambda$-algebra
  $R^{B_2,\square}$, where~$B_2$ denotes
  the Borel subgroup of~$\GL_2$ of upper triangular matrices.
  
We write  \[\rbar=\begin{pmatrix}
      \lambda_{\alphabar}&\varepsilonbar^{-1}\lambda_{\alphabar}^{-1}\eta_{\alpha^2}\\
      0
      &\varepsilonbar^{-1}\lambda_{\alphabar}^{-1}    \end{pmatrix},\]where~$\eta_{\alpha^2}$
  is a cocycle in $Z^1(\Q_p,\varepsilonbar
  \lambda_{\overline{\alpha^2}})$. Rescaling our basis vectors has the
  effect of changing~$\eta_{\alpha^2}$ by a coboundary, so we can and
  do think of~$\eta_{\alpha^2}$ as a class in~$H^1(\Q_p,\varepsilonbar
  \lambda_{\overline{\alpha^2}})$. 

Let~$\mathfrak{b}_2$ be the Lie algebra of~$B_2$, given by the matrices $$ \mathfrak{b}_2 = \begin{pmatrix} \nu + x_{\alpha} &  x_{\alpha^2}   \\   0& \nu-x_{\alpha} \end{pmatrix},$$
where~$\ad^0_{B_2}$ corresponds to~$\nu = 0$.
With respect to the basis given by the matrices corresponding to the
variables~$\{x_{\alpha^2},x_{\alpha}\}$ --- that is, the basis $\left\{
\begin{pmatrix}
  0&1\\0&0
\end{pmatrix},
\begin{pmatrix}
  1&0\\0&-1
\end{pmatrix}\right\}$ --- the Galois
representation~$\ad^0_{B_2}\rbar$ is given explicitly as follows:
$$\left(\begin{matrix}
\varepsilonbar  \lambda_{\alphabar^2} & -2  \eta_{\alpha^2}	\\0
&1							\end{matrix}
\right).$$ Note that if~$M$ is annihilated by~$p$, then~$H^2(\Q_p,M)$
is given by~$H^0(\Q_p,M^*)^{\vee} \simeq \Hom_{G_{\Q_p}}(M,\varepsilonbar)^\vee$. It follows
that~$h^2(\Qp,\ad_{B_2}^0\rbar)=0$ unless~$\alphabar^2=1$
and~$\eta_{\alpha^2}=0$, in which case we have~$h^2(\Qp,\ad_{B_2}^0\rbar)=1$.
We  write~$R^{B_2}$ for~$R^{B_2,\square}$ unless we particularly
want to emphasize the~$\GL_2$-framing variables.
\begin{lem}
  \label{lem: irreducibility and dimension of GL2 ordinary deformation
  ring}The generic fibre~$R^{B_2}[1/p]$ is irreducible, and
has relative dimension~$5$ over~$\Qp$.\end{lem}
\begin{proof}By a standard argument (see \cite[Prop.\ 2]{MR1012172}),
  $R^{B_2}$ has a presentation of the
  form~$\cO[[x_1,\dots,x_r]]/(y_1,\dots,y_s)$,
  where
 \[r=4-h^0(\Qp,\ad_{B_2}\rbar)+h^1(\Qp,\ad_{B_2}^0\rbar) =
  3 - h^0(\Qp,\ad^0_{B_2}\rbar)+h^1(\Qp,\ad_{B_2}^0\rbar),\]
  \[ s=h^2(\Qp,\ad_{B_2}^0\rbar).\]
    Note that, \emph{a priori}, even when~$s > 0$, some of the~$y_i$ may vanish, although one does not expect this to happen.
    By the local Euler characteristic
  formula, $r-s=3+\dim\ad_{B_2}^0\rbar=5$. In particular, if
  $H^2(\Qp,\ad_{B_2}^0\rbar)=0$, then~$R^{B_2}$ is formally smooth over~$\cO$ of relative
  dimension~$5$, and we are done.

If  ~$H^2(\Qp,\ad_{B_2}^0\rbar)\ne 0$, then the  above discussion
  shows that~$\alphabar=\pm 1$, ~$\rbar$ is split,
  and~$s = h^2(\Qp,\ad_{B_2}^0\rbar)=1$. Since any quotient of a formal power
  series ring by a single relation is a local complete intersection,
  it follows from the presentation of the previous paragraph
  that~$R^{B_2}$ is a local complete intersection, and in particular~$\Stwo$.
  Note that, at this point, we don't know if the relation~$y_1$ is non-zero or not, so we do not as yet know the dimension of~$R^{B_2}$.

  Twisting by a quadratic character, we can and do suppose
  that~$\alphabar=1$, so that~$\rbar=1\oplus\varepsilonbar^{-1}$. We
  begin by showing that~$\Spec R^{B_2}[1/p]$ is connected,
  following~\cite[Lem.\ 3.13]{ger}. 
  Note that the map $\Spec R^{B_2}[1/p]\to\Spec\Lambda[1/p]$ admits a section, because we can always find a lift of the
form~$\chi\oplus\chi^{-1}\varepsilon^{-1}$. Since~$\Spec\Lambda[1/p]$ is
connected, it therefore suffices to show that the fibres of this map
over closed points~$x$ of~$\Spec\Lambda[1/p]$ are connected.

By (for example) the proof of~ \cite[Lem.\ 1.2.2]{BLGGT} (see
also~\cite[Lem.\ 3.4.1]{Bellovin}), the irreducible
components of~$\Spec R^{B_2}[1/p]$ are fixed by conjugation by
elements of~$\GL_2(R^{B_2})$ whose image in~$\GL_2(k)$ is
diagonal. It is therefore obviously the case that all the closed points 
which are conjugate to representations
of the
form~$\chi\oplus\chi^{-1}\varepsilon^{-1}$ lie in the same connected
component of the fibre over~$x$, so it suffices to show that each
closed  point of
the form \[r=\begin{pmatrix}
      \chi&*\\
      0
      &\varepsilon^{-1}\chi^{-1}    \end{pmatrix}\]  lies in the same
  connected component as the corresponding point with~$*=0$. To this
  end, we consider the
  representation \[r_t:\diag(t,t^{-1})r\diag(t,t^{-1})^{-1}\to\GL_2(\cO\langle
    t\rangle).\]Note that the specializations of this representation
  at~$t=0$ and~$t=1$ correspond to the two closed points that we are
  considering. 

Letting~$A\subset\cO\langle t\rangle$ be the closed
  subalgebra generated by the matrix entries of the elements of the
  image~$r_t$, one checks exactly as in the proof of~\cite[Lem.\
  1.2.2]{BLGGT} that~$A$ is a complete local Noetherian
  $\cO$-algebra with residue field~$k$. Since~$\rbar$ is split, it
  follows that the representation~$r_t$ arises from a
  map~$R^{B_2}\to A$. Since~$A$ is a domain (being a subring
  of~$\cO\langle t\rangle$), we see that the points corresponding
  to~$t=0$ and~$t=1$ lie on the same irreducible component, as required.

To see that~$R^{B_2}[1/p]$ is moreover irreducible, it is enough to
check that it is normal, or equivalently that it is~$\Rone$ and~$\Stwo$.
We have
already seen that it is~$\Stwo$, and to show that it is~$\Rone$, it suffices to show that there is
an open regular subscheme~$U$ of~$R^{B_2}[1/p]$ whose complement has
codimension at least~$2$. We will in fact show that there is such a subscheme
with the property that the tangent space at any closed point~$u\in U$
has dimension~$5$, 
thus also proving the statement about the
dimension of~$R^{B_2}[1/p]$ (if the one relation in our presentation of~$R^{B_2}$ was trivial,
then~$R^{B_2}$ would be formally smooth of relative dimension~$6$, and there would be no such points).

Over~$R^{B_2}$, we have a universal
lifting~$r^\univ:G_{\Qp}\to\GL_2(R^{B_2})$, and we
let~$H^2_{\ord}:=H^2(G_{\Qp},\ad_{B_2}^0r^{\univ})$, a
finite~$R^{B_2}$-module. Since the cohomology of~$G_{\Qp}$
vanishes in degree greater than~$2$, the formation of~$H^2$ is
compatible with specialization, so that in particular if~$x$ is a
closed point of~$R^{B_2}[1/p]$ with corresponding
representation $r_x:G_{\Qp}\to\GL_2(E_x)$ (with~$E_x$ a finite
extension of~$\Qp$), then
$H^2(G_{\Qp},\ad_{B_2}^0r_x)=H^2_{\ord}\otimes_{R^{B_2}}E_x$. 

We let~$U$ be the complement of the support of~$H^2_{\ord}$ in~$\Spec
R^{B_2}[1/p]$.  (It is not obvious \emph{a priori} that~$U$ is not empty,
but we will prove this below.)
Then at any closed point~$x\in U$, we
have~$H^2(G_{\Qp},\ad_{B_2}^0r_x)=0$, so by a standard Galois cohomology
calculation (essentially identical to the one used in the first
paragraph of this proof), $U$ is
formally smooth over~$\Qp$ at~$x$, with relative tangent
space of dimension~$5$. It follows that~$U$ is regular.

The complement of~$U$ is the support of~$H^2$, so just as above,
its closed points are those~$x$ for which~$\rho_x$ is a direct sum of
two characters whose ratio is the cyclotomic character. But in any
Zariski open neighbourhood of such a point there are points of~$U$
(for example, points which are a direct sum of two characters whose
ratio is not the cyclotomic character, given by twisting the
characters occurring in~$\rho_x$ by unramified characters), so~$U$ is
dense in~$\Spec
R^{B_2}[1/p]$, and $\Spec
R^{B_2}[1/p]$ is equidimensional of relative dimension~$5$
over~$\Qp$.

It remains to show that the complement of~$U$ (that is, the support of~$H^2$) has codimension at
least~$2$, or equivalently that it has relative dimension at most~$3$
over~$\Qp$.  In fact, it has relative dimension at most~$2$
over~$\Qp$: the only
freedom we have is to make twists of the two characters
in~$\rho_x$ (and the determinant is fixed), so
the corresponding dimension is the dimension of~$\GL_2$ minus the dimension
of the centralizer of~$\displaystyle{ \left(\begin{matrix} 1 & 0 \\ 0 & \varepsilon^{-1} \end{matrix} \right)}$,
which equals~$4-2 = 2$,
so we are done.  \end{proof}

Let~$R^{B_2,\tri}$ denote the $B_2$-valued framed deformation ring of~$\rbar$
with fixed determinant. It follows
from the assumption that~$\rbar$ is~$p$-distinguished 
that~$R^{B_2,\square}$ is formally smooth over~$R^{B_2,\tri}$ of
relative dimension
$$\dim \ad^0_{\GL_2} - \dim \ad^0_{B_2} = 3 - 2 = 1$$(see
Lemma~\ref{lemma:getitdone} for the details of an analogous argument
in the symplectic case).

It will prove useful to give (somewhat) explicit descriptions of~$R^{B_2,\tri}$ 
(and thus~$R^{B_2}$) in a number of explicit cases.
Lemma~\ref{lem:explicit} below will also give another proof of Lemma~\ref{lem: irreducibility and dimension of GL2 ordinary deformation
  ring},
although not one we shall generalize to the symplectic context.

Let~$\gamma \in k$, and let
$$\rbar = \left( \begin{matrix} \lambda_{\gamma} & \varepsilonbar^{-1} \lambda^{-1}_{\gamma}
\eta \\ 0 & \varepsilonbar^{-1} \lambda^{-1}_{\gamma} \end{matrix} \right).$$
Via restriction to the character in the upper left hand
corner, the ring~$R^{B_2,\tri}$ is naturally
an algebra over the universal deformation ring~$R^{\GL_1}$ for~$\GL_1$. This gives a map
$$R^{\GL_1} \rightarrow R^{B_2,\tri}.$$
The ring~$R^{\GL_1}$ is formally smooth of relative dimension~$2$
over~$\OL$, and also formally smooth over the Iwasawa
algebra~$\Lambda$ corresponding to restricting the character to
inertia.  Let us choose
isomorphisms~$\Lambda = \OL \llbracket y_2 \rrbracket$
and~$R^{\GL_1} = \OL \llbracket y_1, y_2 \rrbracket$.  In the lemma
below, we shall use~$y_i$ for the variables of~$R^{B_2,\tri}$
corresponding to the algebra structure over~$R^{\GL_1}$, we use~$x_i$
for framing variables, and~$z_i$ for variables related to extensions
(informally corresponding to the upper right corner). More precisely,
by ``framing variables'' we mean the following: the map $c\mapsto (1+\varepsilon c)\rbar$ gives an
isomorphism
\[Z^1(\Q_p,\ad^0_{B_2})\isoto\Hom_{\cO}(R^{B_2,\tri},k[\varepsilon]/(\varepsilon^2))=\Hom_k(\m/\m^2,k),\]where~$\m$
is the maximal ideal of~$R^{B_2,\tri}$ and on the level of reduced
tangent spaces, the framing variables are by definition the
coboundaries $B^1(\Q_p,\ad^0_{B_2})\subset Z^1(\Q_p,\ad^0_{B_2})$.

\begin{lem}  \label{lem:explicit}
Let~$\m$ denote the maximal ideal of~$R^{B_2,\tri}$.
The ring~$R^{B_2,\tri}$  is a complete intersection, is  flat over~$\Lambda$, 
and is irreducible of relative dimension~$4$ over~$\OL$.
The rings~$R^{B_2,\tri}$ as~$R^{\GL_1}$-algebras have the following explicit
presentations. 
\begin{enumerate}
\item  \label{smoothlistening} If~$\gamma^2 \ne 1$ and~$\eta \ne 0$, then~$R^{B_2,\tri} \simeq \OL \llbracket x_{1}, x_{2},  y_{1}, y_{2} \rrbracket$.
\item  \label{jazz} If~$\gamma^2 \ne 1$ and~$\eta = 0$, then~$R^{B_2,\tri} \simeq \OL \llbracket x_{1}, z_1,  y_{1}, y_{2} \rrbracket$.
\item \label{smoothandnot} If~$\gamma^2 = 1$ but~$\eta \ne 0$, then:
\begin{enumerate}
\item  $R^{B_2,\tri}$ is formally smooth over~$\OL$, 
\item  \label{modpcrystalline} $R^{B_2,\tri}$ is formally smooth over~$\Lambda$ unless~$\eta$ is 
peu ramifi\'{e}e,
\item
$R^{B_2,\tri} \simeq \OL\llbracket x_{1}, x_{2}, z_{1}, y_{1}, y_{2} \rrbracket/g_{\eta}$,
where~$g_{\eta} = c_{\eta} y_{1} + d_{\eta} y_{2}  \mod (\lambda,\m^2)$
for~$[c_{\eta}:d_{\eta}] \in \mathbf{P}^1(k)$, and
where~$[c_{\eta}:d_{\eta}]$ depends only on~$\eta \in H^1(\varepsilonbar)$. 
\end{enumerate}
\item \label{nonsmooth} If~$\gamma^2 = 1$ and~$\eta  = 0$, then:
\begin{enumerate}
\item  $R^{B_2,\tri} \simeq \OL\llbracket x_{1}, z_{1}, z_{2}, y_{1}, y_{2} \rrbracket/g$,
where
$$g = z_{1} y_{1} + z_{2} y_{2}  \mod (\lambda,\m^3).$$
\item  \label{lemma:notformallysmooth} The special fibre~$R^{B_2,\tri}/\lambda$ is not formally smooth.
\end{enumerate}
\end{enumerate}
\end{lem}

\begin{remark} Explicit descriptions of ordinary deformation rings (even over general extensions~$K/\Q_p$) for~$\GL_2$ have
been given by B\"{o}ckle in~\cite[\S7]{Boeckle}. However, we require some precise information about these rings as algebras
over~$R^{\GL_1}$ and~$\Lambda$ which is not \emph{explicitly} given in the required form in~\cite{Boeckle}, and
thus we have found it easier to give the argument below. However, all the methods below already appear (in a more complicated setting)
in previous work of B\"{o}ckle and others. 
\end{remark}

\begin{proof} We first note that~$\ad^0_{B_2}$ is simply the~$2$-dimensional
representation given by
$$0 \rightarrow k(\lambda^{2}_{\gamma} \varepsilonbar ) \rightarrow \ad^0_{B_2}
\rightarrow k \rightarrow 0,$$
and where the extension class is given by~$\eta$ (so this is just a twist of~$\rbar$).
The framed tangent space has dimension
$$\dim Z^1(\Q_p,\ad^0_{B_2}) =
\dim H^1(\Q_p,\ad^0_{B_2}) + \dim B^1(\Q_p,\ad^0_{B_2}),$$
with precisely
$$\dim \ad^0_{B_2} - \dim H^0(\Q_p,\ad^{0}_{B_2})$$
framing variables.
The maps from~$R^{\GL_1}$ and from~$\Lambda$ correspond on
tangent spaces to  the maps $Z^1(\Q_p,\ad^0_{B_2}) \rightarrow
H^1(\Q_p,k) $ and $Z^1(\Q_p,\ad^0_{B_2}) \rightarrow H^1(I_{\Qp},k) $ given by the composites of the
maps
$$Z^1(\Q_p,\ad^0_{B_2}) \rightarrow
H^1(\Q_p,\ad^0_{B_2}) \rightarrow H^1(\Q_p,k) \rightarrow  H^1(I_{\Qp},k).$$
We now consider the  four possible cases in turn.

If~$\gamma^2 \ne 1$ and~$\eta \ne 0$, then~$H^0(\Q_p,\ad^0_{B_2})$ is trivial,
and there are two framing variables~$x_1$ and~$x_2$.
The map from~$H^1(\Q_p,\ad^0_{B_2})$
to~$H^1(\Q_p,k)$ is an isomorphism.
Note that~$H^2(\Q_p,\ad^0_{B_2}) = 0$, and so~$R^{B_2,\tri}$
is formally smooth over~$\cO$ and all statements are clear in this case.

If~$\gamma^2 \ne 1$ and~$\eta = 0$, then~$H^0(\Q_p,\ad^0_{B_2}) = k$ and
there is only one framing variable~$x_1$. However, the map
from~$H^1(\Q_p,\ad^0_{B_2})$
to~$H^1(\Q_p,k)$ is now surjective with kernel~$H^1(\Q_p,k(\varepsilonbar \lambda_{\gamma^2}))$,
which is of dimension one. Note that~$H^2(\Q_p,\ad^0_{B_2}) = 0$, and so~$R^{B_2,\tri}$
is formally smooth formally smooth over~$\cO$ and once again all statements are clear.

If~$\gamma^2 = 1$ but~$\eta \ne 0$, then~$H^2(\Q_p,\ad^0_{B_2}) = 0$
and the tangent space has dimension four, exactly two dimensions coming
from framing, one dimension from the image
of~$H^1(\Q_p,\varepsilonbar)$ in $H^1(\Q_p,\ad^0_{B_2})$,
and one dimension coming from the image of~$H^1(\Q_p,\ad^0_{B_2})$
in~$H^1(\Q_p,k)$. To compute the image, it suffices to consider the
(surjective) map
from~$H^1(\Q_p,k)$ to~$H^2(\Q_p,\varepsilonbar)$ and determine the kernel,
or, taking duals, considering the map~$H^0(\Q_p,k) \rightarrow H^1(\Q_p,\varepsilonbar)$
and taking the image. The image of the latter map is precisely given by~$\eta$.

The corresponding ring will fail to be flat over the space of weights~$\Lambda$ precisely when the image  of
$$H^1(\Q_p,\ad^0_{B_2}) \rightarrow H^1(\Q_p,k)$$
maps to zero in~$H^1(I_{\Qp},k)$, or equivalently when the image is unramified.
Under Tate local duality
for~$H^1(\Q_p,k) \times H^1(\Q_p,\varepsilonbar) \rightarrow k$, the unramified
classes are exactly annihilated by the peu ramifi\'{e}e classes. Hence
the failure of formal smoothness over~$\Lambda$ 
 occurs precisely
when~$\eta$ is peu ramifi\'{e}e. 
All the claims follow
except possibly the claim that~$R^{B_2,\tri}$ is flat over~$\Lambda$,
which is also transparent except in the peu ramif\'{e}e case, where~$R^{B_2,\tri}$
is formally smooth over~$\cO$ and is the quotient of a formally smooth $\Lambda$-algebra by the relation
~$g_\eta=y_2 \mod (\m^2,\lambda)$. 
This will be flat over~$\Lambda$ as long as~$y_2 \notin \lambda R^{B_2,\tri}$, which
can be easily ruled out by looking at points in characteristic
zero.
For example,  we see from~\cite[Theorem~2.1.8]{MR3628787} that the
fibre over every point in~$\Lambda[1/p]$ is non-trivial.
In particular, there are points
where the restriction to inertia of the character lifting~$\lambda_\gamma$ is finite of arbitrarily large order, so that~$v(y_2)$ becomes
arbitrarily close to~$0$, which would not be possible if~$y_2 \in \lambda R^{B_2,\tri}$.

It remains to consider the case when~$\rbar = 1 \oplus \varepsilon^{-1}$.
The representation underlying~$\rbar$ decomposes as a direct
sum which induces  corresponding decompositions of~$\ad^0_{B_2}$ and~$\ad^0_{\GL_2}$
respectively.
In particular, the adjoint~$\ad^0_{B_2} = k \oplus \varepsilonbar$ of~$\rbar$ thought of as inside the Borel is naturally a direct summand of~$\ad^0_{\GL_2} = k \oplus \varepsilonbar \oplus \varepsilonbar^{-1}$.
Let~$Z^1(\Q_p,\ad^0_{B_2})$ denote the~$1$-cocycles with values in~$\ad^0_{B_2}$.
There is a natural surjection
$$Z^1(\Q_p,\ad^0_{B_2}) \rightarrow H^1(\Q_p,\ad^0_{B_2}) = H^1(\Q_p,k) \oplus H^1(\Q_p,\varepsilonbar).$$
We now chose a basis for~$Z^1(\Q_p,\ad^0_{B_2})$ as follows:
\begin{enumerate}
\item $r_1$ generates the kernel~$B^1(\Q_p,\ad^0_{B_2})$  of~$Z^1(\Q_p,\ad^0_{B_2})\rightarrow H^1(\Q_p,\ad^0_{B_2})$.
\item $s_1$ and~$s_2$ generate~$H^1(\Q_p,k)$, where~$s_1$ is unramified and~$s_2$ is ramified.
\item $t_1$ and~$t_2$ generate~$H^1(\Q_p,\varepsilonbar)$, where~$t_1$ is
 tr\`{e}s ramifi\'{e}e and~$t_2$ is peu ramifi\'{e}e.
\item Under the alternating cup product pairing
$$H^1(\Q_p,k) \times H^1(\Q_p,\varepsilonbar)
\rightarrow H^2(\Q_p,\varepsilonbar) \simeq k,$$
we have
$s_i \cup t_j = \delta_{ij}$.
\end{enumerate}
We now define the dual basis of
$$\m/(\m^2,\lambda) = \Hom_k(\Hom_k(\m/\m^2,k),k)$$
to be given by~$x_i$, $y_i$, and~$z_i$ for~$i = 1$ for~$x_i$ and~$i = 1,2$ for~$y_i$ and~$z_i$, where
$$x_i(r_j) = y_i(s_j) = z_i(t_j)  = \delta_{ij},$$
and all other combinations vanish.
The representation~$\ad^0_{B_2}(\rhobar)$  has a Lie algebra structure via the map
$$\ad^0_{B_2}(\rhobar) \times \ad^0_{B_2}(\rhobar) \rightarrow \ad^0_{B_2}(\rhobar), \qquad (A,B) \mapsto AB-BA.$$
 The corresponding cup product on cohomology groups composed
with this Lie algebra structure induces a  symmetric bilinear pairing (the bracket cup product)
$$M:Z^1(\Q_p,\ad^0_{B_2})^2  \rightarrow H^1(\Q_p,\ad^0_{B_2})^2  \rightarrow H^2(\Q_p,\ad^0_{B_2}),$$
Writing
$H^1(\Q_p,\ad^0_{B_2}) = H^1(\Q_p,k) \oplus H^1(\Q_p,\varepsilonbar)$, this map can be given explicitly
in our case as follows:
$$M(a,b)  = M((a_1,a_2),(b_1,b_2)) =  2(a_1 \cup b_2 -a_2 \cup b_1).$$
(Note that~$\cup$ is alternating so this map is indeed symmetric.)
As noted by in~\cite[\S1.6]{MR1012172},  the image of the corresponding map
gives the quadratic relations in the deformation ring, which produces
the desired quadratic relation~$g$.

More precisely, note that $H^2(\Q_p,\ad^0_{B_2})=H^2(\Q_p,\varepsilonbar)$ is
1-dimensional. By~\cite[Lem.\ 5.2]{MR3406170}, the relation~$g$ 
can be determined (up to the required order) by the
relation given by the image of the natural map
\[H^2(\Q_p,\ad^0_{B_2})^\vee \to (Z^1(\Q_p,\ad^0_{B_2})^\vee)^2\]
induced by the bracket cup product.  But now the non-zero terms can be read off
from the explicit form of~$M(a,b)$ above and the description of our basis of~$\m/(\m^2,\lambda)$
as a dual basis to the explicit basis of~$Z^1(\Q_p,\ad^0_{B_2})$.
It follows that, after
rescaling, the leading term of~$g$ is given by $y_1z_1+y_2z_2$, as
required.  

Part~(\ref{lemma:notformallysmooth}) is a straightforward  consequences of the presentation 
just determined above. Note that the structure over~$R^{\GL_1}$ and~$\Lambda$
is one again determined by the corresponding map
from~$H^1(\Q_p,\ad^0_{B_2})$ to~$H^1(\Q_p,k)$, and from our explicit
description above this corresponds to our choice of the parameters~$y_1$ and~$y_2$.\end{proof}

We also have:

\begin{lemma} \label{lemma:nonsmoothzero} The points of~$R^{B_2}[1/p]$ which are non-smooth over~$\Lambda$ are --- up to unramified twist ---  crystalline 
extensions of~$\varepsilon^{-1}$ by~$1$.
\end{lemma}

\begin{proof} This is the characteristic zero version of the computation done in the proof of Lemma~\ref{lem:explicit}(\ref{modpcrystalline}), and amounts
to noting that in the Tate duality pairing
$$H^1(\Q_p,\Q_p) \times H^1(\Q_p,\Q_p(1)) \rightarrow \Q_p,$$
the unramified classes in the first group are annihilated exactly by the crystalline extensions in the second.
\end{proof}

We next introduce a class
of partially framed deformation rings, which will allow us to relate framed deformation
rings for different groups.

\subsubsection{Partially framed deformation rings}
Since we are assuming that~$F_v=\Qp$, for the rest of this section we
write~$\Qp$ instead of~$F_v$ and~$\rhobar$ instead
of~$\rhobar|_{G_{F_v}}$.  We shall also henceforth (in this section) write~$R^{B}$
for~$R_v^{B,\alphabetabar_v}$ and~$R^{P}$ for~$R_v^{P}$.
These are framed deformation rings with respect to~$\GSp_4$,
and as such, could also be denoted by~$R^{P,\square}$ and~$R^{P,\square}$ to emphasize the framing. However,
the images of the corresponding Galois representations may always
be conjugated to land in~$B$ or~$P$ respectively. In particular, we may consider deformation
rings in which the image is required to actually land inside these subgroups rather than
land there up to conjugation.

\begin{df} Let~$\cD^{B,\tri}$ and~$\cD^{P,\tri}$
denote the subfunctors consisting of deformations which land inside~$B$
or~$P$ respectively. Let~$R^{B,\tri}$ and~$R^{P,\tri}$
denote the corresponding deformation rings.
\end{df}

(Here the adornment~$\tri$ represents that the framing is all taking place
inside the ``upper right corner'' corresponding to~$B$ or~$P$ respectively.)
The ring~$R^{B,\tri}$ may be identified with the universal framed
deformation of~$\rhobar$ with fixed similitude character thought of
as a \emph{representation to~$B$}. The ring~$R^{P,\tri}$
is not quite the universal deformation ring of~$\rhobar$ to~$P$ (framed in~$P$)
with fixed similitude character, because we are imposing an extra condition
on the restriction of the first two diagonal entries to inertia. On
the other hand, if~$R^{P,\univ} = R^{P,\univ,\square}$ denotes the deformations 
to~$\GSp_4$ of fixed similitude character which may be conjugated
to~$P$ (without imposing this condition on inertia), then there is also a corresponding ring~$R^{P,\univ,\tri}$
which  is the universal~$P$-framed deformation of~$\rhobar$ with
fixed similitude character.
Note that there are tautological maps~$R^{B} \rightarrow R^{B,\tri}$
and~$R^{P} \rightarrow R^{P,\tri}$ respectively.
Let us write~$\ad_{\GSp_4}$, $\ad_{B}$, and~$\ad_{P}$ for the groups~$\ad_{\GSp_4}(\rhobar)$, $\ad_{B}(\rhobar)$, and~$\ad_{P}(\rhobar)$ respectively.
Since~$p$ is odd, there exist corresponding direct factors~$\ad^0_{\GSp_4}$, $\ad^0_{B}$, and~$\ad^0_{P}$ corresponding to deformations
with fixed similitude character.

\begin{lemma} \label{lemma:getitdone} Suppose that~$\rhobar$
  is~$p$-distinguished weight 2 ordinary.
\begin{enumerate}
\item There exists a splitting
$$R^{B,\tri} \rightarrow R^{B,\square}
\rightarrow R^{B,\tri}$$
which realizes~$R^{B} = R^{B,\square}$
as formally smooth over~$R^{B,\tri}$ of relative dimension
$$\dim \ad^{0}_{\GSp_4} - \dim \ad^{0}_{B} = 10 - 6 = 4.$$
\item There exists splittings
$$R^{P,\univ,\tri} \rightarrow R^{P,\univ,\square}
\rightarrow R^{P,\univ,\tri},$$
$$R^{P,\tri} \rightarrow R^{P,\square}
\rightarrow R^{P,\tri}$$
which realize~$R^{P}$ and~$R^{P,\univ}$
as formally smooth over~$R^{P,\tri}$ and~$R^{P,\univ,\tri}$ respectively,
of relative dimension
$$\dim \ad^{0}_{\GSp_4} - \dim \ad^{0}_{P} = 10 - 5 = 5.$$
\end{enumerate}
\end{lemma}

\begin{proof}  As previously noted, the~$p$-distinguished hypothesis implies
the existence (by~\cite[Lem.\ 2.4.6]{cht}) of a \emph{unique} Galois stable filtration~$\Fil^i$ on~$(R^{B,\square})^4$.
In particular, we may choose a splitting of this filtration by a symplectic 
matrix~$M \in \GSp_4(R^{B,\square})$ with the property that~$M \equiv I \mod \m$.
Conjugation by~$M$ induces the desired map from~$\GSp_4$-framed
deformations to~$B$-framed deformations, and thus induces a  splitting from~$R^{B,\tri}$ to~$R^{B}$.
In the~$P$ case, one can additionally choose the splitting such that  the choice
of new vector in~$\Fil^2$ is Galois stable, and then the corresponding conjugate is valued in~$P$.

The~$p$-distinguished hypothesis implies that the maps
$$H^0(\Q_p,\ad^{0}_{P}) \rightarrow H^0(\Q_p,\ad^{0}_{B}) \rightarrow H^0(\Q_p,\ad^{0}_{\GSp_4})$$
are all isomorphisms (see for example the explicit descriptions
of~$\ad^{0}_{P}$ and~$\ad^{0}_{B}$ following the proof of this lemma).  
By construction,  the reduced tangent spaces of~$R^{B,\square}$ and~$R^{B,\tri}$
are given by extensions of~$H^1(\Q_p,\ad^0_{B})$ (in both cases) by~$B^1(\Q_p,\ad^0_{\GSp_4})$
and~$B^1(\Q_p,\ad^0_{B})$ respectively (and analogously with~$P$). 
On the other hand, the map on~$B^1$ groups is precisely dual to the map
$$\ad^0_{B}/H^0(\Q_p,\ad^0_{B}) \rightarrow  \ad^0_{\GSp_4}/H^0(\Q_p,\ad^0_{\GSp_4})$$
(and once more similarly with~$P$).
Hence, from the identification of~$H^0$ groups above, it follows that
the map on
reduced tangent spaces corresponding to $R^{B,\square}\to
R^{B,\tri}$ is an injection whose cokernel has
dimension $$(\dim \ad^{0}_{\GSp_4} - \dim H^0(\Q_p,\ad^{0}_{\GSp_4})) 
- (\dim \ad^{0}_{B} - \dim H^0(\Q_p,\ad^{0}_{B}))$$
$$= \dim \ad^{0}_{\GSp_4} - \dim \ad^{0}_{B}$$
(And similarly in the~$P$ case with~$B$ replaced by~$P$.)

We now prove the maps are formally smooth, which will be a direct consequence of the fact that the
obstruction group is given (for~$R^{B,\square}$ and~$R^{B,\tri}$ or for~$R^{P,\univ,\square}$ and~$R^{P,\univ,\tri}$ and~$R^{P,\square}$ and~$R^{P,\tri}$) by the groups~$H^2(\Q_p,\ad^{0}_{B})$
and~$H^2(\Q_p,\ad^{0}_{P})$ respectively.
We consider first the case of~$B$; for simplicity of notation, we drop
$B$ from the superscripts from now on. Consider a  surjection~$\RRR:=R^{\tri}[[x_1,x_2,x_3,x_4]]
\rightarrow R^{\square}$ 
which induces an isomorphism on reduced
tangent spaces, and let~$J$ denote the kernel (so it suffices to show that~$J = 0$). 
Let~$\wm$ be the radical of~$\RRR$. 
Recall that we have a unique symplectic filtration~$\Fil^i$ on~$(R^{\square})^4 = (\RRR/J)^4$ and a choice of splitting corresponding to the matrix~$M$. Lift
this to a filtration~$\widetilde{\Fil}^i$ for~$\RRR/\wm J$, and consider
 a  corresponding set theoretic
deformation~$\widetilde{\rho}: G_{\Q_p} \rightarrow \GSp_4(\RRR/\wm J)$
which preserves this filtration. (There are no issues lifting filtrations because
the symplectic group is formally smooth.)
The corresponding~$2$-cocycle~$[c] \in H^2(\Q_p,\ad^{0}_{\GSp_4}) \otimes J/\wm J$
then lands in~$H^2(\Q_p,\ad^{0}_{B}) \otimes J/\wm J$. 
Now choose a symplectic splitting of this filtration lifting the one for~$\Fil^i$.
Conjugating~$\widetilde{\rho}$
 by the corresponding matrix $\widetilde{M}$  (lifting~$M$ above) gives a set theoretic map~$\widetilde{M} \widetilde{\rho} \widetilde{M}^{-1}$
 from~$G_{\Q_p}$ to~$B(\RRR/\wm J)$.
But this map lifts~$\rho^{\tri,\univ}: G_{\Q_p} \rightarrow B(R^{\tri})$. By universality of~$\rho^{\tri,\univ}$,
 since~$\RRR/\wm J$ is an~$R^{\tri}$-algebra, there is no obstruction to lifting this to a~$B$-representation of~$G_{\Q_p}$, 
 and hence the class~$[c]$ becomes trivial in the corresponding obstruction
group for the~$B$-deformation problem. Since the obstruction group in this case is~$H^2(\Q_p,\ad^0_{B})$ for both the~$B$-deformation
problem and the ordinary~$\GSp_4$-deformation problem, it follows
that~$[c]$ is trivial and hence that~$J = 0$.

 The same argument applies to~$P$, except now the splitting of~$\widetilde{\Fil}^i$
has to be chosen so that it is preserved by~$G_{\Q_p}$
--- equivalently, an identification of the first two eigenspaces to~$\RRR/J$ to ensure
 that the deformation is of~$P$-type. In the case of~$R^{P}$,
 one additionally requires the set theoretic lift to act diagonally after restriction to inertia
 on~$\widetilde{\Fil}^2$.
\end{proof}

\subsubsection{The $\GSp_4$-deformation rings }We are assuming that~$\rhobar|_{G_{F_v}}$ has image of the form
     \[ \begin{pmatrix}
      \lambda_{\alphabar}&0&  \varepsilonbar^{-1} \lambda^{-1}_{\betabar}  \eta_{\alpha \beta} & \varepsilonbar^{-1}  \lambda^{-1}_{\alphabar}  \eta_{\alpha^2}   \\
      0 &\lambda_{\betabar}& \varepsilonbar^{-1}    \lambda^{-1}_{\betabar}  \eta_{\beta^2} & \varepsilonbar^{-1} \lambda^{-1}_{\alphabar} \eta_{\alpha \beta}  \\
      0&0&\varepsilonbar^{-1}\lambda_{\betabar}^{-1}&0\\
      0&0&0&\varepsilonbar^{-1}\lambda_{\alphabar}^{-1}
    \end{pmatrix}\]
   where $\alphabar\ne \betabar$, and where we write~$\eta_{\delta}$ to denote a  (possibly zero) 
   class in~$H^1(\Q_p,\varepsilonbar \lambda_{\overline{\delta}})$.
    Our analysis of the deformation rings (particularly in the~$B$ case) will depend on which of these classes are equal to zero or not.

    The dimensions of~$B$ and~$P$ are~$7$ and~$6$
    respectively.
    Recall that we  are considering deformations of~$\rhobar$ to~$B$ or~$P$ with fixed similitude character. For~$p > 2$, the
    adjoint representations~$\ppp$ and~$\bbb$ admit a splitting with a canonical one dimensional summand corresponding
    to varying the similitude character. Let~$\ad^0_B(\rhobar) \subset \bbb$ and~$\ad^0_P(\rhobar) \subset \ppp$ denote
the    complementary ~$6$
    and~$5$ dimensional subspaces. Explicitly, ~$\bbb$ is given as follows:
        $$ \mathfrak{b} = \left(\begin{matrix} \nu + x_{\alpha} &  -x_{\alpha/\beta}  &  x_{\alpha \beta} & x_{\alpha^2} \\ 0 & \nu + x_{\beta} & 
        x_{\beta^2} &  x_{\alpha \beta} \\ 0 & 0 &   \nu-x_{\beta}  & x_{\alpha/\beta}  \\  0 & 0 & 0& \nu-x_{\alpha} \end{matrix} \right),$$
where the subspace  with~$x_{\alpha/\beta} = 0$ corresponds to~$\ppp$, and the subspace~$\nu = 0$ corresponds to~$\ad^0_B$.
With respect to the basis given by the matrices corresponding to~$\{x_{\alpha^2},x_{\beta^2},x_{\alpha \beta},x_{\alpha/\beta},x_{\alpha},x_{\beta}\}$, the Galois
representation~$\ad^0_B(\rhobar)$ is given explicitly as follows:
$$\left(\begin{matrix}
\varepsilonbar  \lambda_{\alphabar^2} 	&0			&0				& 2  \lambda_{\alphabar} \lambda^{-1}_{\betabar}  \cdot \eta_{\alpha \beta}		& -2  \eta_{\alpha^2}		&0							\\
0		     	&  \varepsilonbar  \lambda_{\betabar^2}	&0				& 0								&0					&-2 \eta_{\beta^2}		\\
0			&0			&  \varepsilonbar \lambda_{\alphabar \betabar} 	& \lambda_{\alphabar} \lambda^{-1}_{\betabar}   \cdot  \eta_{\beta^2}				& -  \eta_{\alpha \beta}	&- \eta_{\alpha \beta}		\\
0			&0			& 0				& \lambda_{\alphabar} \lambda^{-1}_{\betabar} 						&0					&0							\\
0			&0			& 0				&0								&1					&0							\\
0			&0			& 0				&0								&0					&1							\end{matrix} \right).$$
and,  on the space~$\ad^0_P(\rhobar) \subset \ad^0_B(\rhobar)$ with
respect to the basis
$\{x_{\alpha^2},x_{\beta^2},x_{\alpha \beta},x_{\alpha},x_{\beta}\}$  (not a direct summand!), we have
$$\left(\begin{matrix}
\varepsilonbar  \lambda_{\alphabar^2} 	&0			&0					& -2  \eta_{\alpha^2}		&0							\\
0		     	&  \varepsilonbar  \lambda_{\betabar^2}	&0											&0					&-2 \eta_{\beta^2}		\\
0			&0			&  \varepsilonbar \lambda_{\alphabar \betabar}			& -  \eta_{\alpha \beta}	&- \eta_{\alpha \beta}		\\
0			&0			& 0										&1					&0							\\
0			&0			& 0										&0					&1							\end{matrix} \right).$$

As in the case of~$\GL_2$ above, we may compute~$H^2(\Q_p,\ad^0(\rhobar))$ for~$P$ and~$B$ by counting
whether the subspaces generated by~$\{x_{\alpha^2},x_{\beta^2},x_{\alpha \beta}\}$ generate~$\varepsilonbar$
subspaces and whether these subspaces split. The following lemma is immediate from the explicit description above.

\begin{lem}\label{lem: H2 for B} The dimension of~$H^2(\Q_p,\ad^0_B(\rhobar))$ is zero unless one of the following holds:
\begin{enumerate}
\item  \label{B:caseone} The classes~$\eta_{\alpha \beta}$ and~$\eta_{\alpha^2}$ are both zero, and~$\alphabar^2 = 1$.
\item  \label{B:casetwo}The class~$\eta_{\beta^2}$ is zero, and~$\betabar^2 = 1$. In this case, either:
\begin{enumerate}
\item  \label{B:casetwoa} The conditions of part~$(1)$ also hold, or:
\item  \label{B:casetwob}  The dimension of~$H^2(\Q_p,\ad^0_B(\rhobar))$ is~$1$, and there is a~$\Q_p$-equivariant map from the representation~$V$
underlying~$\rhobar$ to the Borel of~$\GL(2)$ corresponding to the representation
$$W = \left( \begin{matrix}  \lambda_{\betabar} & 0 \\ 0 & \varepsilonbar^{-1} \lambda^{-1}_{\betabar}
\end{matrix} \right) = \lambda_{\betabar} \otimes \left( \begin{matrix} 1 & 0 \\ 0 & \varepsilonbar^{-1}
\end{matrix} \right),$$
and the corresponding map relating~$H^2(\Q_p,\ad^0_B(\rhobar)) = H^2(\Q_p,\ad^0_B(V))$
to~$H^2(\Q_p,\ad^0_{B_2}(W))$ is an isomorphism.
\end{enumerate}
\item  \label{B:casethree} The classes~$\eta_{\alpha \beta}$  and~$\eta_{\beta^2}$ are both zero, and~$\alphabar \betabar = 1$.
\end{enumerate}
Moreover, the dimension of~$H^2(\Q_p,\ad^0_B(\rhobar))$ is~$\ge 2$ only in case~$($\ref{B:casetwoa}$)$, in which
it has dimension~$2$.
\end{lem}

We could give a similar (but easier) computation of~$H^2(\Q_p,\ad^0_P(\rhobar))$, but it is not needed in the sequel so it is omitted.

Recall that~$R^{P,\univ}$ denotes the universal deformation ring for~$P$.
The quotient~$R^{P}$ is given by imposing the condition that the action of inertia on~$\Fil^2$ (given by the upper left~$2 \times 2$ matrix after
changing basis) is through a scalar.  Recall that we also
have corresponding rings~$R^{P,\univ,\tri}$ and~$R^{P,\tri}$
where the image lands in~$P$ directly (rather than up to conjugation).
Any deformation of type~$R^{P,\univ,\tri}$ determines deformations of the three
2-dimensional subquotients of~$\rhobar$, given respectively by
the extension~$\rbar_A$ of~$\varepsilonbar^{-1}\lambda_{\alphabar}^{-1}$
by~$\lambda_{\alphabar}$, by the extension~$\rbar_B$ of ~$\varepsilonbar^{-1}\lambda_{\betabar}^{-1}$
by~$\lambda_{\betabar}$, and the extension~$\rbar_{AB}$
of~$\varepsilonbar^{-1}\lambda_{\betabar}^{-1}$
by~$\lambda_{\alphabar}$.  
Similarly, any triple of such deformations with the appropriate coincidences
of the corresponding characters defines a representation of type~$P$. (These identifications require
that we work with~$\tri$ framings rather than~$\square$ framings, since otherwise there would be superfluous framing variables in this
identification.)

Let~$R_A = R^{B_2,\tri}$ for~$\rbar_A$ and~$R_B = R^{B_2,\tri}$ for~$\rbar_B$.
Let~$R_{AB} = R^{B_2,\tri,\kern+0.15em{\det}}$ for~$\rbar_{AB}$,  where~$R^{B_2,\tri,\kern+0.15em{\det}}$ is the framed~$B_2$
deformation ring in which one does not fix the determinant, so (since~$p>2$)
one has that~$R^{B_2,\tri,\kern+0.15em{\det}} = R^{B_2,\tri} \wotimes_{\cO} R^{\GL_1}$ is formally
smooth over~$R^{B_2,\tri}$ of relative dimension~$2$.
There are natural maps from~$R_A$, $R_B$, and~$R_{AB}$ to~$R^{P,\univ,\tri}$
and~$R^{P,\tri}$ respectively. Write $R^{\GL_1 \times
  \GL_1}=R^{\GL_1}\wotimes_{\cO}R^{\GL_1}$ for the deformation ring
corresponding to the pair of characters~$(\lambda_{\alpha_v},\lambda_{\beta_v})$.  We have the following:

\begin{lemma} \label{reducetotensor} The ring~$R^{P,\univ}$ is formally
smooth over
$$R^{P,\univ,\tri} \simeq (R_A \wotimes_{\OL} R_B) \wotimes_{R^{\GL_1 \times \GL_1}} R_{AB}$$
of relative dimension~$5$.
The ring~$R^{P}$ is formally smooth over
$$R^{P,\tri}  \simeq (R_A \wotimes_{\Lambda} R_B) \wotimes_{R^{\GL_1 \times \GL_1}} R_{AB}$$
of relative dimension~$5$.
\end{lemma}

\begin{proof}  The isomorphisms follow directly from the discussion above,
and the statement about formal smoothness 
is~Lemma~\ref{lemma:getitdone}.
\end{proof}

We now prove the~$P$-part of
Proposition~\ref{prop: local ordinary deformation ring
    dimension}.

\begin{prop}  \label{prop:actuallydiddothis}
Suppose that~$\rhobar$ is~$p$-distinguished. Then
$R^{P,\univ}$ and~$R^{P}$ are both complete intersections.
Moreover, they  are connected
in characteristic zero, and the non-smooth locus in characteristic zero has codimension
at least two. In particular, the generic fibres~$R^{P,\univ}[1/p]$ and~$R^{P}[1/p]$ are irreducible
of dimensions~$15$ and~$14$ over~$\Q_p$ respectively.
\end{prop}

\begin{proof}
The strategy is as follows. By Lemma~\ref{reducetotensor}, we can immediately
reduce to the rings~$R^{P,\univ,\tri}$ and~$R^{P,\tri}$ respectively
and prove that they satisfy the same properties above (with~$15$ and~$14$
replaced by~$10$ and~$9$ respectively).
Given the explicit form of the presentations for the~$2$-dimensional~$B_2$-deformation
rings, in order to show that~$R^{P,\univ,\tri}$
and~$R^{P,\tri}$  are complete intersections,  one can simply
write down enough about the equations for the  tensor products in Lemma~\ref{reducetotensor} and observe
(for the appropriate value~$d = 10$ or~$9$ in either case)
that they are either:
\begin{enumerate}
\item Formally smooth of the relative dimension~$d$ over~$\OL$,
\item Given as a quotient of a power series ring in~$d+1$ variables by one relation,
\item Given as a quotient of a power series ring in~$d+2$ variables by
a~$2$-generator prime ideal which is not contained in~$(\lambda)$.
\end{enumerate}
(The last example occurs only in a single case.)
We say more about this computation below.

For the remaining claims, it suffices to prove that the generic fibre
is connected  and  that our tensor products
are~$\Rone$ and~$\Stwo$ (and thus normal); since they are complete
intersections, it is enough to show that the non-smooth points have codimension at least~$2$.
It is convenient to consider two separate cases.

Suppose that~$\alpha \beta = 1$. In this case, it follows
from the~$p$-distinguishedness hypothesis that~$\alpha^2 \ne 1$ and~$\beta^2 \ne 1$.
In this case, the rings above have a particularly simple form even over~$\OL$.
Namely, $R_A$ and~$R_B$ are formally smooth over~$\OL$ and over~$\Lambda$,
and the resulting tensor product is formally smooth over~$R_{AB}$,
and thus the result follows from Lemma~\ref{lem:explicit}, since the rings~$R^{B_2,\tri}$
satisfy all the required geometric properties above.

Now suppose that~$\alpha \beta \ne 1$. In this case, $R_{AB}$ is formally
smooth over $R^{\GL_1\times\GL_1}$, and by
Lemma~\ref{reducetotensor}, $R^{P,\univ}$ and
~$R^{P,\tri}$ are formally smooth over~$R_{A} \wotimes_{\OL} R_B$
or~$R_{A} \wotimes_{\Lambda} R_B$ respectively. 
Let us now consider the case of~$R^{P,\tri}$, which corresponds to~$R_{A} \wotimes_{\Lambda} R_B$,
the case of~$R_A \wotimes_{\cO} R_B$ being easier and  also following immediately from Lemma~\ref{lem:explicit}.
Since~$R_A$ and~$R_B$ are either smooth or have a non-smooth locus of codimension~$4$
(corresponding to twists of~$1 \oplus \varepsilon^{-1}$ by a (possibly
trivial)
unramified quadratic character), it is certainly the case that the points 
on the generic fibre of~$R_{A} \wotimes_{\Lambda} R_B$  which are non-smooth on~$R_A \wotimes_{\OL} R_B$ have codimension at least~$2$. Hence it suffices to consider the non-smooth
points of~$R_{A} \wotimes_{\Lambda} R_B$ which are smooth on~$R_A \wotimes_{\OL} R_B$.
In particular, such a point must have a tangent space of dimension~$8$, and will be smooth
if and only if it has an infinitesimal deformation which
does not lie on~$R_{A} \wotimes_{\Lambda} R_B$.
Equivalently, given a point~$x= (x_A,x_B)$ on the generic fibre of~$R_{A} \wotimes_{\Lambda} R_B[1/p]$, it will be smooth if it has a deformation in which the weight over~$\Lambda$ varies
for one~$x_A$ or~$x_B$ but remains fixed for the other point.
Equivalently, we can look for a deformation of~$x = (x_A,x_B)$ such that
one point is fixed but the other point varies over~$\Lambda[1/p]$.
For~$x_A$ or~$x_B$,  such a deformation
exists as long as~$x_A$ (or~$x_B$) is  a smooth point  over~$\Lambda$. But the non-smooth points
in characteristic zero over 
the space of weights~$\Lambda$ are (up to unramified twist) 
exactly the crystalline extensions of ~$\varepsilon^{-1}$ by~$1$ (see Lemma~\ref{lemma:nonsmoothzero}), and
hence these non-smooth points certainly have codimension at least~$2$.
To show it is connected, it suffices to note that, for each fibre of~$R_A$
above~$\Lambda$, any~$x_A$ is connected over this fibre to a point which
is smooth over~$\Lambda$. This reduces to showing that any extension of~$\varepsilon^{-1}$
by~$1$ which is crystalline has a deformation to a non-crystalline extension.
But this is trivially achieved by a perturbation of the extension class,
noting that~$H^1(\Q_p,\varepsilon)$ is free of rank~$2$ and the crystalline
subspace is a line of rank~$1$.

It remains to prove the claim that these rings are complete intersections in all
the possible cases.
Almost all the time, the tensor product is either immediately seen to be formally
smooth of the right dimension, or given by a single non-zero equation and of
the right dimension. In fact, the only way in which there
can be two equations is when two of the rings~$R_{A}$, $R_{B}$, and~$R_{AB}$
are not formally smooth. This implies that at least two of~$\alphabar^2$, 
$\betabar^2$, and~$\alphabar \betabar$ are equal to one, and this trivially only
happens when~$\alphabar^2 = 1$ and~$\betabar^2 = 1$, and hence~$\alpha \beta \ne 1$. Thus the only possible
case when there exist at least two equations is 
when~$\alphabar^2 = \betabar^2 = 1$ and~$\eta_{\alphabar^2}
=\eta_{\betabar^2} = 0$.
The corresponding tensor product is then
$$\OL \llbracket x_{A,1}, z_{A,1}, z_{A,2}, y_{A,1}, y_{A,2},
x_{B,1}, z_{B,1}, z_{B,2}, y_{B,1} \rrbracket$$
modulo the ideal (noting tensoring over~$\Lambda$ forces~$y_{A,2} = y_{B,2}$):
$$(z_{A,1} y_{A,1} + z_{A,2} y_{A,2}  + \ldots, z_{B,1} y_{B,1} + z_{B,2} y_{A,2}  + \ldots).$$
This pair of elements is easily seen to generate a height~$2$ prime ideal.

The cases when there are no equations and the rings are formally smooth are trivial.
In the cases when there is an extra generator one has to show that the resulting
equation is non-zero. Essentially the most subtle case of this form occurs
when~$\alphabar^2 = \betabar^2 = 1$ and~$\eta_{\alphabar^2}$ and~$\eta_{\betabar^2}$
are both non-zero and peu ramifi\'{e}e. In that case, there
are na\"{\i}vely two equations which have the following form:
$$y_{A,2} =  h(x_{A,1},z_{A,1},z_{A,2},y_{A,1}),$$
$$y_{A,2}  =   h(x_{B,1},z_{B,1},z_{B,2},y_{B,1}),$$
which immediately reduces to one equation.
(Here~$h$ is the same~$h$ because both extensions generate the same line --- the
other cases are trivial). 
We then need to  show that the resulting equation obtained by taking the difference of the RHS
is non-zero. But this is obviously the case unless the RHS is zero. If this is true, then~$y_{A,2}$
is zero in~$R^{B_2,\tri}/\lambda$, which is impossible since~$R^{B_2,\tri}$ is flat over~$\Lambda$.
\end{proof}

\begin{proof}[Proof of Proposition~\ref{prop: local ordinary deformation ring
    dimension}]By Proposition~\ref{prop:actuallydiddothis}, we only
  need to prove the results for~$R^{B}$.  As in the
  proof of Lemma~\ref{lem: irreducibility and dimension of GL2 ordinary deformation
  ring}, we have a presentation of~$R^{B}$ of the
  form~$\cO[[x_1,\dots,x_r]]/(y_1,\dots,y_s)$,
  where
  \[r=11-h^0(\Qp,\ad_{B}\rhobar)+h^1(\Qp,\ad_{B}^0\rhobar),\
    s=h^2(\Qp,\ad_{B}^0\rhobar),\]so that by the local Euler characteristic
  formula, $r-s=10+\dim\ad_{B}^0\rbar=16$.  In particular, if~$H^2(\Q_p,\ad^0_B(\rhobar))=0$, then~$R^{B}$ is
formally smooth over~$\cO$ of relative dimension~$16$, and there
is nothing to prove. 

It is therefore enough to consider each of the cases of
Lemma~\ref{lem: H2 for B}. In case~(\ref{B:casetwob}), we see that~$R^{B}$
is formally smooth of relative dimension~$11$ over the deformation
ring~$R^{B_2,\tri}$
for~$\rbar=\lambda_{\betabar}\oplus\lambda_{\betabar}^{-1}\varepsilonbar^{-1}$,
so the result follows from Lemma~\ref{lem: irreducibility and dimension of GL2 ordinary deformation
  ring}. From the presentation in the previous paragraph, we see that in cases~(1) and~(3), $R^{B}$ is a complete
intersection, while in case~(\ref{B:casetwoa}), we see that every irreducible
component of~$R^{B}$ has relative dimension at least~$16$
over~$\cO$.  By Lemma~\ref{lemma:getitdone}, we may (and we do)
pass freely between~$R^{B}$ and~$R^{B,\tri}$ when convenient.

Suppose that we are in case~(\ref{B:casethree}), and suppose that~$\eta_{\alpha^2}\ne 0$. Let $\varrhobar$ be the representation
with the same $\alphabar,\betabar$ as~$\rhobar$, but with
$\eta_{\alpha^2}=\eta_{\beta^2}=\eta_{\alpha\beta}=0$. Let
$R_{\varrhobar}^{B,\tri}$ be the corresponding deformation ring. We
claim that in fact~$R_{\varrhobar}^{B,\tri}$ and~$R^{B,\tri}$ are
isomorphic. To see this, note firstly that
since~$\alphabar\betabar=1$, and~$\alphabar\ne\betabar$, we
have~$\alphabar^2\ne
1$. Let~\[\rbar=
  \begin{pmatrix}
    \lambda_{\alphabar}&\eta_{\alpha^2}\\0&\lambda_{\alphabar}^{-1}\varepsilonbar^{-1}
  \end{pmatrix}
.\] We have already shown that, in this case, $R^{B_2,\tri}$ is formally smooth
over~$\Lambda$. In fact, we can be more
explicit. Write~$\widetilde{\Lambda}$ for what we called~$R^{\GL_1}$
above, so that~$\widetilde{\Lambda}$ is the formally
smooth~$\Lambda$-algebra of relative dimension~$1$ which carries the
additional information of the actual lift of~$\lambda_{\alphabar}$
(rather than just its restriction to inertia), so that~$R^{B,\tri}$
and~$R^{B_2,\tri}$ are
naturally
$\widetilde{\Lambda}$-algebras. Let~$\widetilde{\lambda_{\alphabar}}:G_{\Qp}\to\widetilde{\Lambda}^\times$
be the 
universal lift of~$\lambda_{\alphabar}$. Then since~$\alphabar^2\ne
1$, $H^1(G_{\Qp},\widetilde{\lambda_{\alphabar}}^2\varepsilon)$ is a
free~$\widetilde{\Lambda}$-module of rank~$1$, and the universal lift of~$\rbar$
is represented by  \[\begin{pmatrix}
    \widetilde{\lambda_{\alphabar}}&\widetilde{\lambda_{\alphabar}}^{-1}\varepsilon^{-1}\widetilde{\eta_{\alpha^2}}\\0&\widetilde{\lambda_{\alphabar}}^{-1}\varepsilonbar^{-1}
  \end{pmatrix}
\] where~$\widetilde{\eta_{\alpha^2}}$ lifts~$\eta_{\alpha^2}$.
It is then easy to verify that if~$\varrho^{\univ}$ is the universal
upper-triangular lift of~$\varrhobar$ to~$R_{\varrhobar}^{B,\tri}$, then
\[ \varrho^{\univ}
  + \left(\begin{matrix}0 &  0    & 0 &  \widetilde{  \lambda_{\alphabar}}^{-1}\varepsilon^{-1}  \widetilde{\eta_{\alpha^2} }  \\ 0 & 0  & 0 & 0 \\ 0 & 0 &  0      & 0 \\  0 & 0 & 0& 0 \end{matrix} \right)\] is a lift
  of~$\rhobar$. This gives a map $R_{\varrhobar}^{B,\tri}\to
  R^{B,\tri}$, and we can obtain a map $R^{B,\tri}\to
  R_{\varrhobar}^{B,\tri}$ in the same way. It is clear that the
  composites of these maps are the identities, so that~$R_{\varrhobar}^{B,\tri}$ and~$R^{B,\tri}$ are
isomorphic, as claimed.

Accordingly, whenever we are in case~(\ref{B:casethree}), we will assume from now on
that~$\eta_{\alpha^2}=0$. It is now easy to see that in each of the
cases~(\ref{B:caseone}), (\ref{B:casetwoa}), and~(\ref{B:casethree}), $\Spec R^{B,\tri}[1/p]$ is
connected. Indeed, in each case we have
$\eta_{\alpha^2}=\eta_{\alpha\beta}=0$, so by arguing as in the proof of Lemma~\ref{lem: irreducibility and dimension of GL2 ordinary deformation
  ring} (using conjugation by~$\diag(t,1,1,t^{-1})$), we see that
every closed point of~$\Spec R^{B,\tri}[1/p]$ may be path connected
to one which lands in~$P(\Qpbar)$. Now we may immediately conclude by 
knowing the corresponding result
for~$\Spec
R^{P,\univ,\tri}[1/p]$ proved in Proposition~\ref{prop:actuallydiddothis}.

To obtain irreducibility we now argue as in the proof of
Lemma~\ref{lem: irreducibility and dimension of GL2 ordinary
  deformation ring}, by studying the singular locus
of~$\Spec R^{B}[1/p]$. More precisely, we
let~$\rho^{\univ}:G_{\Qp}\to\GSp_4(R^{B})$ be the universal
lifting, and let~$H^2:=H^2(G_{\Qp},\ad^0_B\rho^{\univ})$, a finite
$R^{B}$-module, which is compatible with
specialization. Let~$U$ be the complement of the support of~$H^2$
in~$\Spec R^{B}[1/p]$. At any closed point $x\in U$ with
corresponding representation~$\rho_x:G_{\Qp}\to\GSp_4(E_x)$, we have
$H^2(G_{\Qp},\ad^0_B\rho_x)=0$, so it follows that~$U$ is
formally smooth over~$E_x$ at~$x$ of relative dimension~$16$. In
particular, $U$ is regular.

The points in the complement of~$U$ are those for which
$H^2(G_{\Qp},\ad^0_B\rho_x)\ne 0$. We claim that this has codimension
at least~$2$. We may explicitly
describe this locus as follows (this description follows easily from
the explicit description of~$\ad^0_B$ preceding Lemma~\ref{lem: H2 for B}). In case~(1), we may suppose without loss of
generality (by twisting with a quadratic character if necessary) that~$\alphabar=1$, and then the points in the complement
of~$U$ are those conjugate to representations of the form \[
  \begin{pmatrix}
    1 & 0 & 0 &0\\ 0&\chi& * & 0\\0&0&\varepsilon^{-1}\chi^{-1}& 0 \\0&0&0&\varepsilon^{-1}
  \end{pmatrix}
\]where~$\chi$ lifts~$\lambda_{\betabar}$. The locus of such points
has dimension at most~$13$; 
 indeed, the action of $\PGSp_4$ by
conjugation contributes at most~$10$ to the dimension, and the choice
of~$\chi$ and~$*$ at most~$3$ (there is a two-dimensional family of
choices of~$\chi$, and if~$\chi^2$ is non-trivial then the choice of~$*$
gives one more dimension, while if~$\chi^2$ is trivial then it gives~$2$ dimensions).

In case~(\ref{B:casetwoa}), we have in addition the points of the form 
 \[
  \begin{pmatrix}
    \chi' & * & 0 &0\\ 0&\lambda_{-1}& 0 & 0\\0&0&\varepsilon^{-1}\lambda_{-1}& * \\0&0&0&\varepsilon^{-1}(\chi')^{-1}
  \end{pmatrix}
\]where~$\chi'$ lifts~$\lambda_{\alphabar}=1$. The locus of such
points again has dimension at most~$13$.

In
case~(\ref{B:casethree}), we have the points of the form 
\[
  \begin{pmatrix}
    \chi & * & 0 &*\\ 0&\chi^{-1}& 0 & 0\\0&0&\varepsilon^{-1}\chi& * \\0&0&0&\varepsilon^{-1}\chi^{-1}
  \end{pmatrix},
\]where~$\chi$ lifts~$\lambda_{\alphabar}$. The locus of such points
has dimension at most~$14$ (with~$2$ dimensions for the choice of~$\chi$,
and then generically one dimension each for the choices of the extension class
of~$\chi^{-1}$ by~$\chi$ and of~$\varepsilon^{-1}\chi^{-1}$
by~$\chi$, or two dimensions each if~$\chi^2=1$).

Thus in cases~(\ref{B:caseone}) and~(\ref{B:casethree}), since we know that~$R^{B}$ is a
complete intersection, we see that it is normal (being~$\Rone$ and~$\Stwo$), so
we are done. The case~(\ref{B:casetwob}) having already been dealt
with, we are left with case~(\ref{B:casetwoa}), where we have seen that every irreducible
component of~$\Spec R^{B}[1/p]$ has dimension at least~$16$,
while the complement of~$U$ has dimension at most~$12$. It now suffices to show that~$R^{B}$
is a complete intersection, and thus also normal as above, and to
check that $\Spec R^{B}[1/p]$ has dimension exactly~$16$.

We have a presentation of~$R = R^{B}$ of the
form~$\cO[[x_1,\dots,x_{18}]]/(y_1,y_2)$. This is a complete intersection as long as~$\dim(R) \le 19 - 2=17$, which also
  implies that the relative dimension of~$R$ over~$\cO$ is~$16$, and so the dimension
  of the generic fibre is~$16$.
  Assume otherwise, so that $\dim(R) \ge 18$. Then the support of~$R$
  in~$\Spec \cO[[x_1,\dots,x_{18}]]$
  contains a height one prime~$\mathfrak{p}$ of~$\cO[[x_1,\dots,x_{18}]]$.
Suppose firstly that~$\mathfrak{p}$ has residue characteristic zero, and let~$T$ denote the corresponding
closed subscheme of~$\Spec R[1/p]$, which will have dimension~$17$.
 For any closed point $x\in T$   with
corresponding representation~$\rho_x:G_{\Qp}\to\GSp_4(E_x)$, the tangent space at~$x$ certainly has dimension at least~$\dim(T) = 17$. 
Hence there is an inequality 
$$11-h^0(\Qp,\ad_{B}\rho_x)+h^1(\Qp,\ad_{B}^0\rho_x) \ge 17,$$
and so, by the Euler characteristic formula,
$h^2(G_{\Qp},\ad^0_B\rho_x) \ge  17 - 16 \ge 1$. In particular, it follows that~$x$  lies in the 
support of $H^2$, and hence that~$T \subset U$. But
we have already seen that~$U$ has dimension at most~$12$, and this is a contradiction.

Hence~$R$ can only fail to be a complete intersection if the support of~$R$ contains~$(\lambda)$. 
It follows that~$\dim(R/\lambda) = \dim(k[[x_1,\ldots,x_{18}]])$, and hence that~$R/\lambda = k[[x_1,\ldots,x_{18}]]$.
Twisting, we may without loss of generality assume
that~$\betabar=1$. Let~$\rbar = 1 \oplus \varepsilonbar^{-1}$, and let~$R^{B_2,\tri}$ denote the corresponding fixed determinant
deformation ring
to the Borel of~$\GL(2)$. By realizing~$\rbar$ as the subquotient of
the representation~$\rhobar$ given by the span of the second and third
standard basis vectors, 
there is an induced map
  $$\psi: R^{B_2,\tri} \rightarrow R^{\tri} \rightarrow R \rightarrow R/\lambda = k[[x_1,\dots,x_{18}]].$$
  Let~$W$ denote the representation underlying~$\rbar$, and (as
  previously) $V$ the representation underlying~$\rhobar$.
  Let us now consider the induced map on reduced tangent spaces.
  To compute this, we may look at the corresponding deformation rings, and consider
  the induced map on tangent spaces. 
  For~$R^{B_2,\tri}$, the tangent space is given by~$Z^1(\Q_p,\ad^0_{B_2}(W))$. 
  For~$R^{\tri}$, 
  it is given by~$Z^1(\Q_p,\ad^0_{B}(V))$.
  Note that we are assuming that $\eta_{\alpha
  \beta} = \eta_{\alpha^2} = \eta_{\beta^2} = 0$, and so~$\rhobar$ is
completely split, and so~$\ad^0_{B_2}$ is a direct summand
of~$\ad^0_{B}$. Thus $Z^1(\Q_p,\ad^0_{B_2}(W))\to
  Z^1(\Q_p,\ad^0_{B}(V))$ is injective. 
  On the other hand, 
  Lemma~\ref{lem:explicit}~(\ref{lemma:notformallysmooth}) shows that  (for this~$\rbar$) the ring~$R^{B_2,\tri}/\lambda$
  is \emph{not} formally smooth.  But this is a contradiction; a minimal set of
  generators
  of the maximal ideal of~$R^{B_2,\tri}/\lambda$ satisfy at least one polynomial relation,
  but their images under~$\psi$ do not satisfy any such relation under our assumptions
  because the map on tangent spaces is injective and 
  (as we are currently assuming)~$R^{\tri}/\lambda$ and~$R/\lambda$ is formally smooth.
  Hence~$\lambda$ also cannot be in the support of~$R$,
  and thus~$R$ is a complete intersection.
\end{proof}

\begin{rem}The last argument shows that, in case~(\ref{B:casetwoa}), the ring~$R = R^{B}$ is a complete
intersection. But we certainly expect (in this and in all other cases) the stronger properties
that~$R$ is flat over~$\cO$ and~$R/\lambda$ is also a complete intersection, whereas
the argument only shows that~$\dim(R/\lambda) \le 17$, rather than~$\dim(R) - 1 = 16$,
which would be necessary in order for~$\lambda$ to be a regular element. In general, we 
have
often only attempted to prove exactly enough about the deformation rings that we require for
the argument, rather than giving a fuller account of their geometric properties. We apologize
to readers who examine this argument in closer detail who were
hoping for something more comprehensive. 
\end{rem}
As in~\S\ref{subsec: Galois
  deformation rings}, we say that  a closed point~$x$ of
$R_v^{B,\alphabetabar_v}[1/p]$ 
  (resp.\ $R_v^{P}[1/p]$) is smooth if
  $(R_v^{B,\alphabetabar_v}[1/p])_x$ is regular
  (resp.\ $(R_v^{B,\alphabetabar_v}[1/p])_x$ is
  regular). We say that the corresponding Galois
  representation~$\rho_x$ is pure if it is de Rham, and
  if~$\WD(\rho_x)$ is pure (that is, it arises as the base extension
  of a pure Weil--Deligne representation over a number field).
\begin{lem}\label{lem: smooth points of the big local deformation
    rings}If~$x$ is a closed point of 
  the generic fibre of $R_v^{B,\alphabetabar_v}[1/p]$ or
  $R_v^{P}[1/p]$, and~$x$ is pure, then it is smooth.
  \end{lem}
  \begin{proof}
    We first consider the case of~$B$.
     From the proof of Proposition~\ref{prop: local ordinary deformation ring
    dimension}, we see that it is enough to check that
  $H^2(G_{\Qp},\ad^0_B\rho_x)=0$.
  By Tate local duality, this means
  that it is enough to check that
  $\Hom_{G_{\Qp}}(\rho_x,\rho_x(1))=0$, and therefore it is enough to
  check that~$\Hom_{\WD_{\Qp}}(\WD(\rho_x),\WD(\rho_x(1)))=0$. This
  follows easily from the definition of purity.
  The same argument also applies to~$R^{P,\univ}_v[1/p]$.
  We now consider~$R^{P}_v[1/p]$. 
The non-smooth points~$x$ of~$R^{P}_v[1/p]$ are either
non-smooth in~$R^{P,\univ}_v[1/p]$ (for which the previous argument applies)
or, via the isomorphism of Lemma~\ref{reducetotensor}
and the proof of Proposition~\ref{prop:actuallydiddothis}, 
arise in the following way: the representation~$\rho_x$ 
admits a~$2$-dimensional reducible subquotient~$r_x$
such that the corresponding point on the deformation ring~$R^{B_2}[1/p]$  is not smooth over~$\Lambda_v$.
By Lemma~\ref{lemma:nonsmoothzero}, 
 such 
  representations
  are (up to unramified twist) a crystalline extension 
  of~$\varepsilon^{-1}$ by~$1$. Since these
   are not pure (and purity is preserved by
  taking subquotients), the representation~$\rho_x$ is also not pure, and we are also done in this case.
   \end{proof}

\subsection{Local deformation problems, \texorpdfstring{$l\ne p$}{l ne p}}\label{subsec: l not
  p deformations}\subsubsection{Unobstructed deformations}Assume that $v\nmid p$.  
\begin{prop}\label{prop:smoothlifts}
  If $H^0(F_v,\ad^0\rhobar(1)) = 0$, then~$R_v^\square$ is isomorphic to a
  power series ring over $\cO$ in $10$ variables. If
  furthermore~$\rhobar|_{G_{F_v}}$ is unramified, then so are all of
  its lifts.
\end{prop}
\begin{proof}
By Tate duality, the condition is equivalent to $H^2(F_v,\ad^0\rhobar)
= 0$, and the result follows from a standard calculation in
obstruction theory (see e.g.\ \cite[\S 5.2]{MR1643682}).
\end{proof}

\subsubsection{Taylor--Wiles deformations}\label{subsubsec: TW deformations}Assume that
$q_v \equiv 1 \bmod p$, and that both~$\psi|_{G_{F_v}}$ and
$\rhobar|_{G_{F_v}}$ are unramified.  We take $\Lambda_v = \cO$.  We
assume that $\rhobar(\Frob_v)$ has $4$ distinct eigenvalues in~$k$,
and we fix an ordering of them as $\alphabar_1$, $\alphabar_2$,
$\alphabar_3=\psi(\Frob_v)/\alphabar_2$, $\alphabar_4=\psi(\Frob_v)/\alphabar_1$.
For each $i=1,2$, let
$\overline{\gamma}_i: G_{F_v} \rightarrow k^\times$ be the unramified
character that sends $\Frob_v$ to $\alphabar_i$. 

\begin{lem}\label{lem:anyliftisTW}
  Let $\rho: G_{F_v} \rightarrow \GSp_4(A)$ be any lift of $\rhobar$.
  There are unique continuous characters
  $\gamma_i: G_{F_v} \rightarrow A^\times$ for~$i=1,2$, such that
  $\rho$ is $\GSp_4(A)$-conjugate to a lift of the form
  $\gamma_1 \oplus \gamma_2 \oplus \psi\gamma_2^{-1}
  \oplus\psi\gamma_1^{-1}$, where
  $\gamma_i \bmod \frakm_A = \overline{\gamma}_i$ for each
  $i=1,2$.\end{lem}
\begin{proof}
  This can be proved in exactly the same way as~\cite[Lem.\
  5.1.1]{MR2234862}. 
\end{proof}

Let $\Delta_v = k(v)^\times(p)^2$, where $k(v)^\times(p)$ is the
maximal $p$-power quotient of $k(v)^\times$, and let
$\rho: G_{F_v} \rightarrow \GSp_4(R_v^\square)$ denote the
universal lift. Then $\rho$ is $\GSp_4(R_v^\square)$-conjugate to a
lift of the form    $\gamma_1 \oplus \gamma_2 \oplus \psi\gamma_2^{-1}
  \oplus\psi\gamma_1^{-1}$ as in Lemma~\ref{lem:anyliftisTW}.  For
  $i=1$, $2$, the character
$\gamma_i \circ \Art_{F_v}|_{\cO_{F_v}^\times}$ factors through
$k(v)^\times(p)$, so we obtain a canonical local $\cO$-algebra
morphism $\cO[\Delta_v] \rightarrow R_v^\square$.  Note that this
depends on the choice of ordering $\alphabar_1,\ldots,\alphabar_4$.  It is
straightforward to check that this morphism is formally smooth of
relative dimension $10$.

\subsubsection{Ihara
  avoidance deformations}Let~$v$ be a finite place of~$F$ with $q_v
\equiv 1 \bmod p$. Assume further that $\rhobar|_{G_{F_v}}$ is
trivial, and that~$\psi|_{G_{F_v}}$ is unramified and has trivial
reduction modulo $\lambda$.  We
take $\Lambda_v = \cO$.

Let $\chi = (\chi_1,\chi_2)$ be a pair of continuous characters
$\chi_i: \cO_{F_v}^\times \rightarrow \cO^\times$ that are
trivial modulo $\lambda$.
We let $\cD_v^\chi$ be the functor of lifts $\rho: G_{F_v}
\rightarrow \GSp_4(A)$ such that for all~$\sigma\in I_{F_v}$, the
characteristic polynomial of~$\rho(\sigma)$
is~\[(X-\chi_1(\Art_{F_v}^{-1}(\sigma)))(X-\chi_2(\Art_{F_v}^{-1}(\sigma)))(X-\chi_2(\Art_{F_v}^{-1}(\sigma))^{-1})(X-\chi_1(\Art_{F_v}^{-1}(\sigma))^{-1}).\]

Then $\cD_v^\chi$ is a local deformation problem, and we denote its
representing object by $R_v^\chi$.  
\begin{lem}\label{lem: chi Ihara avoidance deformation ring is smooth}
If $\chi_1,\chi_2\ne1$ and $\chi_1\ne \chi_2^{\pm1}$, then
every closed point of $\Spec R_v^\chi[1/p]$ is smooth.
\end{lem}
\begin{proof} We can choose $\sigma\in I_{F_v}$ with $\chi_1(\Art_{F_v}^{-1}(\sigma))$,
  $\chi_2(\Art_{F_v}^{-1}(\sigma))$, $\chi_1(\Art_{F_v}^{-1}(\sigma))^{-1}$, $\chi_2(\Art_{F_v}^{-1}(\sigma))^{-1}$
  pairwise distinct. As in Lemma~\ref{lem: smooth points of Ihara1}, we
  need to check that for every point~$x$, we have
  $\Hom_{E'[G_{F_v}]}(\rho_x,\rho_x(1))=0$.  Any such homomorphism
  would have to respect the eigenspaces for~$\rho_x(\sigma)$, and must
  therefore be zero.
\end{proof}

The proof of the following two results occupies the rest of this
subsection.
\begin{prop}\label{prop:Ihara1ring}
 Assume that $\chi_1=\chi_2 = 1$.
 Then $R_v^1$ satisfies the following properties: 
 \begin{enumerate}[noitemsep,topsep=3pt]
  \item $\Spec R_v^1$ is equidimensional of dimension $11$ and every generic point has characteristic zero.
  \item Every generic point of $\Spec R_v^1/(\lambda)$ is the
    specialization of a unique generic point of $\Spec R_v^1$.\end{enumerate}
\end{prop}

\begin{prop}\label{prop:Iharachiring}
 Assume that $\chi_1,\chi_2\ne1$ and $\chi_1\ne \chi_2^{\pm1}$. Then $\Spec R_v^\chi$ is irreducible of dimension $11$, and its generic point has characteristic zero. 
\end{prop}

We follow the strategy of~\cite{tay} (which proves the corresponding
results for~$\GL_n$) closely.
A source of minor complications in the case of $\GSp_4$ is that
nilpotent centralizers need not be connected.  Even though we are
interested only in deformation rings with fixed multiplier, we have
found it more convenient to carry out the analysis without fixing
multipliers until the end.  We also take advantage of the fact that we
only care about $\GSp_4$ (rather than, say, $\GSp_{2g}$) to be a bit more ad hoc in our arguments.

Throughout the rest of this section, $q$ will denote an integer which is not a
multiple of $p$. 

\subsubsection{Preliminaries on nilpotent matrices}Let $\cU\subset\GSp_4/\cO$ be the closed subscheme of matrices with characteristic polynomial $(X-1)^4$, and  let $\cN\subset\Lie(\GSp_4)$ be the closed subscheme of matrices with characteristic polynomial $X^4$.

In~\cite{tay}, under the assumption that~$p\ge n$, Taylor uses
truncations to degree~$X^{n-1}$ of the usual exponential and
logarithm maps in order to relate unipotent and nilpotent matrices (see in particular~\cite[Lem.\ 2.4]{tay}). For~$p>3$, we
could in the same way use the
truncations to order~$X^3$ of the usual exponential and logarithmic maps. However, both~$\exp$ and~$\log$
to third order involve terms of the form~$X^3/3!$ and~$(X-1)^3/3$, which
we need to avoid when working in residue characteristic three. 
In the proof of~\cite[Lem.\ 3.15]{jack} an alternative approach
is given (again in the case of~$\GL_n$), using
 the maps~$\exp_1 = 1 + N$ and~$\log_1 = (U-1)$ in order to avoid assumptions on the characteristic. 
However, neither the matrices~$I + N$ for nilpotent~$N$ nor~$U-1$ for unipotent~$U$
will in general be symplectic, and thus our truncated exponential and logarithm
maps 
must be at least quadratic. This motivates the following definitions.

For~$p \ge 3$,
we have the following modified versions of the exponential and logarithm map,
which are the same as the usual definitions up to and including order~$X^2$:
\begin{align*}
\exp_2:\cN&\to\cU\\
N&\mapsto I+N+\frac{N^2}{2}+\frac{N^3}{2}
\end{align*}
and
\begin{align*}
\log_2:\cU&\to\cN\\
U&\mapsto (U-I)-\frac{(U-I)^2}{2}.
\end{align*}
It is easily verified that these maps do indeed have image~$\cU$,
respectively~$\cN$, and that they are in fact inverses to each other,
and in particular are bijective. Additionally, they
 commute with the conjugation action of $\GSp_4$, and for~$m \in \Z$ satisfy
\begin{equation*}
\exp_2\left(mN + m^* N^3\right)=\exp_2(N)^m,\quad\log_2(U^m)=m\log_2(U) + m^* \log_2(U)^3,
\end{equation*}
where~$m^* = \displaystyle{(m - m^3)/3} \in \Z$.

We define the following elements of $\mathcal{N}(\cO)$:
\begin{equation*}
N_0=0, \ N_1=\begin{pmatrix}0&0&0&1\\0&0&0&0\\0&0&0&0\\0&0&0&0\end{pmatrix}\kern-0.3em{,} \  \  N_2=\begin{pmatrix}0&0&1&0\\0&0&0&1\\0&0&0&0\\0&0&0&0\end{pmatrix}\kern-0.3em{,}   \ \ N_3=\begin{pmatrix}0&1&0&0\\0&0&1&0\\0&0&0&-1\\0&0&0&0\end{pmatrix}\kern-0.3em{.}
\end{equation*}

We also let $\cN_i\subset\cN$ be the reduced, locally closed subscheme
consisting of nilpotent matrices of rank $i$, so that $N_i\in\cN_i$.

The following is an analogue of~\cite[Lem.\ 2.5]{tay}.
\begin{prop}\label{prop: nilpotent conjugacy classes} \leavevmode
\begin{enumerate}
\item \label{firstitemGB} $Z_{\GSp_4}(N_i)$ is a smooth group scheme
  over $\Spec \cO$ with fibres of dimensions $11,7,5,3$ for
  $i=0,1,2,3$.  Each connected component of~$Z_{\GSp_4}(N_i)$ is irreducible with irreducible special fibre.  Moreover, $Z_{\GSp_4}(N_i)$~is irreducible except when $i=2$, in which case it has  two components.
\item \label{seconditemGB}  Locally in the \'{e}tale topology, the
  universal nilpotent matrix over $\cN_i$ is conjugate to~$N_i$ by a section of $\GSp_4$.
\item \label{thirditemGB}  $\cN_i$ is smooth over $\Spec \cO$ with irreducible fibres of
  dimensions $0,4,6,8$ for $i=0,1,2,3$. In particular, $\cN_i$ is irreducible.
\end{enumerate}
\end{prop}

\begin{proof} Part~(\ref{firstitemGB}) can be checked by brute force calculation.  For instance in the most interesting case when $i=2$ a direct computation (using that $p>2$) shows that
\begin{equation*}
Z_{\GSp_4}(N_2)\simeq\cO[x,y,z,w,\alpha,\beta,\gamma,\delta,(wx-yz)^{-1}]/(xy,wz,y\gamma-w\alpha-x\delta-z\beta)
\end{equation*}
where the matrix is given by
\begin{equation*}
\begin{pmatrix}
x&y&\alpha&\beta\\
z&w&\gamma&\delta\\
0&0&x&y\\
0&0&z&w
\end{pmatrix}
\end{equation*}
and from this all the properties are clear (for instance, the two components are given by $x=w=0$ and $y=z=0$).

For part~(\ref{seconditemGB}), we explain the case when $i=2$.  The others are similar
but easier. We may view the universal nilpotent $N$ over $\cN_2$ as an endomorphism of $\cO_{\cN_2}^4$ with the ``standard'' symplectic form $\psi$.  Then, by the definition of $\cN_2$, $\ker(N)$ is a local direct summand of rank 2.  Then one checks that
\begin{align*}
\psi':\cO_{\cN_2}^4/\ker N\times\cO_{\cN_2}^4/\ker N&\to\cO_{\cN_2}\\
(v,w)&\mapsto \psi(Nv,w)
\end{align*}
is a well defined non-degenerate symmetric pairing.

\'{E}tale locally, one may trivialize $\psi'$:  For any point $x\in\cN_2$ we may pick a Zariski open neighbourhood $x\in U=\Spec A\subset\cN_2$ over which $\cO_{\cN_2}^4/\ker N$ has a basis $f_1,f_2$ with $\psi'(f_1,f_2)=0$ and $\psi'(f_1,f_1),\psi'(f_2,f_2)\in A^\times$.  Let $A'$ be the \'{e}tale $A$-algebra $A[\sqrt{\psi'(f_1,f_1)},\sqrt{\psi'(f_2,f_2)}]$, so that over $U'=\Spec A'$, $(\cO_{\cN_2}^4/\ker N)_{U'}$ has a basis $f_1'=f_1/\sqrt{\psi'(f_1,f_1)}$, $f_2'=f_2/\sqrt{\psi'(f_2,f_2)}$ with $\psi'(f_1',f_1')=\psi'(f_2',f_2')=1$ and $\psi'(f_1',f_2')=0$.
Now lift $f_1'$ and $f_2'$ to sections $e_1$ and $e_2$ of $\cO_{U'}^4$.  We may further arrange that $\psi(e_1,e_2)=0$ by replacing $e_2$ by $e_2-\psi(e_1,e_2)Ne_1$.  Then $Ne_2,Ne_1,e_1,e_2$ forms a symplectic basis for $\cO_{U'}^4$, and if we let $g\in\GSp_4$ have these elements as columns, then $N_{U'}=gN_2g^{-1}$.

Finally we turn to part~(\ref{thirditemGB}). For each $i$, there is a map
\begin{align*}
\GSp_4&\to\cN_i\\
g&\mapsto gN_ig^{-1}.
\end{align*}
By the first two parts of the proposition, this map is smooth and
surjective.  Indeed, it suffices to check this after base change to a suitable \'{e}tale cover $U\to \cN_i$, over which it becomes isomorphic to $Z_{\GSp_4}(N_i)_U \to U$.  
It follows that $\cN_i$ is smooth over $\cO$.  The fibres of $\cN_i$ are irreducible because those of $\GSp_4$ are, and the statement about dimensions follows from the computation of the dimensions of the fibres of $\GSp_4\to\cN_i$ in part~(\ref{firstitemGB}).
\end{proof}

\begin{rem}
By contrast to the situation for $\GL_n$ considered in~\cite{tay}, it is no longer the case that $Z_{\GSp_4}(N_i)$ is connected, nor is it true that the universal matrix over $\cN_i$ is Zariski locally conjugate to $N_i$ (both fail when $i=2$).
\end{rem}

\subsubsection{Some spaces of polynomials}
Let $\tilde{\cP}=\Gm^3$ be the diagonal torus in~$\GSp_4$; we somewhat
abusively write
\begin{equation*}
\tilde{\cP}=\{(X-\alpha)(X-\beta)(X-\gamma\beta^{-1})(X-\gamma\alpha^{-1})\}
\end{equation*}where the order of the linear factors matters, and we
let
\begin{equation*}
\cP=\tilde{\cP}/W=\{X^4+a_3X^3+a_2X^2+a_1X+a_0\mid a_0\in\Gm,a_3^2a_0=a_1^2\},
\end{equation*}
so that there is a finite map $\pi:\tilde{\cP}\to\cP$, given by
multiplying out the linear factors.  We consider some reduced closed subspaces of $\cP$:
\begin{align*}
\cP_0&=\cP\\
\cP_1&=\pi(\{(X-\alpha)(X-\beta)(X-\gamma\beta^{-1})(X-\gamma\alpha^{-1})\mid \gamma\alpha^{-1}=q\alpha\})\\
\cP_2&=\pi(\{(X-\alpha)(X-q\alpha)(X-\gamma q^{-1}\alpha^{-1})(X-\gamma\alpha^{-1})\})\\
\cP_3&=\pi(\{(X-\alpha)(X-q\alpha)(X-q^2\alpha)(X-q^3\alpha)\})
\end{align*}

We will find it useful to consider some explicit elements of $\GSp_4(R)$, for an $\cO$-algebra $R$.  For $\alpha,\beta,\gamma\in R^\times$ we let
\begin{align*}
\Phi_0(\alpha,\beta,\gamma)&=\diag(\alpha,\beta,\gamma\beta^{-1},\gamma\alpha^{-1})\\
\Phi_1(\alpha,\beta)&=\diag(q\alpha,\beta,q\alpha^2/\beta,\alpha)\\
\Phi_{2,a}(\alpha,\gamma)&=\diag(q\alpha,\gamma\alpha^{-1},\alpha,\gamma q^{-1}\alpha^{-1})\\
\Phi_{2,b}(\alpha,\beta)&=\begin{pmatrix}0&q\alpha&0&0\\ q\beta&0&0&0\\ 0&0&0&\alpha\\0&0&\beta&0\end{pmatrix}\\
\Phi_{3}(\alpha)&=\begin{pmatrix}q^3 \alpha &  0 & \frac{q(1-q^2)}{6}  \alpha &0\\ 0 & q^2  \alpha &0&\frac{(1-q^2)}{6}  \alpha \\ 0&0&q \alpha & 0 \\0&0&0& \alpha \end{pmatrix}
\end{align*}

\subsubsection{Spaces of matrices}
We define $\cN(q)$ to be the closed subscheme of $\GSp_4\times\cN$ consisting of pairs $(\Phi,N)$ satisfying 
$$\Phi  N\Phi^{-1}=\log_2(\exp_2(N)^q) = qN + q^* N^3,$$
where as above we write~$q^* = (q-q^3)/3$. 
This definition is motivated by the following. The actual equation we wish to study has the form
$$\Phi U \Phi^{-1} = U^q$$
 for a unipotent matrix~$U$.
If we let~$N = \log_2(U)$, we have~$U = \exp_2(N)$, and so,
applying~$\log_2$ to the equation above, one finds precisely that
$$\Phi  N\Phi^{-1} = \log_2(U^q) = \log_2(\exp_2(N)^q).$$
Noting that\[N=\frac{1}{q}(qN + q^* N^3) - \frac{q^*}{q^4}(qN + q^*
  N^3)^3,  \]we see that the centralizers of~$N$ and~$qN+q^*N^3$
coincide. It follows that if ~$(\Phi,N)$ is a  point of~$\cN(q)$, then
 ~$(\Psi,N)$ is another point if and only if $\Psi \Phi^{-1}$
 centralizes~$N$ if and only if~$\Psi^{-1}\Phi$ centralizes~$N$.
Note also that if~$N^3=0$, then the condition on~$\Phi$ is
simply that~$\Phi  N\Phi^{-1} = qN$, while, if~$q = 1$, then the equation is simply that~$\Phi N \Phi^{-1} = N$.

Consider the projection
\begin{align*}
\cN(q)&\to \cN\\
(\Phi,N)&\mapsto N
\end{align*}
and let $\cN(q)_i$ denote the locally closed preimage of $\cN_i$.  We let $\cZ_i/\cN_i$ be the centralizer of the universal element over $\cN_i$.  Then there is an action of $\cZ_i$ on $\cN(q)_i$ by $z\cdot (\Phi,N)=(z\Phi,N)$.

\begin{prop}\label{prop: Z torsor}
The above action makes $\cN(q)_i$ into a $\cZ_i$-torsor over $\cN_i$.
\end{prop}
\begin{proof}
By Proposition \ref{prop: nilpotent conjugacy classes}, we may check
the proposition after base change to  a suitable \'{e}tale cover $U\to\cN_i$, over which the universal nilpotent over $U$ is of the form $gN_ig^{-1}$ for some $g\in\GSp_4(U)$.  
Let~$\Phi_i$ be  any of the explicit choices of~$\Phi$ given above for~$N_i$ (for~$i = 2$, take any
specialization of either~$\Phi_{2,a}$ or~$\Phi_{2,b}$). Then one readily checks that~$(\Phi_i,N_i)$ is a point on~$\cN(q)_i$, and that
\begin{align*}
(\cZ_i)_U&\to (\cN(q)_i)_U\\
z&\mapsto (zg\Phi_ig^{-1},gN_ig^{-1})
\end{align*}
is an isomorphism compatible with the $\cZ_i$-action.
\end{proof}

\begin{cor}\label{cor: component properties}
For $i=0,1,2,3$, $\cN(q)_i$ is smooth over $\cO$ with fibres equidimensional of dimension 11.  For $i\not=2$, $\cN(q)_i$ is irreducible with nonempty irreducible special fibre, while $\cN(q)_2$ has two connected components, each of which is irreducible with nonempty irreducible special fibre.
\end{cor}
\begin{proof}
The smoothness and dimension are an immediate consequence of Propositions \ref{prop: nilpotent conjugacy classes} and \ref{prop: Z torsor}.   Moreover, for $i\not=2$, $\cN(q)_i\to\cN_i$ is flat with irreducible fibres, and $\cN_i$ is irreducible, and hence $\cN(q)_i$ is irreducible.  The same argument applies to the special fibre.

Now we explain why $\cN(q)_2$ has two connected components.  As we
explained in the proof of Proposition \ref{prop: nilpotent conjugacy
  classes}, over $\cN_2$ we have the rank 2 non-degenerate quadratic
space $\cO_{\cN_2}^4/\ker(N)$ with quadratic form given by
$v\mapsto\psi(v,Nv)$. Over~$\cN_2$ we have~$N^3=0$, so  the
relation $\Phi N=qN\Phi$ holds on~$\cN(q)_2$, which implies that $\Phi$ preserves $\ker(N)$ and the computation
\begin{equation*}
\psi(\Phi v,N\Phi v)=q^{-1}\psi(\Phi v,\Phi Nv)=q^{-1}\nu(\Phi)\psi(v,Nv)
\end{equation*}
shows that $\Phi$ is an element of the general orthogonal group of
this quadratic space. 

This general orthogonal group has two
components (corresponding to whether the determinant and multiplier
agree or differ by a sign).  As a result we may write
$\cN(q)_2=\cN(q)_{2,a}\coprod\cN(q)_{2,b}$ where $\cN(q)_{2,a}$ is the
locus where $\Phi$ lies in the identity component and $\cN(q)_{2,b}$
is the locus where $\Phi$ lies in the nonidentity component.  Each of
these loci is
in fact nonempty; for example, we can consider points of the form $(\Phi_{2,a}(\alpha,\gamma),N_2)$ and $(\Phi_{2,b}(\alpha,\beta),N_2)$.  As $\cN(q)_{2,a}$ and $\cN(q)_{2,b}$ are unions of connected components, the action of $\cZ_2$ restricts to an action of the identity component $\cZ_2^\circ$ on each of them, and one easily checks that they must each be torsors for $\cZ_2^\circ$, and so the same argument as above shows that $\cN(q)_{2,a}$ and $\cN(q)_{2,b}$ are irreducible with nonempty irreducible special fibre.
\end{proof}

For the rest of this section, we will continue to use the notation
$\cN(q)_{2,a}$ and $\cN(q)_{2,b}$ for the two connected components of
$\cN(q)_2$ as introduced in the proof of Corollary~\ref{cor: component
  properties}.  We also write $\overline{\cN(q)}_{i}$ for the Zariski
closure of $\cN(q)_i$, $\overline{(\cN(q)_{i,\F})}$ for the Zariski
closure of its special fibre, and so on.

\begin{prop}\label{prop: Nq components}
The irreducible components of $\cN(q)$ are $\overline{\cN(q)}_{2,a}$, $\overline{\cN(q)}_{2,b}$,  and $\overline{\cN(q)}_i$ for $i=0,1,3$.  The irreducible components of the special fibre $\cN(q)_\F$ are $\overline{(\cN(q)_{2,a,\F})}$, $\overline{(\cN(q)_{2,b,\F})}$, and $\overline{(\cN(q)_{i,\F})}$ for $i=0,1,3$.  Each irreducible component of $\cN(q)$ has irreducible and generically reduced special fibre.
\end{prop}
\begin{proof}
$\cN(q)$ is set theoretically the disjoint union of the five locally closed subschemes $\cN(q)_{2,a}$, $\cN(q)_{2,b}$, and $\cN(q)_i$ for $i=0,1,3$, which are each irreducible and of the same dimension by Corollary~ \ref{cor: component properties}.  Hence their closures are the irreducible components of $\cN(q)$.  The same argument applies to the special fibre.

To prove the last statement it will suffice to prove that for
$i=0,1,2,3$, $\overline{\cN(q)}_i$ does not contain the generic points
of $\cN(q)_{j,\F}$ for $j\not=i$.  Indeed it already follows from
Corollary~ \ref{cor: component properties} that
$(\overline{\cN(q)_{2,a}})_{\F}$ does not contain the generic point of
$\cN(q)_{2,b,\F}$ and vice versa; and we also see that the special
fibre of each irreducible component of $\cN(q)$ is reduced at the
generic point of the corresponding component of $\cN(q)_\F$. 

In order to do this for $i=0,1,2,3$, let
$\widetilde{\cN(q)}_i\subset\cN(q)$ be the reduced closed subscheme
consisting of pairs $(\Phi,N)$ such that $\mathrm{rank}(N)\leq i$ and the
characteristic polynomial $\mathrm{char}_\Phi(X)$ is in $\cP_i$.  An
easy calculation shows that $\cN(q)_i\subset\widetilde{\cN(q)}_i$, and
hence $\overline{\cN(q)}_i\subset\widetilde{\cN(q)}_i$. (One can either follow the proof
of~\cite[Lem.\ 3.15]{jack}, or observe that we have seen above
that it is enough to check that this holds for the points of the
form~$(z_i\Phi_i,N_i)$ for our explicit choices of~$\Phi_i$ and for
$z_i\in \cZ_i$.) Thus to conclude the proof,
 all we have to do is exhibit a point on each irreducible component of $\cN(q)_\F$ which is only contained in one of the $\widetilde{\cN(q)}_i$'s.  For instance, we may take the following five points:
\begin{itemize}
\item $(\Phi_0(\alpha,\beta,\gamma),0)$ for general values of $\alpha,\beta,\gamma\in\overline{\F}^\times$.
\item $(\Phi_1(\alpha,\beta),N_1)$ for general values of
  $\alpha,\beta\in\overline{\F}^\times$.
\item $(\Phi_{2,a}(\alpha,\gamma),N_2)$ for general values of $\alpha,\gamma\in\overline{\F}^\times$.
\item $(\Phi_{2,b}(\alpha,\beta),N_2)$ for general values of $\alpha,\beta\in\overline{\F}^\times$.
\item $(\Phi_3(1),N_3)$.\qedhere
\end{itemize}
\end{proof}

For $x,y\in\cO^\times$ and $q$ a positive integer which is not a multiple of $p$, we let $\cM(x,y;q)$ be the closed subscheme of $\GSp_4^2/\cO$ consisting of pairs $(\Phi,\Sigma)$ satisfying:
\begin{itemize}
\item The characteristic polynomial of $\Sigma$ is $(X-x)(X-y)(X-y^{-1})(X-x^{-1})$.
\item $\Phi\Sigma\Phi^{-1}=\Sigma^q$.
\end{itemize}
We note that the order of $x$ and $y$ doesn't matter.

There is evidently  an isomorphism
\begin{align*}
\cM(1,1;q)&\to\cN(q)\\
(\Phi,\Sigma)&\mapsto(\Phi,\log_2(\Sigma)).
\end{align*}

We now have the following analogue of~\cite[Lem.\ 3.2]{tay}.
\begin{prop}\label{prop: key properties of components}
Let $q$ be a positive integer with $q\equiv 1\pmod p$.
\begin{enumerate}
\item Let $\cM_i$ be the irreducible components of $\cM(1,1;q)$ with
  their reduced subscheme structure.  Then the special fibres
  $\cM_{i,\F}$ are distinct, generically reduced and irreducible, and their reductions are precisely the irreducible components of $\cM(1,1;q)_\F$.
\item Suppose that either $q\not=1$ and $x,y$ are non trivial
  $(q-1)$st roots of $1$ in $1+\lambda\cO$ with $x\not=y^{\pm1}$; or that $q=1$ and $x,y$ are arbitrary elements of $1+\lambda\cO$.  Then $\cM(x,y;q)^{\mathrm{red}}$ is flat over $\cO$.
\end{enumerate}
\end{prop}
\begin{proof}  \leavevmode
\begin{enumerate}
\item This is an immediate consequence of Proposition \ref{prop: Nq components} and the isomorphism $\cM(1,1;q)\simeq\cN(q)$ above.
\item When $q\not=1$, we observe that, as $x,y,y^{-1},x^{-1}$ are distinct~$(q-1)$st roots of unity,
\begin{equation*}
\mathrm{char}_\Sigma(X)=(X-x)(X-y)(X-x^{-1})(X-y^{-1})|(X^{q-1}-1).
\end{equation*}
Hence, by the Cayley--Hamilton theorem, $\Sigma^q=\Sigma$. This implies that there is an isomorphism~$\cM(x,y;q)=\cM(x,y;1)$.
We are therefore reduced to the case that $q=1$.  

To show that
  $\cM(x,y;1)^{\mathrm{red}}$ is flat over $\cO$, it suffices
  to show that each generic point of its special fibre is the
  specialization of a point of the generic fibre. It suffices in turn
  to show that a Zariski dense set of points of the special fibre lift
  to the generic fibre. Then as $x$ and $y$ reduce to 1, we have
  $\cM(x,y;1)_\F=\cM(1,1;1)_\F\simeq\cN(1)_\F$.  This isomorphism,
  combined with the proof of Proposition~\ref{prop: Nq components},
  shows that the following five kinds of $\overline{\F}$-points are
  Zariski dense in $\cM(x,y;1)_\F$ (because the corresponding points
  are dense in each~$\widetilde{\cN(q)}_i$): \begin{itemize}\item $(g\Phi_0(\alpha,\beta,\gamma)g^{-1},1)$
\item $(g\Phi_1(\alpha,\beta)g^{-1},g\exp_2(N_1)g^{-1})$
\item $(g\Phi_{2,a}(\alpha,\gamma)g^{-1},g\exp_2(N_2)g^{-1})$
\item $(g\Phi_{2,b}(\alpha,\beta)g^{-1},g\exp_2(N_2)g^{-1})$
\item $(g\Phi_3(\alpha)g^{-1},g\exp_2(N_3)g^{-1})$
\end{itemize}where $\alpha,\beta,\gamma\in\overline{\F}^\times$ and $g\in\GSp_4(\overline{\F})$.  Then letting $\tilde{\alpha},\tilde{\beta},\tilde{\gamma}\in W(\overline{\F})$ and $\tilde{g}\in\GSp_4(W(\overline{\F}))$ be lifts, we can lift these to $W(\overline{\F})$ points of $\cM(x,y;1)$ of 
the following form (recall that we are in the case~$q = 1$):
\begin{itemize}
\item $(\tilde{g}\Phi_0(\tilde{\alpha},\tilde{\beta},\tilde{\gamma})\tilde{g}^{-1},\tilde{g}\mathrm{diag}(x,y,y^{-1},x^{-1})\tilde{g}^{-1})$
\item $(\tilde{g}\Phi_1(\tilde{\alpha},\tilde{\beta})\tilde{g}^{-1},\tilde{g}\mathrm{diag}(x,y,y^{-1},x^{-1})\exp_2(N_1)\tilde{g}^{-1})$
\item $(\tilde{g}\Phi_{2,a}(\tilde{\alpha},\tilde{\gamma})\tilde{g}^{-1},\tilde{g}\mathrm{diag}(x,y,y^{-1},x^{-1})\exp_2(N_2)\tilde{g}^{-1})$
\item $(\tilde{g}\Phi_{2,b}(\tilde{\alpha},\tilde{\beta})\tilde{g}^{-1},\tilde{g}\mathrm{diag}(A,X(A^\tau)^{-1}X)\exp_2(N_2)\tilde{g}^{-1})$, where~$X=
  \begin{pmatrix}
    0&1\\1&0
  \end{pmatrix},
$ and $A$ is a 2 by 2 matrix with coefficients in $W(\overline{\F})$ which has
trivial reduction, commutes with
$\begin{pmatrix}0&\tilde{\alpha}\\\tilde{\beta}&0\end{pmatrix}$ and
has eigenvalues $x,y$ (for the existence of such a matrix, use that
$\begin{pmatrix}0&\tilde{\alpha}\\\tilde{\beta}&0\end{pmatrix}$ has
distinct eigenvalues mod $p$, and is therefore diagonalizable).\item $(\tilde{g}\Phi_3(\tilde{\alpha})\tilde{g}^{-1},\tilde{g}\mathrm{diag}(x,y,y^{-1},x^{-1})\exp_2(N_3)\tilde{g}^{-1})$.\qedhere
\end{itemize}
\end{enumerate}
\end{proof}

Next we have an analogue of~\cite[Lem.\ 3.4]{tay}.
\begin{prop}\label{prop:Rchi connected}
Let $q>1$ with $q\equiv 1\pmod p$ and let $x,y$ be non trivial $(q-1)$st roots of $1$ in $1+\lambda\cO$ with $x\not=y^{\pm1}$.  Let $R=\hat{\cO}_{\cM(x,y;q),(1,1)}$ be the complete local ring of $\cM(x,y;q)$ at the point $(1,1)$ of the special fibre.  Then $\Spec R[1/p]$ is connected.
\end{prop}
\begin{proof}
The proof of~\cite[Lem.\ 3.4]{tay} carries over with minor modifications. Let~$\PPP_0$ denote the maximal ideal of~$R[1/p]$ corresponding to~$(\Phi_0,\Sigma_0)$ with~$\Phi_0$ trivial and~$\Sigma_0$
the diagonal matrix~$\diag(x,y,x^{-1},y^{-1})$, and let~$\PPP$ be
another maximal ideal, corresponding to a pair~$(\Phi,\Sigma)$. 
We need to show that~$\PPP$ is in the same connected component as~$\PPP_0$. 
One deduces as in~\cite{tay} that~$\PPP$ is in the same connected component of~$\Spec(R[1/p])$ as the maximal ideal corresponding
to~$(E^{-1} \Phi E,E^{-1} \Sigma  E)$ where ~$E \in \GSp_4(\OL)$ is arbitrary. In order to pass to an upper triangular form, we require the existence
of a filtration~$\Fil^i$ of~$k(\PPP)^4$
such that:
\begin{enumerate}
\item Each~$\Fil^i$ is preserved by~$\Phi$ and~$\Sigma$.
\item The graded pieces~$\gr^i$ are one dimensional and their eigenvalues (in order) are
$\alpha$, $\beta$, $\gamma \beta^{-1}$, $\gamma \alpha^{-1}$, which are the generalized eigenvalues of~$\Phi$. 
\item The orthogonal complement of~$\Fil^i$ is~$\Fil^{4-i}$.
\end{enumerate}
As in the proof of Proposition~\ref{prop: key properties of
  components}, $\Phi$ and~$\Sigma$ commute, so we may choose~$\Fil^1$
to be a common eigenvector of $\Phi$ and~$\Sigma$. We define $\Fil^3$ to be
the orthogonal complement of~$\Fil^1$, and then
choose~$\Fil^2$ to be any lift of a common eigenvector of~$\Phi$
and~$\Sigma$ in~$\Fil^3/\Fil^1$.

The constructions of paths
in~\cite{tay} from upper triangular to diagonal and between diagonal
matrices (eventually to~$(\Phi_0,\Sigma_0)$ and thus connecting~$\PPP$ to~$\PPP_0$) have obvious symplectic modifications.
\end{proof}

\subsubsection{Application to deformation rings}
Now let $\chi=(\chi_1,\chi_2)$ be a pair of continuous characters $\chi_i:\cO_{F_v}^\times\to\cO^\times$ that are trivial mod $\lambda$, let $\tilde{\cD}_v^\chi$ be the functor on $\CNL_\cO$ of continuous homomorphisms $\rho:G_{F_v}\to\GSp_4(A)$ which are trivial mod $\m_A$ and such that for $\sigma\in I_{F_v}$, the characteristic polynomial of $\rho(\sigma)$is
\begin{equation*}
(X-\chi_1(\Art_{F_v}^{-1}(\sigma)))(X-\chi_2(\Art_{F_v}^{-1}(\sigma)))(X-\chi_2(\Art_{F_v}^{-1}(\sigma))^{-1})(X-\chi_1(\Art_{F_v}^{-1}(\sigma))^{-1}).
\end{equation*} As in~\S\ref{subsec: l not p deformations}, we let $\cD_v^\chi\subset\tilde{\cD}_v^\chi$ be the subfunctor of $\rho$ with $\nu\circ\rho=\varepsilon^{-1}$.  The functors $\tilde{\cD}_v^\chi$ and $\cD_v^\chi$ are representable by rings $\tilde{R}_v^\chi$ and $R_v^\chi$.  We also let $\cD_1$ be the functor with $\cD_1(A)$ parameterizing continuous unramified characters $\psi:G_{F_v}\to A^\times$ which are trivial mod $\m_A$.  It is representable by $\cO[[T]]$ with universal object $\chi^{\mathrm{univ}}:G_{F_v}\to\cO[[T]]^\times$ given by $\chi^{\mathrm{univ}}(\Frob_v)=1+T$.

For any $A\in\CNL_\cO$, then as $A$ is complete and $p>2$,
\begin{align*}
1+\m_A\to 1+\m_A\\
t\mapsto t^2
\end{align*}
is a bijection and we denote its inverse by $x\mapsto\sqrt{x}$.  Then we have

\begin{prop}\label{prop: twisting to fix central char}
There is an isomorphism of functors
\begin{align*}
\cD_v^\chi\times\cD_1&\to \tilde{\cD}_v^\chi\\
(\rho,\psi)&\mapsto\rho\otimes\psi
\end{align*}
Consequently there is an isomorphism $R_v^\chi[[T]]\simeq \tilde{R}_v^\chi$.
\end{prop}
\begin{proof}
For the inverse we may take the natural transformation 
\begin{align*}
\tilde{\cD}_v^\chi&\to\cD_v^\chi\times\cD_1\\
\rho&\mapsto(\rho\otimes \sqrt{\varepsilon\cdot (\nu\circ\rho)}^{-1},\sqrt{\varepsilon\cdot(\nu\circ\rho)})
\end{align*}
The only thing that we need to check is that if $\rho\in\tilde{\cD}_v^\chi(A)$ then $\nu\circ\rho$ is trivial on $I_{F_v}$.  For $\sigma\in I_{F_v}$, $(\nu\circ\rho(\sigma))^2$ is the constant term of the characteristic polynomial of $\rho(\sigma)$ which is 1 by definition.  But also $\nu\circ\rho(\sigma)\equiv 1\pmod {\m_A}$, and hence $\nu\circ\rho(\sigma)=1$ as $p>2$.
\end{proof}
 
We may now relate these deformation rings to the spaces of matrices considered in this section.  
\begin{prop}\label{prop: relating deformation ring to M}
Let $\sigma$ be a chosen topological generator of the tame inertia subgroup of $G_{F_v}$.  Let $x=\chi_1(\Art_{F_v}^{-1}(\sigma))$ and $y=\chi_2(\Art_{F_v}^{-1}(\sigma))$.  Then
\begin{equation*}
\tilde{R}_v^\chi\simeq\hat{\cO}_{\cM(x,y;q_v),(1,1)}.
\end{equation*}
\end{prop}\begin{proof}Since~$\rhobar|_{G_{F_v}}$ is trivial, any lifting of it
  factors through the  quotient~$T_v=G_{F_v}/P_{F_v}$, where~$P_{F_v}$
  denotes the maximal pro-prime-to-$p$ subgroup of~$I_{F_v}$ (that is,
  the kernel of any non-trivial homomorphism $I_{F_v}\to\Zp)$. If
  ~$\varphi$ is an arithmetic Frobenius element in~$G_{F_v}$,
  then the group~$T_v$ is topologically generated by~$\varphi$
  and the image of~$\sigma$, subject to the constraints that~$\sigma$
  generates a pro-$p$ group, and that
  $\varphi\sigma\varphi^{-1}=\sigma^{q_v}$. The result then follows
  from the definitions.  
\end{proof}
We can now conclude the proofs of Propositions \ref{prop:Ihara1ring} and \ref{prop:Iharachiring} exactly as in~\cite{tay}.
\begin{proof}[Proof of Proposition \ref{prop:Ihara1ring}]
Combining Propositions~ \ref{prop: key properties of components}~(1) and~\ref{prop: relating deformation ring to M}
with~\cite[Lem.\ 2.7]{tay} proves the corresponding result for $\tilde{R}_v^1$.  The result for $R_v^1$ follows from this and Proposition \ref{prop: twisting to fix central char}.
\end{proof}
 
\begin{proof}[Proof of Proposition \ref{prop:Iharachiring}]
Proposition \ref{prop: key properties of components} implies that
$(\tilde{R}_v^\chi)^{\mathrm{red}}$ is flat over $\cO$.  Proposition
\ref{prop:Rchi connected} implies that $\Spec(\tilde{R}_v^\chi[1/p])$
is connected.  On the other hand, by Lemma~\ref{lem: chi Ihara
  avoidance deformation ring is smooth}, for any closed
point $x\in\Spec(\tilde{R}_v^\chi[1/p])$, the localization
$(\tilde{R}_v^\chi[1/p])_x$ is regular and hence a domain.  Then the
result follows from Propositions~\ref{prop: key properties of
  components}~(2) and~\ref{prop: relating deformation ring to M}, as in the proof of~\cite[Prop.\ 3.1]{tay}.\end{proof}

\subsection{Big  image conditions and vast
  representations}\label{subsec: big image conditions}
\subsubsection{Enormous subgroups}Following~\cite{CG,KT} (which give the analogous definition
for~$\GL_n$) we now define the notion of ``enormous
image,'' with some minor modifications. \begin{defn}\label{defn:enormous image}
We say that a subgroup $H\subset\GSp_4(k)$ is \emph{enormous} if it
satisfies the following conditions:
\begin{enumerate}[label={(E\arabic*)}]
\item   $H^1(H, \ad^0) = 0$ for the~$10$-dimensional representation~$\ad^0$.\item $H$ acts absolutely irreducibly in its natural representation,
in particular,~$H^0(H, \ad^0) = 0$.
 \item For all simple $\overline{k}[H]$-submodules $W \subset \overline{k} \otimes \ad^0$, there is an
  element $h \in H$ such that
  \begin{itemize}
  \item $h\in\GSp_4(k)$ has $4$ distinct eigenvalues, and \item $1$ is an eigenvalue for the action of~$h$ on~$W$.
  \end{itemize}
  \end{enumerate}If~$H$ only satisfies~(E2) and~(E3), then we say that~$H$ is \emph{weakly enormous}.
\end{defn}

\begin{lem} \label{lem:enormoustwisting}
If~$H$ and~$H'$ are subgroups of~$\GSp_4(k)$ with the same image
in~$\PGSp_4(k)$, then~$H$ is enormous \emph{(}resp.\  weakly enormous\emph{)} if and
only if $H'$ is enormous \emph{(}resp.\  weakly enormous\emph{)}.
\end{lem}

\begin{proof} Suppose that~$P$ is the projective image of~$H$ in~$\PGSp_4(k)$
and~$Z$ is the kernel. Then the action of~$H$ on~$\ad^0$ factors through~$P$.
In particular, ~$H^0(P,\ad^0) = H^0(H,\ad^0)$, and there is an inflation--restriction
sequence
$$0 \rightarrow H^1(P,\ad^0) \rightarrow H^1(H,\ad^0) \rightarrow H^1(Z,\ad^0)^P = 0.$$
Hence all the conditions in the definitions of enormousness and weakly enormousness
depend only on the projective representation.
\end{proof}

Note that if~$H' \subset H$ is weakly enormous, then so is~$H$, but if~$H'$ is enormous, then~$H$
is not necessarily enormous.

\begin{rem}  Some ``big image'' conditions in the literature have the additional assumption that~$H$ has no~$p$-power quotient.
In practice, however, that hypothesis is often only used in a very weak way, namely, to ensure that the image of~$\rhobar$
restricted to~$G_{F(\zeta_p)}$ coincides with the restriction
to~$G_{F(\zeta_{p^N})}$ for all~$N\ge 1$.
The stronger hypothesis has the unfortunate side effect of ruling out some perfectly fine Galois representations to
which the Taylor--Wiles method applies, most notably, surjective
representations~$\rhobar:G_{\Q} \rightarrow \GL_2(\F_3)$ with cyclotomic determinant (exactly
the case which arises in the original work of Wiles!).
In order not to rule out some interesting subgroups which occur for~$p = 3$, we therefore
do not assume this hypothesis.
\end{rem}

Let~$p \ge 3$. The cyclotomic character induces a homomorphism:
$$G_F \rightarrow \Z^{\times}_p \simeq (\Z/p\Z)^{\times} \oplus (1 + p \Z_p) \rightarrow (1 + p \Z_p).$$
If~$p$ is unramified in~$F$, then this composite map is surjective. In general, the image contains~$1 + p^{\delta}$
for some integer~$\delta$.
  In order to address the passage from~$F(\zeta_p)$ to~$F(\zeta_{p^N})$ in the Taylor--Wiles
  argument, we have the following lemma:

\begin{lemma} \label{lemma:imageindependence} Suppose that~$p \ge 3$. Let
$$\rhobar:G_{F} \rightarrow \GSp_4(k)$$
be a continuous homomorphism. Then there exists an integer~$\delta$ depending only on~$F$ such that
 the image of~$\rhobar$ restricted to~$G_{F(\zeta_{p^N})}$
is independent of~$N$ for~$N \ge 1 + \delta$ if~$p \ge 5$ or~$N \ge 2 + \delta$ for~$p = 3$.
If~$p$ is unramified in~$F$, then one may take~$\delta = 1$.
\end{lemma}

\begin{proof}  There is a canonical injective
homomorphism
$$\Gal(F(\zeta_{p^N})/F)  \rightarrow (\Z/p^N \Z)^{\times}$$ for
all~$N$, and we will identify~$\Gal(F(\zeta_{p^N})/F)$ with its image
in~$(\Z/p^N \Z)^{\times}$ in the below. We choose~$\delta$ such that for all~$N$,
 the image contains~$1 + p^{\delta}$. In particular, if~$p$ is
 unramified in~$F$, we can take~$\delta=1$.

Let~$M$ denote the fixed field of~$\rhobar$. There are natural maps as follows:
$$
\begin{diagram}
\Gal(M(\zeta_{p^{N+1}})/F)  & \rInto &  \Gal(M/F) \times \Gal(F(\zeta_{p^{N+1}})/F) \\
\dTo & &\dTo \\
\Gal(M(\zeta_{p^{N}})/F)  & \rInto &  \Gal(M/F) \times  \Gal(F(\zeta_{p^N})/F)
\end{diagram}
$$
where the composites of the horizontal maps with the projections to
each factor are surjective.
The images of~$\rhobar$ restricted to~$F(\zeta_{p^N})$ and~$F(\zeta_{p^{N+1}})$ coincide
precisely when the left hand  vertical map has non-trivial kernel (necessarily of order~$p$).
We prove this is so under our assumptions on~$N$.

It suffices to show that the horizontal image of the upper map contains an
element of the form~$(\id_M,1+ m p^N)$ for some~~$m$ with~$(m,p)=1$. 
By the surjectivity onto the second factor, it contains an element of the form~$(g,1+p^{\delta})$.
Let~$m$ be the prime to~$p$ order of~$g$, so that~$h:=g^m$ has~$p$-power order.
Since~$p > 2$, we have~$(g,1+p^\delta)^{m p^{N-\delta}} = (h^{p^{N-\delta}}, 1 + m p^N)$, and hence we are done providing the order of~$h$
divides~$p^{N-\delta}$. Yet all~$p$-power elements of~$\GSp_4(k)$ have order dividing~$p$ if~$p \ge 5$ or order dividing~$p^2$
if~$p =3$.  (The~$p$-Sylow subgroup of~$\GSp_4(k)$ consists of unipotent matrices which 
satisfy~$(\sigma - 1)^4 = 0$, so~$\sigma^{p^k} = 1$ when~$p^k \ge 4$.)
\end{proof}

 In anticipation of Lemma~\ref{lemma:vanishingofcohomology} below,
we make the following definition:
 
 \begin{defn}  \label{defn:vast}
 A representation~$\rhobar: G_{F} \rightarrow \GSp_4(k)$
 is \emph{vast} if one of the following two conditions holds:
 \begin{enumerate}
 \item  The image of~$\rhobar$ restricted to~$G_{F(\zeta_{p^{N}})}$ is
   enormous for all sufficiently large~$N$.
 \item  The image of~$\rhobar$ restricted to~$G_{F(\zeta_{p^{N}})}$ is
   weakly enormous for all sufficiently large ~$N$,
 and the fixed field~$L$ of~$\ad^0 \rhobar$ does not contain~$\zeta_p$.
 \end{enumerate}
 \end{defn}

 \begin{rem}\label{rem: vast unramified N 3} If~$p$ is unramified in~$F$, then, 
 in Definition~\ref{defn:vast},
 one may replace 
 \emph{sufficiently large ~$N$}
 by~$N = 3$, since, by
 Lemma~\ref{lemma:imageindependence}, the image in this case does not depend on~$N$ for~$N \ge 3$.
 \end{rem}
 \begin{rem}
   \label{rem: twisting vast}By Lemma~\ref{lem:enormoustwisting},
   ~$\rhobar$ is vast if and only if any twist of~$\rhobar$ by a
   character is vast.
 \end{rem}

The following lemma will prove useful for constructing Taylor--Wiles primes:
\begin{lemma} \label{lemma:vanishingofcohomology}Suppose that $p\ge 3$.
Let~$\rhobar: G_{F} \rightarrow \GSp_4(k)$ be a continuous
representation. Fix an integer~$N\ge 1$. Suppose either that:
\begin{enumerate}
\item 
The fixed field~$L$  of~$\ad^0 \rhobar$ does not contain~$\zeta_p$, or
\item The restriction of~$\rhobar$ to~$G_{F(\zeta_{p^N})}$ has enormous image.

\end{enumerate}
Then
$$H^1(L(\zeta_{p^N})/F,\ad^0\rhobar(1)) = 0.$$
In particular, if~$\rhobar$ is vast, then the conclusion above holds
for all sufficiently large~$N$.
\end{lemma}

\begin{proof} We first consider the case when~$\zeta_p \notin L$.
By inflation--restriction, it suffices to prove
that the groups
$$H^1(L(\zeta_p)/F,\ad^0\rhobar(1)),  \quad
H^1(L(\zeta_{p^N})/L(\zeta_{p}),  \ad^0\rhobar(1))^{\Gal(L(\zeta_p)/F)}
$$
both vanish. The group~$\Gal(L(\zeta_p)/L) \subset \Gal(L(\zeta_p)/F)$ acts trivially (by conjugation)
 on both the group~$\Gal(L(\zeta_{p^N})/L(\zeta_{p}))$ and the module~$\ad^0$.
However, it acts  by non-trivial scalars on the twist~$\ad^0(1)$ since we are assuming~$\zeta_p \notin L$. Hence the
second group vanishes after taking invariants.  Applying inflation--restriction now to the first group, it suffices to prove that
the groups
$$H^1(L/F, (\ad^0\rhobar(1))^{\Gal(L(\zeta_p)/L)}),  \quad
H^1(L(\zeta_{p})/L,\ad^0\rhobar(1))^{\Gal(L/F)}
$$
both vanish. The second group vanishes because $p\nmid [L(\zeta_p):L]$.  The
first group vanishes because~$\ad^0 \rhobar$ is fixed by~$\Gal(L(\zeta_p)/L)$ and thus has no invariants
after being twisted by the mod-$p$ cyclotomic character (which by assumption
is a non-trivial character of~$\Gal(L(\zeta_p)/L)$).

Now we consider the second case. Let~$M$ denote the splitting field of~$\rhobar$, so that~$M/L$
is a (possibly trivial) cyclic extension of degree prime
to~$p$. Inflation--restriction shows that we have an injection
$$H^1(L(\zeta_{p^N})/F,\ad^0\rhobar(1)) \hookrightarrow H^1(M(\zeta_{p^N})/F,\ad^0\rhobar(1)),$$
so it suffices to show that the latter group vanishes. 
By inflation--restriction, it is enough to show that the cohomology groups $$H^1(F(\zeta_{p^N})/F,H^0(M(\zeta_{p^N}/F(\zeta_{p^N}),\ad^0 \rhobar(1))),$$
$$H^1(M(\zeta_{p^N})/F(\zeta_{p^N}),\ad^0 \rhobar(1))$$ both vanish. We are assuming that
 $$H = \Gal(M(\zeta_{p^N})/F(\zeta_{p^N}))$$ is enormous.  Thus to show that both groups
above vanish, it suffices to note that
$$H^0(H,\ad^0) = H^1(H,\ad^0)=0$$
because~$H$ is enormous.
\end{proof}

\begin{remark} The two parts of this proof are essentially standard
  --- in particular the first part is  exactly
the same as the proof
of Lemma~5.3 of~\cite{MR2783930}.
\end{remark}

We will require a weakly enormous (or in practice vast) image assumption in order to use the
Cebotarev density theorem to guarantee the existence of Taylor--Wiles
primes. Similarly, the following condition will allow us to use
Cebotarev to arrange for our level structures to be neat by increasing
the level at an auxiliary prime.
\begin{defn} \label{defn:tidy}We say that a subgroup $H\subset\GSp_4(k)$ is \emph{tidy} if there
  is an~$h\in H$ with~$\nu(h) \ne 1$, and such that no two eigenvalues of~$h$ have ratio~$\nu(h)$
  (but the eigenvalues need not be distinct). We say that a representation~$\rhobar: G_{F} \rightarrow \GSp_4(k)$
  is tidy if it has tidy image. 
\end{defn}

Note that the property of tidiness is inherited from subgroups.

\begin{lem}\label{lem:bigcentre}
Suppose that~$H \subset \GSp_4(k)$  is absolutely irreducible, and the centre~$Z$ of~$H$ has order at least~$3$.
Then~$H$ is tidy.
\end{lem}

\begin{proof} By Schur's lemma, the centre is cyclic and any element
in the centre is scalar with eigenvalues~$(\zeta,\zeta,\zeta,\zeta)$
for some~$\zeta$. If~$|Z| \ge 3$, there thus exists
such an element~$h$ in the centre with~$\zeta^2 \ne 1$. Since~$\nu(h)  = \zeta^2 \ne 1$, and since the ratio
of every pair of eigenvalues is~$1 \ne \nu(h)$, it follows that~$H$ is tidy.
\end{proof}

\begin{lem}
  \label{lem: criterion for tidiness}Let~$\Delta \subset \GL_2(\F_p) \times \GL_2(\F_p)$ be the subgroup of pairs~$(A,B)$ with~$\det(A) = \det(B)$,
  and consider~$\Delta$ as a subgroup of~$\GSp_4(\F_p)$ via the map of
 ~\S\ref{section:inductions}.
   If~$p \ge 5$
and~$\Delta \subset H$,
then~$H$ is tidy.  \end{lem}
\begin{proof}  The argument is very similar to the proof Lemma~\ref{lem:bigcentre}. The group~$\Delta$
contains a cyclic subgroup of scalar matrices of order~$p - 1 > 2$.
\end{proof}
\begin{lem}
  \label{lem: enormous conditions automatic}If $p\ge 11$ and~$H \subset \GSp_4(k)$ is absolutely irreducible,
   then conditions~\emph{(E1)} and~\emph{(E2)} are satisfied.
  \end{lem}
\begin{proof}
  This is immediate from~\cite[Thm.\ A.9]{jack}.
\end{proof}

\begin{lem}
  \label{lem: GSp4 enormous}If~$p\ge 3$, then~$H=\Sp_4(\F_p)$ is enormous and~$G = \GSp_4(\F_p)$
  is tidy. If~$\rhobar: G_F \rightarrow \GSp_4(\F_p)$ is a surjective representation with similitude character~$\varepsbar^{-1}$,
  then~$\rhobar$ is vast and tidy.
  \end{lem}
  
\begin{proof}  For all such~$p$, the representation~$\ad^0$ is
  absolutely irreducible. Hence for weak enormity it suffices to note that~$H$ contains elements with distinct eigenvalues,
and every such element has at least one eigenvalue~$1$
on~$\ad^0$. Thus for enormity it suffices to
check that~$H^1(\Sp_4(\F_p),\ad^0) = 0$. For~$p \ge 11$, 
this follows from Lemma~\ref{lem: enormous conditions automatic}. For~$p = 3$, $5$, and~$7$,
it can be checked directly using~\texttt{magma}~\cite{MR1484478}. 
(All of the \texttt{magma} code and output for this paper can be found at
the github respository here~\cite{magma}.)

 For tidiness, the centre of~$G$ has order~$p-1$
so the result follows from Lemma~\ref{lem:bigcentre} when~$p > 3$.  (It
also follows from Lemma~\ref{lem: criterion for tidiness}.) When~$p = 3$, the group~$\GSp_4(\F_3)$
contains an element~$g$ of order~$20$  with~$\nu(g) = -1$; more
precisely, its eigenvalues are of the
form~$\zeta,\zeta^3,\zeta^9,\zeta^{27}$ for a~$20$th root of
unity~$\zeta$, and $\nu(g)=\zeta^{10}=\zeta^{30}=-1$. The ratios of
the pairs of eigenvalues
are of the form~$\zeta^{3-1}$, $\zeta^{9-1}$, and~$\zeta^{27-1}$, and
since none of these quantities is equal to~$\zeta^{10} = -1$, we are done.

For vastness, note that the image of~$\rhobar$ restricted to~$G_{F(\zeta_{p})}$ will be~$H = \Sp_4(\F_p)$.
Since this group has no quotients of $p$-power order (indeed~$\PSp_4(\F_p)$ is
simple), the image of the restriction of~$\rhobar$ to~$G_{F(\zeta_{p^N})}$ will also be~$H$
for all~$N$. Hence the image of~$\rhobar$  restricted to~$G_{F(\zeta_{p^N})}$ is always~$H$ and hence enormous;
thus~$\rhobar$ is vast.
\end{proof}

\subsubsection{Representations Induced From Index Two Subgroups}

Suppose that~$G \subset \Sp_4(k)$ is an absolutely irreducible subgroup such that the  underlying
representation~$W$ becomes reducible on an index two subgroup~$H$.
Write~$\chi$ for the quadratic character~$\chi: G \rightarrow G/H \rightarrow k^{\times}$ (we assume
the characteristic of~$k$ is different from~$2$). Write $G/H=\{1,\sigma\}$. Then one may write~$W|_{H} = V \oplus V^{\sigma}$,
and one has the following~$G$-equivariant decompositions (not necessarily into irreducibles):
$$W \otimes W = W \otimes W^{\vee} = k \oplus k(\chi) \oplus \Ind^{G}_{H} (\ad^0(V)) \oplus \Asai(V) \oplus \Asai(V) \otimes \chi,$$
$$\ad^0(W) = \Sym^2(W) =  \Ind^{G}_{H} (\ad^0(V)) \oplus \Asai(V),$$
$$\wedge^2(W) = k \oplus k(\chi) \oplus  \Asai(V) \otimes \chi.$$
Here~$\Asai(V)$ is the Asai representation, which satisfies~$\Asai(V)|_{H} = V \otimes V^{\sigma}$. 
These identifications follow from computing what happens over~$H$ and noting that~$W \simeq W \otimes \chi$.

\begin{lemma} \label{lemma:inductionsbigprep} Suppose that~$\Asai(V)$ and~$\Ind^{G}_{H}(\ad^0(V))$ are absolutely irreducible representations
of~$G$. Suppose that~$G \setminus H$ has an element~$g$ of order neither dividing~$4$ nor divisible by~$p$. Then~$G$ satisfies condition~\emph{(E3)} of enormousness.
\end{lemma}

\begin{proof} Let~$g$ be an element of~$G \setminus H$. Since~$W$ is induced, the eigenvalues of~$g$
are invariant under multiplication by~$-1$. Since~$G \subset \Sp_4(k)$, the eigenvalues are invariant
under inversion. It follows that the eigenvalues are of the form~$(\alpha, \alpha^{-1}, - \alpha, - \alpha^{-1})$
for some~$\alpha$.  If~$g$ has order neither dividing~$4$ nor divisible by~$p$, then~$\alpha^4 \ne 1$ and these eigenvalues
are all distinct.
To show~(E3), it is enough to show that any such element~$g$ has an eigenvalue~$1$ on both~$\Asai(V)$
and~$\Ind^{G}_{H}(\ad^0(V))$. Let~$\Gamma = \langle g \rangle$, and work in the Grothendieck
group of representations of~$\Gamma$. The representation~$\Sym^2$ differs from~$\wedge^2$ by containing
the squares of all the eigenvalues. Hence
$$[\Sym^2] = [\wedge^2] + [\alpha^2,\alpha^{-2},\alpha^{2},\alpha^{-2}].$$
Moreover, since~$\chi(g) = -1$,
$$[\wedge^2] = [1] + [-1] + [-\Asai(V)].$$
It follows by counting eigenvalues in~$W \otimes W$ that
$$
[\wedge^2] = [1] + [-1] + [1,-1,-\alpha^2,-\alpha^{-2}],$$
$$[\Sym^2] = [-1,1,\alpha^2,\alpha^{-2}] + [1,-1,\alpha^2,-\alpha^2,\alpha^{-2},-\alpha^{-2}].$$
from which it follows that
$$[\Asai(V)] = [-1,1,\alpha^2,\alpha^{-2}],$$
$$\Ind^{G}_{H}(\ad^0(V)) = [1,-1,\alpha^2,-\alpha^2,\alpha^{-2},-\alpha^{-2}],$$
both of which have~$1$ as an eigenvalue.
\end{proof}

\begin{lemma} \label{lemma:eq5} \label{lemma:biginduced} Assume~$k$ has characteristic~$p \ge 3$.
Let~$G$ be the group~$\SL_2(k) \wr \Z/2\Z = (\SL_2(k) \times \SL_2(k))
\rtimes \Z/2\Z$, where the semi-direct product swaps the two copies of~$\SL_2(k)$,
considered as a subgroup of~$\Sp_4(k)$ as in~\S\ref{section:inductions}.
Then~$G$ is weakly enormous, and is furthermore enormous if~$\#k \ne 5$.
\end{lemma}

\begin{proof}We begin by checking that property~(E1) holds.
Let~$H = \SL_2(k) \times \SL_2(k)=A\times B$, say, and  
let~$V_A$ and~$V_B$ denote the tautological $2$-dimensional
representations of~$A$ and~$B$, so that~$W|_{H} = V_A \oplus V_B$, 
and~$\Sym^2(W)|_{H} = \ad^0(W)|_{H} = \ad^0(V_A) \oplus \ad^0(V_B) \oplus
V_A \otimes V_B$. 
Since~$H^1(G,\ad^0(W)) = H^1(H,\ad^0(W))^{G/H}$, it suffices to prove
that
$$H^1(H,\ad^0(W)) = H^1(A \times B,\ad^0(W)) = 0.$$
By inflation--restriction, we see that there are
exact sequences:
$$H^1(A,\ad^0(V_A)) \rightarrow H^1(A \times B,\ad^0(V_A))
\rightarrow (H^1(B,k) \otimes \ad^0(V_A))^A=0,$$
$$0 = H^1(A,(V_A \otimes V_B)^{B}) \rightarrow H^1(A \times B,V_A \otimes V_B)
\rightarrow (H^1(B,V_B) \otimes V_A)^{A} = 0.$$
Thus it remains to show that~$H^1(A,\ad^0(V_A)) = 0$. But this is the same as showing
that
$$H^1(\SL_2(k),\Sym^2(k^2)) = 0.$$
This holds for~$\#k \ne 5$ (which we are assuming) by~\cite[Lem.\ 2.48]{MR1605752}.

Property~(E2) is obvious. For property~(E3), it suffices by Lemma~\ref{lemma:inductionsbigprep}  to show that~$G \setminus H$
contains an element~$g$ of order not dividing~$4$ and not divisible by~$p$.
Since~$p^2 - 1$ is always divisible by~$8$, there exists a matrix~$a \in A$ of order exactly~$8$.
The automorphism~$\sigma: A \times B \rightarrow B \times A$ of order~$2$ identifies~$A$ with~$B$,
and with respect to this identification let~$g = \sigma (a,a) = (a,a) \sigma$. Then~$g^2 = (a^2,a^2)$
has order~$4$, so~$g$ has order~$8$ which does not divide~$4$ and is
not divisible by~$p$, as required.
\end{proof}

For~$p = 5$ one has the following substitute: 

\begin{lemma} \label{lemma:neq5} Let~$H/F$ be a quadratic extension, and let~$\rbar: G_{H} \rightarrow \GL_2(\F_5)$
be a surjective representation with  determinant~$\varepsbar^{-1}$. Let~$\rhobar: G_{F} \rightarrow \GSp_4(\F_5)$
be the induction of~$\rbar$ to~$F$, and assume that the image of~$\rhobar|_{G_{F(\zeta_5)}}$
is equal to~$G= \SL_2(\F_5) \wr \Z/2\Z$. Assume furthermore that~$5$ is unramified in~$F$.
Then~$G:=\rhobar(G_{F(\zeta_{5^N})})$ is weakly enormous for all~$N \ge 1$
and~$\zeta_5$ does not lie in the fixed field of~$\ad^0 \rhobar$; in particular, $\rhobar$ is vast.
\end{lemma}

\begin{proof} Since the abelianization of~$G$ has order prime to~$5$, the
image of~$\rhobar$ over~$F(\zeta_5)$ is the same as the image
over~$F(\zeta_{5^N})$ for any~$N$, and is weakly enormous by
Lemma~\ref{lemma:biginduced}.
Let~$\Gamma$ denote the image of~$\rhobar$. Since~$5$ is unramified in~$F$ and the similitude character of~$\Gamma$
is inverse cyclotomic, it follows that the similitude character is
surjective and~$[\Gamma:G] = 4$. In particular, the group~$\Gamma$ is the full pre-image of~$G$ in~$\GSp_4(\F_p)$,
and is generated by pairs~$(A,B)$ in~$\GL_2(\F_p)$ with~$\det(A) =
\det(B)$ together with an involution sending~$(A,B)$ to~$(B,A)$. 
The fixed field~$L$ of~$\ad^0 \rhobar$ is the fixed field of the projective representation.
But one can now observe directly that the image of~$\Gamma$ in~$\PGSp_4(\F_5)$ has 
abelianization~$(\Z/2\Z \oplus \Z/2\Z)$, which does not surject
onto~$\Gal(F(\zeta_5)/F) = \Z/4\Z$. So~$\zeta_p\notin L$,
and~$\rhobar$ is vast, as required.
  \end{proof}

\subsubsection{The enormous subgroups of~$\Sp_4(\F_3)$}

By an exhaustive search, one can determine precisely which of the subgroups of~$\Sp_4(\F_3)$
are enormous. There are~$162$ conjugacy classes of subgroups, and it turns out that
precisely~$11$ of them are enormous, of orders
$40$, $128$, $160$, $192$, $240$,  $320$, $384$, $384$, $1152$, $1920$,
 and $51840$ respectively. Our main interest will be in representations~$\rhobar$ to~$\GSp_4(\F_3)$
which are vast and tidy. In particular, it is of interest to consider subgroups~$G$ of~$\GSp_4(\F_3)$
which are tidy and such that~$H = \Sp_4(\F_3) \cap G$ is enormous. 
Sometimes the tautological~$4$-dimensional representation~$V$ of one of these groups~$G$ fails
to be absolutely irreducible on an index
 two subgroup --- necessarily this subgroup is 
not~$H = G \cap \Sp_4(\F_3)$ because we are assuming that~$H$ is enormous and hence acts absolutely irreducibly on~$V$.
The representation~$V$ underlying~$G$ restricted
to this index two subgroup either becomes reducible over~$\F_3$ or over a non-trivial extension
of~$\F_3$. In the former case, we say that~$G$ is split induced. In this case, the index two subgroup is necessarily
a subgroup of \[\Delta = \{(A,B) \subset \GL_2(\F_p) \times
  \GL_2(\F_p) \text{ where}~\det(A) = \det(B)\}\] and~$G$
is a subgroup of~$\Gamma:= \Delta \rtimes \Z/2\Z$ where~$\Z/2\Z$ swaps the factors. Hence we may write
the index two subgroup in this case as~$G \cap \Delta$.

We collect a number of interesting examples in the following lemma.

\begin{lem} \label{lem:sample} The following groups~$G \subset \GSp_4(\F_3)$ are tidy, and such that the index two subgroup~$H = G \cap \Sp_4(\F_3)$ is enormous. 
\end{lem}

\begin{enumerate}
\item The group~$G = \GSp_4(\F_3)$.
\item A group~$G$ of order~$3840$. The projective image has index~$27$ in~$\PGSp_4(\F_3)$. It may  be identified as the stabilizer of the natural action of~$\PGSp_4(\F_3)$ on the~$27$ lines
of a cubic surface.
\item \label{case:interesting} {\bf Split Inductions\rm}.  The following subgroups of~$\Gamma = \Delta \rtimes \Z/2\Z$:
\begin{enumerate}
\item \label{case:2304} The group~$G = \Gamma$  of order~$2304$.
\item The two groups~$G$ of  index~$3$ in~$\Gamma$. They are the two groups
of order 768  inside~$\GSp_4(\F_3)$ up to conjugacy, and they are distinguished by their
intersections~$H = G \cap \Sp_4(\F_3) \subset \SL_2(\F_3) \wr \Z/2\Z$ and~$H
\cap \Delta \subset \SL_2(\F_3)^2$. Note there is a homomorphism~$\chi: \SL_2(\F_3) \rightarrow A_4
\rightarrow \Z/3\Z$. One intersection~$H \cap \Delta$ is given by pairs~$(A,B)$ with~$\chi(A) = \chi(B)$, 
and the other by pairs with~$\chi(A) = - \chi(B)$. Note that these groups are abstractly
isomorphic (the outer automorphism of~$\SL_2(\F_3)$ sends~$\chi$ to~$-\chi$) but
not conjugate inside~$\GSp_4(\F_3)$.
\end{enumerate}
\item {\bf Other Inductions.\rm} \label{exoticinduction}
 A group~$G$ of order~$480$ with projective image~$S_5 \times \Z/2\Z$.
There is an isomorphism~$\PSL_2(\F_9) = A_6$, and hence a projective~$\F_9$ representation of
the subgroup~$A_5 \subset A_6$. This is not unique --- there are two natural conjugacy
classes of~$A_5$
permuted by the exotic automorphism of~$A_6$. But that automorphism is induced by~$\Frob_3$
acting on the field of coefficients~$\F_9$, so the choice does not matter. 
There is a corresponding lift:
$$\widetilde{A}_5 \rightarrow \GL_2(\F_9)$$
by a group~$\wA_5$ which is a central extension of~$A_5$ by~$\Z/4 \Z$.  The outer automorphism
group of~$\wA_5$ is~$(\Z/2\Z)^2$, and there is a unique such outer automorphism which acts by~$-1$ on the centre
and by an outer automorphism on~$A_5$. Moreover, this lifts to a genuine automorphism~$\sigma$  of~$\wA_5$ of order~$2$.
Then $G:= \widetilde{A}_5 \rtimes \langle \sigma \rangle \subset \GSp_4(\F_3)$ has order~$480$.
This is the only enormous subgroup which both has induced image and is not solvable.
{\bf Warning:\rm} The group~$G$ is not determined up to conjugacy by its order. Indeed,
there exists a second conjugacy class of subgroups~$G'$ of order~$480$ with~$H' = G' \cap \Sp_4(\F_3)$
of order~$240$ such that~$G'$
contains~$\wA_5$ with index two and such that the corresponding outer automorphism is given
by the class of~$\sigma$. The group~$G'$, however, is not a semi-direct product. The groups~$G$
and~$G'$ can be distinguished as follows: the group~$\PGSp_4(\F_3)$ has a natural action on~$40$
points corresponding to the action on~$\mathbf{P}^3(\F_3)$. The orbits of~$G$ are of size~$20$ and~$20$
respectively whereas~$G'$ acts transitively. 
\end{enumerate}

\begin{proof} This can be proved using the computer algebra package
  {\tt magma}~\cite{magma}. We omit the details. Note, however, that case~(\ref{case:2304}) was
proved in Lemma~\ref{lemma:biginduced}. 
\end{proof}

\begin{lem}
  \label{lem: Ind GL2 enormous}Suppose that~$p\ge 3$, that~$K/F$ is a
  quadratic extension such that~$K$ is unramified at~$p$,  and that $\rbar:G_K\to\GL_2(k)$ 
  restricted to~$G_{K(\zeta_p)}$ has image~$\SL_2(k)$.
  Choose~$\sigma\in G_F\setminus G_K$, and assume that
  $\Proj\rbar^\sigma\not\cong \Proj\rbar$, but that
  $\det \rbar^\sigma=\det \rbar$ is equal to~$\varepsbar^{-1}$.
  Let $\rhobar:=\Ind_{G_K}^{G_F}\rbar:G_F\to\GSp_4(k)$. 
  If~$p = 3$, assume that the fixed fields corresponding to the kernels of~$\Proj\rbar$
  and~$\Proj\rbar^{\sigma}$ are disjoint.
 Then~$\rhobar$ is vast and tidy.
  \end{lem}
\begin{proof} Note that~$F$ is necessarily
unramified at~$p$ (since~$K$ is). By Lemmas~\ref{lemma:eq5}
and~\ref{lemma:neq5}, in order to show that~$\rhobar$ is vast, it suffices
to show that for all~$N\ge 1$, the image of~$\rhobar$ restricted to~$F(\zeta_{p^N})$ is~$G =
\SL_2(k) \wr \Z/2\Z$. First assume that~$\#k > 3$. Then~$\PSL_2(k)$ is simple. If the image of~$\rbar^{\sigma}$
is disjoint from the image of~$\rbar$, it would follow by Goursat's Lemma that the image of~$\rhobar|_{G_{K(\zeta_p)}}$ 
is the group~$\SL_2(k)^2$,  and hence 
the image of~$\rhobar|_{G_{K(\zeta_{p^N})}}$ is also~$\SL_2(k)^2$,
and thus the image of~$\rhobar|_{G_{F(\zeta_{p^N})}}$ is~$\SL_2(k) \wr \Z/2\Z$.
 Since the automorphism group of~$\PSL_2(k)$ is~$\PGL_2(k)$, it follows
that  the projective representations associated to~$\rbar$ and~$\rbar^{\sigma}$ have the same
image if and only if they are the same. Since we are assuming
otherwise, we are done unless~$k=\F_3$.

Now assume that~$k = \F_3$, and so the images of~$\rbar$ and~$\rbar^{\sigma}$
restricted to~$K(\zeta_3)$ are both isomorphic to~$\SL_2(\F_3)$, which is a degree two
central extension of~$A_4$.  The  non-trivial
quotients of~$\SL_2(\F_3)$ are given by~$\PSL_2(\F_3) = A_4$ and $\Z/3\Z$. 
By assumption, the fixed fields corresponding to the kernels of~$\Proj\rbar$
  and~$\Proj\rbar^{\sigma}$ are disjoint and both have Galois group~$A_4$.
Thus by Goursat's lemma,
the image of~$\rhobar$ restricted to~$F(\zeta_3)$ is~$\SL_2(\F_3) \wr \Z/2\Z$.
This is enormous, by Lemma~\ref{lem:sample}. The abelianization of  this group
has order prime to~$3$, so
the image of~$\rhobar$ restricted to~$F(\zeta_{3^N})$ is also of this form.

Tidiness follows for~$p \ge 5$ by Lemma~\ref{lem: criterion for tidiness}. For~$p = 3$, the image contains
an element~$g$ of order~$8$ with~$\nu(g) = -1$ and eigenvalues~$(\zeta, - \zeta^{-1}, \zeta, -\zeta^{-1})$
for a primitive~$8$th root of unity~$\zeta$. The ratio of any two eigenvalues is either trivial or is a primitive fourth root
of unity. 
  \end{proof}

\begin{rem} When~$p =3$, we may weaken the hypotheses of this lemma slightly.
By Lemma~\ref{lem:sample} and Lemma~\ref{lemma:imageindependence}, it suffices that the image
of~$\rhobar$ restricted to~$F(\zeta_{27})$ is either~$\SL_2(\F_3) \wr \Z/2\Z$ or
one of the subgroups of~$\SL_2(\F_3) \wr \Z/2\Z$ of index three and order~$384$
considered in Lemma~\ref{lem:sample}. Unfortunately, the hypothesis  that~$\Proj\rbar$
is distinct from~$\Proj\rbar^{\sigma}$ is not \emph{quite} enough to force this.
For example, it is possible that the image of~$\rhobar$
restricted to~$F(\zeta_{3})$ might be the (unique) subgroup of order~$384$ in~$\SL_2(\F_3) \wr \Z/2\Z$ with abelianization~$\Z/6 \Z$,
and this means it is possible that the  image of~$\rhobar$ restricted to~$F(\zeta_{9})$
is the~$2$-Sylow subgroup of~$\SL_2(\F_3) \wr \Z/2\Z$ of order~$128$. However, this latter subgroup is not enormous.
\end{rem}

\subsubsection{Crossing with dihedral extensions} 
The goal of this section is to construct certain representations induced from quadratic fields~$K/F$ which will allow
us to prove modularity results for elliptic curves over~$K$ even when~$K$ is neither totally real nor CM (see 
Theorem~\ref{theorem:funky}).
Suppose that~$F$ is a totally real field in which~$p$ splits completely, and let~$K/F$ be an arbitrary quadratic extension of~$F$ in which~$p$
is unramified.

\begin{lem} \label{lem:di} There exists a Galois extension~$H/F$ containing~$K$ such that:
\begin{enumerate}
\item $D = \Gal(H/F)$ is the dihedral group of order~$8$, and~$\Gal(H/K) = (\Z/2\Z)^2$.
\item $H/F$ is the Galois closure over~$F$ of a quadratic extension~$M/K$.
\item $H/F$ is unramified at each~$v|p$, and~$\langle \Frob_v \rangle \in D$ is not central.
\end{enumerate}Furthermore, $H/F$ may be chosen to be linearly
disjoint from any given fixed finite extension of~$F$ linearly disjoint from~$K/F$. 

\end{lem}

\begin{proof} Let~$L/F$ be a second quadratic extension to be chosen later.
The obstruction to constructing a dihedral extension~$H/F$ containing~$K$ and~$L$ as quadratic subfields with~$\Gal(H/K) \simeq \Gal(H/L) \simeq (\Z/2 \Z)^2$
is the vanishing of the cup product~$\chi_K \cup \chi_L$,  where~$\chi_K, \chi_L \in H^1(F,\F_2)$ are the quadratic characters
corresponding to the fields~$K$ and~$L$. Equivalently, if~$L = F(\sqrt{\beta})$ and~$K = F(\sqrt{\alpha})$, it is the condition
of requiring that~$\beta \in N_{K/F}(K^{\times})$; if~$\beta = N(x + y \sqrt{\alpha})$, then one may take
$$H = F(\sqrt{\alpha},\sqrt{\beta},\sqrt{x + y \sqrt{\alpha}}).$$
The extension~$M = K(\sqrt{x + y \sqrt{\alpha}})$ will have Galois closure~$H$ over~$F$. 
Suppose~$\beta$ can be chosen so that every~$v|p$ is inert in~$L$, and
moreover such that $\beta$ is prime to~$p$. Then~$H/M$ will be
unramified at each~$v|p$,
and~$\Frob_v$ will be non-central, since the fixed field of the
non-trivial central element is the compositum~$K.L$. 

We now construct many such~$\beta$.
Note that~$F_v \simeq \Q_p$ by assumption, and we may assume that~$\alpha$ is a~$v$-adic unit for all~$v|p$.
 Let us consider ~$N_{K/F}(K^{\times})$, which
consists of the non-zero elements of~$F$ of the form~$x^2 - \alpha y^2$ where~$x, y \in F$.
The quadratic form
$$x^2 - \alpha y^2 - \gamma z^2 = 0 \mod v$$ 
for any~$\gamma \ne 0$  always has a non-trivial solution with~$z \ne 0$. Hence, by taking~$\gamma$ to be any quadratic non-residue in~$\F^{\times}_p$,
we may choose~$x$ and~$y$ modulo~$v$ so that~$x^2 - \alpha y^2$ is a non-zero quadratic non-residue. Making such a choice for all~$v|p$, we find that
~$ \beta = x^2 - \alpha y^2$ is a $v$-adic unit and a quadratic non-residue modulo~$v$ for all~$v|p$.
Since~$p > 2$, the resulting extension~$L = F(\sqrt{\beta})$ is thus  inert at all
primes~$v|p$, giving rise to the desired extension~$H/F$. 

Finally, by taking~$x$ and~$y$ sufficiently close to~$1$ and~$0$ respectively  in~$\OL_{F,w}$
for any finite set of auxiliary primes~$w$, can ensure that~$H/K$ splits completely at any such collection of primes, and hence we may ensure~$H/F$ is linearly disjoint from any fixed finite extension of~$F$
which is linearly disjoint from~$K$, as required.
\end{proof}

We may write~$D$ as~$D = \langle a,b | a^2=b^2= (ab)^4=1 \rangle$, where~$[a,b]$ is the order two element of the centre of~$D$.

\begin{lem} \label{lem:double}
Let~$\rbar: G_{F} \rightarrow  \GL_2(\Fbar_p)$ be an absolutely irreducible Galois representation with determinant~$\varepsbar^{-1}$.
Suppose that, for each~$v|p$,
 the restriction~$\rbar|_{G_{F_v}}$ takes the shape
$$\left( \begin{matrix} \chi_v & * \\  0 & \varepsbar^{-1} \chi^{-1}_v \end{matrix} \right)$$
for some unramified character~$\chi_v$.
Let~$K/F$ be an arbitrary quadratic extension linearly disjoint from the
fixed field~$F(\rbar)$ 
of the kernel of~$\rbar$, let~$H/F$ be any corresponding~$D$-extension
as guaranteed by Lemma~\ref{lem:di}, chosen to be linearly disjoint
from~$F(\rbar)$, and let~$M/K$ be a quadratic extension with Galois closure~$H/F$.
Let~$\rhobar$ be the following symplectic induction
$$\rhobar:= \Ind^{G_F}_{G_{K}}(\rbar |_{G_K} \otimes \delta_{M/K}),$$
where~$\delta_{M/K}$ is the quadratic character corresponding to the
extension~$M/K$, and the induction is constructed as
in~\S\ref{section:inductions}.  Let~$\Gamma$ denote the image
of~$\rbar$, and let~$G$ denote the image of~$\rhobar$.
Then:
\begin{enumerate}
\item $\rhobar$ is weight~$2$ ordinary and $p$-distinguished with similitude character~$\varepsbar^{-1}$.
\item If~$V$ denotes the underlying representation of~$G$ given by~$\rbar$, and~$U$ the~$2$-dimensional faithful representation of~$\Gal(H/F) = D$,
then~$\rhobar$ is given by~$V \otimes U$. In particular, $\rhobar$ is absolutely irreducible.
\item If~$\Gamma$ has a central element of order~$2$, then the image of~$\rhobar$ is
$$G = (\Gamma  \times D)/(-1 = [a,b]).$$
 Otherwise, the image is~$G = \Gamma \times D$.
 \end{enumerate}
\end{lem}

\begin{proof} 
The restriction of~$\rbar$ to~$G_{K}$ has determinant~$\varepsbar^{-1}$, which is preserved by the quadratic twist, and hence the induction also has~$\varepsbar^{-1}$
as the similitude character. The induction of~$\delta_{M/K}$
from~$G_{K}$ to~$G_{F}$ is precisely the representation~$U$ of~$D =
\Gal(H/F)$.  

By the construction of Lemma~\ref{lem:di}, for each
place~$v|p$, $\Frob_v \in D$ is not central. It follows
that the  restriction of~$\Gal(H/F)$ acting on~$U$ to the
decomposition group at~$v$ is of the form~$\psi \oplus \chi$ for \emph{distinct} unramified characters~$\psi$ and~$\chi$. Then the representation~$\rhobar|_{G_{F_v}}$ naturally
takes the form~$\rbar|_{G_{F_v}} \otimes\psi \oplus \rbar|_{G_{F_v}} \otimes \chi$.
This is automatically weight~$2$ ordinary and~$p$-distinguished. 

Finally, the image of~$\rhobar$ is the image of the map $\Gamma
\times D\to\GL(V\otimes U)$, and the kernel of this map is given by
the elements of the form~$(z,z^{-1})$ with~$z$ central. \end{proof}

We now show that many of the groups~$G$ occurring as the image of representations~$\rhobar$ as constructed in Lemma~\ref{lem:double} have big image.

\begin{lem}  \label{lem:formixed} Let~$p \ge 5$, and suppose that we are in the setting of
  Lemma~\ref{lem:double}. Suppose either that~$\Gamma = \GL_2(\F_p)$
  or that $p=5$ and~$\Gamma$ is the pre-image in~$\GL_2(\F_5)$ of~$S_4
  \subset S_5\cong \PGL_2(\F_5)$. 
Then~$G \cap \Sp_4(\F_p)$ is enormous unless~$\Gamma = \GL_2(\F_p)$ and~$p = 5$, in which case~$G \cap \Sp_4(\F_p)$ is weakly enormous.
 In any case, ~$\rhobar$ is vast and tidy.
\end{lem}

\begin{proof}By Lemma~\ref{lem:double}, $G$ acts faithfully on~$W =V \otimes U$, where~$V$ is the tautological representation of~$\Gamma$, and~$U$ is the faithful
$2$-dimensional representation of~$D$, with image~$G = (\Gamma \times D)/(-1 = [a,b])$.   We begin by checking condition~(E3). Clearly 
$$W \otimes W^* = (V \otimes V^*) \otimes (U \otimes U^*).$$
The latter factor is the regular representation of the abelianization of~$D$, and is a direct sum of characters of order dividing~$2$.
The first factor is the direct sum of~$\ad^0(V)$ with the trivial character. Both the
trivial representation and the adjoint representation
of~$\GL_2(\F_p)$ have the property
that~$1$ is always an eigenvalue of any element. Hence, for any irreducible summand of~$W \otimes W^*$,
$1$ will always be an eigenvalue on an index two subgroup~$\Sigma \subset G$ which is the kernel of one
of the degree~$2$ characters of~$D$. Yet
given~$g \in \Gamma$, there is an element in~$\Sigma$ with
eigenvalues the roots of~$g$ together with the negatives of the roots of~$g$. Hence it suffices to note that~$\Gamma \cap \SL_2(\F_p)$
has an element with eigenvalues~$\{\alpha,\alpha^{-1}\}$ with~$\alpha \ne \pm \alpha^{-1}$.
(In particular, in the case that~$p=5$ and~$\Gamma$ is the central
cover of~$S_4$, one could take~$g$ to have order~$3$.)

For~$\Gamma=\GL_2(\Fp)$, $p > 5$, since~$D$ has order prime to~$p$, (E1) reduces to the fact
that~$H^1(\SL_2(\F_p),\Sym^{2}(\F^2_p)) = 0$, which is~\cite[Lem.\ 2.48]{MR1605752}.
If~$p = 5$, the group~$G$ is of order~$384 = 4^2 |S_4|$, and therefore satisfies~(E1) automatically because the order is prime to~$p$.

For the final claim, note firstly that for each~$N\ge 1$, we have
$\rhobar(G_{F(\zeta_{p^N})})=G\cap\Sp_4(\Fp)$. Indeed this is clear
for~$N=1$ (as the similitude factor of~$\rhobar$ is~$\varepsilon^{-1}$),
and since~$G$ has no quotients of order~$p$, the same is true for
all~$N>1$. That~$\rhobar$ is
vast  is then an immediate consequence of the previous claims except in the case when~$\Gamma = \GL_2(\F_5)$, where~$G \cap \Sp_4(\F_5)$
is not enormous. But in this case, exactly as in the proof of Lemma~\ref{lemma:neq5}, the image of the projective representation factors through~$\PGL_2(\F_5) \times (\Z/2\Z)^2$
which does not surject onto~$\Z/4\Z$, and hence the fixed field of the adjoint representation cannot contain~$\zeta_5$ when~$E$
is unramified at~$p = 5$. Finally, for tidiness, we note that~$\Gamma$ and hence~$G$ contains a centre of order at least~$p-1$,
and we are done by Lemma~\ref{lem:bigcentre}.
\end{proof}

\subsection{Taylor--Wiles primes}We again fix a global deformation problem
  \[
   \cS = (\rhobar, S, \{\Lambda_v\}_{v\in S}, \psi,\{\cD_v\}_{v\in S}).
  \]
  Then we define a \emph{Taylor--Wiles datum} to be a tuple
  $(Q,(\alphabar_{v,1},\ldots,\alphabar_{v,4})_{v\in Q})$ consisting of:
\begin{itemize}
 \item A finite set of finite places $Q$ of $F$, disjoint from $S$, such that $q_v \equiv 1 \bmod p$ for each $v\in Q$.
 \item For each $v\in Q$, an ordering
   $\alphabar_{v,1},\alphabar_{v,2},\alphabar_{v,3}=\psibar(\Frob_v)\alphabar_{v,2}^{-1},\alphabar_{v,4}=\psibar(\Frob_v)\alphabar_{v,1}^{-1}$
   of the eigenvalues of $\rhobar(\Frob_v)$, which are assumed to be
   $k$-rational and pairwise distinct.
\end{itemize}
Given a Taylor--Wiles datum
$(Q,(\alphabar_{v,1},\ldots,\alphabar_{v,4})_{v\in Q})$, we define the
augmented global deformation problem
  \[
   \cS_Q = (\rhobar, S\cup Q, \{\Lambda_v\}_{v\in S} \cup \{\cO\}_{v\in Q},\psi, \{\cD_v\}_{v\in S} \cup \{\cD_v^\square\}_{v\in Q}).
  \]
  Set $\Delta_Q = \prod_{v\in Q} \Delta_v$. For each $v\in Q$,  the fixed ordering
  $\alphabar_{v,1},\ldots,\alphabar_{v,4}$, determines a
  $\Lambda[\Delta_Q]$-algebra structure on $R^T_{\cS_Q}$ for any
  subset $T$ of $S$ (via the homomorphisms $\cO[\Delta_v]\to
  R_v^\square$ defined in \S\ref{subsubsec: TW deformations}).  Letting
  $\fra_Q = \ker(\Lambda[\Delta_Q] \rightarrow \Lambda)$ be the
  augmentation ideal, the natural surjection
  $R^T_{\cS_Q} \rightarrow R^T_{\cS}$ has kernel $\fra_Q R^T_{\cS_Q}$.

\begin{lem}\label{thm:TWprimes}
 Assume  that $\rhobar$ is vast, that~$p\ge 3$ is unramified in~$F$,
 that~$\psi=\varepsilon^{-1}$, and that~$k$
 contains all of the eigenvalues of all elements of~$\rhobar(G_{F(\zeta_p)})$. 
 Let $q \ge h^1(F_S/F,\ad^0\rhobar(1))$. 
 Then for every $N\ge 1$, there is a choice of Taylor--Wiles datum $(Q_N,(\alphabar_{v,1},\ldots,\alphabar_{v,4})_{v\in Q_N})$ satisfying the following:
 \begin{enumerate}
  \item $\# Q_N = q$.
  \item For each $v\in Q_N$, $q_v \equiv 1 \bmod p^N$.\item $h_{S_{Q_N}^\perp,S}^1(\ad^0\rhobar(1)) = 0$.
 \end{enumerate}
\end{lem}

\begin{proof} Without loss of generality, we may assume that~$N \ge
  3$, and hence (by the definition of vastness and Remark~\ref{rem: vast unramified N 3})
that  ~$\rhobar(G_{F(\zeta_{p^N})})$ is weakly enormous.
By definition, we have \[
   H_{S_{Q_N}^\perp,S}^1(\ad^0\rhobar(1)) = \ker\left(H^1(F_S/F,\ad^0\rhobar(1)) \rightarrow \prod_{v\in Q_N} H^1(F_v,\ad^0\rhobar(1))\right).
  \]
 By induction, it suffices to show that given any cocycle $\kappa$ representing a nonzero element of $H_{\cS^\perp_{Q_N},S}^1(\ad^0\rhobar(1))$, 
 there are infinitely many finite places $v$ of $F$ such that
  \begin{itemize}
  \item $v$ splits in $F(\zeta_{p^N})$;
  \item $\rhobar(\Frob_v)$ has $4$-distinct eigenvalues $\alphabar_{v,1},\ldots,\alphabar_{v,4}$ in $k$;
  \item the image of $\kappa$ in $H^1(F_v,\ad^0\rhobar(1))$ is nonzero.
 \end{itemize}
 By Cebotarev, we are reduced to showing that given any cocycle $\kappa$ representing a nonzero element of $H^1(F_S/F,\ad^0\rhobar(1))$, 
 there is some $\sigma \in G_{F(\zeta_{p^N})}$ such that
  \begin{itemize}
   \item $\rhobar(\sigma)$ has distinct ($k$-rational) eigenvalues;
   \item $p_\sigma \kappa(\sigma) \ne 0$, where $p_\sigma: \ad^0\rhobar \rightarrow (\ad^0\rhobar)^\sigma$ is the $\sigma$-equivariant projection.
  \end{itemize}
  (The latter condition guarantees that the image of~$\kappa$ in
  $H^1(F_v,\ad^0\rhobar(1))$ is not a coboundary.)
 Let $L/F$ be the fixed field of $\ad^0 \rhobar$. 
The kernel of the restriction map 
  \[
   H^1(F_S/F,\ad^0\rhobar(1)) \rightarrow H^1(F_S/L(\zeta_{p^N}),\ad^0\rhobar(1))^{G_F}
  \]
  is, by inflation--restriction, isomorphic to
  $$H^1(\Gal(L(\zeta_{p^N})/F),\ad^0 \rhobar(1)).$$
  The assumption that~$\rhobar$ is vast implies 
  by Lemma~\ref{lemma:vanishingofcohomology} that this group vanishes.
 In particular, the restriction of $\kappa$ defines a nonzero
 $G_{F(\zeta_{p^N})}$-equivariant homomorphism
 $\Gal(F_S/L(\zeta_{p^N})) \rightarrow \ad^0\rhobar$. Let~$W$ be a
 nonzero irreducible sub-$G_{F(\zeta_{p^N})}$-representation of the
 $k$-span of $\kappa(\Gal(F_S/L(\zeta_{p^N}))$.
 Since~$\rhobar(G_{F(\zeta_{p^N})})$ is  weakly enormous and~$k$ is
 sufficiently large, there exists
 $\sigma_0 \in G_{F(\zeta_{p^N})}$ such that $\rhobar(\sigma_0)$ has
 distinct $k$-rational eigenvalues and such that  $W^{\sigma_0}\ne 0$ 
 (this follows from the vastness assumption, in particular,
 by condition~(E3) of~\ref{defn:enormous image}).
 This
 implies that $\kappa(\Gal(F_S/L(\zeta_{p^N}))$ is not contained in
 the kernel of the $\sigma_0$-equivariant projection
 $p_{\sigma_0}: \ad^0\rhobar \rightarrow
 (\ad^0\rhobar)^{\sigma_0}$.  If
 $p_{\sigma_0} \kappa(\sigma_0) \ne 0$, then we take
 $\sigma = \sigma_0$.  Otherwise, we choose
 $\tau \in G_{L(\zeta_{p^N})}$ such that
 $p_{\sigma_0}\kappa(\tau) \ne 0$, and we take
 $\sigma = \tau\sigma_0$;  since
 $\rhobar(\sigma) = \rhobar(\sigma_0)$ and
 $\kappa(\sigma) = \kappa(\sigma_0) + \kappa(\tau)$, we are done.\end{proof}
\begin{defn}\label{defn: odd Galois representation}We say
  that~$\rhobar:G_F\to\GSp_4(\Fpbar)$ is \emph{odd} if the similitude
  character~$\psi$ is odd, i.e.\ if for each place~$v|\infty$ of~$F$
  with corresponding complex conjugation~$c_v$, we have~$\psi(c_v)=-1$.
\end{defn}

\begin{cor}
  \label{cor: existence of TW primes with explicit numbers}Assume that~$\rhobar$ is odd, that
  $\rhobar$ is vast, and that~$k$ contains all of
  the eigenvalues of all elements of~$\rhobar(G_{F(\zeta_p)})$.  Let
  $q \ge h^1(F_S/F,\ad^0\rhobar(1))$.  Then for every
  $N\ge 1$, there is a choice of Taylor--Wiles datum
  $(Q_N,(\alphabar_{v,1},\ldots,\alphabar_{v,4})_{v\in Q_N})$ satisfying the
  following:
 \begin{enumerate}
  \item $\# Q_N = q$.
  \item For each $v\in Q_N$, $q_v \equiv 1 \bmod p^N$.\item There is a local $\Lambda$-algebra surjection  $R_{\cS}^{S,\loc}\llbracket X_1,\ldots,X_g \rrbracket \rightarrow
  R_{\cS_{Q_N}}^S$ with $g=2q-4[F:\Q]+\#S-1$.\end{enumerate}
\end{cor}
\begin{proof}
  By Proposition~\ref{prop: tangent space of local over global} and
  Theorem~\ref{thm:TWprimes}, the claim holds with~$g$ instead equal to
  \[ \#S-1 - \sum_{v| \infty} h^0(F_v,\ad^0\rhobar) + \sum_{v \in Q_N}  h^0(F_v,\ad^0\rhobar(1)).\]
  (Note that 
  the assumption that~$\rhobar$ is vast implies
  that $h^0(F_S/F,\ad^0\rhobar(1))=0$.) For $v\in Q_N$, by the
  assumptions that~$q_v\equiv 1\bmod p$ and that~$\rhobar|_{G_{F_v}}$
  has distinct eigenvalues we have \[h^0(F_v,\ad^0\rhobar(1))=h^0(F_v,\ad^0\rhobar)=2. \]
  For~$v|\infty$ we have $h^0(F_v,\ad^0\rhobar)=4$ by the assumption
  that~$\rhobar$ is odd. It follows that~$g=2q-4[F:\Q]+\#S-1$, as claimed. \end{proof}

\subsection{Global Galois deformation problems} \label{subsec: Galois
  deformation problems}
We now begin to introduce the framework that we need to carry out our Taylor--Wiles patching argument. As always, $F$ is a totally real field in which the prime $p\ge 3$ splits
    completely, and we write~$S_p$ for the set of primes of~$F$ dividing~$p$. Let $\rhobar:=G_F\to\GSp_4(k)$ be an absolutely irreducible
representation. We assume the following hypotheses.
\begin{hypothesis}\label{hyp: hypotheses on rhobar for deformation problems}
\leavevmode
  \begin{enumerate}
  \item The representation $\rhobar$ is vast and tidy.
 
    \item If $v\in S_p$, then  $\rhobar|_{G_{F_v}}$ is $p$-distinguished weight~$2$ ordinary.
 \item There is a set of
    finite places~$R$ of~$F$ which is disjoint from~$S_p$, such that
    \begin{enumerate}
    \item\label{hyp: conditions on R} If $v \in R$, then $\rhobar|_{G_{F_v}}$ is trivial, and $q_v
      \equiv 1 \pmod{p}$. If~$p=3$ then we further insist that~$q_v \equiv 1 \pmod 9$.
    \item If $v\notin S_p\cup R$, then  $\rhobar|_{G_{F_v}}$ is
      unramified.
    \end{enumerate}

\end{enumerate}
\end{hypothesis}

Set~$\psi=\varepsilon^{-1}$, and drop~$\psi$ from our notation for
global deformation problems from now on. Let $I\subset S_p$ be a set of places of cardinality
$\# I$. We will eventually need to assume that $\#I\le 1$,
although the more formal parts of the patching construction can be
carried out without this assumption, so we do not impose it yet. We
write $I^c$ for~$S_p\setminus I$.

By the Cebotarev density theorem and our assumption
that~$\rhobar(G_F)$ is tidy, we can find an  unramified place $v_0\notin R\cup S_p$ of
$F$ with the properties that
\begin{itemize}
\item   $q_{v_0} \not\equiv 1 \pmod{p}$, 
\item no two eigenvalues of $\rhobar(\Frob_{v_0})$ 
  have ratio~$q_{v_0}$, and
\item $v_0$ has residue characteristic greater than~$5$.
\end{itemize}
Then
$H^2(F_{v_0}, \ad \rhobar) = H^0(F_{v_0}, \ad \rhobar(1))^\vee = 0$.
We set $S = R\cup S_p \cup \{ v_0 \}$.

The reason for choosing~$v_0$ is that all liftings
of~$\rhobar|_{G_{F_{v_0}}}$ are automatically unramified by
Proposition~\ref{prop:smoothlifts}, and our choice of level structure
at~$v_0$ will guarantee that our level structures will be neat, by Lemma~\ref{lem:neat_subgroups}.

For each $v\in R$ we
choose a pair of characters $\chi_v = (\chi_{v,1},\chi_{v,2})$, where
$\chi_{v,i}: \cO_{F_v}^\times \rightarrow \cO^\times$  are trivial
modulo $\lambda$. (Note that at this stage the characters~$\chi_{v,i}$ are allowed to be
trivial.) We write $\chi$ for
the tuple $(\chi_v)_{v\in R}$ as well as for the induced character
$\chi = \prod_{v\in R} \chi_v: \prod_{v\in R} \Iw(v) \rightarrow
\cO^\times$. 

For each place~$v|p$, we fix~$\Lambda_v$, (and thus~$\theta_v$) as in~\S\ref{subsec:Ordinary
  lifts}, in the following way: if $v\in I$, then we take
$\Lambda_v=\cO\llbracket\cO_{F_v}^\times(p)\rrbracket$, while if $v\notin I$, then we
take~$\Lambda_v=\cO\llbracket(\cO_{F_v}^\times(p))^2\rrbracket$. We write
$\alphabeta=\{\alphabetabar_v\}_{v\in S_p}$ for a choice of $\alphabar_v$
or~$\betabar_v$ at each $v\in S_p$. 

 We have the
corresponding global deformation problem \begin{multline*} \cS_\chi^{I,\alphabeta} = (\rhobar, S, \{\Lambda_{v,1}\}_{v \in
      I}\cup \{\Lambda_{v,2}\}_{v \in I^c}\cup\{ \cO \}_{v \in
      S\setminus S_p},\\ \{ \cD_v^{P} \}_{v \in I}\cup
    \{\cD_v^{B,\alphabetabar_v} \}_{v \in I^c} \cup \{ \cD_v^\chi
    \}_{v \in R} \cup \{ \cD_{v_0}^\square \}).
  \end{multline*}
Let  $(Q,(\alphabar_{v,1},\ldots,\alphabar_{v,n})_{v\in Q})$ be a
choice of Taylor--Wiles datum. We set $S_{Q} = S \cup Q$ and define
the associated global deformation problem

\begin{multline*}
  \cS_{\chi,Q}^{I,\alphabeta} =(\rhobar, S_{Q},
  \{\Lambda_{v,1}\}_{v \in I}\cup \{\Lambda_{v,2}\}_{v \in I^c}\cup\{
  \cO \}_{v \in S\setminus S_p},\\ \{ \cD_v^{P} \}_{v \in I}\cup
  \{\cD_v^{B,\alphabetabar_v} \}_{v \in I^c} \cup \{ \cD_v^\chi \}_{v
    \in R} \cup \{ \cD_{v}^\square \}_{v\in S_{Q}\setminus (R\cup
    S_p)}).
\end{multline*} Note that by definition~$\cS_{\chi,Q}^{I,\alphabeta}$ does not depend
on the choice of~$\alphabeta_v$ for~$v\in I$.

\subsection{Taylor--Wiles systems: initial construction} \label{subsec: TW systems1}
In the next two sections, we will construct the Taylor--Wiles systems that we
will patch in~\S\ref{subsec: the
          actual patching construction}, using an abstract patching
        criterion explained in~\S\ref{prop: abstract patching}. (\S\ref{subsec: TW systems1}
          is mainly concerned with the construction of the Taylor--Wiles systems,
          whereas~\S\ref{subsec: TW systems2} is mainly concerned
          with proving the required local--global compatibility statements for the corresponding 
          Galois representations.)

Since we are only dealing with the cases that $\#I\le 1$,
we do not need to make use of the full machinery of patching
complexes developed in~\cite{CG,KT,geenew}; rather, we can and do use the notion of ``balanced''
modules introduced in~\cite[\S 2]{CG}, which we recalled in~\S\ref{subsec: balanced}. This has
the advantage that we do not need to consider local global
compatibility at places dividing~$p$ for Galois representations
associated to classes in higher degrees of cohomology, but rather just have to
prove the vanishing of the Euler characteristic of a certain
perfect complex, which follows from a calculation of the cohomology in
terms of automorphic forms. 

We now make the
following hypotheses on a representation $\rhobar:G_F\to\GSp_4(k)$, which include those made in Hypothesis~\ref{hyp:
  hypotheses on rhobar for deformation problems}. \begin{hypothesis}\label{hyp: hypotheses on rhobar for TW primes} 
\leavevmode
  \begin{enumerate}
  \item $F$ is a totally real field in which the prime $p\ge 3$ splits
    completely; we write~$S_p$ for the set of primes of~$F$ dividing~$p$.
  \item The representation~$\rhobar$ is vast and tidy.
    \item For each $v\in S_p$, $\rhobar|_{G_{F_v}}$ is $p$-distinguished weight~$2$ ordinary.
 \item There is a set of
    finite places~$R$ of~$F$ which is disjoint from~$S_p$, such that
    \begin{enumerate}
    \item\label{hyp: conditions on R second time} If $v \in R$, then $\rhobar|_{G_{F_v}}$ is trivial, and $q_v
      \equiv 1 \pmod{p}$. If~$p=3$, then~$q_v\equiv 1 \pmod 9$.
    \item If $v\notin S_p\cup R$, then moreover $\rhobar|_{G_{F_v}}$ is
      unramified.
    \end{enumerate}

  \item There is an  ordinary cuspidal automorphic representation~$\pi$
    of~$\GSp_4(\A_F)$ of parallel weight 2 with central character~$|\cdot|^2$ such that: \begin{enumerate}
    \item     $\rhobar_{\pi,p}\cong \rhobar$.
  \item If $v\in R\cup S_p$, then $\pi_v^{\Iw(v)} \neq 0$.
  \item If $v\notin R\cup S_p$, then $\pi_v^{\GSp_4(\cO_{F_v})}\ne 0$.
  \end{enumerate}

    \end{enumerate}

\end{hypothesis}
 As in~\S\ref{subsec: Galois
  deformation problems}, by the assumption that~$\rhobar(G_F)$
is tidy we can and do choose an  unramified place $v_0\notin R\cup S_p$
with the properties that
\begin{itemize}
\item   $q_{v_0} \not\equiv 1 \pmod{p}$, 
\item  no two eigenvalues of~$\rhobar(\Frob_{v_0})$ 
  have ratio~$q_{v_0}$, and
\item the residue characteristic of~$v_0$ is greater than~$5$.
\end{itemize}

\begin{df} \label{df:levelstructure}
We define an open compact subgroup $K^p = \prod_v K_v$ of $\GSp_4(\A^{\infty,p}_F)$ as follows:
\begin{itemize}
\item If $v \not\in S_p\cup R\cup \{v_0\}$, then $K_v = \GSp_4(\cO_{F_v})$.
\item If $v \in R\cup\{v_0\}$, then $K_v = \Iw_1(v)$.
\end{itemize}
For any Taylor--Wiles datum
$(Q,(\alphabar_{v,1},\ldots,\alphabar_{v,4})_{v\in Q})$, we have open
compact subgroups $K^p_0(Q)$, $K^p_1(Q)$ of~$K^p$ given by
\begin{itemize}
\item If $v\not\in Q$, then $K^p_0(Q)_v=K^p_1(Q)_v=K^p_v$.
\item If $v\in Q$, then $K^p_0(Q)_v=\Iw(v)$, $K^p_1(Q)_v=\Iw_1(v)$.
\end{itemize}
 We  define the open compact subgroup group~$K^p_0(Q, R)$ as follows:
\begin{itemize}
\item If~$v \not\in Q \cup R$, then~$K^p_0(Q,  R)_v = K^p_v$.
\item If~$v \in Q \cup R$, then~$K^p_0(Q, R)_v = \Iw(v)$.
\end{itemize}
Finally, we let~$K^p_1(Q,R) = K^p_1(Q)$. (Note that we already have~$K^p_1(Q)_v = \Iw_1(v)$ for~$v \in R$.)
\end{df} 

The following lemma (applied with $v=v_0$)
guarantees that for any compact open
subgroup~$K_p\subset \GSp_4(F_p)$, ~$K_pK^p_0(Q)$ and ~$K_pK^p_1(Q)$ are neat.
\begin{lemma}\label{lem:neat_subgroups}
  Suppose that $K = \prod_v K_v \subset \GSp_4(\A^\infty_F)$ is an
  open compact subgroup and that there exists a place $v$ of $F$ such
  that $v$ is absolutely unramified of residue characteristic greater than~$5$,
  and $K_v=\Iw_1(v)$. Then $K$ is neat.\end{lemma}\begin{proof}Suppose that there is an element $g_v\in K_v$ which has an
  eigenvalue $\zeta\in\overline{F}_v$ which is a root of unity; by the
  definition of ``neat'' (see Definition~\ref{defn:neat}), it is enough to check that we must
  have~$\zeta=1$. Since the reduction modulo~$v$ of the characteristic
  polynomial of~$g$ is $(X-1)^4$, the $v$-adic valuation of
  $(1-\zeta)$ is at least~$1/4$. On the other hand, if~$v$ has residue
  characteristic~$l$ and~$\zeta\ne 1$ is
  a root of unity, then the $v$-adic valuation of $(1-\zeta)$ is either~$0$, or is at
  most $1/(l-1)$, so we are done, as $l>5$ by assumption.
  \end{proof}

We let \[\TTt=\bigotimes_{v\not\in S_p\cup
R\cup\{v_0\}}\cO[\GSp_4(F_v) \doubleslash \GSp_4(\cO_{F_v})]\]
be the ring of spherical Hecke operators away from the bad places, and similarly we set \[\TTt^Q=\bigotimes_{v\not\in S_p\cup
R\cup\{v_0\}\cup
Q}\cO[\GSp_4(F_v) \doubleslash \GSp_4(\cO_{F_v})].\]
We let
$\widetilde{\m}^{\mathrm{an}}\subset\TTt$ be the maximal ideal
corresponding to $\overline{\rho}$ (the ``an'' stands for
``anaemic''); so by definition~$\m$ containins~$\lambda$, 
and the polynomials $\det(X-\rhobar(\Frob_v))$ and~$Q_v(X)$ are
congruent modulo~$\m$
for each $v\not\in S_p\cup
R\cup\{v_0\}$, where in a slight abuse of notation, if~$v\notin S_p\cup R\cup\{v_0\}$ we
write $Q_v(X)\in \TTT[X]$ for the
polynomial \[X^4-T_{v,1}X^3+(q_vT_{v,2}+(q_v^3+q_v)T_{v,0})X^2-q_v^3T_{v,0}T_{v,1}X+q_v^6T_{v,0}^2\] 
 (cf.\ (\ref{eqn: char poly
  for unramified Hecke})). Similarly we write
$\widetilde{\m}^{\mathrm{an},Q}\subset\TTt^Q$ for the maximal ideal corresponding to $\overline{\rho}$.  For any
choice of $I$ we let \begin{equation*}
\TTt^I=\TTt[\{U_{v,0},U_{\Kli(v),1},U_{v,2}\}_{v\in I},\{U_{v,0},U_{v,1},U_{v,2}\}_{v\in I^c}]
\end{equation*}
and
\begin{equation*}
\TTt^{I,Q}=\TTt^Q[\{U_{v,0},U_{\Kli(v),1},U_{v,2}\}_{v\in I},\{U_{v,0},U_{v,1},U_{v,2}\}_{v\in I^c}]
\end{equation*}
and additionally for any choice of $\alphabeta$ we let $\widetilde{\m}^{I,\alphabeta}\subset
\TTt^I$ be the maximal ideal
\numequation\label{eqn-def-max-ideal}
\widetilde{\m}^{I,\alphabeta}=(\widetilde{\m}^{\mathrm{an}},\{U_{v,0}-1,U_{v,2}-\alpha_v\beta_v\}_{v\in
S_p},\{U_{\Kli(v),1}-\alpha_v-\beta_v\}_{v\in
I},\{U_{v,1}-\alphabeta_v\}_{v\in
I^c})
\end{equation}
and we let $\widetilde{\m}^{I,\alphabeta,Q}\subset \TTt^{I,Q}$ be the maximal ideal
\begin{equation*}
\widetilde{\m}^{I,\alphabeta,Q}=(\widetilde{\m}^{\mathrm{an},Q},\{U_{v,0}-1,U_{v,2}-\alpha_v\beta_v\}_{v\in
S_p},\{U_{\Kli(v),1}-\alpha_v-\beta_v\}_{v\in
I},\{U_{v,1}-\alphabeta_v\}_{v\in
I^c}).
\end{equation*}

Let $\chi=(\chi_{v,1},\chi_{v,2})_{v\in R}$ be any choice of $p$-power
order characters of $I_{F_v}$ for $v\in R$, and also write $\chi_v$
for the corresponding characters of $T(k(v))$ given by
$\chi_{v,1}\circ\Art_{F_v}$, $\chi_{v,2}\circ\Art_{F_v}$.

Then we consider the $\Lambda_I$-module
\begin{equation*}
M^{\chi,I,\alphabeta}=\RHom^0_{\Lambda_I}(M^{\bullet,I}_{K^p},\Lambda_I)_{\widetilde{\m}^{I,\alphabeta},\chi,|\cdot|^2},
\end{equation*}
and the $\Lambda_I[\Delta_Q]$-module
\begin{equation*}
M^{\chi,I,\alphabeta,Q}=\RHom^0_{\Lambda_I}(M^{\bullet,I}_{K^p_1(Q)},\Lambda_I)_{\widetilde{\m}^{I,\alphabeta,Q},\widetilde{\m}_Q,\chi,|\cdot|^2},
\end{equation*}
where:
\begin{itemize}
\item $M^{\bullet,I}_{K^p}$ denotes the complex~$M_I^\bullet$ defined in
  Theorem~\ref{theorem-p-adic-complex}, at tame level~$K^p$.\item The localizations
  $\widetilde{\m}^{I,\alphabeta},\tilde{\m}^{I,\alphabeta,Q}$ are
  defined above.  
\item The localization $\widetilde{\m}_{Q}$ is with respect to the
  maximal ideals~$\widetilde{\m}_v$ of the subalgebras $\cO[T(F_v)/T(\cO_{F_v})_1]$ of the pro-$v$ Iwahori Hecke
  algebras~$\cH_1(v)$ for~$v\in Q$ as considered in \S\ref{subsubsec q
  equals 1}, so that~$\lambda\in\widetilde{\m}_v$,  $U_{v,0}-1\in\widetilde{\m}_v$,
  and $U_{v,1}$ and~$U_{v,2}$ are respectively congruent
  to~$\alphabar_{v,1}$, $\alphabar_{v,1}\alphabar_{v,2}$ modulo~$\widetilde{\m}_v$.
\item The subscript~$\chi$ denotes that we take the
  $\chi$-coinvariants for the action of~$\prod_{v\in R}T(k(v))$.
\item The subscript~$|\cdot|^2$ denotes that we are fixing the central
  character, by taking coinvariants under~$T_{v,0}-q_v^{-2}$ for all
  $v\notin S_p\cup R\cup\{v_0\}$.\end{itemize}

The following lemma motivates our definition using
 ~$\RHom^0_{\Lambda}(M^{\bullet},\Lambda)$, and will be useful for proving various properties
 of $M^{\chi,I,\alphabeta,W}$ below (see also Remark~\ref{rem:
   relationship of RHom to CG}). \begin{lem}\label{lem: justification of RHom} 
Let~$\Lambda\in\CNL_\cO$,
  and let~$M^\bullet$ be a perfect complex of
  $\Lambda$-modules bounded below by 0. Set
  $M:=\RHom^0_\Lambda(M^\bullet,\Lambda)$. Then, writing~$*$ for
  the usual duality of finite-dimensional vector spaces and~${}^\vee$ for Pontryagin duals, we have
  \begin{enumerate}
  \item    \label{bulletone} $M\otimes_\Lambda
    k=(H^0(M^\bullet\otimes_\Lambda^{\mathbf{L}} k))^*$.
    \item For any homomorphism of $\cO$-algebras $\Lambda\to E$,
    $M\otimes_\Lambda E=(H^0(M^\bullet\otimes_\Lambda^{\mathbf{L}} E))^*$.
  \item   \label{bulletCG} For any homomorphism of $\cO$-algebras $\Lambda\to \cO$,
    \[M\otimes_\Lambda\cO=\Hom(H^0(M^\bullet\otimes_\Lambda^{\mathbf{L}}E/\cO),E/\cO) = H^0(M^\bullet\otimes_\Lambda^{\mathbf{L}}E/\cO)^{\vee}.\]\end{enumerate}
\end{lem}\begin{proof}Let $P^\bullet = P^0 \rightarrow P^1 \rightarrow \ldots \rightarrow P^{l_0}$ be a bounded complex of finite projective
  $\Lambda$-modules which is bounded below by~$0$ and is
  quasi-isomorphic to~$M^\bullet$. Then, by definition, we have an exact sequence
  \[\Hom_\Lambda(P^1,\Lambda)\to\Hom_\Lambda(P^0,\Lambda)\to
    M\to 0. \]
In particular it follows that for any $\Lambda$-algebra~$R$, we have
an exact sequence of $R$-modules \[\Hom_R(P^1\otimes_\Lambda
  R,R)\to\Hom_R(P^0\otimes_\Lambda R,R)\to
    M\otimes_\Lambda R\to 0.  \] On the other hand, by definition, we have an exact
  sequence of $R$-modules \[0\to H^0(M^\bullet\otimes_\Lambda^{\mathbf{L}}R)\to
    P^0\otimes_\Lambda R \to P^1\otimes_\Lambda R,\] and therefore, for a field~$R = F$,
   an
  exact sequence 
   \[\Hom_F(P^1\otimes_\Lambda
  F,F)\to\Hom_F(P^0\otimes_\Lambda F,F)\to
  \Hom_F(H^0(M^\bullet\otimes_\Lambda^{\mathbf{L}}F),F)\to
  0.  \]Parts~(1) and~(2) follow immediately with~$F = E$ or~$F = k$. Part~(3) follows from
Lemma~\ref{lem: pairing of O modules} below, applied to the morphism
$P^0\otimes_\Lambda\cO\to P^1\otimes_\Lambda\cO$.
  \end{proof}

\begin{lem} 
  \label{lem: pairing of O modules}If $\phi:M\to N$ is a morphism of
  finite free $\cO$-modules, and~$\phi_{E/\OL} = \phi \otimes E/\OL$ is the map~$M \otimes E/\OL \rightarrow N \otimes E/\OL$, then
the Pontryagin dual~$\phi^{\vee}_{E/\OL}$ of~$\phi_{E/\OL}$ is the map
$$\phi^{\vee}_{E/\OL}: \Hom(N,\OL) \rightarrow \Hom(M,\OL).$$
 In particular, the Pontryagin dual of~$\ker(\phi_{E/\OL})$ is~$\coker(\phi^{\vee}_{E/\OL})$.
  \end{lem}
\begin{proof}  Because~$M$ and~$N$ are  free, the Pontryagin duals of~$M \otimes E/\OL$ and~$N \otimes E/\OL$
are~$\Hom(M,\OL)$ and~$\Hom(N,\OL)$ respectively, and the result follows immediately.
\end{proof}

\begin{rem}
  \label{rem: relationship of RHom to CG} In~\cite{CG} and~\cite{CGGSp4},
  the patched modules are constructed by first taking cohomology
  with coefficients in~$E/\cO$ and then
  taking Pontryagin duals.
 Lemma~\ref{lem: justification
    of RHom}~(\ref{bulletCG}) explains how
    our construction coincides with this in the special case when~$\Lambda = \cO$.
    \end{rem}

\begin{defn}\label{defn: weight kappa}
  For any $I\subset S_p$, a \emph{weight} is a
  homomorphism $\kappa:\Lambda_I\to\cO$; by definition, $\kappa$
  corresponds to a tuple $(\theta_{v,1},\theta_{v,2})_{v\in S_p}$
  where $\theta_{v,i}:I_{F_v}\to \cO^\times$ is a character with
  trivial reduction, and moreover $\theta_{v,1}=\theta_{v,2}$ for
  $v\in I$.  We let $\p_\kappa\subset\Lambda_I$ denote the kernel of
  this homomorphism.

We say that~$\kappa$ is \emph{classical} if there are integers $k_v\ge
l_v\ge 2$ such that~$\theta_{v,1}=\varepsilon^{(k_v+l_v)/2-2}$,
$\theta_{v,2}=\varepsilon^{(k_v-l_v)/2}$ (so that $k_v\equiv l_v\equiv 2$
or~$p+1\pmod{2(p-1)}$, and if $v\in I$, we must
have~$l_v=2$). If~$\kappa$ is classical, then we write~$\omega^\kappa$
for the automorphic vector bundle corresponding to~$(k_v,l_v)_{v\in
  S_p}$, as in~\S\ref{subsec: auto vector bundles}.

  For any $I$ we denote by $\kappa_2$ the classical algebraic weight
  where~$k_v=l_v=2$ for all~$v$.  
 For $I=\emptyset$ we
  pick some sufficiently regular classical algebraic weight,
  $\kappa_{\mathrm{reg}}$; for example, we could choose the one given by
  the characters $\theta_{v,1}=\varepsilon^{2N(p-1)}$ and
  $\theta_{v,2}=\varepsilon^{N(p-1)}$ for all $v\in S_p$, where~$N$ is
  sufficiently large. 

\end{defn}
\begin{rem}
  \label{rem: how regular is regular}In practice we
  choose~$\kappa_{\mathrm{reg}}$ so
  that we can apply Theorems~\ref{thm: cohomology in terms of
    automorphic forms} and~\ref{thm: ordinary classicity for I at most 1} in weight~$\kappa_{\mathrm{reg}}$. We will do
  this without comment from now on. \end{rem}

We will now prove some very important properties of the action of
$\Delta_Q$ on the modules that we patch. It will also be important for
us to understand the action of the diamond operators at the places
in~$R$ (that is, at the places involved in the ``Ihara avoidance''
argument). We can and do treat the places in~$Q$ and in~$R$
simultaneously; recall, by Definition~\ref{df:levelstructure}, we have the
groups~$K^p_1(Q)$ and~$K^p_0(Q, R)$ such that
\begin{itemize}
\item If $v\not\in Q\cup R$, then $K^p_0(Q, R)_v = K^p_v$ which equals~$K^p_1(Q)_v$.
\item If $v\in Q \cup R$, then $K^p_0(Q, R)_v=\Iw(v)$ which contains~$K^p_1(Q)_v = \Iw_1(v)$.
\end{itemize}

In particular, there is an inclusion~$K^p_1(Q) \subset K^p_0(Q, R)$.
In contexts in which we particularly want to emphasize the fact that~$K^p_1(Q)$
has level structure~$\Iw_1(v)$ at~$v \in R$, we write~$K^p_1(Q, R) = K^p_1(Q)$.

Let $K_p$ be any reasonable level structure at $p$ (for example $K_p(I)$). Let  $X_{K_p K_0^p(Q, R), \Sigma}$ be the Shimura variety of the corresponding level $K_p K_0^p(Q, R)$ for a choice $\Sigma$ of good polyhedral cone decomposition. Over the interior $Y_{K_p K_0^p(Q, R)}$ we have for all $v \in Q\cup R$ a flag of subgroups $0 \subset H_v \subset L_v \subset H_v^\bot \subset A[v]$ and all the graded pieces are \'etale $k(v)$-group schemes of rank $1$. 
 We now consider the Shimura variety $X_{K_p K_1^p(Q, R), \Sigma}$ for the same choice of cone decomposition. 

\begin{prop}\label{prop-fixingQ} \leavevmode

 \begin{enumerate}
\item  For all $v \in Q\cup R$, the groups $H_v$, $L_v/H_v$, $H_v^\bot/L_v$ and $A[v]/(H_v^\bot)$ extend to finite \'etale $k(v)$-group schemes of rank $1$ over $X_{K_p K_0^p(Q, R), \Sigma}$.
\item The map $X_{K_p K_1^p(Q, R), \Sigma} \rightarrow X_{K_p K_0^p(Q, R),
    \Sigma}$ is finite \'etale with group $\prod_{v\in Q\cup R}T(k(v))$, and $X_{K_p K_1^p(Q, R), \Sigma}$ identifies with the torsor of trivializations of the groups $H_v$, $L_v/H_v$, $H_v^\bot/L_v$ and $A[v]/(H_v^\bot)$, compatible with duality. 
\end{enumerate}
\end{prop}

\begin{proof} We observe that when $F = \Q$, this is the content of~ \cite[\S2.4.5]{Bord}.
The argument can be adapted to our setting. The
  extension problem is local so let us pick $\sigma \in \Sigma$ and
  consider the completion
$(X_{K_p K_0^p(Q, R),
    \Sigma})^{\wedge}_\sigma 
    \simeq
\Spf~R$ of $X_{K_p K_0^p(Q, R), \Sigma}$ along the $\sigma$-stratum. The
semi-abelian scheme $A$ over $\Spf~R$ is obtained by Mumford's
construction as the quotient of a semi-abelian scheme~$B$ of constant
toric rank by a finite free $\ocal_F$-module $X_\sigma$.  Let
$U_\sigma \hookrightarrow \Spec R$ be the Zariski open complement of
the boundary and let us consider any of the groups $H_v$, $L_v/H_v$,
$H_v^\bot/L_v$ or $A[v]/(H_v^\bot)$. 
If this group is a subquotient of $B[v]$, then since $B$ exists over
all $\Spec R$ and $B[v]$ is a finite \'etale group scheme,  the group
extends  as a subquotient of  $B[v]$. Otherwise, the group maps
isomorphically to its image in $A[v]/B[v] = X_\sigma \otimes_{\ocal_F}
k(v)$ and is constant over $U_\sigma$. Therefore it extends to the
constant group scheme. This proves $(1)$.

We may now define a scheme $X_{K_p K_1^p(Q, R), \Sigma}' \rightarrow
X_{K_p K_0^p(Q, R), \Sigma}$ as the   torsor of trivializations of the
(extended) groups $H_v$, $L_v/H_v$, $H_v^\bot/L_v$ and
$A[v]/(H_v^\bot)$, compatible with duality for all $v|p$. This
scheme is canonically isomorphic to $X_{K_p K_1^p(Q, R), \Sigma}$ because
the two schemes are generically equal, and both are normal, and finite flat  over $X_{K_p K_0^p(Q, R), \Sigma}$.\end{proof}

\begin{prop}\label{prop: free and balanced over Lambda Delta Q} \leavevmode
\begin{enumerate}
\item\label{item: free if I is empty} $M^{\chi,\emptyset,\alphabeta,Q}$ is a finite free $\Lambda_\emptyset[\Delta_Q]$-module.
\item\label{item: balanced if I has size 1} If~$\#I=1$, then $M^{\chi,I,\alphabeta,Q}$ is a balanced
  $\Lambda_{I}[\Delta_Q]$-module. \end{enumerate}
\end{prop}
\begin{proof} The complex $M^{\bullet,I}_{K^p_1(Q)}$  (which is the 
  complex~$M_I^\bullet$ defined in
  Theorem~\ref{theorem-p-adic-complex} for the tame level $K^p_1(Q)$)
  is a perfect complex of $\Lambda_{I}$-modules of amplitude $[0, \#
  I]$. We claim that it is actually a perfect complex of $\Lambda_{I}[
  \prod_{v \in Q\cup R}T(k(v))]$-modules of amplitude $[0, \# I]$.
  
The complex $M^{\bullet,I}_{K^p_1(Q)}$   is obtained by considering
the  cohomology   over $\mathfrak{X}_{K_pK_1^p(Q), \Kli(p^\infty)}$ of
the sheaf of $\Lambda_I$-modules $\Omega^{\kappa_I}(-D)$ and applying
the ordinary idempotent.
Equivalently, it is obtained by  considering the  cohomology   over
$\mathfrak{X}_{K_pK_0^p(Q), \Kli(p^\infty)}$ of the sheaf of
$\Lambda_I$-modules $$\Omega^{\kappa_I}(-D)
\otimes_{\ocal_{\mathfrak{X}_{K_pK_0^p(Q, R), \Kli(p^\infty)}}}
\ocal_{\mathfrak{X}_{K_pK_1^p(Q, R), \Kli(p^\infty)}}$$ and applying  the ordinary idempotent. Using the
independence of  the cohomology with respect to choices of toroidal
compactifications, we may assume that we are in the setting of
Proposition \ref{prop-fixingQ}, so that the morphism
$\mathfrak{X}_{K_pK_1^p(Q, R), \Kli(p^\infty)} \rightarrow
\mathfrak{X}_{K_pK_0^p(Q, R), \Kli(p^\infty)}$ is finite \'etale
with group $\prod_{v \in Q\cup R} T(k(v))$. 
   Therefore, it follows  (by considering a suitable  \'etale covering
   to compute the cohomology) that $M^{\bullet,I}_{K^p_1(Q)}$ is
   represented by a bounded complex of flat complete~$\Lambda_{I}[
   \prod_{v \in Q\cup R}T(k(v))]$-modules.  We can apply Lemma~ \ref{lem:
     finite cohomology implies perfect} (or rather its straightforward
   extension to the semi-local situation; see also~\cite[Prop.\
   2]{MR728135}) to conclude that $M^{\bullet,I}_{K^p_1(Q)}$ is a
   perfect complex of $\Lambda_{I}[ \prod_{v \in Q\cup
     R}T(k(v))]$-modules of amplitude $[0, \# I]$.  It follows that  the corresponding complex~$M^{\bullet,I,\chi}_{K^p_1(Q)}$
   is  a perfect complex of~$\Lambda_I[\Delta_Q]$-modules,
   also of amplitude~$[0,\#I]$.

 Given an ideal~$\widetilde{\m}^{I,\alphabeta}$, 
  we can f localize the complex with respect to the action of a 
  lift of a  suitable idempotent for this ideal in the Hecke algebra,
 and this localization 
 also preserves the property of being perfect of the correct amplitude.
 (The endomorphism ring at the level of derived categories of a perfect complex of~$\Lambda_I$
  modules is a finite~$\Lambda_I$ module. So, if one has a commutative subalgebra, it is a
  semi-local ring. See the discussion following~\cite[Lem.\ 2.12]{KT}
  for a lengthier treatment of such localizations.)

It remains to consider the passage to coinvariants under the centre. To this end, consider the spaces
\begin{equation*}
\MM^{\chi,I,\alphabeta}=\RHom^0_{\Lambda_I}(M^{\bullet,I}_{K^p},\Lambda_I)_{\widetilde{\m}^{I,\alphabeta},\chi},
\end{equation*}
\begin{equation*}
\MM^{\chi,I,\alphabeta,Q}=\RHom^0_{\Lambda_I}(M^{\bullet,I}_{K^p_1(Q)},\Lambda_I)_{\widetilde{\m}^{I,\alphabeta,Q},\widetilde{\m}_Q,\chi},
\end{equation*}
obtained \emph{before} taking coinvariants under the centre. The
component groups of our Shimura varieties are indexed by a finite
abelian (ray) class group~$C = F^{\times} \backslash \A^{\times}_F/U$
for some~$U$. The action
of~$\gamma \in A^{\times}_F$ on components is via the class~$[\gamma]^2$, and the action of the central character on our cohomology groups
is via~$| \cdot |^2$ times a character of~$C$.
  Let~$C = C_p \oplus C^p$, where~$C_p$ is the~$p$-Sylow subgroup of~$C$.
There are always natural isomorphisms of~$\OL[\Delta_Q]$ modules
$$\MM^{\chi,I,\alphabeta} \simeq M^{\chi,I,\alphabeta} \otimes_{\OL} \OL[C_p],$$
$$\MM^{\chi,I,\alphabeta,Q} \simeq M^{\chi,I,\alphabeta,Q} \otimes_{\OL} \OL[C_p],$$
 with~$\Delta_Q$ acting trivially on the second factor. 
The reason for such an isomorphism is that, after localization at a maximal ideal~$\m$ of the Hecke algebra, the elements of the centre which
act through an element of order prime to~$p$ are already determined, because they are fixed modulo~$\m$ and the polynomial~$T^m - 1$ is
separable modulo~$p$ if~$(m,p) =1$. On the other hand, if we consider only the connected components corresponding to the subgroup~$C^p$,
the entire space is canonically isomorphic to~$|C_p|$ copies of this space, and moreover, the action of~$C/C^p = C_p$ on these components
is transitive and fixed point free (this crucially uses that~$p \ne 2$).  Hence working with the~$| \cdot |^2$ part of the cohomology is simply
equivalent to working with the components indexed by~$C^p$ instead of~$C$, and the passage between the cohomology (or complexes)
for either of these two spaces (even before localization) is simply to tensor with~$\OL[C_p]$.

Part~\eqref{item: free if I is empty} follows immediately from these considerations,
  because~$M^{\bullet,I,\chi}_{K^p_1(Q)}$ is perfect of amplitude~$[0,0]$.
By Lemma~\ref{lem:
  presentations of balanced modules}, to prove
  part~\eqref{item: balanced if I has size 1},
  it is enough to prove that the corresponding
perfect complex (of amplitude~$[0,1]$) has Euler characteristic~$0$ after localization at~$\widetilde{\m}^{I,\alphabeta}$. We can check this modulo
any prime ideal of~$\Lambda_I$,
so the result follows from Theorem~\ref{thm: ordinary classicity for I
    at most 1} and Corollary~\ref{cor: dimensions of cohomology in various degrees}.
  \end{proof}

\subsection{Taylor--Wiles systems: local-global compatibility}
\label{subsec: TW systems2}

We write $\TT^{\chi,I,\alphabeta}$ for the $\Lambda_I$-subalgebra of
$\End_{\Lambda_I}(M^{\chi,I,\alphabeta})$ generated by the image of
$\TTt^I$.  Similarly, we write $\TT^{\chi,I,\alphabeta,Q}$ for the
$\Lambda_I$-subalgebra of $\End_{\Lambda_I}(M^{\chi,I,\alphabeta,Q})$
generated by the image of $\TTt^{I,Q}$. We remind the reader that none
of these objects depend on the choice of~$\alphabeta_v$ for~$v\in I$
(but they do depend on the choice of~$\alphabeta_v$ for~$v\notin I$). If~$v\in I$, then by Hensel's lemma and our assumption that
$\alphabar_v\ne\betabar_v$, we can write
\[X^2-U_{\Kli(v),1}X+U_{v,2}=(X-\tilde{\alpha}_v)(X-\tilde{\beta}_v)\]
where $\alphat_v,\betat_v\in\TT^{\chi,I,\alphabeta,Q}$ are
respectively lifts of~$\alphabar_v,\betabar_v$.

If $I\subset I'$, then there is a natural surjective
map~$\Lambda_I\to\Lambda_{I'}$, corresponding to the closed immersion
$\Spec\Lambda_{I'}\to\Spec\Lambda_I$ given by
$\theta_{v,1}=\theta_{v,2}$ for all~$v\in I'$.
 Then we have the following key doubling statement:
\begin{prop}[Doubling]\label{prop: doubling for $I$ to $I'$.}For each
  choice of~$\alphabeta$, and each $I\subset I'$, there are natural surjections
\begin{equation*}
M^{\chi,I,\alphabeta}\otimes_{\Lambda_I}\Lambda_{I'}\to M^{\chi,I',\alphabeta}
\end{equation*}
and
\begin{equation*}
M^{\chi,I,\alphabeta,Q}\otimes_{\Lambda_I}\Lambda_{I'}\to M^{\chi,I',\alphabeta,Q}
\end{equation*}
which commute with all the Hecke operators away from
$I'\setminus I$. Furthermore, if $v\in I'\setminus I$, then these
surjections are equivariant with respect to~$U_{v,0}$ and~$U_{v,2}$,
and intertwine the actions of~$U_{v,1}$ on the source
and~$\widetilde{\alphabeta}_v$ on the target.\end{prop}

\begin{proof}[Proof of Proposition~\ref{prop: doubling for $I$ to
    $I'$.}]
We give the proof for ~$M^{\chi,I,\alphabeta}$, as the argument
for ~$M^{\chi,I,\alphabeta,Q}$ is identical. By induction, it suffices to consider the case that~$I'=I\cup\{v\}$
 for some $v\notin I$. We have 
 a map of complexes 
 \[M^{\bullet,I'}_{K_p}\to
   M^{\bullet,I}_{K_p},\]induced by the restriction map coming from
 the  inclusion
  \[\mathfrak{X}^{G_1,I}_{K, \Kli(p^\infty)}\into
   \mathfrak{X}^{G_1,I'}_{K, \Kli(p^\infty)}, \] together with the natural map
 $\Lambda_I\to\Lambda_{I'}$. This induces a
 map \[M^{\chi,I,\alphabeta}\otimes_{\Lambda_I}\Lambda_{I'}\stackrel{\tau}{\to}
   M^{\chi,I',\alphabeta}, \]and the map that we are seeking is the
 map $\tau\circ U_{v,1}-\widetilde{\alphabeta}'_v\circ\tau$. It is clear that this satisfies all of the claimed properties except
possibly for the surjectivity and the claimed intertwining of~$U_{v,1}$ 
and~$\widetilde{\alphabeta}_v$. 

To see the intertwining, it is convenient to introduce the
module~$M^{\chi,I',\alphabeta,=_v2}$, whose definition is 
\begin{equation*}
M^{\chi,I',\alphabeta,=_v2}=\RHom^0_{\Lambda_{I'}}(M^{\bullet,I}_{K^p}\otimes_{\Lambda_I}\Lambda_{I'},\Lambda_{I'})_{\widetilde{\m}^{I',\alphabeta},\chi,|\cdot|^2};
\end{equation*} that is, it is defined in the same way
as~$M^{\chi,I,\alphabeta}$, but we  are now over the weight
space~$\Lambda_{I'}$, rather than~$\Lambda_I$, and we localize with respect to the
Hecke operator~$(U_{\Kli(v),1}-(\alpha_v+\beta_v))$, rather
than~$(U_{\Iw(v),1}-\alphabeta_v)$. By Lemma~\ref{lem: quadratic relation for U Kli 1 big sheaf}, on~$M^{\chi,I',\alphabeta,=_v2}$ we have the  identity
$$U_{v,1} (U_{\Kli(v),1} - U_{v,1}) = U_{v,2},$$
or equivalently (writing~$\{\alpha_v,\beta_v\}=\{\alphabeta_v,\alphabeta'_v\}$) the identity \numequation\label{eqn: quadratic relation alphabeta yeah}(U_{v,1}-\widetilde{\alphabeta}'_v)(U_{v,1}-\widetilde{\alphabeta}_v)= 0.\end{equation} We need to show that~$U_{v,1}=\widetilde{\alphabeta}'_v$ on
 $M^{\chi,I,\alphabeta}\otimes_{\Lambda_I}\Lambda_{I'}$. Now, noting
 that $M^{\chi,I,\alphabeta}\otimes_{\Lambda_I}\Lambda_{I'}$ is a
 subspace of~$M^{\chi,I',\alphabeta,=_v2}$ (because it is obtained
 from it by
 localizing with respect to $(U_{\Iw(v),1}-\alphabeta_v)$, and
 because~\eqref{eqn: quadratic relation alphabeta yeah} holds
 on~$M^{\chi,I',\alphabeta,=_v2}$), we see that~\eqref{eqn: quadratic relation alphabeta yeah} also holds
 on $M^{\chi,I,\alphabeta}\otimes_{\Lambda_I}\Lambda_{I'}$;
 since~$U_{\Iw(v),1}$ acts via~$\alphabeta_v$ modulo the maximal ideal
 of~$\TT^{\chi,I,\alphabeta}$, it follows from Hensel's lemma
 that~$U_{v,1}=\alphabeta_v$ on
 $M^{\chi,I,\alphabeta}\otimes_{\Lambda_I}\Lambda_{I'}$, as required.

It only remains to check the surjectivity. By Nakamaya's lemma, it is
enough to check surjectivity modulo~$\m_{\Lambda_{I'}}$,
or equivalently (by Lemma~\ref{lem: justification of RHom}) the injectivity of the map
 \numequation\label{eqn: doubling injection special fibre}
 \begin{diagram}
  e(U^{I'})H^0(X^{I',G_1}_{K^pK_p(I'),1},\omega^{2}(-D))_{\tilde{\m}^{I',\alphabeta},\chi,|\cdot|^2}
  \\
  \dTo_{ (U_{v,1}-\alphabeta'_v)} \\
    e(U^{I})
 H^0(X^{I,G_1}_{K_p(I) K^p,1},\omega^{2}(-D))_{\tilde{\m}^{I,\alphabeta},\chi,|\cdot|^2} \end{diagram}
 \end{equation}
on the special fibre. This follows from Theorem~\ref{thm: the doubling
  map is injective}, as in Remark~\ref{rem: doubling as mapping to
  generalized eigenspaces}.
\end{proof}

Recall from~\S\ref{subsec:Ordinary
  lifts} that if $v\in I$, we defined a character~$\theta_v:I_{F_v}\to
\Lambda_v^\times$, and if $v\notin I$ we defined a pair of
characters~$\theta_{v,1},\theta_{v,2}:I_{F_v}\to\Lambda_v^\times$. We
extend all of these characters to~$G_{F_v}$ by sending~$\Art_{F_v}(p)\mapsto
1$. In the following theorem, we allow the Taylor--Wiles datum
$(Q,(\alphabar_{v,1},\ldots,\alphabar_{v,4})_{v\in Q})$ to be
empty. 

\begin{thm}
  \label{thm: existence of Galois representations in big Hecke
    algebras}There is a unique continuous representation
  \[\rho^{\chi,I,\alphabeta,Q}:G_F\to\GSp_4(\TT^{\chi,I,\alphabeta,Q})\]
   which is a deformation of $\overline{\rho}$ of type
   $\cS_{\chi,Q}^{I,\alphabeta}$ such that the induced homomorphism
   $R_{\cS_{\chi,Q}^{I,\alphabeta}}\to\TT^{\chi,I,\alphabeta,Q}$ is a
   homomorphism of $\Lambda_I[\Delta_Q]$-algebras, and moreover such that
\begin{enumerate}
\item If $v\notin S_p\cup R\cup\{v_0\}\cup Q$, then $\det(X-\rho^{\chi,I,\alphabeta,Q}(\Frob_v))=Q_v(X)$.\item If $v\in I$, then
\begin{equation*}
\rho^{\chi,I,\alphabeta,Q}|_{G_{F_v}}\simeq
\begin{pmatrix}
\lambda_{\tilde{\alpha}_v}\theta_v&0&*&*\\
0&\lambda_{\tilde{\beta}_v}\theta_v&*&*\\
0&0& \lambda_{\tilde{\beta}_v}^{-1}\theta_v^{-1}\varepsilon^{-1}&0\\
0&0&0&\lambda_{\tilde{\alpha}_v}^{-1}\theta_v^{-1}\varepsilon^{-1}
\end{pmatrix}.
\end{equation*}
\item If $v\in I^c$, then
\begin{equation*}
\rho^{\chi,I,\alphabeta,Q}|_{G_{F_v}}\simeq
\begin{pmatrix}
\lambda_{U_{v,1}}\theta_{v,1}&*&*&*\\
0&\lambda_{U_{v,2}/U_{v,1}}\theta_{v,2}&*&*\\
0&0& \lambda_{U_{v,2}/U_{v,1}}^{-1}\theta_{v,2}^{-1}\varepsilon^{-1}&*\\
0&0&0&\lambda_{U_{v,1}}^{-1}\theta_{v,1}^{-1}\varepsilon^{-1}
\end{pmatrix}.
\end{equation*}
\end{enumerate}
\end{thm}\begin{proof}First we treat the case $I=\emptyset$.  By Proposition~\ref{prop: free and balanced over Lambda Delta Q},
$M^{\chi,\emptyset,\alphabeta,Q}$ is a finite free
$\Lambda_\emptyset$-module, so there is an injection of $\TTt^{\emptyset,Q}$-modules
\begin{equation*}
M^{\chi,\emptyset,\alphabeta,Q}\to\prod_{\kappa} M^{\chi,\emptyset,\alphabeta,Q}\otimes_{\Lambda_\emptyset,\kappa}E
\end{equation*}
where the product is over all weights $\kappa=(k_v,l_v)_{v|\infty}$
with $k_v\geq l_v\ge 4$, $k_v\equiv l_v\equiv 2$ or $p+1 \pmod{2(p-1)}$. (Note that
these points are scheme-theoretically dense in
$\Spec\Lambda_\emptyset$.)

From the definition of $M^{\chi,\emptyset,\alphabeta,Q}$,
Lemma~\ref{lem: justification of RHom}, and Theorem~\ref{thm: ordinary classicity for I at most 1}, we have
\begin{equation*}
(M^{\chi,\emptyset,\alphabeta,Q}\otimes_{\Lambda,\kappa} E)^\vee=e(\emptyset)H^0(X^{G_1}_{K^p_1(Q)K_p(\emptyset)},\omega^\kappa)_{\widetilde{\m}^{I,\alphabeta,Q}}^{\{T(k(v))=\chi_v\}_{v\in R},|\cdot|^2},
\end{equation*} and by Theorem~\ref{thm: cohomology in terms of automorphic forms}, we have
\begin{equation*}H^0(X^{G_1}_{K^p_1(Q)K_p(\emptyset)},\omega^\kappa)\otimes\barE\simeq\bigoplus_\pi (\pi_f^{K_p^1(Q)K_p(\emptyset)})\otimes\barE,
\end{equation*}
where in the sum, $\pi$ runs over all the cuspidal automorphic representations of weight $(k_v,l_v)$, with $\pi_v$ holomorphic for each $v|\infty$, and $\pi_f$ is the finite part of $\pi$.

Next we observe that for such a $\pi$, if the $\TTt^{\emptyset,Q}$-module
\begin{equation*}
e(\emptyset)(\pi_f^{K^p_1(Q)K_p(\emptyset)}\otimes\barE)_{\widetilde{\m}^{I,\alphabeta,Q}}^{\{T(k(v))=\chi_v\}_{v\in R},|\cdot|^2}
\end{equation*}
is nonzero, then $\pi$ has central character $|\cdot|^2$ and by Proposition \ref{prop: ordinary eigenform in autrep}, $\pi$ is ordinary, and moreover $\TTt^{\emptyset,Q}$ acts on it through a character $\Theta_\pi:\TTt^{\emptyset,Q}\to \barE$, and the ordinary Hecke parameters are $(\Theta_\pi(U_{v,1}),\Theta_\pi(U_{v,2}/U_{v,1}))$.

We now argue as in the proof of~\cite[Prop.\ 3.4.4]{cht}.  By Theorem \ref{thm: existence and properties of Galois
  representations for automorphic representations}, Proposition~\ref{prop: iwahori up eigenvectors converse},
Proposition~\ref{prop: Ihara local global},  Proposition~\ref{prop:
  TW local global},  and Remark~\ref{remark:becausewewantto}, there is a Galois representation $\rho_{\pi,p}:G_F\to\GSp_4(\barE)$ such that
\begin{itemize}
\item If $v\not\in S_p\cup R\cup\{v_0\}\cup Q$, then $\rho_{\pi,p}|_{G_{F_v}}$ is unramified and $\det(X-\rho_{\pi,I}(\Frob_v))=\Theta_\pi(Q_v(X))$.
\item If $v\in S_p$, then
\begin{equation*}
\rho_{\pi,p}|_{G_{F_v}}\simeq
\begin{pmatrix}
\lambda_{\Theta_\pi(U_{v,1})}\theta_{v,1}&*&*&*\\
0&\lambda_{\Theta_\pi(U_{v,2}/U_{v,1})}\theta_{v,2}&*&*\\
0&0& \lambda_{\Theta_\pi(U_{v,2}/U_{v,1})}^{-1}\theta_{v,2}^{-1}\varepsilon^{-1}&*\\
0&0&0&\lambda_{\Theta_\pi(U_{v,1})}^{-1}\theta_{v,1}^{-1}\varepsilon^{-1}
\end{pmatrix}
\end{equation*}
\item If $v\in R$, then for all $\sigma\in I_{F_v}$,  $\det(X-\rho(\sigma))$ is equal to
\begin{equation*}
(X-\chi_{v,1}(\Art_{F_v}^{-1}(\sigma)))(X-\chi_{v,1}(\Art_{F_v}^{-1}(\sigma))^{-1}(X-\chi_{v,2}(\sigma))(X-\chi_{v,2}(\Art_{F_v}^{-1}(\sigma))^{-1}).
\end{equation*}
\item If $v\in Q$, then
\begin{equation*}
\rho|_{G_{F_v}}\simeq\gamma_{v,1}\oplus\gamma_{v,2}\oplus\varepsilon^{-1}\gamma_{v,2}^{-1}\oplus\varepsilon^{-1}\gamma_{v,1}^{-1}
\end{equation*}
for characters $\gamma_{v,i}:G_{F_v}\to \barE^\times$ satisfying
$\overline{\gamma}_i=\lambda_{\overline{\alpha}_{v,i}}$. Furthermore
$T(F_v)$ acts on
$(\pi^{\Iw_1(v)})_{\m_{\overline{\alpha}_1,\overline{\alpha}_2}}$
via the characters $\gamma_{v,i}\circ\Art_{F_v}$.\end{itemize}

After conjugation, we may assume that~$\rho_{\pi,p}$ is valued
in~$\cO_{E_\pi}$ for some finite extension $E_\pi/E$, and
since~$\rhobar_{\pi,p}\cong\rhobar$, we may assume after further
conjugation that $\rhobar_{\pi,p}=\rhobar$. Let~$A$ be the
subring of~$k\oplus\bigoplus_{\pi}\cO_{E_\pi}$ consisting of
those elements~$(a,(a_\pi)_\pi)\in k\oplus\bigoplus_{\pi}\cO_{E_\pi}$ such that for all~$\pi$ the reduction
of~$a_\pi$ modulo the maximal ideal of~$\cO_{E_\pi}$ is equal
to~$a$ (where the direct sum is over the infinitely many~$\pi$ corresponding
to the infinitely many~$\kappa$). Then~$A$ is a local
$\Lambda$-algebra with residue field~$k$ (with the~$\Lambda$-algebra
structure coming from that on~$\cO_{E_\pi}$ given
by~$\kappa$). Set \[\rho_A:=\rhobar\oplus\bigoplus_\pi\rho_{\pi,p}:G_F\to\GSp_4(A).\]There
is a natural injection $\TT^{\chi,\emptyset,\alphabeta,Q}\to A$ (this
map is injective because~$\TT^{\chi,\emptyset,\alphabeta,Q}$ is
reduced, by a standard argument using Proposition~\ref{prop: ordinary eigenform in autrep}). We can
choose (for example, by ordering the~$\kappa$) a decreasing sequence of ideals $I_n$ of~$A$
with~$\cap_nI_n=(0)$ such that each $A/I_n$  is an object of
$\CNL_\Lambda$, and
it follows from~\cite[Lem.\ 7.1.1]{gg} that for each~$n$ the
representation~$\rho_A\otimes_A{A/I_n}$ is
$\ker(\GSp_4(A/I_n)\to\GSp_4(k))$-conjugate to a representation
\[\rho^{\chi,\emptyset,\alphabeta,Q}_n:G_F\to\GSp_4(\TT^{\chi,\emptyset,\alphabeta,Q}/(I_n\cap
\TT^{\chi,\emptyset,\alphabeta,Q})).\] After possibly conjugating again, we can assume
that~$\rho^{\chi,\emptyset,\alphabeta,Q}_{n+1}\pmod{I_n}=\rho^{\chi,\emptyset,\alphabeta,Q}_{n}$,
and we set $\rho^{\chi,\emptyset,\alphabeta,Q}:=\varprojlim_n \rho^{\chi,\emptyset,\alphabeta,Q}_{n}$.
By construction this satisfies the required properties at places~$v\not\in
S_p\cup Q$ (in particular, at the places $v\in R$, the
deformation is of the required type by the definition of~$R_v^\chi$).

It remains to verify the claimed properties
of~$\rho^{\chi,\emptyset,\alphabeta,Q}|_{G_{F_v}}$ for $v\in S_p\cup
Q$. Suppose  that~$v\in S_p$.   We claim firstly that it is enough to show
that there are elements $\nu_{v,1},\nu_{v,2}\in
(\TT^{\chi,\emptyset,\alphabeta,Q})^\times$ such that
\numequation\label{eqn: local global at primes dividing p not in I} \rho^{\chi,\emptyset,\alphabeta,Q}|_{G_{F_v}}\cong\begin{pmatrix}
\lambda_{\nu_{v,1}}\theta_{v,1}&*&*&*\\
0&\lambda_{\nu_{v,2}}\theta_{v,2}&*&*\\
0&0& \lambda_{\nu_{v,2}}^{-1}\theta_{v,2}^{-1}\varepsilon^{-1}&*\\
0&0&0&\lambda_{\nu_{v,1}}^{-1}\theta_{v,1}^{-1}\varepsilon^{-1}
\end{pmatrix}\end{equation}Indeed, if this holds, then the equalities
$\nu_{v,1}=U_{v,1}$ and~$\nu_{v,2}=U_{v,2}/U_{v,1}$  can be checked
after composing with the
injection~$\TT^{\chi,\emptyset,\alphabeta,Q}\into A$, where they
follow from local-global compatibility for
the~$\rho_{\pi,p}|_{G_{F_v}}$. Now, (\ref{eqn: local global at
  primes dividing p not in I}) is equivalent to asking that the
homomorphism~$R_v^\square\to \TT^{\chi,\emptyset,\alphabeta,Q}$
corresponding to~$\rho^{\chi,\emptyset,\alphabeta,Q}|_{G_{F_v}}$
factors through the quotient~$R_v^{B,\alphabetabar_v}$, and this can
again be checked after composing with the
injection~$\TT^{\chi,\emptyset,\alphabeta,Q}\into A$, as required.

Suppose now that~$v\in Q$, so that  we need to check that the morphism
$R_{\cS_{\chi,Q}^{I,\alphabeta}}\to\TT^{\chi,I,\alphabeta,Q}$ is $\Delta_v$-equivariant. By Lemma~\ref{lem:anyliftisTW}, there are unique
characters
$\gamma_{v,1},\gamma_{v,2}:G_{F_v}\to (\TT^{\chi,\emptyset,\alphabeta,Q})^\times$
lifting~$\lambda_{\alphabar_{v,1}}$, $\lambda_{\alphabar_{v,2}}$
respectively such
that
$\rho^{\chi,\emptyset,\alphabeta,Q}|_{G_{F_v}}\cong\gamma_{v,1}\oplus\gamma_{v,2}\oplus\gamma_{v,2}\varepsilon^{-1}\oplus\gamma_{v,1}\varepsilon^{-1}.$
We claim that the action of~$T(F_v)$ on~$M^{\chi,I,\alphabeta,Q}$ is
given by~$\gamma_{v,1}\circ\Art_{F_v}$,
$\gamma_{v,2}\circ\Art_{F_v}$; this can be checked after composing
with the injection ~$\TT^{\chi,\emptyset,\alphabeta,Q}\into A$, so it
follows from the analogous result for~$\rho_{\pi,p}|_{G_{F_v}}$
recalled above. Restricting this claim to~$T(\cO_{F_v})$ gives the
result.

We are done in the case that~$I=\emptyset$. We now prove the result for general~$I$ by induction
on~$\# I$. Accordingly, assume that the result holds for some~$I\ne S_p$,
choose~$w\in I^c$, and set~$I'=I\cup\{w\}$. By Proposition~\ref{prop: doubling for $I$ to $I'$.} we have a natural
surjection of
$\Lambda_{I'}$-algebras \[\TTT^{\chi,I,\alphabeta,Q}\otimes_{\Lambda_I}\Lambda_{I'}\onto
  \TTT^{\chi,I',\alphabeta,Q},  \] and we
let~$\rho^{\chi,I',\alphabeta,Q}$ be the pushforward of $\rho^{\chi,I,\alphabeta,Q}$. It follows from the result
for~$I$ that we need only check that property~(2) holds
for~$v=w$. However, we could equally well have performed the same
construction with~$\alphabeta_w$ replaced with~$\alphabeta_w'$ (the
two candidates for~$\rho^{\chi,I',\alphabeta,Q}$ are conjugate by
  property~(1), the Cebotarev density theorem, and~\cite[Lem.\
  7.1.1]{gg}), so from~(3) and the equivariance properties for Hecke
  operators at~$w$ in Proposition~\ref{prop: doubling for $I$ to
    $I'$.}, we see that $\rho^{\chi,I',\alphabeta,Q}|_{G_{F_w}}$
  admits both~$\lambda_{\alphat_w}\theta_w$
  and~$\lambda_{\betat_w}\theta_w$ as
  subcharacters. Since~$\alphabar_w\ne\betabar_w$, the result follows.
\end{proof}
As a corollary, we have the following result about Galois representations
associated to automorphic representations of parallel weight~$2$. As
ever, some of the hypotheses in this result could be relaxed (in
particular, the assumption that~$\rhobar_{\pi,p}$ is vast and tidy can
presumably easily be relaxed to irreducibility), but in
the interests of brevity we have contented ourselves with this result,
as it is sufficient for our purposes.

\begin{cor}\label{cor: Galois representations for limits of discrete series
    with ordinarity at p}Let~$\pi$ be a cuspidal automorphic
  representation of~$\GSp_4(\A_F)$ of parallel weight~$2$ with central
  character~$|\cdot|^2$. Fix a prime~$p>2$, and assume
that~$\pi$ is ordinary. Then there
is a continuous semisimple representation
$\rho_{\pi,p}:G_F\to\GL_4(\Qpbar)$ such that \begin{enumerate}
\item For each finite place~$v\nmid p$, at which~$\pi_v$ is
  unramified, $\rho_{\pi,p}|_{G_{F_v}}$ is unramified and  \[
    \det(X-\rho_{\pi,p}(\Frob_v))=Q_v(X).\]
  \end{enumerate}
Suppose further that~$\rhobar_{\pi,p}$ is vast and tidy, and that for each $v|p$, the ordinary Hecke parameters $\alpha_v,\beta_v$ of $\pi_v$ satisfy $\alphabar_v\not=\betabar_v$. Then~$\rho_{\pi,p}$ can be conjugated to be valued
in~$\GSp_4(\Qpbar)$, and 
\begin{enumerate}[resume]
\item  $\nu\circ\rho_{\pi,p}=\varepsilon^{-1}$. 

\item For each finite place~$v\nmid p$, we
  have  \[\WD(\rho_{\pi,p}|_{G_{F_v}})^{\semis}\cong
    \recGTp(\pi_v\otimes|\nu|^{-3/2})^{\semis}.\] 

\item For each place~$v|p$,
  then \[\rho_{\pi,p}|_{G_{F_v}}\cong  \begin{pmatrix}
    \lambda_{\alpha_v}&0&*&*\\
0 &\lambda_{\beta_v}&*&*\\
0&0&\lambda_{\beta_v}^{-1}\varepsilon^{-1}&0\\
0&0&0&\lambda_{\alpha_v}^{-1}\varepsilon^{-1}
  \end{pmatrix}.\] 
\end{enumerate}
\end{cor}
\begin{proof}This could be proved by repeating the arguments
  of~\cite[\S4]{MR3200667}, using Theorem~\ref{thm: existence of Galois representations in big Hecke
    algebras} instead of the results of~\cite{MR3361018}. For brevity,
we instead explain how to deduce the result from~\cite[Thm.\
4.14]{MR3200667} and Theorem~\ref{thm: existence of Galois representations in big Hecke
    algebras}.

Firstly, if~$\pi$ is not of general type in the sense
of~\cite{MR2058604}, then the existence of a (unique) semisimple
\emph{reducible} representation~$\rho_{\pi,p}$ satisfying~(1)
is an easy consequence of standard results on Galois representations
for~$\GL_1$ and~$\GL_2$ (see the proof of Lemma~\ref{lem: getting to general type}), and parts~(2)-(4) are then vacuous. 

Accordingly, for the remainder of the proof we assume that~$\pi$ is of
general type, in which case the existence of a
representation~$\rho_{\pi,p}$ satisfying~(1) and~(3) follows
from~\cite[Thm.\ 4.14]{MR3200667}, except that this representation is
only given to be valued in~$\GL_4(\Qpbar)$ rather
than~$\GSp_4(\Qpbar)$.  

Choose a
solvable extension of totally real fields~$F'/F$, linearly disjoint
from $\overline{F}^{\ker\rhobar_{\pi,p}}$ over~$F$, with the properties
that~$p$ splits completely in~$F'$, and that there is an
automorphic representation ~$\Pi$ of
~$\GSp_4(\A_{F'})$ of parallel weight 2 and central
character~$|\cdot|^2$, which is a base change
of~$\pi$ (that is, for each finite place~$w$ of~$F'$, lying over a
place~$v$ of~$F$, we have $\recGTp(\Pi)=\recGTp(\pi)|_{W_{F'_w}}$),
which is holomorphic at all infinite places,  and which satisfies $\Pi_w^{\Iw(w)}\ne 0$
for all finite places~$w$ of~$F'$ (the existence of such an~$F'$ and $\Pi$
follows from~\cite[Prop.\ 4.13]{MR3200667}). 

 We claim that if
~$\rho_{\Pi,p}$ admits a symplectic pairing with
multiplier~$\varepsilon^{-1}$, then so
does~$\rho_{\pi,p}$. Indeed,
since~$\rho_{\Pi,p}=\rho_{\pi,p}|_{G_{F'}}$ is irreducible,
it admits at most one perfect pairing with
multiplier~$\varepsilon^{-1}$; while by~(1), $\rho_{\pi,p}$ admits a
perfect pairing with multiplier~$\varepsilon^{-1}$, which must
therefore also be symplectic.
In addition~(4) holds for~$\rho_{\Pi,p}$ if and only if it holds
for~$\rho_{\pi,p}$. Replacing~$F$
by~$F'$ and~$\pi$ by~$\Pi$, we can and do assume that~$\pi_v^{\Iw(v)}\ne 0$
for all finite places~$v$ of~$F$.  

Taking~$\rhobar:=\rhobar_{\pi,p}$, we see that
Hypothesis~\ref{hyp: hypotheses on rhobar for TW primes} holds, so the required
properties of~$\rho_{\pi,p}$ follow immediately from
Theorem~\ref{thm: existence of Galois representations in big Hecke
  algebras}, taking $I=S_p$, $\chi=1$ and~$Q=\emptyset$. (Note
that as in the proof of Theorem~\ref{thm: existence of Galois representations in big Hecke
  algebras}, it follows from  Theorem~\ref{thm: cohomology in terms of
  automorphic forms} that ~$\pi$ contributes to~$M^{1,S_p,\alphabeta,\emptyset}$.)\end{proof}

We now turn to the final lemmas that we need to prove in order to
construct our Taylor--Wiles systems. \begin{lem}  \label{lem: Nakayama for perfect complexes}Let
  $\Lambda\in \CNL_\cO$, and let
  $f^\bullet:\cC^\bullet\to\cD^\bullet$ be a morphism of bounded
  complexes of $\m_\Lambda$-adically complete and separated flat
  $\Lambda$-modules. Supposed that the induced morphism
  $\cC^\bullet\otimes_\Lambda^{\mathbf{L}}\Lambda/\m_\lambda\to
  \cD^\bullet\otimes_\Lambda^{\mathbf{L}}\Lambda/\m_\lambda$ is a
  quasi-isomorphism. Then~$f^\bullet$ is a quasi-isomorphism.

\end{lem}
\begin{proof}See~\cite[Prop.\ 2.2]{pilloniHidacomplexes}.
\end{proof}

\begin{prop}\label{prop: level Q to level 1}
The natural map $M^{\chi,I,\alphabeta,Q}\to M^{\chi,I,\alphabeta}$ induces an isomorphism $(M^{\chi,I,\alphabeta,Q})_{\Delta_Q}\to M^{\chi,I,\alphabeta}$.
\end{prop}
\begin{proof}We follow the proof of~\cite[Lem.\ 6.25]{KT}. 
We claim that we have natural isomorphisms
\numequation\label{eqn: Q coinvariants from 1 to 0}
(M^{\chi,I,\alphabeta,Q})_{\Delta_Q}\isoto M^{\chi,I,\alphabeta,Q}_{K_0^p(Q)}
\end{equation}
and
\numequation\label{eqn: Q iso to level 0}
M^{\chi,I,\alphabeta,Q}_{K_0^p(Q)}\isoto M^{\chi,I,\alphabeta}
\end{equation}
whose composite is the claimed isomorphism. We begin with~(\ref{eqn:
  Q coinvariants from 1 to 0}). It suffices to show that we have a
natural isomorphism in the derived
category \[(M_{K_1^p(Q)}^{\bullet,I})^{\prod_{v\in Q}T(k(v))}\isoto
  M^{\bullet,I}_{K_0^p(Q)}. \] 
  As in the proof of Proposition~\ref{prop: free and balanced over
  Lambda Delta Q}, the complex on the left (before taking invariants)
  is a perfect complex of~$\Lambda_I[\prod_{v \in Q} T(k(v))]$-modules.
  But now the result is immediate from Proposition~\ref{prop-fixingQ}, as the map~$X_{K_p K_1^p(Q), \Sigma} \rightarrow X_{K_p K_0^p(Q),
    \Sigma}$ is finite \'etale with group $\prod_{v\in Q}T(k(v))$.

We now turn to proving~(\ref{eqn: Q iso to level 0}). Again, we mostly
work
on the level of complexes. We begin by considering the
composite   \[(M^{\bullet,I}_{K_0^p(Q)})_{\tm^{\an,Q},\tm_Q}\to
  (M^{\bullet,I}_{K_0^p(Q)})_{\tm^{\an,Q}}\to (M^{\bullet,I}_{K^p})_{\tm^{\an,Q}}.\]
By Lemma~\ref{lem: Nakayama for perfect complexes}, these maps induce  quasi-isomorphism
of complexes if the following maps are isomorphisms
\[H^*(M^{\bullet,I}_{K_0^p(Q)} \otimes k)_{\tm^{\an,Q},\tm_Q}\to
  H^*(M^{\bullet,I}_{K_0^p(Q)} \otimes k)_{\tm^{\an,Q}}\to
  H^*(M^{\bullet,I}_{K^p} \otimes k)_{\tm^{\an,Q}}.\] 
 This follows formally from Lemmas~\ref{lem: eigenvalues at level q are in the
  support} and~\ref{lem: bijection between localizations at level q and
  level one}, applied at each place in~$Q$, because, for~$K = \GSp_4(\cO_{F_v})$ and~$K' = \Iw(v)$,
   we have the
identities of Hecke
operators
 \[ [{K 1 K'}]   [{K'  1 K }]  =[K:K'] \] 
 \[
[{K' 1 K}]   [{K  1 K' }]
=e_{K} = e_{\GSp_4(\cO_{F_v})},\]and
  we note that
$[{K 1 K'}]$ is the trace from level~$K'$
to~$\GSp_4(\cO_{F_v})$ and $ [{K'  1 K }]$ is the
inclusion from level~$K$ to level~$K'$ (recall
that since~$p>2$, $[K:K'] = [\GSp_4(\cO_{F_v}):\Iw(v)]$ is not divisible
by~$p$).

   Finally, consider the natural map
$$(M^{\bullet,I}_{K^p})_{\tm^{\an,Q}} \rightarrow
(M^{\bullet,I}_{K^p})_{\tm^{\an}}.$$  For our purposes, it suffices to
prove that this map becomes an isomorphism after
applying~$\RHom^0(-,\Lambda)$. Since this map is a localisation, it
suffices to check that it is an isomorphism modulo the maximal ideal
of~$\Lambda$; so by 
Lemma~\ref{lem:
  justification of RHom}~(1), it is in turn enough to prove that 
\[H^0(M^{\bullet,I}_{K^p} \otimes k)_{\tm^{\an,Q}}\to
  H^0(M^{\bullet,I}_{K^p} \otimes k)_{\tm^{\an}}\] is an isomorphism, or in other words, that~$\tm^{\an}$ is the unique maximal ideal~$\mathfrak{n}$ of~$\TTT$
  lying over~$\tm^{\an,Q}$ and
  in the support of $H^0(M^{\bullet,I}_{K^p} \otimes k)$.
  Equivalently, we need to show that the Hecke eigenvalues away from the primes in~$Q$ (which are prime to the level)
  determine the Hecke eigenvalues at~$Q$. This follows from the fact that the Hecke eigenvalues at primes of good
  reduction and residue characteristic different from~$p$ are
  determined by the Galois representation (exactly as in the proof of Theorem~\ref{thm: existence of Galois representations in big Hecke
    algebras}, this local-global compatibility statement for~$H^0$ is
  a  consequence of the corresponding local-global compatibility statement for the
  Galois representations in Theorem~\ref{thm: existence and properties of Galois representations for
    automorphic representations}).
But the Galois representation itself is determined 
 from~$\tm^{\an,Q}$ by the Cebotarev density theorem. Hence~$\mathfrak{n}=\tm^\an$, as required.
  \end{proof}

\subsection{An abstract patching criterion}We have the following slight
variant on~\cite[Prop.\ 2.3, Prop.\ 6.6]{CG} (although our 
formulation is also informed by~\cite[Prop.\ 3.1]{KT}); we leave the details of the proof as an exercise
for the interested reader. 

\begin{prop}\label{prop: abstract patching}Let $l_0$ be equal to either $0$ or~$1$, let~$\Lambda\in\CNL_{\cO}$,
let $S_\infty:=\Lambda\llbracket x_1,\dots,x_q\rrbracket$ for some $q\ge 1$, and set
$\mathbf{a}:=\ker(S_\infty\to\Lambda)$. Let $S_\infty\supset
I_1\supset I_2\supset\dots$ be a decreasing sequence of open
ideals of~$S_\infty$ with $\cap_N I_N=0$. For each $N\ge 1$ we set
$S_N=S_\infty/I_N$.

Suppose that we are given the following data.

\begin{itemize}
\item Objects $R^1_\infty$, $R^\chi_\infty$ of
  $\CNL_\Lambda$.\item Objects $R^1$, $R^\chi$ of~$\CNL_\Lambda$, an $R^1$-module $M^1$, and an  $R^\chi$-module
  $M^\chi$, each of which is finite as a $\Lambda$-module. Furthermore
  if~$l_0=0$, then they are both free as $\Lambda$-modules, and if
  $l_0=1$, then they are balanced $\Lambda$-modules.
\item For each integer~$N\ge 1$, finite $S_N$-modules $M^1_N$,
  $M^\chi_N$, which are free if~$l_0=0$ and balanced if~$l_0=1$,
  together with isomorphisms of $S_N$-modules $M^1_N/\mathbf{a}\isoto
  M^1\otimes_{S_\infty}S_N$, $M^\chi_N/\mathbf{a}\isoto
  M^\chi\otimes_{S_\infty}S_N$ \emph{(}where the action of $S_\infty$ on $M^1$, $M^\chi$ is via the
  augmentation $S_\infty\to\Lambda$\emph{)}.
\item For each $N\ge 1$, objects $R_N^1$, $R_N^\chi$ of $\CNL_{S_N}$,
  and maps of $S_N$-algebras $R_N^1\to R^1/I_N$, $R_N^\chi\to
  R^\chi/I_N$ and $R_N^1\to \End_{S_N}(M_N^1)$,
  $R_N^\chi\to\End_{S_N}(M_N^\chi)$, such that the two following
  diagrams commute.\[\xymatrix{R_N^1\ar[r]\ar[d] &\End_{S_N}(M_N^1)\ar[d]\\
      R^1/I_N\ar[r]&\End_{\Lambda/I_N}(M^1\otimes_{S_\infty} S_N)   } \]
\[\xymatrix{R_N^\chi\ar[r]\ar[d] &\End_{S_N}(M_N^\chi)\ar[d]\\
      R^\chi/I_N\ar[r]&\End_{\Lambda/I_N}(M^\chi\otimes _{S_\infty} S_N)   } \]
\item For each $N\ge 1$, surjections of $\Lambda$-algebras $R_\infty^1\onto R^1_N$,
  $R^\chi_\infty\onto R^\chi_N$.
  \end{itemize}
We suppose also that we are given the following compatibilities
between the data indexed by~$1$ and the data indexed by~$\chi$.
  \begin{itemize}
\item isomorphisms of $\Lambda/\lambda$-algebras
  $R^1_\infty/\lambda\cong R^\chi_\infty/\lambda$, $R^1/\lambda\cong
  R^\chi/\lambda$, and $R^1_N/\lambda\cong R^\chi_N/\lambda$, compatible with the surjections $R^1_\infty\onto
  R^1_N$ and $R^\chi_\infty\onto R^\chi_N$.

\item An isomorphism of $R^1/\lambda\cong
  R^\chi/\lambda$-modules $M^1/\lambda\cong M^\chi/\lambda$.
\item For each $N\ge 1$, isomorphisms of $S_N/\lambda$-modules $M^1_N/\lambda\cong
  M^\chi_N/\lambda$, compatible with all actions, and such that the
  following diagram
  commutes, where we write~$J_N$ for the kernel of the composite
  $\Lambda\to S_\infty\to S_\infty/I_N$.\[\xymatrix{M_N^1/(\lambda,\mathbf{a})\ar[r]\ar[d]&M_N^\chi/(\lambda,\mathbf{a})\ar[d]\\
      M^1/(\lambda,J_N)\ar[r]& M^\chi/(\lambda,J_N)}  \]\end{itemize}

Then we can find the following data.
\begin{itemize}
\item Homomorphisms of $\Lambda$-algebras $S_\infty\to R_\infty^1$,
  $S_\infty\to R_\infty^\chi$.
\item Finite $S_\infty$-modules $M_\infty^1$, $M_\infty^\chi$, which
  are free if $l_0=0$ and balanced if $l_0=1$, together with
  isomorphisms $M_\infty^1\otimes_{S_\infty}\Lambda\isoto M^1$, and
  $M_\infty^\chi\otimes_{S_\infty}\Lambda\isoto M^\chi$.
\item Commutative diagrams of
  $S_\infty$-algebras \[\xymatrix{R^1_\infty\ar[d]\ar[r]&\End_{S_\infty}(M_\infty^1)\ar[d]^{-\otimes_{S_\infty}\Lambda}\\
    R^1\ar[r]&\End_\Lambda(M^1)} \] \[\xymatrix{R^\chi_\infty\ar[d]\ar[r]&\End_{S_\infty}(M_\infty^\chi)\ar[d]^{-\otimes_{S_\infty}\Lambda}\\
    R^\chi\ar[r]&\End_\Lambda(M^\chi)} \]
\item An isomorphism $M^1_\infty/\lambda\isoto
  M^\chi_\infty/\lambda$, compatible with the actions of
  $R^1_\infty/\lambda\isoto R^\chi_\infty/\lambda$, such that the
  following diagram commutes. \[\xymatrix{M_\infty^1/(\lambda,\mathbf{a})\ar[r]\ar[d]&M_\infty^\chi/(\lambda,\mathbf{a})\ar[d]\\    M^1/\lambda\ar[r]& M^\chi/\lambda}  \]\
\end{itemize}

\end{prop}

        \subsection{The patching construction}\label{subsec: the
          actual patching construction}We now apply
        Proposition~\ref{prop: abstract patching} to our spaces of
        $p$-adic automorphic forms. We continue to assume that
        Hypothesis~\ref{hyp: hypotheses on rhobar for TW primes}
        holds. 

Enlarging $E$ if necessary, we can and do assume that $E$ contains a
primitive $p$th root of unity, and a primitive $9$th root of unity if~$p=3$. By Hypothesis~\ref{hyp: hypotheses on rhobar for TW primes}~(\ref{hyp: conditions on R second time}),
for each $v\in R$ we
can and do choose a pair of non-trivial characters $\chi_v = (\chi_{v,1},\chi_{v,2})$, with
$\chi_{v,i}: \cO_{F_v}^\times \rightarrow \cO^\times$ which are trivial
modulo $\lambda$, and such that $\chi_{v,1}\ne \chi_{v,2}^{\pm 1}$. We will now apply the
constructions of the previous sections, simultaneously using both this
choice of~$\chi$, and also the choice~$\chi=1$. In the former case we
will label our objects as we did before, and in the latter we will
replace~$\chi$ by~$1$.

Let
  \[
   q = h^1(F_S/F,\ad^0\rhobar(1)), \quad  \quad g = 2q-4[F:\Q]+\#S-1,
  \]
  and set $\Delta_\infty = \Zp^{2q}$. Let
  $S_\infty = \cT\llbracket \Delta_\infty \rrbracket$, where~$\cT$ is
  as in~\S\ref{subsec: Galois
  deformation rings}.  Viewing
  $S_\infty$ as an augmented $\Lambda$-algebra, we let $\mathbf{a}$ denote
  the augmentation ideal.

  For each $N\ge 1$, we fix a choice of Taylor--Wiles datum
  $(Q_N,(\alphabar_{v,1},\ldots,\alphabar_{v,4})_{v\in Q_N})$ as in
  Corollary~\ref{cor: existence of TW primes with explicit numbers}.
  For $N = 0$, we set $Q_0 = \emptyset$.  For each $N\ge 1$, we let
  $\Delta_N = \Delta_{Q_N} = \prod_{v\in Q_N} k(v)^\times(p)^2$ and
  fix a surjection $\Delta_\infty \onto \Delta_N$.  The kernel
  of this surjection is contained in $(p^N\Zp)^{2q}$, since each
  $v\in Q_N$ satisfies $q_v \equiv 1 \bmod p^N$.  We let $\Delta_0$ be
  the trivial group, viewed as a quotient of $\Delta_\infty$. We
  write~$S_N=\cT[\Delta_N]$. 

For each $N\ge 0$, we  set $R_N^{1,I,\alphabeta} = R_{\cS_{1,Q_N}^{I,\alphabeta}}$ and
$R_N^{\chi,I,\alphabeta} =
R_{\cS_{\chi,Q_N}^{I,\alphabeta}}$. Note that $R_0^{1,I,\alphabeta} = R_{\cS_1}^{I,\alphabeta}$ and $R_0^{\chi,I,\alphabeta} = R_{\cS_\chi}^{I,\alphabeta}$. 
Let  $R^{1,I,\alphabeta,\loc} = R_{\cS_1^{I,\alphabeta}}^{S,\loc}$ and
$R^{\chi,I,\alphabeta,\loc} =
R_{\cS_\chi^{I,\alphabeta}}^{S,\loc}$ denote the corresponding completed tensor
product of local deformation
rings, as in~\S\ref{subsec: Galois
  cohomology and presentations}. By definition we have
  \[
   R^{1,I,\alphabeta,\loc} = (\widehat{\otimes}_{v\in I} R_v^{P})
   \widehat{\otimes}(\widehat{\otimes}_{v\in I^c} R_v^{B,\alphabetabar_v}) \widehat{\otimes} (\widehat{\otimes}_{v\in R} R_v^1) \widehat{\otimes}
    R_{v_0}^\square,
  \]
  \[
   R^{\chi,I,\alphabeta,\loc} =  (\widehat{\otimes}_{v\in I} R_v^{P})
   \widehat{\otimes}(\widehat{\otimes}_{v\in I^c} R_v^{B,\alphabetabar_v})  \widehat{\otimes} (\widehat{\otimes}_{v\in R} R_v^{\chi_v}) \widehat{\otimes}
   R_{v_0}^\square,
  \]
with all completed tensor products being taken over $\cO$.

For any $N\ge 1$, we have $R_{\cS_{1,Q_N}^{I,\alphabeta}}^{S,\loc} = R^{1,I,\alphabeta,\loc}$ and
$R_{\cS_{\chi,Q_N}^{I,\alphabeta}}^{S,\loc} = R^{\chi,I,\alphabeta,\loc}$.  There are canonical
isomorphisms $R^{1,I,\alphabeta,\loc}/(\lambda) \cong R^{\chi,I,\alphabeta,\loc}/(\lambda)$ and
$R_N^{1,I,\alphabeta}/(\lambda) \cong R_N^{\chi,I,\alphabeta}/(\lambda)$ for all $N\ge 0$.  For each
$N\ge 1$, $R_N^{1,I,\alphabeta}$ and $R_N^{\chi,I,\alphabeta}$ are canonically $\Lambda[\Delta_N]$-algebras
and there are canonical isomorphisms
$R_N^{1,I,\alphabeta} \otimes_{\Lambda[\Delta_N]} \Lambda \cong R_0^{1,I,\alphabeta}$ and
$R_N^{\chi,I,\alphabeta} \otimes_{\Lambda[\Delta_N]} \Lambda \cong R_0^{\chi,I,\alphabeta}$, which are compatible
with the isomorphisms modulo $\lambda$.

Fix representatives $\rho_{\cS_\chi^{I,\alphabeta}}$, $\rho_{\cS_1^{I,\alphabeta}}$ of the universal
deformations which are identified modulo $\lambda$ (via the
identifications
$R_{\cS_\chi^{I,\alphabeta}} / (\lambda) \cong R_{\cS_1^{I,\alphabeta}} / (\lambda)$). By
Lemma~\ref{lem: framed deformation ring over deformation ring}, these
give rise to an $R^{1,I,\alphabeta,\loc}$-algebra structure on
$R_N^{1,I,\alphabeta} \widehat{\otimes}_{\Lambda} \cT$ and an $R^{\chi,I,\alphabeta,\loc}$-algebra structure
on $R_N^{\chi,I,\alphabeta} \widehat{\otimes}_{\Lambda} \cT$; the canonical isomorphism
$R^{1,I,\alphabeta,\loc}/(\lambda) \cong R^{\chi,I,\alphabeta,\loc}/(\lambda)$ is compatible with
these algebra structures and with the canonical isomorphisms
$R_N^{1,I,\alphabeta}/(\lambda) \cong R_N^{\chi,I,\alphabeta}/(\lambda)$.  We let $R_\infty^{1,I,\alphabeta}$ and
$R_\infty^{\chi,I,\alphabeta}$ be formal power series rings in~$g$ variables over
$R^{1,I,\alphabeta,\loc}$ and $R^{\chi,I,\alphabeta,\loc}$, respectively. By Proposition~\ref{prop:
  tangent space of local over global} and Corollary~\ref{cor:
  existence of TW primes with explicit numbers}, we can choose local
$\Lambda$-algebra surjections
$R_\infty^{1,I,\alphabeta} \rightarrow R_N^{1,I,\alphabeta}\widehat{\otimes}_{\Lambda} \cT$ and
$R_\infty^{\chi,I,\alphabeta} \rightarrow R_N^{\chi,I,\alphabeta}\widehat{\otimes}_{\Lambda} \cT$ for every
$N\ge 0$.  We can and do assume that these are compatible with our
fixed identifications modulo $\lambda$, and with the natural
isomorphisms $R_N^{1,I,\alphabeta} \otimes_{\Lambda[\Delta_N]} \Lambda \cong R_0^{1,I,\alphabeta}$ and
$R_N^{\chi,I,\alphabeta}\otimes_{\Lambda[\Delta_N]} \Lambda \cong R_0^{\chi,I,\alphabeta}$.

Fix a subset $I\subset S_p$ of
        cardinality $\#I\le 1$, and a choice of~$\alphabeta$. We now
        apply Proposition~\ref{prop: abstract patching}, taking (in
        the notation established in~\S\ref{subsec: Galois
  deformation problems}):\begin{itemize}
        \item $\Lambda$ to be~$\Lambda_I$.
        \item $S_\infty$, $S_N$ to be as above.\item $R^1_\infty:=R^{1,I,\alphabeta}_\infty$, $R^\chi_\infty:=R^{\chi,I,\alphabeta}_\infty$, $R^1:=R_0^{1,I,\alphabeta}$, $R^\chi:=R_0^{\chi,I,\alphabeta}$, $R_N^1:=R_N^{1,\infty,\alphabeta}\wotimes_\Lambda\cT$,
          $R_N^\chi:=R_N^{\chi,\infty,\alphabeta}\wotimes_\Lambda\cT$.
        \item $M^1:=M^{1,I,\alphabeta}$, $M^\chi:=M^{\chi,I,\alphabeta}$, $M^1_N:=M^{1,I,\alphabeta,Q_N}\wotimes_\Lambda\cT$,
          $M^\chi_N:=M^{\chi,I,\alphabeta,Q_N}\wotimes_\Lambda\cT$.\end{itemize}
By Theorem~\ref{thm: existence of Galois representations in big Hecke
    algebras}, Proposition~\ref{prop: level Q to level 1} and Proposition~\ref{prop: free and balanced over Lambda
    Delta Q}, this data satisfies the assumptions of
  Proposition~\ref{prop: abstract patching}. Consequently, we have:
\begin{itemize}
\item $\Lambda_I$-algebra homomorphisms $S_\infty\to
  R^{1,I,\alphabeta}_\infty$ and $S_\infty\to
  R^{\chi,I,\alphabeta}_\infty$.
\item Finite $S_\infty$-modules $M^{1,I,\alphabeta}_\infty$,
  $M^{\chi,I,\alphabeta}_\infty$ which are free if $\#I=0$ and
  balanced if $\#I=1$, together with isomorphisms
  $M^{1,I,\alphabeta}_\infty/\ba\cong M^{1,I,\alphabeta}$,
  $M^{\chi,I,\alphabeta}/\ba\cong M^{\chi,I,\alphabeta}$.
\item Morphisms of $S_\infty$-algebras
  $R_\infty^{1,I,\alphabeta}\to\End_{S_\infty}(M^{1,I,\alphabeta})$,
  $R_\infty^{\chi,I,\alphabeta}\to\End_{S_\infty}(M^{\chi,I,\alphabeta})$,
  which are compatible with the actions of $R^1,R^\chi$ on
  $M^{1,I,\alphabeta}$, $M^{\chi,I,\alphabeta}$ respectively.

\item Isomorphisms
\begin{equation*}
  M_\infty^{1,I,\alphabeta}/\lambda M_\infty^{1,I,\alphabeta}\simeq M_\infty^{\chi,I,\alphabeta}/\lambda M_\infty^{\chi,I,\alphabeta},\quad M^{1,I,\alphabeta}/\lambda M^{1,I,\alphabeta}\simeq M^{\chi,I,\alphabeta}/\lambda M^{\chi,I,\alphabeta}
\end{equation*}
compatible with the actions of
$R_\infty^{1,I,\alphabeta}/(\lambda)\simeq
R_\infty^{\chi,I,\alphabeta}/(\lambda)$ and
$R^{1,I,\alphabeta}/(\lambda)\simeq R^{\chi,I,\alphabeta}/(\lambda)$
and the above isomorphisms.
\end{itemize}

We now briefly pause to introduce some notation that will be in force
throughout the rest of~\S\ref{sec:CG}. We will need to work
with $\cO$-flat modules~$M$ over complete local Noetherian $\cO$-algebras~$R$
which are not necessarily $\cO$-flat, but for which we have good
control of~$R[1/p]$. There are various ways that we
could do this, but we have found it convenient
to reduce to the $\cO$-flat case in the following way. For a
Noetherian complete local $\cO$-algebra $R$ we denote by $R'$ the
maximal $\cO$-flat quotient of $R$ (i.e.\ the image of $R$ in $R[1/p]$,
or equivalently the quotient of $R$ by its ideal of $p$-power
torsion). Note that if~$M$ is an $R$-module that is $\cO$-flat then it
is naturally an $R'$-module.

        Returning to the situation at hand, by definition, $S_\infty$
        is formally smooth over~$\Lambda_I$ of relative dimension
        $2q+11\#S-1$, and 
        $\Lambda_I$ is formally smooth over~$\cO$.
 By Propositions~\ref{prop: local ordinary deformation ring dimension}, \ref{prop:Ihara1ring},
        \ref{prop:Iharachiring}, and~\ref{prop:smoothlifts},  and ~\cite[Lem.\
        3.3]{blght}, $(R_\infty^{1,I,\alphabeta})'$ and
        $(R^{\chi,I,\alphabeta}_\infty)'$ are equidimensional of relative dimension~$g+10\#S+4[F:\Q]-\#I$ over~$\Lambda_I$. By the
        definition of~$g$, we conclude that
\numequation\label{eqn:dimensions of deformation rings}
     \dim (R_\infty^{1,I,\alphabeta})' = \dim (R_\infty^{\chi,I,\alphabeta})' = \dim S_\infty - \#I.
    \end{equation}

    \begin{prop}\label{prop: patched modules are maximal
        CM}$M^{1,I,\alphabeta}_\infty$ is a maximal Cohen--Macaulay
      $(R^{1,I,\alphabeta}_\infty)'$-module, and $M^{\chi,I,\alphabeta}_\infty$ is a maximal Cohen--Macaulay
      $(R^{\chi,I,\alphabeta}_\infty)'$-module.      
    \end{prop}
    \begin{proof}These statements have identical proofs, so we give
      the argument for the first of them. From~(\ref{eqn:dimensions of deformation rings}), we see that the
support of~$M_\infty^{1,I,\alphabeta}$ in~$\Spec S_\infty$ has codimension at
least~$\#I$. By~\cite[Lem.\ 6.2]{CG} (applied to a
resolution~$S_\infty^r\to S_\infty^r$ of~$M_\infty^{1,I,\alphabeta}$ if~$\#I=1$ ---
such a resolution exists, by Lemma~\ref{lem: presentations of balanced modules}  --- and to~$M_\infty^{1,I,\alphabeta}$ itself
if~$\#I=0$), we see that the codimension is precisely~$\#I$, and that
$M_\infty^1$ has depth~$\dim S_\infty-\#I=\dim (R_\infty^{1,I,\alphabeta})'$ over~$S_\infty$. It follows
that the depth of~$M_\infty^{1,I,\alphabeta}$ over~$(R_\infty^{1,I,\alphabeta})'$ is at least~$\dim
(R_\infty^{1,I,\alphabeta})'$, so that $M_\infty^{1,I,\alphabeta}$ is maximal Cohen--Macaulay
over~$(R_\infty^{1,I,\alphabeta})'$, as required. \end{proof}

\subsection{Cycles and modules over products of local deformation
  rings}\label{subsec: cycles and products}
In preparation for our study of the dimensions of certain spaces of
$p$-adic modular forms in the next section, we formalize some
arguments which are at the heart of our version of the ``Ihara
avoidance'' argument
of~\cite{tay}. Following~\cite{emertongeerefinedBM}, we use the
language of cycles on the special fibres of (completed tensor products
of) local deformation rings; our perspective is also informed
by~\cite{shottonGLn}.

We recall some notation for cycles and multiplicities from ~\cite[\S 2]{emertongeerefinedBM}. In
particular, if~$R$ is an equidimensional Noetherian local ring of dimension~$d$ then by a cycle (or a $d$-cycle) on $\Spec R$ we mean simply a formal $\mathbb{Z}$-linear combination of the generic points of $\Spec R$.  We denote the group of cycles on $R$ by $\cZ^d(R)$ (or just $\cZ(R)$, with  the understanding that we will only consider top-dimensional cycles).  If $M$ is a finite $R$-module then the cycle of $M$ is defined
by \[Z(M,R)=\sum_{\eta}\len_{R_\eta}(M_\eta)\cdot\eta\] where the sum
is over the generic points~$\eta$ of~$\Spec R$ and $\len_{R_\eta}(M_\eta)$ denotes the length of $M_\eta$ as a $R_\eta$-module.

If $R$ is an equidimensional, flat, Noetherian $\cO$-algebra of dimension $d+1$ and $\eta$ is a generic point of $\Spec R$ then we write $R^\eta$ for the quotient of $R$ by the minimal prime corresponding to $\eta$, and we let $\overline{\eta}=Z(R^\eta/(\lambda),R/(\lambda))$.  Then \cite[Prop.\ 2.2.13]{emertongeerefinedBM} states that if $M$ is a finite $R$-module which is $\cO$-flat, then
\begin{equation*}
Z(M/\lambda M,R/(\lambda))=\sum_\eta \len_{R_\eta}(M_\eta)\cdot\overline{\eta}
\end{equation*}
where the sum is over the generic points $\eta$ of $R$.

Next we recall several facts about completed tensor products. As in~\S\ref{subsec: the
          actual patching construction}, if $R\in\CNL_{\cO}$, we
        let~$R'$ denote the maximal $p$-torsion free quotient of~$R$. Let $R_1,R_2\in\CNL_\cO$.  First we note that the natural map $R_1\widehat{\otimes} R_2\to
R_1'\widehat{\otimes} R_2'$ induces an isomorphism
$(R_1\widehat{\otimes} R_2)'\simeq R_1'\widehat{\otimes} R_2'$.
(Indeed this follows from the fact that the kernel is
$p$-power torsion and that $R_1'\widehat{\otimes}R_2'$ is $\cO$-flat, see \cite[Lem.\ 1.3]{MR3327536}.)

 Now suppose that $R_1$ and $R_2$ are $\cO$-flat and equidimensional
 of dimensions $d_1+1$ and $d_2+1$ respectively, and further assume that all the
 irreducible components of $\Spec R_i$ and $\Spec R_i/(\lambda)$ for $i=1,2$ are
 geometrically irreducible (for instance by enlarging $\cO$ if
 necessary). Write $R=R_1\widehat{\otimes}R_2$; then  $R$ is $\cO$-flat and
equidimensional of dimension $d_1+d_2+1$. (This, and the other facts
recalled in this paragraph, can be read
off from \cite[Lem 1.4]{MR3327536}.)  Moreover if $\eta_i$ is a generic point of $\Spec R_i$ for $i=1,2$ then the kernel of the natural map
\begin{equation*}
R\to R_1^{\eta_1}\widehat{\otimes}R_2^{\eta_2}
\end{equation*}
is a minimal prime of $R$ which corresponds to a generic point of
$\Spec R$ which we denote by $\eta=(\eta_1,\eta_2)$, and the generic
points of $\Spec R$ are precisely the $(\eta_1,\eta_2)$ as~ $\eta_i$
ranges over the
generic points of $\Spec R_i$ for $i=1,2$.
Similarly if $\p_i\subset R_i/(\lambda)$ is a minimal prime for
$i=1,2$ then \[(\p_1,\p_2)=\ker\left( R/(\lambda)\to R_1/(\lambda,\p_1)\widehat{\otimes}R_2/(\lambda,\p_2)\right)\] is a minimal prime of $R/(\lambda)$, and every minimal prime of $R/(\lambda)$ has this form.  It follows that there is an isomorphism
\begin{align*}
\cZ^{d_1}(R_1/(\lambda))\otimes\cZ^{d_2}(R_2/(\lambda))&\to\cZ^{d_1+d_2}(R/(\lambda)),\\
\eta_1\otimes\eta_2 &\mapsto (\eta_1,\eta_2).
\end{align*}

According to \cite[Lem.\ 2.2.14]{emertongeerefinedBM}, if $M_i$ is a
finite $R_i$-module for $i=1,2$, so that we may form the $R$-module
$M=M_1\widehat{\otimes}_\cO M_2$, then under the above isomorphism we have
\numequation\label{eqn:tensor product of cycles}
Z(M_1/\lambda M_1,R_1/(\lambda))\otimes Z(M_2/\lambda M_2,R_2/(\lambda))=Z(M/\lambda M,R/(\lambda)).
\end{equation}
In particular for a generic point $\eta=(\eta_1,\eta_2)$ of $R$ we have an isomorphism $R^\eta\simeq R^{\eta_1}\widehat{\otimes}R^{\eta_2}$ of $R$-modules and hence, in the notation introduced above, under this isomorphism we have $\overline{\eta}=\overline{\eta}_1\otimes\overline{\eta}_2$.

We wish to apply this discussion to the rings
\begin{equation*}
R^1=\widehat{\otimes}_{v\in R} R_v^1,\quad R^\chi=\widehat{\otimes}_{v\in R}R_v^\chi
\end{equation*}
as well as to $\widetilde{R}^1=R^1\widehat{\otimes}\widetilde{R}$ and
$\widetilde{R}^\chi=R^\chi\widehat{\otimes}\widetilde{R}$, for some
auxiliary  $\widetilde{R}\in\CNL_{\cO}$ with the property that
$\widetilde{R}'$ is irreducible. (In applications $\widetilde{R}^1$
and $\widetilde{R}^\chi$ will be $R_\infty^{1,I,\alphabeta}$ and
$R_\infty^{\chi,I,\alphabeta}$ for some choice of $I$ and
$\alphabeta$; so~$\widetilde{R}$ is formally smooth over a completed
tensor product of the deformation rings considered in Proposition~\ref{prop:
  local ordinary deformation ring dimension}, and $\widetilde{R}'$ is
indeed irreducible.)

We recall that for each $v\in R$ we have $R_v^1/(\lambda)=R_v^\chi/(\lambda)$.  Passing to $p$-torsion free quotients, it is not the case that $(R_v^1)'/(\lambda)$ is identified with $(R_v^\chi)'/(\lambda)$, but Propositions \ref{prop:Ihara1ring} and \ref{prop:Iharachiring} imply that at least the underlying topological spaces of $\Spec (R_v^1)'/(\lambda)$ and $\Spec (R_v^\chi)'/(\lambda)$ coincide with that of $\Spec R_v^1/(\lambda)=\Spec R_v^\chi/(\lambda)$, and so in particular $\cZ((R_v^1)'/(\lambda))=\cZ((R_v^\chi)'/(\lambda))$.  Passing to products we obtain identifications $\cZ((R^1)'/(\lambda))=\cZ((R^\chi)'/(\lambda))$ and $\cZ((\widetilde{R}^1)'/(\lambda))=\cZ((\widetilde{R}^\chi)'/(\lambda))$.

\begin{lem}\label{lem: stupid p torsion free quotient cycle comparison lemma}
Let $M^1$ be a finite $\cO$-flat $\widetilde{R}^1$-module, and let
$M^\chi$ be a finite $\cO$-flat $\widetilde{R}^\chi$-module, such that $M^1/\lambda M^1\simeq M^\chi/\lambda M^\chi$ as $\widetilde{R}^1/(\lambda)=\widetilde{R}^\chi/(\lambda)$-modules.  Then
\begin{equation*}
Z(M^1/\lambda M^1,(\widetilde{R}^1)'/(\lambda))=Z(M^\chi/\lambda M^\chi,(\widetilde{R}^\chi)'/(\lambda))
\end{equation*}
under the identification of $\cZ((\widetilde{R}^1)'/(\lambda))$ with $\cZ((\widetilde{R}^\chi)'/(\lambda))$ from above.
\end{lem}
\begin{proof}
As we explained above we have two quotients
$(\widetilde{R}^1)'/(\lambda)$ and $(\widetilde{R}^\chi)'/(\lambda)$
of $\widetilde{R}^1/(\lambda)=\widetilde{R}^\chi/(\lambda)$ whose
spectra have the same underlying topological space.  Each generic
point of this space corresponds to minimal primes $\p^1$ and $\p^\chi$
of $(\widetilde{R}^1)'/(\lambda)$ and
$(\widetilde{R}^\chi)'/(\lambda)$ as well as to a (not necessarily
minimal) prime $\p$ of
$\widetilde{R}^1/(\lambda)=\widetilde{R}^\chi/(\lambda)$ which is the
preimage of both $\p^1$ and $\p^\chi$.  Then we claim that we have equalities
\begin{equation*}
\len_{((\widetilde{R}^1)'/(\lambda))_{\p^1}}((M^1/\lambda M^1)_{\p^1})=\len_{(\widetilde{R}^1/(\lambda))_\p}((M^1/\lambda M^1)_\p)=\len_{((\widetilde{R}^\chi)'/(\lambda))_{\p^\chi}}((M^\chi/\lambda M^\chi)_{\p^\chi})
\end{equation*}
which exactly gives the statement of the lemma.  Both equalities follow from the fact that if $A\to B$ is a surjective map of rings and $M$ is a finite length $B$-module then $\len_A(M)=\len_B(M)$.
\end{proof}

For the next lemma we need to introduce some more notation.  From a
tuple $\eta=(\eta_v)_{v\in R}$ of generic points $\eta_v$ of $\Spec
(R_v^1)'$ for $v\in R$ we obtain a generic point $\eta$ of $\Spec
(R^1)'$ (resp.\ a generic point also denoted $\eta$ of $\Spec
(\widetilde{R}^1)'$) and moreover these are all of the generic points
of $\Spec (R^1)'$ (resp.\ of $\Spec (\widetilde{R}^1)'$).  By
Proposition \ref{prop:Ihara1ring}, if $\eta_{v,1}$ and $\eta_{v,2}$
are two distinct generic points of $\Spec (R_v^1)'$ for some $v\in R$,
then the cycles $\overline{\eta}_{v,1}$ and $\overline{\eta}_{v,2}$
have disjoint support.

It follows from this and~(\ref{eqn:tensor product of cycles}) that if
$\eta_1$ and $\eta_2$ are two distinct generic points of $\Spec
(R^1)'$ (resp.\ of $\Spec (\widetilde{R}^1)'$) then the supports of
$\overline{\eta}_1$ and $\overline{\eta}_2$ are disjoint.  Finally
recall that by Proposition \ref{prop:Iharachiring}, for each $v\in R$,
$\Spec (R_v^{\chi})'$ is irreducible.  Passing to products, $\Spec
(R^\chi)'$ and $\Spec (\widetilde{R}^\chi)'$ are irreducible as well.
We denote the unique generic point of either by $\eta^\chi$. 

As already indicated, in the statement and proof of the following
lemma, we freely identify the generic points of~$\Spec
(R^1)'$ and~$\Spec
(\widetilde{R}^1)'$ (and we also identify the generic points of ~$\Spec
(R^1)'/(\lambda)$ and~$\Spec
(\widetilde{R}^1)'/(\lambda)$).
\begin{lem}\label{lem: main cycle comparison lemma}
Suppose there exists a finite, $\cO$-flat $\widetilde{R}^1$-module $M^1$, and a finite, $\cO$-flat $\widetilde{R}^\chi$-module $M^\chi$, along with an isomorphism $M^1/\lambda M^1\simeq M^\chi/\lambda M^\chi$ of $\widetilde{R}^1/(\lambda)=\widetilde{R}^\chi/(\lambda)$-modules.  Suppose furthermore that $M^1$ is supported on at least one generic point $\eta$ of $\Spec(\widetilde{R}^1)'$.  Then there exist unique positive integers $d_{\eta}'$ labelled by generic points $\eta=(\eta_v)_{v\in R}$ of $\Spec (R^1)'$ such that
\begin{enumerate}
\item As elements of $\cZ((R^1)'/(\lambda))=\cZ((R^\chi)'/(\lambda))$ we have
\numequation\label{eqn: cycle formula}
\overline{\eta^\chi}=\sum_{\eta}d_\eta'\overline{\eta}
\end{equation}
where the sum is over the generic points $\eta$ of $\Spec (R^1)'$.
\item For each generic point $\eta$ of $\Spec (\widetilde{R}^1)'$ we have
\numequation\label{eqn: length formula}
\len_{(\widetilde{R}^1)'_\eta}(M^1_\eta)=d_\eta'\len_{(\widetilde{R}^\chi)'_{\eta^\chi}}(M^\chi_\eta).
\end{equation}
In particular $M^1$ is supported on every generic point of $\Spec(\widetilde{R}^1)'$.
\end{enumerate}
\end{lem}\begin{proof}As explained above, the cycles $\overline{\eta}$ as
  $\eta$ ranges over the generic points of $\Spec(R^1)'$ have disjoint
  support.  Thus the formula~(\ref{eqn: cycle formula}) uniquely
  determines the integers $d'_\eta$.  Moreover, as the cycle
  $\overline{\eta^\chi}$ is supported on every generic point of $\Spec
  (R^\chi)'/(\lambda)$, (\ref{eqn: cycle formula}) also implies that the integers $d_\eta$ must be positive, if they exist.

Now using \cite[Prop.\ 2.2.13]{emertongeerefinedBM} as recalled above, we have
\begin{equation*}
Z(M^1/\lambda M^1,(\widetilde{R}^1)'/(\lambda))=\sum_\eta \len_{(\widetilde{R}^1)'_\eta}M^1_\eta\cdot\overline{\eta}
\end{equation*}
and
\begin{equation*}
Z(M^\chi/\lambda M^\chi,(\widetilde{R}^\chi)'/(\lambda))=\len_{(\widetilde{R}^\chi)'_{\eta^\chi}}M^\chi_{\eta^\chi}\cdot\overline{\eta^\chi}
\end{equation*}
and moreover these two cycles coincide by Lemma \ref{lem: stupid p torsion free quotient cycle comparison lemma}.

Our hypothesis that $M^1$ is supported on some generic point $\eta$ of $(\widetilde{R}^1)'$ implies that $\len_{(\widetilde{R}^1)'_\eta}M^1_\eta>0$.  Hence by the above equality of cycles, $\len_{(\widetilde{R}^\chi)'_{\eta^\chi}}M^\chi_{\eta^\chi}>0$.  Because the cycles $\overline{\eta}$ have disjoint support, we must have that
\begin{equation*}
d_\eta'=\frac{\len_{(\widetilde{R}^1)'_\eta}(M^1_\eta)}{\len_{(\widetilde{R}^\chi)'_{\eta^\chi}}(M^\chi_{\eta^\chi})}
\end{equation*}
is an integer for each generic point $\eta$ of $\Spec(R^1)'$, and for this choice of $d_\eta'$, the formulas~(\ref{eqn: cycle formula}) and~(\ref{eqn: length formula}) hold.
\end{proof}

\begin{rem}\label{rem: multiplicities are independent of global situation}
We note that the ``multiplicities'' $d_\eta'$ in  Lemma~\ref{lem: main cycle comparison lemma} are
independent of the modules $M^1$ and $M^\chi$ and even of the
auxiliary ring $\widetilde{R}$.  In the next section they will be given a
local representation-theoretic interpretation (see
Proposition~\ref{prop: size of L packets for Ihara avoidance} and Remark~\ref{rem: remark about geometric BM for l not p}). \end{rem}

\begin{rem}\label{rem: Ihara avoidance goes back to RLT}In
 ~\S\ref{subsec: patching padic multiplicities} we will use Lemma~\ref{lem: main cycle comparison lemma} to compute the dimensions of
 spaces of $p$-adic modular forms at Iwahori level.
 The idea of comparing patched modules over ~$R_\infty^{1,I,\alphabeta}$ and~$R_\infty^{\chi,I,\alphabeta}$  goes back to~\cite{tay}; the key point is that $R_\infty^{\chi,I,\alphabeta}[1/p]$ is a domain, which guarantees that the
  support of an appropriate patched module is all
  of~$\Spec R_\infty^{\chi,I,\alphabeta}$, and the isomorphism
  $R_\infty^{1,I,\alphabeta}/(\lambda)=R_\infty^{\chi,I,\alphabeta}/(\lambda)$ which allows us to transfer this information
  to~$R_\infty^{1,I,\alphabeta}$. \end{rem}

\subsection{Multiplicities of patched spaces of \texorpdfstring{$p$}{p}-adic automorphic
  forms}\label{subsec: patching padic multiplicities}

We now make use of our patching constructions to determine the
multiplicities of systems of eigenvalues corresponding to~$\rho$ in
spaces of $p$-adic automorphic forms with~$\#I\le 1$. 

 We begin by introducing some
notation and assumptions. We suppose that we have fixed a representation
$\rho:G_F\to\GSp_4(\cO)$, which satisfies the following
properties. (While this list of properties may appear to be too restrictive to
be useful, we will later use base change to reduce to this situation.)
\begin{hypothesis}\label{hyp: hypotheses on rho}
\leavevmode
  \begin{enumerate}[series=hypothesisconditions]
 \item $F$ is a totally real field in which the prime $p\ge 3$ splits
    completely; we write~$S_p$ for the set of primes of~$F$ dividing~$p$.
  \item  $\nu\circ\rho=\varepsilon^{-1}$.
  \item  For each finite place~$v$ of~$F$, $\rho|_{G_{F_v}}$ is pure.

  \item For each $v\in S_p$, $\rho|_{G_{F_v}}$ is $p$-distinguished weight~$2$ ordinary,
    with unit eigenvalues $\alpha_v,\beta_v\in E$.
  \item There is a finite set $R$ of primes
    of $F$ not dividing $p$ such that if $v\not\in R\cup S_p$, then
    $\rho|_{G_{F_v}}$ is unramified, while if $v\in R$, then:
    \begin{itemize}
    \item $q_v\equiv 1\pmod p$, and if $p=3$ then further $q_v\equiv 1\pmod 9$.
    \item  $\overline{\rho}|_{G_{F_v}}$ is
      trivial.
       \item $\rho|_{G_{F_v}}$ has only unipotent ramification.
\end{itemize}
    \item There exists $\pi=\otimes_v\pi_v$ an ordinary cuspidal automorphic
    representation for $\GSp_4/F$ of parallel weight~$2$ and central
    character $|\cdot|^2$, such that
    $\overline{\rho}_{\pi,p}=\overline{\rho}$, and such that:
    \begin{itemize}
    \item For all $v\not\in  R\cup S_p$, $\pi_v$ is
      unramified.\item For all $v\in R\cup S_p$, $\pi_v^{\Iw(v)}\not=0$.
      \item For each finite place~$v$ of~$F$, $\rho_{\pi,p}|_{G_{F_v}}$ is pure.      \end{itemize}
  \item The representation~$\rhobar$ is vast and tidy.
  \end{enumerate}
\end{hypothesis}
\begin{rem}
  \label{rem: rho hypothesis implies rhobar hypothesis}Note in
  particular that Hypothesis~\ref{hyp: hypotheses on rho} implies that
  Hypothesis~\ref{hyp: hypotheses on rhobar for TW primes} holds for~$\rhobar$.
\end{rem}

As in~\S\ref{subsec: Galois
  deformation problems}, it follows from Hypothesis~\ref{hyp:
  hypotheses on rho}, and in particular from the hypothesis
that~$\rhobar(G_F)$ is tidy, that:
\begin{enumerate}[resume=hypothesisconditions]
\item There exists an absolutely unramified prime
  $v_0\not\in S_p\cup R$ with $q_{v_0}\not\equiv 1\pmod p$ and residue
  characteristic greater than~$5$, such that
  $\rho|_{G_{F_{v_0}}}$ is unramified, and $\rho(\Frob_{v_0})$ has
  (not necessarily distinct)
  eigenvalues with the property that no ratio of these eigenvalues is congruent to
  $q_{v_0}$ modulo~$\lambda$.
\end{enumerate}

Given a closed point $x\in\Spec R^{1,I,\alphabeta}[1/p]$ or $x\in\Spec R^{\chi,I,\alphabeta}[1/p]$, we will always assume that $E$ is
large enough to contain the residue field of $x$, so that in particular $x$
parameterizes a Galois representation $\rho_x:G_F\to \GSp_4(\cO)$.  We
denote by $\p_{x}$ the height one prime ideal which is the kernel
of the corresponding homomorphism $R^{1,I,\alphabeta}\to E$ or $R^{\chi,I,\alphabeta}\to E$, and we
also use the same symbol $\p_{x}$ for the ideals obtained by
pulling back under the homomorphisms
$R^{1,I,\alphabeta,\loc}\to R_\infty^{1,I,\alphabeta}\to
R^{1,I,\alphabeta}$ or under the homomorphisms
$R^{\chi,I,\alphabeta,\loc}\to R_\infty^{\chi,I,\alphabeta}\to
R^{\chi,I,\alphabeta}$.  As in~\S\ref{subsec: l not
  p deformations}, we say that $x$ (or $\rho_x$ or $\p_{x}$)
is \emph{smooth} if $R^{1,I,\alphabeta,\loc}_{\p_{x}}$ (resp.\ $R^{\chi,I,\alphabeta,\loc}_{\p_{x}}$) is a
regular local ring (note that this is equivalent to their completions
being regular). 

For a Galois representation
$\rho':G_F\to\GSp_4(\cO)$ giving rise to a point~$x'$ on one of the
deformation rings $R^{1,I,\alphabeta}$ or $R^{\chi,I,\alphabeta}$, we let
$\tilde{\p}^{\mathrm{an}}_{x'}\subset\TTt$ be the corresponding
prime ideal. Explicitly, this is the prime ideal generated
by 
the coefficients of the polynomials
$Q_v(X)-\det(X-\rho'(\Frob_v))$ for~$v\notin R\cup
S_p\cup\{v_0\}$. (As before, the ``an'' stands for
``anaemic''.)

For any
choice of $I$ and~$\alphabeta$, we 
let $\tilde{\p}_{x'}^{I,\alphabeta}\subset\TTt^I$ denote the prime ideal
\numequation\label{eqn-def-primeideal} 
\tilde{\p}_{x'}^{I,\alphabeta}=(\tilde{\p}_{x'}^{\mathrm{an}},\{U_{v,0}-1,U_{v,2}-\alpha_v'\beta_v'\}_{v\in
S_p},\{U_{\Kli(v),1}-\alpha_v'-\beta_v'\}_{v\in
I},\{U_{v,1}-\alphabeta_v'\}_{v\in
I^c}), \end{equation}where, for $v\in S_p$, $\alpha_v'\equiv \alphabar_v\pmod
\lambda$, $\beta_v'\equiv\betabar_v\pmod \lambda$ and
$\alphabeta_v'\in\{\alpha_v',\beta_v'\}$ are determined by the
local representations $\rho'|_{G_{F_v}}$ as in~\S\ref{subsec:Ordinary
  lifts}.

\begin{defn}\label{defn: Iwahori and optimal level}
  Let $K^{p,\Iw}=\prod_{v\nmid p,\infty}K_v$ and
  $K^{p,\Iw_1}=\prod_{v\nmid p,\infty}K_v'$, where
  \begin{itemize}
  \item $K_v=\Iw(v)$ and $K_v'=\Iw_1(v)$ for $v\in R$.
  \item $K_{v_0}=K_{v_0}'=\Iw_1(v_0)$.
  \item $K_v=K_v'=\GSp_4(\cO_{F_v})$ for $v\not\in R\cup\{v_0\}$.
  \end{itemize}
\end{defn}

We now define some spaces of $p$-adic modular forms.  For any
$I\subset S_p$,  $\alphabeta$, classical algebraic weight~$\kappa$, and choice of~$K^{p, \Iw}$ as in
Definition~\ref{defn: Iwahori and optimal level}, we let \begin{multline*}
S_{\kappa,K^{p, \Iw}K_p(I)}^{I,\alphabeta}=(e(U^I)H^0(\mathfrak{X}^{I,G_1}_{K_p(I) K^p},\omega^{\kappa}(-D))_{\tilde{\m}^{I,\alphabeta}})\otimes_\cO
E[\{U_{v,0}-1\}_{v\in S_p},\\ \{U_{v,0}-q_v^{-2}\}_{v\in
R},\{T_{v,0}-q_v^{-2}\}_{v\not\in S_p\cup R\cup \{v_0\}}]
\end{multline*} We also let
\begin{equation*}
S^{I,\alphabeta}_{\kappa,K^{p,\Iw_1}K_p(I),\chi}:=
(S^{I,\alphabeta}_{\kappa,K^{p,\Iw_1}K_p(I)})^{\prod_{v\in R}\Iw(v)/\Iw_1(v)=\chi}
\end{equation*} be the subspace with ``nebentypus'' corresponding to
$\chi$.  By Lemma~\ref{lem: justification of RHom}~(2), we have isomorphisms

\[(M^{1,I,\alphabeta}/\p_{\kappa}M^{1,I,\alphabeta})[1/p]\simeq
(S^{I,\alphabeta}_{\kappa,K^{p,\Iw}K_p(I)})^\vee,\]\[
(M^{\chi,I,\alphabeta}/\p_{\kappa}M^{\chi,I,\alphabeta})[1/p]\simeq
(S^{I,\alphabeta}_{\kappa,K^{p,\Iw_1}K_p(I),\chi})^\vee.
\]

In particular, with this notation in place, for a Galois
representation $\rho':G_F\to\GSp_4(\cO)$ giving rise to a point on one
of the deformation rings $R^{1,I,\alphabeta}$ or $R^{\chi,I,\alphabeta}$
and of weight $\kappa$ (i.e.\ such that the composition $\Lambda_I\to
R^{1,I,\alphabeta}\to E$ or $\Lambda_I\to R^{\chi,I,\alphabeta}\to E$ is
$\kappa$) we have
\begin{equation*}
(M_\infty^{1,I,\alphabeta}/\p_{x'}M_\infty^{1,I,\alphabeta})[1/p]\simeq
(M^{1,I,\alphabeta}/\p_{x'}M^{1,I,\alphabeta})[1/p]\simeq
(S^{I,\alphabeta}_{\kappa,K^{p,\Iw}K_p(I)}[\tilde{\p}^{I,\alphabeta}_{x'}])^\vee,
\end{equation*}
\begin{equation*}
(M_\infty^{\chi,I,\alphabeta}/\p_{x'}M_\infty^{\chi,I,\alphabeta})[1/p]\simeq
(M^{1,I,\alphabeta}/\p_{x'}M^{1,I,\alphabeta})[1/p]\simeq
(S^{I,\alphabeta}_{\kappa,K^{p,\Iw_1}K_p(I),\chi}[\tilde{\p}^{I,\alphabeta}_{x'}])^\vee.
\end{equation*}

In order to state our results on the dimensions of eigenspaces of
$p$-adic automorphic forms, we need to make a further study of the
local deformation rings at places~$v\in R$. 

\begin{prop}\label{prop: size of L packets for Ihara avoidance}
  Let $\eta_v$ be a generic point of $\Spec R_v^1$ for some $v\in R$.  The set of $y\in (\Spec R_v^{1,\eta_v})(\overline{E})$ such that the
  $L$-packet $L(\rho_y)$ contains a generic representation is
  nonempty, and the number
\begin{equation*}
d_{\eta_v}=\sum_{\pi\in L(\rho_y)}\dim\pi^{\Iw(v)}
\end{equation*}
is independent of such a $y$. More explicitly, the
rank~$n(\eta_v)$ of the monodromy operator~$N$ is generically constant on~$\Spec
R_v^{1,\eta_v}$, and
\begin{itemize}
\item if~$n(\eta_v)=0$, then $d_\eta=8$;
\item if~$n(\eta_v)=1$, then $d_\eta=4$;
\item if~$n(\eta_v)=2$, then $d_\eta=4$;
\item if~$n(\eta_v)=3$, then $d_\eta=1$.
\end{itemize}
\end{prop}\begin{proof}
  This can be read off from~\cite[Tables A.7, A.15]{MR2344630} (note
  that the rank of the monodromy operator is given in the column
  of~\cite[Table A.15]{MR2344630} headed ``$a$''; note also that the unipotent
$L$-packets which contain supercuspidal representations
also contain generic non-supercuspidal representations, namely those
of type Va and XIa, see~\cite[\S1]{schmidt-tables}, and~\cite[Table A.1]{MR2344630}). We see
  that:

\begin{itemize}
\item On the unramified components (those with $n(\eta_v)=0$), the
  $L$-packets containing a generic representation are singletons~$\{\pi\}$ of
  type~I (unramified principal series), so $d_{\eta_v}=8$.\item If $n(\eta_v)=1$, the  $L$-packets containing a generic
  representation are singletons~$\{\pi\}$ of type~ IIa, so $d_{\eta_v}=4$.\item If $n(\eta_v)=2$, the $L$-packets containing a generic
  representation are either singletons~$\{\pi\}$ of type~IIIa, or
  pairs~$\{\pi_a,\pi_b\}$ of respective types VIa
  and VIb, and in either case $d_{\eta_v}=4$. (Note that the
  representations of type Va do not contribute, as they never
  correspond to residually trivial Galois representations.)
\item Finally, if $n(\eta_v)=3$, then the $L$-packets containing a
  generic representation are singletons of type~IVa (unramified twists
  of Steinberg) and $d_{\eta_v}=1$. \qedhere
\end{itemize}
\end{proof}

 We write $\eta=(\eta_v)_{v\in R}$ for a tuple of generic points
 $\eta_v$ of $\Spec R^1_v$ for $v\in R$, which as explained in~\S\ref{subsec: cycles and products} gives rise to a generic point, also denoted $\eta$, of $R^{1,I,\alphabeta,\loc}$ or of $R_\infty^{1,I,\alphabeta}$.  We let
\begin{equation*}
d_\eta=\prod_{v\in R} d_{\eta_v}
\end{equation*}where~$d_{\eta_v}$ is as in Proposition~\ref{prop: size of L packets for Ihara avoidance}.
We also let $d_\rho=d_\eta$ for the generic point $\eta$ of $R^{1,I,\alphabeta,\loc}$
that the local representations of $\rho$ lie on (this point is unique,
as the representations~$\rho|_{G_{F_{v}}}$ are pure by assumption).  Concretely, by Proposition~\ref{prop: size of L packets for Ihara avoidance}, we have
\begin{equation*}
d_\rho=8^{|R_0|}4^{|R_1|+|R_2|}
\end{equation*}
where for $i=0,1,2,3$, $R_i\subset R$ is the set of primes $v\in R$ for which $n(\rho|_{G_{F_v}})=i$.

We can now state our main result about $p$-adic modularity at
Iwahori level. 

\begin{thm}\label{thm: iwahori patching main section}Assume
  Hypothesis~\ref{hyp: hypotheses on rho} for~$\rho = \rho_x$. For
  any~$I$ with $\# I\le 1$, and any
  choice of $\alphabeta$, we have
\begin{equation*} \dim_E
S_{\kappa_2,K^{p,\Iw}K_p(I)}^{I,\alphabeta}[\tilde{\p}_{x}^{I,\alphabeta}]=8d_\rho.
\end{equation*}
\end{thm}
\begin{rem}
  The reason for the factor of~$8=|W|$ on the right hand side is that we
  are working at Iwahori level at the auxiliary place~$v_0$, and not imposing any conditions on the Hecke operators at this place.  It would be possible to impose such
  conditions and remove this factor, but we have found it more
  convenient not to do so (and it makes no difference for our main
  automorphy lifting theorems).
\end{rem}

Before proving the theorem we recall a standard lemma, essentially due
independently to Diamond and Fujiwara (see e.g.\ \cite{MR1440309}) which is the key to proving ``multiplicity one'' (or
``multiplicity~$8d_\rho$'') results in characteristic 0 using the
Taylor--Wiles method.
\begin{lem}\label{lem: smooth implies free main section} Let $R$ be either
$(R_\infty^{1,I,\alphabeta})'$ or $(R_\infty^{\chi,I,\alphabeta})'$ for some
choice of $I$ and $\alphabeta$ and let $M$ be a maximal Cohen--Macaulay
$R$-module.  Let $x\in\Spec R[1/p]$ be a smooth closed point with
residue field $E$, and let $\p_x\subset R$ be the corresponding prime
ideal.  Then $M_{\p_x}$ is a free $R_{\p_x}$-module, and hence if $\eta$ is
the unique generic point of $R$ specializing to $x$, then
\begin{equation*} \dim_{R_\eta} M_\eta=\dim_E (M/\p_{x} M)[1/p].
\end{equation*}
\end{lem}
\begin{proof} The first statement follows from the fact that a maximal
Cohen--Macaulay module over a regular local ring is free, and the
second statement is an immediate consequence of this freeness.
\end{proof}

We also record the following proposition on ``doubling'':
\begin{prop}\label{prop: doubling main section}
Let~$\rho = \rho_x$. For any choice of
$I\subset S_p$, $w\in I^c$, $K^p$ as in Definition~\ref{defn: Iwahori
  and optimal level}, and $\alphabeta$, there is an injection
\begin{equation*}
(U_{w,1}-\alphabeta_w'):S^{I\cup\{w\},\alphabeta}_{\kappa_2,K^{p}K_p(I\cup\{w\})}[\tilde{\p}_{x}^{I\cup\{w\},\alphabeta}]\to
S^{I,\alphabeta}_{\kappa_2,K^{p}K_p(I)}[\tilde{\p}_{x}^{I,\alphabeta}].
\end{equation*}
\end{prop}
\begin{proof}This immediately reduces to the corresponding statement
  with~$\cO$-coefficients, and hence to the injectivity of the map~(\ref{eqn:
    doubling injection special fibre}) (with $w=v$), which we proved in the course
  of the proof of Proposition~\ref{prop: doubling for $I$ to $I'$.}.
  \end{proof}
We are now ready to prove Theorem~\ref{thm: iwahori patching main section}.

\begin{proof}[Proof of Theorem \ref{thm: iwahori patching main section}]  We first
consider the case that $I=\emptyset$.  As $M^{1,\emptyset,\alphabeta}$ is a maximal Cohen--Macaulay $(R^{1,\emptyset,\alphabeta})'$-module, it is supported on some irreducible component of $\Spec (R^{1,\emptyset,\alphabeta})'$ and hence we may apply Lemma~\ref{lem: main cycle comparison lemma}
to the $R_\infty^{1,\emptyset,\alphabeta}$-module
$M^{1,\emptyset,\alphabeta}_\infty$ and the
$R_\infty^{\chi,\emptyset,\alphabeta}$-module
$M^{\chi,\emptyset,\alphabeta}_\infty$.  In particular we conclude that
$M_\infty^{1,\emptyset,\alphabeta}$ is supported on every irreducible
component of $\Spec (R^{1,\emptyset,\alphabeta})'$.

As $M_\infty^{1,\emptyset,\alphabeta}$
is a finite free
$S_\infty$-module,  and $(R_\infty^{1,\emptyset,\alphabeta})^{\red}$
  acts faithfully on $M_\infty^{1,\emptyset,\alphabeta}$ (and is
  therefore finite and torsion free over $S_\infty$), the map $\Spec 
  (R_\infty^{1,\emptyset,\alphabeta}) \rightarrow \Spec S_\infty$ is
  surjective and generalizing by~\cite[\href{http://stacks.math.columbia.edu/tag/080T}{Tag 080T}]{stacks-project}. It follows that we may pick some
$\rho_{\mathrm{reg}}:G_{F}\to\GSp_4(\cO)$ whose corresponding point is
in the support of
$M^{1,\emptyset,\alphabeta}/\p_{\kappa_{\mathrm{reg}}}M^{1,\emptyset,\alphabeta}$
and such that $\rho$ and $\rho_{\mathrm{reg}}$ (or
their corresponding points~$x$ and~$x_{\mathrm{reg}}$) lie on the same component
of $R_\infty^{1,\emptyset,\alphabeta}$; we write~$\eta$ for the
generic point corresponding to this component. Similarly, 
we may pick some
$\rho^\chi_{\mathrm{reg}}:G_{F}\to \GSp_4(\cO)$ whose corresponding 
point~$x^\chi_{\mathrm{reg}}$ is
in the support of
$M^{\chi,\emptyset,\alphabeta}/\p_{\kappa_{\mathrm{reg}}}M^{\chi,\emptyset,\alphabeta}$.

By Proposition~\ref{prop: spaces of invariants in pi regular weight} below,
we have 
\begin{equation*} \dim_E
(M_\infty^{1,\emptyset,\alphabeta}/\p_{x_{\mathrm{reg}}}M_\infty^{1,\emptyset,\alphabeta})[1/p]=\dim_E
S_{\kappa_{\mathrm{reg}},K^{p,\Iw}K_p(\emptyset)}^{\emptyset,\alphabeta}[\tilde{\p}^{\emptyset,\alphabeta}_{x_{\mathrm{reg}}}]=8d_\rho,
\end{equation*} and
\begin{equation*} \dim_E
(M_\infty^{\chi,\emptyset,\alphabeta}/\p_{x^\chi_{\mathrm{reg}}}M_\infty^{\chi,\emptyset,\alphabeta})[1/p]=\dim_E
S_{\kappa_{\mathrm{reg}},K^{p,\Iw_1}K_p(\emptyset),\chi}^{\emptyset,\alphabeta}[\tilde{\p}^{\emptyset,\alphabeta}_{x^\chi_{\mathrm{reg}}}]=8.
\end{equation*}In addition, there are automorphic
representations~$\pi_{\mathrm{reg}}$, $\pi^\chi_{\mathrm{reg}}$
of~$\GSp_4(\A_F)$ of weight~$\kappa_{\mathrm{reg}}$ and central character~$|\cdot|^2$ such that
$\rho_{\pi_{\mathrm{reg}},p}\cong\rho_{\mathrm{reg}}$ and $\rho_{\pi^\chi_{\mathrm{reg}},p}\cong\rho^\chi_{\mathrm{reg}}$.

Applying Lemma~ \ref{lem: smooth implies free main section} to $\p_x$ and
$\p_{x_{\mathrm{reg}}}$ (which we may, by our assumptions on~$\rho$,
and by Theorem~\ref{thm: existence and properties of Galois representations for
    automorphic representations}
  for~$\rho_{\pi_{\mathrm{reg}},p}$, together with
  Lemmas~\ref{lem: smooth points of Ihara1} and~\ref{lem: smooth points of the big local deformation
    rings}),we obtain
\begin{equation*} \dim_E (M_\infty^{1,\emptyset,\alphabeta}/\p_x
M_\infty^{1,\emptyset,\alphabeta})[1/p]=\dim_{(R_{\infty,\eta}^{1,\emptyset,\alphabeta})'}M_{\infty,\eta}^{1,\emptyset,\alphabeta}=
\dim_E
(M_\infty^{1,\emptyset,\alphabeta}/\p_{x_{\mathrm{reg}}}M_\infty^{1,\emptyset,\alphabeta})[1/p]
=8d_\rho.
\end{equation*} As
\begin{equation*}
S^{\emptyset,\alphabeta}_{\kappa_2,K^{p,\Iw},K_p(\emptyset)}[\tilde{\p}_{x}^{\emptyset,\alphabeta}]=(M_\infty^{1,\emptyset,\alphabeta}/\p_{x}
M_\infty^{1,\emptyset,\alphabeta})[1/p]^\vee,
\end{equation*} the theorem is proved for $I=\emptyset$.

Before we go on to the case that $\#I=1$, we note that we may also apply Proposition~\ref{prop: spaces of invariants in pi regular weight} and Lemma~
\ref{lem: smooth implies free main section} to $\p_{x^\chi_{\mathrm{reg}}}$ and conclude that
\begin{equation*}
\dim_{(R_{\infty,\eta^\chi}^{\chi,\emptyset,\alphabeta})'}M_{\infty,\eta^\chi}^{\chi,\emptyset,\alphabeta}=
\dim_E
(M_\infty^{\chi,\emptyset,\alphabeta}/\p_{x^\chi_{\mathrm{reg}}}M_\infty^{\chi,\emptyset,\alphabeta})[1/p]=8.
\end{equation*}

By another application of Lemma~\ref{lem: main cycle comparison lemma}
this implies that $d_\eta'=d_\rho$ (where~$d'_\eta$ is as in
Lemma~\ref{lem: main cycle comparison lemma}). Following
Remark~\ref{rem: multiplicities are independent of global situation}, we will apply this in
the case $\#I=1$ below.

Now consider the case that $\#I=1$.  We consider the automorphic representation~$\pi$
of Hypothesis~\ref{hyp: hypotheses on rho}. By assumption, for all finite places $v$ of~$F$ the
  representation~$\rho_{\pi,p}|_{G_{F_v}}$ is pure, and
  therefore determines a unique component
  of~$R^{1,\emptyset,\alphabeta}_\infty$, which we denote
   by~$\eta_{\pi}$. Arguing as above, we find that $d_{\eta_\pi}'=d_{\rho_{\pi,p}}$ (where~$d'_\eta$ is as in
Lemma~\ref{lem: main cycle comparison lemma}).
  Write~$\tilde{\p}_{\pi}^{I,\alphabeta}$ for the
  height one prime ideal determined by~$\rho_{\pi,p}$.  Then by Proposition \ref{prop: spaces of invariants in pi regular weight} we find that
\numequation\label{eqn: dim is at least 8 rho pi}
\dim_ES^{I,\alphabeta}_{\kappa_2,K^{p,\Iw}K_p(I)}[\tilde{\p}_{\pi}^{I,\alphabeta}]\geq 8d_{\rho_{\pi,p}}.
\end{equation} Again $M_\infty^{1,I,\alphabeta}$ is a maximal
Cohen--Macaulay $(R_\infty^{1,I,\alphabeta})'$-module and so we may apply
Lemmas~\ref{lem: main cycle comparison lemma} and~\ref{lem: smooth implies free main section} to the
$R_\infty^{1,I,\alphabeta}$-module $M_\infty^{1,I,\alphabeta}$ and the
$R_\infty^{\chi,I,\alphabeta}$-module
$M_\infty^{\chi,I,\alphabeta}$. We find that
\begin{align*}\frac{1}{d_{\rho}}\dim_E
  S^{I,\alphabeta}_{\kappa_2,K^{p,\Iw}K_p(I)}[\tilde{\p}_{x}^{I,\alphabeta}]
  &= \frac{1}{d_{\rho}}
\dim_E
(M_\infty^{1,I,\alphabeta}/ \tilde{\p}_{x}^{I,\alphabeta} M_\infty^{1,I,\alphabeta})[1/p]\\&=\frac{1}{d_{\rho}}\dim_{(R_{\infty,\eta}^{1,I,\alphabeta})'}M_{\infty,\eta}^{1,I,\alphabeta}\\
&=\dim_{(R_{\infty,\eta^\chi}^{\chi,I,\alphabeta})'}M_{\infty,\eta^\chi}^{\chi,I,\alphabeta}\\&=\frac{1}{d_{\rho_{\pi,p}}}\dim_{(R_{\infty,\eta_\pi}^{1,I,\alphabeta})'}M_{\infty,\eta_\pi}^{1,I,\alphabeta}\\
&= \frac{1}{d_{\rho_{\pi,p}}}
\dim_E
(M_\infty^{1,I,\alphabeta}/\tilde{\p}_{\pi}^{I,\alphabeta} M_\infty^{1,I,\alphabeta})[1/p]\\&=\frac{1}{d_{\rho_{\pi,p}}}\dim_E
S^{I,\alphabeta}_{\kappa_2,K^{p,\Iw}K_p(I)}[\tilde{\p}_{\pi}^{I,\alphabeta}].
\end{align*}
It follows from~(\ref{eqn: dim is at least 8 rho pi}) that
  \begin{equation*}
  \dim_E S^{I,\alphabeta}_{\kappa_2,K^{p,\Iw}K_p(I)}[\tilde{\p}_{x}^{I,\alphabeta}]\geq 8d_\rho.
  \end{equation*}
On the other hand by Proposition \ref{prop: doubling main section}, we have
that
\begin{equation*} 
\dim_E S^{I,\alphabeta}_{\kappa_2,K^{p,\Iw}K_p(I)}[\tilde{\p}_{x}^{I,\alphabeta}]\leq \dim_E S^{\emptyset,\alphabeta}_{\kappa_2,K^{p,\Iw}K_p(\emptyset)}[\tilde{\p}_{x}^{\emptyset,\alphabeta}]=8d_\rho,
\end{equation*}
and so the theorem is proved.
\end{proof}

\begin{prop}
  \label{prop: spaces of invariants in pi regular weight}In the notation of the proof of Theorem~\ref{thm: iwahori patching main section},
we have 
\begin{equation*} \dim_E S_{\kappa_{\mathrm{reg}},K^{p,\Iw}K_p(\emptyset)}^{\emptyset,\alphabeta}[\tilde{\p}^{\emptyset,\alphabeta}_{x_{\mathrm{reg}}}]=8d_\rho,
\end{equation*} 
\begin{equation*} \dim_E S_{\kappa_{\mathrm{reg}},K^{p,\Iw_1}K_p(\emptyset),\chi}^{\emptyset,\alphabeta}[\tilde{\p}^{\emptyset,\alphabeta}_{x^{\chi}_{\mathrm{reg}}}]=8,
\end{equation*}
\begin{equation*}
  \dim_E
S^{I,\alphabeta}_{\kappa_2,K^{p,\Iw}K_p(I)}[\tilde{\p}_{\pi}^{I,\alphabeta}]\ge
8d_{\rho_{\pi,p}}.
\end{equation*}

In addition, there are automorphic
representations~$\pi_{\mathrm{reg}}$, $\pi^\chi_{\mathrm{reg}}$
of~$\GSp_4(\A_F)$ of weight~$\kappa_{\mathrm{reg}}$ and central character~$|\cdot|^2$ such that
$\rho_{\pi_{\mathrm{reg}},p}\cong\rho_{\mathrm{reg}}$ and $\rho_{\pi^\chi_{\mathrm{reg}},p}\cong\rho^\chi_{\mathrm{reg}}$.
\end{prop}
\begin{proof} By Theorem~\ref{thm: ordinary classicity for I at most 1} and  Theorem~\ref{thm: cohomology in terms of automorphic forms},
  $\dim_E
  S_{\kappa_{\mathrm{reg}},K^{p,\Iw}K_p(\emptyset)}^{\emptyset,\alphabeta}[\tilde{\p}^{\emptyset,\alphabeta}_{x_{\mathrm{reg}}}]$
  is equal to\[ \sum_\pi\dim_{\Ebar}(\pi^\infty)^{K^{p,\Iw}
      K_p(\emptyset),\{U_{v,1}=\alphabeta_v^{\mathrm{reg}},U_{v,2}=\alpha_v^{\mathrm{reg}}\beta_v^{\mathrm{reg}}\}_{v\in
        S_p}}\]where
    the sum is over all the cuspidal automorphic representations~$\pi$
    of weight~$\kappa_{\mathrm{reg}}$ such that $\pi$ has central
      character~$|\cdot|^2$, $\pi_v$ is holomorphic for all
      places~$v|\infty$, and
      $\rho_{\pi,p}\cong \rho_{\mathrm{reg}}$; and we write
      $\alpha_v^{\mathrm{reg}},\beta_v^{\mathrm{reg}}$ for the lifts
      of~$\alphabar_v,\betabar_v$ determined
      by~$\rho_{\mathrm{reg}}|_{G_{F_v}}$. In particular, note that we
      can take~$\pi_{\mathrm{reg}}$ to be any of the automorphic
      representations~$\pi$ contributing to the sum.

Since~$\rho_{\mathrm{reg}}$ is irreducible, such a~$\pi$ is of
      general type in the sense of~\cite{MR2058604} by Lemma~\ref{lem: getting to general type}, and therefore
      corresponds to an essentially self-dual regular algebraic cuspidal
      automorphic representation~$\Pi$ of~$\GL_4(\A_F)$.  By strong
      multiplicity one for~$\GL_4$~\cite{MR623137}, $\Pi$ is uniquely determined by
      the condition that $\rho_{\Pi,p}\cong\rho_{\mathrm{reg}}$, so
      by Theorem~\ref{thm: arthur classification results} we see
      that we can rewrite the above sum as
      \[\left(\sum_{\pi_{v_0}\in
            L(\rho_{\mathrm{reg}}|_{G_{F_{v_0}}})}\dim\pi_{v_0}^{\Iw_1(v_0)}\right)\prod_{v\in R}\left(\sum_{\pi_v\in
            L(\rho_{\mathrm{reg}}|_{G_{F_v}})}\dim\pi_v^{\Iw(v)}\right) \]
      (note that at all places~$v\notin R\cup S_p\cup\{v_0\}$, we are
      taking the space of hyperspecial invariants in an unramified
      representation, which is 1-dimensional; and at the places
      $v\in S_p$, the contribution is 1-dimensional by
      Propositions~\ref{prop: ordclass} and~\ref{prop: ordinary eigenform in autrep}). 

By Proposition~\ref{prop: rec of principal series},
$\pi_{v_0}$ is an irreducible unramified principal series
representation; indeed, by the choice of~$v_0$,
$\rho_{\pi,p}|_{G_{F_{v_0}}}$ is unramified, and no two eigenvalues
of~$\rho_{\pi,p}|_{G_{F_{v_0}}}(\Frob_{v_0})$ can have ratio~$q_{v_0}$. It
follows from Propositions~\ref{prop: Jacquet module iso}
and~\ref{prop: Jacquet module} that we have $\dim\pi_{v_0}^{\Iw_1(v_0)}=8$. The claim then
      follows from Proposition~\ref{prop: size of L packets for Ihara
        avoidance} (which we can apply,
      because for each place~$v\in R$, $\rho_{\mathrm{reg}}|_{G_{F_v}}$ is pure by
      Theorem~\ref{thm: existence and properties of Galois
        representations for automorphic representations}~(4), and
      therefore the corresponding Weil--Deligne representation is generic by~Lemma~\ref{lem: smooth points of Ihara1}, so
      that the corresponding $L$-packet contains a generic
      representation by Proposition~\ref{prop: unipotent generic implies generic}).

The statement
for~$\rho_{\mathrm{reg}}^\chi$ reduces in the same way to the claim that
for each place~$v\in R$, we have \[\sum_{\pi_v\in
            L(\rho^\chi_{\mathrm{reg}}|_{G_{F_v}})}\dim\pi_v^{\Iw_1(v),\chi}=1,\]which
        follows from Proposition~\ref{prop: Ihara local
          global}.
        Finally, in the case of
        $S^{I,\alphabeta}_{\kappa_2,K^{p,\Iw}K_p(I)}[\tilde{\p}_{\pi}^{I,\alphabeta}]$,
        the result follows as above, by computing the contribution of
        the automorphic representation~$\pi$ of Hypothesis~\ref{hyp:
          hypotheses on rho} (note that it contributes by Theorem~~\ref{thm: cohomology
          in terms of automorphic forms}).\end{proof}

\begin{rem}\label{rem: remark about geometric BM for l not p}
In the course of the proof of Theorem~\ref{thm: iwahori patching main
  section}, we showed that for the generic point~$\eta$ corresponding
to~$\rho$, the quantity~$d'_\rho$ of Lemma~\ref{lem: main cycle
  comparison lemma} is equal to~$d_\rho$.
It is presumably
  possible to go further  following~\cite{shottonGLn}, and to use our patched modules to
  show that for each $v\in R$ and each generic point $\eta_v$ of
  $\Spec R_v^1$, if we write
\begin{equation*}
Z(\Spec R_v^{1,\eta_v}/(\lambda))=\overline{\eta}_v
\end{equation*}
then
\begin{equation*}
Z(\Spec (R_v^\chi)^{\mathrm{red}}/(\lambda))=\sum_{\eta_v} d_{\eta_v}\overline{\eta}_v,
\end{equation*}where $d_{\eta_v}$ is as in Proposition~\ref{prop: size of L packets for Ihara avoidance}. \end{rem}

\section{\'Etale descent and the main modularity lifting
  theorem} \label{sec:gluing}\subsection{Introduction}
    Our main goal is to remove the assumption $\#I \leq 1$ of  Theorem~\ref{thm: iwahori patching main
  section}   in order to eventually apply
Theorem~\ref{thm-classicality} with $I = S_p$, and from this  conclude that we have constructed  classical
automorphic representations. 
The starting point is to consider the spaces of $p$-adic automorphic
forms considered in Theorem~\ref{thm: iwahori patching main
  section} for \emph{both}~$\#I = 1$  and~$\#I = 0$. By studying the way in which these
  spaces are related, we will be able to (inductively) determine 
  precise linear combinations
  of such forms which belong to spaces of $p$-adic automorphic forms for larger~$\#I$.
  Our   argument uses the doubling
results of~\S\ref{sec:doubling}, the analytic continuation
results of~\S\ref{sec: higher Coleman theory}, and 
\'etale descent.  Finally, we apply
solvable base change to prove our main modularity lifting theorem.

We briefly indicate some of the main features of our argument. As we
mentioned in the introduction, the analytic continuation arguments
that we are using here are analogous to those used for Hilbert modular
forms of weight at least two, rather than those of weight one -- in
particular, there is no ``gluing'' of the kind used
in~\cite{MR1709306}, and we are simply analytically continuing a
single form at a time (using the method of Kassaei
series~\cite{MR2219265}). This part of the argument is quite standard,
although we have to take some care to show that the regions that we
have analytically continued to are large enough. For this reason, we
ignore the issues of analytic continuation in this introduction.

We show that the conclusion of Theorem~\ref{thm: iwahori patching main section} holds
for all~$I$ by induction on~$\#I$. The key step is to go
from~$\#I\le 1$ to~$\#I\le 2$; indeed, the general inductive step
considers two places $v_1,v_2$ dividing~$p$, and essentially ignores
the other places above~$p$, so for the purpose of exposition we assume
that~$S_p=\{v_1,v_2\}$. Write~$\alpha_i,\beta_i$
for~$\alpha_{v_i},\beta_{v_i}$, $i=1,2$. We denote the various spaces
of forms considered in Theorem~\ref{thm: iwahori patching main
  section} with~$I=\emptyset$ by~$V_{\alpha_1,\alpha_2}$,
$V_{\beta_1,\alpha_2}$, $V_{\alpha_1,\beta_2}$, $V_{\beta_1,\beta_2}$
(so that for example on~$V_{\alpha_1,\alpha_2}$ the eigenvalue
of~$U_{v_1,1}$ is~$\alpha_1$ and the eigenvalue of~$U_{v_2,1}$
is~$\alpha_2$). Each of these spaces has dimension~$d:=8d_\rho$, and
considering the action of~$U_{v_1,1}$ and~$U_{v_2,1}$, we see that
these spaces together span a $4d$-dimensional space of $p$-adic modular forms
of Iwahori level. 

We expect that this space contains a $d$-dimensional subspace of
$p$-adic modular forms which descend to Klingen level (and are
suitably overconvergent in both the $v_1$ and $v_2$ directions). The
difficulty (even if $d=1$) is in identifying this subspace; recall
that there is no obvious relationship between the Hecke eigenvalues
and Fourier coefficients. However, we also have the spaces of forms
for~$I=\{v_1\}$ and~$I=\{v_2\}$, which we denote
by~$V_{\alpha_1,\alpha_2+\beta_2}$, $V_{\beta_1,\alpha_2+\beta_2}$,
$V_{\alpha_1+\beta_1,\alpha_2}$, $V_{\alpha_1+\beta_1,\beta_2}$, where
for example the forms in $V_{\alpha_1,\alpha_2+\beta_2}$ have Klingen
level at~$v_2$ (and are highly overconvergent in the $v_2$ direction),
and are $U_{\Kli(v_2),1}$-eigenforms with
eigenvalue~$\alpha_2+\beta_2$. Again, all of these spaces has
dimension~$d$ by Theorem~\ref{thm: iwahori patching main section}.

Now, the relations between the Hecke operators at Klingen and Iwahori
levels (more precisely, Lemma~\ref{lem: quadratic relation for U Kli 1 big sheaf}) imply that we have a map
\[(U_{v_1,1}-\beta_1):V_{\alpha_1+\beta_1,\alpha_2}\to
  V_{\alpha_1,\alpha_2}.\] Furthermore, this map is injective by
Proposition \ref{prop: doubling main section} (that is, by our main
doubling results), and since the source and target both have
dimension~$d$, this map is in fact an isomorphism. Similarly, we have an
isomorphism \[(U_{v_1,1}-\beta_1):V_{\alpha_1+\beta_1,\beta_2}\isoto
  V_{\alpha_1,\beta_2}\] and thus an
isomorphism of $2d$-dimensional spaces \numequation\label{eqn: intro
  main linear algebra iso}(U_{v_1,1}-\beta_1):V_{\alpha_1+\beta_1,\alpha_2}\oplus
  V_{\alpha_1+\beta_1,\beta_2}\isoto
  V_{\alpha_1,\alpha_2}\oplus V_{\alpha_1,\beta_2}.\end{equation}By pulling back
from Iwahori to Klingen level, we can think
of~$V_{\alpha_1,\alpha_2+\beta_2}$ as a $d$-dimensional subspace of
the target of~\eqref{eqn: intro
  main linear algebra iso}. The inverse image of this space in the
source of~\eqref{eqn: intro
  main linear algebra iso} is the $d$-dimensional space of forms that
we are seeking; considered as living on the right hand side of~\eqref{eqn: intro
  main linear algebra iso}, it comes from Klingen level at~$v_2$, and
on the left hand side of~\eqref{eqn: intro
  main linear algebra iso}, it comes from Klingen level at~$v_1$. We
make this precise using an argument with \'etale descent.

\subsection{\'Etale descent}In this section we carry out the argument
explained above, showing that the conclusion of Theorem~\ref{thm: iwahori patching main section} holds
for all~$I$ by induction on~$\#I$ (in fact, we show slightly more,
keeping track of the overconvergence of our $p$-adic modular forms). Recall that by definition for each choice of~$I,\alphabeta$ we
have \begin{multline*}
S_{\kappa_2,K^{p, \Iw} K_p(I)}^{I,\alphabeta}=(H^0(\mathfrak{X}^{I,G_1}_{K^{p, \Iw}K_p(I)},\omega^{2}(-D))_{\tilde{\m}^{I,\alphabeta}})\otimes_\cO
E[\{U_{v,0}-1\}_{v\in S_p},\\ \{U_{v,0}-q_v^{-2}\}_{v\in
R},\{T_{v,0}-q_v^{-2}\}_{v\not\in S_p\cup R\cup \{v_0\}}]
\end{multline*}  

The maximal ideal $\tilde{\m}^{I,\alphabeta}$ of the Hecke algebra is defined in equation (\ref{eqn-def-max-ideal}). It contains an ordinary projector. 
We have given ourselves (see the beginning of \S\ref{subsec: patching padic multiplicities}) a Galois representation $\rho$ satisfying Hypothesis~\ref{hyp: hypotheses on rho}. We want to prove that it is modular. Associated to this representation is a point $x$ on the deformation space of $\bar{\rho}$ and an ideal ~$\tilde{\p}^{I,\alphabeta}_{{x}}$ (see equation (\ref{eqn-def-primeideal}))  of the Hecke algebra contained in $\tilde{\m}^{I,\alphabeta}$ whose definition we recall here for convenience. It is the ideal of the Hecke 
algebra~$\TTt^I$ given by
  \[\left(\bigotimes_{v\not\in S_p\cup
R\cup\{v_0\}}\cO[\GSp_4(F_v) \doubleslash \GSp_4(\cO_{F_v})]\right)[\{U_{v,0},U_{\Kli(v),1},U_{v,2}\}_{v\in I},\{U_{v,0},U_{v,1},U_{v,2}\}_{v\in I^c}]\]
which is generated by:
\begin{itemize}
\item  the coefficients of $\det(X-\rho(\Frob_v))-Q_v(X)$
for each $v\not\in S_p\cup
R\cup\{v_0\}$, and
\item $\{U_{v,0}-1,U_{v,2}-\alpha_v\beta_v\}_{v\in
S_p},\{U_{\Kli(v),1}-\alpha_v-\beta_v\}_{v\in
I},\{U_{v,1}-\alphabeta_v\}_{v\in
I^c}$,
where, for $v\in S_p$, $\alpha_v$, $\beta_v$ are determined by~$\rho|_{G_{F_v}}$ as in Definition~\ref{defn:
  generic flat ordinary}. 
\end{itemize}

Recall that $\mathcal{X}^{I,G_1}_{K^{p, \Iw} K_p(I)}$ is the analytic adic  space over $\C_p$ associated to $\mathfrak{X}^{I,G_1}_{K^{p, \Iw}K_p(I)}$. By definition, we have: \[
  S^{I,\alphabeta}_{\kappa_2,K^{p,\Iw}K_p(I)}[\tilde{\p}_{{x}}^{I,\alphabeta}] \otimes_{E} \whalingship_p =e(U^I)H^0(\mathcal{X}^{I,G_1}_{K^{p, \Iw}K_p(I)},\omega^{2}(-D))[\tilde{\mathfrak{p}}_x^{I,\alphabeta}].\]

We may also introduce overconvergent versions of these spaces. 
Recall that we defined the dagger space~$\cX^{G_1,\mult,\dag}_{K^{p,\Iw}K_p(I)}$
in~(\ref{eqn: X mult}) (whose associated rigid analytic space is $\cX^{G_1,I}_{K^{p,\Iw}K_p(I)}$).  

 There is a natural injective restriction map: 
 $$e(U^I)H^0(\mathcal{X}^{G_1,\mult,
    \dag}_{K^{p, \Iw}K_p(I)},\omega^{2}(-D)) \rightarrow e(U^I)H^0(\mathcal{X}^{G_1, I}_{K^{p, \Iw}K_p(I)},\omega^{2}(-D)).$$ 
    
    Let $$S^{I,\alphabeta,\dagger}_{\kappa_2,K^{p,\Iw}K_p(I)}  =  e(U^I)H^0(\mathcal{X}^{G_1,\mult,
    \dag}_{K^{p, \Iw}K_p(I)},\omega^{2}(-D)) \cap S^{I,\alphabeta}_{\kappa_2,K^{p,\Iw}K_p(I)} \otimes \C_p$$ 
    where the intersection is taken inside $e(U^I)H^0(\mathcal{X}^{G_1, I}_{K^{p, \Iw}K_p(I)},\omega^{2}(-D))$.  
 
\begin{thm}\label{thm: main thm on gluing}Assume that~$\rho$ satisfies
  Hypothesis~\ref{hyp: hypotheses on rho}. Then for any $I\subset S_p$ and choice of $\alphabeta$,  we have \[\dim_E
S^{I,\alphabeta}_{\kappa_2,K^{p,\Iw}K_p(I)}[\tilde{\mathfrak{p}}_{{x}}^{I,\alphabeta}]=\dim_{\whalingship_p}
S^{I,\alphabeta,\dagger}_{\kappa_2,K^{p,\Iw}K_p(I)}[\tilde{\mathfrak{p}}_x^{I,\alphabeta}]=8d_{\rho}.\] \end{thm}

Before proving this theorem, we record the following important corollary.
\begin{cor}
  \label{cor: rho is automorphic}Suppose that~$\rho$ satisfies
  Hypothesis~\ref{hyp: hypotheses on rho}. Then~$\rho$ is modular.
  More
  precisely, there is an ordinary automorphic
  representation~$\pi'$ of $\GSp_4(\A_F)$ of parallel weight~$2$ and
  central character~$|\cdot|^2$, with
  $\rho_{\pi',p}\cong\rho$, and for every finite
  place~$v$ of~$F$ we have \[\WD(\rho|_{G_{F_v}})^{F-\semis}\cong\recGTp(\pi'_v\otimes|\nu|^{-3/2}).\]
\end{cor}
\begin{proof}
  The existence of~$\pi'$ with~$\rho_{\pi',p}\cong\rho$ is
  immediate from Theorem~\ref{thm: main thm on gluing},
  taking~$I=S_p$, together with Theorem~\ref{thm-classicality} and Theorem~\ref{thm: cohomology in terms
    of automorphic forms}. By Corollary~\ref{cor: Galois representations for limits of
    discrete series with ordinarity at p} we
  have
  \[\WD(\rho|_{G_{F_v}})^{\semis}\cong\recGTp(\pi'_v\otimes|\nu|^{-3/2})^{\semis}\]
  at all finite places~$v$ of~$F$, so we need only prove that the
  monodromy operators agree at the places~$v\in
  R$. Since~$\rho|_{G_{F_v}}$ is pure by assumption, it follows from Lemma~\ref{lem:
    pure implies N is maximal} that it suffices to prove, in the
  notation of Section~\ref{subsec: nonarch LL},
  that $n(\rho|_{G_{F_v}})\le n(\pi'_v) $.

Now, if~$\pi_v$ is any irreducible admissible representation
of~$\GSp_4(F_v)$, then an examination of \cite[Table A.15]{MR2344630} (noting that the column
there headed ``a'' records $n(\pi_v)$) shows that:
\begin{itemize}
\item $n(\pi_v)\ge 1$ if and only if
  $(\pi_v)^{\GSp_4(\cO_{F_v})}=0$.
\item  $n(\pi_v)\ge 2$ if and only if
  $(\pi_v)^{\GSp_4(\cO_{F_v})}=(\pi_v)^{\Par(v)}=0$.
\item $n(\pi_v)=3$ if and only if
  $(\pi_v)^{\Kli(v)}=(\pi_v)^{\Si(v)}=0$.
\end{itemize}
 
Suppose that~$n(\rho|_{G_{F_v}})=1$, so that we need to show that
$(\pi'_v)^{\GSp_4(\cO_{F_v})}=0$. Suppose for the sake of contradiction
that~$(\pi'_v)^{\GSp_4(\cO_{F_v})}\ne 0$; then by Hida theory (more
precisely, by Theorem~\ref{thm: existence of Galois representations in big Hecke
    algebras} and its proof), the Galois
  representation~$\rho_{\pi',p}$ is a $p$-adic limit of Galois
  representations~$\rho_{\pi'',p}$, where~$\pi''$ has regular
  weight and satisfies $(\pi''_v)^{\GSp_4(\cO_{F_v})}\ne 0$. In
  particular, by Theorem~\ref{thm: existence and properties of Galois representations for
    automorphic representations},
  $n(\rho_{\pi'',p}|_{G_{F_v}})=n(\pi''_v)=0$. By the
  semicontinuity of the rank of
  the nilpotent operator~$N$ in such a family, it follows
  that~$n(\rho|_{G_{F_v}})= 0$, a contradiction. We leave the (very
  similar) arguments in the cases~$n(\rho|_{G_{F_v}})=2,3$ to the reader.
\end{proof}

\begin{proof}[Proof of Theorem~\ref{thm: main thm on gluing}]
We prove this by induction on~$\# I$. The result is true for~$I =
\emptyset$ and~$\# I = 1$ by Theorem~\ref{thm: iwahori patching main
  section} and Theorem~\ref{thm: ordinary is overconvergent I equal
  1} (ordinary implies overconvergent if~$\#I\le 1$). For any~$I$, the restriction map \[
  S^{I,\alphabeta,\dagger}_{\kappa_2,K^{p,\Iw}K_p(I)}[\tilde{\p}_x^{I,\alphabeta}] \to
S^{I,\alphabeta}_{\kappa_2,K^{p,\Iw}K_p(I)}[\tilde{\p}_{{x}}^{I,\alphabeta}] \otimes_{E} \whalingship_p \]is
injective, while by  Proposition~\ref{prop: doubling main section} (and a
simple induction) we see that for any~$I$, the dimension
of~$S^{I,\alphabeta}_{\kappa_2,K^{p,\Iw}K_p(I)}[\tilde{\p}_{{x}}^{I,\alphabeta}]$
is at most~$8d_{\rho}$. It therefore suffices to show that
$S^{I,\alphabeta,\dagger}_{\kappa_2,K^{p,\Iw}K_p(I)}[\tilde{\p}_x^{I,\alphabeta}]$
has dimension at
least~$8d_{\rho}$. We may assume that~$\# I \ge 2$, and hence we may write~$I$ as a disjoint union~$J \cup \{v_1,v_2\}$ for two primes~$v_i | p$.
We fix the choice of~$\alphabeta$ at all primes in~$J\cup I^c$.

By the inductive hypothesis applied to~$J$, for each choice of
~$\alphabeta$  at~$v_1$ and~$v_2$ the corresponding eigenspace
$S^{J,\alphabeta,\dagger}_{\kappa_2,K^{p,\Iw}K_p(J)}[\tilde{\p}_x^{J,\alphabeta}]$
is $8d_{\rho}$-dimensional, and we denote these eigenspaces by
~$V_{\alpha_1,\alpha_2}$, $V_{\beta_1,\alpha_2}$,
$V_{\alpha_1,\beta_2}$, $V_{\beta_1,\beta_2}$ (so that for example
on~$V_{\alpha_1,\alpha_2}$ the eigenvalue of~$U_{v_1,1}$
is~$\alpha_1$ and the eigenvalue of~$U_{v_2,1}$
is~$\alpha_2$). Considering the action of~$U_{v_1,1}$ and~$U_{v_2,1}$,
we see that these spaces span a $4 \times 8 d_{\rho} = 32d_\rho$-dimensional subspace~$V$ of $H^0(\mathcal{X}^{J,G_1}_{K^{p,\Iw}K_p(J)},\omega^{2}(-D))$.

By the inductive hypothesis applied to~$J \cup \{v_1\}$ and the two
possible choices of~$\alphabeta$  at~$v_2$ (the choice at~$v_1$ is
irrelevant), we see that the eigenspaces
$S^{J\cup\{v_1\},\alphabeta,\dagger}_{\kappa_2,K^{p,\Iw}K_p(J\cup\{v_1\})}[\tilde{\p}_x^{J\cup\{v_1\},\alphabeta}]$
are both $8d_{\rho}$-dimensional, and we denote the corresponding spaces
by~$V_{\alpha_1+\beta_1,\alpha_2}$
and~$V_{\alpha_1+\beta_1,\beta_2}$. Similarly, the inductive
hypothesis applied to~$J \cup \{v_2\}$ yields $8d_{\rho}$-dimensional
spaces~$V_{\alpha_1,\alpha_2+\beta_2}$
and~$V_{\beta_1,\alpha_2+\beta_2}$.

Recall that our goal is to construct an $8d_{\rho}$-dimensional space of
eigenforms~$V_{\alpha_1 + \beta_1,\alpha_2 + \beta_2}$ which are
eigenforms for the operators~$U_{\Kli(v_1),1}$ and~$U_{\Kli(v_2),1}$
(and for the Hecke operators at all the other places), and
which lie in 
$$H^0(\mathcal{X}^{G_1,\mult, \dag}_{K^{p, \Iw} K_p(I)},\omega^{2}(-D)).$$
We will combine the analytic continuation results of~\S\ref{sec: higher Coleman theory} with a
descent argument to prove the existence of the sought-after
eigenforms.  

We need to introduce some notation in order to be able to describe the adic spaces we are working with. 
Recall that  $\mathcal{X}_{K^{p, \Iw}K_p(I)}$ is the analytic space associated to $X_{K^{p, \Iw} K_p(I)}$.  For each $v \in I$,  $H_v$ refers to the quasi-finite subgroup (of order $p$ over the interior of the moduli space) related to  the Klingen level structure,  and for each $v \in I^c$, ~$L_v \supset H_v$
refers to the quasi-finite  (maximally isotropic rank~$p^2$ over the interior of the moduli space) subgroup  corresponding to  the Iwahori level structure. 

For any tuple $(\epsilon_v) \in [0,1]^I \times [0,2]^{I^c}$ we defined  an analytic adic space $\mathcal{X}_{K^{p, \Iw}
   K_p(I)}((\epsilon_v)_{v \in S_p})$  which is the open subspace of $\mathcal{X}_{K^{p, \Iw}
   K_p(I)}$ where: 
   
   \begin{enumerate} 

\item If $v \in I$,  the degree of the subgroup~$H_v$, which takes values in~$[0,1]$, is greater or equal than $1-\epsilon_v$. 
\item If $v \in I^c$,  the degree of the subgroup~$L_v$ of rank~$p^2$, which takes values in~$[0,2]$, is greater or equal than $2-\epsilon_v$.
Note that we have~$\deg(L_v) = \deg(H_v) + \deg(L_v/H_v)$. 
\end{enumerate}

It will be convenient to adopt the following notation in this proof (note that $I$ is fixed). We write (cf.\ (\ref{eqn: X mult})) \[\cX^{\mult} = \mathcal{X}_{K^{p, \Iw}
   K_p(I)}((0)_{v \in S_p}),\] $$\cX^{\mult, \dagger} = \mathcal{X}^{\mult, \dag}_{K^{p,\Iw}K_p(I)} = \lim_{\epsilon_v \rightarrow 0^+} \mathcal{X}_{K^{p, \Iw}
   K_p(I)}((\epsilon_v)_{v \in S_p}),$$   and    ~$\cX^{\mult, \ddagger}$ for the dagger space
 \[\mathcal{X}^{\mult, \ddagger}:= \lim_{\epsilon_v \rightarrow 0^+} \mathcal{X}_{K^{p, \Iw}
   K_p(I)}((\epsilon_v)_{v \in S_p \setminus \{v_1,v_2\}}),\] 
   where we take the limit over all primes \emph{except}~$v_1$ and~$v_2$. It follows that: 
   \[\mathcal{X}^{\mult, \dagger}  = \lim_{\epsilon_{v_1}, \epsilon_{v_2} \rightarrow 0^+} \mathcal{X}^{\mult, \ddagger}(\epsilon_{v_1}, \epsilon_{v_2}),\]
   and there are maps of locally ringed spaces~$\cX^{\mult} \rightarrow \cX^{\mult, \dagger} \rightarrow \cX^{\mult, \ddagger}$.  
   By adding the subscript~$\Iw(v_1)$ (or~$\Iw(v_2)$, or~$\Iw(v_1,v_2)$) to~$\cX^{\mult}$, $\cX^{\mult, \ddagger}$, $\cX^{\mult, \dagger}$
   we mean the space where one has now added an Iwahori level structure at~$v_1$ (or~$v_2$, or~$v_1$ and~$v_2$) to the relevant space.
  For~$i=1,2$ we write
$d_i^H=\deg H_{v_i}$, $d_i^L=\deg L_{v_i}$, whenever these quantities
are defined. We will adorn~$\cX^{\mult}$ and~$\cX^{\mult,\ddagger}$ with superscripts
indicating the regions (which will typically strictly contain~$\cX^{\mult}$ and~$\cX^{\mult,\ddagger}$)  where various inequalities hold.

Returning to the spaces we defined above, we have
$$\begin{aligned}
V_{\alphabeta_1,\alphabeta_2}  & \subset  \  H^0(\mathcal{X}^{\mult, \dag}_{K^{p, \Iw} K_p(J)},\omega^2(-D)), \\
V_{\alphabeta_1 + \alphabeta'_1,\alphabeta_2} & \subset \ H^0(\mathcal{X}^{\mult, \dag}_{K^{p, \Iw} K_p(J\cup\{v_1\})},\omega^2(-D)), \\
V_{\alphabeta_1,\alphabeta_2 + \alphabeta'_2} & \subset \ H^0(\mathcal{X}^{\mult, \dag}_{K^{p, \Iw}K_p(J\cup\{v_2\})},\omega^2(-D)). \end{aligned}
$$

\begin{lemma}\label{lem: analytic continuation} The elements
  of~$V_{\alphabeta_1,\alphabeta_2}$ extend to~$\cX^{\mult, \ddagger,
    d^{H}_1 \ge 1 - \epsilon, d^{H}_2 \ge 1 - \epsilon,
    d^L_1>1,d^L_2>1}_{\Iw(v_1,v_2)}$ for some~$\epsilon>0$.  Similarly, the elements of~$V_{\alphabeta_1 + \alphabeta'_1,\alphabeta_2}$ and~$V_{\alphabeta_1,\alphabeta_2 + \alphabeta'_2}$
extend to the spaces $\cX^{\mult, \ddagger,d^{H}_1 \ge 1 - \epsilon, d^H_2 > 1- \epsilon, d^L_2>1}_{\Iw(v_2)}$ and  $\cX^{\mult, \ddagger, d^H_1> 1- \epsilon, d^L_1>1,d^{H}_2 \ge 1 - \epsilon}_{\Iw(v_1)}$
 respectively for some~$\epsilon > 0$.
\end{lemma}

\begin{proof} This follows from Lemma~\ref{lem-extension-first} (taking~$I$ there to be~$J$,
  $J\cup\{v_2\}$ and $J\cup\{v_1\}$ respectively). Note that our forms
  are ordinary for~$U_{w,1}$ for the appropriate~$w$, and
  therefore of finite slope for these operators. \end{proof}

By Koecher's principle, 
all of our cohomology groups may be replaced by the cohomology of the corresponding open spaces of ``good reduction''~$\cY^{\mult} \subset \cX^{\mult}$, $\cY^{\mult, \ddagger} \subset \cX^{\mult, \ddagger}$, and~$\cY^{\mult, \dagger} \subset \cX^{\mult, \dagger}$
respectively. (Since the sheaf $\omega^2$ is pulled back from the
minimal compactification, the form of Koecher's principle we are using
is just the following statement: if $\fX$ is a normal formal scheme,
$\fY\subseteq\fX$ is an open formal subscheme whose complement is
codimension $\geq 2$, and $\cL/\fX$ is a line bundle, then
$H^0(\fX,\cL)=H^0(\fY,\cL)$. That the boundary in the minimal
compactification does indeed have codimension $\ge 2$ follows from an analysis of the blowup in the boundary charts.) We now restrict to these spaces to avoid minor technical issues related to the boundary. In particular we will want to use that forgetting the level structure induces finite \'etale maps between our spaces. The reader will check easily that the forms we construct are indeed cuspidal because they  are obtained by ``descent'' of cuspidal forms.

For any $\epsilon>0$ and for $i=1,2$ there is a finite \'etale
map:
$$q_{v_i}: \cY^{\mult, \ddagger, d^H_1\ge 1-\epsilon,d^H_2\ge 1-\epsilon}_{\Iw(v_i)} \rightarrow \cY^{\mult, \ddagger,d^H_1\ge 1-\epsilon,d^H_2\ge 1-\epsilon}.$$
There is a corresponding fibre product map (still finite \'etale):
$$q_{v_1}: \cY^{\mult, \ddagger,d^H_1\ge 1-\epsilon,d^H_2\ge 1-\epsilon }_{\Iw(v_1,v_2)} \rightarrow
\cY^{\mult, \ddagger,d^H_1\ge 1-\epsilon,d^H_2\ge 1-\epsilon}_{\Iw(v_2)}$$(and similarly for~$q_{v_2}$).

We now somewhat abusively also write
$V_{\alpha_1+\beta_1,\alpha_2}$ instead of~$q^*_{v_1}V_{\alpha_1+\beta_1,\alpha_2}$. 
We claim
that the action of~$(U_{v_1,1}-\beta_1)$ induces an isomorphism
\numequation\label{eqn: Z1-alpha1
  iso}(U_{v_1,1}-\beta_1):V_{\alpha_1+\beta_1,\alpha_2}\isoto
V_{\alpha_1,\alpha_2}.\end{equation}To see this, note that 
~$(U_{v_1,1}-\beta_1)$ is injective by Proposition \ref{prop: doubling main section}, and both spaces have the same dimension.  

In the same way, we have an isomorphism
$(U_{v_1,1}-\alpha_1):V_{\alpha_1+\beta_1,\alpha_2}\isoto
V_{\beta_1,\alpha_2}$, so we see that in fact the span
of~$V_{\alpha_1+\beta_1,\alpha_2}$
and~$U_{v_1,1}V_{\alpha_1+\beta_1,\alpha_2}$ is
exactly~$V_{\alpha_1,\alpha_2}\oplus V_{\beta_1,\alpha_2}$. It follows
from Lemma~\ref{lem: analytic continuation} that the forms
in~$V_{\alpha_1,\alpha_2}\oplus V_{\beta_1,\alpha_2}$ extend
to~$\cY^{\mult, \ddagger, d^H_1\geq 1-\epsilon,  d^H_2\geq 1-\epsilon, d^L_2> 1}_{\Iw(v_1,v_2)}$ for some~$\epsilon>0$. 

Set \[U_\epsilon:=\cY^{\mult, \ddagger, d^H_1\geq 1-\epsilon,  d^H_2\geq 1-\epsilon, d^L_2> 1}_{\Iw(v_1,v_2)}. \]The map~$q_{v_2}$ restricts to an  \'etale map:
 \[q_{v_2}:U_\epsilon\to\cY^{\mult, \ddagger, d^H_1\ge
  1-\epsilon,d^H_2\ge 1-\epsilon}_{\Iw(v_1)}. \]
  We claim that for $\epsilon$ sufficiently small, the restriction of~$q_{v_2}$ to~$U_\epsilon$ is surjective. Note that a pre-image of a  point in~$ \cY^{\mult, \ddagger,  d^H_1\geq 1-\epsilon,   d^H_2\geq 1-\epsilon}_{\Iw(v_1,v_2)}$ (without any condition on~$d^L_2$) corresponds to a choice of~$L = L_{v_2}$, which is determined
by a line in~$H^{\perp}/H \subset A[v_2]/H$ for~$H = H_{v_2}$.  We need to show that $\deg(L)>1$ for at least one such $L$.

Let us first assume that $\deg(H) = 1$.  Then we can choose any line $C\subset H^\perp/H$ with $\deg(C)>0$ (such a $C$ exists as $H^\perp/H$ is not \'etale) and the corresponding $L$ has $\deg(L)=\deg(H)+\deg(C)>1$.

We pass from $\deg(H) = 1$ to $\deg(H)>1-\varepsilon$ by a continuity argument.  The function which sends a rank one point $x \in \cY^{\mult, \ddagger, d^H_1\ge1-\epsilon,d^H_2\ge 1-\epsilon}_{\Iw(v_1)}$ to the maximum of $\deg(L)$ is continuous.  It follows that for $\varepsilon$ sufficiently small, we can ensure the existence of a subgroup $L$ such that $\deg(L) >1$. 

Consider the corresponding descent diagram:
$$\begin{diagram}U_\epsilon \times_{\cY^{\mult, \ddagger, d^H_1\ge
  1-\epsilon,d^H_2\ge 1-\epsilon}_{\Iw(v_1)}} U_\epsilon  &   \pile{\rTo^{\ \ \   q_{v_2,1} } \\ \rTo_{\ \ \  q_{v_2,2}}}  & U_\epsilon & \rTo^{q_{v_2} \quad } &  \cY^{\mult, \ddagger, d^H_1\ge
  1-\epsilon,d^H_2\ge 1-\epsilon}_{\Iw(v_1)}. \end{diagram}$$

\begin{lemma}\label{lemma-dificil} After possibly further
  shrinking~$\epsilon>0$, any element of~$V_{\alpha_1,\alpha_2+\beta_2}$ descends to $\cY^{\mult, \ddagger, d^H_1\ge
  1-\epsilon,d^H_2\ge 1-\epsilon}_{\Iw(v_1)}$.
\end{lemma}
\begin{proof}  
Any element of  $V_{\alpha_1,\alpha_2+\beta_2}$ 
tautologically satisfies descent over the (smaller) space~$
\cY^{\mult, \ddagger, d^H_1\geq 1-\epsilon,  d^L_1> 1, d^H_2\geq
  1-\epsilon, d^L_2>1}_{\Iw(v_1,v_2)} \subset U_\epsilon$ to $$q_{v_2} (\cY^{\mult, \ddagger, d^H_1\geq 1-\epsilon,  d^L_1> 1, d^H_2\geq
  1-\epsilon, d^L_2>1}_{\Iw(v_1,v_2)} ) = \cY^{\mult, \ddagger, d^H_1\geq 1-\epsilon,  d^L_1> 1, d^H_2\geq
  1-\epsilon}_{\Iw(v_1)}$$  since it is (by Lemma~\ref{lem: analytic continuation}) obtained simply by pulling back
a form on this space under~$q_{v_2}$. 
Therefore, we deduce that  for any element $G \in V_{\alpha_1, \alpha_2 + \beta_2}$, we have that $q_{v_2,1}^\star G = q_{v_2,2}^\star G$ on 

$$ \cY^{\mult, \ddagger, d^H_1\geq 1-\epsilon,  d^L_1> 1, d^H_2\geq
  1-\epsilon, d^L_2>1}_{\Iw(v_1,v_2)} \times_{\cY^{\mult, \ddagger, d^H_1\geq 1-\epsilon, d^H_2\geq
  1-\epsilon}_{\Iw(v_1)}} \cY^{\mult, \ddagger, d^H_1\geq 1-\epsilon,  d^L_1> 1, d^H_2\geq
  1-\epsilon, d^L_2>1}_{\Iw(v_1,v_2)}.$$

The point is now to show that each connected component of $$U_\epsilon \times_{\cY^{\mult, \ddagger, d^H_1 \ge 1 - \epsilon, d^H_2 \ge 1 - \epsilon}_{\Iw(v_1)}} U_\epsilon$$ intersects 
$$ \cY^{\mult, \ddagger, d^H_1\geq 1-\epsilon,  d^L_1> 1, d^H_2\geq
  1-\epsilon, d^L_2>1}_{\Iw(v_1,v_2)} \times_{\cY^{\mult, \ddagger, d^H_1\geq 1-\epsilon, d^H_2\geq
  1-\epsilon}_{\Iw(v_1)}} \cY^{\mult, \ddagger, d^H_1\geq 1-\epsilon,  d^L_1> 1, d^H_2\geq
  1-\epsilon, d^L_2>1}_{\Iw(v_1,v_2)}$$
  so that  we have that $q_{v_2,1}^\star G = q_{v_2,2}^\star G$ on $U_\epsilon \times_{\cY^{\mult, \ddagger, d^H_1 \ge 1 - \epsilon, d^H_2 \ge 1 - \epsilon}_{\Iw(v_1)}} U_\epsilon$ and can perform the descent of $G$.

It follows from \cite[Thm.\ 2]{MR2388558} that after possibly further shrinking~$\epsilon$, there is a surjective map   
\numequation\label{eqn: components same 0 epsilon}\pi_0(U_0 \times_{\cY^{\mult, \ddagger d^H_1=d^H_2=1}_{\Iw(v_1)}} U_0)
\rightarrow \pi_0(U_\epsilon \times_{\cY^{\mult, \ddagger, d^H_1 \ge 1 - \epsilon, d^H_2 \ge 1 - \epsilon}_{\Iw(v_1)}} U_\epsilon).\end{equation}

We need to see that
$q^*_{v_2,1} G = q^*_{v_2,2} G$ on $$U_\epsilon \times_{\cY^{\mult,
    \ddagger, d^H_1 \ge 1 - \epsilon, d^H_2 \ge 1 -
    \epsilon}_{\Iw(v_1)}} U_\epsilon.$$ By~(\ref{eqn: components same 0 epsilon}), it is enough to  show~$q^*_{v_2,1} G = q^*_{v_2,2} G$ on the subspace
 $$U_0 \times_{\cY^{\mult, \ddagger, d^H_1=d^H_2=1}_{\Iw(v_1)}} U_0.$$
As discussed above, this identity
holds over the region
\[\cY^{\mult, \ddagger, d^H_1=d^H_2=1,d^L_1> 1,d^L_2> 1}_{\Iw(v_1,v_2)} \times_{\cY^{\mult, \ddagger, d^H_1=d^H_2=1}_{\Iw(v_1)}}
\cY^{\mult, \ddagger, d^H_1=d^H_2=1,d^L_1> 1,d^L_2>1}_{\Iw(v_1,v_2)}\] by
definition. It therefore suffices to show that this region
intersects all connected
components of $U_0 \times_{\cY^{\mult, \ddagger,d^H_1=d^H_2=1}_{\Iw(v_1)}} U_0$.

 Accordingly, it is enough to show that every connected component of $U_0
\times_{\cY^{\mult, \ddagger,d^H_1=d^H_2=1}_{\Iw(v_1)}} U_0$ contains
a point which is non-ordinary at~$v_1$. Indeed, if such
a point had~$d^L_1=1$, then we would have ~$\deg(H_{v_1})=\deg(L_{v_1})
= 1$, which implies that~$\deg(L_{v_1}/H_{v_1}) = 0$, so~$L_{v_1}/H_{v_1}$ is
 \'etale, and the point is ordinary at~$v_1$, a contradiction.

 We will prove in  Corollary~ \ref{cor-connected-components} below that
 any connected component of either of the
 spaces~$$\cY_{\Iw(v_1)}^{\mult, \ddagger,
   d^H_1=d^H_2=1,=_{v_2} 1},\ \cY_{\Iw(v_1)}^{\mult, \ddagger,
   d^H_1=d^H_2=1, =_{v_2} 2 }$$ contains a  point which is
 non-ordinary at $v_1$. Recall that the superscripts $=_{v_2} 1$ and $=_{v_2} 2 $ respectively mean the rank $1$ and the ordinary locus at $v_2$.

Now we observe that the maps $U_0^{=_{v_2} 1 } \rightarrow
\cY_{\Iw(v_1)}^{\mult, \ddagger, d^H_1=d^H_2=1,=_{v_2} 1}$   and
$U_0^{=_{v_2} 2 } \rightarrow \cY_{\Iw(v_1)}^{\mult, \ddagger,
  d^H_1=d^H_2=1,=_{v_2} 2}$ are both finite \'etale. It follows  from
Lemma~\ref{lem-con-fin-et} below  that any connected
component of any of the spaces $U_0^{=_{v_2} 1 }$, $U_0^{=_{v_2} 2 }$,
$U_0^{=_{v_2} 1 } \times_{\cY^{\mult, \ddagger,
    d^H_1=d^H_2=1}_{\Iw(v_1)}} U_0^{=_{v_2} 1 }$ or $U_0^{=_{v_2} 2 }
\times_{\cY^{\mult, \ddagger, d^H_1=d^H_2=1}_{\Iw(v_1)}} U_0^{=_{v_2}
  2 }$  contains a  point which is non-ordinary at $v_1$. It  finally
follows that any component of $U_0 \times_{\cY^{\mult, \ddagger,
    d^H_1=d^H_2=1}_{\Iw(v_1)}} U_0$ contains a  point which is
non-ordinary at $v_1$, as required.
\end{proof}

We can now complete the proof of Theorem~\ref{thm: main thm on gluing}. Consider  the diagram:
$$\begin{diagram}
\cY^{\mult, \ddagger,d^H_1\ge 1-\epsilon,d^H_2\ge 1-\epsilon}_{\Iw(v_1,v_2)} & \rTo^{\quad q_{v_1} \quad } & \cY^{\mult, \ddagger,d^H_1\ge 1-\epsilon,d^H_2\ge 1-\epsilon}_{\Iw(v_2)}  \\
\dTo_{q_{v_2}} & & \dTo_{q_{v_2}} \\
\cY^{\mult, \ddagger,d^H_1\ge 1-\epsilon,d^H_2\ge 1-\epsilon}_{\Iw(v_1)} & \rTo^{q_{v_1}} & \cY^{\mult, \ddagger,d^H_1\ge 1-\epsilon,d^H_2\ge 1-\epsilon} \end{diagram}$$

By Lemma~ \ref{lemma-dificil}, we have proved that all elements of our
spaces $V_{\alpha_1 + \beta_1, \alpha_2}$ and $V_{\alpha_1 + \beta_1,
  \beta_2}$ are sections on the whole $\cY^{\mult, \ddagger,d^H_1\ge
  1-\epsilon,d^H_2\ge 1-\epsilon}_{\Iw(v_2)}$ and that all elements of
our spaces $V_{\alpha_1 , \alpha_2+ \beta_2}$ and $V_{ \beta_1, \alpha_2+ \beta_2}$ are sections on the whole $\cY^{\mult, \ddagger,d^H_1\ge 1-\epsilon,d^H_2\ge 1-\epsilon}_{\Iw(v_1)}$. We can pull back these sections to $\cY^{\mult, \ddagger,d^H_1\ge 1-\epsilon,d^H_2\ge 1-\epsilon}_{\Iw(v_1,v_2)}$. 

The isomorphism~(\ref{eqn: Z1-alpha1
  iso}) and the similar
isomorphism $(U_{v_1,1}-\beta_1):V_{\alpha_1+\beta_1,\beta_2}\isoto
V_{\alpha_1,\beta_2}$ induce an
isomorphism \[(U_{v_1,1}-\beta_1):V_{\alpha_1+\beta_1,\alpha_2}\oplus
  V_{\alpha_1+\beta_1,\beta_2}\isoto V_{\alpha_1,\alpha_2}\oplus
  V_{\alpha_1,\beta_2},\] and we
define~$V_{\alpha_1+\beta_1,\alpha_2+\beta_2}$ to be the preimage of
$V_{\alpha_1,\alpha_2+\beta_2}\subset V_{\alpha_1,\alpha_2}\oplus
  V_{\alpha_1,\beta_2}$ under this isomorphism. This is an
  $8d_{\rho}$-dimensional space of eigenforms with the appropriate eigenvalues,
  so we only need to check that all of the elements
  of~$V_{\alpha_1+\beta_1,\alpha_2+\beta_2}$ descend
  to~$\cY^{\mult, \ddagger,d^H_1\ge 1-\epsilon,d^H_2\ge 1-\epsilon}$.

Consider an element~$F$ of this space. By definition,
$F$ has the property that~$(U_{v_1,1}-\beta_1)F$ on~$\cY^{\mult, \ddagger,d^H_1\ge 1-\epsilon,d^H_2\ge 1-\epsilon}_{\Iw(v_1,v_2)}$ is pulled back from~$\cY^{\mult, \ddagger,d^H_1\ge 1-\epsilon,d^H_2\ge 1-\epsilon}_{\Iw(v_1)}$ via~$q^*_{v_2}$. Let~$G = \deg(q_{v_2})^{-1} q_{v_2,*} F$ be the trace of~$F$ to~$\cY^{\mult, \ddagger, d^H_1\ge 1-\epsilon,d^H_2\ge 1-\epsilon}$.
The form~$F$
comes via pullback from~$\cY^{\mult, \ddagger,d^H_1\ge
  1-\epsilon,d^H_2\ge 1-\epsilon}$ if and only if~$F =  q^{*}_{v_2}
G$. Since the trace map at~$v_2$ 
commutes with~$U_{v_1,1}$ (for the usual reasons, ultimately coming
down to Serre--Tate theory and the product structure on the
$p$-divisible group), we deduce (since $q_{v_1}$ is surjective) that~$(U_{v_1,1}-\beta_1)(q^{*}_{v_2}
G-F)=0$, so that~$F = q^{*}_{v_2}G$  (because $(U_{v_1,1}-\beta_1)$ is injective) as required.\end{proof}

We conclude this section with some lemmas that were used above. We first record  the following easy lemma:

\begin{lem}\label{lem-con-fin-et}  If~$S \rightarrow T$ is a finite
  \'{e}tale map of adic spaces of finite type over a field, then the image of any connected component of~$S$ is a connected
component of~$T$. \end{lem}
\begin{proof}Since~$S$ and~$T$ are of finite type, they have only finitely
  many connected components. In particular the connected
   components of $S$ and $T$ are precisely the connected subsets  of $S$ and $T$ which are both open and closed.
  Since finite \'etale morphisms are both open and closed~\cite[Lem.\ 1.4.5,
  Prop.\ 1.7.8]{MR1734903}, the
  result is immediate.
\end{proof}

Next we have the following lemma and its corollary:

\begin{lem}\label{lem: more connected components} Any connected component of $Y^{I, =_{v_2}
    1}_{K^{p, \Iw}K_p(I),1}$ contains a point in $Y^{I,  =_{\{v_1,v_2\}}
    1}_{{K^{p, \Iw} K_p(I)},1}$, and any connected component of $Y^{I, =_{v_2} 2}_{K^{p, \Iw}K_p(I),1}$ contains a point in $Y^{I,  =_{v_1}1, =_{v_2} 2}_{{K_p(I) K^p},1}$.
\end{lem}

\begin{cor}\label{cor-connected-components} Any connected component of
  either of the spaces $$\cY_{\Iw(v_1)}^{\mult, \ddagger,
    d^H_1=d^H_2=1,=_{v_2} 1},\  \cY_{\Iw(v_1)}^{\mult, \ddagger, d^H_1=d^H_2=1, =_{v_2} 2 }$$
  contains a  point which is non-ordinary at $v_1$.
\end{cor}

\begin{proof} The map $\cY_{\Iw(v_1)}^{\mult, \ddagger,
    d^H_1=d^H_2=1,=_{v_2} 1} \rightarrow \cY^{\mult, \ddagger,
    d^H_1=d^H_2=1,=_{v_2} 1}$
    is finite \'etale, so it suffices to prove the claims for
  $\cY^{\mult, \ddagger, d^H_1=d^H_2=1,=_{v_2} 1}$ and $\cY^{\mult,
    \ddagger, d^H_1=d^H_2=1, =_{v_2} 2 }$. Also, the map $\cY^{\mult,
    \ddagger,  d^H_1=d^H_2=1} \rightarrow \cY^{\mult }$ induces an
  isomorphism of $\pi_0$'s, because both spaces have the same rank one
  points and any higher rank point admits a generalization to a rank one
  point. Thus $\cY^{\mult, =_{v_2} 1}$  is the tube of $Y^{I,
    =_{v_2} 1}_{K^{p, \Iw} K_p(I),1}$ and $\cY^{\mult, =_{v_2} 2}$  is the
  tube of $Y^{I, =_{v_2} 2}_{K^{p, \Iw} K_p(I),1}$. Since all these spaces are
  smooth, the tube of a connected component in $Y^{I, =_{v_2}
    i}_{K^{p, \Iw} K_p(I),1}$ is connected for $i=1,2$.  But now 
    by Lemma~\ref{lem: more connected components} these components
    contain points which have rank one at~$v_1$ and  hence are not ordinary at~$v_1$.
\end{proof}

The rest of this section is devoted to proving Lemma \ref{lem: more connected components}.  This statement is a very special case of a general expectation that ``all possible specializations between EKOR strata are realized.''  Unfortunately, as far as we are aware, this exact statement does not yet appear in the literature, but we will explain how it can be deduced from what is available using standard techniques.  This will necessitate a small digression into the theory of stratifications of special fibres of Shimura varieties.

To aid the reader's understanding, we first recall a general strategy
for producing specializations between strata: first one produces a
specialization to a point of a very special stratum, and then one uses
deformation theory at that special point to ``go back up'' to the
desired stratum.  For achieving the first step, there is also a
standard strategy: if one can show that open strata are
(quasi)-affine, while the closures of strata are proper, then it
follows that any component of any stratum must specialize to a point of
a zero dimensional stratum.  This argument becomes a bit more
complicated for non compact Shimura varieties, where one must study
the extension of the stratification to the boundary of the minimal
compactification.  In order to use results readily available in the
literature, we will carry out the first step at spherical level, then
carry out the second step at Iwahori level, and finally explain how
this implies the result that we want at level $K_p(I)$.

First we consider the Ekedahl--Oort stratification at spherical level, see for instance \cite{Viehmann2013}.  Let $K_p = \prod_{v|p} \mathrm{GSp}_4(\ocal_{F_v})$.  Then $Y_{K^{p,Iw}K_p,1}$ has an Ekedahl--Oort stratification into $4^{[F:\Q]}$ strata, according to the four possibilities for each of the finite flat group schemes $\cG_w[p]$ at geometric points.  Let $G_{1,1}=E[p]$ for $E$ a supersingular elliptic curve.  Then these four possibilities are:
\begin{itemize}
\item Ordinary: $\cG_w[p]\simeq \mu_p^2\times(\Z_p/p\Z_p)^2$
\item $p$-rank 1: $\cG_w[p]\simeq \mu_p\times\Z_p/p\Z_p\times G_{1,1}$
\item Supergeneral: $\cG_w[p]$ is connected-connected, but not isomorphic to $G_{1,1}^2$.
\item Superspecial: $\cG_w[p]\simeq G_{1,1}^2$.
\end{itemize}
This stratification refines the $p$-rank stratification, with the last two cases corresponding to $\cG_w$ having $p$-rank 0.  We call a point of $Y_{K^{p,Iw}K_p,1}$ superspecial if $\cG_w[p]$ is superspecial for all $w| p$.  This is the unique zero dimensional stratum.

\begin{lem}\label{lem-superspecial-point} Let $J\subseteq S_p$.  Each irreducible component of $Y_{K^{p,\Iw}K_p,1}^{=_{J^c}2,=_J1}$ contains a point of $Y_{K^{p,\Iw}K_p,1}^{=_{S_p}0}$ in its closure.
\end{lem}

\begin{proof}
It is shown in \cite{2015arXiv150705922B,2015arXiv150705032G} that the Ekedahl--Oort stratification extends to a stratification of the minimal compactification of $Y_{K^{p,\Iw}K_p,1}$, and that each (open) stratum is affine.  Moreover the superspecial locus does not intersect the boundary.  It follows that any component of any Ekedahl--Oort stratum contains a superspecial point in its closure.  By the explicit description of the Ekedahl--Oort stratification recalled above, the $p$-rank strata in the statement of the lemma are also Ekedahl--Oort strata.
\end{proof}

Now we will switch to Iwahori level and consider the Kottwitz--Rapoport stratification, see for instance \cite{MR1952527}.  Let $K_p^{\Iw}=\prod_{v| p}\Iw(v)$.  Then $Y_{K^{p,\Iw}K_p^{\Iw},1}$ and its local model $M^{\loc}_{K_p^{\Iw},1}=\prod_{v| p} M^{\loc}_{\Iw(v),1}$ carry a Kottwitz--Rapoport stratification.  In fact, there is a Kottwitz--Rapoport stratification of $M^{\loc}_{\Iw(v),1}$ and the stratification of $M^{\loc}_{K_p^{\Iw},1}$ is simply the product stratification.  The strata of $M^{\loc}_{\Iw(v)}$ are indexed by a set $\mathrm{Adm}(\mu)$ of cardinality 13.  This set, as well as the partial ordering given by closure, is pictured in \cite[p.\ 1273]{MR2419384}.

We will use below the following argument, which is a consequence of the theory of local models.  If $C$ is an irreducible component of the Kottwitz--Rapoport stratum labeled by $w\in\mathrm{Adm}(\mu)^{S_p}$, then the closure $\overline{C}$ has a decomposition into strata:
$$\overline{C}=\coprod_{w'\leq w}\overline{C}_{w'}.$$
A priori the strata $\overline{C}_{w'}$ might be empty, although it is expected that they are always nonempty.  However, the theory of local models implies that if $\overline{C}_{w'}$ is nonempty, then so is $\overline{C}_{w''}$ for any $w''$ satisfying $w'\leq w''\leq w$.

We will not need to recall in detail the definition of the Kottwitz--Rapoport stratification.  We do recall that, as explained in \cite{MR2419384}, the Kottwitz--Rapoport invariant determines whether the groups of order $p$, $H_w$ and $L_w/H_w$, are \'etale, multiplicative, or connected-connected (and so in particular the Kottwitz--Rapoport invariant determines the $p$-rank of $\cG_w$, a theorem of Genestier--Ng\^o).  Conversely these invariants determine the Kottwitz--Rapoport invariant when the $p$-rank of $\cG_w$ is not 0.  All of this is recorded in the table in \cite[p.\ 1276]{MR2419384}.

We will use the following points:
\begin{itemize}
\item There is a Kottwitz--Rapoport condition, $s_2s_1s_2\tau$ in \cite{MR2419384}, which corresponds to the condition that $\cG_w$ is ordinary and $L_w=\cG[F]$ (equivalently $L_w$ is multiplicative).
\item There is a Kottwitz--Rapoport condition, $s_1s_2\tau$ in \cite{MR2419384}, which corresponds to the condition that $\cG_w$ has $p$-rank 1, $H_w$ is multiplicative, and $L_w=\cG[F]$ (equivalently $H_w$ is multiplicative and $L_w/H_w$ is connected-connected).
\item There are three Kottwitz--Rapoport conditions, $\tau,s_1\tau$, and $s_2\tau$ in \cite{MR2419384}, which have $p$-rank 0 and are in the closure of the first stratum recalled above.  We observe crucially that they are also all in the closure of the second stratum recalled above.  We refer to these three strata as the canonical $p$-rank 0 Kottwitz--Rapoport strata (here ``canonical'' refers to the fact that $L_w=\cG_w[F]$ is the canonical subgroup of $\cG_w$).
\end{itemize}

For $J\subseteq S_p$ we write $Y_{K^{p,\Iw}K_p^{\Iw},1}^{=_{J^c}2,=_J1,m-can}$ for the locus in $Y_{K^{p,\Iw}K_p^{\Iw},1}$ where for $w\in J^c$, $\cG_w$ is ordinary and $L_w=\cG_w[F]$, while for $w\in J$, $\cG_w$ has $p$-rank 1, $H_w$ is multiplicative, and $L_w=\cG_w[F]$.  By what we have just recalled, this is a Kottwitz--Rapoport stratum.

\begin{lem}\label{lem-connected-components-iwahori}
For $J\subseteq J'\subseteq S_p$, any irreducible component of $Y_{K^{p,\Iw}K_p^{\Iw},1}^{=_{J^c}2,=_J1,m-can}$ contains a point of $Y_{K^{p,\Iw}K_p^{\Iw},1}^{=_{(J')^c}2,=_{J'}1,m-can}$ in its closure.
\end{lem}

\begin{proof}
Let $\pi:Y_{K^{p,\Iw}K_p^{\Iw},1}\to Y_{K^{p,\Iw}K_p,1}$ be the projection from Iwahori to spherical level.  It is proper, and the Kottwitz--Rapoport stratum $Y_{K^{p,\Iw}K_p^{\Iw},1}^{=_{J^c}2,=_J1,m-can}$ maps finitely onto the $p$-rank stratum $Y_{K^{p,\Iw}K_p,1}^{=_{J^c}2,=_J1}$ (the fibres correspond to the $p+1$ choices of $H_w$ for $w\in J^c$).  If $C$ is an irreducible component of $Y_{K^{p,\Iw}K_p^{\Iw},1}^{=_{J^c}2,=_J1,m-can}$, then $\pi(C)$ is an irreducible component of $Y_{K^{p,\Iw}K_p,1}^{=_{J^c}2,=_J1}$.  By Lemma \ref{lem-superspecial-point}, the closure $\overline{\pi(C)}$ contains a point which is $p$-rank 0 for all $w\in S_p$.  By the properness of $\pi$ it follows that the closure $\overline{C}$ contains a point which is $p$-rank 0 for all $w\in S_p$.

By what we have shown, in the closure $\overline{C}$, at least one of the canonical $p$-rank 0 Kottwitz--Rapoport strata is nonempty.  Now we apply the argument with local models and the explicit description of the closure relations between the strata recalled above to conclude.
\end{proof}

\begin{rem}
One could give a more direct proof of Lemma \ref{lem-connected-components-iwahori}, avoiding the consideration of the Ekedahl--Oort stratification and the superspecial locus at spherical level, if one knew that the Kottwitz--Rapoport stratification of $Y_{K^{p,\Iw}K_p^\Iw}$ extended to a stratification of the minimal compactification, for which the (open) strata are quasi-affine.  However we lack a reference for these facts.
\end{rem}

\begin{proof}[Proof of Lemma \ref{lem: more connected components}]
Let $\pi:Y_{K^{p,\Iw}K_{p}^\Iw,1}\to Y_{K^{p,\Iw}K_p(I),1}$ be the projection from Iwahori to $K_p(I)$ level.  On $Y_{K^{p,\Iw}K_p(I),1}$, $\pi$ has a ``canonical section'' $s:Y_{K^{p,\Iw}K_p(I),1}\to Y_{K^{p,\Iw}K_{p}^\Iw,1}$, defined by taking $L_w=\cG_w[F]$ for $w\in I$ (recall that on $Y_{K^{p,\Iw}K_p(I),1}$, $H_w$ is multiplicative by definition, and hence $H_w\subseteq\cG_w[F]$).  It follows that for $J\subseteq I$, $s$ and $\pi$ define mutually inverse isomorphisms between $Y_{K^{p,\Iw}K_p^{\Iw},1}^{=_{J^c}2,=_J1,m-can}$ and $Y_{K^{p,\Iw}K_p(I),1}^{I,=_{J^c}2,=_J1}$.  We deduce the following statement from Lemma \ref{lem-connected-components-iwahori}: for $J\subseteq J'\subseteq I$, every irreducible component of $Y_{K^{p,\Iw}K_p(I),1}^{I,=_{J^c}2,=_J1}$ contains a point of $Y_{K^{p,\Iw}K_p(I),1}^{I,=_{(J')^c}2,=_{J'}1}$ in its closure.

Applying this with $J=\{v_2\}$, $J'=\{v_1,v_2\}$ and $J=\emptyset$, $J'=\{v_1\}$ we conclude that any irreducible component of $Y_{K^{p,\Iw}K_p(I),1}^{=_{S_p\setminus \{v_2\}}2,=_{v_2}1}$ contains a point of $Y_{K^{p,\Iw}K_p(I),1}^{I,=_{v_1,v_2}1}$ in its closure, and any irreducible component of $Y_{K^{p,\Iw}K_p(I),1}^{=_{S_p}2}$ contains a point of $Y_{K^{p,\Iw}K_p(I),1}^{I,=_{v_1}1,=_{v_2}2}$ in its closure.  Finally as recalled at the start of section \ref{subsec: Hasse invariants and stratifications}, $Y_{K^{p,\Iw}K_p(I),1}^{=_{S_p\setminus \{v_2\}}2,=_{v_2}1}$ is dense in $Y_{K^{p,\Iw}K_p(I),1}^{=_{v_2}1}$ and $Y_{K^{p,\Iw}K_p(I),1}^{=_{S_p}2}$ is dense in $Y_{K^{p,\Iw}K_p(I),1}^{=_{v_2}2}$, and the lemma follows.
\end{proof}

\subsection{Solvable base change}
  
  We will use solvable
base change to deduce our main modularity lifting theorem from
Corollary~\ref{cor: rho is automorphic}. We firstly prove a couple of
preparatory lemmas, beginning with the following well-known result.

\begin{lemma} \label{lemma:simnotunique} Let~$K$ be a number field, and let $\rho: G_K \rightarrow \GL_4(\Qbar_p)$ be an irreducible representation which preserves a generalized symplectic form with similitude character~$\nu$. Then either~$\nu$ is uniquely determined by~$\rho$, or, if~$\rho$ also admits a similitude character~$\nu \psi$ with~$\psi \ne 1$,
then~$\psi$ has finite order and~$\rho$ is reducible over a quadratic subfield of the fixed field of~$\psi$ and hence also over
the fixed field of~$\psi$.
\end{lemma}

\begin{proof} Let~$V$ denote the underlying representation of~$\rho$, and let~$\nu$ and~$\nu \psi$
denote two possible similitude characters. Then there is an inclusion~$\nu \oplus \nu \psi
\subset \Hom(V^*,V)$, or equivalently, $1 \oplus \psi \subset \Hom(V^*(\nu),V)$. It follows
that~$V \simeq V^*(\nu)$ and~$V \simeq V^*(\nu \psi)$, and thus~$V \simeq V(\psi)$, and also~$V \simeq V(\psi) \simeq V(\psi^2)$. By comparing
determinants, it follows that~$\psi^4$ is trivial, and hence either~$\psi$ or~$\psi^2$ is a quadratic character~$\eta$ such that~$V \simeq V(\eta)$
and hence~$1 \oplus \eta \subset \Hom(V,V)$. 
By Schur's Lemma, $V$ becomes reducible over the fixed field of~$\eta$, which by construction is a quadratic subfield of the fixed field of~$\psi$.
\end{proof}

 We now prove a slightly
technical lemma on solvable base change; it is an analogue
of~\cite[Lem.\ 1.3]{blght} for~$\GSp_4$, but the proof is slightly more involved.

\begin{lem}
  \label{lem: solvable descent for Siegel modular forms}Suppose that
  $p>2$ splits completely in the totally real field~$F/\Q$. Let $F'/F$ be a solvable extension of
  totally real fields. Suppose that $\rho:G_F\to\GSp_4(\Qpbar)$ satisfies:
  \begin{enumerate}
  \item $\nu\circ\rho=\varepsilon^{-1}$.
  \item For all $v|p$, $\rho|_{G_{F_v}}$ is $p$-distinguished weight~$2$ ordinary.
  \item The representation~$\rhobar$ is vast and tidy.
  \item $\rho|_{G_{F'}}$ is irreducible. Furthermore, there is an ordinary
    automorphic representation~$\pi'$ of $\GSp_4(\A_{F'})$ of parallel
    weight~$2$ and central character~$|\cdot|^2$, such for every
    finite place~$w$ of~$F'$ we
    have
    \[\WD(\rho|_{G_{F'_w}})^{F-\semis}\cong\recGTp(\pi'_w\otimes|\nu|^{-3/2}).\]
    \emph{(}So in particular, $\rho_{\pi',p}\cong\rho|_{G_{F'}}$.\emph{)}
\end{enumerate}
Then~$\rho$ is modular. More precisely, there is an ordinary automorphic
representation~$\pi$ of $\GSp_4(\A_F)$ of parallel weight~$2$ and
central character~$|\cdot|^2$, with
$\rho_{\pi,p}\cong\rho$. Furthermore, for every finite place~$v$
of~$F$ we
have
\[\WD(\rho|_{G_{F_v}})^{F-\semis}\cong\recGTp(\pi_v\otimes|\nu|^{-3/2}).\]
\end{lem}
\begin{proof}
Since $\rho_{\pi',p}$ is irreducible, $\pi'$ must be of general type
in the sense of~\cite{MR2058604}, so that it corresponds to a cuspidal
automorphic representation~$\Pi'$ of~$\GL_4(\A_{F'})$. By induction we may reduce to the
  case that~$F'/F$ is cyclic of prime degree, in which case it follows
  from~\cite[Thm.\ 4.2 of \S3]{MR1007299} that there is an automorphic
  representation~$\Pi$ of~$\GL_4(\A_F)$
  with~$\BC_{F'/F}(\Pi)=\Pi'$.

We can
  write
  \[L^S(s,\Pi',\bigwedge^2\otimes |\cdot|^{-2})=\prod_\psi
    L^S(s,\Pi,\bigwedge^2\otimes |\cdot|^{-2}\psi^{-1})\]where the product
  is over the characters~$\psi$ of
  $\A_F^\times/F^\times N_{F'/F}\A_{F'}^\times$. The left hand side
  has a simple pole at~$s=1$ (by the assumption that~$\Pi'$ is the
  transfer of~$\pi'$), while  by the main result of~\cite{MR1610812}, all
  but at most one factor on the right hand side is holomorphic and
  non-vanishing at $s=1$. Thus some factor on the right hand side
  must also have a simple pole at~$s=1$, say
  $L^S(s,\Pi,\bigwedge^2\otimes |\cdot|^{-2}\psi^{-1})$.

  It follows from Theorem~\ref{thm: arthur classification results}
  that ~$\Pi$ is the transfer of a cuspidal automorphic
  representation~$\pi$ of $\GSp_4(\A_F)$ with central
  character~$|\cdot|^2\psi$. Since $\BC_{F'/F}(\Pi)=\Pi'$, we see
  that~$\pi$ is of parallel
  weight~$2$. Letting~$\rho_{\pi,p}:G_F\to\GSp_4(\Qpbar)$ be the
  Galois representation corresponding to~$\pi$ (whose existence
  follows from~\cite[Thm.\ 3.5]{MR3200667} exactly as in the proof
  of Theorem~\ref{thm: existence and properties of Galois representations for
    automorphic representations}), we have
  $\rho_{\pi,p}|_{G_{F'}}\cong\rho|_{G_{F'}}$, so that (since $F'/F$ is
  cyclic of prime degree, and $\rho|_{G_{F'}}$ is irreducible)
  $\rho_{\pi,p}$ differs from~$\rho$ by a twist by a character
  of~$\Gal(F'/F)$. 

Replacing~$\pi$ by the corresponding twist, we may assume
that~$\rho_{\pi,p}$ and~$\rho$ are isomorphic when considered as
representations valued in~$\GL_4(\Qpbar)$. We claim that we
necessarily have
  $\nu\circ\rho=\nu\circ\rho_{\pi,p}=\varepsilon^{-1}$, so that
  $\pi$ has central character~$|\cdot|^2$. Indeed, this follows from
  Lemma~\ref{lemma:simnotunique}, since it holds after
  restriction to~$G_{F'}$, and~$\rho|_{G_{F'}}$ is irreducible by
  assumption. So~$\rho_{\pi,p}\cong\rho$, as required.  

Since we have assumed that~$\rhobar$ is vast and tidy, it follows from Corollary~\ref{cor: Galois representations for limits of discrete series
    with ordinarity at p} that for every finite place~$v$
  of~$F$ we
  have \[\WD(\rho|_{G_{F_v}})^{\semis}\cong\recGTp(\pi_v\otimes|\nu|^{-3/2})^{\semis}.\]
  It remains to check that the monodromy operators agree; but this may
  be checked after base change, and
  since~$\pi'$ is the base change of~$\pi$, it follows from the
  assumption that
  $\WD(\rho|_{G_{F'_w}})^{F-\semis}\cong\recGTp(\pi'_w\otimes|\nu|^{-3/2})$. 
\end{proof}

\subsection{The main modularity lifting
  theorem}\label{subsec: main modularity lifting}

We now prove our main modularity lifting theorem.\begin{thm}\label{thm: final R equals T}Suppose that~$p\ge 3$ splits completely in
  the totally real field~$F/\Q$. Suppose that $\rho:G_F\to\GSp_4(\Qpbar)$ satisfies:
  \begin{enumerate}
  \item $\nu\circ\rho=\varepsilon^{-1}$.
    \item The representation~$\rhobar$ is vast and tidy in the sense of Definitions~\ref{defn:vast} and~\ref{defn:tidy}.
  \item For all $v|p$, $\rho|_{G_{F_v}}$ is $p$-distinguished
     weight~$2$ ordinary in the sense of Definition~\ref{defn: generic flat ordinary}. \item There exists $\pi$ of parallel weight~$2$ and central character~$|\cdot|^2$, which is 
    ordinary at all~$v|p$, such
    that~$\rhobar_{\pi,p}\cong\rhobar$.  \item For all finite places~$v$ of~$F$, $\rho|_{G_{F_v}}$
   and~$\rho_{\pi,p}|_{G_{F_v}}$ are pure.\end{enumerate}Then~$\rho$ is modular.   More
  precisely, there is an ordinary automorphic
  representation~$\pi'$ of $\GSp_4(\A_F)$ of parallel weight~$2$ and
  central character~$|\cdot|^2$ which satisfies
  $\rho_{\pi',p}\cong\rho$. Furthermore, for every finite
  place~$v$ of~$F$ we have \[\WD(\rho|_{G_{F_v}})^{F-\semis}\cong\recGTp(\pi'_v\otimes|\nu|^{-3/2}).\]
\end{thm}
\begin{proof}Choose a solvable extension of totally real
  fields~$F'/F$, linearly disjoint from~$F^{\ker\rhobar}$
  over~$F$, with the following properties:
  \begin{itemize}
  \item $p$ splits completely in~$F'$.
  \item At every place~$w$ of~$F'$ lying over a place~$v\nmid p$
    of~$F$ for which~$\pi_v$ or~$\rho|_{G_{F_v}}$ is ramified,
    $\rhobar|_{G_{F'_w}}$ is trivial, $\rho|_{G_{F'_w}}$ has only
    unipotent ramification, and $q_w\equiv 1\pmod{p^2}$.
  \item There is an automorphic representation~$\pi'$
    of~$\GSp_4(\A_{F'})$ of parallel weight~$2$ which is a base change
    of~$\pi$ (in the sense that for each finite place~$w$ of~$F'$,
    lying over a place~$v$ of~$F$, we have
    $\recGTp(\pi')=\recGTp(\pi)|_{W_{F'_w}}$). Furthermore, for all
    finite places~$w$ of~$F'$ we have $(\pi'_w)^{\Iw(w)}\ne 0$.
  \end{itemize}
(The last property can be arranged by~\cite[Prop.\ 4.13]{MR3200667}.) Then $\rho|_{G_{F'}}$ satisfies Hypothesis~\ref{hyp: hypotheses on
  rho}, so the result follows from Corollary~\ref{cor: rho is
  automorphic} (applied to~$\rho|_{G_{F'}}$) and Lemma~\ref{lem: solvable descent for Siegel modular
  forms}.
\end{proof}

\subsection{Base change and automorphy lifting}\label{subsec: base
  change to relax central character}

Throughout the paper, we have fixed the similitude factor of our
Galois representations to be~$\varepsilon^{-1}$, in order to
streamline both the notation and some arguments. We now explain how to
use base change to relax this condition in our main automorphy lifting
theorem. We do not use this result elsewhere in the paper, so we have
contented ourselves with a slightly ugly statement, and with a sketch
of the proof.

\begin{defn}
  \label{defn: twisted ordinary}We say that a
  representation~$\rho:G_{\Qp}\to\GSp_4(\Qpbar)$ is \emph{twisted $p$-distinguished
     weight~$2$ ordinary} if it is an unramified twist of a
   representation which is $p$-distinguished
     weight~$2$ ordinary in the sense of Definition~\ref{defn: generic
       flat ordinary}. Similarly, we say that an admissible
     representation~$\pi_p$ of~$\GSp_4(\Qp)$ is \emph{twisted
       ordinary} if it is an unramified twist of an
     ordinary representation.
\end{defn}
 \begin{thm}\label{thm: final R equals T relaxed central character}Suppose that~$p\ge 3$ splits completely in
  the totally real field~$F/\Q$. Suppose that $\rho:G_F\to\GSp_4(\Qpbar)$ satisfies:
  \begin{enumerate}
  \item $\nu\circ\rho=\chi\varepsilon^{-1}$, where~$\chi$ is a totally
    even finite order character, which is unramified at all places dividing~$p$.
    \item The representation~$\rhobar$ is vast and tidy in the sense of Definitions~\ref{defn:vast} and~\ref{defn:tidy}.
  \item For all $v|p$, $\rho|_{G_{F_v}}$ is twisted $p$-distinguished
     weight~$2$ ordinary.  \item There exists $\pi$ of parallel weight~$2$, which is  twisted
    ordinary at all~$v|p$, such that~$\rhobar_{\pi,p}\cong\rhobar$.
  \item For all finite places~$v$ of~$F$, $\rho|_{G_{F_v}}$ and~$\rho_{\pi,p}|_{G_{F_v}}$ are pure.\end{enumerate}Then~$\rho$ is modular.   More
  precisely, there is a twisted ordinary automorphic
  representation~$\pi'$ of $\GSp_4(\A_F)$ of parallel weight~$2$ which satisfies
  $\rho_{\pi',p}\cong\rho$. Furthermore, for every finite
  place~$v$ of~$F$ we have \[\WD(\rho|_{G_{F_v}})^{F-\semis}\cong\recGTp(\pi'_v\otimes|\nu|^{-3/2}).\]
\end{thm}
\begin{proof}Let~$\chi'$ be the finite order character
  $G_F\to\Qpbar^\times$ such that
  $\nu\circ\rho_{\pi}=\chi'\varepsilon^{-1}$. Note that~$\chi'$ is
  totally even (since we have~$\chibar'=\chibar$ by
  assumption). We can choose a quadratic extension of totally real
  fields~$F'/F$, linearly disjoint from~$F^{\ker\rhobar}$
  over~$F$, such that:
  \begin{itemize}
  \item $p$ splits completely in~$F'$, and
  \item there are finite order characters
    $\psi,\psi':G_{F'}\to\Qpbar^\times$ such that $\chi|_{G_{F'}}=\psi^2$, $\chi'|_{G_{F'}}=(\psi')^2$.
  \end{itemize}
  Indeed, the obstruction to taking the square root of a character is
  in the 2-torsion of the Brauer group, and there are no obstructions
  to taking a square root of either~$\chi$ or~$\chi'$ at the places
  dividing~$p$ (because both characters are unramified at such places)
  or at the infinite places (because~$\chi,\chi'$ are totally
  even). 

Let~$\pi_{F'}$ be the base change of~$\pi$ to~$F'$. Since
  $\rho|_{G_{F'}}\otimes\psi^{-1}$, $\pi_{F'}\otimes(\psi')^{-1}\circ\Art_{F'}$
  satisfy the hypotheses of Theorem~\ref{thm: final R equals T}, it
  follows that~$\rho|_{G_{F'}}\otimes\psi$ is modular,
  so~$\rho|_{G_F'}$ itself is modular. The result follows from
  Lemma~\ref{lem: solvable descent for Siegel modular forms} (or
  rather, from an obvious generalization of this lemma to the case of
  more general central characters, which may be proved in the same
  way).
\end{proof}

\section{Potential modularity of abelian surfaces}\label{section:
  potential modularity of abelian surfaces}

We now use the potential automorphy methods introduced
in~\cite{MR1954941} to prove the potential modularity of abelian
surfaces. It is presumably possible to follow~\cite[\S1]{MR1954941} quite
closely, but we instead make use of potential modularity results
for~$\GL_2$ and the local to global principle of~\cite[\S3]{frankII}
(see also~\cite[Thm.~1.2]{MBdiditfirst}).

\subsection{Compatible systems and potential automorphy}\label{subsec: potential automorphy}Recall that the notion of a $C$-algebraic automorphic representation
is defined in~\cite{MR3444225}, and in the case of automorphic representations
of~$\GL_n$, this definition agrees with the notion of an algebraic
automorphic representation defined in~\cite{MR1044819}.
\begin{defn}
  \label{defn: automorphic compatible system general number
    field}Let~$K$ be a number field and let~$\cR$ be a strictly
  compatible system of representations of~$G_K$. We say that~$\cR$ is
  \emph{automorphic} if there is an automorphic
representation~$\Pi$ of~$\GL_n(\A_K)$, with the properties that: 
 \begin{enumerate}
\item \label{condn: isobaric algebraic} $\Pi$ is an isobaric direct sum of cuspidal automorphic representations $\boxplus_{i=1}^r\Pi_i$ where
each~$\Pi_i$ is a  $C$-algebraic cuspidal
automorphic 
representation of some~$\GL_{n_i}(\A_K)$. 
\item \label{clozeltotherescue}The fixed field~$M_\Pi$ of the subgroup
of~$\Aut(\C)$ consisting of those~$\sigma\in\Aut(\C)$ with
${}^\sigma\Pi^\infty\cong\Pi^\infty$ is a number field.
\item \label{localglobalfiniteplaces} For each finite place~$v$ of~$K$, $\WD_v(\cR)=\oplus_i\rec(\Pi_{i,v} |\det|_v^{(1-n_i)/2})$.
\end{enumerate}
\end{defn} 

\begin{rem}
  \label{rem: for us coefficient fields are indeed number
    fields}There are many (conjecturally equivalent) variants of Definition~\ref{defn: automorphic compatible system general number
    field} that could be made. The definition is in some sense
  redundant, because condition~\eqref{clozeltotherescue} is implied by
  condition~\eqref{localglobalfiniteplaces}; indeed, by the definition
  of a compatible system, it follows that for all but finitely
  many~$v$, $\Pi_v$ is an unramified principal series representation,
  defined over a number field which may be chosen independently
  of~$v$. Condition~\eqref{clozeltotherescue} then follows from strong
  multiplicity one for isobaric representations~\cite{MR623137}. The reason that we
  have chosen to include the condition separately is that conjecturally (see~\cite{MR1044819})
  condition~\eqref{condn: isobaric
    algebraic} implies condition~\eqref{clozeltotherescue}, and also
  implies the existence of a compatible system~$\cR$ satisfying condition~\eqref{localglobalfiniteplaces}.

In fact, the
  only cases of Definition~\ref{defn: automorphic compatible system general number
    field} that we will need to consider are those where either:
    \begin{enumerate}
    \item
  Each~$\Pi_i$  is regular algebraic, or
  \item\label{item: transferred from GSp4} $K$ is totally real, $\Pi$ is
cuspidal, and~$\Pi$ is the transfer to~$\GL_4$ of a cuspidal
automorphic representation of~$\GSp_4$ of  parallel weight~$2$ and central
  character~$|\cdot|^2$.
  \end{enumerate}
   In either case, condition~(\ref{clozeltotherescue}) is satisfied
 by~\cite[Thm.\
3.13]{MR1044819} and~\cite[Thm.\
  3.2.2]{blasius-harris-ramak}) respectively.
\end{rem}
\begin{rem}
  \label{rem: why we didn't demand local-global compatibility at
    infinity}The reader may wonder why we did not demand an analogue
  of condition~\eqref{localglobalfiniteplaces} of Definition~\ref{defn: automorphic compatible system general number
    field} at the infinite places. One reason is that we do not need
  to do so, as condition~\eqref{localglobalfiniteplaces} already
  determined~$\Pi$ uniquely (indeed, as in Remark~\ref{rem: for us coefficient fields are indeed number
    fields}, this is already true if one only considers
  condition~\eqref{localglobalfiniteplaces} at all but finitely many
  places). The main reason that we do not make a requirement at the
  infinite places is that (in keeping with the literature) our
  definition of a compatible system does not include a requirement
  that the $l$-adic representations are compatible on complex
  conjugations, which makes it harder to formulate a precise
  compatibility. One could certainly ask (as in~\cite{MR3444225}) that the Hodge--Tate weights
  of~$\cR$ correspond to the infinitesimal character of~$\Pi$, but to
  save introducing additional notation and terminology we have not
  done so.
\end{rem}
\begin{rem}
  As explained in Remark~\ref{rem: for us coefficient fields are indeed number
    fields}, condition ~\eqref{localglobalfiniteplaces} of Definition~\ref{defn: automorphic compatible system general number
    field} at all but finitely many places~$v$ determines~$\Pi$
  uniquely. One might ask whether if this condition holds for all but
  finitely many~$v$, it necessarily holds for all~$v$. In general this
  is a hard problem; indeed even if~\eqref{localglobalfiniteplaces} is
  known up to semisimplification,  it is often difficult to show
  that the monodromy operators agree. If however~$\cR$ is pure and~$\Pi$
  is generic then the agreement of monodromy operators is automatic;
  we will use this fact in our arguments below. 
\end{rem}

If~$\cR$ and~$\Pi$ are as in Definition~\ref{defn: automorphic compatible system general number
    field}, we as usual have Gamma factors $L_v(\Pi,s)$ for each
  place~$v|\infty$ of~$K$, and we set \[ \Lambda_{\Pi}(\cR,s)=L(\cR,s)
  \prod_{v|\infty} L_v(\Pi,s).  \] This is of course just the usual completed
$L$-function of~$\Pi$, but we have included~$\cR$ in the notation to
emphasize that the~$L$-functions of~$\Pi$ and~$\cR$ agree -- note that
here it is important that we know Definition~\ref{defn: automorphic compatible system general number
    field}~\eqref{localglobalfiniteplaces} at all finite places, and
  not just at almost all places, or up to semisimplification. As noted above, since we do not \emph{a priori}
  demand any local-global compatibility at~$\infty$ for our compatible system~$\cR$, we use the automorphic representation~$\Pi$
  in this definition mostly as a convenient way to write down the correct Gamma factors at infinity.
For those who find this notation unpleasant, note that --- in the
restrictive context of abelian surfaces over totally real fields
 --- we defined a function~$\Lambda(\cR,s)$ in~\eqref{eqn: completed L function
  motivic version} by explicitly writing down the Gamma factors in question, and then (for
 all the~$A$ and~$\Pi$ that arise in this paper) we indeed have equalities~$\Lambda(\Pi,s) = \Lambda_{\Pi}(\cR_A,s) = \Lambda(\cR_A,s)$.

  We also have an epsilon
factor $\varepsilon(\Pi)$, and a conductor~$N(\Pi)$, and by~\cite[Cor.\ 13.8]{MR0342495}, $\Lambda_\Pi(\cR,s)$ admits a meromorphic continuation to the
entire complex plane, and satisfies the functional equation
\numequation\label{eqn: functional equation for automorphic
  representation} \Lambda_\Pi{}(\cR,s)=\varepsilon(\Pi)N(\Pi)^{-s}
\Lambda_{\Pi}(\cR^\vee, 1-s). \end{equation} 

\begin{defn}
  Let~$A/K$ be an abelian variety. We say that~$A$ is
  \emph{automorphic} if~$\cR_A$ is automorphic in the sense of
  Definition~\ref{defn: automorphic compatible system general number
    field}. We say that it is \emph{potentially automorphic} if there
  is a finite extension of number fields~$L/K$ such that~$\cR_A|_{G_L}$
  is automorphic.
\end{defn}

\begin{rem}\label{rem: automorphic Gamma epsilon N factors agree}
  If~$A/K$ is an abelian variety, and~$A$ is automorphic with the corresponding~$\Pi$ being of the form considered in Remark~\ref{rem:
    for us coefficient fields are indeed number fields}, then the
  Gamma factors~$L_v(\Pi,s)$ 
  and the conductor~$N(\Pi)$  agree with those defined for the
  compatible system~$\cR_A$ in \S\ref{subsec: compatible
    systems}. Indeed, for the Gamma factors this is a direct consequence of
  the definitions, and the conductor respects the
  local Langlands correspondence. In particular, we
  have~$\Lambda_\Pi(\cR_A,s)=\Lambda(\cR_A,s)$, and the functional
  equations~\eqref{eqn: functional equation in automorphic case} and~\eqref{eqn: functional equation for automorphic
  representation} agree; so~$\Lambda(\cR_A,s)$  satisfies the expected
meromorphic continuation and functional equation.
\end{rem}

    \begin{defn}\label{defn: potential modularity for a single
        representation}Let~$F$ be a totally real. We
      say that a representation $\rho:G_F\to\GSp_4(\Qpbar)$ is
      \emph{modular} if there is a cuspidal automorphic representation
      of~$\GSp_4(\A_F)$ of parallel weight~$2$ and central
      character~$|\cdot|^2$ for each finite place~$v$ of~$F$ we have
      $\WD(\rho_{\pi,p}|_{G_{F_v}})^{F-\semis}\cong
      \recGTp(\pi_v\otimes|\nu|^{-3/2})$; in particular,  $\rho\cong\rho_{\pi,p}$.

We say that~$r$ is \emph{potentially modular} if there is a finite
Galois extension $F'/F$ of totally real fields such that~$r|_{G_{F'}}$
is modular.

If~$A/F$ is an abelian surface, we say that~$A$ is modular (resp.\
potentially modular) if~$\rho_{A,p}$ is modular (resp.\ potentially
modular) for some (equivalently, for any) prime~$p$.
    \end{defn}
    \begin{rem}
      \label{rem: relationship between modular and automorphic}The
      relationship between the definitions of what it means for an
      abelian surface~$A/F$ to be (potentially) automorphic or modular is somewhat complicated,
      because of the various possibilities in Arthur's classification
      of the discrete spectrum of~$\GSp_4(\A_F)$. In this paper we
      will only show that~$A$ is (potentially) modular if the
      corresponding automorphic representation of~$\GSp_4$ is of
      general type, in which case~$A$ is also (potentially)
      automorphic, essentially by the definition of ``general type'';
      note that
      this is the case considered in Remark~\ref{rem: for us coefficient fields are
        indeed number fields}~\eqref{item: transferred from GSp4}.
    \end{rem}

We now prove some technical lemmas that we will use in proving our
main potential automorphy/modularity results. A weakly compatible system $\CR$ is defined to be {\em irreducible}
if there is a set $\CL$ of rational primes of Dirichlet density $1$
such that for $\lambda|l \in \CL$ the representation $r_\lambda$ is
irreducible. We say that it is \emph{strongly irreducible} if for all
finite extensions $F'/F$ the compatible system $\CR|_{G_{F'}}$ is
irreducible. If~$n=2$, then we say that~$\cR$ has weight~$0$ if
$H_\tau(\cR)=\{0,1\}$ for each~$\tau$, and we say that~$\cR$ is odd if
$\det\cR(c_v)=-1$ for all $v|\infty$. If~$\boldsymbol{\pi}$ is a cuspidal
automorphic representation of~$\GL_2(\A_F)$ of weight~$0$, then~$\cR(\boldsymbol{\pi})$ is odd and
has weight~$0$.

We have the following standard lemma.

\begin{lem}\label{lem: rank 2 irreducible equivalences} Let~$\cR$ be a
  rank two weakly compatible system.
  \begin{enumerate}
  \item   The following are
  equivalent:
  \begin{enumerate}
  \item $\CR$ is irreducible.
  \item For all $\lambda$ the representation $r_\lambda$ is
    irreducible.
  \item For some $\lambda$ the representation $r_\lambda$ is
    irreducible.
  \end{enumerate}
\item If $\CR$ is irreducible and
  regular, then the following are equivalent:
  \begin{enumerate}
  \item  $\CR$ is strongly irreducible.
  \item  $\Sym^2\CR$ is irreducible.
  \item For all~$\lambda$, $\Sym^2 r_\lambda$ is irreducible.
  \item For some~$\lambda$, $\Sym^2 r_\lambda$ is irreducible.
\end{enumerate}
If these equivalent conditions do not hold, then there is a quadratic
extension $F'/F$ and a weakly \emph{(}equivalently, strongly\emph{)} compatible system $\CX$ of characters of
$G_{F'}$ such that
\[ \CR \cong \Ind_{G_{F'}}^{G_F} \CX. \]
\item If $\CR$ is strongly irreducible and regular, then
  for a density one set of primes~$l$ of~$\Q$, if $\lambda|l$ is a
  place of~$M$, then the image of $\barr_\lambda$ contains
  $\SL_2(\CO_M/\lambda)$.
\end{enumerate}
\end{lem}
\begin{proof}This is well known. Part~(1) is~\cite[Lem.\ 7.1.1]{10author}, and part~(3)
  is~\cite[Lem.\ 7.1.3]{10author}. For part~(2), note that
  by~\cite[Lem.\ 7.1.2]{10author}, either~$\CR$ is strongly
  irreducible, or we can write~$\CR \cong \Ind_{G_{F'}}^{G_F} \CX$. It
  follows that if~$\CR$ is not strongly irreducible, then $\Sym^2
  r_\lambda$ is reducible for every~$\lambda$. 

  Conversely, if $\Sym^2 r_\lambda$ is reducible for some~$\lambda$,
  then there is a nontrivial character~$\psi$ such that
  $r_\lambda\cong r_\lambda\otimes\psi$. Considering determinants,
  $\psi$ is a quadratic character. Letting $F'/F$ be the quadratic
  extension corresponding to~$\psi$, it follows from Schur's lemma
  that~$r_\lambda|_{G_{F'}}$ is reducible, so by part~(1), $\CR$ is
  not strongly irreducible.
\end{proof}

We now use a standard trick with restriction of scalars to give some
slight improvements to some applications of the theorem of Moret-Bailly.

\begin{prop}
  \label{prop: restriction of scalars Snowden}Let $F_1/F$ be a finite
  extension of totally real fields, and let~$p,q>2$ be distinct primes which split
  completely in~$F_1$. Let ~$\rbar:G_{F_1}\to\GL_2(\Fqbar)$ be a
  representation with determinant~$\varepsilonbar^{-1}$.

  Suppose that for each place~$v|q$ of~$F_1$, $\rbar|_{G_{F_{1,v}}}$
  is of the form
$  \begin{pmatrix}
    \lambda_{\alphabar_v}&0\\0&\varepsilon^{-1}\lambda_{\alphabar_v}^{-1}
  \end{pmatrix}$. Suppose also that ~$\rbar$ is unramified at all places
  above~$p$. 

  Let~$\Favoid/F$ be a finite extension. Then there is a finite Galois
  extension $F'/F$ of totally real fields in which~$p$ and~$q$ split
  completely and which is  linearly disjoint from
  $F_1\Favoid/F$, and a $q$-ordinary cuspidal automorphic representation~$\boldsymbol{\pi}$
  of~$\GL_2(\A_{F_1F'})$ of weight~$0$ and trivial central character
  which is unramified at all places dividing~$pq$ and which satisfies
  $\rhobar_{\boldsymbol{\pi},q}\cong\rbar|_{G_{F_1F'}}$.
\end{prop}
\begin{proof}
  In the case $F_1=F$, this is a straightforward consequence
  of~\cite[Thm.\ 8.2.1]{0905.4266}. Indeed, in Snowden's notation, we
  take~$\rhobar=\rbar^\vee$, $\psi=1$, we let~$S$ consist of the
  places dividing~$pq$ and we let~$t$ assign the type~$A$ at places
  lying over~$q$ and type~$AB$ at places lying over~$p$. To prove the general case, one simply replaces the
  scheme~$X$ to which Snowden applies the theorem of Moret-Bailly
  with the restriction of scalars~$\Res_{F_1/F}X_{F_1}$.
\end{proof}

\begin{prop}\label{prop: restriction of scalars of Calegari lifting lemma}Let
  $G$ be a finite group, let $E/\Q$ be a finite extension, and let~$S$
  be a finite set of places of~$E$. Let~$E'/E$ be a finite extension,
  and let $\Favoid/E$ be a finite extension, linearly disjoint
  from~$E'/E$. 

  Let $S'/S$ be the set of places of~$E'$ lying over places
  of~$S$. For each finite place~$v\in S'$, let $H'_v/E'_v$ be a finite
  Galois extension together with a fixed
  inclusion~$\phi_v:\Gal(H'_v/E'_v)\into G$ with image~$D_v$. For each
  real infinite place~$v\in S'$, let~$c_v\in G$  be an element of order
  dividing~$2$. 

Then there exists a number field~$K/E$ and a finite Galois extension
of number fields~$L/K$ such that if we set $K'=KE'$, $L'=LE'$, then 
  \begin{enumerate}
  \item There is an isomorphism~$\Gal(L'/K')=G$.
  \item $L'/E$ is linearly disjoint from~$E'\Favoid/E$.
  \item All places in~$S'$ split completely in~$K'$.
  \item For all finite places~$w$ of~$K'$ above~$v\in S'$, the local
    extension~$L'_w/K'_w$ is equal to~$H'_v/E'_v$. Moreover, there is a commutative diagram:\[ 
\begin{tikzcd} \Gal(L'_w/K'_v)\ar{r}\ar[equals]{d}& D_w\subset
  G\ar[equals]{d}\\
\Gal(H'_v/E'_v)\ar{r}{\phi_v}&D_v\subset G
\end{tikzcd}
\]
\item For all real places~$w|\infty$ of~$K'$ above~$v\in S'$, complex
  conjugation~$c_w\in G$ is conjugate to~$c_v$.
  \end{enumerate}

\end{prop}
\begin{proof}
  The case $E'=E$ is~\cite[Prop.\ 3.2]{frankII} (see also~\cite[Thm.~1.2]{MBdiditfirst}). The general case may
  be proved in exactly the same way, by replacing the application
  of~\cite[Thm.\ 3.1]{frankII} to (an open subscheme of) $X_G/E$ with
  an application of it to~$\Res_{E'/E}X_G$.
\end{proof}

\subsection{Abelian surfaces}\label{subsec:abelian surfaces}
We begin by recalling some results from~\cite{MR2982436}
and~\cite{MR3660222}, which allow us to deal with various
cases where the abelian surfaces have extra endomorphisms, and can be handled with the potential automorphy
theorems of~\cite{BLGGT}. Let $A/F$ be an abelian surface over a
totally real field~$F$,  and let $L/F$ be the minimal extension over
which all its endomorphisms are defined. This is a Galois extension,
and following~\cite{MR2982436} we say that the \emph{Galois type}
of~$A/F$  is the $\Gal(L/F)$-module $\End_L(A)\otimes_\Z\R$. These
possible Galois types are classified in~\cite{MR2982436}, and they are
divided up into 6 families \textbf{A}-\textbf{F}.

The precise classification of monodromy groups in these references is
not actually strictly necessary for our purposes. Write~$\{\rho_{A,l}\}$ for the compatible system of Galois
representations $\{H^1(A_{\overline{F}},\Ql)\}$. In practice, it suffices to know that the $l$-adic representations~$\rho_{A,l}$ fall into  precisely one of the following categories independently of~$l$:
\begin{enumerate}
\item strongly irreducible (type~$\mathbf{A}$),
\item reducible (type~$\mathbf{B}[C_1]$, $\mathbf{C}$, $\mathbf{E}[C_n]$, some~$\mathbf{D}$, some~$\mathbf{F}$),
\item potentially abelian but not reducible (of type the remaining~$\mathbf{D}$ and~$\mathbf{F}$ cases),
\item induced from a quadratic extension~$K/F$ but not potentially abelian, in which case either:
\begin{enumerate}
\item the two~$2$-dimensional representations over~$K$ are equivalent
  up to twist
(type~$\mathbf{E}[D_n]$), or
\item the two~$2$-dimensional representations over~$K$ are not
 equivalent up to twist (type~$\mathbf{B}[C_2]$).
\end{enumerate}
\end{enumerate}

\begin{prop}\label{prop: a bunch of existing cases of potential
    modularity}Suppose that~$A/F$ is not of type~$\mathbf{A}$
  or~$\mathbf{B}[C_2]$. 
  Then~$A$ is potentially automorphic.\end{prop}
\begin{proof}We freely use the discussion of~\cite[\S
4]{MR3660222}. In cases~\textbf{D}, \textbf{F}, $A_L$ is
of CM type, so the compatible system~$\cR_A$ is potentially
abelian, and in particular potentially automorphic.

In cases~$\mathbf{B}[C_1]$, $\mathbf{C}$, and the cases of
type~$\mathbf{E}$ other than those of type~$\mathbf{E}[D_n]$, it
follows from the discussions at the beginnings of~\cite[\S 4.2, 4.4, 4.5]{MR3660222} that we can
write~$\cR_A=\cR_A^1\oplus\cR_A^2$ where each~$\cR_A^i$ is an
irreducible, odd, weight~$0$ weakly compatible system of rank 2
$l$-adic representations of~$G_F$. The potential automorphy of~$\cR_A$
therefore follows from~\cite[Thm.\ 5.4.1]{BLGGT}.

It remains to treat the case that $A$ is of type~$\mathbf{E}[D_n]$. In
this case, as explained in~\cite[\S
4.5-4.6]{MR3660222}, there is a quadratic extension $F'/F$
and a strongly irreducible weakly compatible system~$\cS=\{s_l\}$ of
weight~$0$ representations of~$G_{F'}$
which is defined over~$\Q$~ such that
$\cR_A\cong\Ind_{G_{F'}}^{G_F}\cS$. Furthermore, there is a
finite order character~$\delta$ of~$G_{F'}$ such that if we write
$\Gal(F'/F)=\{1,\sigma\}$, then $s_l^\sigma\cong
s_l\otimes\delta_l$.

It follows from Lemma~\ref{lem: rank 2 irreducible equivalences}~(2)
and~\cite[Prop.\ 5.3.2]{BLGGT} that for a density one set of
primes~$l$, $\Sym^2 \sbar_l|_{G_{F'(\zeta_l)}}$ is irreducible, $l$ is
unramified in~$F'$, and both~$s_l$ and~$\varepsilon_l$ are crystalline at
all primes above~$l$. Fix one such~$l>7$.

Since $\Proj s_l^\sigma\cong \Proj s_l$, it follows from Schur's lemma
that $\Proj s_l$ extends to a representation
$G_F\to\PGL_2(\Qlbar)$. By~\cite[Lem.\ 2.3.17,
2.7.4]{Patrikis}, we may lift this to a representation
$\tilde{r}:G_F\to \GL_2(\Qlbar)$ which is unramified at all but
finitely many places, and is Hodge--Tate at all places dividing~$l$,
with Hodge--Tate weights $(0,1)$. By construction, there is a
character~$\psi:G_{F'}\to\Qlbartimes$ such
that~$\tilde{r}|_{G_{F'}}=s_l\otimes\psi$. Since~$\psi$ is Hodge--Tate
of weight~$0$, it has finite order.

Since $l$ is unramified in~$F'$ and $s_l$ is crystalline at all primes
above~$l$, after possibly replacing ~$\tilde{r}$ by a twist by a finite order
character, we may assume that it is crystalline at all places
dividing~$l$. By~\cite[Prop.\ 2.5]{MR3186511}, $\tilde{r}$ is odd, so
by~\cite[Thm.\ 4.5.1]{BLGGT}, $\tilde{r}$ is potentially
automorphic. Since \[\rho_{A,l}\cong\Ind_{G_{F'}}^{G_F}s_l \cong
  \Ind_{G_{F'}}^{G_F}\bigl(\tilde{r}|_{G_{F'}}\otimes\psi^{-1}\bigr),\]it
follows that~$\cR_A$ is potentially automorphic, as required.
\end{proof}

We say that~$A/F$ is \emph{challenging} if it has type~\textbf{A}
(which is the case that $\End_{\C}A=\Z$) or~$\mathbf{B}[C_2]$. In
the latter case, as explained in~\cite[\S 4.3]{MR3660222},
there is a quadratic extension~$K/F$, and a strongly irreducible
weakly compatible system~$\cS=\{s_l\}$ of rank~$2$, weight~$0$
representations of~$G_K$ with determinant~$\varepsilon_l^{-1}$ such that
$\cR_A\cong\Ind_{G_K}^{G_F}\cS$. Furthermore,
writing~$\Gal(K/F)=\{1,\sigma\}$, $s_l^\sigma$ and $s_l$ do not become
isomorphic after restriction to any finite extension of~$K$.
(The case when~$K/F$ is totally real can be handled using potential
automorphy theorems for~$\GL_2$, but our argument 
(at this point at least) does not
need to distinguish between the various infinity types of~$K$.)

\begin{lem}\label{lem: big image for abelian surfaces}If~$A/F$ is a
  challenging abelian surface, then for a density one set of
  primes~$l$, $\rhobar_{A,l}$ is vast and tidy.
  \end{lem}
  \begin{proof}
    If~$\End_{\C}A=\Z$, then, for all sufficiently large~$l$,
    $\rhobar_{A,l}(G_F)=\GSp_4(\Fl)$ by~\cite{serretovigneras}, so~
    the claim follows from  Lemma~\ref{lem: GSp4 enormous}. 

    If~$A$ is of type~$\mathbf{B}[C_2]$, then writing $\rho_{A,l}\cong\Ind_{G_K}^{G_F}s_l$, we see from Lemma~\ref{lem:
      Ind GL2 enormous} and Lemma~\ref{lem: rank 2 irreducible
      equivalences}~(3) that we need only check that for a density one
    set of primes~$l$, we have $\Proj \sbar_l^\sigma\not\equiv\Proj
    \sbar_l$. (Note that the inverse of the mod~$l$ cyclotomic character is
    surjective for all ~$l$ which are unramified in~$F$.) To see this, note that since $\Proj s_l^\sigma\not\equiv\Proj
    s_l$, we have $\Sym^2 s_l^\sigma\not\equiv\Sym^2 s_l$
    by~\cite[Appendix, Thm.\ B]{MR1728875}. There is therefore some
    finite place~$v$ of~$F$ at which the compatible systems  $\{\Sym^2
    s_l^\sigma\}$, $\{\Sym^2 s_l\}$ are unramified, for which the
    eigenvalues of~$\Frob_v$ differ for the two compatible
    systems. Then the same applies for $\Sym^2
    \sbar_l^\sigma$, $\Sym^2 \sbar_l$ for all sufficiently large~$l$,
    so that in particular $\Proj \sbar_l^\sigma\not\equiv\Proj
    \sbar_l$, as required.
   \end{proof}

\begin{defn}\label{defn: good primes for abelian surface}
  Let $A/F$ be an abelian surface over a totally real field. We say
  that a rational prime~$p\ge 3$ is a \emph{good prime for $A$} if:
  \begin{itemize}
  \item $A$ admits a polarization of degree prime to~$p$.
  \item $p$ splits completely in~$F$.
  \item The representation~$\rhobar_{A,p}$ is vast and tidy.
    \item For each place~$v|p$, $\rho_{A,p}|_{G_{F_v}}$ is $p$-distinguished weight~$2$ ordinary.
  \end{itemize}

\end{defn}

\begin{rem}
  The point of Definition~\ref{defn: good primes for abelian surface}
  is that the good primes~$p$ are the ones for which we can apply our
modularity lifting theorem (Theorem~\ref{thm: final R equals T}) to~$\rho_{A,p}$.
\end{rem}

\begin{lem}\label{lem: existence of good primes}Let~$A/F$ be a challenging abelian surface.
Then the set of rational primes which are
  good primes for~$A$ has relative density one in the set of primes which split completely in~$F$.
\end{lem}
\begin{proof}By Lemma~\ref{lem: big image for abelian surfaces}, it
  suffices to show that $\rhobar_{A,p}|_{G_{F_v}}$ is
  $p$-distinguished weight~$2$ ordinary for a density one set of
  finite places~$v$
  of~$F$ (with residue characteristic~$p$). To do this, we follow the approaches
  of~\cite{MR3494322} and~\cite[Lem.\ A.7]{CGGSp4}. Consider the places~$v$ of~$F$ that are
  split over a prime~$p$ of~$\Q$, for which~$A$ has good reduction;
  the set of such primes has density one. Fix a prime~$l\ne p$. The characteristic
  polynomial of~$\rho_{A,l}(\Frob_v)$ is of the
  form \[x^4-a_1x^3+a_2x^2-pa_1x+p^2\]where $a_1,a_2$ are integers.

  Then~$A$ has good ordinary reduction at~$v$ if and only if
  $p\nmid a_2$. If this holds, then we see that
  $\rhobar_{A,p}|_{G_{F_v}}$  will be  $p$-distinguished weight~$2$ ordinary if and only if $a_1^2-4a_2$ is not divisible by~$p$. By the Weil
  bounds, we have $|a_1|\le 4\sqrt{p}$, $|a_2-2p|\le 4p$, so if~$a_1^2-4a_2$
  is divisible by~$p$, then it is equal to~$pc$  for~$c$ in some
  finite list of integers, independent of~$p$.

  Let~$G$ be the Zariski closure of~$\rho_{A,l}(G_F)$ in~$\GSp_4$, and
  write~$V$ for the standard representation of~$\GSp_4$, and~$\chi$
  for the similitude character. Arguing exactly as in the proof
  of~\cite[Thm.\ 1]{MR3494322}, it follows from the
  Cebotarev density theorem that it is enough to show that  the
  virtual representation $( V^{\otimes 2}-4\wedge^2 V)\otimes\chi^{-1}$ does not have constant trace on any
  connected component of~$G$. 

  By the proof of~\cite[Thm.\ 3]{MR3494322}, we can
  replace~$G$ by the Sato--Tate group of~$A$, which is
  either the connected group~$\mathrm{USp}_4$ (if~$A$ has type~$\mathbf{A}$), or the
  normalizer of~$\mathrm{SU}_2\times\mathrm{SU}_2$
  in~$\mathrm{USp}_4$ (if~$A$ has type~$\mathbf{B}[C_2])$) (which has
  two connected components). The result now follows easily from an
  explicit check.
\end{proof}

\begin{lem}
  \label{lem: existence of auxiliary prime q}Let~$A/F$ be a challenging abelian
  surface. Then there are distinct rational primes $p$, $q$ such
  that  $p$ and~$q$ are both good primes for~$A$, and  for all places~$v|p$ of~$F$, $\rhobar_{A,q}(\Frob_v)$ has distinct
    eigenvalues.

\end{lem}
\begin{proof}By Lemma~\ref{lem: existence of good primes}, a density one subset
  of the set of rational
  primes which split completely in~$F$ are good primes. Let~$p$ be any good prime for~$A$; then, for each place $v|p$ of~$F$,
  $\rho_{A,p}|_{G_{F_v}}$ has $p$-distinguished weight~$2$ ordinary reduction, and in
  particular the eigenvalues of the crystalline Frobenius~$\Frob_v$
  on~$T_pA$ are distinct. Consequently, for all but finitely many
  rational primes~$q$ of good reduction for~$A$,
  $\rhobar_{A,q}(\Frob_v)$ has distinct eigenvalues for all
  places~$v|p$. \end{proof}

If~$q$ is a good prime for an abelian surface~$A/F$, then for each
place~$w|q$ of~$F$ we may write \[\rhobar_{A,q}|_{G_{F_w}}\cong
  \begin{pmatrix}
    \lambda_{\alphabar_w}&0&*&*\\
0 &\lambda_{\betabar_w}&*&*\\
0&0&\varepsilonbar^{-1}\lambda_{\betabar_w}^{-1}&0\\
0&0&0&\varepsilonbar^{-1}\lambda_{\alphabar_w}^{-1}
  \end{pmatrix}
\]Then we write \[(\rhobar_{A,q}|_{G_{F_w}})^\semis:=
  \begin{pmatrix}
    \lambda_{\alphabar_w}&0&0&0\\
0 &\lambda_{\betabar_w}&0&0\\
0&0&\varepsilonbar^{-1}\lambda_{\betabar_w}^{-1}&0\\
0&0&0&\varepsilonbar^{-1}\lambda_{\alphabar_w}^{-1}
  \end{pmatrix}
\]
When reading the proofs of the following two results, it may be
helpful to recall that our convention is that the
representation~$\rhobar_{A,p}$ is the dual of~$A[p]$; this accounts
for the various duals occurring in the proofs.
\begin{lem}
  \label{lem: global realization of GL2 representation}Let $A/F$ be an
  abelian surface over a totally real field, and let~$p$, $q$ be
  primes as in Lemma~\ref{lem: existence of auxiliary prime q}. Fix a
  totally real quadratic extension $F_1/F$ in which~$p$ and~$q$ split
  completely, and which is linearly disjoint from the kernels of the
  actions of~$G_F$ on~$A[p]$ and~$A[q]$.

  Then there is a finite Galois extension of totally real
  fields~$F'/F$, and a representation
  $\rbar_q:G_{F'F_1}\to\GL_2(\Fq)$, with the following
  properties: \begin{enumerate}
  \item $p$ and~$q$ both split completely in~$F'$.
  \item $F'/F$ is linearly disjoint from~$F_1/F$ and from the kernels of the
  actions of~$G_F$ on~$A[p]$ and~$A[q]$.
\item $\det\rbar_q=\varepsilonbar^{-1}$.
\item $\rbar_q(G_{F'F_1})=\GL_2(\Fq)$, and the projective image
  of~$\rbar_q$ is not equal to its conjugate
  under~$\Gal(F'F_1/F')$.\item Set
  $\rhobar_q:=\Ind_{G_{F'F_1}}^{G_{F'}}\rbar_q:G_{F'}\to\GSp_4(\Fq)$
  with similitude factor~$\varepsilonbar^{-1}$. Then
   \begin{itemize}
   \item for any place~$w|q$ of~$F$ and any place~$w'|w$ of~$F'$,
     $\rhobar_q|_{G_{F'_{w'}}}\cong (\rhobar_{A,q}|_{G_{F_w}})^\semis$, and
   \item for any place~$v|p$ of~$F$ and any place~$v'|v$ of~$F'$,
     $\rhobar_q|_{G_{F'_{v'}}}\cong \rhobar_{A,q}|_{G_{F_v}}.$ 
   \end{itemize}
   \item  The representation~$\rhobar_q$ is vast and tidy.
  \end{enumerate}

\end{lem}
\begin{proof}Fix a finite place~$\mathfrak{r}$ of~$F$ not
  dividing~$pq$ and splitting in~$F_1$. We apply Proposition~\ref{prop: restriction of
    scalars of Calegari lifting lemma}, taking $E=F$, $E'=F_1$,
  $G=\GL_2(\Fq)$, $S$ to be the set of places
  dividing~$pq\mathfrak{r}\infty$, and~$\Favoid$ to be the extension
  cut out by the intersection of the kernels of $\rhobar_{A,p}$
  and~$\rhobar_{A,q}$. For each infinite place~$v\in S'$ we
  choose~$c_v$ to have eigenvalues~$\{1,-1\}$. For each place $w\in S$
  dividing~$pq\mathfrak{r}$ we write $w=w_1w_2$ for its decomposition
  in~$F_1$. If $w|p$, then the eigenvalues
  of~$(A[q]^\vee|_{G_{F_w}})(\Frob_w)$ can be written as $\alpha_w$,
  $\beta_w$, $p\beta_w^{-1}$, $p\alpha_w^{-1}$,  while
  if~$w|q$ we use the notation above. In either case, we
  choose~$\phibar_{w_1}$ to correspond to the representation~$
  \begin{pmatrix}
    \lambda_{\alpha_w}&0\\ 0&\varepsilon^{-1}\lambda_{\alpha_w}^{-1}
  \end{pmatrix}
$, and~$\phibar_{w_2}$ to
correspond to~$
  \begin{pmatrix}
    \lambda_{\beta_w}&0\\ 0&\varepsilon^{-1}\lambda_{\beta_w}^{-1}
  \end{pmatrix}
$. Finally, if $w=\mathfrak{r}$, then we choose $\phibar_{w_1}$,
$\phibar_{w_2}$ to have determinant~$\varepsilonbar^{-1}$, in such a way
that~$\phibar_{w_1}$ is unramified, while $\Proj\phibar_{w_2}$ is ramified. 

We obtain an extension~$F'/F$ (the extension~$K/E$ from
Proposition~\ref{prop: restriction of scalars of Calegari lifting
  lemma}, with the ~$\phi_v$ there being our~$\phibar_v$) and a representation $\rbar_q:G_{F_1F'}\to\GL_2(\Fq)$ which
satisfies~(1), and (2). It need not satisfy~(3), but by
construction~$\varepsilonbar\det\rbar_q$ is an even character which is
trivial at all places dividing~$pq\mathfrak{r}$. The obstruction to the existence
of a square root of~$\varepsilonbar\det\rbar_q$ is therefore a class in the
2-torsion of~$\Br_{F_1F'}$ which is trivial at all places
dividing~$pq\mathfrak{r}\infty$. 

We can therefore replace~$F'$ by a quadratic totally real extension in
which~$p$, $q$, $\mathfrak{r}$ split completely, and assume that~$\varepsilonbar\det\rbar_q$
has a square root. By~\cite[Ch.\ X, Thm.\ 5]{MR2467155} we can (by replacing
this square root by a twist by a quadratic character) arrange that
the square root is trivial at all places dividing~$pq\mathfrak{r}$.
 Replacing~$\rbar_q$ by its twist by this square root, we ensure~(3),
at which point~(5) follows (note that for each place~$v|p$ of~$F$,
$\rhobar_{A,q}|_{G_{F_v}}$ is unramified with distinct eigenvalues
of~$\Frob_v$, and is therefore semisimple). Considering the places
lying over~$\mathfrak{r}$, we see that~(4) is satisfied. Finally, (6)
then follows from Lemma~\ref{lem: Ind GL2 enormous}.
\end{proof}

\begin{thm}
  \label{thm: MB GSp4 application}Let~$A/F$ be a challenging abelian
  surface over a totally real field. Then~$A$ is potentially
  modular. More precisely, there is a finite Galois extension of
  totally real fields~$F'/F$ and a prime~$p$ splitting completely in~$F'$ such
  that~$\rho_{A,p}|_{G_{F'}}$ is modular and irreducible.
\end{thm}
\begin{proof}Let~$p$, $q$, $F_1$, $F'$, $\rbar_q$ and~$\rhobar_q$ be as in
  Lemma~\ref{lem: global realization of GL2 representation}. 
  Let~$Y/F'$ denote the moduli space of
  triples~$(B,\imath_p,\imath_q)$ consisting of abelian surfaces~$B$
and symplectic isomorphisms \[\imath_p:B[p]\isoto
  A[p]|_{G_{F'}},\] \[\imath_q:B[q]\isoto \rhobar_q^\vee.\] This is
smooth and geometrically
connected. (Over either~$\C$ or~$\Qbar$, we may identify~$Y$ with the moduli space
of principally polarized abelian surfaces with full level~$pq$ structure.)

We claim that for each place~$v|pq\infty$ of~$F'$,  the
subspace~$\Omega_v:=Y^{{\ord}}(F'_v) \subset Y(F'_v)$ consisting
of points corresponding to abelian surfaces with good ordinary reduction (when~$v$ is finite) is
nonempty. If $v|\infty$, this follows from~$\det r_q^\vee=\varepsilon$,
while if~$v|p$, then~$A$ itself gives a point of~$Y(F'_v)$ (by
point~(5) of Lemma~\ref{lem: global realization of GL2
  representation}). Finally, if~$v|q$, the canonical
lift of~$A$ modulo~$v$ gives a point of~$Y(F'_v)$.  (Since~$A$ has good  ordinary reduction at~$v|p$ and~$v|q$, the
corresponding point on~$Y$ does indeed land in~$\Omega_v$.)

By~\cite[Prop.\ 3.1.1]{BLGGT} (a theorem of Moret-Bailly), we may
find a finite Galois totally real extension~$F''/F'$ in which~$p$
and~$q$ split completely, and which is linearly disjoint from the
compositum of $F_1F'$ and the kernels of the actions of~$G_{F'}$
on~$A[p]$, $A[q]$ and~$\rhobar_q$, with the property
that~$Y(F'') \cap \bigcap_{v|pq}  \Omega_v  \ne 0$. Let~$B/F''$ be a corresponding abelian
surface, which by construction will have good ordinary reduction for all~$v|p$ and~$v|q$.

By Proposition~\ref{prop: restriction of scalars Snowden}, after
replacing~$F''/F'$ with a further totally real extension, we can
maintain all of the above assumptions, and we can further suppose that
there is a ~$q$-ordinary automorphic representation~$\boldsymbol{\pi}$ of~$\GL_2(\A_{F_1F''})$
of weight~$0$ and trivial central character, which is unramified at
all places dividing~$pq$ and which satisfies
$\rhobar_{\boldsymbol{\pi},q}\cong\rbar_q|_{G_{F_1F''}}$. 

It follows from~\cite[Thm.\ 8.6]{MR1871665} that there is an
automorphic representation~$\pi$ of~$\GSp_4(\A_{F''})$ of parallel 
weight~$2$ and trivial central character whose transfer to~$\GL_4(\A_{F''})$ is
the automorphic induction of~$\boldsymbol{\pi}\otimes|\cdot|$, so that in
particular
$\rho_{\pi,q}\cong\Ind_{G_{F_1F''}}^{G_{F''}}\rho_{\boldsymbol{\pi},q}$, so
that $\rhobar_{\pi,q}\cong\rhobar_q|_{G_{F''}}$. In addition, $\pi$ is
ordinary, by construction. The
representation~$\rho_{\pi,p}$ is pure at all finite places
because~$\rho_{\boldsymbol{\pi},q}$ is (for the places away
from~$q$, this is proved in~\cite{MR2327298}, and for the places dividing~$q$ it is
for example a very special case of the main theorem of~\cite{MR3272276}).

We can therefore apply Theorem~\ref{thm: final R equals T} to
$\rho_{B,q}$, and conclude that it is modular. Thus~$\rho_{B,p}$ is
modular, and applying Theorem~\ref{thm: final R equals T} a second
time, we deduce that~$\rho_{A,p}|_{G_{F''}}$ is modular, as
required. (The purity of~$\rho_{B,q}, \rho_{B,p}$ and~$\rho_{A,p}$ at
all finite places is part of Proposition~\ref{prop: pure abelian variety compatible system}.)
\end{proof}

\subsection{Potential modularity and meromorphic
  continuation}\label{subsec: potential modularity and meromorphic continuation}We now
deduce the meromorphic
continuation and functional equation of the $L$-functions associated
to abelian surfaces over totally real fields from our potential
modularity (and automorphy) results.

\begin{thm}
  \label{thm: potential modularity implies meromorphic continuation}
  Let~$F$ be a totally real field, and let~$A/F$ be an abelian
  surface. Then~$\cR_A$ is potentially automorphic, and Conjecture~\ref{conj: L function meromorphic
    continuation} holds for~$A$, for each~$0\le i\le 4$.
\end{thm}
\begin{proof}Since $H^i(A,\Ql)=\wedge^iH^1(A,\Ql)$, it is enough to
  treat the cases $i=1,2$. Note that since for any~$\cR$ we have
  $\varepsilon(\cR)\varepsilon(\cR^\vee)=N(\cR)$
  (see~\cite[(3.4.7)]{MR546607}), and we
  have~$H^i(A,\Ql)^\vee=H^i(A,\Ql)(i)$, the claimed functional
  equation will follow from~(\ref{eqn: functional equation in automorphic
    case}) in the case~$\cR=H^i(A,\Ql)$.

  To see that the meromorphic continuation and the functional
  equation~(\ref{eqn: functional equation in automorphic case}) hold,
  note firstly that if~$A$ has type~$\mathbf{D}$ or~$\mathbf{F}$, then
  the compatible system~$\cR_A$ is potentially abelian, and the result
  follows from a standard argument with Brauer's theorem; more
  precisely, it is immediate from~\cite[Prop.\ 11, Lem.\
  14]{MR3660222}. In the general case, the same argument
  (see e.g.\ the proof of~\cite[Cor.\ 2.2]{MR1954941}) shows that it is
  enough to show that there is a Galois extension of totally real
  fields~$F'/F$ such that for each Galois extension~$F'/F''$ with
  $\Gal(F'/F'')$ solvable, the compatible
  systems~$\cR_A|_{G_{F''}}$ and~$\wedge^2\cR_A|_{G_{F''}}$
  are both automorphic. (Note that the meromorphic continuation and functional
  equations for the compatible systems follow from the functional
  equations~(\ref{eqn: functional equation for automorphic
    representation}) for the corresponding automorphic
  representations.)

  Suppose now that~$A$ has type ~$\mathbf{B}[C_1]$, $\mathbf{C}$, or
  is of type~$\mathbf{E}$ but not of type~$\mathbf{E}[D_n]$. Then as
  we saw in the proof of~\ref{prop: a bunch of existing cases of
    potential modularity}, we can write~$\cR_A=\cR_A^1\oplus\cR_A^2$
  where $\cR_A^1$, $\cR_A^2$ are irreducible, odd, weight~$0$ weakly
  compatible systems of rank 2 $l$-adic representations of~$G_F$. It
  follows from~\cite[Thm.\ 5.4.1]{BLGGT} that there is a Galois
  extension of totally real fields $F'/F$ such that
  $\cR_A^1|_{G_{F'}}$, $\cR_A^2|_{G_{F'}}$ are automorphic and
  irreducible.
It follows from~\cite[Lem.\
  1.3]{blght} that for each
Galois extension
$F'/F''$ with $\Gal(F'/F'')$ solvable,   $\cR_A^1|_{G_{F''}}$,
$\cR_A^2|_{G_{F''}}$ are automorphic. Thus~$\cR_A|_{G_{F''}}$
is automorphic, and since we
have \[\wedge^2\cR_A|_{G_{F''}}=\det\cR_A^1|_{G_{F''}}\oplus
  \det\cR_A^2|_{G_{F''}}\oplus\left(\cR_A^1|_{G_{F''}}\otimes\cR_A^2|_{G_{F''}}\right), \]
it follows from~\cite[Thm.\ M]{MR1792292} that $\wedge^2\cR_A|_{G_{F''}}$ is also
automorphic, as required.

In the remaining cases, namely those of types ~$\mathbf{A}$, $\mathbf{B}[C_2]$, or~$\mathbf{E}[D_n]$, it follows from Theorem~\ref{thm: MB GSp4
  application},  (the proof of)
Proposition~\ref{prop: a bunch of existing cases of potential
  modularity}, Lemma~\ref{lem: getting to general type} and Theorem~\ref{thm: arthur classification
    results}, together with Lemma~\ref{lem: solvable descent for Siegel modular forms}, that there is a Galois
extension of totally real fields $F'/F$ such that for
some~$p$, and all Galois extension $F'/F''$ with $\Gal(F'/F'')$ solvable, $\rho_{A,p}|_{G_{F''}}$ is irreducible and modular (and also automorphic). By the main
result of~\cite{MR2567395} (which is a refinement of the main result
of~\cite{MR1937203}), together with Theorem~\ref{thm: arthur classification
    results}, we see that~$\wedge^2\rho_{A,p}|_{G_{F''}}$ is
  automorphic, as required.
\end{proof}
\begin{rem}
  \label{rem: citing Kim is a bit silly but we do use Arthur}Our use of the results
  of~\cite{MR1937203} and~\cite{MR2567395} in the proof of 
  Theorem~\ref{thm: potential modularity implies meromorphic
    continuation} is almost certainly overkill, and can be  avoided by working with automorphic forms on~$\GSp_4$ rather than~$\GL_4$ as we now explain. 
     Since the four~$4$-dimensional Galois  representations~$H^1(A,\Q_l)$
    are generalized symplectic with respect to the Weil pairing, the exterior square~$\wedge^2 H^1(A,\Q_l) = H^2(A,\Q_l)$ decomposes as
    the direct sum of a $5$-dimensional Galois representation and the one dimensional summand~$\Q_l(-1)$. (The corresponding  Galois
    invariant classes in~$H^2(A,\Q_l(1))$ are generated by the image of a hyperplane section  under the cycle map.)
    Once one knows that the~$4$-dimensional representation~$H^1(A,\Q_l)$
    corresponds (potentially) to an automorphic representation~$\pi$ for~$\GSp_4$, then the~$L$-function associated to the~$5$-dimensional
    summand of~$H^2(A,\Q_l)$ is none other than the degree~$5$ standard  $L$-function, whose analytic properties
    have been known for some time (see~\S6.3 of~\cite{MR892097}).
    On the other hand,  many of our arguments in this paper \emph{do} crucially require passing between~$\GSp_4$ and~$\GL_4$ using Theorem~\ref{thm: arthur classification
    results} and Lemma~\ref{lem: getting to general type}. In particular,   the proof of Theorem~\ref{thm:
    potential modularity implies meromorphic continuation} uses base
  change in the form of Lemma~\ref{lem: solvable descent for Siegel
    modular forms}, and therefore depends directly on Theorem~\ref{thm: arthur classification
    results}; and of course our main modularity lifting theorems also depend on
  these results, in particular to prove that the modules that we patch
  are balanced.
\end{rem}
If~$C/F$ is a curve over a number field, then we can define the
completed $L$-functions $\Lambda_i(C,s)$ and the completed Hasse--Weil
$L$-function~$\Lambda(C,s)$ exactly as for abelian varieties. By
definition we have~$\Lambda_1(C,s)=\Lambda_1(\Jac(C),s)$, where
$\Jac(C)$ is the Jacobian of~$C$.
\begin{cor}\label{cor: genus 2 curves}
  Let $C/F$ be a genus two curve over a totally real field. Then the
  completed Hasse--Weil $L$-function $\Lambda(C,s)$ has a meromorphic
  continuation to the entire complex plane, and satisfies a functional
  equation of the form  $\Lambda(C,s)=\varepsilon N^{-s} \Lambda(C,3-s)$ where~$\varepsilon\in \R$ and~$N\in\Q_{>0}$.
\end{cor}
\begin{proof}
  This follows from Theorem~\ref{thm: potential modularity implies
    meromorphic continuation} with $A=\Jac(C)$.
\end{proof}
Finally, we treat the case of genus one curves over quadratic
extensions of totally real fields.\begin{thm}
  \label{thm: potential modularity of elliptic curves}Let $K/F$ be a
  quadratic extension of a totally real field~$F$, and let~$E/K$ be either a genus one
   curve or an elliptic curve. Then~$E$ is potentially modular. 
More precisely,
  there is a Galois extension of totally real fields $F'/F$ and a
  weight~$0$ cuspidal automorphic representation~$\boldsymbol{\pi}$ of
  $\GL_2(\A_{KF'})$ with trivial central character such that for each
  prime~$l$, we have $\rho_{E,l}|_{G_{KF'}}\cong \rho_{\boldsymbol{\pi},l}$, and in
  fact for each finite place~$v$ of~$KF'$ we have $\WD_v(\rho_{E,l}|_{G_{KF'_v}})^{F-\semis}\cong\rec(\boldsymbol{\pi}_v |\det|_v^{-1/2})$.

  Furthermore,
  Conjecture~\ref{conj: L function meromorphic continuation} holds
  for~$E$.
\end{thm}\begin{proof} We may immediately replace~$E$ by its Jacobian and hence assume that~$E$ is an elliptic curve. 
 If~$E$ is CM, then it is modular, while if~$E$ is isogenous to a
 twist of its
 Galois conjugate over~$F$, then the result follows as in the proof of
Proposition~\ref{prop: a bunch of existing cases of potential
    modularity}. 
 We therefore assume that
 neither of these applies, and set $A=\Res_{K/F}E$. Then~$A$ is an abelian surface of
  type~$\mathbf{B}[C_2]$,
  and~$\cR_A=\Ind_{G_K}^{G_F}\cR_E$. By Theorem~\ref{thm: MB GSp4
    application}, there is a Galois extension of totally real fields
  $F'/F$, linearly disjoint from~$K/F$, and an automorphic
  representation~$\pi$ of $\GSp_4(\A_{F'})$ such that $\rho_{A,p}|_{G_{F'}}\cong \rho_{\pi,p}$.

Let~$\Pi$ be the transfer of~$\pi$ to~$\GL_4(\A_{F'})$.  If~$\kappa$ is the quadratic character of~$G_{F'}$ corresponding
  to~$K':=KF'/F'$, it follows
  that~$\Pi\otimes(\kappa\circ\Art_{F'}\circ\det)\cong\Pi$, so
  by~\cite[Thm.\ 4.2, 5.1 of \S 3]{MR1007299} there is a cuspidal
  automorphic representation~$\boldsymbol{\pi}$ of~$\GL_2(K')$ such that~$\Pi$ is
  the automorphic induction of~$\boldsymbol{\pi}\otimes|\det|$. Write $\Gal(K'/F')=\{1,\tau\}$. We claim that~$\boldsymbol{\pi}$ is of weight~$0$ and
  has trivial central character. Admitting this claim, it
  follows from Theorem~\ref{thm: Galois representations for GL2 over
    quadratic over totally real} that we can write
  \[\rho_{E,p}|_{G_{K'}}\oplus(\rho_{E,p}|_{G_{K'}}
    )^\tau\cong \rho_{\boldsymbol{\pi},p}\oplus \rho_{\boldsymbol{\pi},p}^\tau\]where all four
  $2$-dimensional representations are irreducible. After possibly replacing~$\boldsymbol{\pi}$ by~$\boldsymbol{\pi}^\tau$, we conclude that
  $\rho_{E,p}|_{G_{KF'}}\cong \rho_{\boldsymbol{\pi},p}$, so that by Theorem~\ref{thm: Galois representations for GL2 over
    quadratic over totally real}, for each place
  $v\nmid p$ of~$K'$ we have
  $\WD_v(\rho_{E,l}|_{G_{KF'_v}})^{\semis}\cong\rec(\boldsymbol{\pi}_v
  |\det|_v^{-1/2})^{\semis}$. It follows that in fact \[\WD_v(\rho_{E,l}|_{G_{K'_v}})^{F-\semis}\cong\rec(\boldsymbol{\pi}_v
  |\det|_v^{-1/2})\] (because we know the corresponding statement
  for~$A$). Repeating the argument for a second prime~$p$, we see that
  this holds for all finite places~$v$.

It remains to prove the claim.  By Lemma~\ref{lem: induction of real
    archimedean L parameters from GL2 to GSp4}, for each
  place~$v|\infty$ of~$K'$, either $\boldsymbol{\pi}_v$ corresponds to~$\phi_{0,1}$,
  or $v$ is complex, and the $L$-parameter of $\boldsymbol{\pi}_v$ is scalar, given
  by~$(z/\zbar)^{\pm 1}$; in particular, in either case it is
  algebraic, and so the central character~$\chi_{\boldsymbol{\pi}}$ of~$\boldsymbol{\pi}$ is
  algebraic. Moreover, if $\chi_{\boldsymbol{\pi}}$ is trivial, then
  the second case cannot occur, so that~$\boldsymbol{\pi}$
  automatically has weight~$0$. We therefore assume from now on
  that~$\chi_{\boldsymbol{\pi}}\ne 1$, and derive a contradiction.

Since~$\Pi^\vee\cong\Pi\otimes|\cdot|^{-2}$,
  we
  have \[\boldsymbol{\pi}\boxplus\boldsymbol{\pi}^\tau\cong\boldsymbol{\pi}^\vee\boxplus(\boldsymbol{\pi}^{\tau})^\vee, \]so
  that either~$\boldsymbol{\pi}^\vee\cong\boldsymbol{\pi}$,
  or~$\boldsymbol{\pi}^\vee\cong\boldsymbol{\pi}^\tau$. In the former
  case, we would have~$\boldsymbol{\pi}\cong\boldsymbol{\pi}^\vee\cong\boldsymbol{\pi}\boxtimes\chi_{\boldsymbol{\pi}}^{-1}$,
  from which it follows (if~$\chi_{\boldsymbol{\pi}} \ne 1$) 
  that~$\chi_{\boldsymbol{\pi}}$ is the character of a quadratic
extension~$L'/K'$ and~$\boldsymbol{\pi}$ is induced
from~$\GL(1)/L'$. This implies that~$\rho_{A,p}$ is potentially abelian,
and thus that~$E$ is CM, a contradiction.

We can therefore assume that
$\boldsymbol{\pi}^\vee\cong\boldsymbol{\pi}^\tau$, so that
$\chi_{\boldsymbol{\pi}^\tau}=\chi_{\boldsymbol{\pi}}^{-1}$,
and we shall derive a contradiction from these assumptions. Write~$\chi$
for the $p$-adic character $G_{K'}\to\Qpbartimes$ corresponding
to the algebraic character~$\chi_{\boldsymbol{\pi}}$. Let~$v\nmid p$ be a place of~$K'$ for
which~$\rho_{A,p}|_{G_{K'}}$ is unramified, and let the eigenvalues
of~$\rho_{E,p}(\Frob_v)$ be $\{\alpha_v,q_v/\alpha_v\}$ and those of
$(\rho_{E,p})^\tau(\Frob_v)$ be $\{\beta_v,q_v/\beta_v\}$. Now, 
   ~$\alpha_v q^{-1/2}_v$ is either a Satake parameter
 of~$\boldsymbol{\pi}$ or~$\boldsymbol{\pi}^{\tau}$, so (using
 that~$\chi_{\boldsymbol{\pi}^\tau}=\chi_{\boldsymbol{\pi}}^{-1}$) it
 follows that one of~$\beta_v, q_v/\alpha_v,q_v/\beta_v$ is equal to
 either $q_v\chi(\Frob_v)/\alpha_v$ or
 $q_v\chi^{-1}(\Frob_v)/\alpha_v$.

 Since ~$\chi$ is non-trivial, there is a set of places~$S$ of~$K'$ of
positive density such that~$\chi(\Frob_v) \ne 1$.
Shrinking~$S$ if necessary, we deduce that there exists a set~$S$ of
positive density so that one of the following  equalities holds for
all~$v \in S$:
 $$\begin{aligned}
 q_v \chi(\Frob_v)/\alpha_v = & \ \beta_v, \\
  \ q_v \chi(\Frob_v)/\alpha_v = & \ q_v/\beta_v, \\
 \ q_v \chi^{-1}(\Frob_v)/\alpha_v = & \ \beta_v, \\
  \ q_v \chi^{-1}(\Frob_v)/\alpha_v = & \ q_v/\beta_v.\end{aligned}$$
  By symmetry (replacing~$\chi$ by~$\chi^{-1}$ if necessary and~$\beta_v$ by~$q_v/\beta_v$ if necessary), we may assume that~$\alpha_v \beta_v = q_v \chi(\Frob_v)$
  for all~$v \in S$.

Now, the  representations~$\rho_{E,p}|_{G_{K'}}$
and~$(\rho_{E,p}|_{G_{K'}})^\tau$ have monodromy groups~$\GL(2)$
by~\cite[Thm.\ IV.2.2]{serreabladic} and are not twist equivalent (by
our running assumptions). 
It follows that the monodromy group of their tensor product is the identity component of~$\mathrm{GO}(4)$.
  Since this is connected, we deduce by considering the formal character  of the corresponding
  Lie algebra~$\mathfrak{sl}_2 \times \mathfrak{gl}_2$ that for any fixed character~$\xi$,
 the generic element of the tensor product~$\xi^{-1} \otimes \rho \otimes \rho^{\tau}$ 
  does not have~$1$ has an eigenvalue.

However, for each place~$v \in S$, we have (since~$\chi(\Frob_v) =
\alpha_v \beta_v/q_v$ for $v\in S$):
$$ \varepsilon^{-1}(v) \chi^{-1}(v) \otimes \{\alpha_v,q_v/\alpha_v\} \otimes \{\beta_v,q_v/\beta_v\}  
=  \{1,q^2_v/\alpha^2_v \beta^2_v,q_v/\beta^2_v,q_v \alpha^2_v\}.
  $$
  Since~$S$ has positive density, this is a contradiction, as
  required.
\end{proof}

\subsection{K3 surfaces of large rank}\label{subsec: K3}
If~$A$ is an abelian surface over a totally real field~$F$, then one may define the Kummer surface~$\Km(A)$ to be the resolution of the quotient of~$A$ under the map~$x \mapsto - x$.
The variety~$\Km(A)$ is a smooth projective algebraic K3 surface with (geometric) Picard number~$\ge 17$. (All Picard numbers in this section will
be geometric Picard numbers.)

\begin{prop} \label{prop:kummer} Let~$A$ be an abelian surface over a totally real field~$F$, and let~$X = \Km(A)$. Then Conjecture~\ref{conj:serre} holds for~$X$.
\end{prop}

\begin{proof} The cohomology groups~$H^*(X,\Q_p)$ are trivial in odd degree.
In even degree, they are generated by~$H^*(A,\Q_p)$  plus the~$16$ dimensional space of 
Tate cycles in~$H^2(X,\Q_p)$ spanned
by the~$16$ exceptional divisors in the resolution~$X \rightarrow A/(\pm 1)$. The latter classes are all defined over a finite extension of~$\Q$,
and hence the Galois representation (up to twist) they generate is an Artin representation.
 Hence the result follows from Theorem~\ref{thm: potential modularity implies meromorphic continuation} applied to~$A$, together with the meromorphic continuation 
 of Artin~$L$-functions.
 \end{proof}
 
 More generally, if a K3 surface~$X/F$ admits a Shioda--Inose structure~\cite[\S6]{Morrison} over~$F$,
 then~$H^*(X,\Q_p) \simeq H^*(\Km(A),\Q_p)$ for some abelian surface~$A/F$,
 and Prop~\ref{prop:kummer} implies
 Conjecture~\ref{conj:serre}  for~$X$. 
 It might also happen that~$X/F$ admits a Shioda--Inose
 structure
 over some
 finite extension~$E$.  Recall Ribet's notion of a~$\Q$-curve (\cite{RibetQ})
 as an elliptic curve over~$\Qbar$ all of whose conjugates by~$G_{\Q}$
 are  isogenous:
 
 \begin{df}  \label{Qcurve} An~$F$-abelian variety is an abelian variety~$A$ over a
 Galois extension~$E/F$ all
 of whose $\Gal(E/F)$-conjugates  are isogenous to~$A$ over~$E$.
 \end{df}
 
 Suppose that the conjugates~$A^{\sigma}$ over~$A$ are isogenous
 to (at most) quadratic twists of~$A$ (as necessarily happens if~$\{ \pm 1\}$
 are the only automorphisms in $\End_{\C}(A)$).
 Then the Galois
representations
 associated to~$\wedge^2 H^1(A)$ and thus to~$\Km(A)$ extend
 (even as compatible systems with~$\Q$-coefficients) to~$G_F$.
 Moreover, the (absolutely irreducible) projective Galois representations associated to~$H^1(A)$ also
extend to~$G_F$, and thus, from the vanishing of~$H^2(G_F,\Q/\Z)$ due to Tate, also
 give rise to~$G_F$ representations (now with coefficients).
 If~$F$ is totally real and~$A$ is an~$F$-abelian surface,
  one expects that the methods of this paper will have implications for the potential modularity of~$A$.
   (Note that for primes~$p$ splitting completely in~$E$, the mod~$p$
 representations over~$F$ locally arise from abelian surfaces over~$\Q_p$ --- namely~$A$ itself.)
 We have not endeavored to undertake the task of proving results along these lines, 
 however, since verifying that the Galois
 representations extend in the appropriate manner (especially  when all the different possibilities 
 for~$\End_{\C}(A)$ are taken into account) would necessitate a somewhat involved analysis which we avoid due to issues of time and space.

\subsubsection{General K3 surfaces of Picard rank \texorpdfstring{$\ge 17$}{ge 17}} \label{blather}
An algebraic K3 surface of Picard number~$17$ or~$18$
need not admit a Shioda--Inose structure even over~$\C$ (\cite{Morrison}).
There need not even be a correspondence
between~$X$ and  an abelian surface~$A$ inducing a Hodge isometry
of transcendental lattices~$(T_X \otimes \Q) \simeq (T_A \otimes \Q)$.
The problem, as noted in~\cite{Morrison}, is the following.
Let~$U$  denote the hyperbolic plane --- the lattice of rank two generated
by two isotropic vectors which pair to~$1$.
Then   there are obstructions on the lattices~$T_X$ of signature~$(2,20 - \rho(X)) = (2,3)$ or~$(2,2)$
which arise from~K3 surfaces
to admit an injection of the form~$(T_A \otimes \Q) \hookrightarrow (U \otimes \Q)^3$.
One might still hope to construct abelian varieties from
K3 surfaces of large Picard rank by
directly  lifting the weight two polarized Hodge structure on~$T_X$ to a weight one
Hodge structure of the smallest possible dimension. This amounts to considering the  $\GSpin$ cover of the corresponding orthogonal group
and relating that (in an ad hoc manner) to  weight one Hodge structures via the identification of the associated
 Shimura variety as one of Hodge type. This differs slightly from the
Kuga--Satake construction in which one has a functorial map from weight two
Hodge structures to weight one Hodge structures via the Clifford algebra construction ---
the latter  gives rise  to abelian varieties in a uniform way, but introduces (in general) auxiliary dimensions, and,
for a transcendental lattice~$T_X$ of rank~$5$, would produce an abelian variety of dimension~$2^3 = 8$.
In the case of interest to us,
the corresponding~$\GSpin$ Shimura variety will  now (over~$\C$) be  precisely 
 the moduli of abelian fourfolds with quaternionic multiplication (as considered in~\cite{KR}), where  the
degenerate case~$D = M_2(\Q)$  corresponds to the usual moduli space
of abelian surfaces.
In particular, we arrive at the conclusion that  a K3 surface with~$\rho(X) = 17$ or~$18$  should either admit a correspondence
with an abelian surface~$A$ inducing an isometry~$(T_X \otimes \Q) \simeq (T_A \otimes \Q)$,
or there will exist an abelian fourfold~$A$ with quaternionic multiplication
(a fake abelian surface, see the discussion after the statement of Lemma~\ref{lemma:motives})
and a correspondence inducing an injection~$(T_X \otimes \Q)^4 \hookrightarrow (T_A \otimes \Q)$.
This can also be predicted more arithmetically  by using the Yoga of motives. For convenience, suppose
that~$\rho(X/F) = 17$. Let~$\CR$ be  the compatible system  associated to the transcendental
motive (that is, the motive associated to the transcendental lattice), and assume that
the Galois representations~$r_p$ are strongly irreducible for a density one set of primes~$p$.
One can try to lift~$\CR$ (up to quadratic twist) to a~$4$-dimensional
compatible system~$\CS$ via the isogeny~$\GSp_4 \rightarrow \GO_5$, and then realize~$\CS$
as the motive associated to an abelian surface.  This  happens,  for example, when~$X = \Km(A)$
for some~$A$ over~$F$.
In general, however, one encounters 
two obstructions. The first is that one should expect to have to extend coefficients of the
motive  by a compositum of quadratic fields. This is because
a characteristic polynomial of an element in~$\GO_5$ with coefficients in~$\Q$ lifts to a characteristic
polynomial for~$\GSp_4$ whose coefficients lie either in~$\Q$ or~$\sqrt{D} \cdot \Q$ for some~$D$.
(This is an elementary computation with symmetric polynomials.)
Let~${}^{\sigma}\kern-0.1em{\CS}$ denote the compatible system obtained by applying
an automorphism~$\sigma
\in G_{\Q}$ to the coefficients of~$\CS$. 
Since~$\wedge^2 \CS = \CR = \wedge^2 ({}^{\sigma}\kern-0.1em{\CS})$, the irreducibility assumptions
imply that~$\CS$ is a quadratic twist of~${}^{\sigma}\kern-0.1em{\CS}$ for all~$\sigma \in G_{\Q}$ acting on the coefficients.
The restriction of~$\CS$ will thus have rational coefficients over some field~$E/F$ with~$\Gal(E/F) = (\Z/2\Z)^n$
where all the quadratic twists become trivial. If this restriction
 corresponds to an abelian surface~$A$, this would predict (and even imply,
see the remarks at end of~\cite{Morrison}) that there existed an algebraic cycle
on~$X \times \Km(A)$ which identified the  corresponding transcendental lattices over~$\Q$. 
Moreover, the abelian surface~$A$ would be an~$F$-abelian surface in the sense
of Definition~\ref{Qcurve}.
On the other hand, 
even supposing~$\CS$ has~$\Q$-coefficients over~$E$, it need not be the case that~$\CS$ comes from an abelian surface, even though (for weight reasons) it must be an abelian motive. One also has to allow the possibility that it comes from a fake abelian
surface, that is, a fourfold~$A$ with quaternionic  multiplication (see the proof
of Lemma~\ref{lemma:motives}). 
In summary, given a K3 surface~$X$ of  Picard rank at least~$17$ over a number field~$F$,
one should be able to associate to~$X$ a canonical isogeny class of $F$-abelian  surfaces
or~$F$-fake abelian surfaces.
Under sufficiently big image hypotheses, it should be possible to rigorously justify  the arguments
of this paragraph using the methods and language
of~\cite[\S4]{Patrikis}.

\subsubsection{Fake Kummer surfaces} \label{section:nori} 
This raises the natural question as to whether, given an abelian fourfold with an
 inclusion~$D \hookrightarrow \End^0(A)$ (a fake abelian surface), there are any natural
 geometrical constructions which produce a K3 surface (or, conversely, a construction in the other direction).
 For that matter, one might ask for an explicit geometric construction of either  of these objects. 
 Given six lines in general position in~${\mathbf{P}}^2$, the desingularization~$X$ of the double cover branched over those lines is, in general, a K3
  surface of Picard rank~$16$. If the~$6$ lines are all tangent to a smooth conic, however, then the K3 surface generically has Picard rank~$17$,
  and moreover~$X$ is the Kummer surface associated to the Jacobian of the hyperelliptic curve obtained as the double cover of the conic
  branched at the six tangent points~\cite{MorrisonKummer}. This suggests looking for other degenerations of the six
  lines which could give rise to transcendental lattices with different integral structures.

The following construction, suggested to the authors by Madhav Nori,
 gives a~$3 = 20 -17$ dimensional rational family of such degenerations corresponding to~$D = (-1,3)_{\Q}$. Given a generic point
 in this family of Picard number~$17$, the corresponding~$K3$ cannot be isogenous to a Kummer surface, and so indeed defines
 a genuine false Kummer surface. It is an interesting question to determine whether one can also see the corresponding
 abelian fourfold from this construction --- possibly associated to a generalized Prym variety of some natural cover of curves under
 the map~$\pi: X \rightarrow {\mathbf{P}}^2$.
 Consider five lines~$L_i$ for~$i = 1,\ldots 5$ in~$\mathbf{P}^2$.   These determine a conic~$C$ which
passes through the intersections~$L_1 \cap L_2$, $L_2 \cap L_3$, $L_3 \cap L_4$,  $L_4 \cap L_5$, and~$L_5 \cap L_1$,
which we denote by~$P_i$ for~$i = 1,\ldots 5$. Let~$L_6$ denote a sixth line which is tangent to~$C$ at~$P_6$.
 Note that~$C \cdot \sum_{i=1}^{6} L_i = 2 \left(\sum_{i=1}^{6} P_i\right)$ is divisible
by~$2$.  Let~$Y$ denote the degree~$2$ cover of~${\mathbf{P}}^2$ and~$X$
its desingularization. The lifts of~$P_i$ in~$Y$ for~$i = 1, \ldots, 5$ are ordinary double points, and so the exceptional
divisors~$E_i$ in~$X$ satisfy~$E_i.E_i = -2$. Let~$M = \sum_{i=1}^{5} E_i$.
If~$\pi: X \rightarrow \mathbf{P}^2$ denotes the projection, then~$\pi^{-1}(C) = M + D$, where~$D$ is now an everywhere
unramified double cover of~$C$. But~$C \simeq \mathbf{P}^1$, so~$D$ must decompose into two components~$A+B$
meeting transversally at~$\pi^{-1}(P_6)$.
Note that~$M.M = 5(-2) = -10$, that~$A.B = 1$ (meeting transversally at~$\pi^{-1}(P_6)$), and~$\pi^{-1}(C).\pi^{-1}(C) = 2 (C.C) = 8$.
Moreover, $A.M = B.M = 5$, intersecting in~$E_i$ for~$i = 1,\ldots 5$. 
It follows that, if we let~$E = A - B$, then~$E.E$ is equal to
$$(A + B + M).(A + B + M) - 2 (A+B).M - M.M - 4.A.B = 8 - 20 + 10 - 4 = -6.$$
The class $E$ is transverse to all exceptional classes as well as the pre-image of the hyperplane class, so gives a new class in~$\NS(X)$.
Note that~$U \otimes \Q \simeq (\langle 2k \rangle \oplus \langle -2k \rangle) \otimes \Q$ for any integer~$k$. 
The transcendental lattice of the generic~$X$ is~$(U^2 \oplus \langle -2 \rangle^2)_{\Q} \simeq (\langle 6 \rangle \oplus \langle -6 \rangle
\oplus U \oplus \langle -2 \rangle^2)_{\Q}$, hence the corresponding transcendental lattice of this restricted family is rationally
contained in~$(U \oplus \langle -6 \rangle
\oplus \langle -2 \rangle^2)_{\Q}$. Since this rational family has dimension~$20-17 = 3$, the generic member
will have Picard rank~$17$. As the form of the corresponding orthogonal group does not split, the lattice does not admit
an injection into~$(U \otimes \Q)^3$, and so~$X$ is not isogenous to any Kummer surface. Indeed, from the rational
structure of the resulting lattice, the corresponding fake abelian surface~$A$ will have endomorphisms by~$D = (-1,3)_{\Q}$.
(A related example was also considered in~\cite{Lombardo} --- in particular the divisor denoted in~\cite{Lombardo} by~$X_6$.)

\subsubsection{An Example} Take the conic to be~$y = x^2$, and the points~$P_i$ for~$i = 1,\ldots 5$
to be~$(n,n^2)$ for~$n = -2,\ldots, 2$. Now choose the point of tangency~$P_6$ to be at~$(3,9)$, so~$Y$
can be given by:
$$w^2 = (-x + y) (x + y) (y - 4 z) (-3 x + y + 2 z) (3 x + y + 2 z) (-6 x + y + 9 z).$$
The classes considered above are all defined over~$\Q$, and so~$\Pic(X/\Q) \ge 17$. Let~$\CR = (\Q,S, \{r_{p}\})$ denote the corresponding~$5$-dimensional
compatible system of~$\GO_5(\Q_p)$-representations. One checks that the set~$S$  of places of bad reduction is contained in the set of primes~$\{p \le 11,23,37,\infty\}$. 
Using both the determination of~$S$ and the fact that~$\CR(1)$ is self-dual, one computes that the determinant of~$\CR(1)$ is~$\psi$, the quadratic character
associated to~$K = \Q(\sqrt{-2 \cdot 3 \cdot 7 \cdot 23 \cdot 37})$. The compatible system~$\CR(1) \otimes \psi$ is  valued in~$\SO_5(\Q_p)$.
The  corresponding symplectic compatible system~$\CS$ 
of rank~$4$ need not have coefficients in~$\Q$, since it may come from a~$\Q$-abelian variety
in the sense of Definition~\ref{Qcurve}. Indeed, it has coefficients in~$\Q(\sqrt{3})$, 
and~${}^{\sigma}\kern-0.1em{\CS} \simeq \CS \otimes \psi$, where~$\sigma$ is the non-trivial element of~$\Gal(\Q(\sqrt{3})/\Q)$, and~$\psi$
is the character as above of conductor~$\Delta_K$. In particular,  although~$A$ is a~$\Q$-fourfold, the field of definition will be~$K$.

Over~$\C$, one expects that Nori's construction gives a rational parameterization of a component of the~$\GSpin$ Shimura variety 
associated to the quaternion algebra~$D = (-1,3)_{\Q}$ with some small (possibly trivial) level structure. Over~$\Q$, the~$\Q$-structure appears 
(by examining examples) to be associated
to a twisted form associated to~$\Q$-fourfolds~$A$ over a quadratic extension whose associated rank four motive over~$\Q$ has coefficients in~$\Q(\sqrt{3})$.

\section{Applications to
  modularity} \label{sec:modularityapplications} In this section,  we apply our
main modularity lifting theorem (Theorem~\ref{thm: final R equals T}) to prove modularity theorems for
abelian surfaces. The methods generalize those of~\cite{MR1333035,MR1415322} for
elliptic curves. In \S\ref{subsec: paramodular} and~\S\ref{subsec: paramodular is
false}, we show that our results confirm the paramodular conjecture
of~\cite{MR3165645} in many cases, but that there are counterexamples
to the original formulation of the conjecture (arising from ``fake abelian surfaces'').
\subsection{First modularity results}We begin this section with a proof of Theorem~\ref{theorem:infinite} of the introduction. 

\begin{prop} \label{prop:infinitetwo} Let~$F$ be a totally real
  field in which~$p>2$ splits completely. Let~$A/F$ be an abelian
  surface with good ordinary reduction at all places~$v|p$,
and suppose that at each~$v|p$, the unit root crystalline eigenvalues are distinct modulo~$p$. 
Assume that~$A$ admits a polarization of degree prime to~$p$. Let
$$\rhobar_{A,p}: G_F \rightarrow \GSp_4(\F_p)$$
denote the dual of the mod $p$ Galois representation associated to~$A[p]$, and assume that~$\rhobar_{A,p}$ is vast and tidy.
Assume that~$\rhobar_{A,p}$ is  ordinarily modular, in the sense that
there exists $\pi$ of parallel weight~$2$ and central
character~$|\cdot|^2$ which is unramified and
    ordinary at all~$v|p$, such
    that~$\rhobar_{\pi,p}\cong\rhobar_{A,p}$, and~$\rho_{\pi,p}|_{G_{F_v}}$ is pure for all finite
    places~$v$ of~$F$.
 Then~$A$ is modular.   More
  precisely, there is an ordinary automorphic
  representation~$\pi'$ of $\GSp_4(\A_F)$ of parallel weight~$2$ and
  central character~$|\cdot|^2$ which satisfies
  $\rho_{\pi',p}\cong\rho_{A,p}$. 
 \end{prop}
 
 \begin{proof} As before, we write ~$\rho_{A,p}: G_{F} \rightarrow \GSp_4(\Q_p)$ for the Galois
  representation associated to the dual of the~$p$-adic Tate module of~$A$. 
 The assumption that~$A$ admits a polarization of degree prime to~$p$ implies
 that the image of~$\rho_{A,p}$ lands in~$\GSp_4(\Z_p)$ and~$\rhobar_{A,p}$ lands in~$\GSp_4(\F_p)$.
  By Proposition~\ref{prop: pure abelian variety compatible system}, the
 representation~$\rho_{A,p}$ is pure for all places~$v$ of~$F$. 
 The assumption that~$A$ has good ordinary reduction
 for all~$v|p$ and distinct unit root crystalline eigenvalues for all~$v|p$ implies that the representations~$\rho_{A,p}$
 restricted to~$G_{F_v}$  are~$p$-distinguished weight~$2$ ordinary.
  Prop~\ref{prop:infinitetwo} is then an immediate consequence of
  Theorem~\ref{thm: final R equals T}. \end{proof}

\begin{rem} If~$A$ does \emph{not} have a polarization of order prime to~$p$, then,
by considering the kernel~$A[\lambda]$ of any polarization~$\lambda: A \rightarrow A^t$,
we deduce that the representation~$\rhobar_{A,p}: G_{\Q} \rightarrow \Aut(A[p]) = \GL_4(\F_p)$
is reducible. Hence one could replace the assumption of the existence of a polarization on~$A$
of order prime to~$p$
in Prop.~\ref{prop:infinitetwo} by the assumption that the Galois representation associated 
to~$A[p]$ is irreducible. On the other hand, we do not phrase our theorem
in this way for the following reason: if~$A$ does not have a polarization of order prime to~$p$,
then it need not even be the case that the (necessarily reducible) 
representation~$\rhobar_{A,p}: G_{\Q} \rightarrow \GL_4(\F_p)$
associated to~$A[p]$ lands in any conjugate of~$\GSp_4(\F_p)$. Indeed, let~$E/\Q$ be any elliptic curve such that~$\rbar_{E,3}: G_{\Q} \rightarrow \GL_2(\F_3)$ has surjective image, let~$K/\Q$ be an auxiliary degree~$3$ cyclic extension,
let~$B = \Res_{K/\Q}(E)$, and let~$A$ denote the kernel of the map~$B \rightarrow E$
induced from the trace map~$\Z[\Gal(K/\Q)] \rightarrow \Z$. Then~$A$ is an abelian surface,
and~$\rhobar_{A,3} \simeq \rbar_{E,3} \otimes W$, 
where~$W \in \Ext^1_{G_{\Q}}(\F_3,\F_3)$ is the 
unique non-trivial extension  which splits over~$\Gal(K/\Q)$. The group
theoretic image of~$\rbar_{A,3}$ is
 isomorphic to~$\GL_2(\F_3) \times \Z/3\Z$,
but this is not isomorphic to any subgroup of~$\GSp_4(\F_3)$. 
These examples are also
related to the failure of the Shafarevich--Tate group~$\Sha$ to have square order --- William Stein~\cite{Stein} found  abelian surfaces~$A$ exactly of the form considered above with~$3 \| \Sha(A)[3^{\infty}]$.
\end{rem}

  We now give some examples where one can directly establish the modularity
  of certain residual representations.

\begin{prop}  \label{prop:applications} Let~$F$ be a totally real field in which~$p > 2$ splits completely. Let~$\rhobar_p: G_{F} \rightarrow \GSp_4(\F_p)$ be an absolutely irreducible
representation with similitude factor~$\vareps^{-1}$ which is vast and tidy and
 $p$-distinguished weight~$2$
    ordinary. Suppose furthermore that either:
\begin{enumerate}
\item $p = 3$, and~$\rhobar_{3}$ is induced from a~$2$-dimensional
  representation with inverse cyclotomic determinant over a totally real quadratic extension~$E/F$ in which~$3$ is unramified.
\item $p = 5$, and~$\rhobar_{5}$ is induced from a~$2$-dimensional
  representation valued in~$\GL_2(\F_5)$ with inverse cyclotomic determinant  over a totally real quadratic
  extension~$E/F$ in which~$5$ is unramified.
\item $\rhobar_{p}$ is induced from a character of a quartic CM
  extension~$H/F$ in which~$p$ splits completely. 
  \end{enumerate}
Then~$\rhobar_{p}$ is  ordinarily modular, that is,
 there exists $\pi$ of parallel weight~$2$ and central character~$|\cdot|^2$ which is unramified and
    ordinary at all~$v|p$, such
    that~$\rhobar_{\pi,p}\cong\rhobar$, and~$\rho_{\pi,p}|_{G_{F_v}}$ is pure for all finite
    places~$v$ of~$F$.
 \end{prop}
 
 \begin{proof}
Suppose that we are in one of the first two settings, so that~$p = 3$ or~$5$, and~$\rhobar = \Ind^{G_F}_{G_E} \varrhobar$ for some representation~$\varrhobar: G_{E} \rightarrow \GL_2(\Fbar_p)$ with determinant~$\vareps^{-1}$.
The assumptions on~$\rhobar$ imply that~$\varrhobar|_{G_{E(\zeta_p)}}$ is  absolutely
irreducible, and the restriction of~$\varrhobar$ to the inertia group
at any prime~$w|p$ is an extension of~$\vareps^{-1}$ by~$1$.
If~$p = 5$, the condition on the determinant and the fact that~$E$ is unramified at~$p$ additionally ensures that the projective image of~$\varrhobar$ is not~$A_5$.
The representation~$\varrhobar$ locally has the structure of a representation associated to an ordinary Hilbert modular form of parallel weight two and trivial nebentypus.  Suppose that~$\varrhobar$
 is modular. It follows from~\cite[Thm.\ A]{blggordII} that~$\varrhobar$
 does indeed arise from a Hilbert modular form of this kind, and we may take~$\pi$ to be the automorphic induction of this form from~$E$ to~$F$. Since~$E/F$ is unramified, this will preserve the property
 of being ordinary. As in the proof of Theorem~\ref{thm: MB GSp4
   application}, purity follows from the main results of~\cite{MR2327298,MR3272276}. Hence it suffices to establish the modularity
 of~$\varrhobar$. 

If~$\varrhobar$ has solvable image, 
 then, from a classification of the finite subgroups of~$\GL_2(k)$ for a finite field~$k$ 
 (see for example~\cite{SerreAntwerp}),  we deduce that the projective
  image of~$\varrhobar$ is either~$A_4$, $S_4$, or dihedral, and is in particular a subgroup of 
  $\PGL_2(\C)$.  By a theorem of Tate (see~\cite[Theorem 4]{Serrelifting}),
   this implies that there exists a characteristic zero
  lift of~$\varrhobar$ which is totally odd with finite solvable image,
   and 
the result follows from an application of the theorems of Langlands and Tunnell
 (\cite{MR574808,MR621884})
 as in~\S5 of~\cite{MR1333035}. It remains to consider the
 representations with vast non-solvable image. For~$p = 3$, the only
 non-solvable induced  representations which are vast come from 
 representations (Lemma~\ref{lem:sample}~\eqref{exoticinduction})
 $\varrhobar_3: G_{F} \rightarrow \GL_2(\F_9)$ with projective
 image~$A_5$. The modularity of such a representation follows as in the
 solvable case, except now invoking the odd Artin conjecture
 for totally real fields (\cite[Thm.\ 0.3]{MR3581178}) rather than Langlands--Tunnell. Alternatively, the arguments
of~\cite{MR2180399} over~$\Q$ may be adapted to this setting.

 Thus we are left with the case of  non-solvable representations~$\varrhobar: G_{E} \rightarrow \GL_2(\F_5)$ with determinant~$\vareps^{-1}$, which necessarily are surjective.
 The method of Khare--Wintenberger implies the existence of
 characteristic zero lifts of the required form (for example by~\cite[Thm.\ 7.2.1]{0905.4266} --- the assumption that~$5$ is unramified in~$E$ guarantees
 that~$[E(\zeta_5):E] = 4$). To show that such a lift is modular, it
 suffices (by, for example, the main theorem of~\cite{kis04}) to show
 that~$\varrhobar$ is modular.  However, this follows from a standard argument going back
 to~\cite{MR1415322,MR1981033} by realizing~$\varrhobar$ as the~$5$-torsion of a
 modular elliptic curve over a solvable extension. In our situation,
 we may explicitly invoke~\cite[Prop.\ 2.1.3]{MR3581178}.

Suppose finally that~$\rhobar_p = \Ind^{G_F}_{G_H} \chi$, where~$H/F$
is a quartic CM extension in which~$p$ splits completely.
 Let~$v|p$ be a prime in~$F$. The assumption that~$\rhobar_p$ is ordinary implies that for two of the
 primes~$w|v$ of~$H$ the
 restriction of~$\chi$ to inertia at~$w$ is~$\vareps^{-1}$,  and it is
 trivial at the other two primes above~$v$.
 Let~$\psi$ denote an algebraic Grossencharacter of~$G_{H}$ with conductor prime to~$p$ and CM type corresponding to the mod-$p$ weights of~$\chi$. If~$\psi_p$
 is the~$p$-adic avatar of~$\psi$, then, by construction, the character $\psi_p/\chi \mod p$ is unramified at~$p$, and hence,
 after twisting~$\psi$ by the Teichmuller lift of this character, we may assume that~$\psi_p \equiv \chi \mod p$. Let~$E/F$ denote the intermediate real quadratic field inside~$H$.
 Then the automorphic induction of~$\psi$ to~$\GL_2(\A_E)$ is a Hilbert modular form of parallel weight two which is ordinary at all~$v|p$ and has trivial central character.
Inducing once more to~$F$, we obtain the required form~$\pi$.
 \end{proof}

 For explicit examples of abelian surfaces~$A/\Q$ with~$\End_{\C}(A) = \Z$ whose mod-$3$ or mod-$5$ representations~$\rhobar_{A,p}$
 satisfy Prop.~\ref{prop:applications} --- and hence, by Prop.~\ref{prop:infinitetwo}, are modular --- see~\cite{CGS}.
 In contrast to the examples found in~\cite{BPVY} and~\cite{Berger} of large prime conductor, the examples found in~\cite{CGS}
 have good reduction  outside~$2$, $3$, $5$, and~$7$.
 
 We also have the following application to modularity over number fields which need not be totally real (or even CM).
 
 \begin{theorem} \label{theorem:funky} Let~$F$ be a totally real field in which~$5$ splits completely, and let~$K/F$ be a quadratic extension in which~$5$ is unramified.
 Let~$E/K$ be an elliptic curve which has good ordinary reduction or semistable ordinary reduction for all places~$w|5$
 of~$K$.
 Finally, assume that   the representation~$\varrhobar_{E,5}: G_{K} \rightarrow \GL_2(\F_5)$  has the following properties:
 \begin{enumerate}
 \item The projective image of~$\varrhobar_{E,5}$ is either~$S_5 = \PGL_2(\F_5)$ or~$S_4$.
 \item There exists a representation~$\rbar_5: G_{F} \rightarrow \GL_2(\F_5)$ with determinant~$\vareps^{-1}$ such that ~$\rbar_5|_{G_{K}}\cong\varrhobar_{E,5}$.
 \end{enumerate}
 Then~$E$ is modular.
 In particular, there exist infinitely many modular elliptic curves
 over~$K$ up to twist which are not CM and  do not come from any subfield
 of~$K$.
 \end{theorem}

For example, one could take~$F$ to be~$\Q(\sqrt{d})$ for any~$d  \equiv 1,4 \mod 5$, and then take~$K = \Q(\sqrt[4]{d})$, which is a field of mixed signature.

 \begin{proof} It suffices to prove that the twist of~$E$ by some quadratic character is modular.
 We now apply Lemma~\ref{lem:double} to the representation~$\rbar_5$ of~$G_{F}$ to obtain a representation
 $$\rhobar 
 = \Ind^{G_F}_{G_K}(\rbar_5 |_{G_K} \otimes \delta_{M/K}) = \Ind^{G_F}_{G_K}(\varrhobar_{E,5} \otimes \delta_{M/K}),$$
 where~$\delta_{M/K}$ is  the character of an auxiliary quadratic extension.
 By Lemma~\ref{lem:formixed}, this representation is vast and tidy
 and~$p$-distinguished weight~$2$
    ordinary. 

As in the proof of Proposition~\ref{prop:applications},
 it follows from our hypotheses that~$\rbar_5$ comes from an ordinary Hilbert modular form~$\boldsymbol{\pi}$ for~$F$. By taking the base change
 of this form to~$K/F$, twisting by the quadratic character~$\delta_{M/K}$, and then inducing back to~$F$, we construct a~$\pi$ of parallel weight~$2$ and central character~$|\cdot|^2$ which is unramified and
    ordinary at all~$v|p$ such
    that~$\rhobar_{\pi,p}\cong\rhobar$. Again, the purity
    of~$\rho_{\pi,p}$ follows from the main results of~\cite{MR2327298,MR3272276}.
    It follows from    Theorem~\ref{thm: final
      R equals T} that
    $$\rho = \Ind^{G_F}_{G_K}(\varrho_{E,5} \otimes \delta_{M/K})$$
    is modular,
   and hence (exactly as in the proof of Theorem~\ref{thm: potential
     modularity of elliptic curves}) that~$\varrho_{E,5} \otimes \delta_{M/K}$
   and hence~$E$ is modular. 

It is easy to produce examples of ~$E$ satisfying the hypotheses of
the theorem (starting with an elliptic curve over~$\Q$, for example).
   Using the fact that the genus zero curve~$X(\varrhobar_{E,5})$ is
   isomorphic to~$\mathbf{P}^1$ over~$K$ (there being at least one rational point coming from~$E$), we deduce that there will
   be infinitely many such points. On the other hand, by choosing such points with appropriate local properties (for example, ramified at one prime~$w$ above~$v$
   but not at the other) we may find infinitely many examples which do not arise via base change. Since the mod~$5$-representations associated to these curves are not
   projectively dihedral or cyclic, they also cannot have CM.
    \end{proof}

\subsection{Abelian varieties with fixed~\texorpdfstring{$3$}{3}-torsion}
We have produced a number of residual representations mod~$p$ for small~$p$ which are 
automorphic.
It is natural to ask whether any such representation  (satisfying necessary local conditions) arises  from infinitely
many abelian surfaces over~$F$. The corresponding question for~$2$-dimensional representations 
has a positive answer precisely when~$p = 2$, $3$, or~$5$, where the corresponding moduli space is a smooth curve
of genus zero. We show that for abelian surfaces there is a positive answer for~$p = 2$ and~$p = 3$.
When~$p \ge 5$, the moduli space in question is of general type~\cite{Hulek}, and so one would not expect (in general) that they admit
infinitely many rational points not lying on a  special Shimura
subvariety, although we do not  attempt to address this question.

When~$p = 2$, the problem is pretty much obvious.  The fact that
the corresponding moduli space for the trivial representation~$\rhobar$ is rational
goes back to Igusa (see~\cite[Theorem IV.1.4]{Hulek} and~\cite{Igusa}).
The fact that the corresponding moduli space for non-trivial~$\rhobar$
is unirational is also surely well-known (we shall now give a sketch of this result
although we shall never use this fact).
Fix a  representation~$\rhobar: G_{F} \rightarrow \GSp_4(\F_2)$.
Since~$\GSp_4(\F_2) \simeq S_6$,  one may write
 any~$G$-extension~$L$ of~$F$  for~$G \subseteq S_6$ as the splitting field of a degree~$6$ separable  polynomial~$f(x)$ over~$F$. If one then
 takes~$A$ to be the Jacobian of the curve~$y^2 = g(x)$ for any~$g(x)$
 with~$\Q[x]/g(x) \simeq \Q[x]/f(x)$, then~$\rhobar$
 is the representation associated to the~$2$-torsion of~$A$.
 An elementary computation shows that  this gives  a~$3 = 6 - \dim \PGL_2$ dimensional family of abelian
 surfaces up to isomorphism with fixed~$\rhobar$ for any such~$\rhobar$. Explicitly, one may let~$e_i$ for~$i = 1$ to~$6$ be any basis over~$\Q$ of
 the \'{e}tale $\Q$-algebra~$\Q[x]/f(x)$, and then let~$g(x)$ be the minimal polynomial of~$\sum t_i e_i$.
(The Jacobian~$A$  depends only on~$g(x)$ up to the action of~$\PGL_2$ on~$\mathbf{P}^1$.)

 This leaves the case~$p = 3$. The answer in this case can be extracted from the very extensive literature on the subject,
essentially following the main idea of~\cite{MR1415322}. Fix~$\rhobar: G_F \rightarrow \GSp_4(\F_3)$ with inverse cyclotomic
similitude character, and let~$V$ denote the underling symplectic space over~$\F_3$.
Let~$B(\rhobar)/F$ denote the moduli space of
  pairs~$(A,\imath_3)$ consisting of abelian surfaces~$A$
and symplectic isomorphisms
\[\imath_3:A[3]\isoto 
V^{\vee}.\] 
(The dual is here because our Galois representations have been normalized cohomologically, so it is the dual representations
which actually occur inside the~$p$-adic Tate modules.)
The variety~$B(\rhobar)$ is
smooth and geometrically
connected. Over~$\C$,   we may identify~$B$ with the moduli space
of principally polarized abelian surfaces with full level~$3$ structure.
This space is well-known to be a (geometrically) rational threefold,
and is isomorphic to an open subvariety of the Burkhardt quartic~\cite{Burkhardt,Coble,Baker,Hunt,1705.09006},
specifically, the complement in the Burkhardt quartic of the Hessian hypersurface.

The Burkhardt quartic is exceptional for a number of different reasons, not least
of which is that it admits an action of the group~$\PSp_4(\F_3)$
(tautologically from the description above). 
If~$V = (\mu_3)^2 \oplus (\Z/3\Z)^2$, we write~$B$ for~$B(\rhobar)$. One knows (\cite{1705.09006}) that~$B$ is rational over~$\Q$.
Suppose we knew that~$\rhobar$ actually came from an abelian
surface~$A$, so that~$B(\rhobar)$ admitted a smooth rational point
over~$F$.
One might ask whether this is enough to force the twist~$B(\rhobar)$
to be rational over~$F$; 
this question is resolved in the negative in~\cite{CC}. 
The difficulty in a na\"{\i}ve attempt to replicate the argument of Taylor and Shepherd-Barron (\cite{MR1415322}) in this case is that
the birational map~$B \rightarrow \mathbf{P}^3$ is \emph{not} equivariant with respect to~$\PSp_4(\F_3)$ and any
embedding~$\PSp_4(\F_3) \rightarrow \PGL_4(\Qbar) = \Aut(\mathbf{P}^3)$.  This means that a~$\PSp_4(\F_3)$-twist of~$B$ does not naturally
inherit the structure of a Severi--Brauer variety.

 It turns out, however, that we are lucky. 
There exists a cover~$P(\rhobar) \rightarrow B(\rhobar)$ of degree~$6$ corresponding to an additional choice of level~$2$ structure of~$A$,
namely an odd theta characteristic,
(or, for~$A = \Jac(C)$, a Weierstrass point on the corresponding genus two curve~$C$).
The cover~$P(\rhobar)$ now \emph{does} have the property that it is not only rational, but~$\PSp_4(\F_3)$-equivariantly rational,
which allows us to deduce the rationality of~$P(\rhobar)$ in favourable circumstances, and hence the unirationality of~$B(\rhobar)$.
In particular, this allows us to construct infinitely many rational
points on~$B(\rhobar)$ which correspond to abelian surfaces~$A$ with $\End_{\C}(A)=\Z$, as in the following theorem.

\begin{theorem} \label{boxershorts} Fix~$\rhobar: G_{F} \rightarrow \GSp_4(\F_3)$  with similitude character~$\vareps^{-1}$. 
Then~$B(\rhobar)$ is unirational over~$F$,
and there exist infinitely many principally polarized abelian surfaces~$A/F$ up to twist 
with~$\End_{\C}(A) = \Z$
and such that~$\imath_3: A[3] \simeq V(\rhobar)^{\vee}$.
Moreover, we may additionally assume that these~$A$ are Jacobians of curves  which
have a rational Weierstrass point, and may thus be written in the form~$y^2 = f(x)$ where~$f(x)$ is  a quintic polynomial. Suppose, in addition,  that for all~$v|3$, the representation~$\rhobar|_{G_{F_v}}$ arises as the~$3$-torsion of an abelian surface  over~$F_v$ with good ordinary
reduction. Then we may additionally assume that these~$A$ also all have good ordinary reduction for all~$v|3$.
\end{theorem}

Note that, as with of~\cite[Thm.\ 1.2]{MR1415322}, we do not need to impose any further local hypotheses  at any primes after we impose the global condition
on the similitude character. (In particular, if~$K$ is a local field
of characteristic zero, then by using a globalization argument as
in~\cite[Thm.\ 3.1]{frankII}, this implies that the only requirement on a mod~$3$-representation~$\rhobar: G_{K} \rightarrow \GSp_4(\F_3)$
to arise from a principally polarized abelian
surface over~$K$ is that the similitude character
is~$\varepsilon^{-1}$.)

The remainder of this section is devoted to the proof of 
Theorem~\ref{boxershorts}. We start by defining the non-Galois degree~$6$ cover~$P$ of~$B$ and recalling its basic properties.

\begin{df} Let~$B(2):=\mathcal{A}_2(6)\rightarrow \mathcal{A}_2(3) = B$ denote the cover of~$B$ corresponding to a choice of full level-$2$ structure. It is a Galois cover with
Galois group~$S_6 \simeq \PSp_4(\F_2)$, where we fix this identification up to conjugacy by identifying~$S_6$ generically
with the Galois group of the Weierstrass points on~$C$ with~$\Jac(C) = A$.
(This identification  can be made explicit using the map~$\tau$
below.) 
Then~$P$ denotes the intermediate cover over~$B$ corresponding to the conjugacy class of subgroups~$S_5 \subset S_6 = \PSp_4(\F_2)$ which fix a point.
\end{df}

A more natural definition of~$P$ is given in terms of theta characteristics. Namely, $P$ may be identified with the moduli space~$\mathcal{A}_2(3)^{-}$
of principally polarized abelian surfaces with a symmetric odd  theta structure of level~$3$ (see~\cite[\S2.3]{MR2457735}).
Recall that there is a natural Torelli map~$\tau: \mathcal{M}_2(3) \rightarrow \mathcal{A}_2(3)$ which is a bijection away from the Humbert surface consisting
of principally polarized abelian surfaces which split as the product of two elliptic curves. A rational point on~$P$ in the image of~$\tau$
corresponds to~$A = \Jac(C)$, together with a symplectic isomorphism from~$A[3]$ 
 to~$(\mu_3)^2 \oplus (\Z/3\Z)^2$ and the data of a rational
 Weierstrass point on~$C$, which (after moving this point to infinity) means that~$C$ can
be written in the form~$y^2 = f(x)$ for a quintic polynomial~$f(x)$.

The level~$3$-structure on~$B$ pulled back to~$P$ gives an action over~$\Q(\sqrt{-3})$ of the group~$\PSp_4(\F_3)$.
If one gives~$\PSp_4(\F_3)$ the structure of an \'{e}tale group scheme
over~$\Q$ by viewing it as the group~$G$ of symplectic automorphisms of~$\mu^2_3 \times (\Z/3 \Z)^2$ modulo~$(-1,-1)$, then this action descends to~$\Q$. Equivalently, the twisting of any
automorphism in~$\PSp_4(\F_3)$ by the action of~$\sigma \in G_{\Q}$ is accounted for
by the Galois action on~$\mu^2_3 \times (\Z/3 \Z)^2$ (as in
the formula for~$\sigma \theta$ below).

\begin{prop}\label{prop: P is rational} The variety~$P$ is rational over~$\Q$. Moreover, there exists a birational map~$P \rightarrow \mathbf{P}^3$ over~$\Q$ which is equivariant
with respect to the action of~$G$ for some action of~$G$ on~$\mathbf{P}^3$.
\end{prop}

\begin{proof} The~$G$-equivariant map~$P \rightarrow \mathbf{P}^3_{\Q}$ is the odd theta map denoted~$\mathrm{Th}^{-}$ in~\S2.4 of~\cite{MR2457735}.
The fact that~$\mathrm{Th}^{-}$ is a  birational isomorphism is~\cite[Theorem~0.0.1]{Weddle}. \end{proof}

We now turn to the proof of Theorem~\ref{boxershorts}. From Proposition~\ref{prop: P is rational},
it follows that the rationality of~$P(\rhobar)$ over~$F$ is equivalent to the rationality of~$\mathbf{P}^3(\rhobar)$ over~$F$, where~$\mathbf{P}^3(\rhobar)$ is the twist of~$\mathbf{P}^3$
arising from the projective representation associated to~$\rhobar$.
The action of~$\Sp_4(\F_3)$ on~$\mu^2_3 \times (\Z/3\Z)^2$ over~$\Q(\zeta_3)$ induces a homomorphism~$\theta$
from~$\Sp_4(\F_3)$ to~$\Aut(P)$  and hence to~$\Aut(\mathbf{P}^3)$.
This map satisfies
$$\sigma \theta(\alpha) = \theta \left( \left( \begin{matrix} \vareps (\sigma) & 0 & 0 & 0 \\
0 & \vareps(\sigma) & 0 & 0 \\ 0 & 0 & 1 & 0 \\ 0 & 0 & 0 & 1 \end{matrix} \right)  \alpha
\left( \begin{matrix} \vareps^{-1}(\sigma) & 0 & 0 & 0 \\
0 & \vareps^{-1}(\sigma) & 0 & 0 \\ 0 & 0 & 1 & 0 \\ 0 & 0 & 0 & 1 \end{matrix} \right)\right).$$
Since~$\rhobar$ has similitude factor~$\varepsilonbar^{-1}$, we can
associate to~$\rhobar$ a cocycle
$$\sigma \mapsto \theta \left( \rhobar^{\vee}(\sigma)
\left( \begin{matrix} \vareps^{-1}(\sigma) & 0 & 0 & 0 \\
0 & \vareps^{-1}(\sigma) & 0 & 0 \\ 0 & 0 & 1 & 0 \\ 0 & 0 & 0 & 1 \end{matrix} \right)\right),$$
in~$H^1(F,\PGL_4(\Qbar))$, 
and~$\mathbf{P}^3(\rhobar)$ is the twist of~$\mathbf{P}^3$
by this cocycle.

\begin{lemma} If~$\rhobar: G_{F} \rightarrow \GSp_4(\F_3)$ has
  similitude factor~$\vareps^{-1}$, then~$P(\rhobar)$ is rational.
\end{lemma}

\begin{proof} The proof is very similar to  the proof of~\cite[Lem.\ 1.1]{MR1415322}. We need to show that the cocycle
  in~$H^1(F,\PGL_4(\Qbar))$ corresponding to~$\rhobar$ vanishes, so it
is enough to show that it comes from~$H^1(F,\GL_4(\Qbar))$. It is
therefore
enough (following the argument in~\cite{MR1415322}) to show that we can lift the induced
homomorphism~$\theta:\PSp_4(\F_3)\to\PGL_4(\Qbar)$ to a \emph{unique}
homomorphism~$\wtheta:\Sp_4(\F_3)\to\GL_4(\Qbar)$. 
Since~$\PSp_4(\F_3)$ is perfect (indeed simple) it has a unique Schur
 cover (Darstellungsgruppe).
Since the Schur multiplier of~$\PSp_4(\F_3)$ has order~$2$ (\cite{Atlas}),
the Darstellungsgruppe of~$\PSp_4(\F_3)$ may be identified with~$\Sp_4(\F_3)$, and 
in particular the projective
representation~$\theta$ lifts to a genuine homomorphism~$\wtheta: \Sp_4(\F_3)
\rightarrow \GL_4(\Qbar)$.  It remains to show this lift is unique.

We claim  this follows from the fact that~$\Sp_4(\F_3)$ is perfect.
Indeed, because the group is perfect, every element of~$\Sp_4(\F_3)$ can be written as a product of commutators~$[g,h]$.
Hence it suffices to show that~$\wtheta([g,h])$ is uniquely defined. But~$\wtheta([g,h]) =
[\wtheta(g),\wtheta(h)]$, and  
the commutator of any two elements of~$\GL_n(\Qbar)$  depends only on
their images in~$\PGL_n(\Qbar)$, as required.
We note (although we do not use this fact) that~$\PSp_4(\F_3)$ has no faithful~$4$-dimensional representations
(again by~\cite{Atlas})  and so the representation~$\wtheta$ of~$\Sp_4(\F_3)$ is faithful.
\end{proof}

\begin{rem} \label{rem: soothing Boxer's conscience}
Note that the corresponding facts (that~$\PSL_2(\F_5)$ is simple with 
Darstellungsgruppe~$\SL_2(\F_5)$) lead to a proof 
of~\cite[Lem.\ 1.1]{MR1415322}.
This differs slightly from the original proof in~\cite{MR1415322} as follows: 
Instead of  using the fact that~$\PSL_2(\F_5)$ has Schur multiplier~$\Z/2\Z$ and deducing that
any irreducible projective representation~$\theta$ lifts to a representation~$\wtheta: \SL_2(\F_5)
\rightarrow \GL_2(\Qbar)$, the authors use the fact that the kernel of~$\SL_2(\Qbar) \rightarrow
\PSL_2(\Qbar)$ has order~$2$, and so~$\theta$ \emph{automatically} lifts to a 
representation of a degree~$2$ central extension of~$\PSL_2(\F_5)$ to~$\SL_2(\Qbar)$, and then
argue that the image (and hence source) is~$\SL_2(\F_5)$ (because the split central extension would give
a faithful~$2$-dimensional representation of~$\PSL_2(\F_5)$).
  \end{rem}

We have now proved under the given hypothesis on~$\rhobar$ that~$P(\rhobar)$ is rational,
and hence~$B(\rhobar)$ is unirational. Moreover, twisting~$P$
by~$\rhobar$ leaves the level structure at the odd theta
characteristic unchanged,
so that all the corresponding abelian varieties in the image of the Torelli map are Jacobians of curves~$C$ of the form~$y^2 = f(x)$ where~$f(x)$ is a quintic
(after moving the rational Weierstrass point to infinity).

\begin{proof}[Proof of Theorem~\ref{boxershorts}.]
  We need to
  show that infinitely many of the corresponding points
  of~$P(\rhobar)$ do not admit any extra endomorphisms over~$\C$. We show that we may find
  infinitely many~$A$
  such that the Galois representation associated to~$A[5]$ has
  image containing~$\Sp_4(\F_5)$. If~$A/F$ \emph{did} admit extra
  endomorphisms over~$\overline{F}$, then, from the
  classification of the possible Galois types of endomorphism
  structures on~$A$  recalled at the beginning of~\S\ref{subsec:abelian surfaces}, the Galois
  representation associated to the~$5$-adic Tate module of~$A$ would become
  reducible after making an extension of degree at most~$2$. But the
  action of~$\Sp_4(\F_5)$ on~$\F^4_5$ remains absolutely irreducible
  after restriction to any index two subgroup, which
  forces~$\mathrm{End}_{\Qbar}(A) = \Z$. (In this argument, $5$ could
  have been replaced by any prime~$p$ independent of the level
  structure.)  

  We will now arrange this condition by an application of
  Hilbert irreducibility as in the proof
  of~\cite[Theorem~1.2]{MR1415322}.
  Let~$R(\rhobar) \rightarrow P(\rhobar)$ be the fibre product
  of~$P(\rhobar)$ with~$\mathcal{A}_2(5) \rightarrow
  \mathcal{A}_2$. This is a  Galois cover with Galois group~$\PSp_4(\F_5)$. Recall that~$P(\rhobar)$
  is rational. By Hilbert irreducibility (\cite[\S 9.2, 9.6]{SerreLectures}), we may find infinitely many points~$x \in P(\rhobar)(F)$
  so that the Galois group of the splitting field of any preimage~$y \in R(\rhobar)$
  of~$x$ contains~$\PSp_4(\F_5)$, and moreover, we may restrict~$x$ to
  any non-trivial open subset of~$P(\rhobar)(F_v)$ for all~$v$ in some finite set of primes~$S$.
  If~$A$ denotes the corresponding abelian surface, it follows that the 
  projective Galois representation associated to~$\rhobar_{A,5}$ contains~$\PSp_4(\F_5)$, and thus the image of~$\rhobar_{A,5}$ itself contains~$\Sp_4(\F_5)$.
If  we now use the
  assumption
   that~$P(\rhobar)(F_v)$ has points corresponding
  to abelian surfaces with good ordinary reduction at all~$v|3$, then (since the
  ordinary condition is open) we can choose our~$x \in P(\rhobar)(F)$ so that~$A/F$  has good
  ordinary reduction for all~$v|3$.  
\end{proof}
We obtain the following corollary:
\begin{theorem} \label{thm:bcgpthm}Let~$F$ be a totally real field in which~$3$ splits completely. Then there exist infinitely many abelian surfaces~$A/F$ up to twist with~$\End_{\C}(A) = \Z$
which are modular, and which do not come from any proper subfield
of~$F$ \emph{(}in the sense that, for each~$p$, there is no twist of the corresponding
Galois representation~$\rho_{A,p}$ which extends to the absolute
Galois group of a proper subfield of~$F$\emph{)}.
\end{theorem}

\begin{proof} Let~$H/F$ be a quadratic extension in which every
  prime~$v|3$ is inert, and let~$\sigma$ be the non-trivial element of~$\Gal(H/F)$. Let~$E/H$ be an elliptic curve with good ordinary reduction  for all~$v|3$,
 such that~$\rhobar_{E,3}: G_{H} \rightarrow \GL_2(\F_3)$  is surjective,  and such that the projective images of~$\rhobar_{E,3}$ and~$\rhobar^{\sigma}_{E,3}$
 are totally disjoint. Since~$X(1)$ has genus zero, this can be achieved by choosing a global point which lies over a suitable choice of smooth point
in~$X(1)(k_w)$ (with~$k_w = \OL_H/w$) for suitably chosen primes~$w$
which split over~$F$, for example ensuring that~$\rhobar_{E,3}(\Frob_w)$ and~$\rhobar_{E,3}(\Frob_{w^{\sigma}})$
 give distinct elements of~$\PGL_2(\F_3)$. Similarly, by making choices above primes of~$\Q$ which split completely in~$H$, 
 we may ensure that~$\rhobar = \Ind^{G_F}_{G_H} \rhobar_{E,3}$ does not
 descend after twisting to any proper subfield of~$F$.

By Lemma~\ref{lem: Ind GL2 enormous}, $\rhobar$ is vast and tidy. 
The fact that each prime~$v|3$ in~$H/F$ is inert implies that~$\rhobar$ is~$3$-distinguished and finite flat,
 and hence~$3$-distinguished weight~$2$ ordinary. It follows from
 Proposition~\ref{prop:applications} that~$\rhobar$ is ordinarily modular. The representation~$\rhobar$ at each~$v|3$
arises locally from an abelian variety over~$F_v$ with good ordinary reduction, because it does so globally --- namely,   the restriction of scalars of~$E$ from~$H_v$ to~$F_v$. 
 It follows from Theorem~\ref{boxershorts} that there are infinitely
 many abelian surfaces~$A$ with  up to twist with good ordinary
 reduction at each place~$v|3$ and satisfying~$A[3]\cong\rhobar$.  The choice of~$\rhobar$ ensures that any such~$A$
 does not descend  (even after twist) to any subfield of~$F$. Finally,
 every such~$A$ is modular by Proposition~\ref{prop:infinitetwo}.
\end{proof}

\subsection{The Paramodular Conjecture}\label{subsec: paramodular}We
end this section with a discussion of the relationship between our
results and the ``paramodular conjecture'' of~\cite{MR3165645}
(cf. also the remarks in~\cite[\S8, p.243]{YoshidaOriginal}). Recall
that this conjecture states that there should be a bijection
(determined by the compatibility of Frobenius eigenvalues and Hecke
eigenvalues at unramified places) between isogeny classes of abelian
surfaces~$A/\Q$ of conductor~$N$ with~$\End_\Q A=\Z$, and holomorphic cuspidal
Siegel newforms of weight~$2$ and paramodular level~$N$ which are
``non-lifts'' and have rational Hecke eigenvalues, considered up to scalar
multiplication. Here ``non-lifts'' means that they are orthogonal to
the space of Gritsenko lifts. We explain in this  section why the
paramodular conjecture as originally formulated in~\cite{MR3165645} is
not true. The issue is that Siegel newforms of
weight~$2$  and paramodular level~$N$ with rational eigenvalues will  \emph{not}
always correspond to abelian surfaces. In light of the observations of this paper,
Brumer and Kramer have modified
their conjecture in~\cite{BrumerKramerCorrection} along the lines suggested by the analysis presented here --- we
reproduce their updated conjecture in this paper as
Conjecture~\ref{conj:modified} below. 
  In order to distinguish between the two
versions of this conjecture, we refer to the original formulation (given above) as the
\emph{original} paramodular conjecture, and the modified version (Conjecture~\ref{conj:modified}) as the paramodular conjecture.
Both of these conjectures posit an injective map from isogeny classes of abelian surfaces~$A/\Q$ of conductor~$N$
with~$\End_\Q A=\Z$ to
Siegel newforms of weight~$2$ and paramodular level~$N$, which are
``non-lifts'' and have rational Hecke eigenvalues, and hence, when talking about the implication in this direction,
we do not distinguish between different versions of the conjecture.

We firstly show that all of our examples of modular abelian varieties
verify the paramodular conjecture, before giving a more general
explanation of the relationship between the paramodular conjecture and
the Langlands program, and then explaining some counterexamples to the
original paramodular conjecture.

\begin{lem}
  \label{lem: our things are paramodular}Any abelian variety~$A/\Q$
  satisfying the hypotheses of Proposition~\ref{prop:infinitetwo}
  satisfies the paramodular conjecture; that is, there is a
  corresponding holomorphic cuspidal Siegel newform of weight~$(2,2)$ and
  paramodular level equal to the conductor of~$A$, which is a
  non-lift, has rational Hecke eigenvalues, and is unique up to scalars.
\end{lem}
\begin{proof}
  By Proposition~\ref{prop:infinitetwo} (or more precisely by
  Theorem~\ref{thm: final R equals T}, as applied in the proof of Proposition~\ref{prop:infinitetwo}) there is an $L$-packet of cuspidal automorphic
  representations~$\pi$ of weight~$2$ and general type corresponding to~$A$, whose
  $L$-parameters coincide with those determined by~$A$. The claim that there
  is a unique corresponding newform of level equal to the conductor of~$A$ is
  now a consequence of the theory of newforms due to Roberts and
  Schmidt~\cite{MR2344630} (which assumes that we are working
  with representations of trivial central character, but this is
  harmless, as we can reduce to this case by twisting~$\pi$
  by~$|\cdot|$). This newform is certainly a non-lift, as~$\pi$ is of
  general type (see the discussion following this lemma for a more
  precise description of the non-lifts), and it has rational Hecke
  eigenvalues by local-global compatibility.

  More precisely, by~\cite[Thm.\ 1.1]{SchmidtParamodularPacket}, for each
  prime~$v$ of~$\Q$, there is a unique paramodular representation in
  the $L$-packet at~$v$, namely the unique generic
  representation. Since representations of general type are stable,
  this gives rise to a unique~$\pi$ of weight~$2$ which has a
  paramodular vector at each finite place. Furthermore, for each~$v$ the space of
  paramodular vectors at minimal paramodular level is one-dimensional
  by~\cite[Thm.\ 7.5.1]{MR2344630}, and this minimal paramodular level
  coincides with the conductor of the corresponding $L$-parameter (and
  thus with that of~$A$) by~\cite[Thm.\ 7.5.4(iii)]{MR2344630} and the
  main theorem of~\cite{gantakeda}.
\end{proof}

We now discuss the paramodular conjecture more broadly. Firstly, we discuss
the automorphic side of the conjecture. As explained
in~\cite{SchmidtParamodularPacket}, the space of Siegel modular forms
of weight~$2$ and fixed level can be written as an orthogonal sum of
spaces spanned by eigenforms in automorphic representations of the
various types in Arthur's classification. The Gritsenko lifts are precisely
those of Saito--Kurokawa type, while those of one-dimensional type do
not contribute to the cuspidal spectrum. Since abelian surfaces
with~$\End_\Q A=\Z$  should correspond to automorphic representations
of general type (as their corresponding Galois representations are
irreducible), we see that it is implicit in the statement of the
conjecture that there are no paramodular eigenforms (at least with rational
Hecke eigenvalues) of Yoshida, Soudry, or Howe--Piatetski-Shapiro
type.

This is indeed the case, as is proved
in~\cite{SchmidtParamodularPacket,SchmidtCAP}. The case of Yoshida
type is~\cite[Lem.\ 2.5]{SchmidtParamodularPacket}; in this case, the
parameters are unstable, and the corresponding packet of
representations does not satisfy the required sign condition. Indeed,
at each finite place, the condition that the representation admits a
paramodular vector forces the sign to be trivial, whereas the
condition of being the holomorphic limit of discrete series at
infinity gives a non-trivial sign. Note that the analogous argument
would fail for totally real fields of even degree.

The cases of Soudry and Howe--Piatetski-Shapiro type are~\cite[Prop.\
5.1]{SchmidtCAP}. In these cases the obstructions to the
existence of paramodular vectors are at finite places; it turns out
that at the places where these representations are ramified, there are
no paramodular vectors. In these cases the representations are
parameterized by certain Hecke characters, and the fact that the
representations are ramified at some finite place comes from the fact
that any Hecke character must be ramified. Accordingly, the analogous
argument could fail for totally real fields of class number greater
than~$1$.

It follows from this discussion that the original paramodular conjecture is
equivalent to the claim that there is a bijection between isogeny
classes of abelian surfaces~$A/\Q$ with~$\End_{\Q}A=\Z$, and
cuspidal automorphic representations~$\Pi$ of~$\GL_4(\A_\Q)$ of symplectic
type with multiplier~$|\cdot|^2$, whose infinity type is the one corresponding to the
$L$-parameter~$\phi_{2;1,0}$, and whose Hecke eigenvalues are all
rational. 
In one direction, given~$A$, the
existence of~$\Pi$ is certainly predicted by the
Fontaine--Mazur--Langlands conjecture, the rationality of its Hecke
eigenvalues following from strong multiplicity one. We now explore the converse
direction.

\begin{lem} \label{lemma:motives} Let~$F$ be a totally real field.
Assume the Fontaine--Mazur Conjecture, the Standard Conjectures,  the Hodge Conjecture,
and that the Galois representations associated to any cuspidal automorphic
representation~$\Pi$  for~$\GL_4(\A_{F})$ whose infinity type for
each~$v|\infty$ corresponds to the~$L$-parameter~$\phi_{2;1,0}$ form
an irreducible weakly compatible system. Let~$\Pi$ be such a
representation with the properties that its Hecke eigenvalues are rational, and that~$\Pi$
is of symplectic type with multiplier~$|\cdot|^2$.
Then, associated to~$\Pi$,
there exists a corresponding motive~$A/F$ such that either:
\begin{enumerate}
\item $A/F$ is an abelian surface.
\item \label{caseofdeath} $A/F$ is an abelian~fourfold with endomorphisms over~$F$ by an order in a quaternion
algebra~$D/\Q$.
\end{enumerate}
Moreover, if~$A/F$ is an abelian~fourfold with~$\End_{F}(A) \otimes \Q
= \End_{\C}(A) \otimes \Q$ an indefinite quaternion algebra~$D/\Q$,
and one  assumes only standard automorphy conjectures,
then there exists a corresponding~$\Pi$ of symplectic type with
rational eigenvalues and
multiplier~$|\cdot|^2$.
\end{lem}

One might reasonably (following Serre~\cite[\S0.7, p.\ DeRa-13]{DR}) call an abelian fourfold~$A$ with endomorphisms by an order in a quaternion
algebra~$D/\Q$ a \emph{fake} (or \emph{false}) \emph{abelian surface}
(\emph{fausse surface ab\'{e}lienne}).

\begin{proof}[Sketch of proof] 
(For a more detailed proof of a closely related result, see~\cite[Thm.\ 3.1]{MR3544295}.)
One first obtains from~$\Pi$ a rank~$4$ symplectic 
 weakly compatible irreducible family
of~$p$-adic Galois representations $$\CR=  (\Q,S,\{ Q_v(X) \}, \{r_p \}, \{H_\tau\} )$$
with~$H_{\tau} = (0,0,1,1)$ for all~$\tau | \infty$,
and such that
$$r_p: G_F \rightarrow \GSp_4(\Qbar_p) $$
has inverse cyclotomic similitude character. The Fontaine--Mazur conjecture implies that~$\CR$ arises
 from
a pure irreducible motive~$M$ over $F$ with coefficients in~$\Q$ (we also now assume the standard conjectures~\cite{MR1265519}). Concretely, this means that  $M$ is irreducible and that  for each prime $p$, the $p$-adic \'{e}tale realization of~$M$, $H^i_{et}(M, \Q_p) \otimes_{\Q_p} {\Qbar_p}$ contains $r_p$. By the Brauer--Nesbitt
theorem, all the twists of $r_p$ by automorphisms of the coefficient field $\Qbar_p$ are isomorphic to $r_p$. Therefore,  if we assume the Tate conjecture, we deduce that $H^i_{et}(M, \Q_p) \otimes_{\Q_p} {\Qbar_p}$ is a sum of copies of $r_p$. 
 The rank of~$M$ is therefore~$4d$ for some~$d$. Let~$\End_{\Q}(M) \otimes \Q = D$. Since $M$ is simple,~$D$ is a division algebra.
The centre of~$D$ is a number field~$E$. 
We claim that~$E = \Q$. It suffices to show that, for all~$p$, the centre of~$D \otimes_{\Q} \Qbar_p$ is~$\Qbar_p$.
By the Tate conjecture, however, we can determine~$D \otimes \Qbar_p$ from the endomorphisms
of the $p$-adic \'{e}tale realization of~$M$, which is isomorphic to a direct sum of $d$ copies of $r_p$. 
It follows that~$\End_{\Qbar_p[G_F]}(r_p^d)$ is a matrix algebra over $\Qbar_p$, and thus has centre~$\Qbar_p$. Hence~$E = \Q$.

Let us fix an embedding $F \rightarrow \C$. The Hodge realization of~$M$ is a polarized
Hodge structure of weight one, which gives a polarized torus, and thus (by Riemann) an abelian variety $B$ over $\C$. Since $M$ is defined over $F$, we deduce that, for any automorphism $\sigma$ of $\C$ over $F$, $B^\sigma$ is isogenous to $B$.  Let $2d$ be the dimension of $B$ and $r$ be the degree of its polarization.  Let $\mathcal{A}_{2d, r}$ be the coarse moduli space of  abelian varieties of dimension $2d$ with a degree $r$ polarization. This is a scheme of finite type over $\Q$. We base change it to $F$. Let $[B]$ be the  point on $\mathcal{A}_{2d, r}$ associated to $B$. The set of points on $\mathcal{A}_{2d, r}$ isogenous to $[B]$ is countable. On the other hand, if the residue field of $[B]$ in $\mathcal{A}_{2d, r}$ is transcendental over $F$, we deduce  that its orbit under $\mathrm{Aut}(\C/F)$ is uncountable because there are uncountably many ways to embed a transcendental field of finite type over $F$ into $\C$.  It follows that $[B]$ is defined over a finite algebraic extension $L$ of $F$ and that  $B$ is defined over a finite algebraic extension $L'$ of $L$. Moreover, $M\vert_{L'} = \mathfrak{h}^1(B)$. We now consider $C = \mathrm{Res}_{L'/F} B$, an abelian scheme over $F$. Looking at the $p$-adic \'etale realization $C_p$ of $C$ for a prime $p$, we find that $C_p = M_p \otimes \mathrm{Ind}_{L'/F} \mathbb{1}$ where $M_p$ stands for the $p$-adic realization of $M$. We now let $A$ be the simple factor of $C$ whose $p$-adic \'etale realization contains $r_p$. Then $\cR$ arises from $\mathfrak{h}^1(A)$. It  follows from the Tate conjecture that $\mathfrak{h}^1(A)$ and the motive $M$ we started with are in fact isomorphic. 

Now taking into account that the centre of~$D$ is~$\Q$, we deduce from the Albert classification
(see~\cite[Thm.\ 2, p.\ 201]{mumford}) that~$A$ is one of the following three types:
\begin{enumerate} 
\item Type I: $A/F$ is an abelian surface with~$\End_{\Q}(A) = \Z$.
\item Type II: $A/F$ is an abelian fourfold with~$\End_{\Q}(A) \otimes \Q = D$,
an indefinite quaternion algebra over~$\Q$.
\item\label{item: type III} Type III: $A/F$ is an abelian fourfold with~$\End_{\Q}(A) \otimes \Q = D$,
a definite quaternion algebra over~$\Q$.
\end{enumerate}
(Note that Type IV of the Albert classification cannot occur, because the centre~$F = \Q$ of~$D$ is not a totally imaginary CM field.)

Suppose that~$A/\Q$ is an abelian fourfold with~$\End_{\Q}(A) \otimes \Q = \End_{\C}(A) \otimes \Q = D$ for some indefinite quaternion algebra~$D/\Q$ (and thus of Type II above).
We now construct a suitable compatible system~$\CR$, which (by standard automorphy conjectures)
will give rise to a suitable~$\Pi$.
 The Mumford--Tate conjecture is known for the varieties of type II~\cite{Chi1,Chi2}, and the semisimple
part~$\mathfrak{h}$ of the Lie algebra of the Mumford--Tate group
of~$A$ is (for almost all~$p$)~$\mathfrak{sp}_4$. Let~$p$ be any prime which splits~$D$. Then~$H^1(A,\Q_p)$ has an action of~$D \otimes \Q_p = M_2(\Q_p)$. In particular, it decomposes as~$V_p \oplus V_p$
for an irreducible  4-dimensional representation~$V_p$ whose
monodromy group is contained in~$\GSp_4(\Q_p)$.

If~$Q_v(T)$ denotes the degree~$8$ polynomial in~$\Z[T]$ coming from the characteristic polynomial
of Frobenius at~$v$, then every root of~$Q_v(T)$ has even
multiplicity, and thus~$Q_v(T) = P_v(T)^2$ for a degree~$4$ polynomial~$P_v(T) \in \Z[T]$, which will
be the characteristic polynomial of Frobenius at~$v$ on~$V_p$. By the Weil conjectures, the roots of~$Q_v(T)$ obey the usual symmetry associated
to a weight one motive, and  so the same is true for~$P_v(T)$.
This implies that~$V^{\vee}_p \simeq V_p \otimes \varepsilon^{-1}$. Since the Galois representation has big image in~$V_p$, any isomorphism~$V_p \simeq V_p \otimes \chi$ forces~$\chi$ to be trivial,
and thus from the identification~$V^{\vee}_p \simeq V_p \otimes \varepsilon^{-1}$ above we deduce that the similitude character is
inverse cyclotomic. In particular, 
by standard automorphy conjectures, $V$ will be associated with a~$\Pi$ as in the theorem. 
By strong multiplicity one~\cite{MR623137}, the rationality of Hecke eigenvalues
at almost all primes (in particular primes of good reduction) forces rationality at all primes.
\end{proof}

\subsection{Examples and counterexamples}\label{subsec: paramodular is
false}

In this section, we give some examples of abelian fourfolds~$A/F$ with~$\End_{F}(A) \otimes \Q = D$
for a quaternion algebra~$D/\Q$.
In Lemma~\ref{winteriscoming}, we prove the existence of such~$A$ which also satisfy~$\End_{\C}(A) \otimes \Q = D$
for some indefinite~$D/\Q$.
But first, we construct abelian fourfolds~$A/\Q$ with~$\End_{F}(A) \otimes \Q = D$,
and such that (under standard conjectures, and even unconditionally in some cases), they correspond to a~$\Pi$
as above which comes from a paramodular eigenform with rational Hecke
eigenvalues, and thus contradict the original  paramodular conjecture.

\subsubsection{Abelian fourfolds of type III} \label{section:type3}
We expect that
case~\eqref{item: type III} considered in the proof of
Lemma~\ref{lemma:motives} cannot occur. While we do not show that
here, we instead discuss a minor subtlety which occurs when trying to
construct examples of this kind. 

 Let~$E/\Q$ be an elliptic curve (say without complex multiplication). Let~$F/\Q$ be (say) a totally real field with~$\Gal(F/\Q) = Q$, the quaternion group of order~$8$.
 The group~$Q$ has an irreducible representation~$W/\Q$ of dimension~$4$, which contains a stable integral lattice~$\Lambda \subset W$. 
 Note that, for any prime~$p$, there is a decomposition~$W \otimes \Qbar_p = V \oplus V$ for an irreducible~$2$-dimensional faithful representation~$V$ of~$Q$.
 Now let us define:
 $$A = E^4 = E \otimes_{\Z} \Lambda.$$
 We find that~$A$ is simple over~$\Q$, and~$\End_{\Q}(A) \otimes \Q = \End_{\Q}(W) = D$, where~$D/\Q$ is the Hamilton quaternions.
 The corresponding compatible system~$\CR$ arises from Galois representations
 $$\rho_p:= r_{E,p} \otimes V: G_{\Q} \rightarrow \GL_4(\Qbar_p).$$
 For~$p \ne 2$, the image of this Galois representation lies in~$\GL_4(\Q_p)$.
 On the other hand, the 
 possible symplectic forms associated to~$\rho_p$ are the one dimensional summands of
 $$\wedge^2 \rho = (\Sym^2 \rho_{E,p} \otimes \det(V)) \oplus (\det(\rho_{E,p}) \otimes \Sym^2(V)).$$
 Since~$\det\rho_{E,p} = \varepsilon^{-1}$,
we obtain symplectic representations with inverse cyclotomic similitude character if and only if~$\Sym^2(V)$ contains the trivial representation. But~$\det(V) = 1$ for the faithful complex~$2$-dimensional representation of~$Q$,
so~$\Sym^2(V)$ is the direct sum of the three non-trivial quadratic
characters of~$Q$, and these compatible families do not have the required form. 

More generally, suppose that~$A/\Q$ is an abelian fourfold with~$\End_{\Q}(A) \otimes \Q = D$ for some definite quaternion algebra~$D$.
The corresponding Shimura curves~$X_D$ parametrizing such objects lie in the exceptional
class of Shimura varieties with the property that there is a strict containment~$D \subsetneq \End_{\C}(A)$
for all complex points~$A$ of~$X_D$
(see~\cite[\S9.9]{LB}).
Since~$D/\Q$ is definite, the 
semisimple
part~$\mathfrak{h}$ of the Lie algebra of the
 Mumford--Tate group should (for almost all~$p$) be contained
 in~$\mathfrak{so}_4 = \mathfrak{sl}_2 \times \mathfrak{sl}_2$ rather
 than~$\mathfrak{sp}_4$ (by ~\cite[\S6.1]{Moonen}). 
 This forces the representations~$V_p$
 to decompose (as~$\Qbar_p$-representations) as the tensor product of a representation coming from a modular form with an Artin representation.
 We expect it should be possible to make a careful case by case
analysis to rule out this case occurring, but we have not attempted to
do this.

 \subsubsection{Abelian fourfolds of type II} \label{section:type2}
 One can produce examples of abelian fourfolds with endomorphisms by an  order in  an indefinite quaternion algebra~$D/\Q$ by taking
 the tensor product of a~$2$-dimensional representation with an Artin representation.
 Let~$B/\Q$ be an abelian surface of~$\GL(2)$-type
 with endomorphisms by an order in   a quaternion algebra~$D$  which are defined over a quadratic extension~$K/\Q$, and then take~$L \subset V$ to be a lattice in
  a~$2$-dimensional dihedral representation~$V$ over~$\Q$ which is induced from a quadratic character~$\chi$ of~$K$ which does not extend to~$\Q$ (so the action of~$G_{\Q}$
  on~$V$ is through a dihedral group of order~$8$). Then one can take~$A = B \otimes_{\Z} L$, which may be identified with the restriction of scalars of the quadratic twist~$B \otimes \chi$
 of~$B$ from~$K$ to~$\Q$. 

The action of an order of~$D$ on~$B$ and~$B \otimes \chi$ over~$K$ extends to an action of this order of~$D$ on~$A$. We obtain a compatible system~$\CR$
 of Galois representations
  $$\rho_p:= r_{B,p} \otimes V: G_{\Q} \rightarrow \GL_4(\Qbar_p).$$
 Because~$V$ is induced from~$G_K$, it follows that the characteristic polynomials of Frobenius of this representation all have coefficients in~$\Q$.
It now suffices to show that~$\rho_p$ preserves a symplectic form with
inverse cyclotomic similitude character.  
  The argument proceeds
 exactly as in~\S\ref{section:type3}, except now we have the isomorphism~$\wedge^2 V \simeq \chi_{K/\Q}$. In particular, the trivial character is a summand of~$\Sym^2 V$,
and thus~$\rho_p$ preserves a symplectic form with similitude character~$\varepsilon^{-1}$.
The representations~$\rho_p$  do not arise from abelian surfaces over~$\Q$, since that would
contradict the Tate conjecture for abelian varieties~\cite{Faltings}. Moreover, they are easily seen to be modular. 
Hence these give counterexamples to the original  paramodular conjecture. 

For an explicit example, one could take~$B$
to be the modular abelian surface which is a quotient of~$J_0(243)$ with coefficient field~$\Q(\sqrt{6})$ (see~\cite[Table 3]{CremonaTwist}), which is geometrically
simple and obtains quaternionic multiplication over~$\Q(\sqrt{-3})$. Then take any non-Galois invariant quadratic character~$\chi$ of~$\Q(\sqrt{-3})$, and let~$A = \Res_{K/\Q}(B \otimes \chi)$.

In light of Lemma~\ref{lemma:motives}, 
Brumer and Kramer have formulated the following natural modification of the
original paramodular conjecture (see~\cite{BrumerKramerCorrection}):

\begin{conj}[Paramodular Conjecture of Brumer--Kramer] \label{conj:modified}
Let~$\mathcal{A}_N$ denote the set of
isogeny
classes of abelian surfaces~$A/\Q$ with~$\End_{\Q}A=\Z$ and conductor~$N$,
and~$\mathcal{B}_N$ the set of isogeny classes of
fake abelian surfaces \emph{(}QM abelian fourfolds\emph{)} $B/\Q$ of conductor~$N^2$ with~$\End_{\Q} B$
an order in a non-split quaternion algebra~$D/\Q$.
Let~$\mathcal{P}_N$ denote the set of holomorphic weight 2 paramodular
forms~$f$ of level~$N$ up to nonzero scaling which have rational Hecke
eigenvalues and lie in the orthogonal complement to the space of
Gritsenko lifts.
Then there is a bijection between the set~$\mathcal{A}_N
\cup \mathcal{B}_N$ and~$\mathcal{P}_N$ such that
$$L(C,s) = L(f,s,\mathrm{spin}) \ \text{if $C \in \mathcal{A}_N$
and} \ L(C,s) = L(f,s,\mathrm{spin})^2 \ \text{if $C \in \mathcal{B}_N$.}$$
\end{conj}

We conclude with some remarks on the possible existence of abelian~fourfolds which satisfy
case~(\ref{caseofdeath}) of Lemma~\ref{lemma:motives} 
and additionally have no further endomorphisms over~$\C$ (such varieties
will necessarily be  geometrically simple).

\begin{xlem} There exists a totally real field~$F$, an indefinite quaternion algebra~$D$, and
an abelian~fourfold~$A/F$ with~$\End_{F}(A) \otimes \Q = \End_{\C}(A) \otimes \Q = D$.
\label{winteriscoming}
\end{xlem}

\begin{proof}[Sketch]  Let~$D =  \displaystyle{\left(\frac{-1,3}{\Q} \right)}$ be the unique quaternion algebra over~$\Q$ ramified at (exactly)~$2$ and~$3$.
Let~$\OL_D$ denote the maximal  order in~$D$.
In the standard way, one may also write down an involution~$\dagger$ obtained by conjugating the standard
involution so that~$\tr_{D/\Q}(x x^{\dagger})$ is positive definite, and  write down a non-degenerate alternating form~$\psi$
on~$(\OL^2_D) \otimes \Q$ which satisfies various compatibilities with~$\dagger$. Associated to~$\OL_D$ in the usual way is a Shimura stack (of level one)~$X$ parametrizing tuples~$(A,\lambda,\iota)$ where~$A$ is an~$S$-abelian~fourfold over~$S$,
$\lambda$ is a principal polarization over~$S$, and $\iota: \OL_D \rightarrow \End(A)$ is an injective homomorphism such that the Rosati
involution induced by~$\lambda$ restricts to~$\dagger$ and such that~$\psi$ is compatible with the polarization on homology
as an~$\OL_D$-module. $X$ is a smooth Deligne--Mumford stack over~$\Q$ with a single geometric 
component (cf.~\cite{KR}). The complex points~$X(\C)$
are uniformized by the Siegel upper half space of dimension~$3$, and 
 the generic point of~$X$ over~$\C$
has endomorphisms precisely by~$\OL_D$ (see~\S9.9 of~\cite{LB}).

By~\cite{Milne79}, there exists an abelian surface~$B/\C$ with~$\End_{\C}(B) \otimes \Q = D$, that
the Rosati involution on~$\End B$ is~$x \mapsto x^{\dagger}$, and such that the
restriction of scalars of~$B$ from~$\C$ to~$\R$
gives a point in~$X(\R)$.
(Another way to view this is to consider~$X$ as a~$\mathrm{GSpin}$ Shimura variety associated 
to a~$5$ dimensional quadratic space as in~\cite{KR}, and then  signature~$(1,2)$ subspaces will give Shimura curve
subvarieties.)

 This is not quite sufficient, however, to guarantee a point over~$X(F)$
for a totally real field~$F$ with the correct endomorphisms, nor even a point over~$X(\Qbar)$, since one has to remove from~$X(\C)$
a countable union of proper Shimura subvarieties, which might \emph{a priori} exhaust the~$\Qbar$-points of~$X$.  Moreover, due
to the stackiness of~$X$, there are issues comparing fields of definition versus fields of moduli. We therefore employ a trick already used
in the proof of Theorem~\ref{boxershorts}. Namely, impose level structure by choosing a large prime~$p > 3$ and
fixing a surjective representation~$\rhobar: G_{\Q} \rightarrow
\GSp_4(\F_p)$ with inverse cyclotomic similitude character.
 (Such
representations are abundant --- one source are the duals of the~$p$-torsion of
abelian surfaces over~$\Q$.)
Then~$X$ admits
a geometrically connected cover~$X(\rhobar)$ defined over~$\Q$ with level structure corresponding to~$A[p]^{\vee} = \rhobar \oplus \rhobar$,
with a suitable choice of polarization and compatible action of~$(\OL_D \otimes \Z_p)/p = (M_2(\Z_p))/p = M_2(\F_p)$.

The variety~$X(\rhobar)$ is a fine moduli space which is now a  smooth variety over~$\Q$ with real points, since the
point~$\Res_{\C/\R}(B)$ considered above has the appropriate level structure over~$\R$.
  Employing the theorem of Moret-Bailly~\cite{mb}, we
may deduce the existence of a totally real field~$F$ and a corresponding abelian variety~$A/F$ 
such that~$F$ is disjoint from the splitting field of~$\rhobar$. Because~$X(\rhobar)$ is
a fine moduli space, the variety~$A$ has endomorphisms by~$\OL_D$ over~$F$.
It now suffices to show that it has no further endomorphisms over~$\C$.
The dual of the Tate module of~$A$ decomposes as a Galois
representation as~$\rho \oplus \rho$ where~$\rho: G_F \rightarrow \GSp_4(\Z_p)$
is a lift of~$\rhobar$. Since~$p \ge 5$, the assumption that~$\rhobar$
has surjective image implies that~$\rho$ also has surjective image. 
However, if~$A$ admitted extra endomorphisms over any extension of~$F$, then
the image of~$\rho$ restricted to some open subgroup would lie inside a proper
algebraic subgroup of~$\GSp_4(\Q_p)$, contradicting the fact that image contains an open
subgroup of~$\GSp_4(\Z_p)$.
 \end{proof}

It would be interesting to know whether  (for suitable choices) these
 varieties have points over~$\Q$ which correspond to~$A/\Q$ with~$\End_{\C}(A) \otimes \Q = D$,
but this is not so easy to determine by pure thought. However, the specific~$X$ chosen above (ramified at only~$2$ and~$3$) is possibly
the most likely choice to be rational, since it corresponds to the indefinite quaternion algebra~$D/\Q$ of smallest discriminant. The construction of Nori in~\S\ref{section:nori} suggests that, for this~$D$,
the moduli space is at least geometrically rational. Note, however, that there will be field of moduli issues when one works
at level one, so even the rationality of this space over~$\Q$ does not imply the existence of such~$A$. 
 
 \subsubsection{Cremona's Question}
 We finally consider two~$2$-dimensional irreducible compatible
 systems~$\CS$ of representations of~$G_K$ for some quadratic extension~$K/\Q$,  with inverse cyclotomic
 determinant, Hodge--Tate
 weights~$(0,1)$,
 and coefficients in~$\Q$.
 Note that, for such a family~$\CS$, there is a corresponding family~$\CR = \Ind^{G_{\Q}}_{G_K}\cS$ of~$4$-dimensional
 symplectic representations with inverse cyclotomic similitude character.
 An argument very similar to (but easier than)
 Lemma~\ref{lemma:motives} shows that (assuming all conjectures) either~$\CS$ comes from an elliptic curve, or it arises from a so-called fake elliptic curve,
 namely, an abelian surface~$B/K$ with~$\End_{K}(B) \otimes \Q = D$ for some indefinite quaternion algebra~$D$. The latter can exist only when~$K$ is an imaginary
 quadratic field.  Conjecturally, such compatible systems are in
 bijection with cuspidal cohomological~$\boldsymbol{\pi}$
 for~$\GL(2)/K$ with trivial central character and Hecke eigenvalues
 in~$\Q$.

 One source of such~$B/K$ is to take abelian surfaces over~$\Q$ of~$\GL(2)$-type which acquire quaternionic multiplication over~$K/\Q$.
Assuming the Hodge conjecture and the standard conjectures, it follows
that~\cite[Question
$1$\kern-0.05em{\footnotesize{${}'$}}]{CremonaTwist}  (cf.\ \cite[Conjecture, p.278]{CremonaOld}) is equivalent to asking that~\emph{all} fake elliptic curves~$B$ over~$K$ 
descend to~$\Q$ after twisting by some quadratic character (equivalently, $B$ is  isogenous to a twist of~$B^{\sigma}$  for
the non-trivial element~$\sigma \in \Gal(K/\Q)$). We call such~$B$ \emph{non-autochthonous} because it implies
that the corresponding conjectural~$\boldsymbol{\pi}$ arises via functoriality from a smaller rank group
(cf.\ footnote~$2$ of~\cite{AshTwisted}).
  In this section, we show that the answer to this question is false, namely, we construct autochthonous fake elliptic curves~$B/K$.
If one takes the restriction of scalars~$A = \Res_{K/\Q}(B)$ of such surfaces, then the fourfolds~$A$ give rise to further examples in opposition to the original paramodular conjecture.

 We continue to let~$D$ be the quaternion algebra ramified at
 precisely~$2$ and~$3$. The Shimura curve giving rise to fake elliptic curves with endomorphisms by a maximal order in~$\OL_D$ has genus zero, and is
 well-known (see for example~\cite[Thm.\ 11]{Baba}) to be isomorphic over~$\Q$
 to:
 $$X^2 + Y^2 + 3 Z^2 = 0.$$
 Moreover, more usefully for our purposes,
   Baba and Granath in~\cite{Baba} give  explicit models for genus two curves with endomorphisms by~$\OL_D$. For a certain parameter~$j$,
 they write  down a model (\cite[Thm.\ 15]{Baba}) of a genus two curve~$C$ over~$\Q(\sqrt{-6j})$ such that its endomorphisms
 are all defined (by~\cite[Prop.\ 19]{Baba}) over~$K:=\Q(\sqrt{-6j},\sqrt{j},\sqrt{-27(j+16)})$. With a view to choosing~$K = \Q(\sqrt{-6})$, we let
 $$Z = 3 \sqrt{j}, \qquad  X = \sqrt{-(27j + 16)},  \qquad Y = 4,$$
 and look for solutions to the equation above with~$X,Z \in \Q(\sqrt{-6})$. One such solution is given by~$j = -32/27$, but the corresponding
 surface is not autochthonous. Thus, we parametrize the conic and choose a random such point. Without making too much effort to optimize the height of~$j$,
 one finds that
 \begin{equation} \label{equation:eqrefj}
j = \frac{4 (1 + 2 \sqrt{-6})}{27}
 \end{equation}
 is a suitable point.

 \begin{lem} \label{lemma:explicitgenustwo} Let~$C$ be the following genus two curve:
 $$ y^2 = a_0 + a_1 x + a_2 x^2 + a_3 x^3 + a_4 x^4 + a_5 x^5 + a_6 x^6,$$
 where~$a_i$ are given by the following table, where~$\eta = 1 - \sqrt{-6}$:
 \begin{center}
 \begin{tabular}{|c|l|}
 \hline
 $i$  & $a_i$ \\
 \hline
 $0$  & $ 2^{10} \cdot \eta^{6} \cdot (\sqrt{-6}-4)   $ \\
$1$ &   $2^{10} \cdot \eta^{6}  \cdot 3 $ \\
$2$ &  $2^7 \cdot \eta^4 \cdot (9 \sqrt{-6} + 24)$ \\
$3$ & $- 2^8 \cdot \eta^4 $ \\
$4$ & $2^4 \cdot \eta^2 \cdot  (-9 \sqrt{-6} + 60)$ \\
$5$ & $2^4 \cdot \eta \cdot 3 $ \\
$6$ & $ - 2 \cdot \sqrt{-6}$ \\
  \hline
 \end{tabular}
 \end{center}
The sextic  has discriminant~$2^{90} \cdot 3^6 \cdot (1 - \sqrt{-6})^{30} \cdot (2 - \sqrt{-6})^6$.
 Let~$B = \Jac(C)/\Q(\sqrt{-6})$, and let~$A = \Res_{\Q(\sqrt{-6})/\Q}(B)$. Then~$\End_{\Q}(A) \otimes \Q = D$ 
 and~$\End_{\C}(A) \otimes \Q = M_2(D)$, where~$D/\Q$ is the
 quaternion algebra ramified at precisely~$2$ and~$3$.
 Then~$B$ is autochthonous, and $B$ gives rise to an irreducible~$2$-dimensional compatible system of Galois representations~$\CS$  of~$G_{\Q}$
 with Hodge--Tate weights~$(0,1)$ and inverse cyclotomic determinant,
 and~$A$ gives rise to a~$4$-dimensional compatible system~$\CR$
 of~$4$-dimensional $p$-adic Galois representations of~$G_{\Q}$ with coefficients in~$\Q$
 unramified outside of~$\{2,3,5,7,p\}$, each of which is
 absolutely irreducible and symplectic with inverse cyclotomic
 multiplier.
 \end{lem}

\begin{proof} Let~$K = \Q(\sqrt{-6})$. The curve~$C$ is the
  specialization of the curve in~\cite[Thm.\ 15]{Baba} to the parameter~$j$
as in equation~\ref{equation:eqrefj}, and
$$s =  \sqrt{-6j} =  \frac{-2 \sqrt{-6} + 4}{3}, \qquad
t =  -2 (27 j + 16) = 
 -16 \sqrt{-6} - 40.$$
By~\cite[Prop.\ 19]{Baba}, we deduce that the endomorphisms of~$C$ are defined 
over the field~$K(\sqrt{j},\sqrt{-(27j+16)}) = K$. We now show that~$B$ is autochthonous.
Let~$\p = (11,\sqrt{-6} - 4)$ and~$\q = (11,\sqrt{-6} + 4)$. Then the curves~$X_1 = C(\OL_K/\p)$ and~$X_2 = C(\OL_K/\q)$ 
over~$\F_{11}$ are given explicitly
as follows:
$$X_1: y^2 = 3 x^6 + 3 x^5 + 2 x^4 + 10 x^3 + 8 x^2 + 9 x,$$
$$X_2: y^2 = 8 x^6 + x^5 + 10 x^4 + 6 x^3 + 3 x^2 + 4 x + 4.$$
We compute the zeta functions using {\texttt{magma}} (see~\cite{magma}) to be as follows:
$$Z(X_1,s) = \frac{ (1 -2  p^{-s}+  p^{1-2s})^2}{(1-p^{-s})(1-p^{1-s})},  \quad Z(X_2,s) =  \frac{ (1 - p^{-s} +  p^{1-2s})^2}{(1-p^{-s})(1-p^{1-s})}.$$	
If~$B$ were autochthonous, then in particular the zeta functions of $X_1,X_2$ would differ by a twist by a
finite order character,
but this is impossible since~$2 \ne \pm 1$.

If~$p>3$, then~$p$ splits in~$D$, so there are  
Galois representations
$$r_p: G_K \rightarrow \GL_2(\Q_p)$$
with~$V_p(B)^\vee\cong r_p\oplus r_p$. For all~$p$, there also exist corresponding representations~$\rrr_p: G_K \rightarrow \GL_2(\Qbar_p)$
such that, for~$p > 3$, the representation~$\rrr_p$ is the representation obtained from~$r_p$ by extending scalars.
 We now prove that~$\End_{\C}(B) \otimes \Q = D$. If this were not true, then~$B$ would geometrically have to be isogenous to~$E \times E$ for some elliptic curve~$E$
 with complex multiplication.  This implies that~$B$ itself has complex multiplication over~$\C$,
 which implies that the representations~$r_p$ are potentially reducible. But as the representations~$r_p$ have distinct Hodge--Tate weights,
if they become reducible they do so over a quadratic extension.  This quadratic extension~$L/K$ must be ramified only
 at primes of bad reduction of~$B$, and for~$\p$ which are inert
 in~$L/K$, one must have~$a_{\p} = 0$. But this can be ruled out by
 computation (the only prime~$\p$ of norm less than~$1000$
 with~$a_{\p} = 0$ has norm~$97$).

   Hence~$\End_{\C}(B) \otimes \Q =D$ and~$\End_{\Q}(A) \otimes \Q
   = D$, where~$A = \Res_{K/\Q}(B)$.
It also follows that the representations~$r_p$  and~$\rrr_p$
have inverse cyclotomic
determinant (as otherwise they would be isomorphic to their twists by
a finite order character, and thus potentially reducible). Moreover,  with~$\rho_p:=\Ind^{\Q}_{K} \rrr_p$,
one has
 $$\wedge^2 \rho_p= \Asai(\rrr_p) \oplus \varepsilon^{-1} \oplus \varepsilon^{-1} \cdot \eta_{K/\Q},
 \quad \Sym^2 \Ind^{G_{\Q}}_{G_K} \rho_p = \Ind^{G_{\Q}}_{G_K} \Sym^2 \rrr_p,$$
 and thus~$\rho_p$ is absolutely irreducible and can be chosen to have
 image in~$\GSp_4(\Qbar_p)$ with inverse cyclotomic
 similitude character.
 Finally, the characteristic polynomials of Frobenius will,
 by construction, be degree~$4$ polynomials with coefficients in~$\Q$. \end{proof}

Since one expects the compatible system~$\CS$ to be modular (it is certainly potentially modular, by~\cite{10author}),
it follows that
Cremona's question~\cite[Question $1$\kern-0.05em{\footnotesize{${}'$}}]{CremonaTwist})
is incompatible with standard modularity conjectures. (Similarly, the modularity
of~$\CR = 
\Ind^{G_{\Q}}_{G_{K}} \CS$ is incompatible with the original 
paramodular conjecture, although we have already shown the latter to be false.)
 Of course, from the discussion above, there are natural modifications
 that one could make to Cremona's question (along the lines of
 Conjecture~\ref{conj:modified}) --- namely,
to include all  fake elliptic curves over~$K$, autochthonous or otherwise.

One can presumably show that the~$2$-dimensional~$G_K$-representations~$\rrr_p$
over the field~$K = \Q(\sqrt{-6})$
arising from~$B = \Jac(C)$ are modular
for~$\GL(2)/K$. As in the proof of
Lemma~\ref{lem: our things are paramodular}, we would then obtain a
 cuspidal cohomological automorphic representation~$\boldsymbol{\pi}$ for~$\GL(2)/K$ with trivial central character and rational eigenvalues.
 Since~$\boldsymbol{\pi}$ does not arise (up to twist) from base
 change, this would answer in the negative~\cite[Question $1$\kern-0.05em{\footnotesize{${}'$}}]{CremonaTwist},
 because  the existence of a corresponding elliptic curve~$E/K$ would be
 incompatible with the existence of~$B$ by Faltings' isogeny
 theorem~\cite{Faltings}.

The modularity of the representations~$\rrr_p$ can in principle be established using the Faltings--Serre method (cf.~\cite{BergerBianchi}).  Possibly some computational advantage
would be gained by replacing~$C$ with a curve obtained from a more careful choice of generic point on the Shimura curve 
(in order to work at a manageable level). 
As it turns out, Ciaran Schembri~\cite{Ciaran} has independently found examples of  autochthonous fake elliptic curves which he has verified are modular.

\bibliographystyle{amsalpha}
\bibliography{abeliansurfacesmodular}

\end{document}